\documentclass[11 pt]{article}
\usepackage[english]{babel}
\usepackage[utf8]{inputenc}
\usepackage{amsmath}
\usepackage{amsfonts}
\usepackage{mathrsfs}
\usepackage{enumitem}
\usepackage{graphicx}
\usepackage{url}
\usepackage[colorinlistoftodos]{todonotes}
\usepackage{blindtext}
\usepackage{amssymb}
\usepackage{algorithm}
\usepackage{bbm}
\usepackage{amsthm}
\usepackage{euscript}
\usepackage[margin=1in]{geometry}
\usepackage[pdfborder={0 0 0}]{hyperref}
\hypersetup{
  colorlinks=true, linkcolor=blue, citecolor=blue, urlcolor=blue,
  pdfauthor = {Chenyang Zhong, Marcel Nutz},
  pdfkeywords = {Optimal transport, distance cost, entropic regularization, vanishing regularization},
  pdftitle = {Entropic regularization of Monge’s problem},
  pdfsubject = {Entropic regularization of Monge’s problem},
  pdfpagemode = UseNone
}
\makeatletter
\newcommand*{\rom}[1]{\expandafter\@slowromancap\romannumeral #1@}
\usepackage{tikz}
\usetikzlibrary{calc,arrows.meta,decorations.pathreplacing,positioning,intersections}

\usepackage[capitalize,noabbrev]{cleveref}

\newtheorem{theorem}{Theorem}[section]
\newtheorem{claim}{Claim}
\newtheorem{proposition}[theorem]{Proposition}
\newtheorem{lemma}[theorem]{Lemma}
\newtheorem{corollary}[theorem]{Corollary}
\theoremstyle{definition}
\newtheorem{definition}[theorem]{Definition}
\newtheorem{remark}[theorem]{Remark}

\newtheorem{assumption}[theorem]{Assumption}

\AddToHook{env/theorem/begin}{\crefalias{section}{theorem}}
\AddToHook{env/proposition/begin}{\crefalias{theorem}{proposition}}
\AddToHook{env/lemma/begin}{\crefalias{theorem}{lemma}}
\AddToHook{env/corollary/begin}{\crefalias{theorem}{corollary}}
\AddToHook{env/definition/begin}{\crefalias{theorem}{definition}}
\AddToHook{env/remark/begin}{\crefalias{theorem}{remark}}
\AddToHook{env/example/begin}{\crefalias{theorem}{example}}
\AddToHook{env/assumption/begin}{\crefalias{theorem}{assumption}}
\crefname{theorem}{Theorem}{Theorems}
\crefname{claim}{Claim}{Claims}
\crefname{proposition}{Proposition}{Propositions}
\crefname{lemma}{Lemma}{Lemmas}
\crefname{corollary}{Corollary}{Corollaries}
\crefname{definition}{Definition}{Definitions}
\crefname{remark}{Remark}{Remarks}
\crefname{example}{Example}{Examples}
\crefname{assumption}{Assumption}{Assumptions}

\newcommand{\R}{\mathbb{R}}

\newcommand{\eps}{\varepsilon}
\newcommand{\dist}{d}

\DeclareMathOperator{\OT}{\mathsf{OT}}
\DeclareMathOperator{\EOT}{\mathsf{EOT}}
\DeclareMathOperator{\approxgrad}{\mathrm{ap}\hspace{-1.5pt}\nabla}

\numberwithin{equation}{section}

\begin{document}

\title{Entropic regularization of Monge's problem}

\date{\today}

\author{  
  Marcel Nutz%
  \thanks{Departments of Mathematics and Statistics, Columbia University, mnutz@columbia.edu. Research supported by NSF Grants DMS-2106056, DMS-2407074.} \and Chenyang Zhong%
  \thanks{Department of Statistics, Columbia University, cz2755@columbia.edu.}
  }

\maketitle
\vspace{-1.5em}

\begin{abstract}
We study the vanishing-regularization limit of entropically regularized optimal transport (EOT) for the Euclidean distance cost $c(x,y)=\|x-y\|$ in dimension $d>1$. We develop a comprehensive variational convergence framework that entails two main results. First, we resolve the longstanding entropic selection problem: the EOT minimizer converges to a distinguished optimal transport plan that is characterized explicitly as the solution of a constrained EOT problem on each transport ray. Denoting by $\eps>0$ the regularization parameter, this selection holds for all $o(\eps)$-approximate minimizers, with sharp failure at the $O(\eps)$ scale. Second, we establish an explicit second-order expansion of the entropic transport cost. The second-order term encodes the geometry of the regularization and reveals the optimal asymptotic tradeoff between entropy and transport cost.
\end{abstract}

\tableofcontents
\section{Introduction}\label{Sect.1}
Given two compactly supported Borel probability measures $\mu,\nu$ on $\mathbb{R}^d$ and a continuous function $c:\mathbb{R}^d\times\mathbb{R}^d\to\mathbb{R}$, the Monge--Kantorovich optimal transport (OT) problem is
\begin{equation}\label{OT_Min}
    \OT(\mu,\nu):=\inf_{\gamma \in \Pi(\mu,\nu)} \int_{\mathbb{R}^d\times \mathbb{R}^d}  c(x,y)  d\gamma (x,y),
\end{equation}
where $\Pi(\mu,\nu)$ denotes the set of all  \emph{transport plans} or \emph{couplings} between $\mu$ and $\nu$, that is, all Borel measures on $\mathbb{R}^d\times\mathbb{R}^d$ with marginals $(\mu,\nu)$. 
Any optimizer $\gamma\in \Pi(\mu,\nu)$ of \eqref{OT_Min} is called an \emph{optimal transport plan}. For certain cases, for instance the quadratic cost $c(x,y)=\|x-y\|^2$ and absolutely continuous marginals $(\mu,\nu)$, the optimizer is unique, but in general, the linear nature of~\eqref{OT_Min} entails that there can be many optimal transport plans. In the continuous setting, the most important example with non-unique optimizers is the Euclidean distance cost $c(x,y)=\|x-y\|$, also called the \emph{Monge cost} as it was used in Monge's original formulation \cite{monge1781memoire} of optimal transport. The corresponding value of $\OT(\mu,\nu)$ is known as the 1-Wasserstein distance (or Earthmover's distance, or Kantorovich-Rubinstein distance) between~$\mu$ and~$\nu$. See, e.g., \cite{Santambrogio.15,Villani} for background on optimal transport.

Given a regularization parameter $\eps>0$, the \emph{entropic optimal transport} (EOT) problem adds an entropic penalty to the transport objective and solves
\begin{equation}\label{EOT_Min}
    \EOT_\eps(\mu,\nu):=\inf_{\gamma \in \Pi(\mu,\nu)} \mathcal{C}_{\eps}(\gamma), \qquad
    \mathcal{C}_{\eps}(\gamma):=\int_{\mathbb{R}^d\times \mathbb{R}^d} c(x,y) d\gamma(x,y)+\eps H(\gamma|\mu\otimes \nu).
\end{equation}
Here $\mu\otimes \nu$ is the product of $\mu$ and $\nu$ while $H(P|Q)$ denotes the relative entropy (or Kullback--Leibler divergence) between two probability measures $P$ and $Q$ on the same space, 
\begin{equation}
    H(P|Q):=\begin{cases}
            \int \log\big(\frac{dP}{dQ}\big)dP, & P\ll Q,\\
            \infty, & P\not\ll Q.
        \end{cases}
\end{equation}
Thanks to the strict convexity of the penalty, the EOT problem \eqref{EOT_Min} admits a unique minimizer $\gamma_{\eps}\in \Pi(\mu,\nu)$ for each $\eps>0$, whereas $\eps=0$ recovers the OT problem~\eqref{OT_Min}.

EOT traces back to the Schr\"odinger bridge problem~\cite{Schrodinger.31} on the most likely evolution of a particle system (see~\cite{Follmer.88, Leonard.14} for surveys). However, the recent explosion of the EOT literature mainly stems from its utility as an approximation of OT (and hence of Wasserstein distances) with convenient properties, including fast computation via Sinkhorn's algorithm (e.g.,~\cite{peyre2019computational}) and parametric statistical sample complexity (e.g.,~\cite{ChewiNilesWeedRigollet.25}), among others. 

Therefore, the limit $\eps\to0^+$ of vanishing regularization is of immediate interest. When the OT problem has a unique solution, weak compactness of $\Pi(\mu,
\nu)$ yields that the EOT optimizers~$\gamma_{\eps}$ converge weakly to that optimal transport plan. When the OT problem has more than one solution, it merely follows that any cluster point of $(\gamma_{\eps})$ must be an optimal transport plan, hence one may ask if $\gamma_{\eps}$ converges to some particular plan~$\gamma_0$. From the construction of $\gamma_{\eps}$, such a limit may carry the intuition of having small relative entropy. Indeed, for discrete marginals $(\mu,\nu)$, it is well known that $\gamma_{\eps}$ converges to the optimal transport plan with minimal relative entropy~\cite{CominettiSanMartin.94}, and that result remains valid for general marginals as long as there exists \emph{some} optimal transport plan with finite relative entropy \cite[Theorem~5.5]{Nutz.20}. On the other hand, when the marginals are continuous and the dimension is $d>1$, optimal plans are typically concentrated on lower-dimensional sets. Then, optimal plans are not absolutely continuous wrt.\ $\mu\otimes\nu$ and hence have infinite relative entropy. As a result, it is unclear how to make sense of the aforementioned intuition. Nevertheless, one may hope that $\gamma_{\eps}$ converges to some particular plan~$\gamma_0$. To the best of our knowledge, the first written account of that conjecture appeared in \cite{Leonard.12} where the potential limit is called ``a viscosity solution'' of~\eqref{OT_Min}; however, the conjecture is likely older. From an applied point of view, the relevant question may require an even stronger theorem: Suppose that we have \emph{approximate} EOT solutions $\gamma'_{\eps_n}$ for some sequence $\eps_n\to0^+$, as one would obtain by sampling i.i.d.\ from $(\mu,\nu)$ and running Sinkhorn's algorithm on the resulting discrete EOT problem for finitely many iterations. Then, we would like to assert that $\gamma'_{\eps_n}$ converges to that particular optimal transport plan $\gamma_0$, if the number of samples and iterations is sufficiently large. This problem falls in the broader category of selection problems in optimal transport, which often must be solved on a case-by-case basis (e.g.\ \cite[Section~3]{Santambrogio.15}).

As the Monge cost $c(x,y)=\|x-y\|$ is by far the most important example with non-unique optimal transport plans, it is the canonical setting for our investigation. We also assume throughout that the dimension is $d>1$. The case $d=1$ is special as optimal transports can be absolutely continuous; this case was solved in~\cite{di2018entropic} and, as we will see, is fundamentally different from the higher-dimensional case ($d=1$ is more similar to the aforementioned case of discrete marginals). 

In this work, we develop a variational (Gamma-type) convergence framework for the rescaled EOT cost. This functional convergence, summarized in \cref{th:upperBound,th:lowerBound}, leads to a detailed description of the limit of EOT towards OT, for which we state two main results. The first (Theorem~\ref{Thm1.1}) shows that $o(\eps)$-approximate optimizers $\gamma'_\eps$ of the EOT problem converge to a particular optimal transport plan $\gamma_0$ which we call the \emph{entropic Monge plan} (Definition \ref{Definition1.1}). This plan admits a fairly explicit local description as the solution of a \emph{constrained} EOT problem on each transport ray, with cost function $\tilde{c}(x,y)=-\frac{d-1}{2}\log(\|x-y\|)$. The error tolerance $o(\eps)$ is sharp in that the convergence to~$\gamma_0$ can fail at scale $O(\eps)$; cf.\ \cref{Exa1.1}. The second result (Theorem~\ref{Thm1.2}) is an explicit second-order expansion of the optimal entropic transport cost $\EOT_\eps(\mu,\nu)$ in~$\eps$. While the first-order term is the same for all non-overlapping marginals, the second-order term reflects both the Monge cost and the marginals. It can be interpreted, in analogy with the renormalized energy in the Ginzburg--Landau theory \cite{MR2945619}, as a variational quantity that captures the fine structure of the transport problem.

Altogether, our variational framework provides a comprehensive picture of the vanishing regularization limit of EOT towards Monge's problem. The local description of the entropic Monge plan and the second-order coefficient in the expansion of $\EOT_\eps(\mu,\nu)$ are closely related. Both reflect a tradeoff between the two terms in the EOT problem~\eqref{EOT_Min}: transport cost versus entropy cost, i.e., how far mass is transported versus how much relative entropy it causes. While transporting along transport rays is optimal for the unregularized transport problem, that would cause infinite entropy. The EOT optimizer $\gamma_\eps$ thus spreads mass away from the rays, and the expressions in our main results reflect the optimal tradeoff for small $\eps>0$. Our results show that the optimal tradeoff depends on the precise details of the transport cost. This contrasts sharply with the existing results for discrete marginals (and the continuous Monge problem in $d=1$), where the limiting optimal transport plan is simply the one with minimal relative entropy, regardless of the transport cost. To wit, in those problems, the limit is the result of a two-stage optimization: find the minimizers of the transport cost, then minimize entropy over that set. In the present setting, such a separation does not happen. This insight may also suggest that a generic resolution of the entropic selection conjecture for general costs~$c$ is difficult; arguably, it retroactively justifies our focus on the distance cost.

\subsection{Regularity assumptions and notation regarding transport rays}\label{Sect.1.1}

Let $d\geq 2$ and let $\mathcal{L}^d$ denote Lebesgue measure on $\mathbb{R}^d$. We work with the Euclidean distance cost $c(x,y)=\|x-y\|$ and adopt the conventions $0\log 0:=0$ and $\log 0:=-\infty$. Moreover, we impose the following throughout the paper.

\begin{assumption}\label{assump1}
The marginals $\mu$ and $\nu$ are Borel probability measures on $\mathbb{R}^d$ that are absolutely continuous with respect to $\mathcal{L}^d$, with Lipschitz continuous densities $f$ and $g$. Their supports $\mathcal{X}:=\mathrm{supp}(\mu)$ and $\mathcal{Y}:=\mathrm{supp}(\nu)$ are compact and disjoint. 
\end{assumption}

See \cref{rem:overlap} for the case of overlapping marginal supports. We fix a \textit{Kantorovich potential} $u\in \mathrm{Lip}_1(\mathbb{R}^d)$, the set of $1$-Lipschitz functions on $\mathbb{R}^d$. That is, $u$  is a maximizer of the dual OT problem 
\cite[Remark 6.5]{Villani},
\begin{equation}\label{OTrela}
\sup_{u\in\mathrm{Lip}_1(\mathbb{R}^d)}\int_{\mathbb{R}^d} u(x)(f(x)-g(x))dx =\OT(\mu,\nu).
\end{equation}
The formulation of our main results relies on the geometry of Monge's problem which can be described in terms of \textit{transport rays}. Transport rays are maximal line segments in $\mathbb{R}^d$ along which the Kantorovich potential $u$ decreases at rate one, and any two such rays have disjoint relative interiors. A detailed account of transport rays and related notions is given in Section \ref{Sect.2.1}. Roughly speaking, the connection between optimality and geometry is this: a measure $\gamma\in\Pi(\mu,\nu)$ is an optimal transport plan if and only if it transports mass only along transport rays and, on each ray, only in the direction where~$u$ decreases. This is made precise using the disintegration theorem (e.g., \cite[Theorem~9.1]{MR2011032}) as follows. Denote by $\mathcal{S}$ the set of transport rays, equipped with the natural metric introduced in Definition \ref{TransportRays}. 
Ambrosio \cite{MR2011032} showed that the marginals can be disintegrated along transport rays. Namely, there exist a Borel probability measure $\lambda$ on $\mathcal{S}$, and, for each $T\in\mathcal{S}$, a finite (positive) Radon measure $\tilde{\mu}_T$ on $\mathbb{R}^d$ that is concentrated on $T$, such that
\begin{itemize}
    \item[(a)] for $\lambda$-a.e. $T\in\mathcal{S}$, $\tilde{\mu}_T$ is a probability measure;
    \item[(b)] for any Borel set $A\subseteq\mathbb{R}^d$, $\mu(A) = \int_{\mathcal{S}}\tilde{\mu}_T(A)d\lambda(T)$; 
     \item[(c)] for any Borel set $A\subseteq\mathbb{R}^d$, the mapping  $T\mapsto \tilde{\mu}_T(A)$ is Borel measurable on $\mathcal{S}$.
\end{itemize}
We write this disintegration succinctly as $\mu=\lambda\otimes \tilde{\mu}_T$. Analogously, $\nu=\lambda\otimes\tilde{\nu}_T$ for the second marginal. It is further shown, in the proof of \cite[Theorem 6.2]{MR2011032}, that any optimal transport plan $\gamma\in \Pi(\mu,\nu)$ admits a disintegration along transport rays. Namely, for each $T\in\mathcal{S}$, there exists a finite (positive) Radon measure $\omega_T$ on $\mathbb{R}^d\times\mathbb{R}^d$ such that $\gamma=\lambda\otimes \omega_T$ (meaning that the analogues of (a)--(c) hold), and in addition,
\begin{itemize}
    \item[(i)] $\omega_T\in\Pi(\tilde{\mu}_T,\tilde{\nu}_T)$ for $\lambda$-a.e.\ $T\in\mathcal{S}$;
    \item[(ii)] $\omega_T$ is \emph{monotone,} i.e.,  concentrated on $\{(x,y)\in T\times T: u(x)\geq u(y)\}$, for all $T\in\mathcal{S}$.
\end{itemize}

Returning to the disintegrations of the marginals, Ambrosio \cite[Theorem 6.1]{MR2011032} further showed that $\tilde{\mu}_T$ and $\tilde{\nu}_T$ are absolutely continuous with respect to the Hausdorff measure $\mathcal{H}^1|_T$; we denote their densities on $T$ by $\tilde{f}_T$ and $\tilde{g}_T$, respectively. We can then construct Borel functions $\tilde{f}$ and $\tilde{g}$ on~$\R^d$ such that $\tilde{f}(x)=\tilde{f}_T(x)$ and $\tilde{g}(x)=\tilde{g}_T(x)$ whenever $x$ lies in the relative interior of a transport ray~$T$. While these abstract results were sufficient for the purposes of \cite{MR2011032}---to construct an optimal transport map for Monge's problem---the present study requires a finer analysis. To that end, we derive explicit formulas for $\lambda,\tilde{\mu}_T,\tilde{\nu}_T,\tilde{f},\tilde{g}$ in Section \ref{Sec3} (see Definition~\ref{De3.8} and Lemma~\ref{L3nnn}).

\subsection{Main results}\label{Sect.1.2}

With the assumptions and notation from Section \ref{Sect.1.1} in place, we can now state our main results. We begin by introducing the limiting optimal transport plan for our selection result. Recall the disintegrations $\mu=\lambda\otimes \tilde{\mu}_T$ and $\nu=\lambda\otimes\tilde{\nu}_T$ from Section \ref{Sect.1.1}.

\begin{definition}[Entropic Monge plan $\gamma_0$]\label{Definition1.1}
 For each transport ray $T$, consider the constrained EOT problem
\begin{equation}\label{opt.solver}
\inf_{\substack{\chi_T\in\Pi(\tilde{\mu}_T,\tilde{\nu}_T):\\ \chi_T\text{\;is monotone on $T$}}}\bigg\{-\frac{d-1}{2} \int_{T\times T} \log(\|x-y\|)d\chi_T(x,y) + H(\chi_T|\tilde{\mu}_T\otimes \tilde{\nu}_T)\bigg\}.
\end{equation}
We will show (in \Cref{Lem2}) that for $\lambda$-a.e.\ $T\in\mathcal{S}$, the minimization \eqref{opt.solver} is finite and admits a unique solution $\kappa_T\in\Pi(\tilde{\mu}_T,\tilde{\nu}_T)$. We set
\begin{equation}\label{gamma_0}
    \gamma_0 :=\lambda\otimes\kappa_T \in\Pi(\mu,\nu)
\end{equation}  
and call $\gamma_0$ the \emph{entropic Monge plan.}
\end{definition}

\begin{remark}\label{Rmk1.2}
The monotonicity constraint in \eqref{opt.solver} is crucial to ensure that~$\gamma_0$ is an optimal transport plan, and hence for the validity of the main results below (see \cref{sec:counterex_directional_constraint} for a simple counterexample). The constraint makes~\eqref{opt.solver} a \emph{nonstandard EOT problem}; we will develop its regularity properties from first principles.
\end{remark}

Next, we define the limiting functional for our variational analysis. In addition to the disintegrations $\mu=\lambda\otimes \tilde{\mu}_T$ and $\nu=\lambda\otimes\tilde{\nu}_T$ of the marginals, we recall that $(f,g)$ are the $\mathcal{L}^d$-densities of $(\mu,\nu)$, and $(\tilde {f},\tilde{g})$ are the aggregated $\mathcal{H}^1$-densities of $(\tilde{\mu}_T,\tilde{\nu}_T)$.

\begin{definition}\label{de:Sfunctional}
    Let $\gamma_0'\in\Pi(\mu,\nu)$ be an optimal transport plan and $\gamma_0'=\lambda\otimes\kappa'_T$ its disintegration along transport rays. We define 
    \begin{align}\label{SexpressGeneralPlan}
    \mathtt{S}(\mu,\nu;\gamma_0'):=\, & \int_{\mathcal{S}} \bigg(-\frac{d-1}{2}\int_{T\times T}\log(2\pi\|x-y\|) d\kappa'_T(x,y)+H(\kappa'_T|\tilde{\mu}_T\otimes\tilde{\nu}_T)\bigg)d\lambda(T)\nonumber\\
    &-\frac{1}{2}\int_{\mathcal{X}}f(x)\log(f(x))dx-\frac{1}{2}\int_{\mathcal{Y}}g(y)\log(g(y))dy\nonumber\\
    &+\frac{1}{2}\int_{\mathcal{S}}d\lambda(T)\int_T\tilde{f}(x)\log(\tilde{f}(x))d\mathcal{H}^1(x)+\frac{1}{2}\int_{\mathcal{S}}d\lambda(T)\int_T\tilde{g}(y)\log(\tilde{g}(y))d\mathcal{H}^1(y).
\end{align}
\end{definition}

\begin{lemma}
The right-hand side of \eqref{SexpressGeneralPlan} is well-defined, $\mathtt{S}(\mu,\nu;\gamma_0')\in (-\infty,\infty]$ for any optimal transport plan $\gamma_0'$, and $\mathtt{S}(\mu,\nu;\gamma_0)<\infty$.
\end{lemma}
\begin{proof}

The assertions follow from \Cref{L3.18,L4.11}.
\end{proof}

The next two theorems summarize our core technical results about the convergence of the rescaled EOT functional. Weak convergence of probability measures is understood in the probabilistic sense, i.e., relative to bounded continuous test functions.

\begin{theorem}\label{th:lowerBound}
Consider $\varepsilon_j\rightarrow 0^{+}$ and arbitrary couplings  $\gamma'_{\eps_j}\in\Pi(\mu,\nu)$ converging weakly to an optimal transport plan $\gamma_0'\in\Pi(\mu,\nu)$. Then 
\begin{align*}
   & \liminf_{j\rightarrow\infty}\bigg\{\frac{\mathcal{C}_{\eps_j}(\gamma'_{\eps_j})-\OT(\mu,\nu)}{\eps_j}+\frac{d-1}{2}\log{\eps_j}\bigg\}\nonumber
    \geq \mathtt{S}(\mu,\nu;\gamma_0').
\end{align*}
\end{theorem}

\begin{theorem}\label{th:upperBound}
Let $\gamma_\eps\in\Pi(\mu,\nu)$ be the optimizer of the EOT problem~\eqref{EOT_Min} and $\gamma_0\in\Pi(\mu,\nu)$ the entropic Monge plan. Then
\begin{equation*}
    \limsup_{\eps\rightarrow 0^+}\bigg\{\frac{\mathcal{C}_{\eps}(\gamma_{\eps})-\OT(\mu,\nu)}{\eps}+\frac{d-1}{2}\log{\eps}\bigg\}\leq \mathtt{S}(\mu,\nu;\gamma_0).
\end{equation*}
\end{theorem}

\begin{remark}
    The lower bound (Theorem~\ref{th:lowerBound}) follows the usual definition of Gamma convergence. The upper bound (Theorem~\ref{th:upperBound}) is stated only for the sequence $\gamma_\eps$ of optimizers, which is sufficient for the consequences that we want to draw. In the usual definition of Gamma convergence, the upper bound is stated for a recovery sequence of any point. Our proof uses delicate regularity properties of~$\gamma_0$; it is not clear if an analogous upper bound can be shown for other limiting points.
\end{remark}

Combining \cref{th:upperBound,th:lowerBound} yields the explicit second-order expansion of the entropic transport cost $\EOT_\eps(\mu,\nu)$: the second-order coefficient is $\mathtt{S}(\mu,\nu; \gamma_0)$ of Definition~\ref{de:Sfunctional}, where $\gamma_0$ is the entropic Monge plan. In particular, the second-order coefficient encodes information about the Monge problem and (anticipating on \cref{Cor1.1} below) even the precise limit of the EOT optimizers~$\gamma_\eps$---in contrast to the leading-order coefficient~$\frac{d-1}{2}$ which is independent of $\mu,\nu$.

\begin{theorem}[Second-order expansion of the EOT cost]\label{Thm1.2}
We have
\begin{equation*}
    \EOT_\eps(\mu,\nu)=\OT(\mu,\nu)-\frac{d-1}{2}\eps\log\eps+\mathtt{S}(\mu,\nu; \gamma_0)\eps+o(\eps) \quad \text{as } \varepsilon\rightarrow 0^{+}.
\end{equation*}
\end{theorem}

Next, we state our main result about convergence of the EOT couplings. It shows that $o(\eps)$-approximate minimizers of the EOT problem \eqref{EOT_Min} converge weakly to the entropic Monge plan~$\gamma_0$. More generally, the statement also allows for an approximation error in the marginals, to accommodate the practical scenario where EOT is computed based on statistical samples from the marginals (see~\cite{FournierGuillin.15}). 

\begin{theorem}[Stable entropic selection]\label{Thm1.1}
For each $\eps>0$, consider probability measures~$\mu_\eps$ and~$\nu_\eps$ on $\R^d$ with 1-Wasserstein distance satisfying $W_1(\mu_\eps,\mu)+W_1(\nu_\eps,\nu)=o(\eps)$. Let $\gamma'_\eps\in\Pi(\mu_\eps,\nu_\eps)$ be $o(\eps)$-approximate minimizers of the associated EOT problems, i.e.,
$$
  \int_{\mathbb{R}^d\times \mathbb{R}^d} \|x-y\| d\gamma'_\eps(x,y)+\eps H(\gamma'_\eps|\mu_\eps\otimes \nu_\eps) \leq \EOT_\eps(\mu_\eps,\nu_\eps) + o(\eps).
$$
Then $\gamma'_{\eps}$ converges weakly to the entropic Monge plan $\gamma_0\in\Pi(\mu,\nu)$ as $\eps\rightarrow 0^+$. 
\end{theorem}

\begin{remark}[Sharpness of the $o(\eps)$ error tolerance]\label{Exa1.1}
The error tolerance $o(\eps)$ for the approximate optimizers in Theorem~\ref{Thm1.1} is sharp, even if $\mu_\eps=\mu$ and~$\nu_\eps=\nu$. Indeed, in \Cref{sec:error_tolerance},
we construct $\gamma'_\eps\in\Pi(\mu,\nu)$ which are $O(\eps)$-approximate minimizers for $\EOT_\eps(\mu,\nu)$ but converge to an optimal transport plan different from $\gamma_0$.
\end{remark}

Of course, Theorem~\ref{Thm1.1} shows in particular that the exact EOT minimizer $\gamma_{\eps}\in\Pi(\mu,\nu)$ converges to $\gamma_0$, thereby resolving the entropic selection conjecture for the Euclidean distance cost.

\begin{corollary}[Entropic selection]\label{Cor1.1}
The EOT minimizer $\gamma_{\varepsilon}$ converges weakly to $\gamma_0$ as $\varepsilon\rightarrow 0^{+}$.
\end{corollary}

\subsubsection{Deriving \cref{Thm1.1,Thm1.2} from the variational framework}

\Cref{Thm1.1,Thm1.2} are easily deduced from the variational framework in \cref{th:upperBound,th:lowerBound}.

\begin{proof}[Proofs of \cref{Thm1.1,Thm1.2}]
We first prove  \cref{Thm1.1} in the special case $(\mu_\eps,\nu_\eps)=(\mu,\nu)$. 
Consider a sequence $\eps_j\to0^+$ and $o(\eps_j)$-approximate EOT optimizers $\gamma'_{\eps_j}\in\Pi(\mu,\nu)$. As $\Pi(\mu,\nu)$ is weakly compact, a subsequence (still denoted $\eps_j$) converges weakly to some limit $\gamma_0'\in\Pi(\mu,\nu)$, and it suffices to show that $\gamma_0'=\gamma_0$.

As $\EOT_{\eps_j}(\mu,\nu)\to\OT(\mu,\nu)$, it follows that $\gamma_0'$ is an optimal transport plan. Consider its disintegration $\gamma_0'= \lambda \otimes \kappa_T'$ along transport rays, where $\kappa_T'\in\Pi(\tilde{\mu}_T,\tilde{\nu}_T)$. 
By Theorem~\ref{th:lowerBound},
\begin{align*}
   & \liminf_{j\rightarrow\infty}\bigg\{\frac{\mathcal{C}_{\eps_j}(\gamma'_{\eps_j})-\OT(\mu,\nu)}{\eps_j}+\frac{d-1}{2}\log{\eps_j}\bigg\}\nonumber
    \geq \mathtt{S}(\mu,\nu;\gamma_0').
\end{align*}
Conversely, noting that $\mathcal{C}_{\eps_j}(\gamma'_{\eps_j})\leq \EOT_{\eps_j}(\mu,\nu)+o(\eps_j)= \mathcal{C}_{\eps_j}(\gamma_{\eps_j})+o(\eps_j)$, Theorem~\ref{th:upperBound} yields
\begin{align*}
     & \limsup_{j\rightarrow\infty}\bigg\{\frac{\mathcal{C}_{\eps_j}(\gamma'_{\eps_j})-\OT(\mu,\nu)}{\eps_j}+\frac{d-1}{2}\log{\eps_j}\bigg\}
   \leq \mathtt{S}(\mu,\nu;\gamma_0).
\end{align*}
Combining the two displays gives $\mathtt{S}(\mu,\nu;\gamma_0)\geq \mathtt{S}(\mu,\nu;\gamma_0')$. Recall from Definition~\ref{de:Sfunctional} that
\begin{align*}
    \mathtt{S}(\mu,\nu;\gamma_0')=  \int_{\mathcal{S}} \bigg(-\frac{d-1}{2}\int_{T\times T}\log(2\pi\|x-y\|) d\kappa_T'(x,y)+H(\kappa'_T|\tilde{\mu}_T\otimes\tilde{\nu}_T)\bigg)d\lambda(T) +C,
\end{align*}
where $C$ is independent of $\gamma_0'$, and that $\kappa_T$ is the unique minimizer of the integrand on the right-hand side, for $\lambda$-a.e.\ $T\in\mathcal{S}$. Thus, $\mathtt{S}(\mu,\nu;\gamma_0)\geq \mathtt{S}(\mu,\nu;\gamma_0')$ implies that $\kappa'_T=\kappa_T$ for $\lambda$-a.e.\ $T\in\mathcal{S}$ and hence that $\gamma_0'=\gamma_0$. This completes the proof of Theorem~\ref{Thm1.1} for the special case $(\mu_\eps,\nu_\eps)=(\mu,\nu)$. Moreover, Theorem~\ref{Thm1.2} now follows immediately from Theorems~\ref{th:lowerBound} (taking $\gamma_0':=\gamma_0$) and \ref{th:upperBound}.

In the general case of Theorem~\ref{Thm1.1}, we are given $o(\eps)$-approximate minimizers $\gamma'_\eps\in\Pi(\mu_\eps,\nu_\eps)$ of $\EOT_\eps(\mu_\eps,\nu_\eps)$. It suffices to show that there are  $o(\eps)$-approximate minimizers $\gamma''_\eps\in\Pi(\mu,\nu)$ of $\EOT_\eps(\mu,\nu)$ with $W_1(\gamma'_\eps,\gamma''_\eps)=o(\eps)$, as the already established case then implies that $\gamma''_\eps\to\gamma_0$ and hence also $\gamma'_\eps\to\gamma_0$. 

Indeed, we can take $\gamma''_\eps$ to be the ``shadow'' of $\gamma'_\eps$ onto $\Pi(\mu,\nu)$ in the terminology of \cite{EcksteinNutz.21}. Roughly speaking, the shadow is obtained by concatenating three couplings: a $W_1$-optimal coupling of $\mu$ and $\mu_\eps$, then $\gamma'_\eps$, and finally a $W_1$-optimal coupling of $\nu_\eps$ and~$\nu$. By \cite[Lemma 3.2]{EcksteinNutz.21}, the shadow $\gamma''_\eps\in\Pi(\mu,\nu)$ satisfies  $$W_1(\gamma'_\eps,\gamma''_\eps)=W_1(\mu_\eps,\mu)+W_1(\nu_\eps,\nu) \quad\mbox{and}\quad  H(\gamma''_\eps|\mu\otimes\nu)\leq H(\gamma'_\eps|\mu_\eps\otimes\nu_\eps).$$
Moreover, $|\EOT_\eps(\mu_\eps,\nu_\eps)-\EOT_\eps(\mu,\nu)|\leq W_1(\mu_\eps,\mu)+W_1(\nu_\eps,\nu)$ by \cite[Theorem 3.7]{EcksteinNutz.21}. As $W_1(\mu_\eps,\mu)+W_1(\nu_\eps,\nu)=o(\eps)$, it follows that~$\gamma''_\eps$ has the required properties.
\end{proof}

\subsection{Proof strategy}\label{Sect.1.2.5}

In this subsection, we outline the proof strategy for the variational convergence framework presented in \cref{th:lowerBound,th:upperBound}. 

\subsubsection{Proof strategy for the lower bound}\label{se:lower}
\textbf{Regularity of Kantorovich potential and transport rays.} 
Proving the lower bound in Theorem~\ref{th:lowerBound} hinges on the fine regularity of the Kantorovich potential~$u$. Specifically, we need to estimate differences of the form $u(y)-u(x)$. If $x$ and $y$ lie on the same transport ray, this difference is explicit: $u(y)-u(x)=\pm\|y-x\|$, with the sign determined by their relative order along the ray. The main challenge is to obtain sharp control in directions \emph{orthogonal} to the ray, which requires strengthening the qualitative results of \cite{MR1862796,MR2011032} into quantitative, second-order estimates. Concretely, we develop second-order bounds for $u(x)-u(z+w)$ when $x,z$ lie in the interior of a common ray~$T$ and $w\in O_T$ with small $\|w\|$, where $O_{T}$ denotes the hyperplane orthogonal to $T$. The precise result is presented in \Cref{L2.6}. Roughly speaking, the main conclusion is that for ``most'' $w\in O_T$ in a small ball, 
\begin{equation}\label{eq:strategy2ndOrderGoal}
u(z+w)-u(x)\approx u(z)-u(x)+\frac{1}{2}w^{\top}F(z)w = \pm \|z-x\|+\frac{1}{2}w^{\top}F(z)w,
\end{equation}
where the equality is due to $x$ and $z$ lying on~$T$, there is no first-order term as $w$ is orthogonal to the ray and hence to $\nabla u$, and $F(z)$ is a surrogate for the Hessian of~$u$ in a sense explained next.

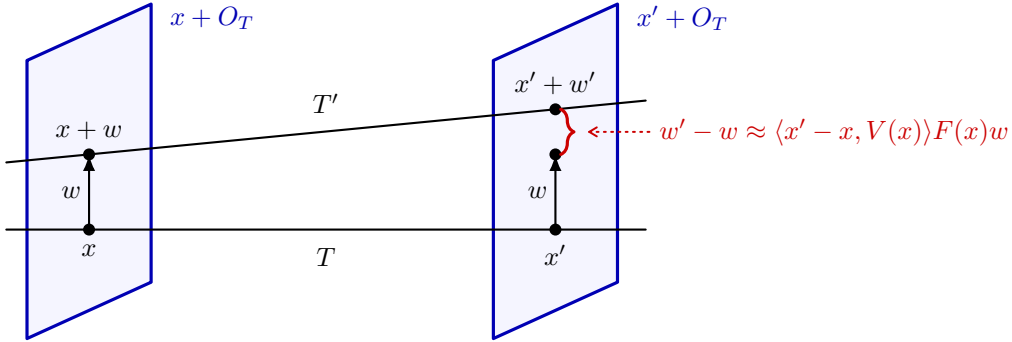
\begin{figure}[tbh]
  \centering
  \begin{tikzpicture}[
      line cap=round, line join=round,
      >={Latex[length=2.2mm]},
      every node/.style={font=\small}
    ]

    \def\x{0.0}      
    \def\xp{6.2}     
    \def\w{1.0}      
    \def\wp{1.6}     
    \def\H{2.1}      

    \def\dx{1.65}
    \def\dy{0.75}

    \pgfmathsetmacro{\yshift}{0.25*\H}

    \tikzset{
      plane/.style={
        draw=blue!70!black, very thick,
        fill=blue!15, fill opacity=0.25
      },
      pt/.style={circle, fill=black, inner sep=1.6pt},
      redbrace/.style={
        decorate,
        draw=red!80!black,
        very thick,
        decoration={brace, amplitude=5pt, raise=2pt, mirror}
      }
    }

    \coordinate (x)     at (\x,0);
    \coordinate (xw)    at (\x,\w);
    \coordinate (xp)    at (\xp,0);
    \coordinate (xpw)   at (\xp,\w);
    \coordinate (xpwp)  at (\xp,\wp);

    \pgfmathsetmacro{\m}{(\wp-\w)/(\xp-\x)} 
    
    \coordinate (Tleft)   at (-1.1,0);
    \coordinate (Tright)  at (\xp+1.2,0);
    \coordinate (TpLeft)  at (-1.1,{\w+\m*(-1.1-\x)});
    \coordinate (TpRight) at (\xp+1.2,{\w+\m*((\xp+1.2)-\x)});

    \coordinate (L1)  at (\x,{-\H+\yshift});
    \coordinate (L2)  at (\x,{ \H+\yshift});
    \coordinate (L1b) at ($(L1)-0.5*(\dx,\dy-1)$);
    \coordinate (L2b) at ($(L2)-0.5*(\dx,\dy)$);
    \coordinate (L3b) at ($(L2)+0.5*(\dx,\dy)$);
    \coordinate (L4b) at ($(L1)+0.5*(\dx,\dy+1)$);
    \draw[plane] (L1b)--(L2b)--(L3b)--(L4b)--cycle;
    \node[blue!70!black, anchor=west] at ($(L3b)+(0.12,-0.18)$) {$x+O_T$};

    \coordinate (R1)  at (\xp,{-\H+\yshift});
    \coordinate (R2)  at (\xp,{ \H+\yshift});
    \coordinate (R1b) at ($(R1)-0.5*(\dx,\dy-1)$);
    \coordinate (R2b) at ($(R2)-0.5*(\dx,\dy)$);
    \coordinate (R3b) at ($(R2)+0.5*(\dx,\dy)$);
    \coordinate (R4b) at ($(R1)+0.5*(\dx,\dy+1)$);
    \draw[plane] (R1b)--(R2b)--(R3b)--(R4b)--cycle;
    \node[blue!70!black, anchor=west] at ($(R3b)+(0.12,-0.18)$) {$x'+O_T$};

    \draw[black, thick] (Tleft) -- (Tright);
    \node[below=3pt] at ($(Tleft)!0.5!(Tright)$) {$T$};

    \draw[black, thick] (TpLeft) -- (TpRight);
    \node[above=3pt] at ($(TpLeft)!0.5!(TpRight)$) {$T'$};

    \node[pt] at (x) {};
    \node[pt] at (xw) {};
    \node[pt] at (xp) {};
    \node[pt] at (xpw) {};
    \node[pt] at (xpwp) {};

    \node[below=2pt] at (x) {$x$};
    \node[above=2pt] at (xw) {$x+w$};

    \node[below=2pt] at (xp) {$x'$};
    \node[above=2pt] at (xpwp) {$x'+w'$};

    \draw[->, thick] (x)  -- (xw)  node[midway, left=0.5pt, xshift=1.5pt] {$w$};
    \draw[->, thick] (xp) -- (xpw) node[midway, left=0.5pt, xshift=1.5pt] {$w$};

    \draw[redbrace] (xpw) -- (xpwp);

    \node[anchor=west, text=red!80!black] (formula) at ($(xpw)!0.5!(xpwp) + (1.25,0)$) {%
      $w'-w \approx \langle x'-x, V(x)\rangle F(x) w$%
    };

    \draw[red!80!black, dotted, line width=0.75pt]

      (formula.west) -- node[midway, above=1pt] {} ($(xpw)!0.5!(xpwp) + (0.45,0)$);

    \draw[red!80!black, line width=0.75pt]
      ($(xpw)!0.5!(xpwp) + (0.45,0)$) ++(0.10,0.06) -- ++(-0.10,-0.06);
    \draw[red!80!black, line width=0.75pt]
      ($(xpw)!0.5!(xpwp) + (0.45,0)$) ++(0.10,-0.06) -- ++(-0.10,0.06);

  \end{tikzpicture}
  \caption{Geometry of nearby transport rays. The ray $T'$ intersects the orthogonal hyperplanes $x+O_T$ and $x'+O_T$ at $x+w$ and $x'+w'$, respectively. We quantify the deviation $w'-w$ between $x'+w'$ and the translation of $x+w$ in the direction of~$T$.}
  \label{fig:ray-geometry-orthogonal}
\end{figure}

\textbf{Step 1: First-order expansion of $V:=\nabla u$.} Indeed, the first derivative $V:=\nabla u$ exists a.e.\ due to the Lipschitz continuity of $u$. To obtain a second-order expansion as in~\eqref{eq:strategy2ndOrderGoal}, we would like to show that $V$ admits a first-order expansion $V(x+w)\approx V(x)+\nabla V(x) w$ for $w\in O_T$ with $\|w\|$ small. However, $\nabla V(x)$ typically does not exist. We leverage the ``countable Lipschitz'' property of \cite{MR1862796} to see that the \emph{approximate derivative} of $V$, denoted $\approxgrad V$, exists in the sense of geometric measure theory \cite{MR3409135}. The approximate derivative $F:=\approxgrad  V$ of the function $V$ is defined via
\begin{equation}\label{eq:strategyApproxDerivClassical}
  \lim_{t\rightarrow 0^+} t^{-d} \mathcal{L}^d\big(\big\{y\in B_d(x,t):\|V(y)-V(x)-F(x)(y-x)\|>\delta\|y-x\|\big\}\big)=0
\end{equation}
for any $\delta>0$, where $B_d(x,t):=\{y\in\mathbb{R}^d:\|y-x\|\leq t\}$. Off the shelf, this tool is not useful, because our quantities of interest are averages over $(d-1)$-dimensional balls in $O_{T}$, not in~$\R^d$. An important insight of our analysis is that, due to the regularity of $u$ in the ray direction, $F=\approxgrad V$ satisfies
  \begin{equation}\label{eq:strategyApproxDerivOrthogonal}
    \lim_{t\rightarrow 0^+} t^{-(d-1)} \mathcal{H}^{d-1}\big(\big\{w\in O_T:\|w\|\leq t, \|\nabla u(x+w)-\nabla u(x)-F(x)w\|>\delta t\big\}\big) =0
  \end{equation}
for a.e.\ $x$ in the interior of a ray $T$. This gives a notion of approximate second derivative of $u$ \emph{w.r.t.\ averages in the orthogonal hyperplane~$O_{T}$} of~$T$, a key tool for our goal~\eqref{eq:strategy2ndOrderGoal}. To sharpen the classical notion~\eqref{eq:strategyApproxDerivClassical} of approximate derivative to our notion~\eqref{eq:strategyApproxDerivOrthogonal} involving averaging only in orthogonal directions, we require control over the geometry of nearby transport rays: for a ray $T'$ close to $T$ and two points $x,x'$ in the interior of $T$, we need to relate the point $x+w$ where $T'$ intersects $x+O_T$ to the point $x'+w'$ where $T'$ intersects $x'+O_T$ (cf.\ \cref{fig:ray-geometry-orthogonal}).
Specifically, we establish that $w'\approx w+\langle x'-x,V(x)\rangle F(x) w$ for ``most'' $w\in O_T$ in a small ball, cf.~\eqref{L2.4.e9n}. A key input for this relation is the Lipschitz regularity of the ray direction, first stated in \cite[Lemma 16]{MR1862796} and slightly strengthened in \Cref{L2.0} for our needs: 
\begin{equation}\label{eq:strategyFlipschitz}
\|V(y)-V(y')\| \leq L(y,y')\|y-y'\|,
\end{equation}
where $L(y,y')$ depends on the distance of $y,y'$ to the ends of their transport rays $T,T'$. Unfortunately, $L(y,y')$ deteriorates when one of the points $y,y'$ approaches an endpoint of its ray (and a counterexample shows that this is unavoidable). Circumventing this difficulty necessitates an additional stability result, stated in \Cref{L2.2}: the distance from $x+w$ to the upper (lower) end of its ray $T'$ is close to the distance from $x$ to the upper (lower) end of its ray $T$, for ``most'' perturbations $w\in O_T$.

These insights also yield regularity and structural properties of the approximate derivative $F$. In particular, for $x,x'$ in the interior of the same ray~$T$, we show a Lipschitz estimate
\begin{equation}\label{eq:F-Lipschitz}
    \|F(x)-F(x')\| \leq C(x,x')\|x-x'\|,
\end{equation}
where $C(x,x')$ depends on the distances of $x,x'$ to the ends of $T$. The main properties of $V$ and $F$ are collected in \Cref{L2.4}; they are used in the proofs of both the lower and upper bounds.

\textbf{Step 2: Second-order expansion of $u$.} With the expansion $V(x+w)\approx V(x)+F(x) w$ established, we use a geometric argument to deduce the second-order expansion~\eqref{eq:strategy2ndOrderGoal}. The key step is to control local increments of the form $u(x+\tilde{w})-u(x+w)$ for nearby $w,\tilde{w}\in O_T$. To this end, we introduce a ``good set'' $\mathcal{G}_t(x)\subseteq O_T\cap B_d(0,t)$ on which both the stability of distances to the ray ends and the first-order expansion of $V(x+w)$ hold. We show in \Cref{prep1.1} that for $w\in \mathcal{G}_t(x)$,
\begin{equation}\label{uwwwne}
    u(x+\tilde{w})-u(x+w)\approx\langle F(x)w,\tilde{w}-w\rangle,
\end{equation}
with an explicit error bound.

\begin{figure}[tbh]
  \centering
  \begin{tikzpicture}[
      line cap=round, line join=round,
      >={Latex[length=2.2mm]},
      every node/.style={font=\small}
    ]

    \tikzset{
      pt/.style={circle, fill=black, inner sep=1.6pt},
      plane/.style={
        draw=blue!70!black, very thick,
        fill=blue!15, fill opacity=0.25
      },
      plabel/.style={text=blue!70!black},
      aux/.style={black!65, line width=0.75pt, densely dashed},
      secant/.style={black, line width=1.0pt},
      ann/.style={
        fill=white, fill opacity=0.95, text opacity=1,
        inner sep=1.0pt
      }
    }

    \def\H{1.78}      
    \def\dx{1.18}     
    \def\dy{0.66}
    \def\hw{0.56}     

    \pgfmathsetmacro{\sx}{\hw*\dx}
    \pgfmathsetmacro{\sy}{\hw*\dy}

    \coordinate (a)   at (0.45,0);
    \coordinate (xpp) at (1.62,0);   
    \coordinate (x)   at (5.40,0);
    \coordinate (xp)  at (9.28,0);   
    \coordinate (b)   at (10.5,0);

    \foreach \P/\name in {xpp/pp, x/c, xp/p}{
      \coordinate (\name LL) at ($ (\P)+(-\sx,-\H-\sy) $);
      \coordinate (\name UL) at ($ (\P)+(-\sx, \H-\sy) $);
      \coordinate (\name UR) at ($ (\P)+( \sx, \H+\sy) $);
      \coordinate (\name LR) at ($ (\P)+( \sx,-\H+\sy) $);
    }

    \draw[plane] (ppLL)--(ppUL)--(ppUR)--(ppLR)--cycle;
    \draw[plane] (cLL)--(cUL)--(cUR)--(cLR)--cycle;
    \draw[plane] (pLL)--(pUL)--(pUR)--(pLR)--cycle;

    \draw[black, thick] (a) -- (b);
    \node[below=2pt] at (a) {$a$};
    \node[below=2pt] at (b) {$b$};
    \node[anchor=west] at ($(b)+(0.22,0)$) {$T$};

    \node[plabel, anchor=west] at ($(ppUR)+(0.14,-0.2)$) {$x''+O_T$};
    \node[plabel, anchor=west] at ($(cUR)+(0.14,-0.2)$) {$x+O_T$};
    \node[plabel, anchor=west] at ($(pUR)+(0.14,-0.2)$) {$x'+O_T$};

    \node[pt] at (xpp) {};
    \node[pt] at (x)   {};
    \node[pt] at (xp)  {};

    \node[below=2pt] at (xpp) {$x''$};
    \node[below=4.3pt] at (x)   {$x$};
    \node[below=2pt] at (xp)  {$x'$};

    \coordinate (xppw) at (1.44,0.86);
    \coordinate (xpw)  at (9.3,0.2);
    \coordinate (xw)   at ($(xppw)!0.50!(xpw)$);
    \coordinate (xtw)  at ($(x)+(0.03,0.9)$);

    \draw[aux] (xtw) -- (xppw);
    \draw[aux] (xtw) -- (xpw);

    \draw[secant] (xppw) -- (xpw);

    \node[pt] at (xppw) {};
    \node[pt] at (xw)   {};
    \node[pt] at (xtw)  {};
    \node[pt] at (xpw)  {};

    \node[above=2.5pt, xshift=5.6pt] at (xppw) {$x''+w''$};
    \node[below=0.45pt, xshift=1pt] at (xw)   {$x+w$};
    \node[above left=2pt, xshift=18pt]               at (xtw)  {$x+\tilde{w}$};
    \node[above=3pt] at (xpw)  {$x'+w'$};

  \end{tikzpicture}
  \caption{Construction for bounding $u(x+\tilde{w})-u(x+w)$ in~\eqref{uwwwne}}
  \label{fig:transport-three-sections}
\end{figure}
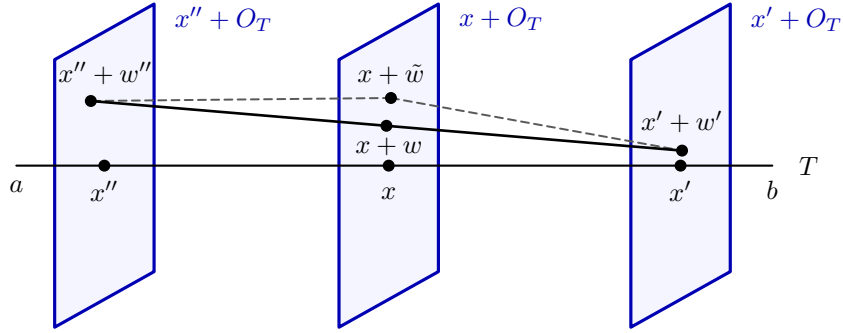

The underlying geometric idea is illustrated in Figure~\ref{fig:transport-three-sections}. Suppose $x$ lies in the interior of the ray $T$ with upper end $a$ and lower end $b$, and fix points $x',x''$ in the interior of $T$ with $x'$ near $b$ and $x''$ near $a$. For $w\in\mathcal{G}_t(x)$ and $t$ sufficiently small, the stability of distances to the ray ends ensures that the ray through $x+w$ intersects $x'+O_T$ and $x''+O_T$; we denote these intersection points by $x'+w'$ and $x''+w''$, respectively. Along the ray, we therefore have
\begin{equation*}
    u(x''+w'')-u(x+w)=\|(x''+w'')-(x+w)\|, \quad u(x+w)-u(x'+w')=\|(x+w)-(x'+w')\|.
\end{equation*}
Since $u\in\mathrm{Lip}_1(\mathbb{R}^d)$, a telescoping decomposition yields the upper bound
\begin{align*}
    u(x+\tilde{w})-u(x+w)&=\big(u(x+\tilde{w})-u(x'+w')\big)-\big(u(x+w)-u(x'+w')\big)\nonumber\\
    & \leq \|(x+\tilde{w})-(x'+w')\|-\|(x+w)-(x'+w')\|.
\end{align*}
A lower bound follows by an analogous argument using $x''+w''$. A detailed analysis, leveraging the first-order expansion $V(x+w)\approx V(x)+F(x)w$, then yields \eqref{uwwwne}. 

To pass from the local estimate \eqref{uwwwne} to the global second-order expansion \eqref{eq:strategy2ndOrderGoal}, we discretize the segment from $0$ to $w$ by selecting $0=\theta_0<\theta_1<\cdots<\theta_m< \theta_{m+1}=1$ such that $\theta_i w\in \mathcal{G}_t(x)$ and $\|\theta_{i+1}w-\theta_i w\|$ is small. Applying the local estimate to consecutive pairs gives
\[u(x+\theta_{i+1}w)-u(x+\theta_i w)\approx \langle F(x)(\theta_i w),\theta_{i+1}w-\theta_iw\rangle=\theta_i(\theta_{i+1}-\theta_i)\langle F(x)w,w\rangle.\] Summing over $i=0,1,\dots,m$ and invoking $u(x+w)-u(x)=\sum_{i=0}^m \big(u(x+\theta_{i+1}w)-u(x+\theta_iw)\big)$ yields the second-order expansion 
for $u(x+w)-u(x)$, and by renaming the variable, also for $u(z+w)-u(z)$. The desired second-order expansion 
\eqref{eq:strategy2ndOrderGoal} now follows by adding $u(z)-u(x)=\pm\|z-x\|$.

\textbf{Detachment and lower bound.} On the strength of these regularity results for $u$, we then proceed to prove \cref{th:lowerBound}. For simplicity of exposition, suppose that $\gamma'_{\eps_j}$ admits a density $\phi_{\eps_j}(x,y)$ wrt.\ $\mathcal{L}^d\otimes\mathcal{L}^d$ that is sufficiently integrable. For $x$ in the interior of a ray $T$ and $t\in\mathbb{R}$, we define the marginal density obtained by integrating out the component of $y$ orthogonal to $T$,
\begin{equation*}
    \psi_{\eps_j}(x,t):=\int_{O_T}\phi_{\eps_j}(x,x+tV(x)+w)d\mathcal{H}^{d-1}(w).
\end{equation*}
Applying Jensen's inequality, \cref{th:lowerBound} can be reduced to two bounds. First, a lower bound on
\begin{equation}\label{lbentropy}
\liminf_{j\rightarrow\infty}\int_{\mathcal{X}\times\mathbb{R}}\psi_{\varepsilon_j}(x,t)\log(\psi_{\varepsilon_j}(x,t))dxdt,
\end{equation}
which is shown in \cref{L3.6} using the variational representation of entropy and properties of Carath\'eodory functions. Second, a precise upper bound on
\begin{equation}\label{ubenergy}
    \limsup_{j\rightarrow\infty}\int_{\mathcal{X}\times\mathbb{R}}G_{\eps_j}(x,x+tV(x))\psi_{\eps_j}(x,t)dxdt,
\end{equation}
where, for $x,z$ in the interior of the same ray~$T$,
\begin{align*}
      G_{\eps_j}(x,z) = \log\bigg(\int_{O_T} e^{-E(x,z+\sqrt{\eps_j}w)/ \eps_j}  d\mathcal{H}^{d-1}(w)\bigg)
\end{align*}
and $E(x,y):=\|x-y\|-u(x)+u(y)$ is the \emph{detachment} of the potential~$u$.
(See \cref{Def2} for the precise expression involving an additional truncation.) Using our second-order expansion~\eqref{eq:strategy2ndOrderGoal}, we obtain, for $x,z$ with $\langle x-z,V(x)\rangle >0$ and $w\in O_T$ in a small ball, the following expansion of the detachment:
\begin{align*}
    E(x,z+w)&=\|x-(z+w)\|- u(x)+u(z+w)\\
    &\approx \sqrt{\|z-x\|^2+\|w\|^2}-\|z-x\|+\frac{1}{2}w^{\top}F(z)w\approx \frac{\|w\|^2}{2\|z-x\|}+\frac{1}{2}w^{\top}F(z)w,
\end{align*}
where the last step is by direct Taylor expansion. Combining this asymptotic description of $E(x,y)$ with complementary uniform estimates on $E(x,y)$ detailed in \cref{L3.2}, we show in \cref{P-lb-1} that $G_{\eps_j}(x,z)\rightarrow G(x,z)$, where 
\begin{equation*}
    G(x,z)=\log\bigg(\int_{O_T}\exp\bigg(-\frac{1}{2}w^{\top}\bigg(\frac{1}{\|x-z\|}\mathbf{I}_d+F(z)\bigg)w\bigg) d\mathcal{H}^{d-1}(w)\bigg).
\end{equation*}

Finally, to bound \eqref{ubenergy}, we write
\begin{align}
    \int_{\mathcal{X}\times\mathbb{R}}G_{\eps_j}(x,x+tV(x))\psi_{\eps_j}(x,t)dxdt   \leq\,& \int_{\mathcal{X}\times\mathbb{R}}G(x,x+tV(x))\psi_{\eps_j}(x,t)dxdt \label{eq:overview-G-split}\\
    &+\int_{\mathcal{X}\times\mathbb{R}}|G_{\eps_j}(x,x+tV(x))-G(x,x+tV(x))|\psi_{\eps_j}(x,t)dxdt. \nonumber
\end{align}
The limit of the first term on the right-hand side is given in \cref{L3.7} using weak convergence arguments exploiting the Lipschitz regularity~\eqref{eq:F-Lipschitz} of~$F$, whereas the second term converges to zero by dominated convergence. (The arguments in the body of the text are more involved because the ``sufficient integrability'' of $\psi_{\eps_j}$ that we assumed for this sketch does not hold in general.)

\subsubsection{Proof strategy for the upper bound}\label{se:upper}

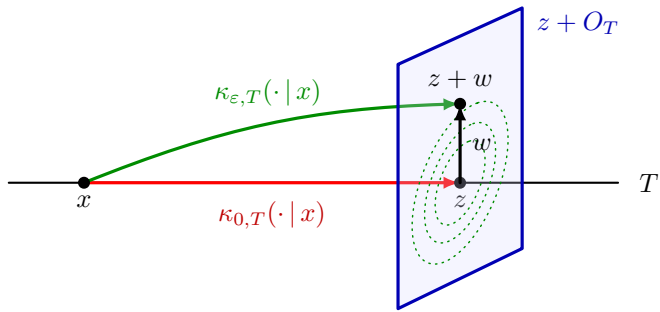
\begin{figure}[tbh]
  \centering
  \begin{tikzpicture}[
      line cap=round, line join=round,
      >={Latex[length=2.2mm]},
      every node/.style={font=\small}
    ]

    \tikzset{
      pt/.style={circle, fill=black, inner sep=1.6pt},
      plane/.style={
        draw=blue!70!black, very thick,
        fill=blue!15, fill opacity=0.25
      },
      gdiff/.style={draw=green!55!black, line width=0.55pt, dotted},
      garrow/.style={draw=green!55!black, very thick},
      glabel/.style={text=green!55!black},
      rlabel/.style={text=red!75!black},
    }

    \def\H{1.95}     
    \def\dx{1.65}    
    \def\dy{0.75}

    \def\wlen{1.05}  

    \def\a{0.78}     
    \def\b{0.36}     
    \def\vstretch{1.20}

    \def\beta{0.18}

    \pgfmathsetmacro{\Axx}{\beta*\dx}
    \pgfmathsetmacro{\Ayy}{\vstretch*(1 + \beta*\dy)}
    \pgfmathsetmacro{\Bxx}{\dx}
    \pgfmathsetmacro{\Byy}{\dy - \beta}

    \coordinate (x)  at (0,0);
    \coordinate (z)  at (5.0,0);
    \coordinate (zw) at ($(z)+(0,\wlen)$);

    \coordinate (Tleft)  at (-1.0,0);
    \coordinate (Tright) at (7.1,0);
    \draw[black, thick] (Tleft) -- (Tright);
    \node[anchor=west] at ($(Tright)+(0.15,0)$) {$T$};

    \draw[->, very thick, red] (x) -- (z);
    \node[below=3pt, rlabel] at ($(x)!0.5!(z)$) {$\kappa_{0,T}(\cdot\,|\,x)$};

    \draw[->, garrow]
      (x) to[bend left=10]
      node[midway, above=2.5pt, glabel] {$\kappa_{\varepsilon, T}(\cdot\,|\,x)$}
      (zw);

    \node[pt] at (x) {};
    \node[pt] at (z) {};
    \node[below=1pt] at (x) {$x$};
    \node[below=1pt] at (z) {$z$};

    \coordinate (P1)  at ($(z)+(0,-\H)$);
    \coordinate (P2)  at ($(z)+(0,\H)$);
    \coordinate (P1b) at ($(P1)-0.5*(\dx,\dy-1.3)$);
    \coordinate (P2b) at ($(P2)-0.5*(\dx,\dy)$);
    \coordinate (P3b) at ($(P2)+0.5*(\dx,\dy)$);
    \coordinate (P4b) at ($(P1)+0.5*(\dx,\dy+1.4)$);
    \draw[plane] (P1b)--(P2b)--(P3b)--(P4b)--cycle;
    \node[blue!70!black, anchor=west] at ($(P3b)+(0.12,-0.20)$) {\!$z+O_T$};

    \begin{scope}
      \clip (P1b)--(P2b)--(P3b)--(P4b)--cycle;
      \foreach \s in {1.00,0.74,0.52}{
        \draw[gdiff]
          plot[samples=180, smooth, domain=0:360]
            ({5.0 + (\s*\a)*cos(\x)*\Axx + (\s*\b)*sin(\x)*\Bxx},
             {0.0 + (\s*\a)*cos(\x)*\Ayy + (\s*\b)*sin(\x)*\Byy});
      }
    \end{scope}

    \draw[->, very thick] (z) -- (zw)
      node[midway, right=0.4pt, xshift=-0.1pt, yshift=-0.4pt] {$w$};

    \node[pt] at (zw) {};
    \node[above=2pt] at (zw) {$z+w$};

  \end{tikzpicture}
  \caption{Heuristic construction for the upper bound. The optimal transport kernel $\kappa_{0,T}(\cdot|x)$ sends~$x$ to a random point $z\in T$. To obtain an absolutely continuous kernel, we diffuse orthogonally to~$T$ by adding $w$, drawn from a $(d-1)$-dimensional Gaussian supported in $O_T$.}
  \label{fig:upperbound-diffusion}
\end{figure}

To show the upper bound in \Cref{th:upperBound}, we need to construct couplings $\tilde{\gamma}_{\eps}\in\Pi(\mu,\nu)$ such that 
\begin{equation}\label{ubb.estimates}
    \limsup_{\eps\rightarrow 0^+}\bigg\{\frac{\mathcal{C}_{\eps}(\tilde{\gamma}_{\eps})-\OT(\mu,\nu)}{\eps}+\frac{d-1}{2}\log{\eps}\bigg\}\leq \mathtt{S}(\mu,\nu;\gamma_0).
\end{equation}
The claim then follows as $\mathcal C_{\eps}(\gamma_{\eps})\leq \mathcal C_{\eps}(\tilde{\gamma}_{\eps})$. Heuristically, suppose that we have solved the constrained EOT problem~\eqref{opt.solver}, giving a kernel $\kappa_{0,T}(\cdot|x)$ that transports mass from $x$ along the transport ray $T$ to a random point $z\in T$, in a way so that $\gamma_{0}:=\mu\otimes \kappa_{0,T}(\cdot|x)$ is a coupling of~$(\mu,\nu)$. This coupling is singular, hence has $\mathcal C_{\eps}(\gamma_{0})=\infty$. We create an absolutely continuous kernel $\kappa_{\eps,T}(\cdot|x)$ by diffusing $\kappa_{0,T}(\cdot|x)$ \emph{orthogonally} to $T$ using a $(d-1)$-dimensional Gaussian distribution with variance depending on $\eps$ and $\approxgrad  V(x)$, chosen to optimally trade off transport and entropy costs. Thus, as illustrated in \cref{fig:upperbound-diffusion}, $x$ is transported to $z+w$, where $O_T$ is the orthogonal complement of $T$ and $w\in O_T$ is the Gaussian perturbation. A formal calculation of the optimal variance and the resulting EOT cost $\mathcal C_{\eps}$ suggests that~\eqref{ubb.estimates} would hold if we take $\tilde{\gamma}_{\eps}$ to be $\mu\otimes \kappa_{\eps,T}(\cdot|x)$. But of course, this is not a coupling of $(\mu,\nu)$ for $\eps>0$: its second marginal is some mixture of Gaussians instead of $\nu$. 

For the rigorous proof, we need a bona fide coupling $\tilde{\gamma}_{\eps}$ of $(\mu,\nu)$. Our strategy is to construct $\tilde{\gamma}_{\eps}$ as a sum $\tilde{\gamma}_{\eps}=\tilde{\gamma}_{\eps,1}+\tilde{\gamma}_{\eps,2}$. The first part, $\tilde{\gamma}_{\eps,1}$, is a sub-probability measure whose marginals are pointwise dominated by (and close to) $\mu,\nu$; we shall call this a \emph{sub-coupling} for brevity. The second part, $\tilde{\gamma}_{\eps,2}$, is a measure with small total mass, designed so that the sum of the two parts is exactly a coupling of $(\mu,\nu)$. The idea is that $\tilde{\gamma}_{\eps,1}$ accounts for most of the mass and approximately satisfies~\eqref{ubb.estimates}. Then, we need $\tilde{\gamma}_{\eps,2}$ so that the sum satisfies the marginal constraint \emph{and} the upper bound~\eqref{ubb.estimates} still goes through. The next paragraphs sketch some of the high-level ideas for this construction.

\textbf{Step 1: An initial sub-coupling mimicking the diffusion heuristic.} We start with a density $p_1$, defined rigorously in \eqref{p1}, whose disintegration along rays mimics a truncated version of the kernel $\kappa_{\eps,T}(\cdot|x)$ from the heuristics above. The truncation serves to avoid the singularity of certain functions near ray ends in the subsequent analysis; we do not enter such technicalities here. The main point is that, by construction, the first marginal of $p_1$ is dominated by $\mu$, and moreover, $p_1$ transports mass in a way similar to the desired heuristics. 

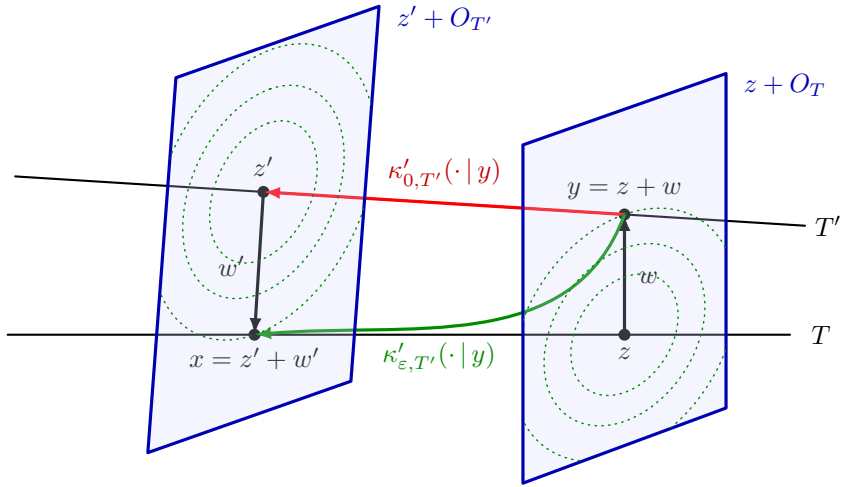
\begin{figure}[tbh]
  \centering
\begin{tikzpicture}[
  line cap=round, line join=round,
  >={Latex[length=2.2mm]},
  every node/.style={font=\small}
]

\tikzset{
  pt/.style={circle, fill=black, inner sep=1.6pt},
  plane/.style={draw=blue!70!black, very thick, fill=blue!15, fill opacity=0.25},
  gdiff/.style={draw=green!55!black, line width=0.55pt, dotted},
  garrow/.style={draw=green!55!black, very thick},
  glabel/.style={text=green!55!black},
  rlabel/.style={text=red!75!black},
}

\def\H{3.00}       
\def\DepthX{2.70}  
\def\DepthY{0.95}
\def\wlen{1.60}   

\def\bR{0.55}
\def\bL{0.55}

\coordinate (Tleft)  at (-1.2,0);
\coordinate (Tright) at (9.2,0);
\draw[black, thick] (Tleft) -- (Tright);
\node[anchor=west] at ($(Tright)+(0.15,0)$) {$T$};

\coordinate (z) at (7.0,0.0);
\coordinate (y) at (7.0,\wlen);

\pgfmathsetmacro{\zpx}{2.20}
\pgfmathsetmacro{\zpy}{1.90}
\coordinate (zp) at (\zpx,\zpy);

\pgfmathsetmacro{\yx}{7.0}
\pgfmathsetmacro{\yy}{\wlen}
\pgfmathsetmacro{\m}{(\yy-\zpy)/(\yx-\zpx)}                
\pgfmathsetmacro{\nxraw}{-\m}                              
\pgfmathsetmacro{\nyraw}{1}
\pgfmathsetmacro{\nlen}{sqrt(\nxraw*\nxraw + \nyraw*\nyraw)}
\pgfmathsetmacro{\nx}{\nxraw/\nlen}                         
\pgfmathsetmacro{\ny}{\nyraw/\nlen}                     
\pgfmathsetmacro{\TpLx}{-1.1}
\pgfmathsetmacro{\TpLy}{\zpy + \m*(\TpLx-\zpx)}
\pgfmathsetmacro{\TpRx}{9.4}
\pgfmathsetmacro{\TpRy}{\zpy + \m*(\TpRx-\zpx)}
\draw[black, thick] (\TpLx,\TpLy) -- (\TpRx,\TpRy);
\node[anchor=west] at (\TpRx,\TpRy) {$T'$};

\pgfmathsetmacro{\tint}{-\zpy/\ny}
\pgfmathsetmacro{\xx}{\zpx + \tint*\nx}
\coordinate (x) at (\xx,0);

\node[pt] at (z)  {};
\node[pt] at (y)  {};
\node[pt] at (zp) {};
\node[pt] at (x)  {};

\node[below=1pt] at (z) {$z$};
\node[above=2pt] at (y) {$y=z+w$};
\node[above=2pt] at (zp) {$z'$};
\node[below=1pt] at (x) {$x=z'+w'$};

\draw[->, very thick] (z) -- (y)
  node[midway, right=1pt, yshift=-1pt] {$w$};

\draw[->, very thick] (zp) -- (x)
  node[midway, left=1pt, yshift=-1pt] {$w'$};

\draw[->, very thick, red] (y) -- (zp);
\node[rlabel, above=1pt] at ($(y)!0.55!(zp)$) {\hspace{5mm}$\kappa_{0,T'}'(\cdot\,|\,y)$};

\draw[->, garrow]
  (y) to[out=-110, in=5]
  node[midway, below=3pt, glabel] {\hspace{-10mm}$\kappa_{\varepsilon,T'}'(\cdot\,|\,y)$}
  (x);

\coordinate (P1) at ($(z)+(0,-\H)$);
\coordinate (P2) at ($(z)+(0,\H)$);
\coordinate (P1b) at ($(P1)-0.5*(\DepthX,\DepthY-3)$);
\coordinate (P2b) at ($(P2)-0.5*(\DepthX,\DepthY)$);
\coordinate (P3b) at ($(P2)+0.5*(\DepthX,\DepthY)$);
\coordinate (P4b) at ($(P1)+0.5*(\DepthX,\DepthY+3.1)$);
\draw[plane] (P1b)--(P2b)--(P3b)--(P4b)--cycle;
\node[blue!70!black, anchor=west] at ($(P3b)+(0.10,-0.15)$) {$z+O_T$};

\pgfmathsetmacro{\aR}{\wlen}
\begin{scope}
  \clip (P1b)--(P2b)--(P3b)--(P4b)--cycle;
  \foreach \scale in {1.00,0.72,0.48}{
    \draw[gdiff]
      plot[samples=180, smooth, variable=\ang, domain=0:360]
      ({7.0 + (\scale*\bR)*sin(\ang)*(\DepthX)},
       {0.0 + (\scale*\aR)*cos(\ang) + (\scale*\bR)*sin(\ang)*(\DepthY)});
  }
\end{scope}

\coordinate (Q1) at ($ (zp) + (-\H*\nx, -\H*\ny) $);
\coordinate (Q2) at ($ (zp) + ( \H*\nx,  \H*\ny) $);
\coordinate (Q1b) at ($(Q1)-0.5*(\DepthX,\DepthY)$);
\coordinate (Q2b) at ($(Q2)-0.5*(\DepthX,\DepthY+2)$);
\coordinate (Q3b) at ($(Q2)+0.5*(\DepthX,\DepthY-2)$);
\coordinate (Q4b) at ($(Q1)+0.5*(\DepthX,\DepthY)$);
\draw[plane] (Q1b)--(Q2b)--(Q3b)--(Q4b)--cycle;
\node[blue!70!black, anchor=west] at ($(Q3b)+(0.10,-0.15)$) {$z'+O_{T'}$};

\pgfmathsetmacro{\aL}{abs(\tint)}
\begin{scope}
  \clip (Q1b)--(Q2b)--(Q3b)--(Q4b)--cycle;
  \foreach \scale in {1.00,0.72,0.48}{
    \draw[gdiff]
      plot[samples=180, smooth, variable=\ang, domain=0:360]
      ({\zpx + (\scale*\aL)*cos(\ang)*(\nx) + (\scale*\bL)*sin(\ang)*(\DepthX)},
       {\zpy + (\scale*\aL)*cos(\ang)*(\ny) + (\scale*\bL)*sin(\ang)*(\DepthY)});
  }
\end{scope}

\end{tikzpicture}
 \caption{Heuristic construction for the density $p_2$. The ``reverse'' optimal transport kernel $\kappa_{0,T'}'(\cdot|y)$ sends~$y$ to a random point $z'\in T'$. We then diffuse orthogonally to~$T'$ by adding $w'$, drawn from a $(d-1)$-dimensional Gaussian supported in $O_{T'}$. The disintegration of $p_2$ mimics a truncated version of this continuous kernel.}
  \label{fig:p2}
\end{figure}

From this construction, only the first marginal is suitable for a sub-coupling: there is no reason why the second marginal of $p_1$ should be dominated by~$\nu$. To build a valid sub-coupling $\tilde{\gamma}_{\eps,1}$, our first key idea is to \emph{reverse the roles of the marginals} and then take a minimum.  Indeed, we introduce a second density $p_2$, obtained by applying the construction of $p_1$ with the roles of $x$ and $y$ reversed. More precisely, as indicated in \cref{fig:p2}, we consider the kernel $\kappa_{0,T'}'(\cdot|y)$ from the constrained EOT problem~\eqref{opt.solver}---which transports mass from $y$ along the ray $T'$ to a random point $z'$---and then diffuse $\kappa_{0,T'}'(\cdot|y)$ orthogonally to $T'$ using another $(d-1)$-dimensional Gaussian distribution with carefully chosen variance (which may differ from the variance used in the construction of $p_1$) to obtain a second continuous kernel $\kappa_{\eps,T'}'(\cdot|y)$. Thus, under $\kappa_{\eps,T'}'(\cdot|y)$, $y$ is transported to $z'+w'$, where $w'\in O_{T'}$ is the Gaussian perturbation. See~\eqref{p2} for the rigorous definition, and note that because the rays $T,T'$ need not be parallel,  $O_{T'}$ is not parallel to $O_{T}$ in general. The disintegration of $p_2$ then mimics a truncated version of $\kappa_{\eps,T'}'(\cdot|y)$. Since the  second marginal of $p_2$ is dominated by~$\nu$, and the first marginal of $p_1$ is dominated by~$\mu$, the density $\min\{p_1, p_2\}$ yields a valid sub-coupling; that is, its marginals $(\mu_1,\nu_1)$ satisfy $\mu_1\leq\mu$ and $\nu_1\leq\nu$ pointwise.

\textbf{Step 2: Balancing the marginals of the sub-coupling on each ray.} The density $p:=\min\{p_1, p_2\}$ from Step~1 yields a valid sub-coupling, meaning that its marginals $(\mu_0,\nu_0)$ satisfy $\mu_0\leq\mu$ and $\nu_0\leq\nu$ pointwise. However, it seems difficult to now construct a remainder measure that couples the missing marginals $(\mu-\mu_0,\nu-\nu_0)$ without invalidating the upper bound~\eqref{ubb.estimates}. To prepare for the construction of the remainder in the Steps 4--5, we modify $p$ by chipping away some of its mass, with the goal of obtaining a measure whose marginals are balanced on each ray. Specifically, we create a density $p_s$ such that the marginals $(\mu_1,\nu_1)$ of $\tilde{\gamma}_{\eps,1}:=\min\{p, p_s\}$ have the following property: for any transport ray $T$, the ``restrictions'' of $\mu_1$ and $\nu_1$ to $T$ have equal mass (where restriction is made rigorous by disintegration). Note that $\tilde{\gamma}_{\eps,1}$ is necessarily a sub-coupling, since we are removing further mass from the sub-coupling $p$. Of course, we must also proceed in a way so that the total mass of $\tilde{\gamma}_{\eps,1}$ remains close to one (cf.\ Step 3).

\begin{figure}[tbh]
  \centering
  \begin{tikzpicture}[
      line cap=round, line join=round,
      >={Latex[length=2.2mm]},
      every node/.style={font=\small}
    ]

    \tikzset{
      pt/.style={circle, fill=black, inner sep=1.6pt},
      ray/.style={black, thick},
      constr/.style={black!60, dashed, line width=0.55pt},
      lavg/.style={black!60, dashed, line width=0.55pt},
      parrow/.style={draw=blue!70!black, very thick},
      psarrow/.style={draw=red!75!black, very thick},
      plabel/.style={text=blue!70!black},
      pslabel/.style={text=red!75!black},
    }

    \pgfmathsetmacro{\mTp}{1.7/7}          
    \pgfmathsetmacro{\bTp}{1.20}          
    \pgfmathsetmacro{\mL}{\mTp/2}         
    \pgfmathsetmacro{\mOL}{-1/\mL}        

    \pgfmathsetmacro{\xmin}{-0.8}
    \pgfmathsetmacro{\xmax}{7.3}

    \pgfmathsetmacro{\Lx}{3.60}
    \pgfmathsetmacro{\Ly}{1.20}

    \pgfmathsetmacro{\Lyleft}{\Ly + \mL*(\xmin-\Lx)}
    \pgfmathsetmacro{\Lyright}{\Ly + \mL*(\xmax-\Lx)}

    \pgfmathsetmacro{\Tpyleft}{\bTp + \mTp*(\xmin)}
    \pgfmathsetmacro{\Tpyright}{\bTp + \mTp*(\xmax)}

    \pgfmathsetmacro{\xX}{2.00}
    \coordinate (x) at (\xX,0.00);

    \pgfmathsetmacro{\yX}{5.00}
    \pgfmathsetmacro{\yY}{\bTp + \mTp*\yX}
    \coordinate (y) at (\yX,\yY);

    \pgfmathsetmacro{\xprimeX}{(\bTp + \xX*\mOL)/(\mOL - \mTp)}
    \pgfmathsetmacro{\xprimeY}{\bTp + \mTp*\xprimeX}
    \coordinate (xp) at (\xprimeX,\xprimeY);

    \pgfmathsetmacro{\yprimeX}{\yX - \yY/\mOL}
    \coordinate (yp) at (\yprimeX,0.00);

    \coordinate (Lleft)  at (\xmin,\Lyleft);
    \coordinate (Lright) at (\xmax,\Lyright);

    \pgfmathsetmacro{\Ix}{(\Ly - \mL*\Lx - \yY + \mOL*\yX)/(\mOL - \mL)}
    \pgfmathsetmacro{\Iy}{\Ly + \mL*(\Ix-\Lx)}
    \coordinate (I) at (\Ix,\Iy);

    \draw[ray] (\xmin,0) -- (\xmax,0);
    \node[anchor=west] at (\xmax+0.18,0) {$T$};

    \draw[ray] (\xmin,\Tpyleft) -- (\xmax,\Tpyright);
    \node[anchor=west] at (\xmax+0.18,\Tpyright) {$T'$};

    \draw[lavg] (Lleft) -- (Lright);
    \node[anchor=west, black!60] at ($(Lright)+(0.18,0.02)$) {$L$};

    \draw[constr] ($(I)+0.15*(1,\mL)$) -- ($(I)-0.15*(1,\mL)$);
    \draw[constr] ($(I)+0.12*(-\mL,1)$) -- ($(I)-0.12*(-\mL,1)$);

    \coordinate (RAa) at ($(I)+0.18*(1,\mL)$);
    \coordinate (RAb) at ($(RAa)+0.18*(-\mL,1)$);
    \coordinate (RAc) at ($(I)+0.18*(-\mL,1)$);
    \draw[black!60, line width=0.55pt] (RAa) -- (RAb) -- (RAc);

    \draw[->, parrow]
      (x) to[bend left=26]
      node[midway, below=-3pt, plabel] {$~~~p$}
      (y);

    \draw[->, psarrow]
      (xp) to[bend right=18]
      node[midway, below=.1pt, pslabel] {$p_s$}
      (yp);

    \draw[constr] (x) -- (xp);
    \node[black!60, anchor=south east] at ($(x)!0.56!(xp)+(0.10,-0.7)$) {$x+O_L$};

    \draw[constr] (y) -- (yp);
    \node[black!60, anchor=west] at ($(y)!0.52!(yp)+(0.05,-0.3)$) {$y+O_L$};

    \node[pt] at (x)  {};
    \node[pt] at (y)  {};
    \node[pt] at (xp) {};
    \node[pt] at (yp) {};

    \node[below=1pt] at (x)  {$x$};
    \node[above=2pt] at (y)  {$y$};
    \node[above=2pt] at (xp) {$x'$};
    \node[below=1pt] at (yp) {$y'$};

  \end{tikzpicture}
  \caption{Symmetry construction on a pair of transport rays. The symmetry points of $x\in T$ and $y\in T'$ are $x'=(x+O_L)\cap T'$ and $y'=(y+O_L)\cap T$, where $O_L$ is the orthogonal complement of the average direction $L$ of $T$ and $T'$.}
  \label{fig:symmetry-points}
\end{figure}
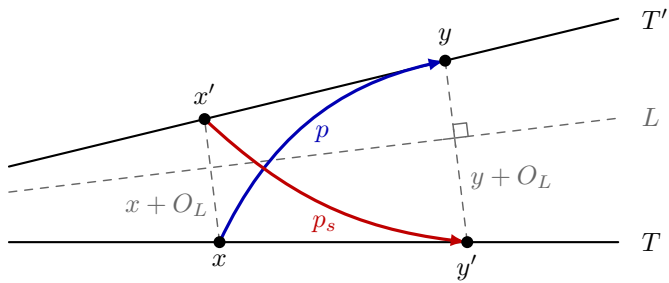

To explain the construction, let us focus on a transport ray $T$. Intuitively, $p$ transports the mass originating on $T$ away from $T$, onto other transport rays. Consider one such ray $T'$, and note that $p$ also transports some mass from $T'$ onto $T$. If for every pair $(T,T')$, the amount of mass flowing from $T$ to $T'$ were exactly the same as the amount from $T'$ to $T$, then the marginals of $p$ would be balanced on~$T$. The density $p$ does not have that property, but if we can obtain sufficient regularity, we may hope that the difference between the two directions is small. Hence, we aim to construct a density $p_s$ such that $\min\{p,p_s\}$ is (exactly) balanced on each ray $T$. This is based on the following geometric idea. Consider two rays $T$ and $T'$, as well as  points $x$ in the interior of $T$ and $y$ in the interior of $T'$. In \eqref{x'd}--\eqref{alt_def}, we then define their ``symmetry points'' $x'\in T'$ and $y'\in T$ as illustrated in \cref{fig:symmetry-points}: Consider the line $L$ whose direction is the average of the directions of $T$ and $T'$, and let $O_L$ be its orthogonal complement. Then $x'$ is the intersection point of $T'$ with $x+O_L$ whereas $y'$ is the intersection point of $T$ with $y+O_L$. Intuitively, if the rays $T$ and $T'$ are close, the mass $p_1(x,y)$ transported from $x$ to $y$ should be similar to the mass $p_1(x',y')$ transported from $x'$ to $y'$. To transport exactly the same amount, we want to remove the difference from the larger value. 

Intuitively, this is achieved by setting $p_s(x,y):=p(x',y')$ and then taking $\min\{p,p_s\}$. That is indeed valid when the rays $T$ and $T'$ are perfectly parallel. In the general case, however, we must account for the change in direction of the rays, which leads to different volume elements along rays when expressed in Euclidean coordinates.  To correct for this, $p_s(x,y)$ is defined as the product of $p(x',y')$ with an explicit Jacobian-type scaling factor. Since the change in direction corresponds to the derivative of $\nabla u$, this term involves the function~$F$ which proxies for the second derivative of~$u$, and the regularity of~$F$ becomes crucial for the downstream analysis (see Step~3). Then $\min\{p, p_s\}$ is the density (w.r.t.\ $\mathcal{L}^d\otimes\mathcal{L}^d$) of a sub-coupling $\tilde{\gamma}_{\eps,1}$ that is balanced on each ray.

In the body of the text, we actually do not apply this idea to $p$ directly.  We apply it separately to~$p_2$ and~$p_1$, yielding the densities $p_3$ and $p_4$ defined rigorously in \eqref{p3}--\eqref{p4}, and then take $\min\{p_1, p_2, p_3, p_4\}$. This is to work with simpler analytic expressions, but ultimately leads to the same outcome.

\textbf{Step 3: The mass of the sub-coupling is close to one.} For our approach to work, we need to show that the sub-coupling constructed so far transports almost all of the mass. This result, stated more precisely in \cref{P4.3}, boils down to proving
\begin{equation}\label{def_re}
    \lim_{\eps\to 0^{+}}\int (p_1(x,y)-p_i(x,y))_{+}dxdy = 0 \quad\text{for}\quad i=2,3,4.
\end{equation}
Consider first $i=2$ and recall that $p_2$ is constructed in a similar way as $p_1$ but starting from the second marginal. Our claim means that for $x\in T$ and $y=z+w\in T'$ as in \Cref{fig:p2}, the values $p_1(x,y)$ and $p_2(x,y)$ are similar. This is surprising at first glance, since for non-parallel rays, the vectors $w$ and $w'$ may have substantially different length and direction. Remarkably, this can be compensated by the difference in the chosen covariance matrices $\Lambda_1,\Lambda_2
$ defining the Gaussian perturbations in the constructions of $p_1,p_2$, albeit in a nontrivial way. Indeed, in the geometry of \Cref{fig:p2}, the quadratic forms governing the Gaussians are $w^{\top} \Lambda_1 w$ and $(w')^{\top}\Lambda_2 w'$. Neither the perturbations $w,w'$ nor the covariances matrices $\Lambda_1,\Lambda_2$ are close to one another. However, on the strength of the structural relation \eqref{Eq2.6.2} for $F=\approxgrad V$ established in \Cref{L2.4}, 
\begin{equation}\label{eq:Lambda-relation-in-overview}
    \Lambda_2\approx \big((\mathbf{I}_d+\|x-z\|F(z))^{-1}\big)^{\top}\Lambda_1 \big(\mathbf{I}_d+\|x-z\|F(z)\big)^{-1}.
\end{equation}
Using the approximate geometric relation in \Cref{fig:ray-geometry-orthogonal}, we further show 
$
    w'\approx -(\mathbf{I}_d+\|x-z\|F(z))w. 
$ 
Combining these two approximations yields the desired result that  $(w')^{\top}\Lambda_2 w' \approx w^{\top} \Lambda_1 w$; cf.~\eqref{Claim.a3.eq1} for the precise statement.

We do not detail the proofs of~\eqref{def_re} for $i=3,4$ here, but confine ourselves to some general remarks about the regularity issues that arise. Novel regularity results are needed for quantities including $F=\approxgrad V$ and the density $\tilde{h}_T$ of $\kappa_T$ (the ray-wise coupling in \cref{Definition1.1}) with respect to $\tilde{\mu}_T\otimes\tilde{\nu}_T$. \emph{Along} each transport ray, these functions enjoy strong regularity. For $\tilde{h}_T$, \cref{Lem2} derives Lipschitz regularity from properties of the constrained EOT problem~\eqref{opt.solver}. For $F$, we have the Lipschitz regularity mentioned in~\eqref{eq:F-Lipschitz} above. \emph{Across} different rays, however, continuity may fail. To overcome this, a recurring tool in our proofs is the use of \emph{Lebesgue points} (in the sense of geometric measure theory). To enable this tool, we develop a parametrization that explicitly separates variation \emph{along} rays from variation \emph{across} rays. Specifically, we localize so that transport rays can be parametrized by their intersection point~$q$ with a fixed hyperplane. Any point $x$ on some ray is then described by $q$ together with the signed distance $s$ from $x$ to the hyperplane; note how a change in $q$ formalizes the idea of variation across rays. When estimating a function $\phi(x)$, we reparametrize as $\phi(x)=\Phi(q,s)$ and decompose 
\[
|\phi(x)-\phi(x')|
=|\Phi(q,s)-\Phi(q',s')|
\leq |\Phi(q,s)-\Phi(q',s)| + |\Phi(q',s)-\Phi(q',s')|.
\]
The first term captures the variation across rays, which we often address by invoking Lebesgue points of $\Phi(\cdot,s)$ for fixed $s$. For the second term, we use the established regularity along rays. Such decompositions, and variations thereof, are used in integral estimates throughout the paper.

\textbf{Step 4: Transporting the remaining marginals on each ray.} So far, we have constructed a sub-coupling $\tilde{\gamma}_{\eps,1}$ whose marginals $(\mu_1,\nu_1)$ are balanced on each transport ray, i.e., the marginals $(\mu_{1,T},\nu_{1,T})$ induced on any ray $T$ have the same mass. It remains to couple the remaining masses $(\mu_2,\nu_2):=(\mu-\mu_1,\nu-\nu_1)$ in a way that creates sufficiently small transport cost and entropy cost, such as not to invalidate the upper bound~\eqref{ubb.estimates} when combined with $\tilde{\gamma}_{\eps,1}$. We start by focusing on the transport cost and recall that this cost will be optimized if we transport mass only along rays (and only downwards on each ray). Since both $(\mu,\nu)$ and $(\mu_1,\nu_1)$ are balanced on each ray, so are $(\mu_2,\nu_2)$, hence they can indeed be coupled by moving mass only along rays. For simplicity, let us assume that $\mu_{2,T}$ and $\nu_{2,T}$ are concentrated on disjoint intervals, so that the directional constraint can be ignored. (In the actual proof, the marginal supports $\mathcal{X}$ and $\mathcal{Y}$ are partitioned into $\{\mathcal{X}_k\}_{k=1}^K$ and $\{\mathcal{Y}_{k'}\}_{k'=1}^{K'}$, and we work separately on each pair $(\mathcal{X}_k,\mathcal{Y}_{k'})$. The construction is so that, for each pair, either the disjointness assumption holds, or no mass needs to be transported at all.) There are many choices to couple $(\mu_{2,T},\nu_{2,T})$ for each $T$, and they are equivalent in terms of transport cost. But with a view towards the entropy cost, a more diffuse choice is preferable. Since an explicit formula will be helpful for the fine estimates needed, we choose the product coupling of~$\mu_{2,T}$ and~$\nu_{2,T}$---it turns out to be sufficiently close to optimal for the entropy cost. Call $\bar{\gamma}_{\eps,2}$ the coupling of $(\mu_2,\nu_2)$ obtained by integrating these product measures over all transport rays.

\textbf{Step 5: Transporting the remaining marginals with controlled entropy.} In Step 4, we have coupled the remainder marginals $(\mu_2,\nu_2)$ in a way that is optimal for the transport cost, but of course any transport concentrated on rays is singular and hence causes an infinite entropy cost. By contrast, we require a tight control of the entropy to guarantee the desired upper bound \eqref{ubb.estimates}. Specifically, the density $r_\eps$ of $\tilde{\gamma}_{\eps,2}$ w.r.t.\ $\mathcal{L}^d\otimes\mathcal{L}^d$ needs to satisfy
\begin{equation}\label{entropy.control}
   \lim_{\eps\to0^+} \int_{\mathbb{R}^d\times\mathbb{R}^d}r_\eps(x,y)\max\big\{\!\log\bigl(\eps^{(d-1)/ 2}r_\eps(x,y)\bigr),0\big\}dxdy =0.
\end{equation}
To construct $r_\eps$, the basic idea is to spread the mass of the ray-wise transport $\bar{\gamma}_{\eps,2}$ in the proximity of each transport ray, without changing the marginals $(\mu_2,\nu_2)$. The standard smoothing tool for couplings is the \emph{block approximation technique} popularized by \cite{CarlierDuvalPeyreSchmitzer.17}; in general, this technique approximates a given (often singular) coupling $\gamma\in\Pi(\mu,\nu)$ by an absolutely continuous one $\hat\gamma$ as follows. Partition $\mathcal{X}$ into blocks $(A_i)$ and $\mathcal{Y}$ into blocks $(B_j)$. Next, replace each restriction $\gamma|_{A_i\times B_j}$ by the product of its own marginals, and define $\hat\gamma$ as the sum of these product measures. Then $\hat\gamma$ has the same marginals $(\mu,\nu)$ as $\gamma$ and $\hat\gamma\ll\mu\otimes\nu$. Typically, $A_i$ and $B_j$ are squares with a common side length $s$ which is optimized to achieve a particular tradeoff.

However, such a block approximation is not sharp enough to achieve~\eqref{entropy.control}. To attain this tight control, smoothening must be implemented at different scales in the direction of the ray and the orthogonal ones (namely, $\eps$ and $\sqrt{\eps}$). Of course, different rays have different directions, putting this outside the scope of the usual block approximation. Instead, we propose a bi-level, \emph{anisotropic} refinement of the block approximation technique: First, we partition $\mathcal{X}$ into squares $(A_i)$ and $\mathcal{Y}$ into squares $(B_j)$, all of sidelength $\sqrt{\eps}$. Next, for any fixed pair $(i,j)$, we further partition both $A_i$ and $B_j$ into slices $A_{i,l}$ and $B_{j,k}$ of thickness $\eps$, where the slices are cut orthogonally to the line joining the centers of the squares $A_i$ and $B_j$ (see \cref{{fig:two-layer-anisotropic-blocks}}). Now, we replace each $\gamma|_{A_{i,l}\times B_{j,k}}$ by the product of its own marginals. The idea behind this approximation is that the \emph{line joining the centers acts as a proxy for a transport ray:} since the initial kernel transports mass along rays, the approximation roughly achieves our goal of diffusing at scale $\eps$ in the parallel direction and scale $\sqrt{\eps}$ in the orthogonal direction, \emph{simultaneously} for all rays and exactly preserving the marginals~$(\mu_2,\nu_2)$. The rigorous version of this construction can be found in \Cref{Sect.4.2.2} and the proof of \eqref{entropy.control} is completed in Section~\ref{Sect.4.5}.

\begin{figure}[t]
\centering
\begin{tikzpicture}[
  line cap=round, line join=round,
  >={Latex[length=2.1mm]},
  every node/.style={font=\normalsize},
  gridline/.style={black!45, line width=0.35pt},
  gridborder/.style={black!60, line width=0.80pt},
  hiA/.style={fill=green!55!black, fill opacity=0.22, draw=none},
  hiB/.style={fill=cyan!55!blue,  fill opacity=0.22, draw=none},
  hiR/.style={fill=red!75!black,  fill opacity=0.22, draw=none},
  centerline/.style={green!55!black, dashed, line width=0.95pt},
  pt/.style={circle, fill=black, inner sep=1.6pt},
  slice/.style={draw=red!75!black, line width=1.05pt}
]

\def\cell{0.62}        
\def\gap{1.35}         
\def\sepLR{1.60}       

\node at (3.3,4.5) {\bf{Level 1}};

\coordinate (X0) at (0.00,0.70); 

\def\ai{1}\def\aj{1}
\coordinate (AiSW) at ($(X0)+(\ai*\cell,\aj*\cell)$);
\path[hiA] (AiSW) rectangle ++(\cell,\cell);
\coordinate (AiC) at ($(AiSW)+(0.5*\cell,0.5*\cell)$);

\draw[gridborder] (X0) rectangle ++(4*\cell,4*\cell);
\foreach \i in {1,2,3}{
  \draw[gridline] ($(X0)+(\i*\cell,0)$) -- ($(X0)+(\i*\cell,4*\cell)$);
}
\foreach \j in {1,2,3}{
  \draw[gridline] ($(X0)+(0,\j*\cell)$) -- ($(X0)+(4*\cell,\j*\cell)$);
}

\node at ($(X0)+(2*\cell,-0.40)$) {$\mathcal X$};

\node[pt] at (AiC) {};

\draw[<->, black!70, line width=0.8pt]
  ($(X0)+(-0.30,1*\cell)$) -- ($(X0)+(-0.30,2*\cell)$)
  node[midway, left=2pt] {$\sqrt{\varepsilon}$};

\coordinate (Y0) at ($(X0)+(4*\cell+\gap,0.31)$);

\def\bi{2}\def\bj{1}
\coordinate (BjSW) at ($(Y0)+(\bi*\cell,\bj*\cell)$);
\path[hiB] (BjSW) rectangle ++(\cell,\cell);
\coordinate (BjC) at ($(BjSW)+(0.5*\cell,0.5*\cell)$);

\draw[gridborder] (Y0) rectangle ++(5*\cell,3*\cell);
\foreach \i in {1,2,3,4}{
  \draw[gridline] ($(Y0)+(\i*\cell,0)$) -- ($(Y0)+(\i*\cell,3*\cell)$);
}
\foreach \j in {1,2}{
  \draw[gridline] ($(Y0)+(0,\j*\cell)$) -- ($(Y0)+(5*\cell,\j*\cell)$);
}

\node at ($(Y0)+(2.5*\cell,-0.40)$) {$\mathcal Y$};

\node[pt] at (BjC) {};

\draw[centerline] (AiC) -- (BjC);

\begin{scope}[xshift=\sepLR cm]

\node at (10.6,4.5) {\bf{Level 2}};

\coordinate (ZA0) at (8.05,0.55);
\coordinate (ZB0) at (11.15,1.55);
\def\Z{2.10}

\draw[gridborder, line width=0.85pt] (ZA0) rectangle ++(\Z,\Z);
\path[hiA] (ZA0) rectangle ++(\Z,\Z);
\node at ($(ZA0)+(\Z/2,-0.35)$) {$A_i$};

\draw[gridborder, line width=0.85pt] (ZB0) rectangle ++(\Z,\Z);
\path[hiB] (ZB0) rectangle ++(\Z,\Z);
\node at ($(ZB0)+(\Z/2,-0.35)$) {$B_j$};

\coordinate (ZAcenter) at ($(ZA0)+(\Z/2,\Z/2)$);
\coordinate (ZBcenter) at ($(ZB0)+(\Z/2,\Z/2)$);
\node[pt] at (ZAcenter) {};
\node[pt] at (ZBcenter) {};

\pgfmathsetmacro{\dxL}{3.10}
\pgfmathsetmacro{\dyL}{1.00}
\pgfmathsetmacro{\lenL}{sqrt(\dxL*\dxL + \dyL*\dyL)}
\pgfmathsetmacro{\ux}{\dxL/\lenL}
\pgfmathsetmacro{\uy}{\dyL/\lenL}
\pgfmathsetmacro{\nx}{-\uy}
\pgfmathsetmacro{\ny}{\ux}

\def\breakGap{0.55}   
\def\capLen{0.14}     

\coordinate (midLine) at ($(ZAcenter)!0.5!(ZBcenter)$);
\coordinate (brL) at ($(midLine)+(-0.5*\breakGap*\ux,-0.5*\breakGap*\uy)$);
\coordinate (brR) at ($(midLine)+( 0.5*\breakGap*\ux, 0.5*\breakGap*\uy)$);

\draw[centerline] (ZAcenter) -- (brL);
\draw[centerline] (brR) -- (ZBcenter);

\draw[green!55!black, line width=0.95pt]
  ($(brL)+(\capLen*\nx,\capLen*\ny)$) -- ($(brL)+(-\capLen*\nx,-\capLen*\ny)$);
\draw[green!55!black, line width=0.95pt]
  ($(brR)+(\capLen*\nx,\capLen*\ny)$) -- ($(brR)+(-\capLen*\nx,-\capLen*\ny)$);

\def\sep{0.51}

\pgfmathsetmacro{\tSecondRight}{1*\sep} 
\pgfmathsetmacro{\tThirdRight}{0*\sep}  

\begin{scope}
  \clip (ZB0) rectangle ++(\Z,\Z);

  \coordinate (pB2) at ($(ZBcenter)+(\tSecondRight*\ux,\tSecondRight*\uy)$);
  \coordinate (pB3) at ($(ZBcenter)+(\tThirdRight*\ux,\tThirdRight*\uy)$);

  \path[hiB, even odd rule]
    (ZB0) rectangle ++(\Z,\Z)
    ($(pB3)+(-6*\nx,-6*\ny)$) -- ($(pB3)+(6*\nx,6*\ny)$) --
    ($(pB2)+(6*\nx,6*\ny)$) -- ($(pB2)+(-6*\nx,-6*\ny)$) -- cycle;

  \path[hiR]
    ($(pB3)+(-6*\nx,-6*\ny)$) -- ($(pB3)+(6*\nx,6*\ny)$) --
    ($(pB2)+(6*\nx,6*\ny)$) -- ($(pB2)+(-6*\nx,-6*\ny)$) -- cycle;
\end{scope}

\foreach \t in {-1.53,-1.02,-0.51,0.00,0.51,1.02,1.53}{
  \begin{scope}
    \clip (ZA0) rectangle ++(\Z,\Z);
    \coordinate (pA) at ($(ZAcenter)+(\t*\ux,\t*\uy)$);
    \draw[slice] ($(pA)+(-6*\nx,-6*\ny)$) -- ($(pA)+(6*\nx,6*\ny)$);
  \end{scope}
  \begin{scope}
    \clip (ZB0) rectangle ++(\Z,\Z);
    \coordinate (pB) at ($(ZBcenter)+(\t*\ux,\t*\uy)$);
    \draw[slice] ($(pB)+(-6*\nx,-6*\ny)$) -- ($(pB)+(6*\nx,6*\ny)$);
  \end{scope}
}

\pgfmathsetmacro{\outN}{1.08} 
\coordinate (e0base) at ($(ZBcenter)+(\outN*\nx,\outN*\ny)$);
\coordinate (e1base) at ($(e0base)+(\sep*\ux,\sep*\uy)$);

\coordinate (epsShift) at (-.55mm,2mm);
\coordinate (e0) at ($(e0base)+(epsShift)$);
\coordinate (e1) at ($(e1base)+(epsShift)$);

\draw[<->, red!75!black, line width=0.9pt] (e0) -- (e1)
  node[midway, above=1pt] {$\!\!\!\varepsilon$};

\end{scope}

\end{tikzpicture}
\caption{Two-level anisotropic block approximation. Level 1 partitions $\mathcal X$ and $\mathcal Y$ into squares $A_i,B_j$ of side length $\sqrt{\varepsilon}$. For each fixed $(i,j)$, level 2 slices $A_i$ and $B_j$ orthogonally to the line joining their centers, producing slabs of thickness~$\varepsilon$. The line serves as a proxy for a transport ray.}
\label{fig:two-layer-anisotropic-blocks}
\end{figure}
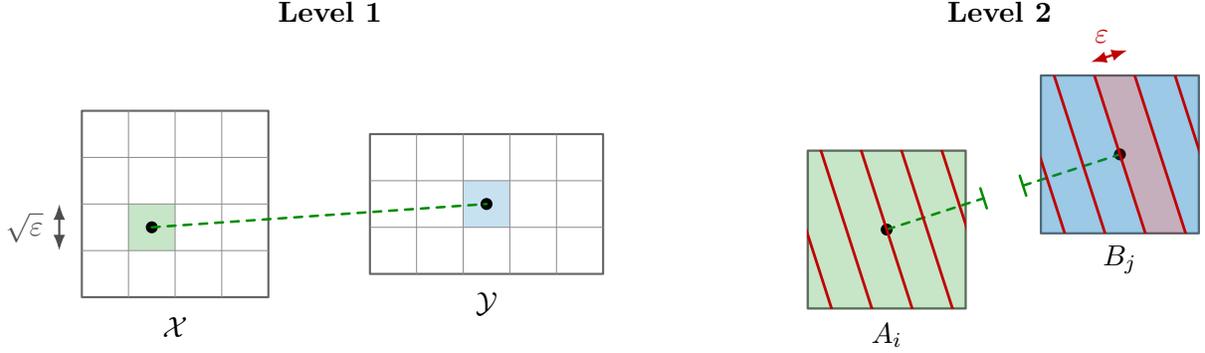

\begin{remark}[Overlapping marginals]\label{rem:overlap}
We expect that adaptations of the techniques developed here can also treat the case of overlapping marginal supports (i.e., $\mathcal X\cap\mathcal Y\neq\emptyset$) which we have excluded in \cref{assump1}. This case gives rise to several new features. Most notably, transport rays may degenerate to singletons, called \emph{zero-length rays}, so that optimal transports do not move the mass on those points. Hence, we define the entropic Monge coupling $\gamma_{0}$ as before on rays of positive length, and as the identity map on rays of zero length. The cost asymptotics are changed, even at the first order: for example, for $\mu=\nu=\mathrm{Unif}[0,1]^d$, direct computation gives leading coefficient $d$ rather than $(d-1)/2$. Heuristically, $d-1$ arises from Gaussian fluctuations in the $(d-1)$-dimensional orthogonal hyperplane $O_T$ along a positive-length ray, while $d$ corresponds to full-space fluctuations around zero-length rays, and the missing factor $1/2$ reflects the kink of $\|x-y\|$ at $x=y$. To extend our results to the overlapping case, we expect two additional ingredients to be essential. First, a regularity theory describing the behavior of the Kantorovich potential $u$ near zero-length rays, in the spirit of \eqref{eq:strategy2ndOrderGoal}; and second, a multi-level refinement of our current bi-level block approximation, capable of handling the lack of a uniform lower bound on $\|x-y\|$ when $x\in\mathcal{X}$ and $y\in\mathcal{Y}$. In view of these significant technical complications and the length of the present document, we plan to address the overlapping case in future work. 
\end{remark}

\subsection{Related literature}\label{Sect.1.3}

In view of the enormous literature on entropic optimal transport, our review focuses on studies about the limit $\eps\to0$ of vanishing regularization, as well as transport rays.

For discrete problems, the study of this limit goes back to~\cite{CominettiSanMartin.94}; see also~\cite{Weed.18} for a non-asymptotic result and \cite{AltschulerNilesWeedStromme.21} for a semi-discrete problem. In those settings, $\EOT_{\eps} - \OT = O(\eps)$. This linear rate is characteristic of problems that admit an optimal transport with finite relative entropy \cite[Proposition~A.1]{EcksteinNutz.22}. Here, we are mainly interested in continuous problems, where no such optimal transport exists and hence the rates are slower than linear. The link between $\EOT_{\eps}$ and $\OT$ goes back to~\cite{Mikami.02, Mikami.04} in the Schr\"odinger bridge setting, which can be seen as a dynamic version of entropic optimal transport with quadratic cost. Gamma convergence was shown in~\cite{Leonard.12}  and early quantitative results for quadratic cost, from a large deviations viewpoint, are \cite{AdamsDirrPeletierZimmer.11,DuongLaschosRenger.13,ErbarMaasRenger.15}. A stochastic control viewpoint is presented in~\cite{ChenGeorgiouPavon.16}. In a setting closer to ours, \cite{CarlierDuvalPeyreSchmitzer.17} showed Gamma convergence using the now-standard block approximation technique. 

Next, we review in more detail the existing \emph{first-order} results for EOT. For $C^{1,1}$ cost functions satisfying a nondegeneracy condition (invertibility of the mixed derivative $D^{2}_{xy}c(x,y)$), \cite{CarlierPegonTamanini.22} showed that $\EOT_{\eps} - \OT = -\frac{d}{2} \eps\log\eps + O(\varepsilon)$. A similar result, under different technical conditions, was shown in \cite{EcksteinNutz.22}. The same rate was later found in \cite{MalamutSylvestre.25}, which additionally disentangled the contributions of the transport and the divergence costs, and showed convergence at rate $\sqrt\eps$ of the optimal couplings in $W_2$ metric. Clearly, the distance cost considered in the present work does \emph{not} satisfy the nondegeneracy condition. For general Lipschitz costs, \cite{CarlierPegonTamanini.22,EcksteinNutz.22} showed $\EOT_{\eps} - \OT \leq -d \eps\log\eps + O(\varepsilon)$. This rate is sharp for the distance cost when the two marginals are identical ($\mu=\nu$). By contrast, the present work shows that the exact first-order term is $-\frac{d-1}{2}\eps\log\eps$ for the distance cost and non-overlapping marginals. Thus, the convergence is faster than in either of the previously known cases. Given the previous results, a possible intuition for the factor $\frac{d-1}{2}$ is the following: On the one hand, the  factor $1/2$ relates to the regularity of $\|x-y\|$ for $x\neq y$, similarly to the $C^{1,1}$ case in \cite{CarlierPegonTamanini.22}. On the other hand, the factor $d-1$ (instead of $d$) relates to the fact that diffusing mass along the direction of a transport ray does not lead to additional transport costs; additional costs only occur for the $d-1$ orthogonal dimensions.

A \emph{second-order} expansion of the optimal value was obtained in \cite{ConfortiTamanini.19} for the Schr\"odinger bridge setting. Building on that result, \cite{Chizat2020Faster} provided a second-order expansion for EOT with quadratic cost. Independently, \cite{Pal.19} obtained a second-order expansion for EOT with a more general class of cost functions modeled on the quadratic. In all these works, the optimal transport is unique (and given by a map), so the question of selection is trivial. However, the proof idea in \cite{Pal.19} is related to the present work: the key step therein is that for small regularization parameter, the conditional distribution of the optimal coupling is well approximated by a Gaussian distribution. 

A related literature studied the convergence of the optimal dual potentials (which however play a minor role in the present work). Based on compactness arguments, \cite{GigliTamanini.21, NutzWiesel.21} showed (subsequential) convergence to Kantorovich potentials in the uniform and $L^1$ sense respectively, and \cite{ChiariniConfortiGrecoTamanini.22} showed the convergence of the gradients in $L^2$ for quadratic cost. Again for quadratic cost, quantitative results for the convergence of the potentials and their gradients in $L^2$ were recently established in \cite{PooladianWeed.21,CarlierPegonTamanini.22,MalamutSylvestre.25}, and a rate in $L^\infty$ was recently established by \cite{LopezRivera.25}. Closely related to the convergence of the potentials is the large deviations principle of \cite{BerntonGhosalNutz.21} describing the convergence of the optimal couplings, again limited to situations where the optimal transport is unique.

Next, we review the previous literature on the entropic selection problem with distance cost, again focusing on continuous marginals. The first result is due to \cite{di2018entropic} and considers marginals in dimension $d=1$ with compact support and regular densities. A distinct feature of $d=1$ is that optimal transport plans, despite being concentrated on transport rays, can have support of dimension $d+d$, meaning that they are not singular. The authors show that the domain can be decomposed into components, and on each component, the limiting optimal transport plan behaves like a product measure. In particular, the limit is not singular (and has finite relative entropy) for non-overlapping marginal supports. When the marginals overlap, there can be mass which is not moved by any optimal transport plan, and that mass will remain fixed also under the selected limit (which is therefore singular). After removing this fixed mass, however, the remaining part is nonsingular. The proof proceeds by gamma convergence arguments. Recently, \cite{ley2025decomposition} generalized the result of \cite{di2018entropic} to less regular marginals, using a different approach based on stochastic orders, a theorem of Kellerer, and geometric properties related to the cyclical invariance of EOT \cite{BerntonGhosalNutz.21}. 

Turning to the present case $d>1$, we mention that the ideas and results of the present work have been disseminated for several years in increasing generality, for instance in \cite{Nutz2023EntropicSelectionOWPS} or \cite{Nutz2024EntropicSelectionIMSI}. More recently, the preprint \cite{aryan2025entropic} took up the same problem of entropic selection that is solved in \cref{Cor1.1}. The result in \cite{aryan2025entropic} has  similar assumptions and provides a (partial) local description of the selected limit, which is however different from the one in \cref{Cor1.1}. First, the local optimization problem in \cite[Theorem~1.1]{aryan2025entropic} involves a cost function that has a factor $e^{c\|x-y\|}$ with an unidentified constant $c\geq0$ that ``depends on $\mu,\nu$ and the transport ray'' \cite{aryan2025entropic}. We observe that the cost function would coincide with the one in \cref{Cor1.1} if $c=0$ for all rays. On the other hand, it seems to be (implicitly) assumed in \cite{aryan2025entropic} that $c\neq 0$ for all rays, for otherwise the key step in the proof of \cite[Theorem~1.1]{aryan2025entropic}, namely the application of \cite[Lemma~4.2]{aryan2025entropic}, seems to fail. Second, the optimization problem in \cite[Theorem~1.1]{aryan2025entropic} is an unconstrained EOT problem over a transport ray, in contrast to the constrained problem in \cref{Cor1.1}. As seen in \cref{sec:counterex_directional_constraint}, the solution of the unconstrained problem need not be an optimal transport, and that would  preclude it from being the limit of the EOT optimizers. In terms of arguments, \cite{aryan2025entropic} follows a substantially different approach, working primarily with the dual potentials and appealing to the stability result of \cite{ghosal2022stability} based on cyclical invariance.

Our work makes extensive use of transport rays, which originate in the existence theory for optimal maps in Monge's problem. The seminal work \cite{MR530375} outlined a strategy based on disintegrating $\mu$ and $\nu$ along the affine regions singled out by a Kantorovich potential, thus anticipating the modern decomposition along transport rays. However, making this fully rigorous required further advances; see \cite{BianchiniDaneri2018} for a comprehensive and modern treatment. Using partial differential equations, \cite{MR1464149} proved existence of an optimal transport map assuming Lipschitz marginal densities with disjoint supports. Their analysis introduces the transport density  and interprets transport rays as characteristics for limiting PDEs obtained through $p$-Laplace-type approximations. Another approach is approximation of the cost function. Indeed, \cite{MR1862796} constructed optimal maps for strictly convex perturbations of the distance cost and passed to the limit, thereby producing an optimal map for the original cost under minimal assumptions on the marginal densities. Independently, \cite{TW} obtained a similar existence result via approximation and PDE methods. These works highlighted a selection problem different from the one in the present work: the limiting Monge map can depend on the chosen approximation. The lecture notes \cite{MR2011032} (and the subsequent CIME notes \cite{ambrosio2003existence}) provide an alternative proof strategy emphasizing Gamma convergence and stability of minimizers together with a detailed analysis of the measure-theoretic properties of transport rays. The latter results are essential for the present work, which also contributes new regularity results where needed for our construction of the EOT limit.

\subsection{Notation}\label{Sect.1.6}

We denote by $\langle \cdot,\cdot\rangle$ the inner product on $\mathbb{R}^d$ and by $e_1,e_2,\cdots,e_d$ the standard basis of $\mathbb{R}^d$. For any $v\in\mathbb{R}^d\backslash\{0\}$, set $O(v):=\{x\in\mathbb{R}^d:\langle x,v\rangle=0\}$. For any $m\in\mathbb{N}^{*}$, $x_0\in \mathbb{R}^m$, and $r>0$, we define $B_m(x_0,r):=\{x\in\mathbb{R}^m:\|x-x_0\|\leq r\}$. For any $A,B\subseteq\mathbb{R}^d$, we define 
\begin{equation}\label{dfrakdef}
    \dist(A,B):=\inf_{x\in A,y\in B}\|x-y\|.
\end{equation} 
For any $x,x'\in\mathbb{R}^d$, we define $[x,x']:=\{tx+(1-t)x':0\leq t\leq 1\}$. If $x\neq x'$, we define $[x,x'):=\{tx+(1-t)x':0< t\leq 1\}$, $(x,x']:=\{tx+(1-t)x':0\leq t< 1\}$, and $(x,x'):=\{tx+(1-t)x':0< t<1\}$; if $x=x'$, we define $[x,x')=(x,x']=(x,x'):=\emptyset$. We also define $\mathrm{int}([x,x'])=\mathrm{int}([x,x'))=\mathrm{int}((x,x']):=(x,x')$.

We set $[m]:=\{1,2,\cdots,m\}$ for any $m\in\mathbb{N}^*$ and $[0]:=\emptyset$. For any $m\in\mathbb{N}^{*}$, we denote by $\mathbf{I}_m$ the $m\times m$ identity matrix. For a matrix $A$, we denote by $\|A\|_2$ the spectral norm of $A$. For any $s\in\mathbb{R}$, we define $s_{+}:=\max\{s,0\}$ and $s_{-}:=\max\{-s,0\}$. For any $s,t\in\mathbb{R}$, we define $s\wedge t:=\min\{s,t\}$ and $s\vee t:=\max\{s,t\}$. For any $x\geq 0$, we define
\begin{equation}\label{Phi_def}
    \Phi(x):= x\log x, \qquad \tilde{\Phi}(x):=x(\log x)_{+}.
\end{equation}  
For any $m\in\mathbb{N}^{*}$, we denote by $\mathcal{H}^m$ the $m$-dimensional Hausdorff measure on $\mathbb{R}^d$. For any two measures $\mu_1,\mu_2$ on the same measurable space, we write $\mu_1\ll \mu_2$ if $\mu_1$ is absolutely continuous with respect to $\mu_2$. 

As $\mathcal{X},\mathcal{Y}$ are compact and disjoint, there exist $D\geq 1$ and $d_0\in (0,1]$ such that
\begin{equation}\label{domainconst}
    \mathcal{X},\mathcal{Y}\subseteq B_d(0,D), \qquad \dist(\mathcal{X},\mathcal{Y}) \geq 10 d_0.
\end{equation}
As $f,g$ are Lipschitz functions with compact supports, there exists $M_0>0$ such that for every $x,x'\in\mathbb{R}^d$,
\begin{equation}\label{bddLip}
    f(x)\leq M_0, \quad g(x)\leq M_0, \quad |f(x)-f(x')|\leq M_0\|x-x'\|, \quad |g(x)-g(x')|\leq M_0\|x-x'\|.
\end{equation}

We denote by $C,c$ positive constants that depend only on the dimension $d$ and the fixed marginals $\mu,\nu$. The values of these constants may change from line to line.

\subsection{Organization of the paper}\label{Sect.1.5}

The remainder of this paper is organized as follows. Section \ref{Sect.2.1} reviews notation and background, and introduces a decomposition of the domains $\mathcal{X}$ and $\mathcal{Y}$. Proofs for Section \ref{Sect.2.1} are deferred to \cref{proof_Section2}. In Section \ref{Sect.2}, we establish key properties of the Kantorovich potential and its transport rays. Building on these results, we present the proofs of Theorems \ref{th:lowerBound} and \ref{th:upperBound} in Sections \ref{Sec3} and \ref{Sec4}, respectively. \Cref{sec:counterex} concludes with counterexamples. 

\section{Background on transport rays and domain decomposition}\label{Sect.2.1}

This section has two parts. Section~\ref{background_transport} collects notation and background on transport rays. Section~\ref{sectdom} introduces a decomposition of the domains $\mathcal{X}$,  $\mathcal{Y}$ that will be used in the proofs of the main results and establishes several related technical lemmas. The proofs of all results in this section are deferred to Appendix~\ref{proof_Section2}.

\subsection{Background on transport rays}\label{background_transport}

In this subsection, we review notation and background on transport rays, largely following~\cite{MR2011032, MR1862796}. 
Let $u\in \mathrm{Lip}_1(\mathbb{R}^d)$ be the Kantorovich potential fixed in Section \ref{Sect.1.1}.

\begin{definition}[Transport rays]\label{TransportRays}
A \emph{transport ray} is a line segment with endpoints $a,b\in \mathbb{R}^d$ such that the following properties hold:
\begin{itemize}
    \item[(a)] $a\in \mathcal{X},b\in \mathcal{Y}, a\neq b$.
    \item[(b)] $u(a)-u(b)=\|a-b\|$.
    \item[(c)] The segment $[a,b]$ is maximal for properties (a) and (b), i.e., for any $t>0$ such that $a+t(a-b)\in \mathcal{X}$, we have $|u(a+t(a-b))-u(b)|<\|a+t(a-b)-b\|$, and for any $t>0$ such that $b+t(b-a) \in \mathcal{Y}$, we have $|u(b+t(b-a))-u(a)|<\|b+t(b-a)-a\|$.
\end{itemize}
We then call $a$ the \emph{upper end} of the transport ray and $b$ the \emph{lower end} of the transport ray. For any transport ray $T$, we denote by $a_T,b_T$ the upper and lower ends of $T$, by $\mathrm{int}(T)$ the relative interior, and by $L_T$ the line containing $T$. The \emph{direction} of~$T$ and its orthogonal complement are 
\begin{equation*}
    V_T:=\frac{a_T-b_T}{\|a_T-b_T\|},\qquad O_T:=\{x\in\mathbb{R}^d:\langle x, V_T\rangle=0 \}.
\end{equation*}
We denote by $\mathcal{S}$ the set of transport rays. For any two transport rays $T,T'$, we define the distance between $T$ and $T'$ to be $\sqrt{\|a_T-a_{T'}\|^2+\|b_T-b_{T'}\|^2}$. Hereafter, we equip $\mathcal{S}$ with this metric. 

We denote by $\mathcal{T}_1$ the set of points that lie on some transport ray, by $\mathcal{T}_1^{*}$ the set of points that lie in the relative interior of some transport ray, and by $\mathcal{E}\subseteq \mathcal{T}_1$ the set of endpoints of transport rays. For $x\in\mathcal{T}_1^{*}$, we denote by $T(x)$ the (unique, by \cite[Lemma 10]{MR1862796}) transport ray containing~$x$, and by $L(x)$ the line containing $T(x)$.
\end{definition}

\begin{definition}[Transport set]\label{Def2.7}
A set $A\subseteq\mathbb{R}^d$ is called a \textit{transport set} if for any $x\in A\cap\mathcal{T}_1^{*}$, we have $\mathrm{int}(T(x))\subseteq A$. 
\end{definition}

\begin{definition}[Rays of length zero]\label{Def1.3}
We define 
\begin{equation*}
    \mathcal{T}_0:=\big\{z\in\mathcal{X}\cap\mathcal{Y}: |u(z)-u(z')|<\|z-z'\|\text{ for any } z'\in\mathcal{X}\cup\mathcal{Y}, z'\neq z\big\},
\end{equation*}
and refer to this set as the ``rays of length zero''. 
\end{definition}

By \cite[Lemma 25]{MR1862796} and the comments immediately following that lemma, the sets $\mathcal{T}_0$, $\mathcal{T}_1$, and $\mathcal{E}$ are Borel sets. Hence $\mathcal{T}_1^{*}=\mathcal{T}_1\backslash\mathcal{E}$ is also a Borel set. Moreover, as we assume $\mathcal{X}\cap\mathcal{Y}=\emptyset$, we have $\mathcal{T}_0=\emptyset$. Hence by \cite[Lemmas 9 and 25]{MR1862796},
\begin{equation}\label{new1.1}
    \mathcal{X}\cup\mathcal{Y}\subseteq \mathcal{T}_0\cup\mathcal{T}_1=\mathcal{T}_1^{*}\cup\mathcal{E} \quad \mbox{and} \quad  \mathcal{L}^d(\mathcal{E})=0.
\end{equation}

\begin{definition}[The functions $\alpha$ and $\beta$]\label{def_alphabeta}
For any $z\in \mathbb{R}^d$, we define
\begin{equation*}
     \alpha(z):=\sup\{\|z-x\|:x\in \mathcal{X}, u(x)-u(z)=\|z-x\|\},
\end{equation*}
\begin{equation*}
    \beta(z):=\sup\{\|z-y\|:y\in \mathcal{Y}, u(z)-u(y)=\|z-y\|\},
\end{equation*}
where $\sup\emptyset:=-\infty$. By \cite[Lemma 24]{MR1862796}, $\alpha,\beta:\mathbb{R}^d\rightarrow \mathbb{R}\cup\{-\infty\}$ are upper semicontinuous, hence Borel measurable. 
\end{definition}

By the discussion immediately following \cite[Lemma 24]{MR1862796} and the fact that $\mathcal{T}_0=\emptyset$, we have the following interpretation of the functions $\alpha$ and $\beta$: for $x\in\mathcal{T}_1^*$, $\alpha(x)$ and $\beta(x)$ represent the distances from $x$ to the upper and lower ends of $T(x)$, respectively; for $x\in\mathcal{E}$, we have $\alpha(x)=0$ or $\beta(x)=0$; for $x\in\mathbb{R}^d\backslash \mathcal{T}_1$, we have $\alpha(x)=-\infty$ or $\beta(x)=-\infty$.

In the following definition, we use the notion of \textit{approximate derivative} as discussed in \cite[Section~6.1.3]{MR3409135}. For a function $\varphi:\mathbb{R}^d\rightarrow\mathbb{R}^m$ (where $m\in\mathbb{N}^{*}$), we denote by $\approxgrad \varphi(x)$ the approximate derivative of $\varphi$ at $x$ (if it exists).

\begin{definition}[The functions $V$ and $F$]\label{DefVF}
For any $x\in\mathbb{R}^d$, we define
\begin{equation*}
    V(x):=\begin{cases}
        \nabla u(x) & \text{if }\nabla u(x)\text{ exists }\\
        0 & \text{otherwise}.
    \end{cases} 
\end{equation*}
By \cite[Lemma 12]{MR1862796} (see also \cite[Proposition 4.2]{MR2011032}), for any $x\in \mathcal{T}_1^{*}$, $\nabla u(x)$ exists and is equal to the direction of the transport ray containing $x$ (in particular, $\|V(x)\|=1$). By the discussion immediately above \cite[Theorem 4.3]{MR2011032}, $\approxgrad V(x)$ exists for $\mathcal{L}^d$-a.e.\ $x\in\mathcal{T}_1^{*}$. For any $x\in\mathbb{R}^d$, we define
\begin{equation*}
    F(x):=\begin{cases}
        \approxgrad V(x) & \text{if }\approxgrad V(x) \text{ exists }\\
        0 & \text{otherwise}.
    \end{cases}
\end{equation*}
Both $V:\mathbb{R}^d\rightarrow \mathbb{R}^d$ and $F:\mathbb{R}^d\rightarrow\mathbb{R}^{d\times d}$ are Borel measurable (see \cite[p.\,177]{MR2135805}).
\end{definition}

\begin{definition}[The functions $a$ and $b$]\label{Def2.6}
For any $x\in\mathbb{R}^d$, we define
\begin{equation*}
    a(x):=\begin{cases}
        x+\alpha(x) V(x) & \text{if }x\in\mathcal{T}_1^{*}\\
        x & \text{otherwise},
    \end{cases}\qquad b(x):=\begin{cases}
        x-\beta(x)V(x) & \text{if }x\in\mathcal{T}_1^{*}\\
        x & \text{otherwise}.
    \end{cases}
\end{equation*}
Note that for $x\in\mathcal{T}_1^{*}$, $a(x)$ and $b(x)$ are the upper and lower ends of $T(x)$. The functions $a(\cdot),b(\cdot)$ are Borel measurable as $\alpha(\cdot),\beta(\cdot),V(\cdot)$ are Borel measurable and $\mathcal{T}_1^{*}$ is a Borel set.
\end{definition}

We record two straightforward measurability lemmas for ease of reference.

\begin{lemma}\label{L2.0nn}
The mapping $T:\mathcal{T}_1^{*}\rightarrow \mathcal{S}$ is Borel measurable.
\end{lemma}

\begin{lemma}\label{Le2.2}
Define $\EuScript{I}:\mathcal{S}\rightarrow\mathbb{R}^d\times\mathbb{R}^d$ by $\EuScript{I}(T)=(a_T,b_T)$. Then $\EuScript{I}(\mathcal{S})$ is a Borel subset of $\mathbb{R}^d\times\mathbb{R}^d$. Consequently, $\mathcal{S}$ equipped with its Borel $\sigma$-algebra is a standard Borel space.
\end{lemma}

The following lemma is a strengthening of \cite[Lemma 16]{MR1862796}.

\begin{lemma}\label{L2.0}
For any $x_1,x_2\in\mathcal{T}_1^*$, we have 
\begin{equation*}
    \|V(x_1)-V(x_2)\|\leq \frac{4\|x_1-x_2\|}{\min\{\alpha(x_1),\beta(x_1),\alpha(x_2),\beta(x_2)\}}.
\end{equation*}
\end{lemma}

The following result is from \cite[Proposition 5.2]{MR1464149}.

\begin{lemma}\label{L2.1}
Let $B\subseteq\mathbb{R}^d$ be Lebesgue measurable with $\mathcal{L}^d(B)=0$. Then for $\mathcal{L}^d$-a.e.\ $z_0\in\mathcal{T}_1^{*}$, we have $\mathcal{H}^1(T(z_0)\cap B) =0$. 
\end{lemma}

\subsection{Domain decomposition and related properties}\label{sectdom}

In this subsection, we introduce a decomposition scheme for the domains $\mathcal{X}$ and $\mathcal{Y}$ and establish several related lemmas. The main goal of the decomposition is to localize so that source and target domains can be separated by a hyperplane. This enables a parametrization of any $x\in T$ by the intersection point $q$ of the ray $T$ with the hyperplane and the signed distance from $x$ to the hyperplane along the ray, which will allow us to decompose various quantities of interest into transverse and parallel components.

As $\inf_{x\in\mathcal{X},y\in\mathcal{Y}}\|x-y\|\geq 10 d_0$  by \eqref{domainconst}, there exist disjoint Borel sets $\{\mathcal{X}_k\}_{k\in [K]}$ and $\{\mathcal{Y}_{k'}\}_{k'\in [K']}$ (where $K,K'\in\mathbb{N}^{*}$ are of constant order) as follows:
\begin{itemize}
    \item[(a)] There exist $\mathtt{h}\in (0,d_0/ (10\sqrt{d}) )$ (of constant order), and points $\mathtt{x}_k=(\mathtt{x}_{k,1},\cdots,\mathtt{x}_{k,d})$, $\mathtt{y}_{k'}=(\mathtt{y}_{k',1},\cdots,\mathtt{y}_{k',d})$ in $\mathbb{R}^d$, such that for every $k\in[K]$ and $k'\in[K']$,
    \begin{equation}\label{eq:mathtthDefn}
    \mathcal{X}_k=\prod_{\ell=1}^d(\mathtt{x}_{k,\ell}-\mathtt{h},\mathtt{x}_{k,\ell}+\mathtt{h}], \qquad \mathcal{Y}_{k'}=\prod_{\ell=1}^d(\mathtt{y}_{k',\ell}-\mathtt{h},\mathtt{y}_{k',\ell}+\mathtt{h}].
    \end{equation}
    \item[(b)] We have $\mathcal{X}\subseteq\bigcup_{k=1}^K\mathcal{X}_k$ and $\mathcal{Y}\subseteq\bigcup_{k'=1}^{K'}\mathcal{Y}_{k'}$.
    \item[(c)] For every $k\in [K]$, $\mathcal{X}_k\cap\mathcal{X}\neq\emptyset$; for every $k'\in[K']$, $\mathcal{Y}_{k'}\cap\mathcal{Y}\neq\emptyset$. 
    \item[(d)] For every $k\in [K],k'\in[K']$, $\mathcal{X}_{k}\cap\mathcal{Y}_{k'}=\emptyset$, and there exists a hyperplane
    \begin{equation}\label{Hkkdef}
        H_{k,k'}=\{z\in\mathbb{R}^d:\langle z, \mathtt{a}_{k,k'}\rangle = \mathtt{b}_{k,k'}\}
    \end{equation}
    (where $\mathtt{a}_{k,k'}\in\mathbb{R}^d$, $\|\mathtt{a}_{k,k'}\|=1$, and $\mathtt{b}_{k,k'}\in\mathbb{R}$), such that 
    \begin{equation}\label{bddakk}
        \langle z,\mathtt{a}_{k,k'}\rangle \geq \mathtt{b}_{k,k'}+2d_0,\text{ for every }z\in\mathcal{X}_k; \langle z,\mathtt{a}_{k,k'}\rangle \leq \mathtt{b}_{k,k'}-2d_0\text{ for every } z\in\mathcal{Y}_{k'}.
    \end{equation}
\end{itemize}
We fix such a choice of $\{\mathcal{X}_k\}_{k\in[K]}$, $\{\mathcal{Y}_{k'}\}_{k'\in [K']}$, and $\{H_{k,k'}\}_{k\in[K],k'\in[K']}$ throughout the rest of this paper. Note that by \eqref{domainconst}, the above conditions imply that for every $k\in[K], k'\in[K']$,
\begin{equation}\label{distcondi}
   \mathcal{X}_k,\mathcal{Y}_{k'}\subseteq B_d(0,D+1), \quad 
    \dist(\mathcal{X}_k,H_{k,k'})\geq 2d_0, \quad \dist(\mathcal{Y}_{k'},H_{k,k'})\geq 2d_0.
\end{equation}
For every $k\in[K]$, we denote by $\mathcal{X}_k^{\circ}$ the interior of $\mathcal{X}_k$; for every $k'\in[K']$, we denote by $\mathcal{Y}_{k'}^{\circ}$ the interior of $\mathcal{Y}_{k'}$. For each $k\in[K],k'\in[K']$, we extend $\mathtt{a}_{k,k'}$ to an orthonormal basis of $\mathbb{R}^d$: $\big\{\mathtt{a}_{k,k'},\mathtt{e}_{k,k'}^{(1)},\cdots,\mathtt{e}_{k,k'}^{(d-1)}\big\}$. Note that any $q\in H_{k,k'}$ satisfies $q=\mathtt{b}_{k,k'}\mathtt{a}_{k,k'}+\sum_{\ell=1}^{d-1}\langle q, \mathtt{e}_{k,k'}^{(\ell)}\rangle \mathtt{e}_{k,k'}^{(\ell)}$. For any $q\in H_{k,k'}$, we set
\begin{equation}\label{qexpre}
    \Omega_{k,k'}(q):=\Big(\langle q, \mathtt{e}_{k,k'}^{(1)}\rangle,\cdots,\langle q, \mathtt{e}_{k,k'}^{(d-1)}\rangle\Big).
\end{equation}
Note that $\Omega_{k,k'}:H_{k,k'}\rightarrow \mathbb{R}^{d-1}$ is a bijective mapping, with
\begin{equation}\label{inverse_map}
    \Omega_{k,k'}^{-1}(x)=\mathtt{b}_{k,k'}\mathtt{a}_{k,k'}+\sum_{\ell=1}^{d-1}x_{\ell}\mathtt{e}_{k,k'}^{(\ell)},\quad\text{ for all }\quad  x\in\mathbb{R}^{d-1}.
\end{equation}

Now for every $k\in[K],k'\in[K']$, we define
\begin{equation}\label{def_T11}
    \mathcal{T}_{1;k,k'}^*:=\big\{x\in\mathcal{T}_1^*: a(x)\in\mathcal{X}_k,b(x)\in\mathcal{Y}_{k'}\big\}, 
\end{equation}
where $a(\cdot),b(\cdot),T(\cdot)$ are as in Definitions~\ref{TransportRays} and~\ref{Def2.6}. As $a(\cdot),b(\cdot)$ are Borel measurable, $\mathcal{T}_{1;k,k'}^*$ is a Borel set. Moreover,
\begin{equation}\label{disj}
    \{\mathcal{T}_{1;k,k'}^*\}_{k\in[K],k'\in[K']}\text{ are mutually disjoint}, \quad \text{ and }\quad \mathcal{T}_1^{*}=\bigcup_{k\in[K],k'\in[K']}\mathcal{T}_{1;k,k'
    }^{*}.
\end{equation}
For every $k\in[K],k'\in[K']$, we define
\begin{equation}\label{deftildehat}
   \bar{H}_{k,k'}:=H_{k,k'}\cap \mathcal{T}_{1;k,k'}^{*}  ,  \quad \tilde{H}_{k,k'}:=\big\{q\in H_{k,k'}\cap\mathcal{T}_1^*:\langle V(q),\mathtt{a}_{k,k'}\rangle\neq 0,\min\{\alpha(q),\beta(q)\}\geq d_0\big\}. 
\end{equation}
Note that $\bar{H}_{k,k'}$ and $\tilde{H}_{k,k'}$ are Borel subsets of $H_{k,k'}$, and $\bar{H}_{k,k'}\subseteq \tilde{H}_{k,k'}$ (using \eqref{distcondi}). As $\mathcal{X},\mathcal{Y}\subseteq B_d(0,D)$ (see \eqref{domainconst}), from Definitions \ref{TransportRays} and \ref{def_alphabeta} we have
\begin{equation}\label{bddsalphabeta}
    \mathcal{T}_1\subseteq B_d(0,D),\quad\text{and}\quad \max\{\alpha(x),\beta(x)\}\leq 2D \quad\text{for every }x\in\mathcal{T}_1. 
\end{equation}
Note that this implies  $\bar{H}_{k,k'}\subseteq \tilde{H}_{k,k'}\subseteq B_d(0,D)$, hence
\begin{equation}\label{BddHtildekk}
    \mathcal{H}^{d-1}(\bar{H}_{k,k'})\leq \mathcal{H}^{d-1}(\tilde{H}_{k,k'})\leq C. 
\end{equation}

We record two measurability lemmas for later reference. Given a subset $A$ of the separating hyperplane, \cref{L2.1n} states that the set of points on rays passing through $A$ is measurable, and \cref{L2.4n} states that the set of rays themselves is measurable. Recall that all proofs are deferred to Appendix~\ref{proof_Section2}.

\begin{lemma}\label{L2.1n}
For any $k\in[K],k'\in[K']$ and any Borel set $A\subseteq H_{k,k'}$, the set
\begin{equation}\label{defWkk}
    \mathfrak{W}_{k,k';A}:=\big\{q+tV(q):q\in A\cap\mathcal{T}_1^*,\langle V(q),\mathtt{a}_{k,k'}\rangle\neq 0,-\beta(q)<t<\alpha(q)\big\}\subseteq \mathcal{T}_1^{*}
\end{equation}
is Borel measurable.
\end{lemma}

\begin{lemma}\label{L2.4n}
For any $k\in[K],k'\in[K']$ and any Borel set $A\subseteq H_{k,k'}$, the set
\begin{equation}\label{defskkA}
    \mathcal{S}_{k,k';A}:=\big\{T\in\mathcal{S}:a_T\in\mathcal{X}_k,b_T\in\mathcal{Y}_{k'},T \text{ intersects }H_{k,k'}\text{ at some }q\in A\big\} 
\end{equation}
is a Borel subset of $\mathcal{S}$ (recall that the metric on $\mathcal{S}$ is given as in Definition \ref{TransportRays}).
\end{lemma}

Next, we fix $k\in[K],k'\in[K']$. Define (recall \eqref{deftildehat} and Definition \ref{TransportRays})
\begin{equation}\label{deftildeTkk}
    \tilde{\mathcal{T}}^{*}_{1;k,k'}:=\big\{x\in\mathcal{T}_1^{*}:T(x)\text{ intersects }H_{k,k'}\text{ at some }q\in\tilde{H}_{k,k'}\big\}.
\end{equation}
From the fact that $\bar{H}_{k,k'}\subseteq \tilde{H}_{k,k'}$, we have 
\begin{equation}\label{subs}
    \mathcal{T}^{*}_{1;k,k'}\subseteq \tilde{\mathcal{T}}^{*}_{1;k,k'}. 
\end{equation}

Our next goal is a disintegration of the Lebesgue measure along transport rays. To that end, some additional notation is necessary. 
As $\tilde{H}_{k,k'}$ is a Borel subset of $H_{k,k'}$, by Lemma \ref{L2.1n}, $\tilde{\mathcal{T}}^{*}_{1;k,k'}=\mathfrak{W}_{k,k';\tilde{H}_{k,k'}}$ is a Borel subset of $\mathbb{R}^d$. Note that for any $q\in \tilde{H}_{k,k'}$, $\min\{\alpha(q),\beta(q)\}\geq d_0$. Hence by Lemma \ref{L2.0}, for any $q,q'\in \tilde{H}_{k,k'}$, we have $\|V(q)-V(q')\|\leq 4d_0^{-1}\|q-q'\|$. By Kirszbraun's theorem, $V|_{\tilde{H}_{k,k'}}$ can be extended to $\tilde{V}_{k,k'}:H_{k,k'}\rightarrow\mathbb{R}^d$ such that $\|\tilde{V}_{k,k'}\|_{\mathrm{Lip}}\leq 4d_0^{-1}$ (where $\|\cdot\|_{\mathrm{Lip}}$ is the Lipschitz norm). Without loss of generality, we may assume that $\|\tilde{V}_{k,k'}(q)\|\leq 1$ for all $q\in H_{k,k'}$ (since otherwise we can replace $\tilde{V}_{k,k'}$ by its composition with the projection onto the closed unit ball without increasing the Lipschitz norm). Now for any $q\in H_{k,k'}$ and $t\in\mathbb{R}$, we define
\begin{equation}\label{Psikkp}
    \Psi_{k,k'}(q,t):=q+t\tilde{V}_{k,k'}(q). 
\end{equation}
Note that for any $q,q'\in H_{k,k'}$ and $t,t'\in\mathbb{R}$ such that $\max\{|t|,|t'|\}\leq 8D$,
\begin{align*}
   &\|\Psi_{k,k'}(q,t)-\Psi_{k,k'}(q',t')\|\leq \|q-q'\|+|t|\|\tilde{V}_{k,k'}(q)-\tilde{V}_{k,k'}(q')\|+\|\tilde{V}_{k,k'}(q')\||t-t'|\nonumber\\
   \leq\, & (1+32Dd_0^{-1})\|q-q'\|+|t-t'|\leq 2(1+32Dd_0^{-1})\sqrt{\|q-q'\|^2+(t-t')^2}.
\end{align*}
Hence by Kirszbraun's theorem, $\Psi_{k,k'}|_{H_{k,k'}\times [-8D,8D]}$ can be extended to a Lipschitz mapping $\tilde{\Psi}_{k,k'}:H_{k,k'}\times \mathbb{R}\rightarrow\mathbb{R}^d$ with $\|\tilde{\Psi}_{k,k'}\|_{\mathrm{Lip}}\leq 2(1+32Dd_0^{-1})\leq C$. For any $x\in\mathbb{R}^{d-1}$ and $t\in\mathbb{R}$, we set (recall \eqref{qexpre})
\begin{equation}\label{defEP}
    \EuScript{P}(x,t):=\tilde{\Psi}_{k,k'}\big(\Omega_{k,k'}^{-1}(x),t\big).
\end{equation} 
For any $(q,t)\in H_{k,k'}\times\mathbb{R}$, we set 
\begin{equation}\label{defRk}
    R_{k,k'}(q,t):=\begin{cases}
          J\EuScript{P}(\Omega_{k,k'}(q),t), & \text{ if } J\EuScript{P}(\Omega_{k,k'}(q),t)\text{ exists}\\
          0 & \text{ otherwise},
    \end{cases}
\end{equation}
where $J\EuScript{P}(\cdot)$ is the Jacobian of $\EuScript{P}(\cdot)$ (see \cite[Section 3.2.2]{MR3409135}). Note that $\EuScript{P}:\mathbb{R}^d\rightarrow\mathbb{R}^d$ is Lipschitz continuous with $\|\EuScript{P}\|_{\mathrm{Lip}}=\|\tilde{\Psi}_{k,k'}\|_{\mathrm{Lip}}\leq  C$, which implies that for every $(q,t)\in H_{k,k'}\times\mathbb{R}$, 
\begin{equation}\label{Rkkb}
    0\leq R_{k,k'}(q,t) \leq C.
\end{equation}
Moreover, by Rademacher's theorem (see, e.g., \cite[Theorem 3.2]{MR3409135}), $J\EuScript{P}(\Omega_{k,k'}(q),t)$ exists for $\mathcal{H}^{d-1}\otimes\mathcal{L}^1$-a.e.\ $(q,t)\in H_{k,k'}\times\mathbb{R}$.    

With this notation in place, the following lemma gives a disintegration of the Lebesgue measure restricted to $\tilde{\mathcal{T}}^{*}_{1;k,k'}$. An integral in $\R^d$ can be performed by first integrating along each ray and then summing over rays by integrating over the separating hyperplane.

\begin{lemma}\label{L2.3n}
For any $k\in[K],k'\in[K']$ and any nonnegative Lebesgue measurable function $\varphi:\mathbb{R}^d\rightarrow [0,\infty)$, we have 
\begin{equation}\label{E2.3n.1}
\int_{\tilde{\mathcal{T}}^{*}_{1;k,k'}}\varphi(x)dx =\int_{\tilde{H}_{k,k'}}d\mathcal{H}^{d-1}(q) \int_{\mathbb{R}} \mathbbm{1}_{(-\beta(q),\alpha(q))}(t) 
    \varphi(q+tV(q)) R_{k,k'}(q,t) dt.
\end{equation}
\end{lemma}

The next lemma shows that the marginals are balanced (induce the same mass) along each ray.

\begin{lemma}\label{L2.5n}
For any $k\in[K],k'\in[K']$, the following holds. For $\mathcal{H}^{d-1}$-a.e.\ $q\in \tilde{H}_{k,k'}$, we have
    \begin{equation}\label{L2.5.e1}
        \int_{\mathbb{R}}\mathbbm{1}_{(-\beta(q),\alpha(q))}(t)f(q+tV(q))R_{k,k'}(q,t)dt=\int_{\mathbb{R}}\mathbbm{1}_{(-\beta(q),\alpha(q))}(t)g(q+tV(q))R_{k,k'}(q,t)dt.
    \end{equation}
In particular, for $\mathcal{L}^d$-a.e.\ $z_0\in \tilde{\mathcal{T}}_{1;k,k'}^*$, if we set $q$ to be the intersection point of $T(z_0)$ and $H_{k,k'}$ (note that $q\in\tilde{H}_{k,k'}$), then \eqref{L2.5.e1} holds.
\end{lemma}

For future reference, we record the integrability of the log-distance to the ray end.

\begin{lemma}\label{Lemma2.7}
We have
\begin{equation*}
    \bigg|\int_{\mathcal{T}_1^{*}} f(x)\log(\alpha(x))dx\bigg|\leq C, \qquad \bigg|\int_{\mathcal{T}_1^{*}} f(x)\log(\beta(x))dx\bigg|\leq C,
\end{equation*}
and similarly for $g$ instead of $f$.
\end{lemma}

Lastly, we also bound the left tail behavior of the distance to the ray end. 

\begin{lemma}\label{Lem4.4}
For any $s\geq 0$, we have
\begin{equation*}
    \mathcal{L}^d(\{x\in\mathcal{T}_1^{*}:\alpha(x)\leq s\})\leq Cs, \qquad \mathcal{L}^d(\{x\in\mathcal{T}_1^{*}:\beta(x)\leq s\})\leq Cs.
\end{equation*}
\end{lemma}

\section{Properties of the Kantorovich potential and its transport rays}\label{Sect.2}

This section establishes differentiability properties of the Kantorovich potential~$u$ as well as continuity properties of the functions $\alpha,\beta$ which measure the distance of a point to the end of its transport ray (Definition~\ref{def_alphabeta}). We recall the functions $V,F$ which proxy for the first and second derivatives of~$u$ (Definition~\ref{DefVF}). To state the results, we introduce the following sets capturing the regularity of the functions $\alpha,\beta,V$ in directions orthogonal to the transport ray containing $x$.

\begin{definition}\label{Defd}
For any $x\in \mathcal{T}_1^{*}$ and $t,\delta>0$, we define
\begin{equation*}
    \mathcal{D}(x;t,\delta):=\big\{w\in\mathbb{R}^d:\langle w,V(x) \rangle=0, \|w\|\leq t,\max\{|\alpha(x+w)-\alpha(x)|,|\beta(x+w)-\beta(x)|\}>\delta\big\},
\end{equation*}
\begin{equation*}
    \mathcal{D}'(x;t,\delta):=\big\{w\in\mathbb{R}^d:\langle w,V(x) \rangle=0, \|w\|\leq t,\max\{|\alpha(x+w)-\alpha(x)|,|\beta(x+w)-\beta(x)|\}\leq\delta\big\}, 
\end{equation*}
\begin{equation*}
    \mathcal{W}(x;t,\delta):=\big\{w\in\mathbb{R}^d:\langle w,V(x) \rangle=0, \|w\|\leq t,\|V(x+w)-V(x)-F(x) w\|>\delta t\big\},
\end{equation*}
\begin{equation*}
    \mathcal{W}'(x;t,\delta):=\big\{w\in\mathbb{R}^d:\langle w,V(x) \rangle=0, \|w\|\leq t,\|V(x+w)-V(x)-F(x) w\|\leq\delta t\big\},
\end{equation*}
where the functions $\alpha,\beta,V,F$ are as in Definitions~\ref{def_alphabeta} and~\ref{DefVF}.
\end{definition}

Intuitively, the set $\mathcal{D}'(x;t,\delta)$ consists of the points in directions orthogonal to the transport ray containing $x$ where $\alpha,\beta$ are close to being continuous, and the set $\mathcal{W}'(x;t,\delta)$ consists of those points where the first-order approximation $V(x+w)\approx V(x)+F(x) w$ holds.  

We also recall the standard notion of approximate continuity (see \cite[Definition~1.29]{MR3409135}). A function $\varphi:\mathbb{R}^d\rightarrow[-\infty,\infty)$ is said to be \textit{approximately continuous} at $x\in\mathbb{R}^d$ if for any $\delta>0$,
\begin{equation*}
    \lim_{t\rightarrow 0^+}\frac{\mathcal{L}^d(\{y\in B_d(x,t):|\varphi(y)-\varphi(x)|\geq \delta\})}{t^d}=0. 
\end{equation*}

Our first theorem states two results. First, the distances $\alpha$ and $\beta$ to the ray ends are approximately continuous at every interior point of a ray. The novelty here, which will be crucial in subsequent proofs, is the validity at \emph{every} point, rather than almost every point (which would be a straightforward result). Second, the theorem establishes \emph{transverse} approximate continuity in the sense that the average is taken only over $(d-1)$-dimensional balls orthogonal to the ray.

\begin{theorem}[Continuity of $\alpha$ and $\beta$]\label{L2.2}
For $\mathcal{L}^d$-a.e.\ $z_0\in\mathcal{T}_1^*$, the following holds for $T=T(z_0)$: for every $x\in\mathrm{int}(T)$, the functions $\alpha,\beta$ are approximately continuous at $x$, and 
\begin{equation}\label{Eq2.1n}
    \lim_{t\rightarrow 0^{+}}\frac{\mathcal{H}^{d-1}(\mathcal{D}(x;t,\delta))}{t^{d-1}} = 0,\quad\text{for all }\delta>0.
\end{equation}
\end{theorem}

The next theorem collects several key properties of the Kantorovich potential~$u$ and its (approximate) derivatives at all points in the interior of a transport ray. Part~(a) states the existence of the approximate derivative~$F$ of~$V$ (corresponding to the second derivative of $u$) in the usual sense as well as in the transverse sense. Part~(b) derives a bound and Lipschitz regularity of $F$ along rays. In addition, it contains two structural properties of~$F$: orthogonality of its rows (and columns) to the ray direction $V$, and the transformation formula~\eqref{Eq2.6.4} that will be crucial to derive~\eqref{eq:Lambda-relation-in-overview}. Finally, (c) states the symmetry of~$F$. (Symmetry will be derived as a by-product of the other arguments as the applicability of standard results is not obvious.)

\begin{theorem}[Properties of $V$ and $F$]\label{L2.4}
For $\mathcal{L}^d$-a.e.\ $z_0\in\mathcal{T}_1^*$, the following holds for $T=T(z_0)$:
\begin{itemize}
    \item[(a)] For every $x\in\mathrm{int}(T)$, 
    \begin{equation}\label{Eq2.11n}
        \approxgrad V(x)\text{ exists } (\text{hence }F(x)=\approxgrad V(x)),
    \end{equation}
\begin{equation}\label{Eq2.12}
   \lim_{t\rightarrow 0^{+}}\frac{\mathcal{H}^{d-1}(\mathcal{W}(x;t,\delta))}{t^{d-1}}=0,\quad \text{for all } \delta>0.
\end{equation} 
\item[(b)] For every $x\in\mathrm{int}(T)$,
\begin{equation}\label{L2.4.enew3}
     \|F(x)\|_2 \leq \frac{4}{\min\{\alpha(x),\beta(x)\}},
\end{equation}
\begin{equation}\label{Eq2.6.5}
    F(x)V(x)=F(x)^{\top}V(x)=0.
\end{equation}
For every $x,x'\in\mathrm{int}(T)$,
\begin{equation}\label{L2.4.enew4}
    \|F(x)-F(x')\|_2\leq \frac{128\|x-x'\|}{\min\{\alpha(x),\beta(x)\}\min\{\alpha(x'),\beta(x')\}},
\end{equation}
\begin{equation}\label{Eq2.6.2}
    F(x') =  \big(\mathbf{I}_d-\langle x'-x, V_T\rangle F(x')\big)F(x).
\end{equation} 
For every $x,x',x_0\in\mathrm{int}(T)$,
\begin{equation}\label{Eq2.6.4}
    \mathbf{I}_d-\langle x'-x_0, V_T\rangle F(x')=\big(\mathbf{I}_d-\langle x'-x, V_T\rangle F(x')\big)\big(\mathbf{I}_d-\langle x-x_0, V_T\rangle F(x)\big). 
\end{equation}
\item[(c)] For every $x\in\mathrm{int}(T)$, $F(x)^{\top} = F(x)$.  
\end{itemize}
\end{theorem}
\begin{remark}\label{Rmk.L2.4}
For ease of reference, we record the special case of~\eqref{Eq2.6.4} where $x_0=x'$:
\begin{equation*}
    \big(\mathbf{I}_d-\langle x'-x, V_T\rangle F(x')\big)\big(\mathbf{I}_d-\langle x-x', V_T\rangle F(x)\big)=\mathbf{I}_d \quad \text{for all $x,x'\in\mathrm{int}(T)$.}
\end{equation*}
\end{remark}

The last theorem is the approximate second-order expansion of $u$, making rigorous the heuristic~\eqref{eq:strategy2ndOrderGoal}. It is a key ingredient for the proof of the variational lower bound in \cref{Sec3}.

\begin{theorem}[Second-order expansion of $u$]\label{L2.6}
For $\mathcal{L}^d$-a.e.\ $z_0\in\mathcal{T}_1^*$, the following holds for $T=T(z_0)$: for every $x\in \mathrm{int}(T)$, 
    \begin{equation*}
        \lim_{t\rightarrow 0^{+}} \frac{\mathcal{H}^{d-1}\big(\big\{w\in O_{T}:\|w\|\leq t,\big|u(x+w)-u(x)-\frac{1}{2}w^{\top} F(x) w\big|>\delta t^2\big\}\big)}{t^{d-1}}=0,\quad \text{ for all }\delta>0.
   \end{equation*}
    
\end{theorem}

The remainder of this section is devoted to the proofs of the three theorems.  

\subsection{Proof of Theorem \ref{L2.2}}\label{Sect.3.1.2}

For any transport ray $T\in\mathcal{S}$, define $\mathcal{A}(T)$ to be the set of $x\in\mathrm{int}(T)$ such that \eqref{Eq2.1n} holds. We define $\tilde{\alpha},\tilde{\beta}:\mathbb{R}^d\rightarrow\mathbb{R}$ by
\begin{equation}\label{defa}
    \tilde{\alpha}(x):=\begin{cases}
        \alpha(x), & \text{ if }x\in\mathcal{T}_1, \\
        0, & \text{ otherwise},  
    \end{cases} \qquad \tilde{\beta}(x):=\begin{cases}
        \beta(x), & \text{ if }x\in\mathcal{T}_1, \\
        0, & \text{ otherwise}. 
    \end{cases}
\end{equation}
Since $\mathcal{T}_1$ is a Borel set and $\alpha,\beta$ are Borel measurable (see Definition~\ref{def_alphabeta}), the functions $\tilde{\alpha},\tilde{\beta},\mathbbm{1}_{\mathcal{T}_1}(\cdot)$ are Borel measurable. Let $A$ be the set of $x\in\mathbb{R}^d$ such that $\tilde{\alpha},\tilde{\beta},\mathbbm{1}_{\mathcal{T}_1}(\cdot)$ are approximately continuous at $x$. By \cite[Theorem 1.37]{MR3409135}, we have $\mathcal{L}^d(\mathbb{R}^d\backslash A)=0$. Hence by Lemma~\ref{L2.1}, for $\mathcal{L}^d$-a.e.\ $z_0\in\mathcal{T}_1^*$, $T=T(z_0)$ satisfies
\begin{equation}\label{ESq2.2}
   \mathcal{H}^1(\mathrm{int}(T)\cap(\mathbb{R}^d\backslash A))= \mathcal{H}^1(T\cap (\mathbb{R}^d\backslash A))=0.
\end{equation}
In the following, we fix any $T\in\mathcal{S}$ such that \eqref{ESq2.2} holds.

\paragraph{Step 1.}

We first show that $\mathrm{int}(T)\cap A\subseteq \mathcal{A}(T)$, which by \eqref{ESq2.2} implies that
\begin{equation}\label{Esq2.3}
    \mathcal{H}^1(\mathrm{int}(T)\backslash \mathcal{A}(T))=0.
\end{equation}

In the following, we consider any $x\in \mathrm{int}(T)\cap  A$ and $\delta\in (0,\min\{\alpha(x),\beta(x),1\}\slash 100)$. As $\tilde{\alpha},\tilde{\beta},\mathbbm{1}_{\mathcal{T}_1}(\cdot)$ are approximately continuous at $x$, for any $\delta'>0$, there exists $t_0\in (0,\delta\slash 10)$ depending on $x,\delta,\delta'$, such that for any $t\in (0,t_0)$,
\begin{equation*}
    \mathcal{L}^d\big(B_d(x,2t)\cap (\mathbb{R}^d\backslash\mathcal{T}_1)\big)\leq \frac{1}{2}\delta' t^d,
\end{equation*}
\begin{equation*}
    \mathcal{L}^d\big(B_d(x,2t)\cap \big\{y\in\mathbb{R}^d:\max\{|\tilde{\alpha}(y)-\tilde{\alpha}(x)|,|\tilde{\beta}(y)-\tilde{\beta}(x)|\}\geq \delta\slash 2\big\}\big)\leq \frac{1}{2}\delta' t^d. 
\end{equation*}
Consequently, noting \eqref{defa}, for any $t\in (0,t_0)$, we have
\begin{align*}
&\mathcal{L}^d\big(B_d(x,2t)\cap \big\{y\in\mathbb{R}^d:\max\{|\alpha(y)-\alpha(x)|,|\beta(y)-\beta(x)|\}\geq \delta\slash 2\big\}\big)\nonumber\\
\leq\,& \mathcal{L}^d\big(B_d(x,2t)\cap \big\{y\in\mathbb{R}^d:\max\{|\tilde{\alpha}(y)-\tilde{\alpha}(x)|,|\tilde{\beta}(y)-\tilde{\beta}(x)|\}\geq \delta\slash 2\big\}\big)\nonumber\\
&+\mathcal{L}^d\big(B_d(x,2t)\cap (\mathbb{R}^d\backslash\mathcal{T}_1)\big)\leq \delta' t^d.
\end{align*}
Below, we fix any $t\in (0,t_0)$ and take $t' := t+\frac{32t^2}{\min\{\alpha(x),\beta(x)\}}\in [t,2t]$ (note that $t<t_0<\delta\slash 10$). By the above display, there exists $s\in [-t,t]$ such that
\begin{align}\label{Eq2.3ol}
  &\mathcal{H}^{d-1}\big(\big\{w\in O_T:\|w\|\leq t,\max\{|\alpha(x+sV_T+w)-\alpha(x)|,|\beta(x+sV_T+w)-\beta(x)|\}\geq \delta\slash 2\big\}\big)\nonumber\\
  \leq\,& \frac{\mathcal{L}^d\big(B_d(x,2t)\cap \big\{y\in\mathbb{R}^d:\max\{|\alpha(y)-\alpha(x)|,|\beta(y)-\beta(x)|\}\geq \delta\slash 2\big\}\big)}{2t}\leq \delta' t^{d-1}. 
\end{align}
Set $x':=x+s V_T\in\mathrm{int}(T)$. Note that
\begin{equation}\label{boundss}
    \|x'-x\|=|s|\leq t.
\end{equation}
As 
\begin{equation}\label{diff}
    \max\{|\alpha(x')-\alpha(x)|,|\beta(x')-\beta(x)|\}\leq |s|\leq t \leq \frac{\delta}{2},
\end{equation}
by the above display and recalling Definition~\ref{Defd}, we have $\mathcal{H}^{d-1}(\mathcal{D}(x';t,\delta))\leq \delta' t^{d-1}$. Hence
\begin{align}\label{Eq2.4}
    \mathcal{H}^{d-1}(\mathcal{D}'(x';t,\delta))&\geq \mathcal{H}^{d-1}(\{w\in O_T:\|w\|\leq t\})-\delta' t^{d-1}\nonumber\\
    &=\mathcal{L}^{d-1}(B_{d-1}(0,t))-\delta' t^{d-1} = \bigg(\frac{\pi^{(d-1)\slash 2}}{\Gamma((d+1) \slash 2)}-\delta'\bigg)t^{d-1}. 
\end{align}
For any $w\in\mathcal{D}'(x;t',2\delta)$, $\|w\|\leq t'\leq 2t$ and 
\begin{equation}\label{Eq2.6}
   \min\{\alpha(x+w),\beta(x+w)\}\geq \min\{\alpha(x),\beta(x)\}-2\delta\geq \frac{1}{2}\min\{\alpha(x),\beta(x)\}>0.
\end{equation}
Note that this implies $x+w\in\mathcal{T}_1^*$. Hence by Lemma \ref{L2.0}, we have
\begin{equation}\label{Eq2.8}
   \sqrt{2-2\langle V(x+w),V_T\rangle}=\|V(x+w)-V(x)\|\leq\frac{8\|w\|}{\min\{\alpha(x),\beta(x)\}}\leq \frac{16t}{\min\{\alpha(x),\beta(x)\}}\leq \frac{1}{5},
\end{equation}
\begin{equation}\label{eqnnn1.1}
  \text{and } \qquad \langle V(x+w),V_T\rangle\geq\frac{1}{2}.
\end{equation}
Now we define the mapping $\varphi:\mathcal{D}'(x;t', 2\delta)\rightarrow O_T$ by setting
\begin{equation}\label{defva}
      \varphi(w):=x+w-x'-\frac{\langle x-x', V_T \rangle}{\langle V(x+w), V_T\rangle}V(x+w) 
\end{equation}
for any $w\in\mathcal{D}'(x;t', 2\delta)$. Note that for any $w\in\mathcal{D}'(x;t',2\delta)$, by \eqref{boundss}, we have 
\begin{equation}\label{Eq2.7}
    \|(x'+\varphi(w))-(x+w)\|=\frac{|\langle x'-x, V_T \rangle|}{\langle V(x+w), V_T\rangle}\leq 2\|x'-x\|\leq 2t\leq \delta<\frac{\min\{\alpha(x),\beta(x)\}}{100},
\end{equation}
hence by \eqref{Eq2.6}, $\min\{\alpha(x+w),\beta(x+w)\}> \|(x'+\varphi(w))-(x+w)\|$. Therefore, noting \eqref{defva} and the fact that $\varphi(w)\in O_T$, we conclude that
\begin{equation}\label{conclun1}
    x'+\varphi(w)\text{ is the intersection of }x'+O_T\text{ and }T(x+w),\quad \text{for every }w\in \mathcal{D}'(x;t',2\delta).
\end{equation}

By the definition of $\mathcal{D}'(x';t,\delta)$, for any $w'\in\mathcal{D}'(x';t,\delta)$, we have 
\begin{align}\label{Eq2.6n}
     \min\{\alpha(x'+w'),\beta(x'+w')\}\geq\,& \min\{\alpha(x'),\beta(x')\}-\delta\nonumber\\
    \geq\,&\min\{\alpha(x),\beta(x)\}-\frac{3}{2}\delta \geq \frac{1}{2}\min\{\alpha(x),\beta(x)\}>0,
\end{align}
where we use \eqref{diff} in the second inequality. Note that this implies $x'+w'\in\mathcal{T}_1^*$. Hence by Lemma \ref{L2.0}, \eqref{Eq2.6n}, and the fact that $\min\{\alpha(x'),\beta(x')\}\geq \min\{\alpha(x),\beta(x)\}\slash 2$ (by \eqref{diff}), we have
\begin{equation}\label{Eq2.8n}
   \sqrt{2-2\langle V(x'+w'),V_T\rangle}=\|V(x'+w')-V(x')\|\leq\frac{8\|w'\|}{\min\{\alpha(x),\beta(x)\}}\leq \frac{8t}{\min\{\alpha(x),\beta(x)\}}\leq \frac{1}{5}, 
\end{equation}
and so 
\begin{equation}\label{lbddd}
    \langle V(x'+w'),V_T\rangle\geq \frac{1}{2}.
\end{equation}
Now we define $\psi:\mathcal{D}'(x';t,\delta)\rightarrow O_T$ by setting
\begin{equation}\label{pside}
       \psi(w'):=x'+w'-x-\frac{\langle x'-x, V_T \rangle}{\langle V(x'+w'), V_T\rangle}V(x'+w') 
\end{equation}
for any $w'\in \mathcal{D}'(x';t,\delta)$. Note that for any $w'\in \mathcal{D}'(x';t,\delta)$, by \eqref{boundss}, we have  
\begin{equation}\label{Eq2.7n}
     \|(x+\psi(w'))-(x'+w')\|=\frac{\|x'-x\|}{\langle V(x'+w'), V_T\rangle}\leq 2\|x'-x\|\leq 2t\leq \delta<\frac{\min\{\alpha(x),\beta(x)\}}{100}.
\end{equation}
Hence by \eqref{Eq2.6n}, $\min\{\alpha(x'+w'),\beta(x'+w')\}> \|(x+\psi(w'))-(x'+w')\|$, and so $x+\psi(w')$ is the intersection of $x+O_T$ and $T(x'+w')$. Using the fact that $x,x'$ are on the same transport ray and $x'+w',x+\psi(w')$ are on the same transport ray, we obtain
\begin{align*}
    & \max\{|\alpha(x+\psi(w'))-\alpha(x)|,|\beta(x+\psi(w'))-\beta(x)|\}\nonumber\\
    \leq\,& \max\{|\alpha(x'+w')-\alpha(x')|,|\beta(x'+w')-\beta(x')|\}+\big|\|(x+\psi(w'))-(x'+w')\|-\|x-x'\|\big|\nonumber\\
    \leq\,& \delta+\|x-x'\|\bigg|1-\frac{1}{\langle V(x'+w'), V_T\rangle}\bigg|
    =\delta + \frac{\|x-x'\||\langle V(x'+w')-V_T,V_T\rangle|}{\langle V(x'+w'), V_T\rangle}\nonumber\\\leq\,&\delta+2t\|V(x'+w')-V_T\|=\delta+2t\|V(x'+w')-V(x')\|\nonumber\\
    \leq\,&\delta+\frac{16 t^2}{\min\{\alpha(x),\beta(x)\}}\leq\delta+\frac{\delta^2}{\min\{\alpha(x),\beta(x)\}}\leq 2\delta,
\end{align*}
where the second inequality follows from the fact that $w'\in\mathcal{D}'(x';t,\delta)$ and the equality in \eqref{Eq2.7n}, and the first inequality in the last line uses \eqref{Eq2.8n}. Moreover, using \eqref{pside}, we obtain
\begin{align*}
   \|\psi(w')\|\leq\,& \|w'\|+\bigg\|x'-x-\frac{\langle x'-x, V_T \rangle}{\langle V(x'+w'), V_T\rangle}V(x'+w')\bigg\| \nonumber\\
   \leq\,& t+\bigg\|\langle x'-x,V_T\rangle V_T-\frac{\langle x'-x, V_T \rangle}{\langle V(x'+w'), V_T\rangle}V_T\bigg\|\nonumber\\
   & +\bigg\|\frac{\langle x'-x, V_T \rangle}{\langle V(x'+w'), V_T\rangle}V_T-\frac{\langle x'-x, V_T \rangle}{\langle V(x'+w'), V_T\rangle}V(x'+w')\bigg\|\nonumber\\
   \leq\,& t + \|x'-x\|\bigg|1-\frac{1}{\langle V(x'+w'), V_T\rangle}\bigg|+\frac{|\langle x'-x, V_T \rangle|}{\langle V(x'+w'), V_T\rangle}\|V(x'+w')-V_T\|\nonumber\\
   \leq\,& t+\|x'-x\|\frac{|\langle V(x'+w')-V_T,V_T\rangle|}{\langle V(x'+w'), V_T\rangle}+\|x'-x\|\frac{\|V(x'+w')-V_T\|}{\langle V(x'+w'),V_T\rangle}\nonumber\\
   \leq\,& t+4t\|V(x'+w')-V_T\|=t+4t\|V(x'+w')-V(x')\|\leq t+\frac{32t^2}{\min\{\alpha(x),\beta(x)\}}=t',
\end{align*}
where the two inequalities in the last line use \eqref{boundss}, \eqref{lbddd}, and \eqref{Eq2.8n}. By the above two displays, we have $\psi(w')\in \mathcal{D}'(x;t',2\delta)$. Consequently, by \eqref{conclun1}, $x'+\varphi(\psi(w'))$ is the intersection of $x'+O_T$ and $T(x+\psi(w'))=T(x'+w')$ (as $x'+w',x+\psi(w')$ are on the same transport ray), and so $w'=\varphi(\psi(w'))\in\varphi(\mathcal{D}'(x;t',2\delta))$. Therefore, we conclude that
\begin{equation}\label{Eq2.10}
   \mathcal{D}'(x';t,\delta) \subseteq \varphi(\mathcal{D}'(x;t',2\delta)).
\end{equation}

By \eqref{Eq2.6} and Lemma \ref{L2.0}, for any $w,w'\in\mathcal{D}'(x;t',2\delta)$, we have
\begin{equation*}
    \|V(x+w')-V(x+w)\|\leq \frac{8}{\min\{\alpha(x),\beta(x)\}}\|w-w'\|.
\end{equation*}
Hence
\begin{align*}
       &  \bigg\|\frac{V(x+w)}{\langle V(x+w), V_T\rangle}-\frac{V(x+w')}{\langle V(x+w'), V_T\rangle}\bigg\|\nonumber\\
       \leq\,& \bigg\|\frac{V(x+w)}{\langle V(x+w), V_T\rangle}- \frac{V(x+w')}{\langle V(x+w), V_T\rangle} \bigg\|+\bigg\|\frac{V(x+w')}{\langle V(x+w), V_T\rangle}-\frac{V(x+w')}{\langle V(x+w'), V_T\rangle}\bigg\|    \nonumber\\
       \leq\,&  \frac{\|V(x+w)-V(x+w')\|}{\langle V(x+w), V_T\rangle}
       +\frac{|\langle V(x+w')-V(x+w), V_T\rangle|}{\langle V(x+w), V_T\rangle  \langle V(x+w'), V_T\rangle}\nonumber\\
       \leq\,& 8\|V(x+w')-V(x+w)\|\leq \frac{64}{\min\{\alpha(x),\beta(x)\}}\|w-w'\|,
\end{align*}
where the second and third inequalities use \eqref{eqnnn1.1}. Consequently, by \eqref{defva} and \eqref{boundss}, we have  
\begin{align*}
    \|\varphi(w)-\varphi(w')\| 
   \leq\,& \|w-w'\|+|\langle x'-x, V_T \rangle|\bigg\|\frac{V(x+w)}{\langle V(x+w), V_T\rangle}-\frac{V(x+w')}{\langle V(x+w'), V_T\rangle}\bigg\|\nonumber\\
   \leq\,& \|w-w'\|+\frac{64t}{\min\{\alpha(x),\beta(x)\}}\|w-w'\|. 
\end{align*}
By Kirszbraun's theorem (see, e.g., \cite[Section 2.10.43]{MR0257325}), there exists a Lipschitz mapping $\tilde{\varphi}:O_T\rightarrow O_T$ such that $\tilde{\varphi}(w)=\varphi(w)$ for all $w\in\mathcal{D}'(x;t',2\delta)$ and $\|\tilde{\varphi}\|_{\mathrm{Lip}}\leq 1+\frac{64 t}{\min\{\alpha(x),\beta(x)\}}$. By \cite[Theorem 2.8]{MR3409135}, 
\begin{equation}\label{Eq2.11}
    \mathcal{H}^{d-1}(\varphi(\mathcal{D}'(x;t',2\delta)))=\mathcal{H}^{d-1}(\tilde{\varphi}(\mathcal{D}'(x;t',2\delta)))\leq\bigg(1+\frac{64 t}{\min\{\alpha(x),\beta(x)\}}\bigg)^{d-1}\mathcal{H}^{d-1}(\mathcal{D}'(x;t',2\delta)).
\end{equation}

By \eqref{Eq2.4}, \eqref{Eq2.10}, and \eqref{Eq2.11}, we have 
\begin{align*}
    \mathcal{H}^{d-1}(\mathcal{D}'(x;t',2\delta))\geq\, &\bigg(1+\frac{64t}{\min\{\alpha(x),\beta(x)\}}\bigg)^{-(d-1)}\mathcal{H}^{d-1}(\mathcal{D}'(x';t,\delta))\nonumber\\
    \geq\, & \bigg(1+\frac{64t}{\min\{\alpha(x),\beta(x)\}}\bigg)^{-(d-1)}\bigg(\frac{\pi^{(d-1)\slash 2}}{\Gamma((d+1) \slash 2)}-\delta'\bigg)t^{d-1}. 
\end{align*}
Hence
\begin{align*}
   &   \mathcal{H}^{d-1}(\mathcal{D}(x; t, 2\delta))\leq \mathcal{H}^{d-1}(\mathcal{D}(x;t',2\delta))=\frac{\pi^{(d-1)\slash 2}}{\Gamma((d+1) \slash 2)}(t')^{d-1}
   -\mathcal{H}^{d-1}(\mathcal{D}'(x;t',2\delta))\nonumber\\
   \leq\,& \frac{\pi^{(d-1)\slash 2}}{\Gamma((d+1) \slash 2)}(t')^{d-1}-\bigg(1+\frac{64t}{\min\{\alpha(x),\beta(x)\}}\bigg)^{-(d-1)}\bigg(\frac{\pi^{(d-1)\slash 2}}{\Gamma((d+1) \slash 2)}-\delta'\bigg)t^{d-1}.
\end{align*}
As $\delta'>0$ can be arbitrarily small, by the above display, we get
\begin{equation*}
    \lim_{t\rightarrow 0^{+}}\frac{\mathcal{H}^{d-1}(\mathcal{D}(x;t,2\delta))}{t^{d-1}}=0.
\end{equation*}
Hence $x\in\mathcal{A}(T)$. We conclude that $\mathrm{int}(T)\cap A\subseteq\mathcal{A}(T)$, as desired. 

\paragraph{Step 2.} 

Next, we show that $\mathcal{A}(T)$ is a closed subset of $\mathrm{int}(T)$ under the subspace topology. 

Let $x_n\in \mathcal{A}(T)$ for each $n\in\mathbb{N}^{*}$, and suppose that $x_n\rightarrow x$ for some $x\in\mathrm{int}(T)$. Consider any $\delta\in (0,\min\{\alpha(x),\beta(x),1\}\slash 100)$ and $\delta'\in (0,\delta)$. Note that there exists $n_0\in\mathbb{N}^{*}$ such that $\|x_n-x\|\leq \delta'$ for all $n\geq n_0$. Below we fix $n\geq n_0$. As $x_n\in\mathcal{A}(T)$, there exists $t_0\in (0,\delta\slash 10)$ such that $\mathcal{H}^{d-1}(\mathcal{D}(x_n;t,\delta))\leq \delta' t^{d-1}$ for all $t\in (0,t_0)$. Following a similar argument as in Step~1 (with $x_n$ playing the role of $x'$), we obtain that  
\begin{align*}
    \mathcal{H}^{d-1}(\mathcal{D}(x;t,2\delta))\leq\,& \mathcal{H}^{d-1}(\mathcal{D}(x;t'',2\delta))\nonumber\\
    \leq\,&\frac{\pi^{(d-1)\slash 2}}{\Gamma((d+1) \slash 2)}(t'')^{d-1}-\bigg(1+\frac{64\delta'}{\min\{\alpha(x),\beta(x)\}}\bigg)^{-(d-1)}\bigg(\frac{\pi^{(d-1)\slash 2}}{\Gamma((d+1) \slash 2)}-\delta'\bigg)t^{d-1},
\end{align*}
where $t'':=t+\frac{64 \delta' t}{\min\{\alpha(x),\beta(x)\}}\in [t,2t]$. As $\delta'>0$ can be arbitrarily small, by the above display, we get
\begin{equation*}
    \lim_{t\rightarrow 0^{+}}\frac{\mathcal{H}^{d-1}(\mathcal{D}(x;t,2\delta))}{t^{d-1}}=0.
\end{equation*}
Hence $x\in\mathcal{A}(T)$. We conclude that $\mathcal{A}(T)$ is a closed subset of $\mathrm{int}(T)$ under the subspace topology. By \eqref{Esq2.3}, $\mathrm{int}(T)\backslash \mathcal{A}(T)$ is an open subset of $\mathrm{int}(T)$ (under the subspace topology) with zero $\mathcal{H}^1$ measure, which implies $\mathrm{int}(T)\backslash \mathcal{A}(T)=\emptyset$ and $\mathcal{A}(T)=\mathrm{int}(T)$. Therefore, for any $x\in\mathrm{int}(T)$, \eqref{Eq2.1n} holds. 

\paragraph{Step 3.} Finally, we show that $\alpha,\beta$ are approximately continuous at every $x\in\mathrm{int}(T)$. 

Fix any $x\in\mathrm{int}(T)$ and $\delta\in (0,\min\{\alpha(x),\beta(x),1\}\slash 100)$. From the result of Step 2, \eqref{Eq2.1n} holds for $x$. Therefore, for any $\delta'>0$, there exists $t_0\in (0,\delta\slash 10)$ such that for all $t\in (0,t_0)$,
\begin{equation}
    \mathcal{H}^{d-1}(\mathcal{D}(x;t,\delta))\leq \delta' t^{d-1}. 
\end{equation}

Below we consider any $t\in (0,t_0)$. For any $x'\in\mathrm{int}(T)$ such that $\|x'-x\|\leq t$, following a similar argument as in Step~1 (with $x$ playing the role of $x'$ and $x'$ playing the role of $x$), we can deduce that 
\begin{align}\label{L2.2.e2}
   & \mathcal{H}^{d-1}(\mathcal{D}(x'; t, 2\delta))\nonumber\\
   \leq\,&\frac{\pi^{(d-1)\slash 2}}{\Gamma((d+1) \slash 2)}(t')^{d-1} -\bigg(1+\frac{64 t}{\min\{\alpha(x),\beta(x)\}}\bigg)^{-(d-1)}\bigg(\frac{\pi^{(d-1)\slash 2}}{\Gamma((d+1) \slash 2)}-\delta'\bigg)t^{d-1},
\end{align}
where $t' := t+\frac{32 t^2}{\min\{\alpha(x),\beta(x)\}}$. Moreover, note that for any $y\in B_d(x,t)$ such that $\max\{|\alpha(y)-\alpha(x)|,|\beta(y)-\beta(x)|\}> 3\delta$, if $y=x'+w$ (where $x'\in [x-tV(x),x+t V(x)]$ and $w\in O_T$), then $\|w\|\leq t$ and 
\begin{align}\label{L2.2.e1}
   & \max\{|\alpha(x'+w)-\alpha(x')|,|\beta(x'+w)-\beta(x')|\}\nonumber\\
   \geq\,&\max\{|\alpha(y)-\alpha(x)|,|\beta(y)-\beta(x)|\}-\max\{|\alpha(x')-\alpha(x)|,|\beta(x')-\beta(x)|\}\nonumber\\
   \geq\,&\max\{|\alpha(y)-\alpha(x)|,|\beta(y)-\beta(x)|\}-t\geq 3\delta-t> 2\delta. 
\end{align}
By \eqref{L2.2.e2} and \eqref{L2.2.e1}, we have  
\begin{align*}
   & \mathcal{L}^d\big(B_d(x,t)\cap\big\{y\in\mathbb{R}^d:\max\{|\alpha(y)-\alpha(x)|,|\beta(y)-\beta(x)|\}> 3\delta\big\}\big)\nonumber\\
   \leq\,& \int_{[x-tV(x),x+tV(x)]}\mathcal{H}^{d-1}(\mathcal{D}(x';t,2\delta))d\mathcal{H}^1(x')\nonumber\\
   \leq\,& 2t\Bigg(\frac{\pi^{(d-1)\slash 2}}{\Gamma((d+1) \slash 2)}(t')^{d-1}-\bigg(1+\frac{64 t}{\min\{\alpha(x),\beta(x)\}}\bigg)^{-(d-1)}\bigg(\frac{\pi^{(d-1)\slash 2}}{\Gamma((d+1) \slash 2)}-\delta'\bigg)t^{d-1}\Bigg).
\end{align*}
Taking $t\rightarrow 0^{+}$, and then taking $\delta'\rightarrow 0^{+}$, we obtain that
\begin{equation*}
    \lim_{t\rightarrow 0^{+}}\frac{\mathcal{L}^d\big(B_d(x,t)\cap\big\{y\in\mathbb{R}^d:\max\{|\alpha(y)-\alpha(x)|,|\beta(y)-\beta(x)|\}> 3\delta\big\}\big)}{t^d} = 0.
\end{equation*}
Hence $\alpha,\beta$ are approximately continuous at $x$.

\subsection{Proof of Theorem \ref{L2.4}(a)--(b)}\label{Sect.3.1.3}

We begin with the following lemma.

\begin{lemma}\label{L2.3}
For any $x\in\mathcal{T}_1^{*}$, if $\approxgrad V(x)$ exists and $\alpha,\beta$ are approximately continuous at $x$, then $F(x) V(x) = 0$. 
\end{lemma}
\begin{proof}

Let $x\in\mathcal{T}_1^*$ be such that $\approxgrad V(x)$ exists and $\alpha,\beta$ are approximately continuous at $x$. Fix an arbitrary $\delta_0\in (0,1)$, and for any $t>0$ set  
\begin{equation*}
    U_{\delta_0,t}:=\big\{x+s V(x)+w: t\slash 2\leq s\leq t, w\in O(V(x)), \|w\|\leq \delta_0 t   \big\}.
\end{equation*}
Note that
\begin{equation}\label{Usubl}
    U_{\delta_0,t}\subseteq B_d(x,2t), \qquad  \mathcal{L}^d(U_{\delta_0,t})\geq c_{\delta_0} t^d,
\end{equation}
where $c_{\delta_0}>0$ depends only on $\delta_0,d$. 

Now consider any $\delta\in (0,\min\{\alpha(x),\beta(x),1\}\slash 10)$. As $\approxgrad V(x)$ exists (note that this implies $\approxgrad V(x)=F(x)$; see Definition~\ref{DefVF}) and $\alpha,\beta$ are approximately continuous at $x$, there exists $t_0\in (0,\delta\slash 10)$, such that for any $t\in (0,t_0)$, we have
\begin{equation*}
   \mathcal{L}^d\big(B_d(x,2t)\cap\big\{y\in\mathbb{R}^d:\max\{|\alpha(y)-\alpha(x)|,|\beta(y)-\beta(x)|\}\geq\delta\big\}\big)\leq c_{\delta_0} t^d\slash 4, 
\end{equation*}
\begin{equation*}
    \mathcal{L}^d\big(B_d(x,2t)\cap\big\{y\in\mathbb{R}^d:\|V(y)-V(x)-F(x)(y-x)\|\geq\delta\|y-x\|\big\}\big)\leq c_{\delta_0} t^d\slash 4.
\end{equation*}
Below we fix any $t\in (0,t_0)$. By the above three displays, there exists $y_0=x+s_0V(x)+w_0\in U_{\delta_0,t}$ (where $t\slash 2\leq s_0\leq t$, $w_0\in O(V(x))$, and $\|w_0\|\leq\delta_0 t$) such that 
\begin{equation}\label{Eq2.17}
    \max\{|\alpha(y_0)-\alpha(x)|,|\beta(y_0)-\beta(x)|\}<\delta, \quad \|V(y_0)-V(x)-F(x)(y_0-x)\|<\delta\|y_0-x\|.
\end{equation}
As $\delta<\min\{\alpha(x),\beta(x),1\}\slash 10$, \eqref{Eq2.17} implies that $\min\{\alpha(y_0),\beta(y_0)\}\geq \min\{\alpha(x),\beta(x)\}\slash 2>0$, and in particular, $y_0\in \mathcal{T}_1^*$. Moreover, from $s_0\leq t<t_0<\delta\slash 10<\min\{\alpha(x),\beta(x)\}\slash 100$, we get $x+s_0V(x)\in \mathrm{int}(T(x))$,  and so $V(x+s_0V(x))=V(x)$ and
\begin{equation*}
   \min\{\alpha(x+s_0V(x)),\beta(x+s_0V(x))\}\geq\min\{\alpha(x),\beta(x)\}-s_0\geq \frac{\min\{\alpha(x),\beta(x)\}}{2}.
\end{equation*}
Hence by Lemma \ref{L2.0}, we have 
\begin{equation}\label{Eq2.16}
    \|V(y_0)-V(x)\|=\|V(y_0)-V(x+s_0V(x))\|\leq\frac{8}{\min\{\alpha(x),\beta(x)\}}\|w_0\|\leq \frac{8 \delta_0 t}{\min\{\alpha(x),\beta(x)\}}.
\end{equation}
By \eqref{Eq2.17} and \eqref{Eq2.16}, noting that $y_0\in B_d(x,2t)$ (by \eqref{Usubl}),
we have
\begin{align*}
    \|F(x)(y_0-x)\|\leq\,& \|V(y_0)-V(x)-F(x)(y_0-x)\|+\|V(y_0)-V(x)\|\nonumber\\
    \leq\,&\delta\|y_0-x\|+\frac{8\delta_0 t}{\min\{\alpha(x),\beta(x)\}}\leq 2\delta t+\frac{8\delta_0 t}{\min\{\alpha(x),\beta(x)\}}.
\end{align*}
Now note that
\begin{align*}
    \|F(x)(y_0-x)\|=\|F(x)(s_0V(x)+w_0)\|\geq\,& s_0\|F(x)V(x)\|-\|F(x)w_0\|\nonumber\\
    \geq\,& \frac{t}{2}\|F(x)V(x)\|-\delta_0 t \|F(x)\|_2,
\end{align*}
where we recall that $s_0\geq t\slash 2$ and $\|w_0\|\leq \delta_0 t$. Combining the above two displays, we obtain
\begin{equation*}
    \|F(x)V(x)\|\leq 2\delta_0\|F(x)\|_2+4\delta+\frac{16\delta_0}{\min\{\alpha(x),\beta(x)\}}.
\end{equation*}
Therefore, taking first $\delta\rightarrow 0^+$ and then $\delta_0\rightarrow 0^{+}$, we conclude that $F(x)V(x)=0$.
\end{proof}

\begin{proof}[Proof of Theorem~\ref{L2.4}(a)--(b)] For any $T\in\mathcal{S}$, we define
\begin{equation}\label{defB_T}
    \mathcal{B}(T):=\big\{x\in\mathrm{int}(T): \eqref{Eq2.12}\text{ holds}\big\}.
\end{equation}
Let $B$ be the set of $x\in\mathbb{R}^d$ such that $\approxgrad V(x)$ exists. Note that $\mathcal{L}^d(\mathcal{T}_1^{*}\backslash B)=0$ (see Definition~\ref{DefVF}). By Lemma~\ref{L2.1}, for $\mathcal{L}^d$-a.e.\ $z_0\in\mathcal{T}_1^*$, $T=T(z_0)$ satisfies 
\begin{equation}\label{Eq2.14}
   \mathcal{H}^1(\mathrm{int}(T)\cap (\mathbb{R}^d\backslash B))=\mathcal{H}^1(T\cap (\mathcal{T}_1^{*}\backslash B))=0. 
\end{equation}
In the following, we consider any $T\in\mathcal{S}$ such that \eqref{Eq2.14} and the conclusion of Theorem~\ref{L2.2} hold. 

\paragraph{Step 1.} In this step, we will show that $\mathrm{int}(T)\cap B\subseteq \mathcal{B}(T)$, which by \eqref{Eq2.14} implies that 
\begin{equation}\label{Eq2.15}
     \mathcal{H}^1(\mathrm{int}(T)\backslash \mathcal{B}(T))=0. 
\end{equation}
In the following, we consider any
\begin{equation}\label{L2.4.e1}
    x\in \mathrm{int}(T)\cap B,  \quad \delta\in (0,\min\{\alpha(x),\beta(x),1\}\slash 10).
\end{equation}

Since $\approxgrad V(x)$ exists, we have $\approxgrad V(x)=F(x)$, and for any $\delta'>0$, there exists $t_0\in (0,\delta\slash 10)$ such that for all $t\in (0,t_0)$, 
\begin{equation*}
    \mathcal{L}^d\big(\big\{y\in \mathbb{R}^d:\|y-x\|\leq 2t, \|V(y)-V(x)-F(x)(y-x)\|\geq \delta\|y-x\|\slash 10 \big\}\big)\leq \delta' t^d.
\end{equation*}
Below we fix any $t\in (0,t_0)$. By the above display and the fact that $\sqrt{s^2+\|w\|^2}\leq 2t$ for any $s\in[-t,t]$ and $w\in O_T$ with $\|w\|\leq t$, there exists $s\in [-t,t]$ such that 
\begin{align}\label{condis}
&\mathcal{H}^{d-1}\big(\big\{ w \in O_T: \|w\|\leq t,\|V(x+sV_T+w)-V(x)-F(x)(sV_T+w)\| \geq \delta t\slash 5 \big\}\big)\nonumber\\
    \leq\,&\mathcal{H}^{d-1}\big(\big\{ w \in O_T: \|w\|\leq t, \|V(x+sV_T+w)-V(x)-F(x)(sV_T+w)\|\geq \delta\sqrt{s^2+\|w\|^2}\slash 10 \big\}\big)\nonumber\\
    \leq\,& \frac{\mathcal{L}^d\big(\big\{y\in \mathbb{R}^d:\|y-x\|\leq 2t, \|V(y)-V(x)-F(x)(y-x)\|\geq \delta\|y-x\|\slash 10 \big\}\big)}{2t}
\leq \delta' t^{d-1}. 
\end{align}
Fix such an $s\in [-t,t]$, and set $x':=x+s V_T$. As the conclusion of Theorem~\ref{L2.2} holds for $T$ and $x\in\mathrm{int}(T)$, the functions $\alpha$ and $\beta$ are approximately continuous at $x$. Moreover, since $x\in B$, we have that $\approxgrad V(x)$ exists. Hence by Lemma~\ref{L2.3}, we have $F(x)V_T=F(x)V(x)=0$. Using this and \eqref{condis}, we obtain
\begin{equation}\label{L2.4.eq10}
     \mathcal{H}^{d-1}\big(\big\{ w \in O_T: \|w\|\leq t,\|V(x'+w)-V(x)-F(x) w\|\geq\delta t \slash 5 \big\}\big)\leq\delta' t^{d-1}. 
\end{equation}
Below, we denote
\begin{equation}\label{L2.4.e6}
    \mathscr{E}:=\big\{ w \in O_T: \|w\|\leq t,\|V(x'+w)-V(x)-F(x) w\|<\delta t \slash 5 \big\}.
\end{equation}

Note that for any $w\in \mathcal{D}'(x;t,\delta)$ (recall Definition~\ref{Defd}), it follows from \eqref{L2.4.e1} that
\begin{equation}\label{alphabetan}
    \min\{\alpha(x+w),\beta(x+w)\}\geq \min\{\alpha(x),\beta(x)\}-\delta\geq \min\{\alpha(x),\beta(x)\}\slash 2>0.
\end{equation}
In particular, $x+w\in\mathcal{T}_1^*$. Hence by Lemma~\ref{L2.0}, we have  
\begin{equation}\label{L2.4.e2}
    \|V(x+w)-V(x)\|\leq\frac{8}{\min\{\alpha(x),\beta(x)\}}\|w\|\leq \frac{8t}{\min\{\alpha(x),\beta(x)\}} \leq \frac{1}{5},
\end{equation}
where the last inequality uses the fact that $t<t_0<\delta\slash 10<\min\{\alpha(x),\beta(x)\}\slash 100$. Consequently, 
\begin{align}\label{L2.4.e4}
   &|\langle V(x+w),V_T\rangle-1|=|\langle V(x+w)-V(x),V(x)\rangle|\nonumber\\
   \leq\,&\|V(x+w)-V(x)\|\leq \frac{8}{\min\{\alpha(x),\beta(x)\}}\|w\| \leq \frac{8t}{\min\{\alpha(x),\beta(x)\}},
\end{align}
\begin{equation}\label{L2.4.e3}
    \langle V(x+w),V_T\rangle  \geq 1-\frac{8t}{\min\{\alpha(x),\beta(x)\}}\geq \frac{4}{5}>0.
\end{equation}
Now for any $w\in \mathcal{D}'(x;t,\delta)$, we define 
\begin{equation}\label{varphi.def}
    \varphi(w):=x+w-x'-\frac{\langle x-x',V_T \rangle}{\langle V(x+w),V_T\rangle}V(x+w) \in O_T.
\end{equation}
Note that for any $w\in \mathcal{D}'(x;t,\delta)$, we have $x'+\varphi(w)\in L(x+w)$. Moreover, by \eqref{L2.4.e3}, we have 
\begin{equation*}
    \|(x'+\varphi(w))-(x+w)\|=\frac{|\langle x-x',V_T\rangle|}{|\langle V(x+w),V_T\rangle|}\leq 2\|x'-x\|\leq 2t\leq \frac{\delta}{5}. 
\end{equation*}
Hence by \eqref{alphabetan} and \eqref{L2.4.e1}, we have
\begin{equation}\label{L2.4.e10}
    \min\{\alpha(x+w),\beta(x+w)\}\geq \min\{\alpha(x),\beta(x)\}-\delta\geq 9\delta>\|(x'+\varphi(w))-(x+w)\|,
\end{equation}
and so $x'+\varphi(w)\in\mathrm{int}(T(x+w))$. Moreover, by \eqref{L2.4.e2}--\eqref{varphi.def}, we have
\begin{align}\label{L2.4.e5}
    \|\varphi(w)-w\|=\, & \Big\|\langle x-x', V(x)\rangle V(x)-\frac{\langle x-x', V(x)\rangle}{\langle V(x+w),V(x)\rangle} V(x+w)\Big\|\nonumber\\
   \leq\,&\|x-x'\|\Big\|V(x)-\frac{V(x+w)}{\langle V(x+w),V(x)\rangle}\Big\|\nonumber\\
   =\, &\|x-x'\|\frac{\|(\langle V(x+w),V(x)\rangle-1) V(x)-(V(x+w)-V(x))\|}{|\langle V(x+w),V(x)\rangle|}\nonumber\\
   \leq\, & 2\|x-x'\|\big(|\langle V(x+w),V(x)\rangle-1|+\|V(x+w)-V(x)\|\big)\nonumber\\
   \leq\, &  \frac{32\|x-x'\|}{\min\{\alpha(x),\beta(x)\}}\|w\|\leq \frac{32t}{\min\{\alpha(x),\beta(x)\}}\|w\|.
\end{align}

Take 
\begin{equation}\label{L2.4.e7}
    t':=\bigg(1+\frac{32  t}{\min\{\alpha(x),\beta(x)\}}\bigg)^{-1}t.
\end{equation}
By \eqref{L2.4.e5} and the definition of $\mathcal{D}'(x;t',\delta)$, for any $w\in \mathcal{D}'(x;t',\delta)\subseteq \mathcal{D}'(x;t,\delta)$, we have  
\begin{equation*}
    \|\varphi(w)\|\leq \bigg(1+\frac{32t}{\min\{\alpha(x),\beta(x)\}}\bigg)\|w\|\leq \bigg(1+\frac{32t}{\min\{\alpha(x),\beta(x)\}}\bigg)t'= t,
\end{equation*}
\begin{align*}
   & \max\{|\alpha(x'+\varphi(w))-\alpha(x')|, |\beta(x'+\varphi(w))-\beta(x')|\}\nonumber\\
   \leq\, & \max\{|\alpha(x+w)-\alpha(x)|,|\beta(x+w)-\beta(x)|\}+\big|\|(x'+\varphi(w))-(x+w)\|-\|x-x'\|\big|\nonumber\\
   \leq\,& \delta+\|\varphi(w)-w\|\leq \delta+\frac{32t}{\min\{\alpha(x),\beta(x)\}}\|w\|\leq \delta+\frac{32t^2}{\min\{\alpha(x),\beta(x)\}}\leq 2\delta,
\end{align*}
Hence we have
\begin{equation}\label{L2.4.e9}
    \varphi(w)\in\mathcal{D}'(x';t,2\delta),\quad\text{for every }w\in \mathcal{D}'(x;t',\delta).
\end{equation}

For any $w\in \mathcal{D}'(x;t',\delta)$, if $\varphi(w)\in\mathscr{E}$ (see \eqref{L2.4.e6}), then using the fact that $V(x+w)=V(x'+\varphi(w))$ (as $x'+\varphi(w)\in\mathrm{int}(T(x+w))$) and \eqref{L2.4.e5}, we obtain
\begin{align*}
   &\|V(x+w)-V(x)-F(x)w\|= \|V(x'+\varphi(w))-V(x)-F(x)w\| \nonumber\\
   \leq\,& \|V(x'+\varphi(w))-V(x)-F(x)\varphi(w)\|+\|F(x)(\varphi(w)-w)\|\nonumber\\
   \leq\,& \frac{\delta t}{5}+\frac{32t}{\min\{\alpha(x),\beta(x)\}}\|F(x)\|_2\|w\|\leq \frac{\delta t}{5}+\frac{32t^2\|F(x)\|_2}{\min\{\alpha(x),\beta(x)\}}.
\end{align*}
Hence, for $t$ sufficiently small (depending on $x$ and $\delta$), which we assume for the rest of the proof, for any $w\in \mathcal{D}'(x;t',\delta)$ such that $\varphi(w)\in \mathscr{E}$, we have $\|V(x+w)-V(x)-F(x)w\|\leq \delta t$. Therefore, recalling Definition~\ref{Defd}, we obtain
\begin{equation*}
    \{w\in \mathcal{D}'(x;t',\delta):\varphi(w)\in\mathscr{E}\}\subseteq \mathcal{W}'(x;t,\delta),
\end{equation*}
and so 
\begin{align}\label{L2.4.e8}
    \mathcal{W}(x;t,\delta)=\,&\{w\in O_T:\|w\|\leq t\}\backslash\mathcal{W}'(x;t,\delta)\nonumber\\
    \subseteq\,& \{w\in O_T:t'< \|w\|\leq t\}\cup \big(\{w\in O_T:\|w\|\leq t'\}\backslash \mathcal{W}'(x;t,\delta)\big) \nonumber\\
    \subseteq\,& \{w\in O_T:t'< \|w\|\leq t\}\cup \mathcal{D}(x;t',\delta)\cup  \{w\in\mathcal{D}'(x;t',\delta):\varphi(w)\notin\mathscr{E}\}.
\end{align}

We now proceed to bound $\mathcal{H}^{d-1}\big(\big\{w\in\mathcal{D}'(x;t',\delta):\varphi(w)\notin \mathscr{E}\big\}\big)$. For any $w'\in\mathcal{D}'(x';t,2\delta)$, we define
\begin{equation}\label{psi_def}
    \psi(w'):=x'+w'-x-\frac{\langle x'-x,V_T \rangle}{\langle V(x'+w'),V_T\rangle}V(x'+w')\in O_T  .
\end{equation}
For any $w'\in \mathcal{D}'(x';t,2\delta)$, using the definition of $\mathcal{D}'(x';t,2\delta)$ and the fact that $\min\{\alpha(x'),\beta(x')\}\geq\min\{\alpha(x),\beta(x)\}-t$, we obtain
\begin{equation}\label{x'w'pro}
    \min\{\alpha(x'+w'),\beta(x'+w')\}\geq \min\{\alpha(x'),\beta(x')\}-2\delta\geq \min\{\alpha(x),\beta(x)\}-t-2\delta\geq \frac{1}{2}\min\{\alpha(x),\beta(x)\},
\end{equation}
where the last inequality uses \eqref{L2.4.e1} and the fact that $t<t_0<\delta\slash 10$. Using \eqref{x'w'pro} and arguing similarly as in \eqref{L2.4.e2}--\eqref{L2.4.e10}, we obtain that for any $w'\in\mathcal{D}'(x';t,2\delta)$,
\begin{equation}\label{eqnVx'}
  \langle V(x'+w'), V_T\rangle \geq \frac{4}{5},\qquad  x+\psi(w')\in\mathrm{int}(T(x'+w')).
\end{equation}
Moreover, for any $w\in \mathcal{D}'(x;t',\delta)$, noting \eqref{L2.4.e9}, we have $\psi(\varphi(w))=w$. Hence by \eqref{L2.4.e9}, for any $w\in \mathcal{D}'(x;t',\delta)$ such that $\varphi(w)\notin\mathscr{E}$, we have $w=\psi(\varphi(w))\in \psi\big(\mathcal{D}'(x';t,2\delta)\backslash \mathscr{E}\big)$, and so 
\begin{equation}\label{L2.4.eq11}
    \big\{w\in\mathcal{D}'(x;t',\delta):\varphi(w)\notin\mathscr{E}\big\}\subseteq \psi\big(\mathcal{D}'(x';t,2\delta)\backslash \mathscr{E}\big).
\end{equation}
For any $w_1',w_2'\in \mathcal{D}'(x';t,2\delta)$,  by Lemma~\ref{L2.0} and \eqref{x'w'pro}, we have
\begin{equation}\label{L2.4.e14}
    \|V(x'+w_1')-V(x'+w_2')\|\leq\frac{8}{\min\{\alpha(x),\beta(x)\}}\|w_1'-w_2'\|. 
\end{equation}
Consequently, using \eqref{eqnVx'}, we obtain
\begin{align}
   &\bigg\|\frac{V(x'+w_1')}{\langle V(x'+w_1'),V_T\rangle}-\frac{V(x'+w_2')}{\langle V(x'+w_2'),V_T\rangle}\bigg\|  
 \nonumber\\
 \leq\,& \bigg\|\frac{V(x'+w_1')}{\langle V(x'+w_1'),V_T\rangle}-\frac{V(x'+w_2')}{\langle V(x'+w_1'),V_T\rangle}\bigg\|+\bigg\|\frac{V(x'+w_2')}{\langle V(x'+w_1'),V_T\rangle}-\frac{V(x'+w_2')}{\langle V(x'+w_2'),V_T\rangle}\bigg\|\nonumber\\
   \leq\,&  \frac{\|V(x'+w_1')-V(x'+w_2')\|}{\langle V(x'+w_1'), V_T\rangle}+\frac{|\langle V(x'+w_1')-V(x'+w_2'), V_T\rangle|}{\langle V(x'+w_1'), V_T\rangle  \langle V(x'+w_2'), V_T\rangle}\nonumber\\
       \leq\,&8\|V(x'+w_1')-V(x'+w_2')\|\leq \frac{64}{\min\{\alpha(x), 
 \beta(x)\}}\|w_1'-w_2'\|,
\end{align}
and so
\begin{align}\label{L2.4.e15}
   \|\psi(w_1')-\psi(w_2')\|\leq\,&\|w_1'-w_2'\|+|\langle x'-x,V_T \rangle|\bigg\|\frac{V(x'+w_1')}{\langle V(x'+w_1'),V_T\rangle}-\frac{V(x'+w_2')}{\langle V(x'+w_2'),V_T\rangle}\bigg\|\nonumber\\
   \leq\,&\bigg(1+\frac{64 t}{\min\{\alpha(x),\beta(x)\}}\bigg)\|w_1'-w_2'\|.
\end{align}
By Kirszbraun's theorem, there exists a Lipschitz mapping $\tilde{\psi}:O_T\rightarrow O_T$ such that $\tilde{\psi}(w')=\psi(w')$ for all $w'\in\mathcal{D}'(x';t,2\delta)$ and $\|\tilde{\psi}\|_{\mathrm{Lip}}\leq 1+\frac{64t}{\min\{\alpha(x),\beta(x)\}}$. 
Hence by \eqref{L2.4.eq11} and \cite[Theorem 2.8]{MR3409135}, we have 
\begin{align}\label{L2.4.eq14}
   &\mathcal{H}^{d-1}\big(\big\{w\in\mathcal{D}'(x;t',\delta):\varphi(w)\notin\mathscr{E}\big\}\big)\leq \mathcal{H}^{d-1}\big(\psi\big(\mathcal{D}'(x';t,2\delta)\backslash \mathscr{E}\big)\big)=\mathcal{H}^{d-1}\big(\tilde{\psi}\big(\mathcal{D}'(x';t,2\delta)\backslash \mathscr{E}\big)\big)\nonumber\\
   \leq\,& \bigg(1+\frac{64t}{\min\{\alpha(x),\beta(x)\}}\bigg)^{d-1}\mathcal{H}^{d-1}\big(\mathcal{D}'(x';t,2\delta)\backslash \mathscr{E}\big)\leq \bigg(1+\frac{64t}{\min\{\alpha(x),\beta(x)\}}\bigg)^{d-1}\delta' t^{d-1},
\end{align}
where the last inequality uses $\mathcal{H}^{d-1}\big(\mathcal{D}'(x';t,2\delta)\backslash \mathscr{E}\big)\leq \mathcal{H}^{d-1}(\{w\in O_T:\|w\|\leq t\}\backslash\mathscr{E})\leq \delta' t^{d-1}$ (see \eqref{L2.4.eq10} and \eqref{L2.4.e6}).

By \eqref{L2.4.e8} and \eqref{L2.4.eq14}, we have (recall \eqref{L2.4.e7})
\begin{align*}
     \mathcal{H}^{d-1}(\mathcal{W}(x;t,\delta)) 
    \leq\,&\mathcal{H}^{d-1}\big(\big\{w\in O_T:t'< \|w\|\leq t\big\}\big)+\mathcal{H}^{d-1}(\mathcal{D}(x;t',\delta))\nonumber\\
    & +\mathcal{H}^{d-1}\big(\big\{w\in\mathcal{D}'(x;t',\delta):\varphi(w)\notin\mathscr{E}\big\}\big)\nonumber\\
    \leq\,& \frac{\pi^{(d-1)\slash 2}}{\Gamma((d+1) \slash 2)}\bigg(1-\bigg(1+\frac{32t}{\min\{\alpha(x),\beta(x)\}}\bigg)^{-(d-1)}\bigg)t^{d-1}\nonumber\\
    &+\mathcal{H}^{d-1}(\mathcal{D}(x;t',\delta))+\bigg(1+\frac{64t}{\min\{\alpha(x),\beta(x)\}}\bigg)^{d-1}\delta' t^{d-1}.
\end{align*}
Hence by Theorem~\ref{L2.2}, we have 
\begin{equation*}
    \limsup_{t\rightarrow 0^{+}}\frac{\mathcal{H}^{d-1}(\mathcal{W}(x;t,\delta))}{t^{d-1}}\leq \delta'. 
\end{equation*}
As $\delta'>0$ can be arbitrarily small, we have
\begin{equation*}
     \lim_{t\rightarrow 0^{+}}\frac{\mathcal{H}^{d-1}(\mathcal{W}(x;t,\delta))}{t^{d-1}}=0.
\end{equation*}
Hence $x\in \mathcal{B}(T)$. We conclude that $\mathrm{int}(T)\cap B\subseteq \mathcal{B}(T)$.

\paragraph{Step 2.} In this step, we show that for any $x\in\mathrm{int}(T)$ such that both \eqref{Eq2.11n} and \eqref{Eq2.12} hold, 
\begin{equation}\label{L2.4.e7n}
    \|F(x)\|_2  \leq  \frac{4}{\min\{\alpha(x),\beta(x)\}}.
\end{equation} 

Fix any $x\in\mathrm{int}(T)$ such that both \eqref{Eq2.11n} and \eqref{Eq2.12} hold. By \eqref{Eq2.1n} and \eqref{Eq2.12}, for any $\delta>0$,
\begin{equation*}
    \lim_{t\rightarrow 0^{+}}\frac{\mathcal{H}^{d-1}\big(\big\{w\in O_T:\|w\|\leq t, \max\{|\alpha(x+w)-\alpha(x)|, |\beta(x+w)-\beta(x)|\}>\delta\big\}\big)}{t^{d-1}}=0,
\end{equation*}
\begin{equation*}
    \lim_{t\rightarrow 0^{+}}\frac{\mathcal{H}^{d-1}\big(\big\{w\in O_T:\|w\|\leq t, \|V(x+w)-V(x)-F(x)w\|>\delta t\big\}\big)}{t^{d-1}}=0. 
\end{equation*}
Below, we fix any $w_0\in O_T$ with $\|w_0\|=1$ and any $\delta\in (0,\min\{\alpha(x),\beta(x),1\}\slash 10)$. By the above two displays, there exist $(t_n)_{n=1}^{\infty}\subseteq (0,1)$ and $(w_n)_{n=1}^{\infty}\subseteq O_T$ such that as $n\rightarrow \infty$,
\begin{equation}\label{L2.4.eq21}
    t_n\rightarrow 0,\quad (1-\delta)t_n   \leq \|w_n\|\leq t_n, \quad \frac{w_n}{\|w_n\|}\rightarrow w_0,
\end{equation}
\begin{equation}\label{L2.4.eq24}
  \max\{|\alpha(x+w_n)-\alpha(x)|,|\beta(x+w_n)-\beta(x)|\}\leq \delta, 
\end{equation}
\begin{equation}\label{L2.4.eq22}
    \|V(x+w_n)-V(x)-F(x)w_n\|\leq \delta t_n.
\end{equation}
By \eqref{L2.4.eq24}, we have
\begin{equation*}
    \min\{\alpha(x+w_n),\beta(x+w_n)\}\geq \min\{\alpha(x),\beta(x)\}-\delta>0.
\end{equation*}
Hence by Lemma~\ref{L2.0},  
\begin{equation}\label{L2.4.eq23}
    \|V(x+w_n)-V(x)\|\leq \frac{4}{\min\{\alpha(x),\beta(x)\}- \delta }\|w_n\|.
\end{equation}
Combining \eqref{L2.4.eq21}, \eqref{L2.4.eq22}, and \eqref{L2.4.eq23}, we obtain
\begin{equation*}
    \|F(x) w_n\|\leq \frac{4}{\min\{\alpha(x),\beta(x)\}-\delta}\|w_n\|+\delta t_n\leq \bigg(\frac{4}{\min\{\alpha(x),\beta(x)\}-\delta}+\delta\bigg) t_n.
\end{equation*}
Hence
\begin{align}\label{L2.4.eq25}
   \|F(x) w_0\|\leq\,&  \Big\|F(x)\frac{w_n}{\|w_n\|}\Big\|+\Big\|F(x)\Big(\frac{w_n}{\|w_n\|}-w_0\Big)\Big\|\nonumber\\
   \leq\,& \bigg(\frac{4}{\min\{\alpha(x),\beta(x)\}-\delta}+\delta\bigg)\frac{t_n}{\|w_n\|}+\|F(x)\|_2\Big\|\frac{w_n}{\|w_n\|}-w_0\Big\|\nonumber\\
   \leq\,& \frac{4}{\min\{\alpha(x),\beta(x)\}-\delta}\cdot\frac{1}{1-\delta}+\frac{\delta}{1-\delta}+\|F(x)\|_2\Big\|\frac{w_n}{\|w_n\|}-w_0\Big\|.
\end{align}
First taking $n\rightarrow\infty$, and then taking $\delta\rightarrow 0^{+}$, we obtain that for any $w_0\in O_T$ such that $\|w_0\|=1$,
\begin{equation*}
    \|F(x)w_0\|\leq \frac{4}{\min\{\alpha(x),\beta(x)\}}.
\end{equation*}
By Theorem~\ref{L2.2} and Lemma~\ref{L2.3}, we have $F(x)V(x)=0$ (note that \eqref{Eq2.11n} holds). Hence for any $w\in \mathbb{R}^d$ such that $\|w\|=1$, we have $\|F(x)w\|\leq \frac{4}{\min\{\alpha(x),\beta(x)\}}$, and so
\begin{equation*}
    \|F(x)\|_2\leq  \frac{4}{\min\{\alpha(x),\beta(x)\}}.
\end{equation*}

\paragraph{Step 3.} In this step, we show that for any $x,x'\in\mathrm{int}(T)$ such that \eqref{Eq2.11n} and \eqref{Eq2.12} hold, we have 
\begin{equation}\label{L2.4.eq18}
    \|F(x)-F(x')\|_2 \leq \frac{128\|x-x'\|}{\min\{\alpha(x),\beta(x)\}\min\{\alpha(x'),\beta(x')\}},
\end{equation}
\begin{equation}\label{L2.4.e15n}
    F(x') =  \big(\mathbf{I}_d-\langle x'-x, V_T\rangle F(x')\big)F(x).
\end{equation} 

In the following, we fix any $x,x'\in\mathrm{int}(T)$ such that \eqref{Eq2.11n} and \eqref{Eq2.12} hold. We also fix
\begin{equation}\label{deltacon}
    \delta\in \bigg(0,\frac{\min\{\alpha(x),\beta(x),\alpha(x'),\beta(x'),1\}}{100\max\{\|x-x'\|,1\}}\bigg),
\end{equation}
and assume that $t\in (0,\delta^3\slash 10)$. 

For any $w\in\mathcal{D}'(x;t,\delta)$, define $\varphi(w)$ as in \eqref{varphi.def}. Arguing similarly as in \eqref{alphabetan}--\eqref{L2.4.e3} and \eqref{L2.4.e5}, we obtain
\begin{equation}\label{L2.4.e8n}
    \|V(x+w)-V(x)\|  \leq  \frac{8}{\min\{\alpha(x),\beta(x)\}}\|w\|,\qquad \langle V(x+w),V_T\rangle\geq \frac{4}{5},  
\end{equation}
\begin{equation}\label{L2.4.e5n}
    \big|\|(x'+\varphi(w))-(x+w)\|-\|x-x'\|\big|\leq\|\varphi(w)-w\|\leq\frac{32\|x-x'\|\|w\|}{\min\{\alpha(x),\beta(x)\}}\leq\frac{32\|x-x'\|t}{\min\{\alpha(x),\beta(x)\}}\leq \delta,
\end{equation}
where the last inequality in \eqref{L2.4.e5n} follows from \eqref{deltacon}:
\begin{equation*}
    \frac{32\|x-x'\|t}{\min\{\alpha(x),\beta(x)\}}\leq \frac{32\|x-x'\|}{\min\{\alpha(x),\beta(x)\}}\cdot\frac{\delta^3}{10}\leq \frac{32\|x-x'\|}{\min\{\alpha(x),\beta(x)\}}\cdot\frac{\min\{\alpha(x),\beta(x)\}}{100\max\{\|x-x'\|,1\}}\cdot\frac{\delta}{10}\leq \delta.
\end{equation*}
By \eqref{varphi.def}, we have $\langle (x'+\varphi(w))-(x+w),V(x+w)\rangle=\frac{\langle x'-x,V_T \rangle}{\langle V(x+w),V_T\rangle}$; since $\langle V(x+w),V_T\rangle>0$ (by \eqref{L2.4.e8n}), $\langle (x'+\varphi(w))-(x+w),V(x+w)\rangle$ has the same sign as $\langle x'-x, V_T\rangle$. Hence by the definition of $\mathcal{D}'(x;t,\delta)$ and \eqref{L2.4.e5n}, we have
{\small
\begin{align}\label{intargu}
   &\min\big\{\alpha(x+w)-\langle (x'+\varphi(w))-(x+w),V(x+w)\rangle,\beta(x+w)+\langle (x'+\varphi(w))-(x+w),V(x+w)\rangle\big\} \nonumber\\
   \geq\,& \min\big\{\alpha(x)-\delta-\langle x'-x,V_T\rangle, \beta(x)-\delta+\langle x'-x,V_T\rangle\big\}\nonumber\\
   &-\big|\langle (x'+\varphi(w))-(x+w),V(x+w)\rangle-\langle x'-x,V_T\rangle\big|\nonumber\\
   =\, & \min\{\alpha(x'),\beta(x')\}-\delta-\big|\|(x'+\varphi(w))-(x+w)\|-\|x-x'\|\big|\geq \min\{\alpha(x'),\beta(x')\}-2\delta>0,
\end{align}
}and so
\begin{equation}\label{L2.4.e38}
   x'+\varphi(w)\in \mathrm{int}(T(x+w)).
\end{equation}
Moreover, by the definition of $\mathcal{D}'(x;t,\delta)$ and \eqref{L2.4.e5n}, we have 
\begin{align}\label{L2.4.e28}
     &\max\{|\alpha(x'+\varphi(w))-\alpha(x')|, |\beta(x'+\varphi(w))-\beta(x')|\}\nonumber\\
   \leq\,& \max\{|\alpha(x+w)-\alpha(x)|,|\beta(x+w)-\beta(x)|\}+\big|\|(x'+\varphi(w))-(x+w)\|-\|x-x'\|\big|\leq 2\delta,
\end{align}
\begin{equation}\label{L2.4.e29}
    \|\varphi(w)\|\leq \|w\|+\|\varphi(w)-w\|\leq \bigg(1+\frac{32\|x-x'\|}{\min\{\alpha(x),\beta(x)\}}\bigg)t.
\end{equation}
Below, we denote
\begin{equation}\label{defini}
    M:=1+\frac{32\|x-x'\|}{\min\{\alpha(x),\beta(x)\}}.
\end{equation}
By \eqref{L2.4.e28} and \eqref{L2.4.e29},
\begin{equation}\label{Ln}
    \varphi(w)  \in \mathcal{D}'(x';Mt,2\delta),\quad\text{for any }  w\in\mathcal{D}'(x;t,\delta).
\end{equation}

Moreover, for any $w\in\mathcal{D}'(x;t,\delta)$, as $\|V(x+w)\|=\|V(x)\|=1$, we have
\begin{equation*}
    \|V(x+w)-V(x)\|^2=2-2\langle V(x+w),V(x)\rangle.
\end{equation*}
Hence by \eqref{L2.4.e8n}, we have
\begin{align}\label{eq2.4.new}
    &  \langle V(x+w),V(x)\rangle\leq \|V(x+w)\|\|V(x)\|=1,\nonumber\\
    &  \langle V(x+w),V(x)\rangle=1-\frac{1}{2}\|V(x+w)-V(x)\|^2\geq 1-\frac{32\|w\|^2}{(\min\{\alpha(x),\beta(x)\})^2} \geq \frac{1}{2},
\end{align}
where the last inequality uses $\|w\|\leq t\leq \delta\slash 10\leq \min\{\alpha(x),\beta(x)\}\slash 100$. Therefore, by \eqref{varphi.def},
\begin{align}\label{L2.4.e9n}
   & \|\varphi(w)-w+\langle x-x',V(x)\rangle F(x)w\|\nonumber\\
   =\,&\bigg\|x-x'-\frac{\langle x-x',V(x) \rangle}{\langle V(x+w),V(x)\rangle}V(x+w)+\langle x-x',V(x)\rangle F(x)w\bigg\|\nonumber\\
   =\,& \big|\langle x-x',V(x)\rangle\big|\bigg\| V(x)-\frac{V(x+w)}{\langle V(x+w),V(x)\rangle}+F(x)w\bigg\|\nonumber\\
   \leq\,& \|x-x'\|\bigg(\|V(x)-V(x+w)+F(x)w\|+\bigg|1-\frac{1}{\langle V(x+w),V(x)\rangle}\bigg|\bigg)\nonumber\\
   \leq\,& \|x-x'\|\big(\|V(x)-V(x+w)+F(x)w\|+2|\langle V(x+w),V(x)\rangle-1|\big)\nonumber\\
   \leq\,& \|x-x'\|\bigg(\|V(x)-V(x+w)+F(x)w\|+\frac{64\|w\|^2}{(\min\{\alpha(x),\beta(x)\})^2}\bigg).
\end{align}

Now for any $w'\in\mathcal{D}'(x';Mt,2\delta)$, define $\psi(w')$ as in \eqref{psi_def}. As
\begin{equation*}
    \min\{\alpha(x'+w'),\beta(x'+w')\}\geq \min\{\alpha(x'),\beta(x')\}-2\delta\geq \min\{\alpha(x'),\beta(x')\}\slash 2,
\end{equation*}
by Lemma~\ref{L2.0}, we have
\begin{align*}
    \|V(x'+w')-V(x')\|\leq\,& \frac{8}{\min\{\alpha(x'),\beta(x')\}}\|w'\|\leq \frac{8Mt}{\min\{\alpha(x'),\beta(x')\}}\nonumber\\
    \leq\,& \frac{8}{\min\{\alpha(x'),\beta(x')\}}\cdot\bigg(1+\frac{32\|x-x'\|}{\min\{\alpha(x),\beta(x)\}}\bigg)\cdot\frac{\delta^3}{10} \leq \frac{1}{5},
\end{align*}
where the last two inequalities use \eqref{defini} and \eqref{deltacon}, respectively. Consequently,
\begin{equation*}
    |\langle V(x'+w'),V(x')\rangle-1|=|\langle V(x'+w')-V(x'),V(x')\rangle|\leq \|V(x'+w')-V(x')\|\leq \frac{1}{5}.
\end{equation*}
By the above two displays, noting \eqref{deltacon}, we have
\begin{align}
    & \big|\|(x+\psi(w'))-(x'+w')\|-\|x-x'\|\big|\leq \|\psi(w')-w'\|\nonumber\\
    \leq\,&\Big\|\langle x'-x, V(x')\rangle V(x')-\frac{\langle x'-x, V(x')\rangle}{\langle V(x'+w'),V(x')\rangle} V(x'+w')\Big\|\nonumber\\
   \leq\,&\|x-x'\|\Big\|V(x')-\frac{V(x'+w')}{\langle V(x'+w'),V(x')\rangle}\Big\|\nonumber\\
   =\,&\|x-x'\|\frac{\|(\langle V(x'+w'),V(x')\rangle-1) V(x')-(V(x'+w')-V(x'))\|}{|\langle V(x'+w'),V(x')\rangle|}\nonumber\\
   \leq\,& 2\|x-x'\|\big(|\langle V(x'+w'),V(x')\rangle-1|+\|V(x'+w')-V(x')\|\big)\nonumber\\
   \leq\,& 4\|x-x'\|\|V(x'+w')-V(x')\|\leq\frac{32\|x-x'\|}{\min\{\alpha(x'),\beta(x')\}}\|w'\|\leq \frac{32\|x-x'\|Mt}{\min\{\alpha(x'),\beta(x')\}}\nonumber\\
   \leq\,& \frac{32\|x-x'\|}{\min\{\alpha(x'),\beta(x')\}}\cdot \bigg(1+\frac{32\|x-x'\|}{\min\{\alpha(x),\beta(x)\}}\bigg)\cdot\frac{\delta^3}{10}\leq \delta.
\end{align}
Therefore, arguing similarly as in \eqref{intargu}, we obtain
\begin{equation*}
    x+\psi(w')\in\mathrm{int}(T(x'+w')),\quad\text{for any }w'\in\mathcal{D}'(x';Mt,2\delta).
\end{equation*}
Consequently, noting \eqref{Ln}, we have
\begin{equation}\label{psiphi}
    \psi(\varphi(w))=w,\quad\text{for any }w\in\mathcal{D}'(x;t,\delta).
\end{equation}

We now define
\begin{equation}\label{L2.4.eq32}
    \mathscr{G}(x,x';t,\delta):=\big\{  w\in\mathcal{D}'(x;t,\delta)\cap\mathcal{W}'(x;t,\delta): \varphi(w)\notin \mathcal{W}'(x';Mt,2\delta)\big\},
\end{equation}
where $M$ is as in \eqref{defini}. By \eqref{Ln} and \eqref{psiphi}, we have 
\begin{align}\label{L2.4.eq11nn}
 \mathscr{G}(x,x';t,\delta)
    =\,&\big\{  w\in\mathcal{D}'(x;t,\delta)\cap\mathcal{W}'(x;t,\delta): \varphi(w)\in\mathcal{D}'(x';Mt,2\delta)\backslash \mathcal{W}'(x';Mt,2\delta)\big\}\nonumber\\
    \subseteq\,&\psi\big(\mathcal{D}'(x';Mt,2\delta)\backslash \mathcal{W}'(x';Mt,2\delta)\big). 
\end{align}
Arguing similarly as in \eqref{L2.4.e14}--\eqref{L2.4.e15}, we obtain $\psi$ is $\Big(1+\frac{64\|x-x'\|}{\min\{\alpha(x'),\beta(x')\}}\Big)$-Lipschitz on $\mathcal{D}'(x';Mt,2\delta)$. By Kirszbraun's theorem, there exists a Lipschitz mapping $\tilde{\psi}:O_T\rightarrow O_T$ such that $\tilde{\psi}(w')=\psi(w')$ for all $w'\in\mathcal{D}'(x';Mt,2\delta)$ and $\|\tilde{\psi}\|_{\mathrm{Lip}}\leq 1+\frac{64\|x-x'\|}{\min\{\alpha(x'),\beta(x')\}}$. 
Hence by \eqref{L2.4.eq11nn} and \cite[Theorem 2.8]{MR3409135}, 
\begin{align}\label{L2.4.eq35}
  \mathcal{H}^{d-1}(\mathscr{G}(x,x';t,\delta))\leq\,&  \mathcal{H}^{d-1}\big(\tilde{\psi}\big(\mathcal{D}'(x';Mt,2\delta)\backslash \mathcal{W}'(x';Mt,2\delta)\big)\big)  \nonumber\\
  \leq\,& \bigg(1+\frac{64\|x-x'\|}{\min\{\alpha(x'),\beta(x')\}}\bigg)^{d-1} \mathcal{H}^{d-1}\big(\mathcal{D}'(x';Mt,2\delta)\backslash \mathcal{W}'(x';Mt,2\delta)\big)\nonumber\\
  \leq\,& \bigg(1+\frac{64\|x-x'\|}{\min\{\alpha(x'),\beta(x')\}}\bigg)^{d-1} \mathcal{H}^{d-1}\big(\mathcal{W}(x';Mt,2\delta)\big).
\end{align}

Define
\begin{equation}
    \mathscr{G}'(x,x';t,\delta):=\big\{  w\in\mathcal{D}'(x;t,\delta)\cap\mathcal{W}'(x;t,\delta): \varphi(w)\in \mathcal{D}'(x';Mt,2\delta)\cap\mathcal{W}'(x';Mt,2\delta)\big\}.
\end{equation}
By \eqref{Ln} and \eqref{L2.4.eq32}, we have 
\begin{align*}
    \mathscr{G}'(x,x';t,\delta) =\,&\big\{  w\in\mathcal{D}'(x;t,\delta)\cap\mathcal{W}'(x;t,\delta): \varphi(w)\in \mathcal{W}'(x';Mt,2\delta)\big\}\nonumber\\
    =\,& (\mathcal{D}'(x;t,\delta)\cap\mathcal{W}'(x;t,\delta))\backslash \mathscr{G}(x,x';t,\delta)\nonumber\\
    =\,& \{w\in O_T:\|w\|\leq t\}\backslash (\mathcal{D}(x;t,\delta)\cup\mathcal{W}(x;t,\delta)\cup\mathscr{G}(x,x';t,\delta)).
\end{align*}
Hence by \eqref{L2.4.eq35},
\begin{align*}
& \mathcal{H}^{d-1}(\{w\in O_T:\|w\|\leq t\})\nonumber\\
  \geq\,& \mathcal{H}^{d-1}(\mathscr{G}'(x,x';t,\delta))\nonumber\\
   \geq\,& \mathcal{H}^{d-1}(\{w\in O_T:\|w\|\leq t\})-\mathcal{H}^{d-1}(\mathcal{D}(x;t,\delta))-\mathcal{H}^{d-1}(\mathcal{W}(x;t,\delta))-\mathcal{H}^{d-1}(\mathscr{G}(x,x';t,\delta))\nonumber\\
   \geq\,& \mathcal{H}^{d-1}(\{w\in O_T:\|w\|\leq t\})-\mathcal{H}^{d-1}(\mathcal{D}(x;t,\delta))-\mathcal{H}^{d-1}(\mathcal{W}(x;t,\delta))\nonumber\\
   &-\bigg(1+\frac{64\|x-x'\|}{\min\{\alpha(x'),\beta(x')\}}\bigg)^{d-1} \mathcal{H}^{d-1}(\mathcal{W}(x';Mt,2\delta)).
\end{align*}
As the conclusion of Theorem~\ref{L2.2} holds for $T$ and \eqref{Eq2.11n}--\eqref{Eq2.12} hold for $x$ and $x'$, we have
\begin{equation*}
     \lim_{t\rightarrow 0^{+}}\frac{\mathcal{H}^{d-1}(\mathcal{D}(x;t,\delta))}{t^{d-1}}=\lim_{t\rightarrow 0^{+}}\frac{\mathcal{H}^{d-1}(\mathcal{W}(x;t,\delta))}{t^{d-1}}=\lim_{t\rightarrow 0^+}\frac{\mathcal{H}^{d-1}(\mathcal{W}(x';Mt,2\delta))}{t^{d-1}}=0.
\end{equation*}
Combining the above two displays, we get
\begin{equation}\label{L2.4.e37}
    \lim_{t\rightarrow 0^{+}}\frac{\mathcal{H}^{d-1}(\mathscr{G}'(x,x';t,\delta))}{t^{d-1}} = \lim_{t\rightarrow 0^{+}}\frac{\mathcal{H}^{d-1}(\{w\in O_T:\|w\|\leq t\})}{t^{d-1}}=\frac{\pi^{(d-1)\slash 2}}{\Gamma((d+1)\slash 2)}.
\end{equation}

By \eqref{L2.4.e37}, for any $w_0\in O_T$ with $\|w_0\|=1$, there exist $(t_n)_{n=1}^{\infty}\subseteq (0,1)$ and
$(w_n)_{n=1}^{\infty}\subseteq O_T$ such that $w_n\in \mathscr{G}'(x,x';t_n,\delta)$ for every $n\in\mathbb{N}^*$ and, as $n\rightarrow \infty$, 
\begin{equation}\label{L2.4.e1n}
    t_n\rightarrow 0,\quad t_n \slash 2\leq \|w_n\|\leq t_n,\quad \frac{w_n}{\|w_n\|}\rightarrow w_0.
\end{equation}
In what follows, we assume $n$ is sufficiently large so that $t_n<\delta^3\slash 10$. As $w_n\in\mathcal{D}'(x;t_n,\delta)\cap\mathcal{W}'(x;t_n,\delta)$ and $\varphi(w_n)\in \mathcal{W}'(x';Mt_n,2\delta)$, by \eqref{L2.4.e5n} and \eqref{L2.4.e9n}, we have 
\begin{equation}\label{L2.4.e6n}
    \|\varphi(w_n)-w_n\|\leq \frac{32\|x-x'\|\|w_n\|}{\min\{\alpha(x),\beta(x)\}},
\end{equation}
\begin{align}\label{L2.4.e10n}
   & \|\varphi(w_n)-w_n+\langle x-x',V(x)\rangle F(x)w_n\|\nonumber\\
   \leq\,& \|x-x'\|\bigg(\|V(x)-V(x+w_n)+F(x)w_n\|+\frac{64\|w_n\|^2}{(\min\{\alpha(x),\beta(x)\})^2}\bigg)\nonumber\\
   \leq\,& \|x-x'\|\bigg( \delta t_n +\frac{64\|w_n\|^2}{(\min\{\alpha(x),\beta(x)\})^2}\bigg),
\end{align}
\begin{equation}\label{L2.4.e2n}
    \|V(x+w_n)-V(x)-F(x)w_n\|\leq \delta t_n,
\end{equation}
\begin{equation}
    \|V(x'+\varphi(w_n))-V(x')-F(x')\varphi(w_n)\|\leq 2M\delta t_n.
\end{equation}
By \eqref{L2.4.e38}, we have
\begin{equation}\label{L2.4.e39}
    V(x+w_n)=V(x'+\varphi(w_n)).
\end{equation}
By \eqref{L2.4.e2n}--\eqref{L2.4.e39} and the fact that $V(x)=V(x')$, we have
\begin{equation}\label{L2.4.e3n}
    \|F(x)w_n-F(x')\varphi(w_n)\|\leq (2M+1)\delta t_n. 
\end{equation}
Using \eqref{L2.4.e3n}, \eqref{L2.4.e1n}, and \eqref{L2.4.e6n}, we obtain that
\begin{align}\label{L2.4.e4n}
   & \Big\|(F(x)-F(x'))\frac{w_n}{\|w_n\|}\Big\| \leq \frac{\|F(x)w_n-F(x')\varphi(w_n)\|+\|F(x')(\varphi(w_n)-w_n)\|}{\|w_n\|}\nonumber\\
   \leq\,& 2(2M+1)\delta+\|F(x')\|_2\cdot\frac{\|\varphi(w_n)-w_n\|}{\|w_n\|}\leq 2(2M+1)\delta+\frac{32\|x-x'\|\|F(x')\|_2}{\min\{\alpha(x),\beta(x)\}}\nonumber\\
   \leq\,& 2(2M+1)\delta+\frac{128\|x-x'\|}{\min\{\alpha(x),\beta(x)\}\min\{\alpha(x'),\beta(x')\}},
\end{align}
where the last inequality uses \eqref{L2.4.e7n}. Moreover, by \eqref{L2.4.e3n}, \eqref{L2.4.e1n}, and \eqref{L2.4.e10n}, we have
\begin{align}\label{L2.4.e12n}
    &\Big\|\Big(F(x)-F(x')\big(\mathbf{I}_d-\langle x-x',V(x)\rangle F(x)\big)\Big)\frac{w_n}{\|w_n\|}\Big\|\nonumber\\
    \leq\,&\frac{\|F(x)w_n-F(x')\varphi(w_n)\|+\big\|F(x')\big(\varphi(w_n)-\big(\mathbf{I}_d-\langle x-x',V(x)\rangle F(x)\big)w_n\big)\big\|}{\|w_n\|}\nonumber\\
    \leq\,& 2(2M+1)\delta+\frac{\|F(x')\|_2}{\|w_n\|}\cdot\big\|\varphi(w_n)-\big(\mathbf{I}_d-\langle x-x',V(x)\rangle F(x)\big)w_n\big\|\nonumber\\
    \leq\,& 2(2M+1)\delta+\|F(x')\|_2\|x-x'\|\bigg(2\delta+\frac{64\|w_n\|}{(\min\{\alpha(x),\beta(x)\})^2}\bigg).
\end{align}

In \eqref{L2.4.e4n}, taking $n\rightarrow\infty$, and then taking $\delta\rightarrow 0^{+}$, we obtain that for any $w_0\in O_T$ such that $\|w_0\|=1$, 
\begin{equation*}
    \|(F(x)-F(x'))w_0\|\leq  
    \frac{128\|x-x'\|}{\min\{\alpha(x),\beta(x)\}\min\{\alpha(x'),\beta(x')\}}.
\end{equation*}
By Theorem~\ref{L2.2} and Lemma~\ref{L2.3}, we have $F(x)V(x)=F(x')V(x')=0$, and so
$(F(x)-F(x'))V(x)=0$ (note that $V(x)=V(x')$). Hence
\begin{equation*}
    \|F(x)-F(x')\|_2\leq \frac{128\|x-x'\|}{\min\{\alpha(x),\beta(x)\}\min\{\alpha(x'),\beta(x')\}}.
\end{equation*}

In \eqref{L2.4.e12n}, taking $n\rightarrow\infty$, we get
\begin{equation*}
    \big\|\big(F(x)-F(x')\big(\mathbf{I}_d-\langle x-x',V(x)\rangle F(x)\big)\big)w_0\big\|\leq 2(2M+1)\delta+2\|F(x')\|_2\|x-x'\|\delta.
\end{equation*}
Taking $\delta\rightarrow 0^{+}$, we obtain that for any $w_0\in O_T$ with $\|w_0\|=1$,
\begin{equation*}
    \big(F(x)-F(x')\big(\mathbf{I}_d-\langle x-x',V(x)\rangle F(x)\big)\big)w_0=0.
\end{equation*}
By Theorem~\ref{L2.2} and Lemma~\ref{L2.3}, we have $F(x)V(x)=F(x')V(x)=0$. Hence
\begin{equation*}
    \big(F(x)-F(x')\big(\mathbf{I}_d-\langle x-x',V(x)\rangle F(x)\big)\big)V(x)=F(x)V(x)-F(x')V(x)=0.
\end{equation*}
Therefore,
\begin{equation*}
    F(x)=F(x')\big(\mathbf{I}_d-\langle x-x',V(x)\rangle F(x)\big),
\end{equation*}
and so
\begin{equation*}
    F(x') =  \big(\mathbf{I}_d-\langle x'-x, V_T\rangle F(x')\big)F(x).
\end{equation*} 

\paragraph{Step 4.} In this step, we will show that for any $x\in\mathrm{int}(T)$, if for some $F_0\in \mathbb{R}^{d\times d}$ with $F_0 V_T=0$, 
\begin{equation}\label{L2.4.eq17}
    \lim_{t\rightarrow 0^{+}}\frac{\mathcal{H}^{d-1}\big(\big\{w\in O_T:\|w\|\leq t, \|V(x+w)-V(x)-F_0 w\|>\delta t\big\}\big)}{t^{d-1}}=0
\end{equation}
holds for all $\delta>0$, then $\approxgrad V(x)$ exists and $\approxgrad V(x)=F(x)=F_0$.

In the following, we fix any $\delta\in (0,\min\{\alpha(x),\beta(x),1\}\slash 20)$. By \eqref{L2.4.eq17}, for any $\delta'>0$, there exists $t_0\in (0,\delta\slash 10)$ such that for any $t\in (0,t_0)$, we have 
\begin{equation}
    \mathcal{H}^{d-1}\big(\big\{w\in O_T:\|w\|\leq t, \|V(x+w)-V(x)-F_0 w\|\geq \delta t\slash 5\big\}\big)\leq\delta' t^{d-1}. 
\end{equation}

Below, we consider any $t\in (0,t_0)$ and any $x'\in\mathrm{int}(T)$ such that $\|x'-x\|\leq t$. Note that
\begin{equation}\label{alphalow}
    \min\{\alpha(x'),\beta(x')\}\geq \min\{\alpha(x),\beta(x)\}-t\geq \frac{19}{20}\min\{\alpha(x),\beta(x)\},
\end{equation}
where the last inequality uses $t<t_0<\delta<\min\{\alpha(x),\beta(x)\}\slash 20$. Arguing similarly as in Step~1 (with $x$ playing the role of $x'$, $x'$ playing the role of $x$, and $F_0$ playing the role of $F(x)$) and using \eqref{alphalow}, we obtain that for $t$ sufficiently small (depending on $x,\delta,\delta',F_0$),
\begin{align}\label{L2.4.enn1}
    & \mathcal{H}^{d-1}\big(\big\{w\in O_T:\|w\|\leq t,\|V(x'+w)-V(x')-F_0w\|>\delta t\big\}\big) \nonumber\\
    \leq\,& \frac{\pi^{(d-1)\slash 2}}{\Gamma((d+1) \slash 2)}\bigg(1-\bigg(1+\frac{35 t}{\min\{\alpha(x),\beta(x)\}}\bigg)^{-(d-1)}\bigg)t^{d-1}\nonumber\\
    &+\mathcal{H}^{d-1}(\mathcal{D}(x';t,\delta))+\bigg(1+\frac{70 t}{\min\{\alpha(x),\beta(x)\}}\bigg)^{d-1}\delta' t^{d-1}.
\end{align}
Moreover, we have
\begin{align}\label{L2.4.enn3}
    &\mathcal{H}^{d-1}\big(\big\{w\in O_T:\|w\|\leq t,\|V(x'+w)-V(x')-F_0 w\|>\sqrt{\delta}\|w\|\big\}\big)\nonumber\\
\leq\,&\mathcal{H}^{d-1}\big(\big\{w\in O_T:\|w\|\leq \sqrt{\delta}t\big\}\big)+\mathcal{H}^{d-1}\big(\big\{w\in O_T:\|w\|\leq t,\|V(x'+w)-V(x')-F_0 w\|>\delta t \big\}\big)\nonumber\\
=\,& \frac{\pi^{(d-1)\slash 2}}{\Gamma((d+1) \slash 2)}\delta^{(d-1)\slash 2} t^{d-1}+\mathcal{H}^{d-1}\big(\big\{w\in O_T:\|w\|\leq t,\|V(x'+w)-V(x')-F_0 w\|>\delta t \big\}\big).
\end{align}

Now note that for any $y\in B_d(x,t)$ such that $\|V(y)-V(x)-F_0(y-x)\|>\sqrt{\delta}\|y-x\|$, if $y=x'+w$ (where $x'\in  [x-tV(x), x+tV(x)]$ and $w\in O_T$), then $\|w\|\leq \|y-x\|\leq t$ and 
\begin{align}\label{L2.4.enn2}
     \|V(x'+w)-V(x')-F_0 w\|=\,&\|V(x'+w)-V(x)-F_0(x'-x+w)\|\nonumber\\
    =\,& \|V(y)-V(x)-F_0(y-x)\|>\sqrt{\delta}\|y-x\|\geq \sqrt{\delta}\|w\|,
\end{align}
where the first equality uses $F_0V(x)= F_0 V_T=0$. Combining \eqref{L2.4.enn1}--\eqref{L2.4.enn2}, we obtain 
\begin{align}\label{L2.4.enn4}
    &\mathcal{L}^d\big(\big\{y\in B_d(x,t): \|V(y)-V(x)-F_0(y-x)\|>\sqrt{\delta}\|y-x\|\big\}\big)\nonumber\\
    \leq\,&\int_{[x-tV(x),x+tV(x)]}
\mathcal{H}^{d-1}\big(\big\{w\in O_T:\|w\|\leq t,\|V(x'+w)-V(x')-F_0 w\|>\sqrt{\delta}\|w\|\big\}\big)d\mathcal{H}^1(x') \nonumber\\
\leq\,& \frac{2\pi^{(d-1)\slash 2}}{\Gamma((d+1) \slash 2)}\delta^{(d-1)\slash 2}t^d+2\bigg(1+\frac{70 t}{\min\{\alpha(x),\beta(x)\}}\bigg)^{d-1}\delta' t^d\nonumber\\
&+\frac{2\pi^{(d-1)\slash 2}}{\Gamma((d+1) \slash 2)}\bigg(1-\bigg(1+\frac{35 t}{\min\{\alpha(x),\beta(x)\}}\bigg)^{-(d-1)}\bigg)t^{d}\nonumber\\
&+\int_{[x-tV(x),x+tV(x)]}\mathcal{H}^{d-1}(\mathcal{D}(x';t,\delta))d\mathcal{H}^1(x').
\end{align}

For any $y\in\mathbb{R}^d$ such that $y=x'+w$ with $x'\in[x-tV(x),x+tV(x)]$ and $w\in\mathcal{D}(x';t,\delta)$, we have $\|w\|\leq t$, $\|y-x\|\leq\sqrt{2}t$, and 
\begin{align*}
  &\max\{|\alpha(y)-\alpha(x)|,|\beta(y)-\beta(x)|\}=\max\{|\alpha(x'+w)-\alpha(x)|,|\beta(x'+w)-\beta(x)|\}\nonumber\\
  \geq\,& \max\{|\alpha(x'+w)-\alpha(x')|,|\beta(x'+w)-\beta(x')|\}-t>\delta-t\geq \frac{\delta}{2}.  
\end{align*}
Hence we have
\begin{align}\label{L2.4.enn5}
   &\int_{[x-tV(x),x+tV(x)]}\mathcal{H}^{d-1}(\mathcal{D}(x';t,\delta))d\mathcal{H}^1(x')\nonumber\\
   =\,& \mathcal{L}^d\big(\big\{y\in\mathbb{R}^d: y=x'+w, x'\in [x-tV(x),x+tV(x)], w\in \mathcal{D}(x';t,\delta)\big\}\big)\nonumber\\
   \leq\,& \mathcal{L}^d\big(\big\{y\in B_d(x,\sqrt{2}t):\max\{|\alpha(y)-\alpha(x)|,|\beta(y)-\beta(x)|\}\geq \delta\slash 2\big\}\big).
\end{align}
Since the conclusion of Theorem~\ref{L2.2} holds for $T$, $\alpha$ and $\beta$ are approximately continuous at $x$, hence
\begin{equation}\label{L2.4.enn6}
    \lim_{t\rightarrow 0^{+}}\frac{\mathcal{L}^d\big(\big\{y\in B_d(x,\sqrt{2}t):\max\{|\alpha(y)-\alpha(x)|,|\beta(y)-\beta(x)|\}\geq \delta\slash 2\big\}\big)}{t^d} = 0.
\end{equation}
By \eqref{L2.4.enn5} and \eqref{L2.4.enn6}, we have
\begin{equation}\label{L2.4.enn7}
    \lim_{t\rightarrow 0^{+}} t^{-d}\int_{[x-tV(x),x+tV(x)]}\mathcal{H}^{d-1}(\mathcal{D}(x';t,\delta))d\mathcal{H}^1(x')=0. 
\end{equation}

By \eqref{L2.4.enn4} and \eqref{L2.4.enn7}, taking first $t\rightarrow 0^+$ and then $\delta'\rightarrow 0^{+}$, we obtain
\begin{equation*}
    \limsup_{t\rightarrow 0^{+}} t^{-d}\mathcal{L}^d\big(\big\{y\in B_d(x,t): \|V(y)-V(x)-F_0(y-x)\|>\sqrt{\delta}\|y-x\|\big\}\big)
    \leq \frac{2\pi^{(d-1)\slash 2}}{\Gamma((d+1) \slash 2)}\delta^{(d-1)\slash 2}. 
\end{equation*}
Hence for any $\delta_0>0$, when $\delta\in (0,\delta_0^2)$, we have
\begin{align*}
    &\limsup_{t\rightarrow 0^{+}} t^{-d}\mathcal{L}^d\big(\big\{y\in B_d(x,t): \|V(y)-V(x)-F_0(y-x)\|>\delta_0\|y-x\|\big\}\big)\nonumber\\
    \leq\,& \limsup_{t\rightarrow 0^{+}} t^{-d}\mathcal{L}^d\big(\big\{y\in B_d(x,t): \|V(y)-V(x)-F_0(y-x)\|>\sqrt{\delta}\|y-x\|\big\}\big)\nonumber\\
    \leq\,& \frac{2\pi^{(d-1)\slash 2}}{\Gamma((d+1) \slash 2)}\delta^{(d-1)\slash 2}.
\end{align*}
Taking $\delta\rightarrow 0^{+}$, we obtain that
\begin{equation*}
    \lim_{t\rightarrow 0^{+}} t^{-d}\mathcal{L}^d\big(\big\{y\in B_d(x,t): \|V(y)-V(x)-F_0(y-x)\|>\delta_0\|y-x\|\big\}\big)=0. 
\end{equation*}
Hence $\approxgrad V(x)$ exists and $\approxgrad V(x)=F(x)=F_0$. 

\paragraph{Step 5.} In this step, we will show that $\mathcal{B}(T)$ (recall \eqref{defB_T}) is a closed subset of $\mathrm{int}(T)$ under the subspace topology.

Suppose that $(x_n)_{n=1}^{\infty}\subseteq\mathcal{B}(T)$, $x\in\mathrm{int}(T)$, and $x_n\rightarrow x$ as $n\rightarrow\infty$. For any $n\in\mathbb{N}^*$, if $\approxgrad V(x_n)$ does not exist, then $F(x_n)=0$ by Definition~\ref{DefVF}, and so $F(x_n)V_T=0$; since $x_n\in\mathcal{B}(T)$, \eqref{Eq2.12} holds for $x_n$, and the result of Step~4 implies that $\approxgrad V(x_n)$ exists, which is a contradiction. Therefore, for every $n\in\mathbb{N}^*$, $\approxgrad V(x_n)$ exists and $\approxgrad V(x_n)=F(x_n)$. By \eqref{L2.4.eq18}, $\{F(x_n)\}_{n=1}^{\infty}$ is a Cauchy sequence, hence $\tilde{F}(x):=\lim_{n\rightarrow\infty}F(x_n)$ exists. 

Consider any $\delta\in(0, \min\{\alpha(x),\beta(x),1\}\slash 50)$, and take $\delta'\in (0,\delta\slash 10)$ sufficiently small. There exists $n_0\in\mathbb{N}^{*}$, such that for all $n\geq n_0$, $\|x_n-x\|\leq \delta'$ and $\|F(x_n)-\tilde{F}(x)\|_2\leq \delta'$. Below we fix $n\geq n_0$. As $x_n\in\mathcal{B}(T)$, there exists $t_0\in (0,\delta\slash 10)$ such that for all $t\in (0,t_0)$, $\mathcal{H}^{d-1}(\mathcal{W}(x_n;t,\delta))\leq \delta' t^{d-1}$. Following a similar argument as in Step~1 (with $x_n$ playing the role of $x'$), we can deduce that for $t\in (0,t_0)$, 
\begin{align}\label{L2.4.new1}
&\mathcal{H}^{d-1}\big(\big\{w\in O_T:\|w\|\leq t, \|V(x+w)-V(x)-\tilde{F}(x)w\|>2\delta t\big\}\big)\nonumber\\
\leq\,& \frac{\pi^{(d-1)\slash 2}}{\Gamma((d+1) \slash 2)}\bigg(1-\bigg(1+\frac{32\delta'}{\min\{\alpha(x),\beta(x)\}}\bigg)^{-(d-1)}\bigg)t^{d-1}\nonumber\\
    &+\mathcal{H}^{d-1}(\mathcal{D}(x;t,\delta))+\bigg(1+\frac{64 \delta'}{\min\{\alpha(x),\beta(x)\}}\bigg)^{d-1}\delta' t^{d-1}.
\end{align}
As the conclusion of Theorem~\ref{L2.2} holds for $T$, we have 
\begin{equation}\label{L2.4.new2}
    \lim_{t\rightarrow 0^{+}}\frac{\mathcal{H}^{d-1}(\mathcal{D}(x;t,\delta))}{t^{d-1}} = 0.
\end{equation}
By \eqref{L2.4.new1} and \eqref{L2.4.new2}, we have
\begin{align*}
   & \limsup_{t\rightarrow 0^{+}}\frac{\mathcal{H}^{d-1}\big(\big\{w\in O_T:\|w\|\leq t, \|V(x+w)-V(x)-\tilde{F}(x)w\|>2\delta t\big\}\big)}{t^{d-1}}\nonumber\\
    \leq\,& \frac{\pi^{(d-1)\slash 2}}{\Gamma((d+1) \slash 2)}\bigg(1-\bigg(1+\frac{32\delta'}{\min\{\alpha(x),\beta(x)\}}\bigg)^{-(d-1)}\bigg)+\bigg(1+\frac{64 \delta'}{\min\{\alpha(x),\beta(x)\}}\bigg)^{d-1}\delta'.
\end{align*}
Taking $\delta'\rightarrow 0^{+}$, we obtain that
\begin{equation}\label{limto0}
    \lim_{t\rightarrow 0^{+}}  \frac{\mathcal{H}^{d-1}\big(\big\{w\in O_T:\|w\|\leq t, \|V(x+w)-V(x)-\tilde{F}(x)w\|>2\delta t\big\}\big)}{t^{d-1}}=0.
\end{equation}

For every $n\in\mathbb{N}^*$, using \cref{L2.3} and noting the fact that $\approxgrad V(x_n)$ exists and $\alpha,\beta$ are approximately continuous at $x_n$ (since the conclusion of \cref{L2.2} holds for $T$), we obtain that $F(x_n)V_T=0$. As $\tilde{F}(x)=\lim_{n\rightarrow\infty}F(x_n)$, we have $\tilde{F}(x)V_T=0$. Consequently, using \eqref{limto0} and the result of Step 4, we conclude that $\approxgrad V(x)$ exists and $\approxgrad V(x)=F(x)=\tilde{F}(x)$. Using this and noting \eqref{limto0}, we conclude that \eqref{Eq2.12} holds for $x$, and so $x\in\mathcal{B}(T)$. Therefore, $\mathcal{B}(T)$ is a closed subset of $\mathrm{int}(T)$ under the subspace topology. 

\paragraph{Step 6.} Using the results from Steps~1 and~5, we conclude that $\mathrm{int}(T)\backslash \mathcal{B}(T)$ is an open subset of $\mathrm{int}(T)$ (under the subspace topology) with zero $\mathcal{H}^1$ measure. Hence $\mathrm{int}(T)\backslash \mathcal{B}(T)=\emptyset$ and $\mathcal{B}(T)=\mathrm{int}(T)$. Noting \eqref{defB_T}, we conclude that for every $x\in\mathrm{int}(T)$, \eqref{Eq2.12} holds. 

Now for any $x\in\mathrm{int}(T)$, if $\approxgrad V(x)$ does not exist, then $F(x)=0$ by \cref{DefVF}, and so $F(x)V_T=0$; as \eqref{Eq2.12} holds for $x$, the result of Step~4 implies that $\approxgrad V(x)$ exists, which is a contradiction. Hence $\approxgrad V(x)$ exists (and so $F(x)=\approxgrad V(x)$) for every $x\in \mathrm{int}(T)$. We have established part (a).

By part (a) and the results from Steps~2 and~3, we obtain \eqref{L2.4.enew3}, \eqref{L2.4.enew4}, and \eqref{Eq2.6.2}. Now for any $x,x',x_0\in\mathrm{int}(T)$, using \eqref{Eq2.6.2}, we obtain
\begin{align*}
   & \big(\mathbf{I}_d-\langle x'-x_0, V_T\rangle F(x')\big)-\big(\mathbf{I}_d-\langle x'-x, V_T\rangle F(x')\big)\big(\mathbf{I}_d-\langle x-x_0, V_T\rangle F(x)\big) \nonumber\\
   =\,& -\langle x-x_0, V_T\rangle F(x')+\langle x-x_0, V_T\rangle F(x)-\langle x'-x,V_T\rangle\langle x-x_0,V_T\rangle F(x') F(x)\nonumber\\
   =\,& -\langle x-x_0, V_T\rangle \big(F(x')-F(x)+\langle x'-x, V_T\rangle F(x')
   F(x) \big)=0. 
\end{align*}
This establishes \eqref{Eq2.6.4}. 

By part (a), Theorem~\ref{L2.2}, and Lemma~\ref{L2.3}, we have $F(x)V(x)=0$ for every $x\in\mathrm{int}(T)$. Now for any $x\in\mathrm{int}(T)$, $\delta\in (0,\min\{\alpha(x),\beta(x),1\}\slash 10)$, $t\in (0,\delta^3\slash 10)$, and $w\in \mathcal{D}'(x;t,\delta)\cap\mathcal{W}'(x;t,\delta)$, by \eqref{eq2.4.new} and Definition~\ref{Defd}, we have  
\begin{equation*}
    |\langle V(x+w)-V(x), V(x)\rangle|=|\langle V(x+w), V(x)\rangle-1|\leq \frac{32\|w\|^2}{(\min\{\alpha(x), \beta(x)\})^2},
\end{equation*}
\begin{equation*}
    \|V(x+w)-V(x)-F(x)w\|\leq \delta t. 
\end{equation*}
Hence
\begin{align}\label{L2.4.new1.1}
    |\langle F(x) w, V(x)\rangle|
    \leq\,& |\langle V(x+w)-V(x), V(x)\rangle|+|\langle V(x+w)-V(x)-F(x)w, V(x)\rangle|\nonumber\\
    \leq\,& |\langle V(x+w)-V(x), V(x)\rangle|+\|V(x+w)-V(x)-F(x)w\|\nonumber\\
    \leq\,& \frac{32\|w\|^2}{(\min\{\alpha(x), \beta(x)\})^2}+\delta t.
\end{align}
By Theorem~\ref{L2.2} and part (a),
\begin{equation*}
    \lim_{t\rightarrow 0^{+}} \frac{\mathcal{H}^{d-1}(\mathcal{D}(x;t,\delta)\cup \mathcal{W}(x;t,\delta))}{t^{d-1}} = 0.
\end{equation*}
Hence for any $w_0\in O_T$ with $\|w_0\|=1$, there exist $(t_n)_{n=1}^{\infty}\subseteq (0,1)$ and $(w_n)_{n=1}^{\infty}\subseteq O_T$ such that $w_n\in \mathcal{D}'(x;t_n,\delta)\cap\mathcal{W}'(x;t_n,\delta)$ for every $n\in\mathbb{N}^*$ and, as $n\rightarrow\infty$, 
\begin{equation*}
    t_n\rightarrow 0, \quad t_n  \slash  2\leq \|w_n\| \leq t_n, \quad \frac{w_n}{\|w_n\|}\rightarrow w_0. 
\end{equation*}
By \eqref{L2.4.new1.1}, for $n$ sufficiently large, we have
\begin{equation*}
    \Big|\Big\langle F(x)\frac{w_n}{\|w_n\|},V(x)\Big\rangle\Big|\leq \frac{32\|w_n\|}{(\min\{\alpha(x),\beta(x)\})^2}+\frac{\delta t_n}{\|w_n\|}\leq \frac{32\|w_n\|}{(\min\{\alpha(x),\beta(x)\})^2}+2\delta. 
\end{equation*}
Taking $n\rightarrow\infty$ and noting \eqref{L2.4.enew3}, we obtain $|\langle F(x)w_0, V(x)\rangle|\leq 2\delta$. Taking $\delta\rightarrow 0^{+}$, we conclude that $w_0^{\top} F(x)^{\top}V(x)=\langle F(x) w_0, V(x)\rangle=0$ for any $w_0\in O_T$ with $\|w_0\|=1$. Now for any $z\in \mathbb{R}^d$, we can write $z=\lambda_1 V(x)+\lambda_2 w_0$ for some $\lambda_1,\lambda_2\in\mathbb{R}$ and $w_0\in O_T$ with $\|w_0\|=1$. Thus using the fact that $F(x)V(x)=0$, we obtain
\begin{equation*}
    z^{\top} F(x)^{\top}V(x)=\lambda_1 (F(x)V(x))^{\top} V(x)+\lambda_2 w_0^{\top}F(x)^{\top}V(x)=0.
\end{equation*}
Hence $F(x)^{\top}V(x)=0$ for every $x\in\mathrm{int}(T)$. This establishes \eqref{Eq2.6.5}. 
\end{proof}

\subsection{Proof of \Cref{L2.6}}\label{Sect.3.1.5}

Recall \Cref{Defd}. For $x$ in the interior of a transport ray $T$, the following lemma bounds the difference $u(x+\tilde{w})-u(x+w)$ for $w,\tilde{w}\in O_T$.  

\begin{lemma}\label{prep1.1}
Fix any transport ray $T\in\mathcal{S}$ and any $x\in\mathrm{int}(T)$. For any $\delta\in \Big(0,\frac{\min\{\alpha(x),\beta(x)\}}{6}\Big)$, $t\in \Big(0,\frac{\min\{\alpha(x),\beta(x)\}}{16} \cdot \min\Big\{1,\frac{\delta}{\max\{\alpha(x),\beta(x)\}}\Big\}\Big)$, and $w, \tilde{w}\in O_T$ such that $\|\tilde{w}-w\|\leq\frac{1}{2}\min\{\alpha(x),\beta(x)\}$, if $w\in\mathcal{D}'(x;t,\delta)\cap\mathcal{W}'(x;t,\delta)$, then
\begin{align*}
    \big|u(x+\tilde{w})-u(x+w)-\langle F(x)w,\tilde{w}-w\rangle\big|\leq\,& 2|\langle  F(x)w,\tilde{w}-w\rangle|\cdot\frac{\|\tilde{w}-w\|}{\min\{\alpha(x),\beta(x)\}}\nonumber\\
    &+2\delta t\|\tilde{w}-w\|+\frac{2\|\tilde{w}-w\|^2}{\min\{\alpha(x),\beta(x)\}}.
\end{align*}
\end{lemma}
\begin{proof}

\noindent\textbf{Step 1.} We first derive an upper bound on $u(x+\tilde{w})-u(x+w)$. As $w\in\mathcal{D}'(x;t,\delta)$, we have
\begin{equation*}
    \min\{\alpha(x+w),\beta(x+w)\}\geq \min\{\alpha(x),\beta(x)\}-\delta\geq \frac{1}{2}\min\{\alpha(x),\beta(x)\}.
\end{equation*}
Hence by \Cref{L2.0}, 
\begin{equation}\label{Vdiff}
    \|V(x+w)-V(x)\|\leq \frac{8}{\min\{\alpha(x),\beta(x)\}}\|w\|\leq \frac{8t}{\min\{\alpha(x),\beta(x)\}}\leq \frac{1}{2},
\end{equation}
which implies 
\begin{equation}\label{psoti2}
   \langle V(x+w),V(x)\rangle =1-\frac{\|V(x+w)-V(x)\|^2}{2}\geq \frac{1}{2}.     
\end{equation}
Now let $x':=b_T+3\delta V(x)\in\mathrm{int}(T)$ (recall that $b_T$ is the lower end of $T$). Note that
\begin{equation}\label{eqneww}
    \langle x-x',V(x)\rangle=\beta(x)-3\delta\geq \frac{1}{2}\min\{\alpha(x),\beta(x)\}>0, \text{ so } \|x-x'\|=\langle x-x',V(x)\rangle.
\end{equation}
We define (note \eqref{psoti2})
\begin{equation}\label{def_wp}
    w':=x-x'+w-\frac{\langle x-x', V(x)\rangle}{\langle V(x+w),V(x)\rangle}V(x+w).
\end{equation}
Note that $w'\in O_T$ and $x'+w'\in L(x+w)$. Moreover, by \eqref{Vdiff}--\eqref{eqneww},
\begin{align*}
   \bigg|\frac{\langle x-x', V(x)\rangle}{\langle V(x+w),V(x)\rangle}-\langle x-x',V(x)\rangle\bigg|=\,& \frac{|\langle x-x', V(x)\rangle||\langle V(x+w),V(x)\rangle-1|}{\langle V(x+w),V(x)\rangle}\nonumber\\
   \leq\,& 2\beta(x)\|V(x+w)-V(x)\|\leq \frac{16 \beta(x) t}{\min\{\alpha(x),\beta(x)\}}\leq \delta.
\end{align*}
Hence using \eqref{eqneww} and the fact that  $w\in\mathcal{D}'(x;t,\delta)$, we obtain
\begin{equation*}
    \langle (x'+w')-(x+w) , V(x+w)\rangle=-\frac{\|x-x'\|}{\langle V(x+w),V(x)\rangle}\leq 0,
\end{equation*}
\begin{align*}
    \beta(x+w)+\langle (x'+w')-(x+w), V(x+w)\rangle =\,&\beta(x+w)-  \frac{\langle x-x', V(x)\rangle}{\langle V(x+w),V(x)\rangle}\nonumber\\
    \geq\,&(\beta(x)-\delta)-(\langle x-x',V(x)\rangle+\delta)=\delta>0.
\end{align*}
Consequently, $x'+w'\in\mathrm{int}(T(x+w))$ and $u(x+w)-u(x'+w')=\|(x+w)-(x'+w')\|$. Hence
\begin{align*}\label{u_dif}
    u(x+\tilde{w})-u(x+w)&=(u(x+\tilde{w})-u(x'+w'))-(u(x+w)-u(x'+w'))\nonumber\\
    & \leq \|(x+\tilde{w})-(x'+w')\|-\|(x+w)-(x'+w')\|.
\end{align*}
By \eqref{def_wp} and \eqref{psoti2}, we have
\begin{equation}\label{diff2}
    \|(x+w)-(x'+w')\|=\frac{\|x-x'\|}{\langle V(x+w),V(x)\rangle}.
\end{equation}
Moreover, 
\begin{align*}
  &\|(x+\tilde{w})-(x'+w')\|^2  =\bigg\| \tilde{w} -w + \frac{\langle x-x',V(x)\rangle}{\langle V(x+w),V(x)\rangle}V(x+w)\bigg\|^2\nonumber\\
    =\, & \|\tilde{w}-w\|^2+\bigg(\frac{\|x-x'\|}{\langle V(x+w), V(x)\rangle}\bigg)^2+\frac{2\|x-x'\|}{\langle V(x+w),V(x)\rangle}\langle V(x+w),\tilde{w}-w\rangle.
\end{align*}
Combining the above three displays and noting \eqref{eqneww}, we get
\begin{align*}
    u(x+\tilde{w})-u(x+w)&\leq \frac{\|(x+\tilde{w})-(x'+w')\|^2-\|(x+w)-(x'+w')\|^2}{\|(x+\tilde{w})-(x'+w')\|+\|(x+w)-(x'+w')\|}\nonumber\\
  &= \frac{\frac{2\|x-x'\|}{\langle V(x+w),V(x)\rangle}\langle V(x+w),\tilde{w}-w\rangle+\|\tilde{w}-w\|^2}{\|(x+\tilde{w})-(x'+w')\|+\|(x+w)-(x'+w')\|}.
\end{align*}
As $w\in \mathcal{W}'(x;t,\delta)$, we have $\|V(x+w)-V(x)-F(x)w\|\leq \delta t$; thus, noting that $\langle V(x),\tilde{w}-w\rangle =0$ (as $w,\tilde{w}\in O_T$), we obtain
\begin{align}\label{vwx2n}
     &|\langle V(x+w),\tilde{w}-w\rangle-\langle F(x)w,\tilde{w}-w\rangle|=|\langle V(x+w)-V(x)-F(x)w,\tilde{w}-w\rangle|\nonumber\\
     \leq\,&\|V(x+w)-V(x)-F(x)w\|\|\tilde{w}-w\|\leq \delta t \|\tilde{w}-w\|.
\end{align}
By the above two displays, we obtain
\begin{equation}\label{uwww}
    u(x+\tilde{w})-u(x+w)\leq \frac{\frac{2\|x-x'\|}{\langle V(x+w),V(x)\rangle}\big(\langle F(x)w,\tilde{w}-w\rangle+\delta t \|\tilde{w}-w\|\big)+\|\tilde{w}-w\|^2}{\|(x+\tilde{w})-(x'+w')\|+\|(x+w)-(x'+w')\|}. 
\end{equation}

By the triangle inequality and \eqref{diff2}, 
\begin{equation*}
    \|(x+\tilde{w})-(x'+w')\|
    \begin{cases}
        \geq \|(x+w)-(x'+w')\|-\|\tilde{w}-w\|=\frac{\|x-x'\|}{\langle V(x+w),V(x)\rangle}-\|\tilde{w}-w\|,\\
        \leq  \|(x+w)-(x'+w')\|+\|\tilde{w}-w\|=\frac{\|x-x'\|}{\langle V(x+w),V(x)\rangle}+\|\tilde{w}-w\|.
    \end{cases}
\end{equation*}
By \eqref{eqneww}, $\|\tilde{w}-w\|\leq \frac{1}{2}\min\{\alpha(x),\beta(x)\}\leq\|x-x'\|$. Hence by the above display, 
\begin{align*}
 & \frac{1}{\|(x+\tilde{w})-(x'+w')\|+\|(x+w)-(x'+w')\|}\nonumber\\
 \leq\,&\frac{1}{\frac{2\|x-x'\|}{\langle V(x+w),V(x)\rangle}-\|\tilde{w}-w\|}
 =\frac{\langle V(x+w), V(x) \rangle}{2\|x-x'\|}\bigg(1-\frac{\|\tilde{w}-w\|\langle V(x+w), V(x) \rangle}{2\|x-x'\|}\bigg)^{-1}\nonumber\\
 \leq\,&\frac{\langle V(x+w), V(x)\rangle}{2\|x-x'\|}\bigg(1+\frac{\|\tilde{w}-w\|}{\|x-x'\|}\bigg)\leq\frac{\langle V(x+w), V(x)\rangle}{2\|x-x'\|}\bigg(1+\frac{2\|\tilde{w}-w\|}{\min\{\alpha(x),\beta(x)\}}\bigg),
\end{align*} 
where the second to last inequality uses $(1-s)^{-1}\leq 1+2s$ for $s\in [0,1/2]$. Similarly,
\begin{align*}
 & \frac{1}{\|(x+\tilde{w})-(x'+w')\|+\|(x+w)-(x'+w')\|}\nonumber\\
 \geq\,&\frac{1}{\frac{2\|x-x'\|}{\langle V(x+w),V(x)\rangle}+\|\tilde{w}-w\|}
 =\frac{\langle V(x+w), V(x) \rangle}{2\|x-x'\|}\bigg(1+\frac{\|\tilde{w}-w\|\langle V(x+w), V(x) \rangle}{2\|x-x'\|}\bigg)^{-1}\nonumber\\
 \geq\,&\frac{\langle V(x+w), V(x)\rangle}{2\|x-x'\|}\bigg(1-\frac{\|\tilde{w}-w\|}{\|x-x'\|}\bigg)\geq \frac{\langle V(x+w), V(x)\rangle}{2\|x-x'\|}\bigg(1-\frac{2\|\tilde{w}-w\|}{\min\{\alpha(x),\beta(x)\}}\bigg),
\end{align*} 
where the second to last inequality uses $(1+s)^{-1}\geq 1-s$ for $s\geq 0$. By \eqref{diff2} and \eqref{eqneww}, 
\begin{align*}
    &\frac{\frac{2\|x-x'\|}{\langle V(x+w),V(x)\rangle}\cdot\delta t \|\tilde{w}-w\|+\|\tilde{w}-w\|^2}{\|(x+\tilde{w})-(x'+w')\|+\|(x+w)-(x'+w')\|}\leq\frac{\frac{2\|x-x'\|}{\langle V(x+w),V(x)\rangle}\cdot\delta t \|\tilde{w}-w\|+\|\tilde{w}-w\|^2}{\frac{\|x-x'\|}{\langle V(x+w),V(x)\rangle}}\nonumber\\
    \leq\,&2\delta t\|\tilde{w}-w\|+\frac{\|\tilde{w}-w\|^2}{\|x-x'\|}\leq 2\delta t\|\tilde{w}-w\|+\frac{2\|\tilde{w}-w\|^2}{\min\{\alpha(x),\beta(x)\}}.
\end{align*} 
By \eqref{uwww} and the above three displays, we obtain 
\begin{align}\label{ubbww}
    u(x+\tilde{w})-u(x+w)\leq\,&\langle F(x)w,\tilde{w}-w\rangle +2|\langle  F(x)w,\tilde{w}-w\rangle|\cdot\frac{\|\tilde{w}-w\|}{\min\{\alpha(x),\beta(x)\}}\nonumber\\
    &+2\delta t\|\tilde{w}-w\|+\frac{2\|\tilde{w}-w\|^2}{\min\{\alpha(x),\beta(x)\}}.
\end{align}

\bigskip

\noindent \textbf{Step 2.} Now we derive a lower bound on $u(x+\tilde{w})-u(x+w)$. Let $x'':=a_T-3\delta V(x)\in\mathrm{int}(T)$ (recall that $a_T$ is the upper end of $T$). Note that
\begin{equation}\label{eqneww2}
    \langle x''-x,V(x)\rangle=\alpha(x)-3\delta\geq \frac{1}{2}\min\{\alpha(x),\beta(x)\}>0, \text{ so } \|x-x''\|=\langle x''-x,V(x)\rangle.
\end{equation}
We define (note \eqref{psoti2})
\begin{equation}\label{def_wp2}
    w'':=x-x''+w-\frac{\langle x-x'', V(x)\rangle}{\langle V(x+w),V(x)\rangle}V(x+w).
\end{equation}
Note that $w''\in O_T$ and $x''+w''\in L(x+w)$. Arguing similarly as in Step 1, we can deduce that $x''+w''\in\mathrm{int}(T(x+w))$ and $u(x''+w'')-u(x+w)=\|(x''+w'')-(x+w)\|$. Hence 
\begin{align*}
   u(x+\tilde{w})-u(x+w)&=(u(x''+w'')-u(x+w))-(u(x''+w'')-u(x+\tilde{w}))\nonumber\\
    &\geq\|(x''+w'')-(x+w)\|-\|(x''+w'')-(x+\tilde{w})\|.
\end{align*}
By \eqref{def_wp2} and \eqref{psoti2}, we have
\begin{equation*}
    \|(x''+w'')-(x+w)\|=\frac{\|x-x''\|}{\langle V(x+w),V(x)\rangle}. 
\end{equation*}
Moreover, 
\begin{align*}
  &\|(x''+w'')-(x+\tilde{w})\|^2  =\bigg\|\frac{\langle x''-x,V(x)\rangle}{\langle V(x+w),V(x)\rangle}V(x+w)-(\tilde{w}-w)\bigg\|^2\nonumber\\
    =\, & \|\tilde{w}-w\|^2+\bigg(\frac{\|x-x''\|}{\langle V(x+w),V(x)\rangle}\bigg)^2-\frac{2\|x-x''\|}{\langle V(x+w),V(x)\rangle}\langle V(x+w),\tilde{w}-w\rangle.
\end{align*}
Combining the above three displays and noting \eqref{eqneww2} and \eqref{vwx2n}, we get
\begin{align*}
    u(x+\tilde{w})-u(x+w)&\geq \frac{\|(x''+w'')-(x+w)\|^2-\|(x''+w'')-(x+\tilde{w})\|^2}{\|(x''+w'')-(x+\tilde{w})\|+\|(x''+w'')-(x+w)\|}\nonumber\\
  & =  \frac{\frac{2\|x-x''\|}{\langle V(x+w),V(x)\rangle}\langle V(x+w),\tilde{w}-w\rangle-\|\tilde{w}-w\|^2}{\|(x''+w'')-(x+\tilde{w})\|+\|(x''+w'')-(x+w)\|}\nonumber\\
  &\geq \frac{\frac{2\|x-x''\|}{\langle V(x+w),V(x)\rangle}\big(\langle F(x)w,\tilde{w}-w\rangle-\delta t \|\tilde{w}-w\|\big)-\|\tilde{w}-w\|^2}{\|(x''+w'')-(x+\tilde{w})\|+\|(x''+w'')-(x+w)\|}.
\end{align*}
Using this and arguing similarly as in Step 1, we obtain 
\begin{align}\label{lbbww}
    u(x+\tilde{w})-u(x+w)\geq\,&\langle F(x)w,\tilde{w}-w\rangle-2|\langle  F(x)w,\tilde{w}-w\rangle|\cdot\frac{\|\tilde{w}-w\|}{\min\{\alpha(x),\beta(x)\}}\nonumber\\
    &-2\delta t\|\tilde{w}-w\|-\frac{2\|\tilde{w}-w\|^2}{\min\{\alpha(x),\beta(x)\}}.
\end{align}
The conclusion of the lemma follows from \eqref{ubbww} and \eqref{lbbww}.
\end{proof}

Using \Cref{prep1.1}, we derive the following lemma bounding the error in the approximation $u(x+w)\approx u(x)+\frac{1}{2}w^{\top}F(x) w$ for $x$ in the interior of a transport ray $T$ and $w\in O_T$. 

\begin{lemma}\label{prep1.2}
Fix any transport ray $T\in\mathcal{S}$ and any $x\in\mathrm{int}(T)$ such that $\|F(x)\|_2\leq \frac{4}{\min\{\alpha(x),\beta(x)\}}$. For any $\delta\in \Big(0,\frac{\min\{\alpha(x),\beta(x)\}}{6}\Big)$, $t\in \Big(0,\frac{\min\{\alpha(x),\beta(x)\}}{16} \cdot \min\Big\{1,\frac{\delta}{\max\{\alpha(x),\beta(x)\}}\Big\}\Big)$, $w\in\mathcal{D}'(x;t,\delta)\cap\mathcal{W}'(x;t,\delta)$, and $\delta'\in \big(0,\frac{1}{4}\big)$ such that
\begin{equation}\label{criterion2}
\mathcal{L}^1\big(\big\{\theta\in [0,1]:\theta w\in\mathcal{D}(x;t,\delta)\cup\mathcal{W}(x;t,\delta)\big\}\big)\leq \delta',
\end{equation}
we have
\begin{equation}\label{desired2}
    \bigg|u(x+w)-u(x)-\frac{1}{2}w^{\top} F(x) w\bigg|  \leq  \frac{40\delta'\|w\|^2}{\min\{\alpha(x),\beta(x)\}}+2\delta t\|w\|.
\end{equation}
\end{lemma}
\begin{proof}
Let $I_0:=\lceil  (2\delta')^{-1} \rceil-1$. For each $i\in [I_0-1]$, we set $B_i:=\big[2(i-1)\delta',2i\delta'\big)$. We also set $B_{I_0}:=\big[2(I_0-1)\delta',1\big]$. Note that $\mathcal{L}^1(B_i)\in [2\delta',4\delta']$ for every $i\in [I_0]$. Consequently, by \eqref{criterion2}, for every $i\in [I_0]$, there exists $\theta_i\in B_i$ such that $\theta_i w\notin \mathcal{D}(x;t,\delta)\cup\mathcal{W}(x;t,\delta)$; as $\|\theta _i w\|\leq \|w\|\leq t$, we have $\theta_i w\in \mathcal{D}'(x;t,\delta)\cap\mathcal{W}'(x;t,\delta)$. We also let $\theta_0:=0$ and $\theta_{I_0+1}:=1$; note that $\theta_i w\in \mathcal{D}'(x;t,\delta)\cap\mathcal{W}'(x;t,\delta)$ for $i\in \{0,I_0+1\}$. For any $i\in \{0\}\cup [I_0]$, we have
\begin{equation*}
    \|\theta_i w-\theta_{i+1}w\|=|\theta_i-\theta_{i+1}|\|w\|\leq 8\delta'\|w\|\leq 8\delta' t \leq \frac{1}{8}\min\{\alpha(x),\beta(x)\}.
\end{equation*}
Hence by \Cref{prep1.1} (with $w,\tilde{w}$ replaced by $\theta_i w,\theta_{i+1} w$),
\begin{align*}
&\big|u(x+\theta_{i+1}w)-u(x+\theta_i w)-\langle F(x)(\theta_i w), (\theta_{i+1}-\theta_i)w\rangle \big|\nonumber\\
\leq\,& 2 |\langle F(x) (\theta_i w),(\theta_{i+1}-\theta_i)w\rangle|\cdot\frac{\|\theta_{i+1}w-\theta_iw\|}{\min\{\alpha(x),\beta(x)\}}
+2\delta t \|\theta_{i+1}w-\theta_i w\|+\frac{2\|\theta_{i+1}w-\theta_i w\|^2}{\min\{\alpha(x),\beta(x)\}}\nonumber\\
\leq\,& \frac{2\|F(x)\|_2\|w\|^3}{\min\{\alpha(x),\beta(x)\}}\cdot(\theta_{i+1}-\theta_i)^2+2\delta t\|w\|(\theta_{i+1}-\theta_i)+\frac{2\|w\|^2}{\min\{\alpha(x),\beta(x)\}}\cdot (\theta_{i+1}-\theta_i)^2\nonumber\\
\leq\,& \frac{3\|w\|^2}{\min\{\alpha(x),\beta(x)\}}\cdot (\theta_{i+1}-\theta_i)^2+2\delta t\|w\|(\theta_{i+1}-\theta_i),
\end{align*}
where the last inequality uses $\|F(x)\|_2\|w\|\leq \frac{4}{\min\{\alpha(x),\beta(x)\}}\cdot t\leq \frac{1}{4}$. Using the above display along with the observation that $\sum_{i=0}^{I_0}(\theta_{i+1}-\theta_i)^2\leq 8\delta'\sum_{i=0}^{I_0}(\theta_{i+1}-\theta_i)=8\delta'$, we get
\begin{align*}
   & \bigg|u(x+w)-u(x)-\sum_{i=0}^{I_0}\langle F(x)(\theta_i w),(\theta_{i+1}-\theta_i)w\rangle\bigg|\nonumber\\
   \leq\,& \sum_{i=0}^{I_0}\big|u(x+\theta_{i+1} w)-u(x+\theta_i w)-\langle F(x)(\theta_i w),(\theta_{i+1}-\theta_i)w\rangle\big|\nonumber\\
   \leq\,& \frac{3\|w\|^2}{\min\{\alpha(x),\beta(x)\}}\cdot\sum_{i=0}^{I_0}(\theta_{i+1}-\theta_i)^2+2\delta t\|w\|\sum_{i=0}^{I_0}(\theta_{i+1}-\theta_i)\leq\frac{24\delta'\|w\|^2}{\min\{\alpha(x),\beta(x)\}}+2\delta t\|w\|.
\end{align*}
Since $\sum_{i=0}^{I_0}(\theta_{i+1}+\theta_i)(\theta_{i+1}-\theta_i)=\sum_{i=0}^{I_0}(\theta_{i+1}^2-\theta_i^2)=1$, we have
\begin{align*}
    &\bigg|\sum_{i=0}^{I_0}\langle F(x)(\theta_i w),(\theta_{i+1}-\theta_i)w\rangle-       \frac{1}{2}w^{\top}F(x)w\bigg|=|w^{\top}F(x)w|\bigg|\sum_{i=0}^{I_0}\theta_i (\theta_{i+1}-\theta_i)-\frac{1}{2}\bigg|\nonumber\\
    \leq\,& \|F(x)\|_2\|w\|^2\bigg|\sum_{i=0}^{I_0}\theta_i(\theta_{i+1}-\theta_i)-\frac{1}{2}\sum_{i=0}^{I_0}(\theta_{i+1}+\theta_i)(\theta_{i+1}-\theta_i)\bigg|\nonumber\\
    =\,&\frac{1}{2}\|F(x)\|_2\|w\|^2\sum_{i=0}^{I_0}(\theta_{i+1}-\theta_i)^2\leq \frac{16\delta' \|w\|^2}{\min\{\alpha(x),\beta(x)\}}.
\end{align*}
The conclusion of the lemma follows from the above two displays. 
\end{proof}

\begin{proof}[Proof of \Cref{L2.6}]

Consider $T\in\mathcal{S}$ for which the conclusions of \cref{L2.2} and \cref{L2.4}(a)--(b) hold (note that by Sections~\ref{Sect.3.1.2} and~\ref{Sect.3.1.3}, this is the case for $T=T(z_0)$ for $\mathcal{L}^d$-a.e.\ $z_0\in\mathcal{T}_1^*$). Fix any $x\in\mathrm{int}(T)$. By \eqref{Eq2.1n} and \eqref{Eq2.12}, for any $\delta>0$ and $\delta'\in (0,1)$, there exists $t_0\in (0,1)$ such that for any $t\in (0,t_0)$, 
\begin{equation}\label{L3.4.eq}
    \mathcal{H}^{d-1}\big(\mathcal{D}(x;t,\delta)\cup\mathcal{W}(x;t,\delta)\big)\leq \delta' t^{d-1}. 
\end{equation}
For any $t>0$, denote $O_{T,t}:=\{w\in O_T:\|w\|\leq t\}$. For any $t\in (0,t_0)$ and $\delta_0\in (0,1)$, we have
\begin{align}\label{Lem3.4.eq1}
  & \mathcal{H}^{d-1}\otimes\mathcal{L}^1\big(\big\{(w,\theta)\in O_{T,t}\times [0,1]:\theta w\in \mathcal{D}(x;t,\delta)\cup\mathcal{W}(x;t,\delta)\big\}\big)\nonumber\\
  \leq\,& \mathcal{H}^{d-1}\otimes\mathcal{L}^1\big(\big\{(w,\theta)\in O_{T,t}\times (\delta_0,1]:\theta w\in \mathcal{D}(x;t,\delta)\cup\mathcal{W}(x;t,\delta)\big\}\big)+\mathcal{H}^{d-1}\otimes\mathcal{L}^1(O_{T,t}\times [0,\delta_0])\nonumber\\
  \leq\,& \int_{(\delta_0,1]}\mathcal{H}^{d-1}\big(\theta^{-1}(\mathcal{D}(x;t,\delta)\cup\mathcal{W}(x;t,\delta))\big)d\theta + \frac{\pi^{(d-1)\slash 2}}{\Gamma((d+1)\slash 2)}\delta_0 t^{d-1}\nonumber\\
  \leq\,&\delta_0^{-(d-1)}\delta' t^{d-1}+\frac{\pi^{(d-1)\slash 2}}{\Gamma((d+1)\slash 2)}\delta_0 t^{d-1},
\end{align}
where the last inequality uses the fact that for any $\theta\in (\delta_0,1]$, by \eqref{L3.4.eq},
\begin{equation*}
    \mathcal{H}^{d-1}\big(\theta^{-1}(\mathcal{D}(x;t,\delta)\cup\mathcal{W}(x;t,\delta))\big)=\theta^{-(d-1)}\mathcal{H}^{d-1}\big(\mathcal{D}(x;t,\delta)\cup\mathcal{W}(x;t,\delta)\big)\leq \delta_0^{-(d-1)}\delta' t^{d-1}.
\end{equation*}
Taking $\delta_0=(\delta')^{1\slash d}$ in \eqref{Lem3.4.eq1}, we obtain
\begin{equation*}
    \mathcal{H}^{d-1}\otimes\mathcal{L}^1\big(\big\{(w,\theta)\in O_{T,t}\times [0,1]:\theta w\in \mathcal{D}(x;t,\delta)\cup\mathcal{W}(x;t,\delta)\big\}\big)\leq C_d(\delta')^{1\slash d} t^{d-1},
\end{equation*}
where $C_d$ depends only on $d$. Hence
\begin{equation*}
    \limsup_{t\rightarrow 0^{+}} \frac{\mathcal{H}^{d-1}\otimes\mathcal{L}^1\big(\big\{(w,\theta)\in O_{T,t}\times [0,1]:\theta w\in \mathcal{D}(x;t,\delta)\cup \mathcal{W}(x;t,\delta)\big\}\big)}{t^{d-1}}\leq C_d(\delta')^{1\slash d}.
\end{equation*}
Taking $\delta'\rightarrow 0^{+}$, we conclude that for any $\delta>0$, 
\begin{equation*}
    \lim_{t\rightarrow 0^{+}} \frac{\mathcal{H}^{d-1}\otimes\mathcal{L}^1\big(\big\{(w,\theta)\in O_{T,t}\times [0,1]:\theta w\in \mathcal{D}(x;t,\delta)\cup \mathcal{W}(x;t,\delta)\big\}\big)}{t^{d-1}}=0.
\end{equation*}

Below, we fix any $\delta\in\Big(0,\frac{\min\{\alpha(x),\beta(x)\}}{6}\Big)$ and $\delta'\in \big(0,\frac{1}{4}\big)$. For any $t>0$, let $\mathscr{Q}_{\delta,\delta',t}$ be the set of $w\in O_{T,t}$ such that $\mathcal{L}^1\big(\big\{\theta\in [0,1]:\theta w\in\mathcal{D}(x;t,\delta)\cup\mathcal{W}(x;t,\delta)\big\}\big)\leq \delta'$. Note that
\begin{equation*}
    \mathcal{H}^{d-1}\big(O_{T,t}\backslash \mathscr{Q}_{\delta,\delta',t}\big)\leq (\delta')^{-1}
    \mathcal{H}^{d-1}\otimes\mathcal{L}^1\big(\big\{(w,\theta)\in O_{T,t}\times [0,1]:\theta w\in\mathcal{D}(x;t,\delta)\cup\mathcal{W}(x;t,\delta)\big\}\big).
\end{equation*}
Using the above two displays, along with \eqref{Eq2.1n} and \eqref{Eq2.12}, we obtain that 
\begin{equation*}
    \lim_{t\to 0^+}\frac{\mathcal{H}^{d-1}\big(O_{T,t}\backslash \mathscr{Q}_{\delta,\delta',t}\big)}{t^{d-1}}=0,\qquad \lim_{t\to 0^+}\frac{\mathcal{H}^{d-1}\big(\mathcal{D}(x;t,\delta)\cup\mathcal{W}(x;t,\delta)\big)}{t^{d-1}}=0.
\end{equation*}
Fix any $\delta''>0$. Using the above display and noting \Cref{Defd}, we obtain the existence of $t_1\in \Big(0,\frac{\min\{\alpha(x),\beta(x)\}}{16} \cdot \min\Big\{1,\frac{\delta}{\max\{\alpha(x),\beta(x)\}}\Big\}\Big)$ such that for any $t\in (0,t_1)$,
\begin{equation*}
    \mathcal{H}^{d-1}\Big(O_{T,t}\big\backslash \big(\mathscr{Q}_{\delta,\delta',t}\cap\mathcal{D}'(x;t,\delta)\cap\mathcal{W}'(x;t,\delta)\big)\Big)\leq \delta'' t^{d-1}.
\end{equation*}
For any $t\in (0,t_1)$ and $w\in\mathscr{Q}_{\delta,\delta',t}\cap\mathcal{D}'(x;t,\delta)\cap\mathcal{W}'(x;t,\delta)$, by \Cref{prep1.2}, we have 
\begin{equation*}
    \bigg|u(x+w)-u(x)-\frac{1}{2}w^{\top} F(x) w\bigg|  \leq  \frac{40\delta'\|w\|^2}{\min\{\alpha(x),\beta(x)\}}+2\delta t\|w\|\leq \bigg(\frac{40\delta'}{\min\{\alpha(x),\beta(x)\}}+2\delta\bigg)t^2.
\end{equation*}
Now for any $\delta_0>0$, take $\delta,\delta'$ sufficiently small such that $\frac{40\delta'}{\min\{\alpha(x),\beta(x)\}}+2\delta\leq \delta_0$. By the above two displays, for any $t\in (0,t_1)$, 
\begin{align*}
   &\mathcal{H}^{d-1}\Big(  \Big\{ w\in O_T:\|w\|\leq t,\Big|u(x+w)-u(x)-\frac{1}{2}w^{\top} F(x) w\Big|>\delta_0 t^2 \Big\}\Big)\nonumber\\
   \leq\,& \mathcal{H}^{d-1}\Big(O_{T,t}\big\backslash \big(\mathscr{Q}_{\delta,\delta',t}\cap\mathcal{D}'(x;t,\delta)\cap\mathcal{W}'(x;t,\delta)\big)\Big)\leq \delta'' t^{d-1}.
\end{align*}
Hence
\begin{equation*}
    \limsup_{t\rightarrow 0^{+}} \frac{\mathcal{H}^{d-1}\big(\big\{w\in O_{T}:\|w\|\leq t,\big|u(x+w)-u(x)-\frac{1}{2}w^{\top} F(x) w\big|>\delta_0 t^2\big\}\big)}{t^{d-1}}\leq \delta''.
\end{equation*}
Taking $\delta''\to 0^+$, we conclude that
\begin{equation*}
    \lim_{t\rightarrow 0^{+}} \frac{\mathcal{H}^{d-1}\big(\big\{w\in O_{T}:\|w\|\leq t,\big|u(x+w)-u(x)-\frac{1}{2}w^{\top} F(x) w\big|>\delta_0 t^2\big\}\big)}{t^{d-1}}=0.\qedhere
\end{equation*}
\end{proof}

\subsection{Proof of Theorem~\ref{L2.4}(c)}\label{Sect.3.1.4}

We begin with the following lemma. For $x$ in the interior of a transport ray $T$ and $w,w'\in O_T$, the lemma bounds the error in the approximation $u(x+w')-u(x+w)\approx \frac{1}{2}\langle F(x)(w+w'),w'-w\rangle$.  

\begin{lemma}\label{Lemmaww2}
Fix any transport ray $T\in\mathcal{S}$ and any $x\in\mathrm{int}(T)$ such that $\|F(x)\|_2\leq \frac{4}{\min\{\alpha(x),\beta(x)\}}$. For any $\delta\in \Big(0,\frac{\min\{\alpha(x),\beta(x)\}}{6}\Big)$, $t\in \Big(0,\frac{\min\{\alpha(x),\beta(x)\}}{16} \cdot \min\Big\{1,\frac{\delta}{\max\{\alpha(x),\beta(x)\}}\Big\}\Big)$, $w,w'\in\mathcal{D}'(x;t,\delta)\cap\mathcal{W}'(x;t,\delta)$, and $\delta'\in \big(0,\frac{1}{4}\big)$ such that
\begin{equation}\label{ca2}
\mathcal{L}^1\big(\big\{\theta\in [0,1]:  w+\theta(w'-w)   \in\mathcal{D}(x;t,\delta)\cup\mathcal{W}(x;t,\delta)\big\}\big)\leq \delta',
\end{equation}
we have
\begin{equation}\label{desired_conclusion}
    \bigg|u(x+w')-u(x+w)-\frac{1}{2}
    \langle F(x) (w+w'),w'-w\rangle    
\bigg|  \leq  \frac{40\delta'\|w'-w\|^2}{\min\{\alpha(x),\beta(x)\}}+2\delta t\|w'-w\|.
\end{equation}
\end{lemma}
\begin{proof}
Let $I_0:=\lceil  (2\delta')^{-1} \rceil-1$. For each $i\in [I_0-1]$, we set $B_i:=\big[2(i-1)\delta',2i\delta'\big)$. We also set $B_{I_0}:=\big[2(I_0-1)\delta',1\big]$. Note that $\mathcal{L}^1(B_i)\in [2\delta',4\delta']$ for every $i\in [I_0]$. Consequently, by \eqref{ca2}, for every $i\in [I_0]$, there exists $\theta_i\in B_i$ such that $w+\theta_i(w'-w)\notin \mathcal{D}(x;t,\delta)\cup\mathcal{W}(x;t,\delta)$; as $\|w+\theta_i(w'-w)\|=\|(1-\theta_i)w+\theta_i w'\|\leq \max\{\|w\|,\|w'\|\}\leq t$, we have $w+\theta_i(w'-w)\in \mathcal{D}'(x;t,\delta)\cap\mathcal{W}'(x;t,\delta)$. We also let $\theta_0:=0$ and $\theta_{I_0+1}:=1$; note that $w+\theta_i(w'-w)\in \mathcal{D}'(x;t,\delta)\cap\mathcal{W}'(x;t,\delta)$ for $i\in \{0,I_0+1\}$. For any $i\in \{0\}\cup [I_0]$, we have
\begin{equation*}
    \big\|(w+\theta_i(w'-w))-(w+\theta_{i+1}(w'-w))\big\|=|\theta_i-\theta_{i+1}|\|w'-w\| \leq (8\delta')\cdot (2t) \leq \frac{1}{4}\min\{\alpha(x),\beta(x)\}.
\end{equation*}
Hence by \Cref{prep1.1} (with $w,\tilde{w}$ replaced by $w+\theta_i(w'-w),w+\theta_{i+1}(w'-w)$),
\begin{align*}
&\big|u(x+w+\theta_{i+1}(w'-w))-u(x+w+\theta_i(w'-w))\nonumber\\
&\hspace{0.05in}-\langle F(x)(w+\theta_i(w'-w)), (\theta_{i+1}-\theta_i)(w'-w)\rangle \big|\nonumber\\
\leq\,& 2 |\langle F(x) (w + \theta_i (w'-w)),(\theta_{i+1}-\theta_i)(w'-w)\rangle|\cdot\frac{\|(\theta_{i+1}-\theta_i)(w'-w)\|}{\min\{\alpha(x),\beta(x)\}}\nonumber\\
&+2\delta t \|(\theta_{i+1}-\theta_i)(w'-w)\|+\frac{2\|(\theta_{i+1}-\theta_i)(w'-w)\|^2}{\min\{\alpha(x),\beta(x)\}}\nonumber\\
\leq\,& \frac{2\|F(x)\|_2\max\{\|w\|,\|w'\|\}\|w'-w\|^2}{\min\{\alpha(x),\beta(x)\}}\cdot(\theta_{i+1}-\theta_i)^2\nonumber\\
&+2\delta t\|w'-w\|(\theta_{i+1}-\theta_i)+\frac{2\|w'-w\|^2}{\min\{\alpha(x),\beta(x)\}}\cdot (\theta_{i+1}-\theta_i)^2\nonumber\\
\leq\,& \frac{3\|w'-w\|^2}{\min\{\alpha(x),\beta(x)\}}\cdot (\theta_{i+1}-\theta_i)^2+2\delta t\|w'-w\|(\theta_{i+1}-\theta_i),
\end{align*}
where the last inequality uses $\|F(x)\|_2\max\{\|w\|,\|w'\|\}\leq\frac{4}{\min\{\alpha(x),\beta(x)\}}\cdot t\leq \frac{1}{4}$. Using the above display along with the observation that $\sum_{i=0}^{I_0}(\theta_{i+1}-\theta_i)^2\leq 8\delta'\sum_{i=0}^{I_0}(\theta_{i+1}-\theta_i)=  8\delta'$, we get
\begin{align}\label{Ldisplay1}
&\bigg|u(x+w')-u(x+w)-\sum_{i=0}^{I_0}\langle F(x)(w+\theta_i(w'-w)),(\theta_{i+1}-\theta_i)(w'-w)\rangle\bigg|\nonumber\\
   \leq\,& \sum_{i=0}^{I_0}\big|u(x+w+\theta_{i+1}(w'-w))-u(x+w+\theta_i(w'-w))\nonumber\\
&\hspace{0.3in}-\langle F(x)(w+\theta_i(w'-w)), (\theta_{i+1}-\theta_i)(w'-w)\rangle \big|\nonumber\\
   \leq\,& \frac{3\|w'-w\|^2}{\min\{\alpha(x),\beta(x)\}}\cdot\sum_{i=0}^{I_0}(\theta_{i+1}-\theta_i)^2+2\delta t\|w'-w\|\sum_{i=0}^{I_0}(\theta_{i+1}-\theta_i)\nonumber\\
   \leq\,& \frac{24\delta'\|w'-w\|^2}{\min\{\alpha(x),\beta(x)\}}+2\delta t\|w'-w\|,
\end{align}
Since $\sum_{i=0}^{I_0}(\theta_{i+1}+\theta_i)(\theta_{i+1}-\theta_i)=\sum_{i=0}^{I_0}(\theta_{i+1}^2-\theta_i^2)=1$, we have
\begin{equation*}
    \bigg|\sum_{i=0}^{I_0}\theta_i (\theta_{i+1}-\theta_i)-\frac{1}{2}\bigg|=\bigg|\sum_{i=0}^{I_0}\theta_i(\theta_{i+1}-\theta_i)-\frac{1}{2}\sum_{i=0}^{I_0}(\theta_{i+1}+\theta_i)(\theta_{i+1}-\theta_i)\bigg|\leq\frac{1}{2}\sum_{i=0}^{I_0}(\theta_{i+1}-\theta_i)^2\leq 4\delta'.
\end{equation*}
Hence
\begin{align}\label{Ldisplay2}
    &\bigg|\sum_{i=0}^{I_0}\langle F(x)(w+\theta_i(w'-w)),(\theta_{i+1}-\theta_i)(w'-w)\rangle -       \frac{1}{2}\langle F(x)(w+w'),w'-w\rangle\bigg|\nonumber\\
    =\,& \bigg|\sum_{i=0}^{I_0}\langle F(x)(w'-w),w'-w\rangle\cdot\theta_i(\theta_{i+1}-\theta_i)-       \frac{1}{2}\langle F(x)(w'-w),w'-w\rangle\bigg|\nonumber\\
    =\,& |\langle F(x)(w'-w),w'-w\rangle|\bigg|\sum_{i=0}^{I_0}\theta_i (\theta_{i+1}-\theta_i)-\frac{1}{2}\bigg|\nonumber\\
    \leq\,&  4\delta'\|F(x)\|_2\|w'-w\|^2\leq\frac{16\delta'\|w'-w\|^2}{\min\{\alpha(x),\beta(x)\}}.
\end{align}
The conclusion of the lemma follows from \eqref{Ldisplay1} and \eqref{Ldisplay2}. 
\end{proof}

\begin{proof}[Proof of \Cref{L2.4}(c)]

Fix any $T\in\mathcal{S}$ for which the conclusions of \Cref{L2.2}, \Cref{L2.4}(a)--(b), and \Cref{L2.6} hold (note that by Sections~\ref{Sect.3.1.2}--\ref{Sect.3.1.5}, this is the case for $T=T(z_0)$ for $\mathcal{L}^d$-a.e.\ $z_0\in\mathcal{T}_1^*$). Fix any $x\in  \mathrm{int}(T)$. By \eqref{Eq2.1n} and \eqref{Eq2.12}, for any $\delta>0$ and $\delta'\in (0,1)$, there exists $t_0\in (0,1)$ such that for any $t\in (0,t_0)$, 
\begin{equation}\label{L3.1.eq7n}
    \mathcal{H}^{d-1}\big(\mathcal{D}(x;t,\delta)\cup\mathcal{W}(x;t,\delta)\big)\leq \delta' t^{d-1}. 
\end{equation}
For any $t>0$, denote $O_{T,t}:=\{w\in O_T:\|w\|\leq t\}$. Define
\begin{equation*}
    \mathscr{E}_{\delta,t}:=\big\{(w,w',\theta)\in O_{T,t}\times O_{T,t}\times [0,1]:(1-\theta)w+\theta w'\in\mathcal{D}(x;t,\delta)\cup\mathcal{W}(x;t,\delta)\big\}.
\end{equation*}
For any $t\in (0, t_0)$ and $\delta_0\in (0,1)$, we have 
\begin{align}\label{L3.1.eq8}
   & \mathcal{H}^{d-1}\otimes \mathcal{H}^{d-1}\otimes \mathcal{L}^1(\mathscr{E}_{\delta,t})\nonumber\\
   \leq\,& \int_{O_{T,t}}d\mathcal{H}^{d-1}(w)\int_{(\delta_0,1]}d\theta \int_{O_{T,t}} \mathbbm{1}_{(1-\theta) w+ \theta w'\in \mathcal{D}(x;t,\delta)\cup\mathcal{W}(x;t,\delta)}d\mathcal{H}^{d-1}(w')\nonumber\\
   &+\mathcal{H}^{d-1}\otimes\mathcal{H}^{d-1}\otimes \mathcal{L}^1\big(O_{T,t}\times O_{T,t}\times [0,\delta_0]\big)\nonumber\\
   \leq\,&\int_{O_{T,t}}d\mathcal{H}^{d-1}(w)\int_{(\delta_0,1]}
   \mathcal{H}^{d-1}\Big(\theta^{-1}\big((\mathcal{D}(x;t,\delta)\cup\mathcal{W}(x;t,\delta))-(1-\theta)w\big)\Big) 
   d\theta \nonumber\\
   &+\bigg(\frac{\pi^{(d-1)\slash 2}}{\Gamma((d+1)\slash 2)}\bigg)^2\delta_0 t^{2(d-1)}\leq \frac{\pi^{(d-1)\slash 2}}{\Gamma((d+1)\slash 2)}\delta_0^{-(d-1)}\delta' t^{2(d-1)}+ \frac{\pi^{d-1}\delta_0}{\Gamma((d+1)\slash 2)^2}t^{2(d-1)},
\end{align}
where the last inequality uses the fact that by \eqref{L3.1.eq7n}, for any $w\in O_{T,t}$ and $\theta\in (\delta_0,1]$, 
\begin{align*}
   & \mathcal{H}^{d-1}\Big(\theta^{-1}\big((\mathcal{D}(x;t,\delta)\cup\mathcal{W}(x;t,\delta))-(1-\theta) w\big)\Big)\nonumber\\
   =\,& \theta^{-(d-1)}\mathcal{H}^{d-1}\big(\mathcal{D}(x;t,\delta)\cup\mathcal{W}(x;t,\delta)\big) \leq \delta_0^{-(d-1)}\delta' t^{d-1}.
\end{align*}
Taking $\delta_0=(\delta')^{1\slash d}$ in \eqref{L3.1.eq8}, we obtain that $\mathcal{H}^{d-1}\otimes \mathcal{H}^{d-1}\otimes \mathcal{L}^1(\mathscr{E}_{\delta,t})\leq C_d(\delta')^{1\slash d}t^{2(d-1)}$, where $C_d>0$ depends only on $d$. Hence
\begin{align*}
    \limsup_{t\rightarrow 0^{+}} \frac{\mathcal{H}^{d-1}\otimes \mathcal{H}^{d-1}\otimes \mathcal{L}^1(\mathscr{E}_{\delta,t})}{t^{2(d-1)}}\leq   C_d(\delta')^{1\slash d}. 
\end{align*}
Taking $\delta'\rightarrow 0^{+}$, we conclude that for any $\delta>0$, 
\begin{equation}\label{Eqn1.1n}
  \lim_{t\rightarrow 0^{+}} \frac{\mathcal{H}^{d-1}\otimes \mathcal{H}^{d-1}\otimes \mathcal{L}^1(\mathscr{E}_{\delta,t})}{t^{2(d-1)}}=0.
\end{equation}

Below, we fix any $\delta\in\Big(0,\frac{\min\{\alpha(x),\beta(x)\}}{6}\Big)$ and $\delta'\in \big(0,\frac{1}{4}\big)$. For any $t>0$, let $\mathscr{R}_{\delta,\delta',t}$ be the set of $(w,w')\in O_{T,t}\times O_{T,t}$ such that $\mathcal{L}^1\big(\big\{\theta\in [0,1]:(1-\theta) w+\theta w'\in\mathcal{D}(x;t,\delta)\cup\mathcal{W}(x;t,\delta)\big\}\big)\leq \delta'$. Note that
\begin{equation*}
    \mathcal{H}^{d-1}\otimes\mathcal{H}^{d-1}\big((O_{T,t}\times O_{T,t})\backslash \mathscr{R}_{\delta,\delta',t}\big)\leq (\delta')^{-1}
    \mathcal{H}^{d-1}\otimes \mathcal{H}^{d-1}\otimes\mathcal{L}^1(\mathscr{E}_{\delta,t}).
\end{equation*}
Define $\tilde{\mathscr{R}}_{\delta,\delta',t}:=\mathscr{R}_{\delta,\delta',t}\cap\big\{(w,w')\in O_{T,t}\times O_{T,t}: w,w'\in\mathcal{D}'(x;t,\delta)\cap\mathcal{W}'(x;t,\delta)\big\}$. We have
\begin{align}\label{Eqn1.2n}
&\mathcal{H}^{d-1}\otimes\mathcal{H}^{d-1}\big((O_{T,t}\times O_{T,t})\backslash \tilde{\mathscr{R}}_{\delta,\delta',t}\big)\nonumber\\
\leq\,&\mathcal{H}^{d-1}\otimes\mathcal{H}^{d-1}\big((O_{T,t}\times O_{T,t})\backslash \mathscr{R}_{\delta,\delta',t}\big)+\mathcal{H}^{d-1}\otimes\mathcal{H}^{d-1}\big(O_{T,t}\times (\mathcal{D}(x;t,\delta)\cup\mathcal{W}(x;t,\delta))\big)\nonumber\\
&+\mathcal{H}^{d-1}\otimes\mathcal{H}^{d-1}\big((\mathcal{D}(x;t,\delta)\cup\mathcal{W}(x;t,\delta))\times O_{T,t}\big)\nonumber\\
\leq\,&(\delta')^{-1}
    \mathcal{H}^{d-1}\otimes \mathcal{H}^{d-1}\otimes\mathcal{L}^1(\mathscr{E}_{\delta,t})+\frac{2\pi^{(d-1)/2}}{\Gamma((d+1)/2)}t^{d-1} \mathcal{H}^{d-1}\big(\mathcal{D}(x;t,\delta)\cup\mathcal{W}(x;t,\delta)\big).
\end{align}
Using \eqref{Eqn1.1n} and \eqref{Eqn1.2n}, along with \eqref{Eq2.1n} and \eqref{Eq2.12}, we obtain that 
\begin{equation*}
    \lim_{t\to 0^+}\frac{\mathcal{H}^{d-1}\otimes\mathcal{H}^{d-1}\big((O_{T,t}\times O_{T,t})\backslash \tilde{\mathscr{R}}_{\delta,\delta',t}\big)}{t^{2(d-1)}}=0.
\end{equation*}
Fix any $\delta''>0$. Using the above display, we obtain the existence of $t_1 \in \Big(0,\frac{\min\{\alpha(x),\beta(x)\}}{16} \cdot \min\Big\{1,\frac{\delta}{\max\{\alpha(x),\beta(x)\}}\Big\}\Big)$ such that for any $t\in (0,t_1)$,
\begin{equation*}
    \mathcal{H}^{d-1}\otimes\mathcal{H}^{d-1}\big((O_{T,t}\times O_{T,t})\backslash \tilde{\mathscr{R}}_{\delta,\delta',t}\big)\leq \delta'' t^{2(d-1)}.
\end{equation*}
For any $t\in (0,t_1)$ and $(w,w')\in \tilde{\mathscr{R}}_{\delta,\delta',t}$, by \Cref{Lemmaww2}, we have 
\begin{align*}
    \bigg|u(x+w')-u(x+w)-\frac{1}{2}
    \langle F(x) (w+w'),w'-w\rangle    
\bigg|  &\leq  \frac{40\delta'\|w'-w\|^2}{\min\{\alpha(x),\beta(x)\}}+2\delta t\|w'-w\|\nonumber\\
 &\leq \bigg(\frac{160\delta'}{\min\{\alpha(x),\beta(x)\}}+4\delta\bigg)t^2.
\end{align*}
Now for any $\delta_0>0$, take $\delta,\delta'$ sufficiently small such that $\frac{160\delta'}{\min\{\alpha(x),\beta(x)\}}+4\delta\leq \delta_0$. By the above two displays, for any $t\in (0,t_1)$,
\begin{align*}
   &\mathcal{H}^{d-1}\otimes\mathcal{H}^{d-1}\Big(\Big\{(w,w')\in O_{T,t}\times O_{T,t}:\Big|u(x+w')-u(x+w)-\frac{1}{2}
    \langle F(x) (w+w'),w'-w\rangle\Big|>\delta_0 t^2\Big\}\Big)\nonumber\\
   &\leq\mathcal{H}^{d-1}\otimes\mathcal{H}^{d-1}\big((O_{T,t}\times O_{T,t})\backslash \tilde{\mathscr{R}}_{\delta,\delta',t}\big)\leq \delta'' t^{2(d-1)}.
\end{align*}
Hence
{\small
\begin{equation*}
    \limsup_{t\rightarrow 0^{+}} \frac{\mathcal{H}^{d-1}\otimes\mathcal{H}^{d-1}\big(\big\{(w,w')\in O_{T,t}\times O_{T,t}:\big|u(x+w')-u(x+w)-\frac{1}{2}
    \langle F(x) (w+w'),w'-w\rangle\big|>\delta_0 t^2\big\}\big)}{t^{2(d-1)}}\leq \delta''.
\end{equation*}
}Taking $\delta''\to 0^+$, we conclude that for any $\delta_0>0$,
{\small
\begin{equation}\label{res_ww}
   \lim_{t\rightarrow 0^{+}} \frac{\mathcal{H}^{d-1}\otimes\mathcal{H}^{d-1}\big(\big\{(w,w')\in O_{T,t}\times O_{T,t}:\big|u(x+w')-u(x+w)-\frac{1}{2}
    \langle F(x) (w+w'),w'-w\rangle\big|>\delta_0 t^2\big\}\big)}{t^{2(d-1)}}=0.
\end{equation}}

Now by \Cref{L2.6}, for any $\delta_0>0$,
\begin{equation*}
        \lim_{t\rightarrow 0^{+}} \frac{\mathcal{H}^{d-1}\big(\big\{w\in O_{T,t}:\big|u(x+w)-u(x)-\frac{1}{2}w^{\top} F(x) w\big|>\delta_0 t^2/2 \big\}\big)}{t^{d-1}}=0,
\end{equation*}
and an analogous equation with $w$ replaced by $w'$ holds. Hence
{\small
\begin{equation}\label{res_ww2}
    \lim_{t\rightarrow 0^+}\frac{\mathcal{H}^{d-1}\otimes\mathcal{H}^{d-1}\big(\big\{(w,w')\in O_{T,t}\times O_{T,t}:\big|u(x+w')-u(x+w)-\frac{1}{2}(w')^{\top}F(x)w'+\frac{1}{2}w^{\top}F(x)w\big|>\delta_0 t^2\big\}\big)}{t^{2(d-1)}}=0.
\end{equation}
}Combining \eqref{res_ww} and \eqref{res_ww2}, and noting that
\begin{equation*}
    \langle F(x) (w+w'),w'-w\rangle-\big((w')^{\top}F(x)w'-w^{\top}F(x)w\big)=\langle F(x)w,w'\rangle-\langle F(x)w',w\rangle,
\end{equation*}
we get for any $\delta_0>0$, 
\begin{equation}\label{L3.1.eq22}
    \lim_{t\rightarrow 0^{+}} \frac{\mathcal{H}^{d-1}\otimes\mathcal{H}^{d-1}\big(\big\{(w,w')\in O_{T,t}\times O_{T,t}:\big|\langle F(x)w,w'\rangle-\langle F(x)w',w\rangle\big|>4\delta_0 t^2\big\}\big)}{t^{2(d-1)}}=0.
\end{equation}

By \eqref{L3.1.eq22}, for any $w_0,w_0'\in O_T$ such that $\|w_0\|=\|w_0'\|=1$ and any $\delta_0>0$, there exist $(t_n)_{n=1}^{\infty}\subseteq (0,1)$ and $(w_n)_{n=1}^{\infty},(w_n')_{n=1}^{\infty}\subseteq O_T$ such that as $n\rightarrow\infty$,
\begin{equation*}
    t_n\rightarrow 0,\quad t_n\slash 2\leq \|w_n\|,\|w_n'\|\leq t_n,\quad \frac{w_n}{\|w_n\|}\rightarrow w_0,\quad \frac{w_n'}{\|w_n'\|}\rightarrow w_0',  
\end{equation*}
\begin{equation*}
    \big|\langle F(x)w_n,w_n'\rangle -\langle F(x)w_n',w_n\rangle\big|\leq 4\delta_0 t_n^2.
\end{equation*}
Consequently,
\begin{equation*}
    \bigg|\bigg\langle F(x)\frac{w_n}{\|w_n\|},\frac{w_n'}{\|w_n'\|}\bigg\rangle-\bigg\langle F(x)\frac{w_n'}{\|w_n'\|},\frac{w_n}{\|w_n\|}\bigg\rangle\bigg|\leq \frac{4\delta_0t_n^2}{\|w_n\|\|w_n'\|}\leq 16\delta_0. 
\end{equation*}
Taking $n\rightarrow\infty$, we obtain (note \eqref{L2.4.enew3}) $\big|\langle F(x)w_0,w_0'\rangle-\langle F(x)w_0',w_0\rangle\big|\leq 16\delta_0$. Since $\delta_0>0$ is arbitrary, we have $\langle F(x)w_0,w_0'\rangle=\langle F(x)w_0',w_0\rangle$. Combining this with \eqref{Eq2.6.5} yields $F(x)^{\top}=F(x)$ for all $x\in\mathrm{int}(T)$. This completes the proof of Theorem~\ref{L2.4}(c). 
\end{proof}

\section{Proof of Theorem \ref{th:lowerBound}}\label{Sec3}

In this section, we present the proof of the variational lower bound. Section~\ref{Sec3:1} introduces the detachment $E(x,y)$ of the Kantorovich potential as well as the functions $G_{\eps}$ and $G$ (cf.\ Section~\ref{se:lower}), and establishes the main properties of $G_{\eps}$ and $G$. Section~\ref{Sec3:2} derives explicit formulas for $\lambda,\tilde{\mu}_T,\tilde{\nu}_T,\tilde{f},\tilde{g}$ (cf.\ Section~\ref{Sect.1.1}) and analyzes their key properties. Finally, Section~\ref{Sec3:3} completes the proof of Theorem~\ref{th:lowerBound}.

\subsection{Convergence of $G_{\eps}$ to $G$}\label{Sec3:1}

We introduce the detachment $E(x,y)$ of the potential $u$ and the functions $G_{\eps},G$ as follows. 

\begin{definition}[Detachment]\label{Def1}
For any $x,y\in\mathbb{R}^d$, we define
\begin{equation*}
    E(x,y):=\|x-y\|-u(x)+u(y). 
\end{equation*}
Since $|u(x)-u(y)|\leq \|x-y\|$ for all $x,y\in\mathbb{R}^d$, we have $E(x,y)\geq 0$.  
\end{definition}

\begin{definition}[The functions $G_{\eps}$ and $G$]\label{Def2}
For any $\eps>0$, $x\in\mathcal{T}_1^{*}$, and $z\in L(x)$, we define
\begin{align}
     G_{\eps}(x,z):=\,&\log\bigg(\varepsilon^{-(d-1)\slash 2}\int_{O(V(x))}e^{-E(x,z+w)\slash\varepsilon}\mathbbm{1}_{\|w\|\leq 2D} d\mathcal{H}^{d-1}(w)\bigg)\nonumber\\
     =\,& \log\bigg(\int_{O(V(x))} e^{-E(x,z+\sqrt{\eps}w)\slash \eps} \mathbbm{1}_{\|w\|\leq 2D\eps^{-1\slash 2}} d\mathcal{H}^{d-1}(w)\bigg). 
\end{align}
For any $x\in\mathcal{T}_1^{*}$ and $z\in \mathrm{int}(T(x))$ such that $\langle V(x),x-z\rangle>0$, we define
\begin{equation}\label{Eq3.33}
    G(x,z):=\log\bigg(\int_{O(V(x))}\exp\bigg(-\frac{1}{2}w^{\top}\bigg(\frac{1}{\|x-z\|}\mathbf{I}_d+F(z)\bigg)w\bigg) d\mathcal{H}^{d-1}(w)\bigg)
\end{equation}
if $\mathbf{I}_d+\|x-z\|F(z)$ is positive definite, and set $G(x,z):=0$ if $\mathbf{I}_d+\|x-z\|F(z)$ is not positive definite. For any $x\in\mathcal{T}_1^{*}$ and $z\in L(x)$, we set $G(x,z):=-\infty$ if $z\notin\mathrm{int}(T(x))$, or if $z\in\mathrm{int}(T(x))$ and
$\langle V(x),x-z\rangle\leq 0$.
\end{definition}

The main results of this subsection are summarized in the following two propositions. \Cref{P-lb-1} shows that $G$ is the pointwise limit of $G_{\eps}$ as $\eps\to 0^+$ under suitable assumptions. 

\begin{proposition}[Convergence of $G_{\eps}$ to $G$]\label{P-lb-1}
For $\mathcal{L}^d$-a.e.\ $x\in\mathcal{T}_1^{*}$, the following holds:
\begin{itemize}
    \item[(a)] For any $z\in \mathrm{int}(T(x))$ such that $z\neq x$, $\lim_{\eps\rightarrow 0^{+}}G_{\eps}(x,z)=G(x,z)$.
    \item[(b)] For any $z\in L(x)\backslash T(x)$, if either $\langle z-a(x), V(x)\rangle>0$, or $\langle z-b(x),V(x)\rangle<0$ and $z\in\mathcal{Y}$, then $\lim_{\eps\rightarrow 0^{+}}G_{\eps}(x,z)=G(x,z)$. 
\end{itemize}
\end{proposition}

\Cref{P-lb-2} collects key properties of $G_{\eps},G$ that will be used later in the proof of Theorem~\ref{th:lowerBound}: part (a) provides uniform upper bounds on $G_{\eps},G$ in terms of the log-distance to the ray end; part (b) gives an explicit formula for $G$; and part (c) describes continuity properties of $G$.

\begin{proposition}[Properties of $G_{\eps}$ and $G$]\label{P-lb-2}
For $\mathcal{L}^d$-a.e.\ $x\in\mathcal{T}_1^{*}$, the following holds:
\begin{itemize}
    \item[(a)] For any $\eps>0$ and $z\in L(x)$ such that $\max\{\|x-z\|,\|z-a(x)\|\}\leq 2D$, we have
    \begin{equation*}
        G_{\eps}(x,z)\leq -\frac{d-1}{2}\log(\alpha(x))+C;
    \end{equation*}
    for any $z\in L(x)$, we have 
    \begin{equation*}
        G(x,z)\leq -\frac{d-1}{2}\log(\alpha(x))+C,
    \end{equation*}
   where $C>0$ is a constant that only depends on $D,d$.    
    \item[(b)] For any $z\in\mathrm{int}(T(x))$ such that $\langle x-z,V(x)\rangle>0$, we have 
\begin{equation}\label{L3.6.2.eq1}
    G(x,z)= \frac{d-1}{2}\log(2\pi\|x-z\|)-\frac{1}{2}\log\det\big(\mathbf{I}_d+\|x-z\|F(z)\big),
\end{equation}
where $\mathbf{I}_d+\|x-z\|F(z)$ is positive definite.
\item[(c)] For any $t\in (-\beta(x),0)$, we have $\lim_{s\rightarrow t}G(x,x+sV(x))=G(x,x+tV(x))$. Moreover, $\lim\limits_{t\rightarrow 0^{-}}G(x,x+tV(x))=-\infty$ and $\lim\limits_{t\rightarrow (-\beta(x))^{+}}G(x,x+tV(x))$ exists in $[-\infty,\infty)$. 
\end{itemize}    
\end{proposition}

The remainder of this subsection is devoted to the proofs of \Cref{P-lb-1,P-lb-2}.

\subsubsection{Proof of \Cref{P-lb-1}}

We begin with uniform estimates on $E(x,y)$. 
\begin{lemma}\label{L3.2}
Let $x\in\mathcal{T}_1^{*}$ and $y\in\mathbb{R}^d$, and let $a$ denote the upper end of $T(x)$. Then 
\begin{equation}\label{eq:oldL3.1}
    E(x,y)\geq   \|y-x\|+\|x-a\|-\|y-a\|.
\end{equation} 
Moreover, let $z\in L(x)$ and $w\in O_{T(x)}$. If $\langle V(x), a-z\rangle\leq 0$, then $E(x,z+w)\geq \|x-a\|$, whereas if $\langle V(x), a-z\rangle\geq 0$, then
\begin{equation*}
    E(x,z+w)\geq \frac{\|x-a\|\|w\|^2}{4(\max\{\|z-a\|,\|z-x\|\}+\|w\|)^2}.
\end{equation*}
\end{lemma}

\begin{proof}
As $|u(y)-u(a)|\leq \|y-a\|$, we have $u(y)\geq u(a)-\|y-a\|$. Using $u(a)-u(x)=\|a-x\|$,
\begin{equation*}
    E(x,y)=\|y-x\| -u(x)+u(y)
    \geq\|y-x\|-\|y-a\|+u(a)-u(x)= \|y-x\|- \|y-a\|+\|x-a\|.
\end{equation*}
This shows~\eqref{eq:oldL3.1}. If $\langle V(x), x-z\rangle\geq 0$, by \eqref{eq:oldL3.1}, we have 
\begin{align*}
     E(x,z+w)\geq\,&\|z+w-x\|+\|x-a\|-\|z+w-a\|\nonumber\\
    =\,&\sqrt{\|z-x\|^2+\|w\|^2}+\|x-a\|-\sqrt{\|z-a\|^2+\|w\|^2}\nonumber\\
    =\,& \|x-a\|-\frac{\|x-a\|(\|z-x\|+\|z-a\|)}{\sqrt{\|z-x\|^2+\|w\|^2}+\sqrt{\|z-a\|^2+\|w\|^2}}\nonumber\\
    =\,&\|x-a\|\|w\|^2\frac{\frac{1}{\sqrt{\|z-x\|^2+\|w\|^2}+\|z-x\|}+\frac{1}{\sqrt{\|z-a\|^2+\|w\|^2}+\|z-a\|}}{\sqrt{\|z-x\|^2+\|w\|^2}+\sqrt{\|z-a\|^2+\|w\|^2}}\geq\frac{\|x-a\|\|w\|^2}{4(\|z-a\|+\|w\|)^2}.
\end{align*}
If $\langle V(x), x-z\rangle<0$ and $\langle V(x), a-z\rangle\geq 0$, by \eqref{eq:oldL3.1}, we have 
\begin{align*}
     E(x,z+w)\geq\,&\|z+w-x\|+\|x-a\|-\|z+w-a\|\nonumber\\
    =\,&\sqrt{\|z-x\|^2+\|w\|^2}+\|x-a\|-\sqrt{\|z-a\|^2+\|w\|^2}\nonumber\\
    =\,& \|x-a\|+\frac{\|x-a\|(\|z-x\|-\|z-a\|)}{\sqrt{\|z-x\|^2+\|w\|^2}+\sqrt{\|z-a\|^2+\|w\|^2}}\nonumber\\
    =\,& \|x-a\|\frac{\sqrt{\|z-x\|^2+\|w\|^2}+\|z-x\|+\frac{\|w\|^2}{\sqrt{\|z-a\|^2+\|w\|^2}+\|z-a\|}}{\sqrt{\|z-x\|^2+\|w\|^2}+\sqrt{\|z-a\|^2+\|w\|^2}}\nonumber\\
    \geq\,& \frac{\|x-a\|\|w\|^2}{4(\max\{\|z-x\|,\|z-a\|\}+\|w\|)^2}.
\end{align*}
If $\langle V(x), a-z\rangle\leq 0$, by \eqref{eq:oldL3.1}, we have
\begin{align*}
    E(x,z+w)\geq\,&\|z+w-x\|+\|x-a\|-\|z+w-a\|\nonumber\\
    =\,&\sqrt{\|z-x\|^2+\|w\|^2}+\|x-a\|-\sqrt{\|z-a\|^2+\|w\|^2}\geq \|x-a\|. \qedhere
\end{align*}
\end{proof}

The next lemma shows that if $x$ and $z$ lie in the interior of a common transport ray, then the matrix $\mathbf{I}_d+\langle x-z, V(x)\rangle F(z)$ is positive definite. 

\begin{lemma}\label{L2.7}
For $\mathcal{L}^d$-a.e.\ $z_0\in\mathcal{T}_1^*$, the following holds for $T=T(z_0)$: for any $x,z\in\mathrm{int}(T)$, $\mathbf{I}_d+\langle x-z, V(x)\rangle F(z)$ is a positive definite matrix. 
\end{lemma}

\begin{proof}
We fix any $T\in\mathcal{S}$ for which the conclusions of \cref{L2.4,L2.6} hold. 

We first consider $x,z\in \mathrm{int}(T)$ such that $\langle V(x),x-z\rangle>0$. From the conclusion of \cref{L2.6} (with $x$ replaced by $z$), for any $w_0\in  O_T$ such that $\|w_0\|=1$ and any $\delta>0$, there exist $(t_n)_{n=1}^{\infty}\subseteq (0,1)$ and $(w_n)_{n=1}^{\infty}\subseteq O_T$ such that as $n\rightarrow\infty$,
\begin{equation*}
     t_n\rightarrow 0, \quad t_n\slash 2\leq \|w_n\|\leq t_n, \quad \frac{w_n}{\|w_n\|}\rightarrow w_0, \quad \Big|u(z+w_n)-u(z)-\frac{1}{2}w_n^{\top}F(z) w_n\Big|\leq \delta t_n^2.
\end{equation*}
By \cref{L3.2}, we have  
\begin{equation*}
    E(x,z+w_n)\geq \frac{\|x-a\|\|w_n\|^2}{4(\|z-a\|+\|w_n\|)^2}.
\end{equation*}
Since $E(x,z+w_n)=\|z+w_n-x\|-u(x)+u(z+w_n)$ and $u(x)-u(z)=\|x-z\|$ (as $\langle V(x),x-z\rangle>0$), we have 
\begin{align*}
     E(x,z+w_n)
     =\,& \sqrt{\|x-z\|^2+\|w_n\|^2}-\|x-z\|+u(z+w_n)-u(z)\nonumber\\
     =\,& \frac{\|w_n\|^2}{\sqrt{\|x-z\|^2+\|w_n\|^2}+\|x-z\|}+u(z+w_n)-u(z).
\end{align*}
Hence 
\begin{align*}
    & \frac{\|w_n\|^2}{\sqrt{\|x-z\|^2+\|w_n\|^2}+\|x-z\|}+\frac{1}{2}w_n^{\top}F(z)w_n + 4 \delta\|w_n\|^2\nonumber\\
    \geq\,& \frac{\|w_n\|^2}{\sqrt{\|x-z\|^2+\|w_n\|^2}+\|x-z\|}+\frac{1}{2}w_n^{\top}F(z)w_n+\delta t_n^2\nonumber\\
    \geq\,&\frac{\|w_n\|^2}{\sqrt{\|x-z\|^2+\|w_n\|^2}+\|x-z\|}+u(z+w_n)-u(z)
    =E(x,z+w_n)\geq \frac{\|x-a\|\|w_n\|^2}{4(\|z-a\|+\|w_n\|)^2}.
\end{align*}
Dividing both sides by $\|w_n\|^2$ and taking $n\rightarrow\infty$ (note \eqref{L2.4.enew3}), we obtain
\begin{equation}
     \frac{1}{2\|x-z\|}+\frac{1}{2}w_0^{\top}F(z) w_0+4\delta\geq \frac{\|x-a\|}{4\|z-a\|^2}. 
\end{equation}
As $\delta>0$ is arbitrary and $\|w_0\|=1$, we conclude that
\begin{equation}
     w_0^{\top}\bigg(\frac{1}{\|x-z\|}\mathbf{I}_d+F(z)\bigg)w_0\geq \frac{\|x-a\|}{2\|z-a\|^2}. 
\end{equation}
Now for any $y\in\mathbb{R}^d$, we can write $y=s_1 V(x)+s_2 w_0$, where $s_1,s_2\in\mathbb{R}$, $w_0\in O_T$, and $\|w_0\|=1$. Then by \eqref{Eq2.6.5}, we have \begin{align}\label{Lem3.4.eq14}      y^{\top}\big(\mathbf{I}_d+\|x-z\|F(z)\big)y=\,&s_1^2+s_2^2 w_0^{\top}\big(\mathbf{I}_d+\|x-z\|F(z)\big)w_0\nonumber\\
     \geq\,&s_1^2+\frac{s_2^2\|x-z\|\|x-a\|}{2\|z-a\|^2}\geq  \frac{\|x-z\|\|x-a\|}{2\|z-a\|^2}\|y\|^2,
\end{align}
where we note that $\|y\|^2=s_1^2+s_2^2$ and $\max\{\|x-z\|,\|x-a\|\}\leq \|z-a\|$. Hence for any $x,z\in\mathrm{int}(T)$ such that $\langle V(x),x-z\rangle>0$, $\mathbf{I}_d+\langle x-z, V(x)\rangle F(z)=\mathbf{I}_d+\|x-z\| F(z)$ is positive definite.

Now if $x=z\in\mathrm{int}(T)$, then $\mathbf{I}_d+\langle x-z, V(x)\rangle F(z)=\mathbf{I}_d$ is positive definite. If $x,z\in\mathrm{int}(T)$ and $\langle V(x), x-z\rangle <0$, then by \Cref{Rmk.L2.4} (with $x'=z$),
\begin{equation*}
     \big(\mathbf{I}_d+\langle x-z, V(x)\rangle F(z)\big)\big(\mathbf{I}_d+\langle z-x, V(z)\rangle F(x)\big)=\mathbf{I}_d.
\end{equation*}
As $\langle V(z),z-x\rangle=\langle V(x),z-x\rangle>0$, $\mathbf{I}_d+\langle z-x, V(z)\rangle F(x)$ is positive definite, hence $\mathbf{I}_d+\langle x-z, V(x)\rangle F(z)$ is also positive definite. 
\end{proof}

We now complete the proof of \Cref{P-lb-1}.

\begin{proof}[Proof of \Cref{P-lb-1}(a)]
\noindent\textbf{Case 1.} We first consider the case where $x\in\mathcal{T}_1^{*}$, $z\in \mathrm{int}(T(x))$, and $\langle V(x), x-z\rangle >0$. Note that in this case, $u(x)-u(z)=\|x-z\|$, and for any $w\in O_{T(x)}$,
\begin{align*}
    E(x,z+w) =\,&\|x-z-w\|-u(x)+u(z+w)\nonumber\\
    =\,& \sqrt{\|x-z\|^2+\|w\|^2}-\|x-z\|+u(z+w)-u(z)\nonumber\\
    =\,& \frac{\|w\|^2}{\sqrt{\|x-z\|^2+\|w\|^2}+\|x-z\|}+u(z+w)-u(z).
\end{align*}
Hence
\begin{equation*}
    E(x,z+w)\begin{cases}
        \geq \frac{\|w\|^2}{2\|x-z\|+\|w\|}+u(z+w)-u(z)\\
        \leq \frac{\|w\|^2}{2\|x-z\|}+u(z+w)-u(z).
    \end{cases}
\end{equation*}
Replacing $w$ by $\sqrt{\eps}w$, we obtain that
\begin{equation}\label{Eq2.26}
    E(x,z+\sqrt{\eps}w)\begin{cases}
        \geq \frac{\eps\|w\|^2}{2\|x-z\|+\sqrt{\eps}\|w\|}+u(z+\sqrt{\eps}w)-u(z)\\
        \leq \frac{\eps\|w\|^2}{2\|x-z\|}+u(z+\sqrt{\eps}w)-u(z).
    \end{cases}
\end{equation}

Below, we assume that the conclusions of \cref{L2.4,L2.6} hold for $T=T(x)$ (note that this is the case for $\mathcal{L}^d$-a.e.\ $x\in\mathcal{T}_1^{*}$).

We fix any $M\geq 1$ and $\delta,\delta'\in (0,1)$. By \cref{L2.6} (with $t=M\sqrt{\eps}$), for $\eps>0$ sufficiently small, we have
\begin{equation*}
   \mathcal{H}^{d-1}\Big(\Big\{w\in O_{T(x)}:\|w\|\leq M\sqrt{\eps},\Big|u(z+w)-u(z)-\frac{1}{2}w^{\top} F(z) w\Big|>\delta \eps M^2\Big\}\Big) \leq \delta' (M\sqrt{\eps})^{d-1}.
\end{equation*}
Thus, denoting 
\begin{equation*}
    \mathscr{D}:=\Big\{w\in O_{T(x)}:\|w\|\leq M,\Big|u(z+\sqrt{\eps}w)-u(z)-\frac{\eps}{2}w^{\top} F(z) w\Big|\leq \delta \eps M^2\Big\},
\end{equation*}
we have
\begin{equation}\label{Lem3.8.eq1}
    \mathcal{H}^{d-1}\big(\big\{w\in O_{T(x)}:\|w\|\leq M\big\}\backslash\mathscr{D}\big)\leq \delta' M^{d-1}.
\end{equation}
Hence by \eqref{Eq2.26} and the definition of $\mathscr{D}$, for $\eps>0$ sufficiently small,
\begin{align}\label{Eq2.31}
 &\int_{\{w\in O(V(x)):\|w\|\leq M\}} 
 e^{-E(x,z+\sqrt{\eps}w)\slash \eps}d\mathcal{H}^{d-1}(w)\leq \int_{\mathscr{D}} e^{-E(x,z+\sqrt{\eps}w)\slash \eps}d\mathcal{H}^{d-1}(w)+\delta' M^{d-1}\nonumber\\
 \leq\,& e^{\delta M^2}\int_{\mathscr{D}} e^{-\frac{\|w\|^2}{2\|x-z\|+\sqrt{\eps}\|w\|}-\frac{1}{2}w^{\top}F(z)w}d\mathcal{H}^{d-1}(w)+\delta' M^{d-1}\nonumber\\
 \leq\,& e^{\delta M^2}\int_{\{w\in O(V(x)):\|w\|\leq M\}} e^{-\frac{\|w\|^2}{2\|x-z\|+\sqrt{\eps}\|w\|}-\frac{1}{2}w^{\top}F(z)w}d\mathcal{H}^{d-1}(w)+\delta' M^{d-1},
\end{align}
\begin{align}\label{Eq2.32}
&\int_{O(V(x))} 
 e^{-E(x,z+\sqrt{\eps}w)\slash \eps}\mathbbm{1}_{\|w\|\leq 2D\eps^{-1\slash 2}}d\mathcal{H}^{d-1}(w)\geq\int_{\{w\in O(V(x)):\|w\|\leq M\}} 
 e^{-E(x,z+\sqrt{\eps}w)\slash \eps}d\mathcal{H}^{d-1}(w)\nonumber\\
 \geq\,& \int_{\mathscr{D}} e^{-E(x,z+\sqrt{\eps}w)\slash \eps}d\mathcal{H}^{d-1}(w) \geq  e^{-\delta M^2}\int_{\mathscr{D}}e^{-\frac{\|w\|^2}{2\|x-z\|}-\frac{1}{2}w^{\top}F(z)w}d\mathcal{H}^{d-1}(w).
\end{align}

By the dominated convergence theorem (note \eqref{L2.4.enew3}), sending $\eps\rightarrow 0^{+}$ and then $\delta,\delta'\rightarrow 0^{+}$ in \eqref{Eq2.31} yields
\begin{align*}
    &\limsup_{\eps\rightarrow 0^{+}} \int_{\{w\in O(V(x)):\|w\|\leq M\}} 
 e^{-E(x,z+\sqrt{\eps}w)\slash \eps}d\mathcal{H}^{d-1}(w)\nonumber\\
 \leq\,&\int_{\{w\in O(V(x)):\|w\|\leq M\}} e^{-\frac{1}{2}w^{\top}\big(\frac{1}{\|x-z\|}\mathbf{I}_d+F(z)\big)w}d\mathcal{H}^{d-1}(w).
\end{align*}
By \eqref{bddsalphabeta}, we have $\|x-z\|\leq \alpha(z)\leq 2D$. Hence by \Cref{L3.2}, for any $w\in O_{T(x)}$ such that $\|w\|\leq 2D\eps^{-1\slash 2}$, we have $E(x,z+\sqrt{\eps}w)\geq c\eps \alpha(x)\|w\|^2$, and so
\begin{align*}
    &\int_{\{w\in O(V(x)):\|w\|> M\}} 
 e^{-E(x,z+\sqrt{\eps}w)\slash \eps} \mathbbm{1}_{\|w\|\leq 2D\eps^{-1\slash 2}} d\mathcal{H}^{d-1}(w)\nonumber\\
 \leq\,&\int_{\{w\in O(V(x)):\|w\|> M\}}\exp(-c\alpha(x)\|w\|^2)d\mathcal{H}^{d-1}(w). 
\end{align*}
Consequently, we have
\begin{align*}
   & \limsup_{\eps\rightarrow 0^{+}} \int_{O(V(x))} 
 e^{-E(x,z+\sqrt{\eps}w)\slash \eps} \mathbbm{1}_{\|w\|\leq 2D\eps^{-1\slash 2}} d\mathcal{H}^{d-1}(w)\nonumber\\
 \leq\,& \int_{O(V(x))} e^{-\frac{1}{2}w^{\top}\big(\frac{1}{\|x-z\|}\mathbf{I}_d+F(z)\big)w}d\mathcal{H}^{d-1}(w)+\int_{\{w\in O(V(x)):\|w\|> M\}}\exp(-c\alpha(x)\|w\|^2)d\mathcal{H}^{d-1}(w).
\end{align*}
By the dominated convergence theorem, taking $M\rightarrow\infty$ in the above display, we obtain
\begin{align}\label{Eq2.36}
    &\limsup_{\eps\rightarrow 0^{+}} \int_{O(V(x))} 
 e^{-E(x,z+\sqrt{\eps}w)\slash \eps} \mathbbm{1}_{\|w\|\leq 2D\eps^{-1\slash 2}} d\mathcal{H}^{d-1}(w)\nonumber\\
 \leq\,& \int_{O(V(x))} e^{-\frac{1}{2}w^{\top}\big(\frac{1}{\|x-z\|}\mathbf{I}_d+F(z)\big)w}d\mathcal{H}^{d-1}(w).
\end{align}

By \eqref{L2.4.enew3} and \eqref{Lem3.8.eq1},
\begin{align*}
    \int_{\{w\in O_{T(x)}:\|w\|\leq M\}\backslash\mathscr{D}}e^{-\frac{\|w\|^2}{2\|x-z\|}-\frac{1}{2}w^{\top}F(z)w}d\mathcal{H}^{d-1}(w)\leq\,  & \int_{\{w\in O_{T(x)}   :     \|w\|\leq M\}\backslash\mathscr{D}}e^{\frac{1}{2}\|F(z)\|_{2}\|w\|^2}d\mathcal{H}^{d-1}(w)\nonumber\\
    \leq \, &  e^{2M^2/\min\{\alpha(z),\beta(z)\}}\cdot \delta' M^{d-1}.
\end{align*}
Sending $\eps\rightarrow 0^{+}$ in \eqref{Eq2.32} and noting the above display, we obtain that
\begin{align*}
&\liminf_{\eps\rightarrow 0^{+}}\int_{O(V(x))} 
 e^{-E(x,z+\sqrt{\eps}w)\slash \eps} \mathbbm{1}_{\|w\|\leq 2D\eps^{-1\slash 2}} d\mathcal{H}^{d-1}(w)\nonumber\\
 \geq\,& e^{-\delta M^2} \bigg(\int_{\{w\in O(V(x)):\|w\|\leq M\}}e^{-\frac{\|w\|^2}{2\|x-z\|}-\frac{1}{2}w^{\top}F(z)w}d\mathcal{H}^{d-1}(w)-e^{2M^2/\min\{\alpha(z),\beta(z)\}}\cdot \delta' M^{d-1}\bigg).
\end{align*}
Sending $\delta,\delta'\rightarrow 0^{+}$, we obtain
\begin{align*}
    &\liminf_{\eps\rightarrow 0^{+}}\int_{O(V(x))} 
 e^{-E(x,z+\sqrt{\eps}w)\slash \eps} \mathbbm{1}_{\|w\|\leq 2D\eps^{-1\slash 2}} d\mathcal{H}^{d-1}(w)\nonumber\\
 \geq\,& \int_{\{w\in O(V(x)):\|w\|\leq M\}}e^{-\frac{\|w\|^2}{2\|x-z\|}-\frac{1}{2}w^{\top}F(z)w}d\mathcal{H}^{d-1}(w).
\end{align*}
Finally, taking $M\rightarrow\infty$ and using the monotone convergence theorem, we get
\begin{align}\label{Eq2.37}
    &\liminf_{\eps\rightarrow 0^{+}}\int_{O(V(x))} 
 e^{-E(x,z+\sqrt{\eps}w)\slash \eps} \mathbbm{1}_{\|w\|\leq 2D\eps^{-1\slash 2}} d\mathcal{H}^{d-1}(w)\nonumber\\
 \geq\,&\int_{O(V(x))} e^{-\frac{1}{2}w^{\top}\big(\frac{1}{\|x-z\|}\mathbf{I}_d+F(z)\big)w}d\mathcal{H}^{d-1}(w).
\end{align}

By \eqref{Eq2.36} and \eqref{Eq2.37}, and noting \Cref{Def2} and \Cref{L2.7}, for $\mathcal{L}^d$-a.e.\ $x\in\mathcal{T}_1^{*}$ and any $z\in \mathrm{int}(T(x))$ such that $\langle V(x),x-z\rangle>0$, we have $\lim_{\eps\rightarrow 0^{+}} G_{\eps}(x,z) = G(x,z)$.

\bigskip

\noindent\textbf{Case 2.} We now consider the case where $x\in\mathcal{T}_1^{*}$, $z\in \mathrm{int}(T(x))$, and $\langle V(x), x-z  \rangle< 0$. In this case, we have $u(z)-u(x)=\|x-z\|$, and for any $w\in O(V(x))$,
\begin{align}
    E(x,z+w) =\,&\|x-z-w\|-u(x)+u(z+w)\nonumber\\
    =\,& \sqrt{\|x-z\|^2+\|w\|^2}+\|x-z\|+u(z+w)-u(z)\nonumber\\
    \geq\,& \sqrt{\|x-z\|^2+\|w\|^2}+\|x-z\|-\|w\|\geq \|x-z\|,
\end{align}
where we use the fact that $\|u\|_{\mathrm{Lip}}\leq 1$. Hence
\begin{equation*}
    G_{\eps}(x,z)=\log\bigg(\int_{O(V(x))} 
 e^{-E(x,z+\sqrt{\eps}w)\slash \eps}\mathbbm{1}_{\|w\|\leq 2D\eps^{-1\slash 2}} d\mathcal{H}^{d-1}(w)\bigg)
 \leq\log\big(C \eps^{-(d-1)\slash 2} e^{-\|x-z\|\slash \eps}\big).
\end{equation*}
Therefore, $\lim_{\eps\rightarrow 0^{+}} G_{\eps}(x,z) = -\infty = G(x,z)$. 
\end{proof}

\begin{proof}[Proof of \Cref{P-lb-1}(b)]

\noindent\textbf{Case 1.} We first consider the case where $\langle z-a(x), V(x)\rangle >0$. By \cref{L3.2}, 
\begin{align*}
    &\int_{O(V(x))} e^{-E(x,z+\sqrt{\eps}w)\slash \eps} \mathbbm{1}_{\|w\|\leq 2D\eps^{-1\slash 2}} d\mathcal{H}^{d-1}(w)\nonumber\\
    \leq\,& \int_{O(V(x))} e^{-\alpha(x)\slash\eps}\mathbbm{1}_{\|w\|\leq 2D\eps^{-1\slash 2}}d\mathcal{H}^{d-1}(w)\leq C e^{-\alpha(x)\slash \eps}\eps^{-(d-1)\slash 2}. 
\end{align*}
As $\alpha(x)>0$, we have $\lim_{\eps\rightarrow 0^{+}}G_{\eps}(x,z)=-\infty =G(x,z)$.

\noindent\textbf{Case 2.} Now consider the case where $\langle z-b(x), V(x)\rangle <0$ and $z\in \mathcal{Y}$. For any $w\in O(V(x))$,
\begin{align}\label{Eq3.21}
    & E(x,z+w)=\|z+w-x\|-u(x)+u(z+w)\nonumber\\
    =\,& \sqrt{\|z-x\|^2+\|w\|^2}+(u(z)-u(x))+(u(z+w)-u(z))\nonumber\\
    \geq\,& \sqrt{\|z-x\|^2+\|w\|^2}-\|w\|+u(z)-u(x)  \nonumber\\
    \geq\,&   \|z-x\|+u(z)-u(x)-\|w\|=E(x,z)-\|w\|,
\end{align}
where we use the fact that $\|u\|_{\mathrm{Lip}}\leq 1$. By \cref{TransportRays} and the fact that $z\in\mathcal{Y}$,
\begin{equation}\label{Eq3.23}
   E(x,z)= \|z-x\|+u(z)-u(x)>0.
\end{equation}
Moreover, by \cref{L3.2}, for any $w\in O(V(x))$ such that $\|w\|\leq 2D$,  
\begin{equation}\label{Eq3.24}
    E(x,z+w)\geq c' \|w\|^2,
\end{equation}
where $c'>0$ only depends on $x,z,D$. Now for any $M\geq 1$, by \eqref{Eq3.21}--\eqref{Eq3.24}, we have
\begin{align*}
    &\int_{O(V(x))} e^{-E(x,z+\sqrt{\eps}w)\slash \eps} \mathbbm{1}_{\|w\|\leq 2D\eps^{-1\slash 2}} d\mathcal{H}^{d-1}(w)\nonumber\\
    \leq\,& \int_{O(V(x))} e^{-E(x,z+\sqrt{\eps}w)\slash \eps} \mathbbm{1}_{\|w\|\leq M} d\mathcal{H}^{d-1}(w) +   \int_{O(V(x))} e^{-E(x,z+\sqrt{\eps}w)\slash \eps} \mathbbm{1}_{M\leq\|w\|\leq 2D\eps^{-1\slash 2}} d\mathcal{H}^{d-1}(w)\nonumber\\
    \leq\,& \int_{O(V(x))} e^{-(E(x,z)-M\sqrt{\eps})\slash \eps} \mathbbm{1}_{\|w\|\leq M} d\mathcal{H}^{d-1}(w)+\int_{O(V(x))} e^{-c'\|w\|^2} \mathbbm{1}_{M\leq\|w\|\leq 2D\eps^{-1\slash 2}} d\mathcal{H}^{d-1}(w)\nonumber\\
    \leq\,& CM^{d-1}e^{-(E(x,z)-M\sqrt{\eps})\slash \eps}+\int_{O(V(x))} e^{-c'\|w\|^2} \mathbbm{1}_{\|w\|\geq M} d\mathcal{H}^{d-1}(w).
\end{align*}
Taking $\eps\rightarrow 0^{+}$ and noting \eqref{Eq3.23}, we obtain
\begin{equation*}
   \limsup_{\eps\rightarrow 0^{+}} \int_{O(V(x))} e^{-E(x,z+\sqrt{\eps}w)\slash \eps} \mathbbm{1}_{\|w\|\leq 2D\eps^{-1\slash 2}} d\mathcal{H}^{d-1}(w)\leq \int_{O(V(x))} e^{-c'\|w\|^2} \mathbbm{1}_{\|w\|\geq M} d\mathcal{H}^{d-1}(w).
\end{equation*}
Taking $M\rightarrow\infty$, by the dominated convergence theorem, we obtain
\begin{eqnarray*}
    \lim_{\eps\rightarrow 0^{+}} \int_{O(V(x))} e^{-E(x,z+\sqrt{\eps}w)\slash \eps} \mathbbm{1}_{\|w\|\leq 2D\eps^{-1\slash 2}} d\mathcal{H}^{d-1}(w) = 0.
\end{eqnarray*}
Hence $\lim_{\eps\rightarrow 0^{+}}G_{\eps}(x,z)=-\infty=G(x,z)$.
\end{proof}

\subsubsection{Proof of \Cref{P-lb-2}}

\begin{proof}[Proof of \Cref{P-lb-2}(a)]
Fix any $x\in\mathcal{T}_1^*$ and $z\in L(x)$ with $\max\{\|x-z\|,\|z-a(x)\|\}\leq 2D$. If $\langle V(x),a(x)-z\rangle\geq 0$, for any $w\in O_{T(x)}$ such that $\|w\|\leq 2D\eps^{-1\slash 2}$, using \cref{L3.2} and the fact that $\max\{\|x-z\|,  \|z-a(x)\|\}\leq 2D$, we obtain $E(x,z+\sqrt{\eps}w)\geq c\eps\alpha(x)\|w\|^2$. Hence for any $\eps>0$,
\begin{align*}
    G_{\eps}(x,z)=\,&\log\bigg(\int_{O(V(x))}e^{-E(x,z+\sqrt{\eps}w)\slash \eps}\mathbbm{1}_{\|w\|\leq 2D\eps^{-1\slash 2}} 
 d\mathcal{H}^{d-1}(w)\bigg)\nonumber\\
 \leq\,& \log\bigg(\int_{O(V(x))}\exp\big(-c\alpha(x)\|w\|^2\big)d\mathcal{H}^{d-1}(w)\bigg) \leq -\frac{d-1}{2}\log(\alpha(x))+C.
\end{align*}

Now if $\langle V(x), a(x) - z\rangle<0$, by \cref{L3.2}, we have 
\begin{align*}
     G_{\eps}(x,z)=\,&\log\bigg(\int_{O(V(x))}e^{-E(x,z+\sqrt{\eps}w)\slash \eps}\mathbbm{1}_{\|w\|\leq 2D\eps^{-1\slash 2}} 
 d\mathcal{H}^{d-1}(w)\bigg)\nonumber\\
 \leq\,& \log\bigg(\int_{O(V(x))}e^{-\alpha(x) \slash \eps}\mathbbm{1}_{\|w\|\leq 2D\eps^{-1\slash 2}} 
 d\mathcal{H}^{d-1}(w)\bigg)\leq \log\big(C e^{-\alpha(x)\slash\eps}\eps^{-(d-1)\slash 2}\big)\nonumber\\
 \leq\,& -\frac{\alpha(x)}{\eps}-\frac{d-1}{2}\log{\eps}+C \leq -\frac{d-1}{2}\log(\alpha(x)) + C.
\end{align*}

For any $x\in\mathcal{T}_1^*$ and $z\in L(x)$, if $z\notin\mathrm{int}(T(x))$, or if $z\in\mathrm{int}(T(x))$ and
$\langle V(x),x-z\rangle\leq 0$, we have $G(x,z)=-\infty$. For $\mathcal{L}^d$-a.e.\ $x\in\mathcal{T}_1^*$ and any $z\in \mathrm{int}(T(x))\backslash\{x\}$, by \cref{P-lb-1}(a), we have $\lim_{\eps\rightarrow 0^{+}}G_{\eps}(x,z)=G(x,z)$. Hence using the established upper bound on $G_{\eps}(x,z)$, we obtain that $G(x,z)\leq -\frac{d-1}{2}\log(\alpha(x))+C$. 
\end{proof}

\begin{proof}[Proof of \Cref{P-lb-2}(b)]

By \cref{L2.4,L2.7}, for $\mathcal{L}^d$-a.e.\ $x\in \mathcal{T}_1^{*}$ and any $z\in\mathrm{int}(T(x))$ such that $\langle V(x), x-z\rangle>0$, we have $F(z)^{\top}=F(z)$, $F(z)V(z)=F(z)^{\top}V(z)=0$, and $\mathbf{I}_d+\langle x-z, V(x)\rangle F(z)=\mathbf{I}_d+\|x-z\| F(z)$ is positive definite. Hence by \eqref{Eq3.33}, we have 
\begin{align*}
 G(x,z)=\,&\log\bigg(\int_{O(V(x))}e^{-\frac{1}{2}w^{\top}\big(\frac{1}{\|x-z\|}\mathbf{I}_d+F(z)\big)w}d\mathcal{H}^{d-1}(w)\bigg)\nonumber\\
  =\,& \frac{d-1}{2}\log(\|x-z\|)+\log\bigg(\int_{O(V(x))}e^{-\frac{1}{2}w^{\top}\big(\mathbf{I}_d+\|x-z\|F(z)\big)w}d\mathcal{H}^{d-1}(w)\bigg)\nonumber\\
  =\,& \frac{d-1}{2}\log(2\pi\|x-z\|)-\frac{1}{2}\log\det\big(\mathbf{I}_d+\|x-z\|F(z)\big). \qedhere
\end{align*} 
\end{proof}

\begin{proof}[Proof of \Cref{P-lb-2}(c)]
Fix any $x\in\mathcal{T}_1^{*}$ such that the conclusions of \Cref{L2.4,L2.7} hold for $T=T(x)$ and the conclusion of \Cref{P-lb-2}(b) holds for $x$ (note that this is the case for $\mathcal{L}^d$-a.e.\ $x\in\mathcal{T}_1^{*}$). 

By \cref{L2.4} (see \eqref{L2.4.enew4}) and \Cref{P-lb-2}(b), $F(\cdot)$ is continuous on $\mathrm{int}(T(x))$, and
\begin{equation*}
    G(x,x+sV(x))=\frac{d-1}{2}\log(2\pi|s|)-\frac{1}{2}\log\det\big(\mathbf{I}_d-s F(x+s V(x))\big)
\end{equation*}
for any $s\in (-\beta(x),0)$. Hence for any $t\in (-\beta(x),0)$,
\begin{equation*}
    \lim_{s\rightarrow t}\big(\mathbf{I}_d-s F(x+s V(x))\big)=\mathbf{I}_d-t F(x+t V(x)).
\end{equation*}
By \cref{L2.7} (with $z=x+tV(x)$), $\det\big(\mathbf{I}_d-t F(x+t V(x))\big)>0$. Hence
\begin{equation*}
     \lim_{s\rightarrow t}\log\det\big(\mathbf{I}_d-s F(x+s V(x))\big) = \log\det\big(\mathbf{I}_d-t F(x+t V(x))\big).
\end{equation*}
We conclude that $\lim_{s\rightarrow t}G(x,x+sV(x))=G(x,x+tV(x))$.

By \Cref{P-lb-2}(b), for any $t\in(-\beta(x),0)$, 
\begin{equation*}
    G(x,x+tV(x))=\frac{d-1}{2}\log(2\pi|t|)-\frac{1}{2}\log\det\big(\mathbf{I}_d-tF(x+tV(x))\big).
\end{equation*}
Hence $\lim\limits_{t\rightarrow 0^{-}} G(x,x+tV(x))=-\infty$. By \cref{Rmk.L2.4} (with $x'=x+tV(x)$), for any $t\in (-\beta(x),0)$, we have $\big(\mathbf{I}_d-tF(x+tV(x))\big)\big(\mathbf{I}_d+tF(x)\big)=\mathbf{I}_d$, hence
\begin{equation}\label{L.3.10.e1}
    G(x,x+tV(x))=\frac{d-1}{2}\log(2\pi|t|)+\frac{1}{2}\log\det\big(\mathbf{I}_d+tF(x)\big).
\end{equation}
By \cref{L2.7} (with $z=x+tV(x)$), for any $t\in (-\beta(x),0)$, $\mathbf{I}_d-tF(x+tV(x))$ is positive definite, hence $\mathbf{I}_d+t F(x)$ is positive definite. Thus $\mathbf{I}_d-\beta(x) F(x)=\lim\limits_{t\rightarrow (-\beta(x))^{+}}\big(\mathbf{I}_d+t F(x)\big)$ is positive semidefinite. Hence by \eqref{L.3.10.e1}, if $\det\big(\mathbf{I}_d-\beta(x) F(x)\big)=0$, then $\lim\limits_{t\rightarrow (-\beta(x))^{+}}G(x,x+tV(x))=-\infty$; if $\det\big(\mathbf{I}_d-\beta(x) F(x)\big)>0$, then 
\begin{equation*}
    \lim\limits_{t\rightarrow (-\beta(x))^{+}}G(x,x+tV(x))=\frac{d-1}{2}\log(2\pi\beta(x))+\frac{1}{2}\log\det\big(\mathbf{I}_d-\beta(x)F(x)\big)\in (-\infty,\infty).  \qedhere
\end{equation*}
\end{proof}

\subsection{Conditional marginal densities on transport rays}\label{Sec3:2}

This subsection derives explicit formulas for $\lambda,\tilde{\mu}_T,\tilde{\nu}_T,\tilde{f},\tilde{g}$ (cf.\ Section~\ref{Sect.1.1}) and establishes their key properties. 

The following lemma derives an explicit formula for the Jacobian $R_{k,k'}(q,t)$ defined in \eqref{defRk}. Recall the sets $\bar{H}_{k,k'}$, $\tilde{H}_{k,k'}$, and $\tilde{\mathcal{T}}_{1;k,k'}^*$ from \eqref{deftildehat} and \eqref{deftildeTkk}. 

\begin{lemma}\label{Lem3.10n}
For any $k\in[K],k'\in[K']$, the following holds for $\mathcal{L}^d$-a.e.\ $z_0\in\tilde{\mathcal{T}}_{1;k,k'}^*$. Let $q$ be the intersection point of $T(z_0)$ and $H_{k,k'}$ (note that $q\in\tilde{H}_{k,k'}$). Then for every $t\in (-\beta(q),\alpha(q))$, 
\begin{equation}\label{L3.10E1}
    R_{k,k'}(q,t)= |\langle V(q), \mathtt{a}_{k,k'}\rangle| \det(\mathbf{I}_d+tF(q)),
\end{equation}
where $\mathtt{a}_{k,k'}$ is as in \eqref{Hkkdef}.  
\end{lemma}
\begin{proof}

Recall the discussion below \eqref{subs}. Note that by \eqref{Psikkp}, for any $x=(x_1,\cdots,x_{d-1})\in\mathbb{R}^{d-1}$ and $t\in\mathbb{R}$ such that $|t|\leq 2D$, we have
\begin{equation*}
    \EuScript{P}(x,t)=\Psi_{k,k'}\big(\Omega_{k,k'}^{-1}(x),t\big)=\Omega_{k,k'}^{-1}(x)+t(\tilde{V}_{k,k'}\circ\Omega_{k,k'}^{-1})(x).
\end{equation*}
Note that
\begin{equation}\label{pxt0.1}
    \frac{\partial \EuScript{P}(x,t)}{\partial t}=(\tilde{V}_{k,k'}\circ\Omega_{k,k'}^{-1})(x),
\end{equation}
and if $\tilde{V}_{k,k'}\circ \Omega_{k,k'}^{-1}$ is differentiable at $x$, then for every $\ell\in [d-1]$, we have (by \eqref{inverse_map}) 
\begin{equation}\label{pxt1}
    \frac{\partial \EuScript{P}(x,t)}{\partial x_{\ell}}=\mathtt{e}_{k,k'}^{(\ell)}+t\frac{\partial (\tilde{V}_{k,k'}\circ\Omega_{k,k'}^{-1})(x)}{\partial x_{\ell}}.
\end{equation}

By \cite[Theorem 1.35]{MR3409135} and Rademacher's theorem (note that $\tilde{V}_{k,k'}\circ \Omega_{k,k'}^{-1}$ is Lipschitz continuous), for $\mathcal{H}^{d-1}$-a.e.\ $q\in\tilde{H}_{k,k'}$, we have
\begin{equation}\label{niceco}
    \lim_{s\rightarrow 0^+}\frac{\mathcal{H}^{d-1}(\{q'\in \tilde{H}_{k,k'}:\|q'-q\|\leq s\})}{\mathcal{H}^{d-1}(\{q'\in H_{k,k'}:\|q'-q\|\leq s\})}=1, \quad \tilde{V}_{k,k'}\circ \Omega_{k,k'}^{-1}\text{ is differentiable at }\Omega_{k,k'}(q).
\end{equation}
By \Cref{L2.3n,L2.4,L2.7}, for $\mathcal{L}^d$-a.e.\ $z_0\in\tilde{\mathcal{T}}_{1;k,k'}^*$, if $q$ denotes the intersection point of $T(z_0)$ and $H_{k,k'}$ (note that $q\in\tilde{H}_{k,k'}$, and so $\min\{\alpha(q),\beta(q)\}\geq d_0$), then \eqref{niceco} holds, and
\begin{align}\label{niceco.1}
  &  \|F(q)\|_2\leq \frac{4}{d_0}, \quad F(q)V(q)=0, \quad
 \lim_{s\rightarrow 0^{+}}\frac{\mathcal{H}^{d-1}(\mathcal{W}(q;s,\delta))}{s^{d-1}}=0\quad\text{for all } \delta>0,\nonumber\\
 & \hspace{1in} \det(\mathbf{I}_d+tF(q))>0 \quad\text{for all }t\in (-\beta(q),\alpha(q)).
\end{align}
Below, we consider any $q\in \tilde{H}_{k,k'}$ such that \eqref{niceco} and \eqref{niceco.1} hold. Note that by the definition of $\tilde{H}_{k,k'}$ in \eqref{deftildehat}, $\langle V(q),\mathtt{a}_{k,k'}\rangle\neq 0$; without loss of generality, we assume that $\langle V(q),\mathtt{a}_{k,k'}\rangle>0$. We denote $x_0:=\Omega_{k,k'}(q)$. 

In the following, we consider any $\delta\in (0,\langle V(q),\mathtt{a}_{k,k'}\rangle d_0\slash 100)$. By \eqref{niceco} and \eqref{niceco.1}, for any $\delta'\in (0,1)$, there exists $s_0\in (0,\langle V(q),\mathtt{a}_{k,k'}\rangle d_0^2\delta\slash 100)$ such that the following holds:
\begin{itemize}
    \item[(a)] for all $s\in (0,s_0)$,
\begin{equation}\label{bddWqs}
     \mathcal{H}^{d-1}(\{q'\in H_{k,k'} \backslash \tilde{H}_{k,k'}:\|q'-q\|\leq s\})\leq \delta' s^{d-1},\quad
    \mathcal{H}^{d-1}(\mathcal{W}(q;s,\delta)) \leq \delta' s^{d-1};
\end{equation}
    \item[(b)] for any $q'\in\tilde{H}_{k,k'}$ with $\|q'-q\|\leq s_0$ (which implies $\|\Omega_{k,k'}(q')-x_0\|=\|q'-q\| \leq s_0$), we have 
\begin{align}\label{E3.1nn}
   & \big\|V(q')-V(q)-\nabla (\tilde{V}_{k,k'}\circ \Omega_{k,k'}^{-1}) (x_0)\cdot(\Omega_{k,k'}(q')-x_0)\big\|\nonumber\\
   =\,&\big\|(\tilde{V}_{k,k'}\circ\Omega_{k,k'}^{-1})(\Omega_{k,k'}(q'))-(\tilde{V}_{k,k'}\circ\Omega_{k,k'}^{-1})(x_0)-\nabla (\tilde{V}_{k,k'}\circ \Omega_{k,k'}^{-1}) (x_0)\cdot(\Omega_{k,k'}(q') - x_0)\big\|\nonumber\\
   \leq\,& \delta\|\Omega_{k,k'}(q')-x_0\|=\delta\|q'-q\|.
\end{align} 
\end{itemize}

Below, we fix any $s\in (0,s_0\slash 4)$. For any $q'\in \tilde{H}_{k,k'}$ such that $\|q'-q\|\leq s$, we define
\begin{equation}\label{defpsiq}
   \psi(q'):=q'-q-\frac{\langle q'-q, V(q)\rangle}{\langle V(q), V(q')\rangle}V(q')\in O(V(q)),
\end{equation}
where we note that by \cref{L2.0} and the fact that $\min\{\alpha(q),\beta(q),\alpha(q'),\beta(q')\}\geq d_0$ (recall \eqref{deftildehat}), 
\begin{equation}\label{lbdine}
    \|V(q')-V(q)\|\leq \frac{4 \|q'-q\|}{d_0}\leq \frac{4s}{d_0}\leq \frac{4s_0}{d_0} \leq \frac{1}{25} \,  \Rightarrow \, \langle V(q),V(q')\rangle=\frac{2-\|V(q')-V(q)\|^2}{2}\geq \frac{1}{2}. 
\end{equation}
Note that for any $q'\in \tilde{H}_{k,k'}$ with $\|q'-q\|\leq s$,  
\begin{equation*}
    \|q+\psi(q')- q'\|=\bigg|\frac{\langle q'-q, V(q)\rangle}{\langle V(q), V(q')\rangle}\bigg|\leq 2\|q'-q\|\leq 2s\leq 2s_0 \leq \frac{d_0}{50},
\end{equation*}
hence $\|q+\psi(q')- q'\|<\min\{\alpha(q'),\beta(q')\}$. Consequently,
\begin{equation}\label{insp}
    q+\psi(q')\in\mathrm{int}(T(q')),\quad V(q')=V(q+\psi(q')).
\end{equation}
Moreover,
\begin{equation}\label{lbdalpgaq}
    \min\{\alpha(q+\psi(q')),\beta(q+\psi(q'))\} \geq \min\{\alpha(q'),\beta(q')\}-\|q+\psi(q')- q'\|\geq d_0-\frac{d_0}{50}\geq \frac{d_0}{2}.
\end{equation}
By \eqref{defpsiq} and \eqref{lbdine},
\begin{align}\label{E3.48}
  &  \big\|\psi(q')-(\mathbf{I}_d-V(q) V(q)^{\top})(q'-q)\big\|=\bigg\|\frac{\langle q'-q, V(q)\rangle}{\langle V(q), V(q')\rangle}V(q')-\langle  q'-q,V(q)\rangle V(q)
\bigg\|\nonumber\\
\leq\, & \|q'-q\|\cdot \frac{\|V(q')-V(q)\|+|\langle V(q),V(q')\rangle  -1|}{\langle V(q),V(q')\rangle}\nonumber\\
\leq\,& 2\|q'-q\|(\|V(q')-V(q)\|+|\langle V(q),V(q')-V(q)\rangle|)\nonumber\\
\leq\,& 4\|q'-q\|\|V(q')-V(q)\|\leq 16sd_0^{-1}\|q'-q\|,
\end{align}
which implies
\begin{equation}\label{upsi0.1}
    \|\psi(q')\|\leq \|(\mathbf{I}_d-V(q) V(q)^{\top})(q'-q)\|+16sd_0^{-1}\|q'-q\|\leq (1+16sd_0^{-1})\|q'-q\|\leq 2s.
\end{equation}

Define $\mathscr{D}_s:=\big\{w\in O(V(q)):\|w\|\leq 2s, \min\{\alpha(q+w),\beta(q+w)\}\geq d_0\slash 2\big\}$. By Lemma \ref{L2.0}, for any $w\in\mathscr{D}_s$, we have
\begin{align*}
    \|V(q+w)-V(q)\|\leq \frac{4\|w\|}{d_0\slash 2}\leq \frac{16 s}{d_0},
\end{align*}
which implies
\begin{align}\label{lbdq}
   & |\langle V(q+w),\mathtt{a}_{k,k'}\rangle-\langle V(q),\mathtt{a}_{k,k'}\rangle|\leq \|V(q+w)-V(q)\| \leq \frac{16 s}{d_0}\nonumber\\
   \Rightarrow \,& \langle V(q+w),\mathtt{a}_{k,k'}\rangle\geq \langle V(q),\mathtt{a}_{k,k'}\rangle-\frac{16 s}{d_0}\geq \langle V(q),\mathtt{a}_{k,k'}\rangle-\frac{16\delta}{d_0}\geq \frac{1}{2}\langle V(q),\mathtt{a}_{k,k'}\rangle>0.
\end{align}
For any $w\in \mathscr{D}_s$, we define  
\begin{equation}
    \varphi(w):=q+w+\frac{\mathtt{b}_{k,k'}-\langle q+w,\mathtt{a}_{k,k'}\rangle}{\langle V(q+w),\mathtt{a}_{k,k'}\rangle}V(q+w)=q+w-\frac{\langle w,\mathtt{a}_{k,k'}\rangle}{\langle V(q+w),\mathtt{a}_{k,k'}\rangle}V(q+w),
\end{equation}
where we use $\langle q,\mathtt{a}_{k,k'}\rangle =\mathtt{b}_{k,k'}$ (recall \eqref{Hkkdef}). Note that $\varphi(w)\in H_{k,k'}$. Moreover, by \eqref{lbdq},
\begin{equation*}
    \|\varphi(w)-(q+w)\|=\bigg|\frac{\langle w,\mathtt{a}_{k,k'}\rangle}{\langle V(q+w),\mathtt{a}_{k,k'}\rangle}\bigg|\leq \frac{2\|w\|}{\langle V(q),\mathtt{a}_{k,k'}\rangle} \leq \frac{s_0}{\langle V(q),\mathtt{a}_{k,k'}\rangle}\leq\frac{\delta}{10}\leq \frac{d_0}{1000}, 
\end{equation*}
hence $\|\varphi(w)-(q+w)\|<d_0\slash   2\leq  \min\{\alpha(q+w),\beta(q+w)\}$, and consequently,
\begin{equation}\label{upsphi0.1}
     \varphi(w)\in\mathrm{int}(T(q+w)).
\end{equation}
By \eqref{lbdalpgaq} and \eqref{upsi0.1}, for any $q'\in\tilde{H}_{k,k'}$ such that $\|q'-q\|\leq s$, we have $\psi(q')\in \mathscr{D}_s$, hence noting \eqref{insp} and \eqref{upsphi0.1}, we get $\varphi(\psi(q'))=q'$. Consequently, taking 
\begin{equation}\label{Esdd}
    \mathscr{E}_s:=\big\{q'\in \tilde{H}_{k,k'}: \|q'-q\|\leq s, \psi(q')\notin \mathcal{W}'(q;2s,\delta)\big\},
\end{equation}
we have
\begin{equation}\label{Esde}
    \mathscr{E}_s \subseteq \varphi\big(\mathscr{D}_s\backslash \mathcal{W}'(q;2s,\delta)\big). 
\end{equation}
For any $w_1,w_2\in\mathscr{D}_s$, by \cref{L2.0}, 
\begin{equation*}
    \|V(q+w_1)-V(q+w_2)\|\leq \frac{4}{d_0\slash 2}\|w_1-w_2\|=8d_0^{-1}\|w_1-w_2\|, 
\end{equation*}
hence by \eqref{lbdq},
\begin{align*}
   & \bigg\|\frac{\langle w_1,\mathtt{a}_{k,k'}\rangle}{\langle V(q+w_1),\mathtt{a}_{k,k'}\rangle}V(q+w_1)-\frac{\langle w_2,\mathtt{a}_{k,k'}\rangle}{\langle V(q+w_2),\mathtt{a}_{k,k'}\rangle}V(q+w_2)\bigg\|\nonumber\\
   \leq\,& \bigg|\frac{\langle w_1,\mathtt{a}_{k,k'}\rangle}{\langle V(q+w_1),\mathtt{a}_{k,k'}\rangle}-\frac{\langle w_2,\mathtt{a}_{k,k'}\rangle}{\langle V(q+w_2),\mathtt{a}_{k,k'}\rangle}\bigg|+\bigg|\frac{\langle w_2,\mathtt{a}_{k,k'}\rangle}{\langle V(q+w_2),\mathtt{a}_{k,k'}\rangle}\bigg|\|V(q+w_2)-V(q+w_1)\|\nonumber\\
   \leq\,& \frac{|\langle w_1,\mathtt{a}_{k,k'}\rangle| |\langle V(q+w_2)-V(q+w_1),\mathtt{a}_{k,k'}\rangle|+|\langle V(q+w_1),\mathtt{a}_{k,k'}\rangle||\langle w_1-w_2,\mathtt{a}_{k,k'}\rangle|}{\langle V(q+w_1),\mathtt{a}_{k,k'}\rangle \langle V(q+w_2),\mathtt{a}_{k,k'}\rangle}\nonumber\\
   &+\frac{\|w_2\|}{\langle V(q+w_2),\mathtt{a}_{k,k'}\rangle}\cdot 8d_0^{-1}\|w_1-w_2\|\nonumber\\
   \leq\,& \frac{4(\|w_1\|\cdot 8d_0^{-1}\|w_1-w_2\|+\|w_1-w_2\|)}{\langle V(q),\mathtt{a}_{k,k'}\rangle^2}+\frac{2\|w_2\|}{\langle V(q),\mathtt{a}_{k,k'}\rangle}\cdot 8d_0^{-1}\|w_1-w_2\|\nonumber\\
   \leq\,& \frac{(100sd_0^{-1}+4)\|w_1-w_2\|}{\langle V(q),\mathtt{a}_{k,k'}\rangle^2}\leq \frac{5\|w_1-w_2\|}{\langle V(q),\mathtt{a}_{k,k'}\rangle^2},
\end{align*}
which implies
\begin{equation*}
    \|\varphi(w_1)-\varphi(w_2)\|\leq \|w_1-w_2\|+\frac{5\|w_1-w_2\|}{\langle V(q),\mathtt{a}_{k,k'}\rangle^2}\leq \frac{6\|w_1-w_2\|}{\langle V(q),\mathtt{a}_{k,k'}\rangle^2}.
\end{equation*}
By Kirszbraun's theorem, there exists a Lipschitz mapping $\tilde{\varphi}:O(V(q))\rightarrow H_{k,k'}$ such that $\tilde{\varphi}(w)=\varphi(w)$ for all $w \in  \mathscr{D}_s$ and $\|\tilde{\varphi}\|_{\mathrm{Lip}}\leq \frac{6}{\langle V(q),\mathtt{a}_{k,k'}\rangle^2}$. Hence by \eqref{Esde}, \cite[Theorem 2.8]{MR3409135}, and \eqref{bddWqs}, we have (note that $2s\in (0,s_0)$)
\begin{align*}
   \mathcal{H}^{d-1}(\mathscr{E}_s)\leq\,& \mathcal{H}^{d-1}\big(\varphi\big(\mathscr{D}_s\backslash \mathcal{W}'(q;2s,\delta)\big)\big)=\mathcal{H}^{d-1}\big(\tilde{\varphi}\big(\mathscr{D}_s\backslash \mathcal{W}'(q;2s,\delta)\big)\big)\nonumber\\
   \leq\,& \bigg(\frac{6}{\langle V(q),\mathtt{a}_{k,k'}\rangle^2}\bigg)^{d-1}\mathcal{H}^{d-1}(\mathcal{W}(q;2s,\delta))\nonumber\\
   \leq\,&  \bigg(\frac{6}{\langle V(q),\mathtt{a}_{k,k'}\rangle^2}\bigg)^{d-1}\cdot \delta' (2s)^{d-1}=\delta'\bigg(\frac{12s}{\langle V(q),\mathtt{a}_{k,k'}\rangle^2}\bigg)^{d-1}.
\end{align*}
Combining \eqref{bddWqs} and the above display, and noting \eqref{Esdd}, we get
\begin{align}\label{E3.1new}
   & \mathcal{H}^{d-1}\big(\big\{q'\in H_{k,k'}:\|q'-q\|\leq s  \big\}  \big\backslash \big\{q'\in\tilde{H}_{k,k'}:\|q'-q\|\leq s,  \psi(q')\in \mathcal{W}'(q;2s,\delta)\big\}\big)\nonumber\\
   \leq  \,& \mathcal{H}^{d-1}(\{q'\in H_{k,k'} \backslash \tilde{H}_{k,k'}:\|q'-q\|\leq s\})+\mathcal{H}^{d-1}(\mathscr{E}_s)\leq  2\delta'\bigg(\frac{12s}{\langle V(q),\mathtt{a}_{k,k'}\rangle^2}\bigg)^{d-1}.
\end{align}

Below, we consider any $q'\in\tilde{H}_{k,k'}$ such that $\|q'-q\|\leq s$ and $\psi(q')\in \mathcal{W}'(q;2s,\delta)$. By \eqref{insp} and the definition of $\mathcal{W}'(q;2s,\delta)$ (recall \cref{Defd}), we have     
\begin{equation*}
   \|V(q')-V(q)-F(q)\psi(q')\|= \|V(q+\psi(q'))-V(q)-F(q)\psi(q')\|\leq 2\delta  s,
\end{equation*}
which combined with \eqref{E3.1nn} implies
\begin{equation*}
    \|\nabla (\tilde{V}_{k,k'}\circ \Omega_{k,k'}^{-1}) (x_0)\cdot(\Omega_{k,k'}(q')-x_0)-F(q)\psi(q')\|\leq 3\delta s.
\end{equation*}
By \eqref{niceco.1} and \eqref{E3.48}, we have
\begin{align*}
 &\|F(q)\psi(q')-F(q)(q'-q)\|=\|F(q)\psi(q')-F(q)(\mathbf{I}_d-V(q) V(q)^{\top})(q'-q)\|\nonumber\\
    \leq\,& \|F(q)\|_2\|\psi(q')-(\mathbf{I}_d-V(q) V(q)^{\top})(q'-q)\|
    \leq\frac{4}{d_0}\cdot 16sd_0^{-1}\|q'-q\|\leq 64s^2d_0^{-2}.
\end{align*}
By the above two displays,
\begin{equation*}
    \|\nabla (\tilde{V}_{k,k'}\circ \Omega_{k,k'}^{-1}) (x_0)\cdot(\Omega_{k,k'}(q')-x_0)-F(q)(q'-q)\|\leq 3\delta s+64s^2 d_0^{-2}\leq 4\delta s.
\end{equation*}
Hence
\begin{align*}
   & \big\{q'\in\tilde{H}_{k,k'}:\|q'-q\|\leq s,  \psi(q')\in \mathcal{W}'(q;2s,\delta)\big\}\nonumber\\
   \subseteq\, & \big\{q'\in H_{k,k'}:\|q'-q\|\leq s, \|\nabla (\tilde{V}_{k,k'}\circ \Omega_{k,k'}^{-1}) (x_0)\cdot (\Omega_{k,k'}(q')-x_0) -F(q)(q'-q)\|\leq 4\delta s \big\} .
\end{align*}
Noting \eqref{E3.1new}, we thus have
\begin{align*}
    & \mathcal{H}^{d-1}\big(\big\{q'\in H_{k,k'}:\|q'-q\|\leq s,\|\nabla (\tilde{V}_{k,k'}\circ \Omega_{k,k'}^{-1}) (x_0)\cdot(\Omega_{k,k'}(q')-x_0)-F(q)(q'-q)\|> 4\delta s \big\}\big)\nonumber\\
    &  \leq 2\delta'\bigg(\frac{12s}{\langle V(q),\mathtt{a}_{k,k'}\rangle^2}\bigg)^{d-1}.
\end{align*}
Consequently,
\begin{align*}
   & \limsup_{s\rightarrow 0^+}\frac{\mathcal{H}^{d-1}\big(\big\{q'\in H_{k,k'}: \|q'-q\|\leq s,\|\nabla (\tilde{V}_{k,k'}\circ \Omega_{k,k'}^{-1}) (x_0)\cdot(\Omega_{k,k'}(q')-x_0)-F(q)(q'-q)\|> 4\delta s \big\}\big)}{s^{d-1}}\nonumber\\
   & \leq 2\delta'\bigg(\frac{12}{\langle V(q),\mathtt{a}_{k,k'}\rangle^2}\bigg)^{d-1}.
\end{align*}
Taking $\delta'\rightarrow 0^+$ in the above display, we conclude that for any $\delta>0$, 
\begin{equation}\label{Enewco}
    \lim_{s\rightarrow 0^+}\frac{\mathcal{H}^{d-1}\big(\big\{q'\in H_{k,k'}: \|q'-q\|\leq s, \|\nabla (\tilde{V}_{k,k'}\circ \Omega_{k,k'}^{-1}) (x_0)\cdot(\Omega_{k,k'}(q')-x_0)-F(q)(q'-q)\|> 4\delta s \big\}\big)}{s^{d-1}}=0.
\end{equation}

Let $\Delta_{k,k'}:\mathbb{R}^{d-1}\rightarrow\{z\in\mathbb{R}^{d}:\langle z,\mathtt{a}_{k,k'}\rangle =0\}$ be such that
\begin{equation*}
    \Delta_{k,k'}(x):=\sum_{\ell=1}^{d-1}x_{\ell}\mathtt{e}_{k,k'}^{(\ell)},\quad\text{ for all } x=(x_1,\cdots,x_{d-1})\in\mathbb{R}^{d-1}.
\end{equation*}
By \eqref{Enewco}, for any $x^*\in \mathbb{R}^{d-1}$ such that $\|x^*\|=1$ and any $\delta>0$, there exist $(s_n)_{n=1}^{\infty}\subseteq (0,1)$ and $(x_n)_{n=1}^{\infty}\subseteq \mathbb{R}^{d-1}$ such that as $n\rightarrow\infty$,
\begin{equation*}
    s_n\rightarrow 0, \quad s_n\slash 2\leq \|x_n\|\leq s_n,  \quad \frac{x_n}{\|x_n\|}\rightarrow x^*, 
\end{equation*}
\begin{equation*}
    \|\nabla (\tilde{V}_{k,k'}\circ \Omega_{k,k'}^{-1}) (x_0)x_n-F(q)\cdot\Delta_{k,k'}(x_n)\|\leq 4\delta s_n.
\end{equation*}
Note that the above implies that
\begin{equation*}
    \bigg\|\nabla (\tilde{V}_{k,k'}\circ \Omega_{k,k'}^{-1}) (x_0)\frac{x_n}{\|x_n\|}-F(q)\cdot\Delta_{k,k'}\bigg(\frac{x_n}{\|x_n\|}\bigg)  \bigg\|\leq \frac{4\delta s_n}{\|x_n\|}\leq 8\delta.
\end{equation*}
Taking $n\rightarrow\infty$, we get 
\begin{equation*}
    \|\nabla (\tilde{V}_{k,k'}\circ \Omega_{k,k'}^{-1}) (x_0)x^*-F(q)\cdot\Delta_{k,k'}(x^*)\|\leq 8\delta. 
\end{equation*}
Since $\delta>0$ is arbitrary, we get $\nabla (\tilde{V}_{k,k'}\circ \Omega_{k,k'}^{-1}) (x_0)x^*=F(q)\cdot\Delta_{k,k'}(x^*)$. Hence for any $x^*\in\mathbb{R}^{d-1}$, $\nabla (\tilde{V}_{k,k'}\circ \Omega_{k,k'}^{-1}) (x_0)x^*=F(q)\cdot\Delta_{k,k'}(x^*)$. Consequently, for any $\ell\in[d-1]$, we have 
\begin{equation*}
    \frac{\partial (\tilde{V}_{k,k'}\circ \Omega_{k,k'}^{-1})(x_0)}{\partial x_{\ell}} = F(q)\mathtt{e}_{k,k'}^{(\ell)}.
\end{equation*}
By \eqref{pxt1}, the above display implies that for any $t\in (-\beta(q),\alpha(q))$ (note that by \eqref{bddsalphabeta}, $|t|< 2D$),
\begin{equation*}
    \frac{\partial \EuScript{P}(x_0,t)}{\partial x_{\ell}}=\mathtt{e}_{k,k'}^{(\ell)}+t\frac{\partial (\tilde{V}_{k,k'}\circ\Omega_{k,k'}^{-1})(x_0)}{\partial x_{\ell}}=(\mathbf{I}_d+tF(q))\mathtt{e}_{k,k'}^{(\ell)}.
\end{equation*}
Moreover, by \eqref{pxt0.1}, we have
\begin{equation*}
    \frac{\partial \EuScript{P}(x_0,t)}{\partial t} = (\tilde{V}_{k,k'}\circ\Omega_{k,k'}^{-1})(x_0) = \tilde{V}_{k,k'}(q) = V(q).
\end{equation*}
Hence by \eqref{defRk}, noting that $F(q)V(q)=0$ (recall \eqref{niceco.1}), we conclude that
\begin{align*}
    R_{k,k'}(q,t)=\, & \Big|\det\Big[V(q),(\mathbf{I}_d+tF(q))\mathtt{e}_{k,k'}^{(1)},\cdots,(\mathbf{I}_d+tF(q))\mathtt{e}_{k,k'}^{(d-1)}\Big]\Big| \nonumber\\
    =\, & \Big|\det \Big((\mathbf{I}_d+tF(q)) \Big[V(q),\mathtt{e}_{k,k'}^{(1)},\cdots, \mathtt{e}_{k,k'}^{(d-1)}\Big]\Big)\Big|\nonumber\\
    =\, & |\det(\mathbf{I}_d+tF(q))| \cdot \Big|\det\Big[V(q),\mathtt{e}_{k,k'}^{(1)},\cdots, \mathtt{e}_{k,k'}^{(d-1)}\Big]\Big|
    =\det(\mathbf{I}_d+tF(q)) |\langle\mathtt{a}_{k,k'},V(q)\rangle|,
\end{align*}
where for the last equality, we use \eqref{niceco.1} and the fact that
\begin{align*}
&\Big[V(q),\mathtt{e}_{k,k'}^{(1)},\cdots, \mathtt{e}_{k,k'}^{(d-1)}\Big]=\Big[\mathtt{a}_{k,k'},\mathtt{e}_{k,k'}^{(1)},\cdots, \mathtt{e}_{k,k'}^{(d-1)}\Big]\begin{bmatrix}
        \langle  \mathtt{a}_{k,k'}, V(q)\rangle & 0 & \cdots & 0\\
        \langle \mathtt{e}_{k,k'}^{(1)}, V(q)\rangle & 1  & \cdots & 0\\
       &  \cdots  &&&  \\
        \langle  \mathtt{e}_{k,k'}^{(d-1)}, V(q)\rangle & 0 & \cdots & 1\\
    \end{bmatrix} \\
    \Rightarrow\, & \left|\det\Big[V(q),\mathtt{e}_{k,k'}^{(1)},\cdots, \mathtt{e}_{k,k'}^{(d-1)}\Big]\right|=\left|\det\begin{bmatrix}
        \langle  \mathtt{a}_{k,k'}, V(q)\rangle & 0 & 0 & \cdots & 0\\
        \langle \mathtt{e}_{k,k'}^{(1)}, V(q)\rangle & 1 & 0 & \cdots & 0\\
       &  \cdots  &&&  \\
        \langle  \mathtt{e}_{k,k'}^{(d-1)}, V(q)\rangle & 0 & 0 & \cdots & 1\\
    \end{bmatrix}\right|=\left|\langle\mathtt{a}_{k,k'}, V(q)\rangle\right|.  \qedhere
\end{align*}
\end{proof}

The next lemma introduces Borel sets $\mathfrak{H}_{k,k'}\subseteq \tilde{H}_{k,k'}$ for $k\in[K],k'\in [K']$ and $\mathfrak{T}\subseteq \mathcal{T}_1^{*}$ with certain desirable properties. 

\begin{lemma}\label{Defn4.1.1n}
For each $k\in[K],k'\in[K']$, there exists a Borel set $\mathfrak{H}_{k,k'}\subseteq \tilde{H}_{k,k'}$ such that
\begin{equation}\label{new.con}
    \int_{\mathbb{R}} \mathbbm{1}_{(-\beta(q),\alpha(q))}(t) 
R_{k,k'}(q,t) dt=0 \quad\text{ for }\mathcal{H}^{d-1}\text{-a.e.\ }q\in \tilde{H}_{k,k'}\backslash \mathfrak{H}_{k,k'}, 
\end{equation}
and for every $q\in\mathfrak{H}_{k,k'}$, we have
\begin{equation}\label{qrelations}
    \int_{T(q)}f(x)\det\big(\mathbf{I}_d+\langle x-q, V(q)\rangle F(q)\big)d\mathcal{H}^1(x)=\int_{T(q)}g(x)\det\big(\mathbf{I}_d+\langle x-q, V(q)\rangle F(q)\big)d\mathcal{H}^1(x),
\end{equation}
\eqref{L3.10E1} holds for all $t\in (-\beta(q),\alpha(q))$, and the conclusions of Theorems~\ref{L2.2}--\ref{L2.4} and \cref{L2.7} hold for all $z_0\in\mathrm{int}(T(q))$. Define (recall \eqref{defWkk})
\begin{equation}\label{defT}
    \mathfrak{T}:=\bigcup_{k\in[K],k'\in[K']}\mathfrak{W}_{k,k';\bar{H}_{k,k'}\cap \mathfrak{H}_{k,k'}} \subseteq \mathcal{T}_1^{*}.
\end{equation}
Then $\mathfrak{T}$ is a Borel transport set, and the conclusions of Theorems~\ref{L2.2}--\ref{L2.4} and \cref{L2.7} hold for every $z_0\in\mathfrak{T}$. 
\end{lemma}
\begin{proof}
Let $\mathfrak{T}_0$ be the set of $z_0\in\mathcal{T}_1^*$ for which the conclusions of Theorems~\ref{L2.2}--\ref{L2.4} and \cref{L2.7} hold. Note that $\mathfrak{T}_0$ is a transport set. By Theorems~\ref{L2.2}--\ref{L2.4} and \cref{L2.7}, $\mathfrak{T}_0$ is Lebesgue measurable and $\mathcal{L}^d(\mathcal{T}_1^*\backslash\mathfrak{T}_0)=0$. For each $k\in [K],k'\in[K']$, let $\mathfrak{T}_{1;k,k'}$ be the set of $z_0\in\tilde{\mathcal{T}}_{1;k,k'}^*$ for which the conclusions of \cref{L2.5n,Lem3.10n} hold. Note that $\mathfrak{T}_{1;k,k'}$ is a transport set. By \cref{L2.5n,Lem3.10n}, $\mathfrak{T}_{1;k,k'}$ is Lebesgue measurable and $\mathcal{L}^d(\tilde{\mathcal{T}}_{1;k,k'}^*\backslash \mathfrak{T}_{1;k,k'})=0$.

Now for each $k\in[K],k'\in[K']$, by \cref{L2.3n} and the fact that $\mathcal{L}^d(\tilde{\mathcal{T}}_{1;k,k'}^*\backslash (\mathfrak{T}_0\cap\mathfrak{T}_{1;k,k'}))=0$, there exists a Borel set $\bar{\mathfrak{H}}_{k,k'}\subseteq \tilde{H}_{k,k'}$, such that $\mathcal{H}^{d-1}(\tilde{H}_{k,k'}\backslash \bar{\mathfrak{H}}_{k,k'})=0$, and for every $q\in  \bar{\mathfrak{H}}_{k,k'}$, 
\begin{equation}\label{realtion1}
    \int_{\mathbb{R}} \mathbbm{1}_{(-\beta(q),\alpha(q))}(t) 
    \mathbbm{1}_{q+tV(q)\in \tilde{\mathcal{T}}_{1;k,k'}^*\backslash (\mathfrak{T}_0\cap\mathfrak{T}_{1;k,k'})} R_{k,k'}(q,t) dt=0.
\end{equation}
As $\tilde{\mathcal{T}}_{1;k,k'}^*\backslash (\mathfrak{T}_0\cap\mathfrak{T}_{1;k,k'})$ is a transport set, we have $\mathbbm{1}_{q+tV(q)\in \tilde{\mathcal{T}}_{1;k,k'}^*\backslash (\mathfrak{T}_0\cap\mathfrak{T}_{1;k,k'})}=\mathbbm{1}_{q\in \tilde{\mathcal{T}}_{1;k,k'}^*\backslash (\mathfrak{T}_0\cap\mathfrak{T}_{1;k,k'})}$ for any $q\in\tilde{H}_{k,k'}$ and $t\in (-\beta(q),\alpha(q))$. Hence \eqref{realtion1} implies that either $q\in \mathfrak{T}_0\cap\mathfrak{T}_{1;k,k'}$ or $\int_{\mathbb{R}} \mathbbm{1}_{(-\beta(q),\alpha(q))}(t) 
R_{k,k'}(q,t) dt=0$. We define $\mathfrak{H}_{k,k'}$ to be the set of $q\in \bar{\mathfrak{H}}_{k,k'}$ such that 
\[\int_{\mathbb{R}} \mathbbm{1}_{(-\beta(q),\alpha(q))}(t) 
R_{k,k'}(q,t) dt>0\]
(note that this implies $q\in\mathfrak{T}_0\cap\mathfrak{T}_{1;k,k'}$). Note that $\mathfrak{H}_{k,k'}\subseteq \tilde{H}_{k,k'}$ is a Borel set and \eqref{new.con} holds. As $\mathfrak{T}_0\cap\mathfrak{T}_{1;k,k'}$ is a transport set, we have $\mathrm{int}(T(q))\subseteq \mathfrak{T}_0\cap\mathfrak{T}_{1;k,k'}$ for every $q\in \mathfrak{H}_{k,k'}$. In particular, \eqref{L3.10E1} (for all $t\in (-\beta(q),\alpha(q))$) and \eqref{qrelations} hold for every $q\in\mathfrak{H}_{k,k'}$, and the conclusions of Theorems~\ref{L2.2}--\ref{L2.4} and \cref{L2.7} hold for every $z_0\in\mathrm{int}(T(q))$ when $q\in\mathfrak{H}_{k,k'}$. 

By \cref{L2.1n}, $\mathfrak{T}$ as defined in \eqref{defT} is a Borel transport set. Note that from the previous paragraph, for every $z_0\in\mathfrak{T}$, the conclusions of Theorems~\ref{L2.2}--\ref{L2.4} and \cref{L2.7} hold.
\end{proof}

The following lemma shows that the ratio \eqref{ratio_eq} is independent of the choice of the base point $x_0$, and that the two marginal densities are balanced on each transport ray. 

\begin{lemma}\label{Ln2.2}
Fix any $x\in\mathfrak{T}$.
\begin{itemize}
    \item[(a)] The strict positivity of $\int_{T(x)}f(z)\det\big(\mathbf{I}_d+\langle z-x_0, V(x)\rangle F(x_0)\big)d\mathcal{H}^1(z)$ and the value of 
\begin{equation}\label{ratio_eq}
    \frac{\det\big(\mathbf{I}_d+\langle x-x_0, V(x)\rangle F(x_0)\big)}{\int_{T(x)}f(z)\det\big(\mathbf{I}_d+\langle z-x_0, V(x)\rangle F(x_0)\big)d\mathcal{H}^1(z)}
\end{equation}
are independent of the choice of $x_0\in\mathrm{int}(T(x))$.
   \item[(b)] For any $x_0\in\mathrm{int}(T(x))$, we have
\begin{equation*}
\int_{T(x)}\!\!f(z)\det\!\big(\mathbf{I}_d+\langle z-x_0, V(x) \rangle F(x_0)\big)d\mathcal{H}^1(z) = \int_{T(x)}\!\!g(z)\det\!\big(\mathbf{I}_d+\langle z-x_0, V(x) \rangle F(x_0)\big)d\mathcal{H}^1(z).
\end{equation*}
\end{itemize}
\end{lemma}
\begin{proof}

As $x\in\mathfrak{T}$, for any $x_0,x_0',z\in\mathrm{int}(T(x))$, by \cref{L2.4} (using \eqref{Eq2.6.4} with $x,x',x_0$ replaced by $x_0,x_0',z$), we have 
\begin{equation*}
    \det\big(\mathbf{I}_d+\langle z-x_0', V(x)\rangle F(x_0')\big)=\det\big(\mathbf{I}_d+\langle x_0-x_0', V(x)\rangle F(x_0')\big)\det\big(\mathbf{I}_d+\langle z-x_0, V(x)\rangle F(x_0)\big). 
\end{equation*}
Taking $z=x$, we obtain that
\begin{equation*}
    \det\big(\mathbf{I}_d+\langle x-x_0', V(x)\rangle F(x_0')\big)=\det\big(\mathbf{I}_d+\langle x_0-x_0', V(x)\rangle F(x_0')\big)\det\big(\mathbf{I}_d+\langle x-x_0, V(x)\rangle F(x_0)\big).
\end{equation*}
Moreover,
\begin{align*}
&\int_{T(x)}f(z)\det\big(\mathbf{I}_d+\langle z-x_0', V(x)\rangle F(x_0')\big)d\mathcal{H}^1(z)\nonumber\\
=\,&\det\big(\mathbf{I}_d+\langle x_0-x_0', V(x)\rangle F(x_0')\big) \int_{T(x)}f(z)\det\big(\mathbf{I}_d+\langle z-x_0, V(x)\rangle F(x_0)\big)d\mathcal{H}^1(z).
\end{align*}
By \cref{L2.7}, $\det\big(\mathbf{I}_d+\langle x_0-x_0', V(x)\rangle F(x_0')\big)>0$. Thus part (a) follows from the two displays.

We proceed to the proof of part (b). By \eqref{defT}, as $x\in\mathfrak{T}$, there exist $k\in[K],k'\in[K']$ and $q\in \bar{H}_{k,k'}\cap\mathfrak{H}_{k,k'}$ such that $x\in \mathrm{int}(T(q))$. By \eqref{qrelations}, noting that $T(x)=T(q)$ and $V(x)=V(q)$, 
\begin{equation*}
    \int_{T(x)}f(z)\det\big(\mathbf{I}_d+\langle z-q, V(x)\rangle F(q)\big)d\mathcal{H}^1(z)=\int_{T(x)}g(z)\det\big(\mathbf{I}_d+\langle z-q, V(x)\rangle F(q)\big)d\mathcal{H}^1(z).
\end{equation*}
For any $x_0,z\in\mathrm{int}(T(x))$, by \cref{L2.4} (using \eqref{Eq2.6.4} with $x,x',x_0$ replaced by $x_0,q,z$), we have
\begin{equation*}
    \det\big(\mathbf{I}_d+\langle z-q, V(x)\rangle F(q)\big)=\det\big(\mathbf{I}_d+\langle x_0-q, V(x)\rangle F(q)\big)  \det\big(\mathbf{I}_d+\langle z-x_0, V(x)\rangle F(x_0)\big).
\end{equation*}
By \cref{L2.7}, $\det\big(\mathbf{I}_d+\langle x_0-q, V(x)\rangle F(q)\big)>0$. Therefore, the claim follows by combining the two displays.
\end{proof}

\begin{definition}[The sets $\mathfrak{T}^{*}$, $\mathcal{S}^{*}$ and the function $\mathfrak{F}$]\label{factors}
We define $\mathfrak{T}^{*}$ to be the set of $x\in\mathfrak{T}$ such that $\int_{T(x)}f(z)\det\big(\mathbf{I}_d+\langle z-x_0, V(x)\rangle F(x_0)\big)d\mathcal{H}^1(z)>0$
for some (hence all) $x_0\in\mathrm{int}(T(x))$. Note that $\mathfrak{T}^{*}$ is a transport set. We further write $\mathcal{S}^{*}:=\{T(x):x\in\mathfrak{T}^*\}$.

For any $x\in\mathfrak{T}^*$, we define
\begin{equation*}
    \mathfrak{F}(x):=\frac{\det(\mathbf{I}_d+\langle x-x_0, V(x)\rangle F(x_0))}{\int_{T(x)}f(z)\det(\mathbf{I}_d+\langle z-x_0, V(x)\rangle F(x_0))d\mathcal{H}^1(z)}.
\end{equation*}
For any $x\in\mathbb{R}^d\backslash\mathfrak{T}^*$, we define $\mathfrak{F}(x):=0$. Note that by \cref{Ln2.2}, $\mathfrak{T}^{*}$ and $\mathfrak{F}(x)$ are independent of the choice of $x_0\in\mathrm{int}(T(x))$. 
\end{definition}
\begin{remark}\label{positivityF}
By \cref{L2.7}, $\mathfrak{F}(x)>0$ for $x \in \mathfrak{T}^{*}$, and $\mathfrak{F}(x)= 0$ for $x\in\mathbb{R}^d\backslash \mathfrak{T}^{*}$. 
\end{remark}

\begin{definition}\label{defkk}
For any $k\in[K],k'\in[K']$ and $x\in \tilde{\mathcal{T}}^{*}_{1;k,k'}$, we define $\mathfrak{q}_{k,k'}(x)\in\tilde{H}_{k,k'}$ to be the intersection point of $T(x)$ and $H_{k,k'}$. Note that by \eqref{Hkkdef}, for any $x\in\tilde{\mathcal{T}}^{*}_{1;k,k'}$, we have
\begin{equation*}
    \mathfrak{q}_{k,k'}(x)=x-\frac{\langle x, \mathtt{a}_{k,k'}\rangle-\mathtt{b}_{k,k'}}{\langle V(x), \mathtt{a}_{k,k'}\rangle}V(x). 
\end{equation*}
In particular, the mapping $\mathfrak{q}_{k,k'}:\tilde{\mathcal{T}}^{*}_{1;k,k'}\rightarrow H_{k,k'}$ is Borel measurable.
\end{definition}

\begin{lemma}\label{L3n}
$\mathfrak{T}^{*}$ is a Borel subset of $\mathbb{R}^d$ and $\mathcal{S}^*$ is a Borel subset of $\mathcal{S}$. Moreover, the mapping $\mathfrak{F}:\mathbb{R}^d\rightarrow [0,\infty)$ is Borel measurable.
\end{lemma}
\begin{proof}

Fix any $k\in[K],k'\in[K']$. For any $q\in \bar{H}_{k,k'}\cap\mathfrak{H}_{k,k'}$, we define 
\begin{align*}
    \mathfrak{I}_{k,k'}(q):=&\int_{T(q)}f(z)\det\big(\mathbf{I}_d+\langle z-q, V(q)\rangle F(q)\big)d\mathcal{H}^1(z) \nonumber\\
    =&\int_{\mathbb{R}}\mathbbm{1}_{(-\beta(q),\alpha(q))}(t)f(q+tV(q))\det(\mathbf{I}_d+t F(q))dt.
\end{align*}
As the function $(q,t)\mapsto\mathbbm{1}_{(-\beta(q),\alpha(q))}(t)f(q+tV(q))\det(\mathbf{I}_d+t F(q))$ on the domain $(\bar{H}_{k,k'}\cap\mathfrak{H}_{k,k'})\times\mathbb{R}$ is Borel measurable, by Fubini's theorem (\cite[Section 8.8]{MR924157}), $\mathfrak{I}_{k,k'}:\bar{H}_{k,k'}\cap\mathfrak{H}_{k,k'}\rightarrow [0,\infty)$ is Borel measurable. Hence $\{q\in \bar{H}_{k,k'}\cap\mathfrak{H}_{k,k'}:\mathfrak{I}_{k,k'}(q)>0\}$ is a Borel set. By \eqref{defT} (also note \cref{factors}), we have
\begin{equation}\label{Tkkk}
    \mathfrak{T}^{*}\cap\mathcal{T}^*_{1;k,k'}=\mathfrak{W}_{k,k';\{q\in \bar{H}_{k,k'}\cap\mathfrak{H}_{k,k'}:\mathfrak{I}_{k,k'}(q)>0\}}.
\end{equation}
Hence by \cref{L2.1n}, $\mathfrak{T}^{*}\cap\mathcal{T}^*_{1;k,k'}$ is Borel measurable. By \eqref{disj},
\begin{equation}\label{Tkkkk}
    \mathfrak{T}^{*}=\bigcup_{k\in[K],k'\in[K']}(\mathfrak{T}^{*}\cap\mathcal{T}^*_{1;k,k'}).
\end{equation}
Hence $\mathfrak{T}^*$ is Borel subset of $\mathbb{R}^d$. Moreover, 
\begin{equation*}
   \mathcal{S}^*= \{T(x):x\in\mathfrak{T}^*\}=\bigcup_{k\in[K],k'\in[K']}\mathcal{S}_{k,k';\{q\in \bar{H}_{k,k'}\cap\mathfrak{H}_{k,k'}:\mathfrak{I}_{k,k'}(q)>0\}}.
\end{equation*}
Hence by \cref{L2.4n}, $\mathcal{S}^*$ is a Borel subset of $\mathcal{S}$. 

For any $k\in[K], k'\in [K']$ and $x\in \mathfrak{T}^{*}\cap\mathcal{T}^*_{1;k,k'}$, by \eqref{Tkkk}, we have $\mathfrak{q}_{k,k'}(x)\in \bar{H}_{k,k'}\cap\mathfrak{H}_{k,k'}$ and $\mathfrak{I}_{k,k'}(\mathfrak{q}_{k,k'}(x))>0$ (recall \cref{defkk}). As $\mathfrak{I}_{k,k'}:\bar{H}_{k,k'}\cap\mathfrak{H}_{k,k'}\rightarrow [0,\infty)$ is Borel measurable, the function $x\mapsto \mathfrak{I}_{k,k'}(\mathfrak{q}_{k,k'}(x))$ on the domain $\mathfrak{T}^{*}\cap\mathcal{T}^*_{1;k,k'}$ is also Borel measurable. Note that for any $x\in \mathfrak{T}^{*}\cap\mathcal{T}^*_{1;k,k'}$, we have (recall \cref{factors})
\begin{equation*}
    \mathfrak{F}(x)=\frac{\det(\mathbf{I}_d+\langle x-\mathfrak{q}_{k,k'}(x), V(x)\rangle F(\mathfrak{q}_{k,k'}(x)))}{\mathfrak{I}_{k,k'}(\mathfrak{q}_{k,k'}(x))},
\end{equation*}
hence $\mathfrak{F}|_{\mathfrak{T}^{*}\cap\mathcal{T}^*_{1;k,k'}}$ is Borel measurable for every $k\in[K],k'\in[K']$. Moreover, by \cref{factors}, $\mathfrak{F}(x)=0$ for every $x\in\mathbb{R}^d\backslash\mathfrak{T}^*$. Hence by \eqref{Tkkkk}, the mapping $\mathfrak{F}:\mathbb{R}^d\rightarrow [0,\infty)$ is Borel measurable.
\end{proof}

\begin{definition}[The measures $\lambda,\tilde{\mu}_T,\tilde{\nu}_T$ and the functions $\tilde{f},\tilde{g}$]\label{De3.8}
With $T:\mathcal{T}_1^*\rightarrow \mathcal{S}$ as in \cref{TransportRays} (note that $T$ is Borel measurable by \cref{L2.0nn}), we define $\lambda$ to be the Borel measure on $\mathcal{S}$ such that for any Borel set $A\subseteq \mathcal{S}$, $\lambda(A)=\mu(\{x\in\mathcal{T}_1^*:T(x)\in A\})$. Note that by \eqref{new1.1}, $\mu(\mathcal{T}_1^*)=1$, hence $\lambda(\mathcal{S})=1$ and $\lambda$ is a probability measure. 

For any $x\in\mathbb{R}^d$, we define
\begin{equation}\label{defft}
    \tilde{f}(x):=f(x)\mathfrak{F}(x), \qquad \tilde{g}(x):=g(x)\mathfrak{F}(x).
\end{equation}
Note that by \cref{L3n}, $\tilde{f},\tilde{g}:\mathbb{R}^d\rightarrow [0,\infty)$ are Borel measurable. For any $T\in\mathcal{S}$, we define $\tilde{\mu}_T$ and $\tilde{\nu}_T$ to be Borel measures on $\mathbb{R}^d$ such that for any Borel set $A\subseteq\mathbb{R}^d$, 
\begin{equation}\label{ress}
    \tilde{\mu}_T(A)=\int_{T} \tilde{f}(x)\mathbbm{1}_A(x) d\mathcal{H}^1(x), \qquad \tilde{\nu}_T(A)=\int_{T} \tilde{g}(x)\mathbbm{1}_A(x) d\mathcal{H}^1(x). 
\end{equation}
\end{definition}

The following lemma derives explicit formulas for $\lambda$.

\begin{lemma}\label{L3.14n}
For any Borel set $A\subseteq \mathcal{S}$, we have 
\begin{align}\label{L3.14n1}
    &\lambda(A)=\nu(\{x\in\mathcal{T}_1^*:T(x)\in A\})\nonumber\\
    =\,& \sum_{k=1}^K\sum_{k'=1}^{K'}\int_{\bar{H}_{k,k'}\cap\mathfrak{H}_{k,k'}\cap T^{-1}(A)} |\langle V(q), \mathtt{a}_{k,k'}\rangle| d\mathcal{H}^{d-1}(q)\int_{T(q)}f(x)\det\big(\mathbf{I}_d+\langle x-q, V(q)\rangle F(q)\big)d\mathcal{H}^1(x)\nonumber\\
    =\,& \sum_{k=1}^K\sum_{k'=1}^{K'}\int_{\bar{H}_{k,k'}\cap\mathfrak{H}_{k,k'}\cap T^{-1}(A)}|\langle V(q), \mathtt{a}_{k,k'}\rangle| d\mathcal{H}^{d-1}(q)\int_{T(q)}g(x)\det\big(\mathbf{I}_d+\langle x-q, V(q)\rangle F(q)\big)d\mathcal{H}^1(x),
\end{align}  
where $\mathtt{a}_{k,k'},\bar{H}_{k,k'},\mathfrak{H}_{k,k'}$ are as in \eqref{Hkkdef}, \eqref{deftildehat}, and \cref{Defn4.1.1n}, respectively. Moreover, for any $k\in[K],k'\in[K']$ and any Borel set $A\subseteq\mathcal{S}$ such that $T^{-1}(A)\subseteq\tilde{\mathcal{T}}_{1;k,k'}^*$, we have
\begin{equation*}
\lambda(A)=\int_{\mathfrak{H}_{k,k'}\cap T^{-1}(A)}|\langle V(q), \mathtt{a}_{k,k'}\rangle| d\mathcal{H}^{d-1}(q)\int_{T(q)}f(x)\det\big(\mathbf{I}_d+\langle x-q, V(q)\rangle F(q)\big)d\mathcal{H}^1(x).
\end{equation*}
\end{lemma}
\begin{proof}

By \eqref{disj} and \cref{L2.3n}, for any Borel set $A\subseteq\mathcal{S}$, we have
\begin{align*}
    \lambda(A)&=\mu(\{x\in\mathcal{T}_1^{*}:T(x)\in A\})=\sum_{k=1}^K\sum_{k'=1}^{K'}\mu(\{x\in\mathcal{T}_{1;k,k'}^*:T(x)\in A\})\nonumber\\
    &=\sum_{k=1}^K\sum_{k'=1}^{K'}\int_{\mathcal{T}_{1;k,k'}^*}\mathbbm{1}_{T(x)\in A} f(x)  dx \nonumber\\
&=\sum_{k=1}^K\sum_{k'=1}^{K'}\int_{\bar{H}_{k,k'}\cap T^{-1}(A)}d\mathcal{H}^{d-1}(q)\int_{\mathbb{R}}\mathbbm{1}_{(-\beta(q),\alpha(q))}(t)f(q+tV(q))
R_{k,k'}(q,t)dt.
\end{align*}
For any $k\in[K],k'\in[K']$, with $\mathfrak{H}_{k,k'}$ as defined in \cref{Defn4.1.1n}, we have
\begin{align}\label{L3.14n3}
   & \int_{\bar{H}_{k,k'}\cap T^{-1}(A)}d\mathcal{H}^{d-1}(q)\int_{\mathbb{R}}\mathbbm{1}_{(-\beta(q),\alpha(q))}(t)f(q+tV(q))R_{k,k'}(q,t)dt\nonumber\\
   =\,& \int_{\bar{H}_{k,k'}\cap\mathfrak{H}_{k,k'}\cap T^{-1}(A)}d\mathcal{H}^{d-1}(q)\int_{\mathbb{R}}\mathbbm{1}_{(-\beta(q),\alpha(q))}(t)f(q+tV(q))R_{k,k'}(q,t)dt\nonumber\\
   =\,& \int_{\bar{H}_{k,k'}\cap\mathfrak{H}_{k,k'}\cap T^{-1}(A)}|\langle V(q), \mathtt{a}_{k,k'}\rangle|d\mathcal{H}^{d-1}(q)\int_{\mathbb{R}}\mathbbm{1}_{(-\beta(q),\alpha(q))}(t)f(q+tV(q))\det(\mathbf{I}_d+tF(q))dt\nonumber\\
   =\,& \int_{\bar{H}_{k,k'}\cap\mathfrak{H}_{k,k'}\cap T^{-1}(A)}|\langle V(q), \mathtt{a}_{k,k'}\rangle|d\mathcal{H}^{d-1}(q)\int_{T(q)}f(x)\det\big(\mathbf{I}_d+\langle x-q, V(q)\rangle F(q)\big)d\mathcal{H}^1(x),
\end{align}
where we use \eqref{new.con} and \eqref{L3.10E1} in the first and second equalities, respectively. Hence we have 
\begin{align*}
    &\lambda(A)\nonumber\\
    =\,& \sum_{k=1}^K\sum_{k'=1}^{K'}\int_{\bar{H}_{k,k'}\cap\mathfrak{H}_{k,k'}\cap T^{-1}(A)}|\langle V(q), \mathtt{a}_{k,k'}\rangle|d\mathcal{H}^{d-1}(q)\int_{T(q)}f(x)\det\big(\mathbf{I}_d+\langle x-q, V(q)\rangle F(q)\big)d\mathcal{H}^1(x).
\end{align*}
Similarly, we can deduce that
\begin{align*}
    &\nu(\{x\in\mathcal{T}_1^{*}:T(x)\in A\})\nonumber\\
    =\,& \sum_{k=1}^K\sum_{k'=1}^{K'}\int_{\bar{H}_{k,k'}\cap\mathfrak{H}_{k,k'}\cap T^{-1}(A)}|\langle V(q), \mathtt{a}_{k,k'}\rangle|d\mathcal{H}^{d-1}(q)\int_{T(q)}g(x)\det\big(\mathbf{I}_d+\langle x-q, V(q)\rangle F(q)\big)d\mathcal{H}^1(x).
\end{align*}
As \eqref{qrelations} holds for any $k\in[K],k'\in[K']$ and $q\in \mathfrak{H}_{k,k'}$, \eqref{L3.14n1} follows from the above two displays. 

Now for any $k\in[K],k'\in[K']$ and any Borel set $A$ such that $T^{-1}(A)\subseteq\tilde{\mathcal{T}}_{1;k,k'}^*$, by Lemma \ref{L2.3n} and arguing similarly as in \eqref{L3.14n3}, we obtain that  
\begin{align*}
    \lambda(A)=\,&\mu(\{x\in\mathcal{T}_1^*:T(x)\in A\})=\mu(\{x\in\tilde{\mathcal{T}}_{1;k,k'}^*:T(x)\in A\})\nonumber\\
    =\,&\int_{\tilde{\mathcal{T}}_{1;k,k'}^*}\mathbbm{1}_{T(x)\in A}f(x)dx=\int_{\tilde{H}_{k,k'}\cap T^{-1}(A)}d\mathcal{H}^{d-1}(q)\int_{\mathbb{R}}\mathbbm{1}_{(-\beta(q),\alpha(q))}(t)f(q+tV(q))
R_{k,k'}(q,t)dt\nonumber\\
=\,& \int_{\mathfrak{H}_{k,k'}\cap T^{-1}(A)}|\langle V(q), \mathtt{a}_{k,k'}\rangle|d\mathcal{H}^{d-1}(q)\int_{T(q)}f(x)\det\big(\mathbf{I}_d+\langle x-q, V(q)\rangle F(q)\big)d\mathcal{H}^1(x). \qedhere
\end{align*}
\end{proof}

Next, we show that the marginals $\mu$ and $\nu$ are concentrated on the set $\mathfrak{T}^*$. 

\begin{lemma}\label{Lem3.15}
We have $\mu(\mathbb{R}^d\backslash\mathfrak{T}^{*})=0$ and $\nu(\mathbb{R}^d\backslash \mathfrak{T}^*)=0$.
\end{lemma}
\begin{proof}

By \eqref{new1.1}, $\mu(\mathbb{R}^d\backslash\mathfrak{T}^{*})=\mu(\mathcal{T}_1^*\backslash\mathfrak{T}^{*})$. Hence by \eqref{disj}, we have
\begin{equation}\label{L3.15n1}
    \mu(\mathbb{R}^d\backslash\mathfrak{T}^{*})= \sum_{k=1}^K\sum_{k'=1}^{K'}\mu(\mathcal{T}_{1;k,k'}^{*}\backslash \mathfrak{T}^{*}).
\end{equation}
Now consider any $k\in[K],k'\in[K']$. By \cref{L2.3n}, 
\begin{align*}
&\mu(\mathcal{T}_{1;k,k'}^{*}\backslash \mathfrak{T}^{*})=\int_{\mathcal{T}_{1;k,k'}^{*}}\mathbbm{1}_{x\notin \mathfrak{T}^{*}}f(x)dx\nonumber\\
=\,&\int_{\bar{H}_{k,k'}}d\mathcal{H}^{d-1}(q)\int_{\mathbb{R}}\mathbbm{1}_{(-\beta(q),\alpha(q))}(t)\mathbbm{1}_{q+tV(q)\notin\mathfrak{T}^*}  f(q+tV(q))R_{k,k'}(q,t)dt\nonumber\\
=\,&\int_{\bar{H}_{k,k'}\cap\mathfrak{H}_{k,k'}}d\mathcal{H}^{d-1}(q)\int_{\mathbb{R}}\mathbbm{1}_{(-\beta(q),\alpha(q))}(t)\mathbbm{1}_{q+tV(q)\notin\mathfrak{T}^*}
f(q+tV(q))R_{k,k'}(q,t)dt\nonumber\\
=\,& \int_{\bar{H}_{k,k'}\cap\mathfrak{H}_{k,k'}}|\langle V(q), \mathtt{a}_{k,k'}\rangle| d\mathcal{H}^{d-1}(q)\int_{\mathbb{R}}\mathbbm{1}_{(-\beta(q),\alpha(q))}(t)\mathbbm{1}_{q+tV(q)\notin\mathfrak{T}^*}
f(q+tV(q))\det(\mathbf{I}_d+tF(q))dt\nonumber\\
=\,&
\int_{\bar{H}_{k,k'}\cap\mathfrak{H}_{k,k'}}|\langle V(q), \mathtt{a}_{k,k'}\rangle| d\mathcal{H}^{d-1}(q)\int_{T(q)}\mathbbm{1}_{x\notin\mathfrak{T}^*}
f(x)\det\big(\mathbf{I}_d+\langle x-q,V(q)\rangle F(q)\big)d\mathcal{H}^1(x),
\end{align*}
where we use \eqref{new.con} and \eqref{L3.10E1} in the third and fourth equalities, respectively. Now noting \eqref{defT}, for any $q\in \bar{H}_{k,k'}\cap\mathfrak{H}_{k,k'}$ and $x\in \mathrm{int}(T(q))$, we have $q,x\in\mathfrak{T}$, and $x\notin\mathfrak{T}^*$ if and only if $q\notin \mathfrak{T}^*$. Hence for any $q\in \bar{H}_{k,k'}\cap\mathfrak{H}_{k,k'}$,
\begin{align*}
&\int_{T(q)}\mathbbm{1}_{x\notin\mathfrak{T}^*}
f(x)\det\big(\mathbf{I}_d+\langle x-q,V(q)\rangle F(q)\big)d\mathcal{H}^1(x)\nonumber\\
=\, &\mathbbm{1}_{q\notin\mathfrak{T}^*}\int_{T(q)}f(x)\det\big(\mathbf{I}_d+\langle x-q,V(q)\rangle F(q)\big)d\mathcal{H}^1(x)=0,
\end{align*}
where we use the definition of $\mathfrak{T}^*$ in the last equality. Combining the above two displays, we get $\mu(\mathcal{T}_{1;k,k'}^{*}\backslash \mathfrak{T}^{*})=0$ for every $k\in[K],k'\in[K']$. Hence by \eqref{L3.15n1}, we conclude that $\mu(\mathbb{R}^d\backslash\mathfrak{T}^{*})=0$. As $\mathbb{R}^d\backslash\mathfrak{T}^*$ is a Borel transport set, by \cite[Lemma 27]{MR1862796}, we have $\nu(\mathbb{R}^d\backslash\mathfrak{T}^*)=\mu(\mathbb{R}^d\backslash\mathfrak{T}^*)=0$. 
\end{proof}

We record the following measurability lemma for later use.

\begin{lemma}\label{L3.15n}
For any nonnegative Borel measurable function $\varphi:\mathbb{R}^d\rightarrow [0,\infty)$, the mapping $\mathcal{S}\ni T\mapsto \int_T \varphi(x)d\mathcal{H}^1(x)$ is Borel measurable. 
\end{lemma}
\begin{proof}

For any $T\in\mathcal{S}$, with $a_T,b_T$ as defined in \cref{TransportRays}, we have 
\begin{equation*}
    \int_T \varphi(x)d\mathcal{H}^1(x)=\|b_T-a_T\|\int_{[0,1]}\varphi(a_T+t(b_T-a_T))dt.
\end{equation*}
As the mapping $(a,b,t)\mapsto \|b-a\|\varphi(a+t(b-a))$ (where $(a,b,t)\in\mathbb{R}^d\times\mathbb{R}^d\times [0,1]$) is Borel measurable, by Fubini's theorem and noting the metric on $\mathcal{S}$ as given in \cref{TransportRays}, the mapping $T\mapsto \|b_T-a_T\|\int_{[0,1]}\varphi(a_T+t(b_T-a_T))dt$ (where $T\in\mathcal{S}$) is Borel measurable. 
\end{proof}

The following lemma gives the disintegration of the marginals $\mu$ and $\nu$ along transport rays, expressed in terms of $\lambda,\tilde{\mu}_T,\tilde{\nu}_T$ (cf.\ \Cref{Sect.1.1}).  

\begin{lemma}\label{L3nnn}
For any Borel set $A\subseteq\mathbb{R}^d$, the mappings $T\mapsto\tilde{\mu}_T(A)$ and $T\mapsto\tilde{\nu}_T(A)$ (where $T\in\mathcal{S}$) are Borel measurable, and 
\begin{equation*}
     \mu(A)=\int_{\mathcal{S}}\tilde{\mu}_T(A)d\lambda(T), \qquad \nu(A)=\int_{\mathcal{S}}\tilde{\nu}_T(A)d\lambda(T).
\end{equation*}
Moreover, $\lambda(\mathcal{S}\backslash\mathcal{S}^*)=0$, and for any $T\in\mathcal{S}^*$, $\tilde{\mu}_T$ and $\tilde{\nu}_T$ are Borel probability measures on $\mathbb{R}^d$. For any $T\in\mathcal{S}$, we have $\tilde{\mu}_T(\mathbb{R}^d\backslash \mathrm{int}(T))=\tilde{\nu}_T(\mathbb{R}^d\backslash \mathrm{int}(T))=0$.
\end{lemma}
\begin{proof}

Fix any Borel set $A\subseteq\mathbb{R}^d$. By \eqref{ress} and \cref{L3.15n}, the mappings $T\mapsto\tilde{\mu}_T(A)$ and $T\mapsto\tilde{\nu}_T(A)$ (where $T\in\mathcal{S}$) are Borel measurable. By \eqref{new1.1}, \eqref{disj}, and \cref{L2.3n}, we have
\begin{align}\label{L3.16ne.1}
    &\mu(A)=\int_A f(x)dx=\int_{\mathcal{T}_1^{*}} f(x)\mathbbm{1}_A(x)dx=\sum_{k=1}^K\sum_{k'=1}^{K'}\int_{\mathcal{T}^{*}_{1;k,k'}}f(x)\mathbbm{1}_A(x)dx\nonumber\\
    =\, &\sum_{k=1}^K\sum_{k'=1}^{K'}\int_{\bar{H}_{k,k'}}d\mathcal{H}^{d-1}(q)\int_{\mathbb{R}}\mathbbm{1}_{(-\beta(q),\alpha(q))}(t)f(q+tV(q))\mathbbm{1}_A(q+tV(q))R_{k,k'}(q,t)dt\nonumber\\
    =\, &\sum_{k=1}^K\sum_{k'=1}^{K'}\int_{\bar{H}_{k,k'}\cap\mathfrak{H}_{k,k'}}d\mathcal{H}^{d-1}(q)\int_{\mathbb{R}}\mathbbm{1}_{(-\beta(q),\alpha(q))}(t)f(q+tV(q))\mathbbm{1}_A(q+tV(q))R_{k,k'}(q,t)dt\nonumber\\
    =\,&\sum_{k=1}^K\sum_{k'=1}^{K'}\int_{\bar{H}_{k,k'}\cap\mathfrak{H}_{k,k'}}|\langle V(q), \mathtt{a}_{k,k'}\rangle|d\mathcal{H}^{d-1}(q)\nonumber\\
    &\hspace{1in}\int_{\mathbb{R}}\mathbbm{1}_{(-\beta(q),\alpha(q))}(t)f(q+tV(q))\mathbbm{1}_A(q+tV(q))\det(\mathbf{I}_d+tF(q))dt\nonumber\\
    =\,& \sum_{k=1}^K\sum_{k'=1}^{K'} \int_{\bar{H}_{k,k'}\cap\mathfrak{H}_{k,k'}} |\langle V(q), \mathtt{a}_{k,k'}\rangle|
    d\mathcal{H}^{d-1}(q)\int_{T(q)}f(x)\mathbbm{1}_A(x)\det\big(\mathbf{I}_d+\langle x-q, V(q)\rangle F(q)\big)d\mathcal{H}^1(x),
\end{align}
where we use \eqref{new.con} and \eqref{L3.10E1} in the third and fourth lines, respectively. Moreover, by \cref{L3.14n}, we have
\begin{align}\label{L3.16n1}
    & \int_{\mathcal{S}}\tilde{\mu}_T(A)d\lambda(T)\nonumber\\
=\,&\sum_{k=1}^K\sum_{k'=1}^{K'}\int_{\bar{H}_{k,k'}\cap\mathfrak{H}_{k,k'}}|\langle V(q), \mathtt{a}_{k,k'}\rangle|\bigg(\int_{T(q)}f(x)\det\big(\mathbf{I}_d+\langle x-q, V(q)\rangle F(q)\big)d\mathcal{H}^1(x)\bigg)\tilde{\mu}_{T(q)}(A)d\mathcal{H}^{d-1}(q).
\end{align}
Note that for any $k\in[K],k'\in[K']$ and $q\in \bar{H}_{k,k'}\cap\mathfrak{H}_{k,k'}$, by \cref{factors,De3.8}, we have
\begin{align*}
    \tilde{\mu}_{T(q)}(A)=\,&\int_{T(q)} f(x)\mathfrak{F}(x) \mathbbm{1}_A(x) d\mathcal{H}^1(x)\nonumber\\
    =\,& \mathbbm{1}_{\int_{T(q)}f(x)\det(\mathbf{I}_d+\langle x-q, V(q)\rangle F(q))d\mathcal{H}^1(x)>0}\cdot\frac{\int_{T(q)}f(x)\det(\mathbf{I}_d+\langle x-q, V(q)\rangle F(q))\mathbbm{1}_A(x)d\mathcal{H}^1(x)}{\int_{T(q)}f(x)\det(\mathbf{I}_d+\langle x-q, V(q)\rangle F(q))d\mathcal{H}^1(x)},
\end{align*}
hence
\begin{align}\label{L3.16n1.1}
 &\bigg(\int_{T(q)}f(x)\det\big(\mathbf{I}_d+\langle x-q, V(q)\rangle F(q)\big)d\mathcal{H}^1(x)\bigg)\tilde{\mu}_{T(q)}(A)\nonumber \\
 =\,& \int_{T(q)}f(x)\det\big(\mathbf{I}_d+\langle x-q, V(q)\rangle F(q)\big)\mathbbm{1}_A(x)d\mathcal{H}^1(x).
\end{align}
By \eqref{L3.16n1} and \eqref{L3.16n1.1}, we have
\begin{align}\label{L3.16ne.2}
&\int_{\mathcal{S}}\tilde{\mu}_T(A)d\lambda(T)\nonumber\\
=\,&\sum_{k=1}^K\sum_{k'=1}^{K'}\int_{\bar{H}_{k,k'}\cap\mathfrak{H}_{k,k'}}|\langle V(q), \mathtt{a}_{k,k'}\rangle|d\mathcal{H}^{d-1}(q)\int_{T(q)}f(x)\det\big(\mathbf{I}_d+\langle x-q, V(q)\rangle F(q)\big)\mathbbm{1}_A(x) d\mathcal{H}^1(x).
\end{align}
Combining \eqref{L3.16ne.1} and \eqref{L3.16ne.2}, we get $\mu(A)=\int_{\mathcal{S}}\tilde{\mu}_T(A)d\lambda(T)$. Similarly, we can deduce that $\nu(A)=\int_{\mathcal{S}}\tilde{\nu}_T(A)d\lambda(T)$.  

By \eqref{ress}, for every $T\in\mathcal{S}$, we have $\tilde{\mu}_T\big(\mathbb{R}^d\backslash \mathrm{int}(T)\big)=\tilde{\nu}_T\big(\mathbb{R}^d\backslash \mathrm{int}(T)\big)=0$. By \cref{L3n,Lem3.15}, $\mathcal{S}^*$ is a Borel subset of $\mathcal{S}$, and (recall \cref{factors,De3.8})
\begin{equation*}
   \lambda(\mathcal{S}\backslash\mathcal{S}^*)= \lambda(\mathcal{S}\backslash\{T(x):x\in\mathfrak{T}^{*}\})=\mu\big(\big\{z\in\mathcal{T}_1^{*}:T(z)\notin \{T(x):x\in\mathfrak{T}^{*}\}\big\}\big)\leq\mu(\mathcal{T}_1^*\backslash\mathfrak{T}^*)=0.
\end{equation*}
Now for any $T\in\mathcal{S}^*$, by \cref{factors,De3.8}, we have
\begin{equation*}
    \tilde{\mu}_T(\mathbb{R}^d)=\int_T \tilde{f}(x)d\mathcal{H}^1(x)=\int_T f(x)\mathfrak{F}(x)d\mathcal{H}^1(x)=\frac{\int_T f(x)\det(\mathbf{I}_d+\langle x-x_0, V_T\rangle F(x_0))d\mathcal{H}^1(x)}{\int_T f(z)\det(\mathbf{I}_d+\langle z-x_0, V_T\rangle F(x_0))d\mathcal{H}^1(z)}=1,
\end{equation*}
\begin{equation*}
    \tilde{\nu}_T(\mathbb{R}^d)=\int_T \tilde{g}(x)d\mathcal{H}^1(x)=\int_T g(x)\mathfrak{F}(x)d\mathcal{H}^1(x)=\frac{\int_T g(x)\det(\mathbf{I}_d+\langle x-x_0, V_T\rangle F(x_0))d\mathcal{H}^1(x)}{\int_T f(z)\det(\mathbf{I}_d+\langle z-x_0, V_T\rangle F(x_0))d\mathcal{H}^1(z)}=1,
\end{equation*}
where $x_0$ is an arbitrary point in $\mathrm{int}(T)$, and we use \cref{Ln2.2}(b) in the second display. Hence for any $T\in\mathcal{S}^*$, $\tilde{\mu}_T$ and $\tilde{\nu}_T$ are Borel probability measures on $\mathbb{R}^d$.
\end{proof}

The next lemma provides a disintegration of the Lebesgue measure restricted to $\mathfrak{T}^*$. 

\begin{lemma}\label{Lem3.18nn}
For any Borel measurable function $\varphi:\mathbb{R}^d\rightarrow [0,\infty)$, we have
\begin{equation*}
    \int_{\mathfrak{T}^*}\varphi(x)dx=\int_{\mathcal{S}}d\lambda(T)\int_{T}\varphi(x)\mathfrak{F}(x)d\mathcal{H}^1(x).
\end{equation*}
\end{lemma}
\begin{proof}

By \eqref{disj} and \cref{L2.3n}, noting \eqref{defT} and the fact that $\mathfrak{T}^*\subseteq\mathfrak{T}$, we have 
\begin{align*}
    &\int_{\mathfrak{T}^*}\varphi(x)dx=\sum_{k=1}^K\sum_{k'=1}^{K'}\int_{\mathfrak{T}^*\cap\mathcal{T}^*_{1;k,k'}}\varphi(x)dx\nonumber\\
    =\,&\sum_{k=1}^K\sum_{k'=1}^{K'}\int_{\bar{H}_{k,k'}\cap\mathfrak{H}_{k,k'}}d\mathcal{H}^{d-1}(q)\int_{\mathbb{R}}\mathbbm{1}_{(-\beta(q),\alpha(q))}(t)\varphi(q+tV(q))\mathbbm{1}_{\mathfrak{T}^*}(q+tV(q)) R_{k,k'}(q,t)dt\nonumber\\
    =\,&\sum_{k=1}^K\sum_{k'=1}^{K'}\int_{\bar{H}_{k,k'}\cap\mathfrak{H}_{k,k'}}|\langle V(q), \mathtt{a}_{k,k'}\rangle|d\mathcal{H}^{d-1}(q)\nonumber\\
    &\hspace{1in}\int_{\mathbb{R}}\mathbbm{1}_{(-\beta(q),\alpha(q))}(t)\varphi(q+tV(q))\mathbbm{1}_{\mathfrak{T}^*}(q+tV(q)) \det(\mathbf{I}_d+tF(q))dt\nonumber\\
    =\,&\sum_{k=1}^K\sum_{k'=1}^{K'}\int_{\bar{H}_{k,k'}\cap\mathfrak{H}_{k,k'}}|\langle V(q), \mathtt{a}_{k,k'}\rangle|\mathbbm{1}_{\int_{T(q)}f(x)\det(\mathbf{I}_d+\langle x-q, V(q)\rangle F(q))d\mathcal{H}^1(x)>0}\nonumber\\
    &\hspace{1.1in}\times \bigg(\int_{T(q)}\varphi(x)\det\big(\mathbf{I}_d+\langle x-q,V(q)\rangle F(q)\big)d\mathcal{H}^1(x)\bigg) d\mathcal{H}^{d-1}(q),
\end{align*}
where we use \eqref{L3.10E1} in the second to last equality. Using \cref{L3.14n}, the partition of the domain, and \cref{factors},
\begin{align*}
 &\int_{\mathcal{S}}d\lambda(T)\int_{T}\varphi(x)\mathfrak{F}(x)d\mathcal{H}^1(x)\nonumber\\
 =\,&\sum_{k=1}^K\sum_{k'=1}^{K'}\int_{\bar{H}_{k,k'}\cap\mathfrak{H}_{k,k'}}|\langle V(q), \mathtt{a}_{k,k'}\rangle|\mathbbm{1}_{\int_{T(q)}f(x)\det(\mathbf{I}_d+\langle x-q, V(q)\rangle F(q))d\mathcal{H}^1(x)>0}\nonumber\\
    &\hspace{1.1in}\times \bigg(\int_{T(q)}\varphi(x)\det\big(\mathbf{I}_d+\langle x-q,V(q)\rangle F(q)\big)d\mathcal{H}^1(x)\bigg) d\mathcal{H}^{d-1}(q),
\end{align*}
The conclusion of the lemma follows from the above two displays.
\end{proof}

We close this subsection with the following integrability result.

\begin{lemma}\label{L3.18}
We have 
\begin{equation*}
    \int_{\mathcal{S}}d\lambda(T)\int_{T}\tilde{f}(x)|\log(\tilde{f}(x))|d\mathcal{H}^1(x) \leq C, \qquad \int_{\mathcal{S}}d\lambda(T)\int_{T}\tilde{g}(x)|\log(\tilde{g}(x))|d\mathcal{H}^1(x) \leq C.
\end{equation*}
\end{lemma}
\begin{proof}

By \cref{factors,De3.8}, for any $k\in[K],k'\in[K']$, $q\in \bar{H}_{k,k'}\cap\mathfrak{H}_{k,k'}$, and $x\in\mathrm{int}(T(q))$, we have
\begin{equation*}
    \tilde{f}(x)=\mathbbm{1}_{\int_{T(q)}f(z)\det(\mathbf{I}_d+\langle z-q, V(q)\rangle F(q))d\mathcal{H}^1(z)>0}\frac{f(x)\det(\mathbf{I}_d+\langle x- q, V(q)\rangle F(q))}{\int_{T(q)}f(z)\det(\mathbf{I}_d+\langle z-q, V(q)\rangle F(q))d\mathcal{H}^1(z)}.
\end{equation*}
Hence if $\int_{T(q)}f(z)\det\big(\mathbf{I}_d+\langle z-q, V(q)\rangle F(q)\big)d\mathcal{H}^1(z)>0$, then 
\begin{align*}
  \int_{T(q)}\tilde{f}(x)|\log(\tilde{f}(x))|d\mathcal{H}^1(x)\leq\,&    \frac{\int_{T(q)}\big|\Phi\big(f(z)\det(\mathbf{I}_d+\langle z-q, V(q)\rangle F(q))\big)\big|d\mathcal{H}^1(z)}{\int_{T(q)}f(z)\det(\mathbf{I}_d+\langle z-q, V(q)\rangle F(q))d\mathcal{H}^1(z)}\nonumber\\
  &+\bigg|\log\bigg(\int_{T(q)}f(z)\det\big(\mathbf{I}_d+\langle z-q, V(q)\rangle F(q)\big)d\mathcal{H}^1(z)\bigg)\bigg|,
\end{align*}
where $\Phi(\cdot)$ is as in \eqref{Phi_def}. If $\int_{T(q)}f(z)\det\big(\mathbf{I}_d+\langle z-q, V(q)\rangle F(q)\big)d\mathcal{H}^1(z)=0$, then
\begin{equation*}
    \int_{T(q)}\tilde{f}(x)|\log(\tilde{f}(x))|d\mathcal{H}^1(x)=0.
\end{equation*}
Therefore, by \cref{L3.14n}, 
\begin{align}\label{L3.18nn}
    &\int_{\mathcal{S}}d\lambda(T)\int_{T}\tilde{f}(x)|\log(\tilde{f}(x))|d\mathcal{H}^1(x)\nonumber\\
    =\,& \sum_{k=1}^K\sum_{k'=1}^{K'}\int_{\bar{H}_{k,k'}\cap\mathfrak{H}_{k,k'}}|\langle V(q), \mathtt{a}_{k,k'}\rangle|\bigg(\int_{T(q)}f(x)\det\big(\mathbf{I}_d+\langle x-q, V(q)\rangle F(q)\big)d\mathcal{H}^1(x)\bigg)\nonumber\\
    &\hspace{1.2in}\times \bigg( \int_{T(q)}\tilde{f}(x)|\log(\tilde{f}(x))|d\mathcal{H}^1(x)\bigg)d\mathcal{H}^{d-1}(q)\nonumber\\
    \leq\,& \sum_{k=1}^K\sum_{k'=1}^{K'}\int_{\bar{H}_{k,k'}\cap\mathfrak{H}_{k,k'}}\bigg(\int_{T(q)}\big|\Phi\big(f(x)\det\big(\mathbf{I}_d+\langle x-q, V(q)\rangle F(q)\big)\big)\big|d\mathcal{H}^1(x)\bigg)d\mathcal{H}^{d-1}(q)\nonumber\\
    &+\sum_{k=1}^K\sum_{k'=1}^{K'}\int_{\bar{H}_{k,k'}\cap\mathfrak{H}_{k,k'}}\bigg|\Phi\bigg(\int_{T(q)}f(x)\det\big(\mathbf{I}_d+\langle x-q, V(q)\rangle F(q)\big)d\mathcal{H}^1(x)\bigg)\bigg|d\mathcal{H}^{d-1}(q),
\end{align}
where we use $|\langle V(q), \mathtt{a}_{k,k'}\rangle|\leq \|V(q)\|\|\mathtt{a}_{k,k'}\|=1$ in the fourth line. 

For any $k\in [K],k'\in[K']$ and $q\in \bar{H}_{k,k'}\cap\mathfrak{H}_{k,k'}$, by \eqref{distcondi}, we have $\min\{\alpha(q),\beta(q)\}\geq 2d_0$. Hence by \cref{L2.4} (see \eqref{L2.4.enew3}), we have $\|F(q)\|_2\leq \frac{4}{\min\{\alpha(q),\beta(q)\}}\leq 2d_0^{-1}\leq C$. Using Hadamard's inequality (see, e.g., \cite[Lemma 2.5]{ipsen2008perturbation}), for any $x\in\mathrm{int}(T(q))$, we obtain
\begin{align*}
    \big|\det\big(\mathbf{I}_d+\langle x-q, V(q)\rangle F(q)\big)\big| &\leq \big\|\mathbf{I}_d+\langle x-q, V(q)\rangle F(q)\big\|_2^d\leq \big(1+\|x-q\|\|F(q)\|_2\big)^d\nonumber\\
    &\leq (1+2D\cdot C)^d\leq C,
\end{align*}
where we use \eqref{bddsalphabeta} in the second line. Hence noting \eqref{bddLip} and \eqref{bddsalphabeta}, we get
\begin{equation*}
0\leq f(x)\det\big(\mathbf{I}_d+\langle x-q, V(q)\rangle F(q)\big)\leq C,\quad\text{for all }x\in\mathrm{int}(T(q)),
\end{equation*}
\begin{equation*}
    0\leq \int_{T(q)}f(x)\det\big(\mathbf{I}_d+\langle x-q, V(q)\rangle F(q)\big) d\mathcal{H}^1(x) \leq C.
\end{equation*}
Therefore, we have
\begin{equation*}
    \int_{T(q)}\big|\Phi\big(f(x)\det\big(\mathbf{I}_d+\langle x-q, V(q)\rangle F(q)\big)\big)\big|d\mathcal{H}^1(x)\leq C,
\end{equation*}
\begin{equation*}
    \bigg|\Phi\bigg(\int_{T(q)}f(x)\det\big(\mathbf{I}_d+\langle x-q, V(q)\rangle F(q)\big)d\mathcal{H}^1(x)\bigg)\bigg|\leq C.
\end{equation*}
Combining these with \eqref{L3.18nn}, we obtain that
\begin{equation*}
    \int_{\mathcal{S}}d\lambda(T)\int_{T}\tilde{f}(x)|\log(\tilde{f}(x))|d\mathcal{H}^1(x)\leq C\sum_{k=1}^K\sum_{k'=1}^{K'}\mathcal{H}^{d-1}(\bar{H}_{k,k'}\cap\mathfrak{H}_{k,k'})\leq C,
\end{equation*}
where we use \eqref{BddHtildekk} in the second inequality. Similarly, noting \eqref{qrelations}, we deduce the claim for~$\tilde{g}$.
\end{proof}

\subsection{Proof of \cref{th:lowerBound}}\label{Sec3:3}

In this subsection, we complete the proof of \cref{th:lowerBound}. Throughout, we fix $\eps_j\rightarrow 0^+$ and arbitrary couplings  $\gamma'_{\eps_j}\in\Pi(\mu,\nu)$ that converge weakly to an optimal transport plan $\gamma_0'\in\Pi(\mu,\nu)$. Let $\gamma_0'=\lambda\otimes\kappa'_T$ be the disintegration of $\gamma_0'$ along transport rays. Note that, for each $j\in  \mathbb{N}^*$, if $\gamma'_{\eps_j}$ is not absolutely continuous with respect to $\mathcal{L}^d\otimes \mathcal{L}^d$, then $\mathcal{C}_{\eps_j}(\gamma'_{\eps_j})=\infty$. Hence, without loss of generality, we assume hereafter that $\gamma'_{\eps_j}$ is absolutely continuous with respect to $\mathcal{L}^d\otimes \mathcal{L}^d$ for all $j\in\mathbb{N}^*$. 

For each $j\in\mathbb{N}^*$, we denote by $\phi_{\eps_j}(\cdot,\cdot)$ the density of $\gamma'_{\eps_j}$ with respect to $\mathcal{L}^d\otimes\mathcal{L}^d$, and define 
\begin{equation}\label{P3.1.eq1}
    \psi_{\varepsilon_j}(x,t):=\begin{cases}
  \int_{O(V(x))}\phi_{\varepsilon_j}(x,x+tV(x)+w)d\mathcal{H}^{d-1}(w) , & \text{for }x\in\mathcal{T}_1^{*},t\in\mathbb{R},\\
     0,   & \text{for }x\in\mathbb{R}^d\backslash\mathcal{T}_1^{*},t\in\mathbb{R}.
    \end{cases}
\end{equation}
Since $\gamma_{\eps_j}'\in\Pi(\mu,\nu)$, we may assume without loss of generality that
\begin{equation}\label{phij_0}
    \phi_{\eps_j}(x,y)=0,\quad \text{ for any }(x,y)\in (\mathbb{R}^d\times\mathbb{R}^d)\backslash (\mathcal{X}\times \mathcal{Y}).
\end{equation}
Note that by \eqref{new1.1}, for any Lebesgue measurable set $A \subseteq \mathbb{R}^d$, we have
\begin{align}\label{L3.11.e1}
\int_{\mathbb{R}^d\times\mathbb{R}}\psi_{\eps_j}(x,t)\mathbbm{1}_{A}(x)dxdt
=\,& \int_{A\cap\mathcal{T}_1^{*}}dx\int_{\mathbb{R}}dt \int_{O(V(x))}\phi_{\varepsilon_j}(x,x+tV(x)+w)d\mathcal{H}^{d-1}(w)\nonumber\\
=\,& \int_{A\cap\mathcal{T}_1^{*}}dx \int_{\mathbb{R}^d}\phi_{\eps_j}(x,y)dy=\int_{A\cap\mathcal{T}_1^{*}}f(x)dx=\mu(A). 
\end{align}

We begin with the following elementary lemma.

\begin{lemma}\label{Lemma3.3}
Assume that $x\in \mathcal{X}\cap\mathcal{T}_1^{*}$, $z\in L(x)$, $w\in O(V(x))$, and $z+w\in\mathcal{Y}$. Let $a$ be the upper end of $T(x)$. Then we have $\max\{\|x-z\|, \|z-a\|, \|w\|\} \leq 2D$.
\end{lemma}
\begin{proof}
As $a, x,z+w\in\mathcal{X}\cup\mathcal{Y}\subseteq B_d(0,D)$, we have
\begin{equation*}
    \|z-a\|\leq \sqrt{\|w\|^2+\|z-a\|^2} =\|(z+w)-a\|\leq 2D,
\end{equation*}
\begin{equation*}
    \max\{\|x-z\|, \|w\|\}\leq\sqrt{\|w\|^2+\|x-z\|^2}=\|(z+w)-x\|\leq 2D. \qedhere
\end{equation*}
\end{proof}
\begin{remark}
By \eqref{P3.1.eq1}--\eqref{phij_0} and \Cref{Lemma3.3}, and noting \eqref{domainconst}, we have
\begin{equation}\label{compact_supports}
\{(x,t)\in\mathbb{R}^{d+1}:\psi_{\eps_j}(x,t)>0\}\subseteq (\mathcal{X}\cap\mathcal{T}_1^*)\times [-2D,2D]\subseteq  B_d(0,D)\times [-2D,2D].
\end{equation}
\end{remark}

The following result was proved in \cite{JacodMemin.81} (see also \cite[Lemma A.3]{lacker2015mean} for a reference in English).

\begin{lemma}\label{L3.24}
Suppose that $X, Y$ are complete separable metric spaces. Let $\varphi:X\times Y\rightarrow \mathbb{R}$ be a bounded function such that
\begin{itemize}
    \item[(a)] $\varphi(\cdot,y)$ is Borel measurable for each $y\in Y$,
    \item[(b)] $\varphi(x,\cdot)$ is continuous for each $x\in X$.
\end{itemize}
Then, for any Borel probability measures $(P_n)_{n\in\mathbb{N}},P$ on $X\times Y$ such that $P_n(\cdot\times Y)=P(\cdot\times Y)$ for every $n\in\mathbb{N}$ and $P_n$ converges weakly to $P$ as $n\rightarrow\infty$, we have $\int_{X\times Y}\varphi dP_n\rightarrow \int_{X\times Y}\varphi dP$. 
\end{lemma}

The following lemma bounds the entropy of the marginal density $\psi_{\varepsilon_j}(x,t)$. 

\begin{lemma}\label{L3.6}
We have
\begin{align*}
&\liminf_{j\rightarrow\infty}\int_{(\mathcal{X}\cap\mathcal{T}_1^*)\times\mathbb{R}}\psi_{\varepsilon_j}(x,t)\log(\psi_{\varepsilon_j}(x,t))dxdt\nonumber\\
\geq\,& \int_{\mathcal{S}}H(\kappa'_T|\tilde{\mu}_T\otimes\tilde{\nu}_T)d\lambda(T) + \int_{\mathcal{X}}f(x)\log(f(x))dx+\int_{\mathcal{S}}d\lambda(T)\int_T \tilde{g}(y)\log(\tilde{g}(y))d\mathcal{H}^1(y).
\end{align*}
\end{lemma}
\begin{proof}

\noindent\textbf{Step 1.}
Denote by $C(\mathbb{R}^{d+1})$ the set of continuous functions on $\mathbb{R}^{d+1}$ and by $C_b(\mathbb{R}^{d+1})$ the set of continuous bounded functions on $\mathbb{R}^{d+1}$. Following the proof of \cite[Lemma 1.3]{Nutz.20} and noting \eqref{compact_supports}, we obtain that for each $j\in\mathbb{N}^*$,
\begin{align}\label{variation}
    &\int_{(\mathcal{X}\cap\mathcal{T}_1^*) \times\mathbb{R}}\psi_{\varepsilon_j}(x,t)\log(\psi_{\varepsilon_j}(x,t))dx dt  \nonumber\\
    =\,& \sup_{\substack{\varphi\in C(\mathbb{R}^{d+1}):\\\int_{\mathbb{R}^{d+1}}e^{\varphi(x,t)}dxdt<\infty}}\bigg\{\int_{(\mathcal{X}\cap\mathcal{T}_1^*)\times\mathbb{R}}\varphi(x,t)\psi_{\varepsilon_j}(x,t)dxdt-\log\bigg(\int_{\mathbb{R}^{d+1}}e^{\varphi(x,t)}dxdt\bigg)\bigg\}.
\end{align} 
For any $\varphi\in C(\mathbb{R}^{d+1})$ such that $\int_{\mathbb{R}^{d+1}}e^{\varphi(x,t)}dxdt<\infty$, there exists $\tilde{\varphi}\in C_b(\mathbb{R}^{d+1})$ such that $\tilde{\varphi}|_{B_d(0,D)\times [-2D,2D]}=\varphi|_{B_d(0,D)\times [-2D,2D]}$. Using \eqref{new1.1} and \eqref{compact_supports}, we obtain 
\begin{align*}
 &\int_{(\mathcal{X}\cap\mathcal{T}_1^*)\times\mathbb{R}}\varphi(x,t)\psi_{\varepsilon_j}(x,t)dxdt=\int_{(\mathcal{X}\cap\mathcal{T}_1^*)\times\mathbb{R}}\tilde{\varphi}(x,t)\psi_{\varepsilon_j}(x,t)dxdt\nonumber\\
    =\,& \int_{(\mathcal{X}\cap\mathcal{T}_1^{*})\times\mathbb{R}} \tilde{\varphi}(x,t) dx dt \int_{O(V(x))} \phi_{\varepsilon_j}(x,x+tV(x)+w)d\mathcal{H}^{d-1}(w)\nonumber\\
    =\,& \int_{(\mathcal{X}\cap\mathcal{T}_1^*)\times\mathbb{R}^d} \tilde{\varphi}(x,\langle y-x, V(x)\rangle)\phi_{\varepsilon_j}(x,y)dxdy.
\end{align*}
For any $x,y\in\mathbb{R}^d$, define $\Psi(x,y):=\tilde{\varphi}(x,\langle y-x, V(x)\rangle)\mathbbm{1}_{\mathcal{X}\cap\mathcal{T}_1^{*}}(x)$. Note that $\Psi(\cdot,\cdot)$ is bounded and Borel measurable, and for each $x\in\mathbb{R}^d$, $\Psi(x,\cdot)$ is continuous. Hence by \cref{L3.24} (recall that $\phi_{\eps_j}$ is the density of $\gamma'_{\eps_j}\in\Pi(\mu,\nu)$) and the fact that $\mu(\mathcal{X}\cap\mathcal{T}_1^*)=1$ (see \eqref{new1.1}), 
\begin{align*}
&\lim_{j\rightarrow\infty}\int_{(\mathcal{X}\cap\mathcal{T}_1^*)\times\mathbb{R}}\varphi(x,t)\psi_{\eps_j}(x,t) dx dt
= \lim_{j\rightarrow\infty}\int_{(\mathcal{X}\cap\mathcal{T}_1^{*})\times\mathbb{R}^d}\tilde{\varphi}(x,\langle y-x, V(x)\rangle)\phi_{\eps_j}(x,y)dxdy\nonumber\\
=\,&\lim_{j\rightarrow\infty}\int_{\mathbb{R}^d\times\mathbb{R}^d} \Psi(x,y) \phi_{\eps_j}(x,y) dxdy=\int_{\mathbb{R}^d\times\mathbb{R}^d} \Psi(x,y)d\gamma_0'(x,y)\nonumber\\
=\,&\int_{(\mathcal{X}\cap\mathcal{T}_1^*)\times\mathbb{R}^d} \tilde{\varphi}(x,\langle y-x, V(x)\rangle) d\gamma_0'(x,y)=\int_{\mathbb{R}^d\times\mathbb{R}^d} \varphi(x,\langle y-x, V(x)\rangle) d\gamma_0'(x,y).
\end{align*}
Hence by \eqref{variation},  
\begin{align}\label{Eqnnew}
 & \liminf_{j\rightarrow\infty}\int_{(\mathcal{X}\cap\mathcal{T}_1^*)\times\mathbb{R}}\psi_{\varepsilon_j}(x,t)\log(\psi_{\varepsilon_j}(x,t))dx dt\nonumber\\
 \geq\,& \sup_{\substack{\varphi\in C(\mathbb{R}^{d+1}):\\\int_{\mathbb{R}^{d+1}}e^{\varphi(x,t)}dxdt<\infty}}\bigg\{\int_{\mathbb{R}^d\times\mathbb{R}^d} \varphi(x,\langle y-x, V(x)\rangle)d\gamma_0'(x,y)-\log\bigg(\int_{\mathbb{R}^{d+1}}e^{\varphi(x,t)}dxdt\bigg)\bigg\}.
\end{align}

\noindent\textbf{Step 2.} For any $(x,y)\in\mathbb{R}^d\times\mathbb{R}^d$, we define $\mathfrak{r}(x,y):=(x,\langle y-x, V(x)\rangle)\in\mathbb{R}^{d}\times\mathbb{R}$. Define $\EuScript{Q}$ to be the pushforward of $\gamma_0'$ by $\mathfrak{r}$. Note that $\EuScript{Q}$ is a Borel probability measure on $\mathbb{R}^{d}\times\mathbb{R}$ whose first marginal is given by $\mu$. By \cref{Lem3.15}, $\EuScript{Q}$ is concentrated on $\mathfrak{T}^*\times\mathbb{R}$. For any $\varphi\in C(\mathbb{R}^{d+1})$ such that $\int_{\mathbb{R}^{d+1}}e^{\varphi(x,t)}dxdt<\infty$, we have
\begin{equation*}
    \int_{\mathbb{R}^{d+1}}\varphi d\EuScript{Q}=\int_{\mathbb{R}^d\times\mathbb{R}^d}\varphi(x,\langle y-x,V(x)\rangle)d\gamma_0'(x,y).
\end{equation*}
Arguing similarly as in the proof of \cite[Lemma 1.3]{Nutz.20}, we can deduce that if $\EuScript{Q}\not\ll \mathcal{L}^{d+1}$, then  
\begin{equation*}
    \sup_{\substack{\varphi\in C(\mathbb{R}^{d+1}):\\\int_{\mathbb{R}^{d+1}}e^{\varphi(x,t)}dxdt<\infty}}\bigg\{\int_{\mathbb{R}^{d}\times\mathbb{R}^{d}}\varphi(x,\langle y-x,V(x)\rangle)d\gamma_0'(x,y)-\log\bigg(\int_{\mathbb{R}^{d+1}}e^{\varphi(x,t)}dxdt\bigg)\bigg\}=\infty.
\end{equation*}
Combined with \eqref{Eqnnew}, this completes the proof of the lemma when $\EuScript{Q}\not\ll \mathcal{L}^{d+1}$.

Throughout the rest of the proof, we assume that $\EuScript{Q}\ll \mathcal{L}^{d+1}$, and denote by $\mathfrak{h}(\cdot,\cdot)$ the density of $\EuScript{Q}$ with respect to $\mathcal{L}^{d+1}$. Arguing similarly as in the proof of \cite[Lemma 1.3]{Nutz.20} and noting that $\EuScript{Q}$ is concentrated on $\mathfrak{T}^*\times\mathbb{R}$, we get
\begin{align}\label{Eq3.60}
   & \sup_{\substack{\varphi\in C(\mathbb{R}^{d+1}):\\\int_{\mathbb{R}^{d+1}}e^{\varphi(x,t)}dxdt<\infty}}\bigg\{\int_{\mathbb{R}^{d}\times\mathbb{R}^{d}}\varphi(x,\langle y-x,V(x)\rangle)d\gamma_0'(x,y)-\log\bigg(\int_{\mathbb{R}^{d+1}}e^{\varphi(x,t)}dxdt\bigg)\bigg\}\nonumber\\
   &=\int_{\mathbb{R}^{d+1}}\mathfrak{h}(x,t)\log(\mathfrak{h}(x,t))dxdt=\int_{\mathfrak{T}^*\times\mathbb{R}}\mathfrak{h}(x,t)\log(\mathfrak{h}(x,t))dxdt.
\end{align}

Now note that for any Borel set $A\subseteq\mathbb{R}^{d+1}$,   
\begin{align*}
   \EuScript{Q}(A)=\,&\EuScript{Q}(A\cap (\mathfrak{T}^*\times\mathbb{R}))=\int_{\mathfrak{T}^*}dx\int_{\mathbb{R}}\mathbbm{1}_{(x,t)\in A} \mathfrak{h}(x,t)dt
 \nonumber\\
    =\,& \int_{\mathcal{S}}d\lambda(T)\int_{T}\mathfrak{F}(x)d\mathcal{H}^1(x)\int_{\mathbb{R}}\mathbbm{1}_{(x,t)\in A} \mathfrak{h}(x,t)dt=\int_{\mathcal{S}}d\lambda(T)\int_{T\times\mathbb{R}}\mathbbm{1}_{(x,t)\in A}\mathfrak{F}(x)\mathfrak{h}(x,t)d\mathcal{H}^1(x)dt,
\end{align*}
where we use \cref{Lem3.18nn} in the third equality. As $\gamma_0'=\lambda\otimes\kappa_T'$, we have
\begin{align*}
    \EuScript{Q}(A)=\,& \gamma_0'(\mathfrak{r}^{-1}(A))=
    \gamma_0'(\{(x,y)\in\mathbb{R}^d\times\mathbb{R}^d:(x,\langle y-x,V(x)\rangle)\in A\})\nonumber\\
    =\,&\int_{\mathcal{S}}\kappa'_T(\{(x,y)\in T\times T:(x,\langle y-x, V(x) \rangle)\in A\})d\lambda(T).
\end{align*}
Using the above two displays and
arguing similarly as in the proof of \cite[Theorem 9.2]{MR2011032}, we get for $\lambda$-a.e.\ $T\in\mathcal{S}$, 
\begin{equation}\label{Enq3.61}
    \int_{T\times\mathbb{R}}\mathbbm{1}_{(x,t)\in A}\mathfrak{F}(x)\mathfrak{h}(x,t)d\mathcal{H}^1(x)dt=\kappa'_T(\{(x,y)\in T\times T:(x,\langle V(x),y-x\rangle)  \in  A\})
\end{equation}
for every Borel set $A\subseteq\mathbb{R}^{d+1}$. 

Now consider any $T\in\mathcal{S}$ such that \eqref{Enq3.61} holds and any Borel set $B  \subseteq T\times T$. Take $A=\{(x,t)\in T\times \mathbb{R}:(x,x+tV(x))\in B\}$, and note that $A$ is a Borel subset of $\mathbb{R}^{d+1}$. We have
\begin{equation}\label{Enq.3.1.1}
    \kappa_T'(B)=\kappa'_T(\{(x,y)\in T\times T:(x,\langle V(x),y-x\rangle)  \in  A\})
\end{equation}
and
\begin{align}\label{Enq.3.1.2}
    \int_{T\times\mathbb{R}}\mathbbm{1}_{(x,t)\in  A}  \mathfrak{F}(x)  \mathfrak{h}(x,t)d\mathcal{H}^1(x)dt=\,& \int_{T\times\mathbb{R}}\mathbbm{1}_{(x,x+tV(x))\in  B} \mathfrak{F}(x) \mathfrak{h}(x,t)d\mathcal{H}^1(x)dt\nonumber\\
    =\,&\int_{B} \mathfrak{F}(x)\mathfrak{h}(x,\langle y-x,V(x)\rangle)d\mathcal{H}^1(x)d\mathcal{H}^1(y).
\end{align}
By \eqref{Enq3.61}--\eqref{Enq.3.1.2}, for $\lambda$-a.e.\ $T\in\mathcal{S}$ and any Borel set $B\subseteq T\times T$, we have
\begin{equation*}
    \kappa_T'(B)=\int_{B} \mathfrak{F}(x)\mathfrak{h}(x,\langle y-x,V(x)\rangle)d\mathcal{H}^1(x)d\mathcal{H}^1(y).
\end{equation*}
Hence for $\lambda$-a.e.\ $T\in \mathcal{S}$, $\kappa_T'\ll \mathcal{H}^1|_T\otimes \mathcal{H}^1|_T$, and the density of $\kappa_T'$ with respect to $\mathcal{H}^1|_T\otimes \mathcal{H}^1|_T$ is given by $\mathfrak{F}(x)\mathfrak{h}(x,\langle y-x,V(x)\rangle)$ for $(x,y)\in T\times T$. Further, for $\lambda$-a.e.\ $T\in\mathcal{S}$, as $\kappa_T'\in \Pi(\tilde{\mu}_T,\tilde{\nu}_T)$, for $\mathcal{H}^1\otimes\mathcal{H}^1$-a.e.\ $(x,y)\in  T\times T$, if $\tilde{f}(x)=0$ or $\tilde{g}(y)=0$, then $\mathfrak{F}(x)\mathfrak{h}(x,\langle y-x,V(x)\rangle)= 0$. Consequently, 
\begin{align*}
&\int_{\mathcal{S}}H(\kappa_T'|\tilde{\mu}_T\otimes\tilde{\nu}_T)d\lambda(T)=\int_{\mathcal{S}}d\lambda(T)\int_{T\times T}\log\bigg(\frac{\mathfrak{F}(x)\mathfrak{h}(x,\langle y-x,V(x)\rangle)}{\tilde{f}(x)\tilde{g}(y)}\bigg) d\kappa_T'(x,y)\nonumber\\
    =\,& \int_{\mathcal{S}}d\lambda(T)\int_{T\times T}\log(\mathfrak{h}(x,\langle y-x,V(x)\rangle)) d\kappa_T'(x,y)-\int_{\mathcal{S}}d\lambda(T)\int_{T\times T}\log(f(x))d\kappa_T'(x,y)\nonumber\\
    &-\int_{\mathcal{S}}d\lambda(T)\int_{T\times T}\log(\tilde{g}(y))d\kappa_T'(x,y)\nonumber\\
    =\,& \int_{\mathcal{S}}d\lambda(T)\int_{T\times T}\log(\mathfrak{h}(x,\langle y-x,V(x)\rangle))\mathfrak{F}(x)\mathfrak{h}(x,\langle y-x,V(x)\rangle)d\mathcal{H}^1(x)d\mathcal{H}^1(y)\nonumber\\
    & -\int_{\mathbb{R}^d\times\mathbb{R}^d}\log(f(x))d\gamma_0'(x,y)-\int_{\mathcal{S}}d\lambda(T)\int_{T}\tilde{g}(y)\log(\tilde{g}(y))d\mathcal{H}^1(y) \nonumber\\
    =\,& \int_{\mathcal{S}}d\lambda(T)\int_{T\times\mathbb{R}} \mathfrak{F}(x) \mathfrak{h}(x,t) \log(\mathfrak{h}(x,t))d\mathcal{H}^1(x)dt\nonumber\\
    &-\int_{\mathcal{X}}f(x)\log(f(x))dx-\int_{\mathcal{S}}d\lambda(T)\int_{T}\tilde{g}(y)\log(\tilde{g}(y))d\mathcal{H}^1(y),
\end{align*}
where we use the fact that $\gamma_0'=\lambda\otimes\kappa_T'$ in the third equality. Note that by Lemma \ref{Lem3.18nn},
\begin{equation*}
    \int_{\mathcal{S}}d\lambda(T)\int_{T\times\mathbb{R}} \mathfrak{F}(x) \mathfrak{h}(x,t) \log(\mathfrak{h}(x,t))d\mathcal{H}^1(x)dt=\int_{\mathfrak{T}^*\times\mathbb{R}}\mathfrak{h}(x,t)\log(\mathfrak{h}(x,t))dxdt.
\end{equation*}
Combining the above two displays, we obtain
\begin{align}\label{Eqnnnnew}
    \int_{\mathfrak{T}^*\times\mathbb{R}}\mathfrak{h}(x,t)\log(\mathfrak{h}(x,t))dxdt=\,&\int_{\mathcal{S}}H(\kappa_T'|\tilde{\mu}_T\otimes\tilde{\nu}_T)d\lambda(T)+\int_{\mathcal{X}}f(x)\log(f(x))dx\nonumber\\
    &+  \int_{\mathcal{S}}d\lambda(T)\int_{T}\tilde{g}(y)\log(\tilde{g}(y))d\mathcal{H}^1(y).
\end{align}
The conclusion of the lemma follows from \eqref{Eqnnew}, \eqref{Eq3.60}, and \eqref{Eqnnnnew}. 
\end{proof}

The next lemma will be used to bound the first term of the decomposition~\eqref{eq:overview-G-split} discussed in the proof overview.

\begin{lemma}\label{L3.7}
For any $M\geq 0$, $\delta\in (0,1)$, and Borel measurable set $\mathcal{X}_0\subseteq \{x\in\mathbb{R}^d:\alpha(x)\geq \delta\}$, 
\begin{equation*}
    \limsup_{j\rightarrow\infty} \int_{\mathcal{X}_0\times\mathbb{R}} \big(G(x,x+tV(x))\vee (-M)\big) \psi_{\eps_j}(x,t)dxdt\leq \int_{\mathcal{X}_0\times\mathbb{R}^d}G(x,y)\vee (-M)d\gamma_0'(x,y).
\end{equation*}
\end{lemma}
\begin{proof}

Let $\mathscr{X}\subseteq \mathcal{T}_1^{*}$ be the set of $x\in \mathcal{T}_1^{*}$ for which the conclusion of \cref{P-lb-2} holds. Note that $\mathcal{L}^d(\mathcal{T}_1^{*}\backslash \mathscr{X})=0$. For any $x\in\mathbb{R}^d$ and $t\in\mathbb{R}$, we define
\begin{equation*}
    \tilde{G}_M(x,t):=\begin{cases}
        G(x,x+tV(x))\vee (-M) & \text{ if }x\in\mathscr{X},t>-\beta(x),\\
        \Big(\lim\limits_{s\rightarrow (-\beta(x))^{+}}G(x,x+sV(x))\Big)\vee (-M) & \text{ if }x\in\mathscr{X}, t\leq -\beta(x),\\
        0 & \text{ if }x\in\mathbb{R}^d\backslash\mathscr{X}.
    \end{cases} 
\end{equation*}
Note that by \cref{P-lb-2}(a) and the fact that $\mathcal{X}_0\subseteq \{x\in\mathbb{R}^d:\alpha(x)\geq \delta\}$, for any $x\in\mathcal{X}_0$ and $t\in\mathbb{R}$, we have $-M\leq \tilde{G}_M(x,t)\leq \frac{d-1}{2}\log(\delta^{-1})+C$, where $C>0$ only depends on $D,d$. By \cref{P-lb-2}(c), for any $x\in\mathbb{R}^d$, $\tilde{G}_M(x,\cdot)$ is continuous. Consequently, for any $x\in\mathbb{R}^d$, the function $\mathbb{R}^d\ni y \mapsto\mathbbm{1}_{\mathcal{X}_0\cap\mathcal{X}\cap\mathcal{T}_1^*\cap\mathscr{X}}(x)\tilde{G}_M(x,\langle y-x,V(x)\rangle)$ is continuous. By \eqref{compact_supports}, for any $j\in\mathbb{N}^*$, 
\begin{align*}
     &  \int_{\mathcal{X}_0\times\mathbb{R}} \big(G(x,x+tV(x))\vee (-M)\big) \psi_{\eps_j}(x,t)dxdt\nonumber\\
   =\,& \int_{(\mathcal{X}_0\cap\mathcal{X}\cap\mathcal{T}_1^{*})\times\mathbb{R}} \big(G(x,x+tV(x))\vee (-M)\big) \psi_{\eps_j}(x,t)dxdt\nonumber\\
    =\,& \int_{(\mathcal{X}_0\cap\mathcal{X}\cap\mathcal{T}_1^{*}\cap\mathscr{X})\times\mathbb{R}} \big(G(x,x+tV(x))\vee (-M)\big) \psi_{\eps_j}(x,t)dxdt.
\end{align*}
Moreover, for any $x\in\mathscr{X}$ and $t\leq -\beta(x)$, by definition, $G(x,x+tV(x))=-\infty$, hence $G(x,x+tV(x))\vee (-M)=-M\leq \tilde{G}_M(x,t)$. Hence by \cref{L3.24}, we have
\begin{align*}
   & \limsup_{j\rightarrow\infty} \int_{\mathcal{X}_0\times\mathbb{R}} \big(G(x,x+tV(x))\vee (-M)\big) \psi_{\eps_j}(x,t)dxdt\nonumber\\
    =\,& \limsup_{j\rightarrow\infty} \int_{(\mathcal{X}_0\cap\mathcal{X}\cap\mathcal{T}_1^{*}\cap\mathscr{X})\times\mathbb{R}} \big(G(x,x+tV(x))\vee (-M)\big) \psi_{\eps_j}(x,t)dxdt\nonumber\\
    \leq\, & \limsup_{j\rightarrow\infty} \int_{(\mathcal{X}_0\cap\mathcal{X}\cap\mathcal{T}_1^{*}\cap\mathscr{X})\times\mathbb{R}} \tilde{G}_M(x,t) \psi_{\eps_j}(x,t)dxdt\nonumber\\
=\, &\limsup_{j\rightarrow\infty} \int_{(\mathcal{X}_0\cap\mathcal{X}\cap\mathcal{T}_1^{*}\cap\mathscr{X})\times\mathbb{R}^d} \tilde{G}_M(x,\langle y-x,V(x)\rangle)\phi_{\eps_j}(x,y)dxdy\nonumber\\
   =\,& \int_{(\mathcal{X}_0\cap\mathcal{X}\cap\mathcal{T}_1^{*}\cap\mathscr{X})\times\mathbb{R}^d} \tilde{G}_M(x,\langle y-x,V(x)\rangle)d\gamma_0'(x,y)\nonumber\\
   =\,& \int_{(\mathcal{X}_0\cap\mathcal{X}\cap\mathcal{T}_1^{*}\cap\mathscr{X})\times\mathbb{R}^d} \big(G(x,y)\vee (-M)\big) d\gamma_0'(x,y)
   = \int_{\mathcal{X}_0\times\mathbb{R}^d} \big(G(x,y)\vee (-M)\big) d\gamma_0'(x,y),
\end{align*}
where the last equality uses \eqref{new1.1} and the fact that the first marginal of $\gamma_0'$ is $\mu\ll\mathcal{L}^d$. 
\end{proof}

The next lemma bounds and computes certain integrals that will appear in the main proof.
\begin{lemma}\label{Lem3.11}
We have
\begin{equation*}
    \int_{\mathcal{S}}d\lambda(T)\int_{T\times T}\big|\log\det\big(\mathbf{I}_d+\|x-y\|F(y)\big)\big|d\kappa_T'(x,y)\leq C,
\end{equation*}
\begin{equation*}
    \int_{\mathcal{S}^*}d\lambda(T)\int_{T} |\log(\mathfrak{F}(x))|\tilde{f}(x)d\mathcal{H}^1(x)\leq C,
\end{equation*}
as well as the identities
\begin{align*}
    &\int_{\mathcal{S}}d\lambda(T)\int_{T\times T}\log\det\big(\mathbf{I}_d+\|x-y\|F(y)\big)d\kappa_T'(x,y)\nonumber\\
    =\,& -\int_{\mathcal{X}}f(x)\log(f(x))dx+\int_{\mathcal{Y}}g(y)\log(g(y))dy\nonumber\\
&+\int_{\mathcal{S}}d\lambda(T)\int_T\tilde{f}(x)\log(\tilde{f}(x))d\mathcal{H}^1(x)-\int_{\mathcal{S}}d\lambda(T)\int_T\tilde{g}(y)\log(\tilde{g}(y))d\mathcal{H}^1(y),
\end{align*}
\begin{equation*}
    \int_{\mathcal{S}^*}d\lambda(T)\int_{T}\log(\mathfrak{F}(x))\tilde{f}(x)d\mathcal{H}^1(x)=-\int_{\mathcal{X}}f(x)\log(f(x))dx+\int_{\mathcal{S}}d\lambda(T)\int_T\tilde{f}(x)\log(\tilde{f}(x))d\mathcal{H}^1(x).
\end{equation*}
\end{lemma}
\begin{proof}

Recall from \cref{L3n,L3nnn} that $\mathcal{S}^*$ is a Borel subset of $\mathcal{S}$ and $\lambda(\mathcal{S}\backslash\mathcal{S}^*)=0$. Now fix any $T\in\mathcal{S}^*$ and $x_0\in\mathrm{int}(T)$. By \cref{L2.4} (see \eqref{Eq2.6.4}; also note \cref{Defn4.1.1n}), for any $x,y\in\mathrm{int}(T)$,
\begin{align*}
    \det\big(\mathbf{I}_d+\langle x-x_0, V_T\rangle F(x_0)\big)&=\det\big(\mathbf{I}_d+\langle y-x_0, V_T\rangle F(x_0)\big)\det\big(\mathbf{I}_d-\langle y-x, V_T\rangle F(y)\big).
\end{align*}
By \cref{L2.7}, we have
\begin{equation*}
    \det\big(\mathbf{I}_d+\langle x-x_0, V_T\rangle F(x_0)\big),\det\big(\mathbf{I}_d+\langle y-x_0, V_T\rangle F(x_0)\big),\det\big(\mathbf{I}_d-\langle y-x, V_T\rangle F(y)\big)>0.
\end{equation*}
Hence noting \cref{factors}, we have (note that $\mathfrak{F}(x),\mathfrak{F}(y)>0$)
\begin{align*}
    &\log\det\big(\mathbf{I}_d-\langle y-x, V_T\rangle F(y)\big)\nonumber\\
    =\,&\log\det\big(\mathbf{I}_d+\langle x-x_0, V_T\rangle F(x_0)\big)-\log\det\big(\mathbf{I}_d+\langle y-x_0, V_T\rangle F(x_0)\big)\nonumber\\
    =\,&\log(\mathfrak{F}(x))-\log(\mathfrak{F}(y))=\mathbbm{1}_{\mathfrak{F}(x)>0}\log(\mathfrak{F}(x))-\mathbbm{1}_{\mathfrak{F}(y)>0}\log(\mathfrak{F}(y)).
\end{align*}
Therefore, we have
\begin{align}\label{E3.54n}
&\int_{\mathcal{S}}d\lambda(T)\int_{T\times T}\big|\log\det\big(\mathbf{I}_d+\|x-y\|F(y)\big)\big|d\kappa_T'(x,y)\nonumber\\
\leq\,&\int_{\mathcal{S}}d\lambda(T)\int_{T\times T}\big(\mathbbm{1}_{\mathfrak{F}(x)>0}|\log(\mathfrak{F}(x))| +\mathbbm{1}_{\mathfrak{F}(y)>0}|\log(\mathfrak{F}(y))|\big)d\kappa_T'(x,y)
\end{align}
and
\begin{align}\label{E3.54}
&\int_{\mathcal{S}}d\lambda(T)\int_{T\times T}\log\det\big(\mathbf{I}_d+\|x-y\|F(y)\big)d\kappa_T'(x,y)\nonumber\\
=\,&\int_{\mathcal{S}}d\lambda(T)\int_{T\times T}\big(\mathbbm{1}_{\mathfrak{F}(x)>0}\log(\mathfrak{F}(x))-\mathbbm{1}_{\mathfrak{F}(y)>0}\log(\mathfrak{F}(y))\big)d\kappa_T'(x,y).
\end{align}

Since $\kappa_T'\in\Pi(\tilde{\mu}_T,\tilde{\nu}_T)$ for $\lambda$-a.e.\ $T\in\mathcal{S}$, we have 
\begin{align*}
    &\int_{\mathcal{S}}d\lambda(T)\int_{T\times T} \mathbbm{1}_{\mathfrak{F}(x)>0}|\log(\mathfrak{F}(x))|d\kappa_T'(x,y)=\int_{\mathcal{S}}d\lambda(T)\int_{T} \mathbbm{1}_{\mathfrak{F}(x)>0}|\log(\mathfrak{F}(x))|d\tilde{\mu}_T(x)\nonumber\\
    =& \int_{\mathbb{R}^d} \mathbbm{1}_{\mathfrak{F}(x)>0}|\log(\mathfrak{F}(x))|f(x)dx= \int_{\mathbb{R}^d}\mathbbm{1}_{\mathfrak{F}(x)>0}|\log(\tilde{f}(x))-\log(f(x))|f(x)dx\nonumber\\
    \leq& \int_{\mathbb{R}^d}f(x)|\log(\tilde{f}(x))|dx+\int_{\mathcal{X}}f(x)|\log(f(x))|dx  \nonumber\\
    =&\int_{\mathcal{S}}d\lambda(T)\int_{T}\tilde{f}(x)|\log(\tilde{f}(x))|d\mathcal{H}^1(x)+\int_{\mathcal{X}}f(x)|\log(f(x))|dx\leq C,
\end{align*}
where we use \cref{L3nnn} in the second and fourth lines, and use \eqref{domainconst}, \eqref{bddLip}, and \cref{L3.18} in the last inequality. Note that by \cref{positivityF}, the above display implies 
\begin{equation*}
    \int_{\mathcal{S}^*}d\lambda(T)\int_{T} |\log(\mathfrak{F}(x))|\tilde{f}(x)d\mathcal{H}^1(x)=\int_{\mathcal{S}}d\lambda(T)\int_{T} \mathbbm{1}_{\mathfrak{F}(x)>0}|\log(\mathfrak{F}(x))|d\tilde{\mu}_T(x)\leq C.
\end{equation*}
Similarly, we can deduce that $\int_{\mathcal{S}}d\lambda(T)\int_{T\times T} \mathbbm{1}_{\mathfrak{F}(y)>0}|\log(\mathfrak{F}(y))|d\kappa_T'(x,y)\leq C$. Hence by \eqref{E3.54n} and \eqref{E3.54}, we have
\begin{equation*}
    \int_{\mathcal{S}}d\lambda(T)\int_{T\times T}\big|\log\det\big(\mathbf{I}_d+\|x-y\|F(y)\big)\big|d\kappa_T'(x,y)\leq C,
\end{equation*}
\begin{align*}
    &\int_{\mathcal{S}}d\lambda(T)\int_{T\times T}\log\det\big(\mathbf{I}_d+\|x-y\|F(y)\big)d\kappa_T'(x,y)\nonumber\\
    =\,& \int_{\mathcal{S}}d\lambda(T)\int_{T\times T}\mathbbm{1}_{\mathfrak{F}(x)>0}\log(\mathfrak{F}(x))d\kappa_T'(x,y)-\int_{\mathcal{S}}d\lambda(T)\int_{T\times T}\mathbbm{1}_{\mathfrak{F}(y)>0}\log(\mathfrak{F}(y))d\kappa_T'(x,y)\nonumber\\
    =\,& \int_{\mathcal{S}}d\lambda(T)\int_{T}\mathbbm{1}_{\mathfrak{F}(x)>0}\log(\mathfrak{F}(x))d\tilde{\mu}_T(x)-\int_{\mathcal{S}}d\lambda(T)\int_{T}\mathbbm{1}_{\mathfrak{F}(y)>0}\log(\mathfrak{F}(y))d\tilde{\nu}_T(y)\nonumber\\
    =\,& \int_{\mathbb{R}^d}\mathbbm{1}_{\mathfrak{F}(x)>0}\log(\mathfrak{F}(x))f(x)dx-\int_{\mathbb{R}^d}\mathbbm{1}_{\mathfrak{F}(y)>0}\log(\mathfrak{F}(y))g(y)dy\nonumber\\
    =\,& \int_{\mathbb{R}^d}\mathbbm{1}_{\mathfrak{F}(x)>0}(\log(\tilde{f}(x))-\log(f(x)))f(x)dx-\int_{\mathbb{R}^d}\mathbbm{1}_{\mathfrak{F}(y)>0}(\log(\tilde{g}(y))-\log(g(y)))g(y)dy\nonumber\\
    =\,& \int_{\mathbb{R}^d}(\log(\tilde{f}(x))-\log(f(x)))f(x)dx-\int_{\mathbb{R}^d}(\log(\tilde{g}(y))-\log(g(y)))g(y)dy\nonumber\\
    =\,& -\int_{\mathcal{X}}f(x)\log(f(x))dx+\int_{\mathcal{Y}}g(y)\log(g(y))dy\nonumber\\
    &+\int_{\mathcal{S}}d\lambda(T)\int_T\tilde{f}(x)\log(\tilde{f}(x))d\mathcal{H}^1(x)-\int_{\mathcal{S}}d\lambda(T)\int_T\tilde{g}(y)\log(\tilde{g}(y))d\mathcal{H}^1(y),
\end{align*}
where in the second display, we use \cref{L3nnn} in the fourth and seventh lines, and use \cref{positivityF} and \cref{Lem3.15} in the sixth line. Moreover, we have
\begin{align*}
&\int_{\mathcal{S}^*}d\lambda(T)\int_{T}\log(\mathfrak{F}(x))\tilde{f}(x)d\mathcal{H}^1(x)=\int_{\mathcal{S}}d\lambda(T)\int_{T} \mathbbm{1}_{\mathfrak{F}(x)>0}\log(\mathfrak{F}(x))d\tilde{\mu}_T(x)\nonumber\\
=\,& \int_{\mathbb{R}^d}\mathbbm{1}_{\mathfrak{F}(x)>0}\log(\mathfrak{F}(x))f(x)dx=\int_{\mathbb{R}^d}\mathbbm{1}_{\mathfrak{F}(x)>0}(\log(\tilde{f}(x))-\log(f(x)))f(x)dx\nonumber\\
=\,&\int_{\mathbb{R}^d}(\log(\tilde{f}(x))-\log(f(x)))f(x)dx=-\int_{\mathcal{X}}f(x)\log(f(x))dx+\int_{\mathcal{S}}d\lambda(T)\int_T\tilde{f}(x)\log(\tilde{f}(x))d\mathcal{H}^1(x).
\end{align*}
\end{proof}

Lastly, we introduce and bound a variant of the function $G_{\eps}$ introduced in \Cref{Def2}. 

\begin{definition}[The function $\bar{G}_{\eps}$]\label{Def2n}
For any $\eps>0$, $x\in \mathcal{T}_1^*$, and $z\in L(x)$, we define
\begin{equation*}
     \bar{G}_{\eps}(x,z):=\log\bigg(\varepsilon^{-(d-1)\slash 2}\int_{O(V(x))}e^{-E(x,z+w)\slash\varepsilon}\mathbbm{1}_{\|w\|\leq 2D,z+w\in\mathcal{Y}} d\mathcal{H}^{d-1}(w)\bigg). 
\end{equation*}
\end{definition}
\begin{remark}\label{barGequation}
Note that $\bar{G}_{\eps}(x,z)\leq G_{\eps}(x,z)$ for all $\eps>0$, $x\in \mathcal{T}_1^*$, and $z\in L(x)$. 
\end{remark}

\begin{lemma}\label{NewGepsconvergence}
For any $x\in\mathcal{T}_1^*$ and $z\in L(x)$ such that $\langle z-b(x),V(x)\rangle<0$ and $z\notin \mathcal{Y}$, we have $\lim_{\eps\rightarrow 0^{+}}\bar{G}_{\eps}(x, z)=-\infty$.
\end{lemma}
\begin{proof}

By \Cref{L3.2}, for any $w\in O(V(x))$ such that $\|w\|\leq 2D$,
\begin{equation*}
    E(x,z+w)\geq \frac{\|x-a(x)\|\|w\|^2}{4(\max\{\|z-a(x)\|,\|z-x\|\}+2D)^2}\geq c'\|w\|^2, 
\end{equation*}
where $c'>0$ only depends on $x,z,D$. As $\mathcal{Y}$ is a closed set and $z\notin\mathcal{Y}$, there exists $\delta_0\in (0,1)$ such that $B_d(z,\delta_0)\cap \mathcal{Y}=\emptyset$. Hence for any $w\in O(V(x))$ such that $\|w\|<\delta_0$, we have $z+w\notin\mathcal{Y}$. Consequently, 
\begin{align*}
      \int_{O(V(x))}e^{-E(x,z+w)\slash\eps}\mathbbm{1}_{\|w\|\leq 2D,z+w\in\mathcal{Y}} d\mathcal{H}^{d-1}(w)\leq\int_{\{w\in O(V(x)):\delta_0\leq \|w\|\leq 2D\}} e^{-E(x,z+w)\slash\eps}d\mathcal{H}^{d-1}(w) &\nonumber\\
    \leq \int_{\{w\in O(V(x)):\delta_0\leq \|w\|\leq 2D\}} e^{-c'\|w\|^2\slash\eps}d\mathcal{H}^{d-1}(w)\leq C e^{-c'\delta_0^2\slash\eps}& ,
\end{align*}
which implies $\bar{G}_{\eps}(x,z)\leq \log\big(C\eps^{-(d-1)\slash 2}e^{-c'\delta_0^2\slash \eps}\big)$. Therefore, $\lim_{\eps\rightarrow 0^{+}}\bar{G}_{\eps}(x, z)=-\infty$. 
\end{proof}

We now complete the proof of Theorem \ref{th:lowerBound}.

\begin{proof}[Proof of \cref{th:lowerBound}]
We denote $\mathcal{X}^{\circ}:=\{x\in\mathbb{R}^d:f(x)>0\}$ and $\mathcal{Y}^{\circ}:=\{y\in\mathbb{R}^d:g(y)>0\}$. For every $j\in\mathbb{N}^*$, by \eqref{OTrela}, we have 
\begin{align*}
  \mathcal{C}_{\eps_j}(\gamma'_{\eps_j})-\OT(\mu,\nu)
   =\,& \int_{\mathbb{R}^d\times\mathbb{R}^d}\|x-y\|\phi_{\eps_j}(x,y)dx dy-\int_{\mathbb{R}^d}u(x)(f(x)-g(x))dx\nonumber\\
   &+\eps_j\int_{\mathcal{X}^{\circ}\times\mathcal{Y}^{\circ}}\phi_{\eps_j}(x,y)\log\bigg(\frac{\phi_{\eps_j}(x,y)}{f(x)g(y)}\bigg)dxdy.
\end{align*} 
Hence, noting that $\gamma_{\eps_j}'\in \Pi(\mu,\nu)$ and $E(x,y)=\|x-y\|-u(x)+u(y)$ for all $x,y\in\mathbb{R}^d$, we obtain
\begin{align}\label{Eq2.20}
    &\mathcal{C}_{\eps_j}(\gamma'_{\eps_j})-\OT(\mu,\nu)\nonumber\\
   =\,& \int_{\mathbb{R}^d\times\mathbb{R}^d}(\|x-y\|-u(x)+u(y))\phi_{\eps_j}(x,y)dxdy+\eps_j\int_{\mathcal{X}^{\circ}\times\mathcal{Y}^{\circ}}\phi_{\eps_j}(x,y)\log (\phi_{\eps_j}(x,y))dxdy\nonumber\\
   & -\eps_j\int_{\mathcal{X}^{\circ}}f(x)\log(f(x))dx-\eps_j\int_{\mathcal{Y}^{\circ}}g(y)\log(g(y))dy\nonumber\\
   =\,& \int_{\mathcal{X}\times \mathcal{Y}} \phi_{\eps_j}(x,y)\big(E(x,y)+\eps_j\log (\phi_{\eps_j}(x,y))\big)dxdy\nonumber\\
   &  -\eps_j\int_{\mathcal{X}}f(x)\log(f(x))dx-\eps_j\int_{\mathcal{Y}}g(y)\log(g(y))dy.
\end{align} 

In the following, we lower bound $\int_{\mathcal{X}\times \mathcal{Y}} \phi_{\eps_j}(x,y)\big(E(x,y)+\eps_j\log (\phi_{\eps_j}(x,y))\big)dxdy$. With $\Phi(\cdot)$ defined as in \eqref{Phi_def}, using \eqref{new1.1}, we have
\begin{align*}
 &\int_{\mathcal{X}\times\mathcal{Y}} \phi_{\eps_j}(x,y)\big(E(x,y)+\eps_j\log(\phi_{\eps_j}(x,y))\big) dxdy\nonumber\\
=\,&\int_{(\mathcal{X}\cap\mathcal{T}_1^{*})\times\mathcal{Y}} \phi_{\eps_j}(x,y)\big(E(x,y)+\varepsilon_j\log(\phi_{\eps_j}(x,y))\big) dxdy\nonumber\\
=\,& \int_{\mathcal{X}\cap\mathcal{T}_1^{*}} dx\int_{L(x)}d\mathcal{H}^1(z)\nonumber\\
&\hspace{0.9in}\int_{O(V(x))}\mathbbm{1}_{z+w\in\mathcal{Y}}\cdot\phi_{\eps_j}(x,z+w)\big(E(x,z+w)+\varepsilon_j\log(\phi_{\eps_j}(x,z+w))\big)
d\mathcal{H}^{d-1}(w)\nonumber\\
=\,& \varepsilon_j\int_{\mathcal{X}\cap\mathcal{T}_1^{*}} dx\int_{L(x)} d\mathcal{H}^1(z)\int_{O(V(x))} \mathbbm{1}_{z+w\in\mathcal{Y}}\cdot e^{-E(x,z+w)\slash\varepsilon_j}\Phi\big(\phi_{\eps_j}(x,z+w)e^{E(x,z+w)\slash\varepsilon_j}\big)
d\mathcal{H}^{d-1}(w).
\end{align*}
Hence by Jensen's inequality and the convexity of $\Phi(\cdot)$ on $[0,\infty)$, 
\begin{align}\label{Eq2.19}
 &\int_{\mathcal{X}\times\mathcal{Y}} \phi_{\eps_j}(x,y)\big(E(x,y)+\eps_j\log(\phi_{\eps_j}(x,y))\big) dxdy\nonumber\\
\geq\,& \varepsilon_j\int_{\mathcal{X}\cap\mathcal{T}_1^{*}} dx\int_{L(x)} \bigg(\int_{O(V(x))}\mathbbm{1}_{z+w\in\mathcal{Y}} \cdot e^{-E(x,z+w)\slash\varepsilon_j}d\mathcal{H}^{d-1}(w)\bigg)\nonumber\\
&\hspace{1.1in}\times\Phi\bigg(\frac{\int_{O(V(x))}\mathbbm{1}_{z+w\in\mathcal{Y}}\cdot\phi_{\eps_j}(x,z+w)d\mathcal{H}^{d-1}(w)}{\int_{O(V(x))}\mathbbm{1}_{z+w\in\mathcal{Y}}\cdot e^{-E(x,z+w)\slash\varepsilon_j}d\mathcal{H}^{d-1}(w)}\bigg)d\mathcal{H}^1(z)\nonumber\\
=\,& \varepsilon_j\int_{\mathcal{X}\cap\mathcal{T}_1^{*}} dx\int_{L(x)}\bigg(\int_{O(V(x))}\phi_{\eps_j}(x,z+w)d\mathcal{H}^{d-1}(w)\bigg)\nonumber\\
&\hspace{1.1in}\times\log\bigg(\frac{\int_{O(V(x))}\phi_{\eps_j}(x,z+w)d\mathcal{H}^{d-1}(w)}{\int_{O(V(x))}e^{-E(x,z+w)\slash\varepsilon_j} \mathbbm{1}_{\|w\|\leq 2D, z+w\in\mathcal{Y}} d\mathcal{H}^{d-1}(w)}\bigg)d\mathcal{H}^1(z)\nonumber\\
=\,& \eps_j\int_{(\mathcal{X}\cap\mathcal{T}_1^{*})\times\mathbb{R}}\psi_{\eps_j}(x,t)\log\bigg(\frac{\psi_{\eps_j}(x,t)}{\eps_j^{(d-1)\slash 2}\exp(\bar{G}_{\eps_j}(x,x+tV(x)))}\bigg)dxdt\nonumber\\
=\,& -\frac{d-1}{2}\eps_j\log\eps_j+\eps_j\int_{(\mathcal{X}\cap\mathcal{T}_1^{*})\times\mathbb{R}} \psi_{\eps_j}(x,t)\log(\psi_{\eps_j}(x,t))dxdt\nonumber\\
&-\eps_j\int_{(\mathcal{X}\cap\mathcal{T}_1^{*})\times\mathbb{R}}\bar{G}_{\eps_j}(x,x+tV(x))\psi_{\eps_j}(x,t)dxdt,
\end{align}
where the first equality uses \eqref{phij_0} and \cref{Lemma3.3}, and the second equality uses the definitions of $\psi_{\eps_j}(\cdot,\cdot)$ and $\bar{G}_{\eps_j}(\cdot,\cdot)$ as in \eqref{P3.1.eq1} and \cref{Def2n}. 


Note that by \eqref{P3.1.eq1}, for any $(x,t)\in\mathbb{R}^d\times \mathbb{R}$, if $\psi_{\eps_j}(x,t)>0$, then $x\in\mathcal{T}_1^{*}$, and there exists $w\in O(V(x))$ such that $\phi_{\varepsilon_j}(x,x+tV(x)+w)>0$. Hence by \eqref{phij_0}, $x\in\mathcal{X}\cap\mathcal{T}_1^{*}$ and $x+tV(x)+w\in\mathcal{Y}$, which by Lemma \ref{Lemma3.3} (with $z=x+tV(x)$) implies $\max\{|t|,\|(x+tV(x))-a(x)\|\}\leq 2D$. Therefore, for any $(x,t)\in\mathbb{R}^d\times\mathbb{R}$,
\begin{equation}\label{P3.1.eq2}
  \psi_{\eps_j}(x,t)>0 \,\Rightarrow\,\max\{|t|,\|(x+tV(x))-a(x)\|\}\leq 2D.
\end{equation}

Let $\mathscr{G}$ denote the set of $x\in\mathcal{T}_1^{*}$ such that the conclusions of \Cref{P-lb-1,P-lb-2} hold for $x$. Note that $\mathcal{L}^d(\mathcal{T}_1^{*}\backslash \mathscr{G})=0$. In the following, we fix any $M, M'\geq 10$ and $\delta\in (0,1)$, and set $\mathcal{A}_{\delta}:=\big\{x\in\mathbb{R}^d:\alpha(x)\geq \delta\big\}\cap\mathscr{G}$. By \eqref{L3.11.e1}, we have  
\begin{align}\label{P3.1.e4}
    &\int_{(\mathcal{X}\cap\mathcal{T}_1^{*}\cap\mathcal{A}_{\delta})\times\mathbb{R}}\bar{G}_{\eps_j}(x,x+tV(x))\psi_{\eps_j}(x,t)dxdt\nonumber\\
    \leq\,& \int_{(\mathcal{X}\cap\mathcal{T}_1^{*}\cap\mathcal{A}_{\delta})\times\mathbb{R}}\big(\bar{G}_{\eps_j}(x,x+tV(x))\vee (-M)\big) \psi_{\eps_j}(x,t)dxdt\nonumber\\
    =\,& \int_{(\mathcal{X}\cap\mathcal{T}_1^{*}\cap\mathcal{A}_{\delta})\times\mathbb{R}}\big(\bar{G}_{\eps_j}(x,x+tV(x))\vee (-M)+M\big)\psi_{\eps_j}(x,t)dxdt- M \mu(\mathcal{X}\cap\mathcal{T}_1^{*}\cap\mathcal{A}_{\delta}). 
\end{align}
For any $(x,t)\in (\mathcal{X}\cap\mathcal{T}_1^*\cap\mathcal{A}_{\delta})\times\mathbb{R}$ such that $\psi_{\eps_j}(x,t)>0$, by \eqref{P3.1.eq2}, $\max\{|t|,\|(x+tV(x))-a(x)\|\}\leq 2D$. Hence by \Cref{barGequation} and \cref{P-lb-2}(a) (note that $\alpha(x)\geq \delta$), we have
\begin{equation}\label{P3.1.e1}
   \bar{G}_{\eps_j}(x,x+tV(x))  \leq  G_{\eps_j}(x,x+tV(x))\leq -\frac{d-1}{2}\log(\alpha(x))+C\leq \frac{d-1}{2}\log(\delta^{-1})+C\leq C\log(e\delta^{-1}), 
\end{equation}
\begin{align}\label{P3.1.e2}
    \psi_{\eps_j}(x,t)\mathbbm{1}_{ \psi_{\eps_j}(x,t)> M'}\leq\,& \frac{\max\big\{\psi_{\eps_j}(x,t) \log(\psi_{\eps_j}(x,t)),0\big\}}{\log(M')}\nonumber\\
  \leq\,&\frac{\psi_{\eps_j}(x,t) \log(\psi_{\eps_j}(x,t))+1}{\log(M')}=\frac{\psi_{\eps_j}(x,t) \log(\psi_{\eps_j}(x,t))+\mathbbm{1}_{|t|\leq 2D}}{\log(M')},
\end{align}
where the second inequality in \eqref{P3.1.e2} uses $s\log s\geq -e^{-1}$ for any $s\geq 0$, and the equality uses $|t|\leq 2D$. By \eqref{P3.1.e1} and \eqref{P3.1.e2},
\begin{align}\label{P3.1.e5}
& \int_{(\mathcal{X}\cap\mathcal{T}_1^{*}\cap\mathcal{A}_{\delta})\times\mathbb{R}}\big(\bar{G}_{\eps_j}(x,x+tV(x))\vee(-M)+M\big) \psi_{\eps_j}(x,t)\mathbbm{1}_{\psi_{\eps_j}(x,t)> M'}dxdt\nonumber\\
\leq\,& (C\log(e\delta^{-1})+M)\int_{(\mathcal{X}\cap\mathcal{T}_1^{*}\cap\mathcal{A}_{\delta})\times\mathbb{R}}\psi_{\eps_j}(x,t)\mathbbm{1}_{\psi_{\eps_j}(x,t)> M'}dxdt\nonumber\\
\leq\,& \frac{C\log\big(e\delta^{-1}\big)+M}{\log(M')}\int_{(\mathcal{X}\cap\mathcal{T}_1^{*}\cap\mathcal{A}_{\delta})\times\mathbb{R}}\big(\psi_{\eps_j}(x,t) \log(\psi_{\eps_j}(x,t))+\mathbbm{1}_{|t|\leq 2D}\big) dxdt\nonumber\\
\leq\,& \frac{C(\log\big(e\delta^{-1}\big)+M)}{\log(M')}\int_{(\mathcal{X}\cap\mathcal{T}_1^{*}\cap\mathcal{A}_{\delta})\times\mathbb{R}}\psi_{\eps_j}(x,t) \log(\psi_{\eps_j}(x,t))dxdt+\frac{C(\log(e\delta^{-1})+M)\mathcal{L}^d(\mathcal{X})}{\log(M')}\nonumber\\
\leq\,& \frac{C(\log\big(e\delta^{-1}\big)+M)}{\log(M')}\int_{(\mathcal{X}\cap\mathcal{T}_1^{*})\times\mathbb{R}}\psi_{\eps_j}(x,t) \log(\psi_{\eps_j}(x,t))dxdt+\frac{C(\log\big(e\delta^{-1}\big)+M)}{\log(M')},
\end{align}
where the last inequality uses (note \eqref{P3.1.eq2})
\begin{equation*}
    \int_{((\mathcal{X}\cap\mathcal{T}_1^{*})\backslash\mathcal{A}_{\delta})\times\mathbb{R}}\psi_{\eps_j}(x,t) \log(\psi_{\eps_j}(x,t))dxdt\geq -e^{-1}\mathcal{L}^d\otimes\mathcal{L}^1(\mathcal{X}\times [-2D,2D])\geq -C.
\end{equation*}

Combining \eqref{Eq2.19}, \eqref{P3.1.e4}, and \eqref{P3.1.e5}, we obtain 
\begin{align}\label{Eqq3.1}
&\eps_j^{-1}\int_{\mathcal{X}\times\mathcal{Y}} \phi_{\eps_j}(x,y)\big(E(x,y)+\eps_j\log(\phi_{\eps_j}(x,y))\big)dxdy+\frac{d-1}{2}\log{\eps_j}\nonumber\\
\geq\,&
-\int_{(\mathcal{X}\cap\mathcal{T}_1^{*}\cap\mathcal{A}_{\delta})\times\mathbb{R}}\bar{G}_{\eps_j}(x,x+tV(x))\psi_{\eps_j}(x,t)dxdt -\int_{((\mathcal{X}\cap\mathcal{T}_1^{*})\backslash\mathcal{A}_{\delta})\times\mathbb{R}}\bar{G}_{\eps_j}(x,x+tV(x))\psi_{\eps_j}(x,t)dxdt\nonumber\\
&+\int_{(\mathcal{X}\cap\mathcal{T}_1^{*})\times\mathbb{R}} \psi_{\eps_j}(x,t)\log(\psi_{\eps_j}(x,t))dxdt\nonumber\\
\geq\,& - \int_{(\mathcal{X}\cap\mathcal{T}_1^{*}\cap\mathcal{A}_{\delta})\times\mathbb{R}}\big(\bar{G}_{\eps_j}(x,x+tV(x))\vee (-M)+M\big)\psi_{\eps_j}(x,t)\mathbbm{1}_{\psi_{\eps_j}(x,t)\leq M'}dxdt\nonumber\\
&+\bigg(1-\frac{C(\log\big(e\delta^{-1}\big)+M)}{\log(M')}\bigg)\int_{(\mathcal{X}\cap\mathcal{T}_1^{*})\times\mathbb{R}}\psi_{\eps_j}(x,t) \log(\psi_{\eps_j}(x,t))dxdt\nonumber\\
&-\sup_{j\in \mathbb{N}^*}\bigg\{\int_{((\mathcal{X}\cap\mathcal{T}_1^{*})\backslash\mathcal{A}_{\delta})\times\mathbb{R}}\bar{G}_{\eps_j}(x,x+tV(x))\psi_{\eps_j}(x,t)dxdt\bigg\}\nonumber\\
&+M \mu(\mathcal{X}\cap\mathcal{T}_1^{*}\cap\mathcal{A}_{\delta})-\frac{C(\log(e\delta^{-1})+M)}{\log(M')}.
\end{align}

Sequentially applying \Cref{barGequation}, \eqref{P3.1.eq2}, \cref{P-lb-2}(a), and \eqref{L3.11.e1}, we obtain that for any $j\in \mathbb{N}^*$, 
\begin{align}\label{P3.1.eeq1}
\int_{((\mathcal{X}\cap\mathcal{T}_1^{*})\backslash\mathcal{A}_{\delta})\times\mathbb{R}}\bar{G}_{\eps_j}(x,x+tV(x))\psi_{\eps_j}(x,t)dxdt\leq \int_{((\mathcal{X}\cap\mathcal{T}_1^{*})\backslash\mathcal{A}_{\delta})\times\mathbb{R}}G_{\eps_j}(x,x+tV(x))\psi_{\eps_j}(x,t)dxdt &\nonumber\\
\leq C\int_{((\mathcal{X}\cap\mathcal{T}_1^{*})\backslash\mathcal{A}_{\delta})\times\mathbb{R}}\log(e\max\{\alpha(x)^{-1},1\})\psi_{\eps_j}(x,t)dxdt&\nonumber\\
\leq C\int_{(\mathcal{X}\cap\mathcal{T}_1^{*})\backslash\mathcal{A}_{\delta}}\log(e\max\{\alpha(x)^{-1},1\})f(x)dx.\hspace{0.6in}& \raisetag{1.42\baselineskip}
\end{align}
By \cref{Lemma2.7},
\begin{equation}\label{P3.1.eeq2}
    \int_{\mathcal{X}\cap\mathcal{T}_1^{*}}\log(e\max\{\alpha(x)^{-1},1\})f(x)dx<\infty.
\end{equation}
Let $\tilde{A}_{\delta}:=\{x\in \mathbb{R}^d:\alpha(x)\geq \delta\}$. As $\mathcal{L}^d(\mathcal{X})<\infty$ and 
\begin{equation}\label{P3.1.eqnew8}
    \lim_{\delta\rightarrow 0^{+}}\mathcal{T}_1^{*}\backslash\tilde{\mathcal{A}}_{\delta}=\emptyset,
\end{equation}
we have $\lim_{\delta\rightarrow 0^{+}} \mathcal{L}^d((\mathcal{X}\cap\mathcal{T}_1^{*})\backslash\tilde{\mathcal{A}}_{\delta})=0$. Noting $\mathcal{A}_{\delta}=\tilde{\mathcal{A}}_{\delta}\cap\mathscr{G}$ and $\mathcal{L}^d(\mathcal{T}_1^{*}\backslash \mathscr{G})=0$, we thus have 
\begin{equation}\label{P3.1.eeq3}
    \lim_{\delta\rightarrow 0^{+}} \mathcal{L}^d((\mathcal{X}\cap\mathcal{T}_1^{*})\backslash\mathcal{A}_{\delta})=\lim_{\delta\rightarrow 0^{+}} \mathcal{L}^d((\mathcal{X}\cap\mathcal{T}_1^{*})\backslash\tilde{\mathcal{A}}_{\delta})=0.
\end{equation}
By \eqref{P3.1.eeq2} and \eqref{P3.1.eeq3}, we have
\begin{equation*}
    \lim_{\delta\rightarrow 0^{+}}\int_{(\mathcal{X}\cap\mathcal{T}_1^{*})\backslash\mathcal{A}_{\delta}}\log(e\max\{\alpha(x)^{-1},1\})f(x)dx = 0.
\end{equation*}
Hence by \eqref{P3.1.eeq1},
\begin{equation}\label{P3.1.eq8}
    \limsup_{\delta\rightarrow 0^{+}}\sup_{j\in \mathbb{N}^*}\bigg\{\int_{((\mathcal{X}\cap\mathcal{T}_1^{*})\backslash\mathcal{A}_{\delta})\times\mathbb{R}}\bar{G}_{\eps_j}(x,x+tV(x))\psi_{\eps_j}(x,t)dxdt\bigg\}\leq 0.
\end{equation}

By \cref{L3.6}, we obtain 
\begin{align}\label{P3.1.eq7}
&\liminf_{j\rightarrow\infty}\bigg\{\int_{(\mathcal{X}\cap\mathcal{T}_1^{*})\times\mathbb{R}}\psi_{\varepsilon_j}(x,t)\log(\psi_{\varepsilon_j}(x,t))dxdt\bigg\}\nonumber\\
\geq\,& \int_{\mathcal{S}}H(\kappa'_T|\tilde{\mu}_T\otimes\tilde{\nu}_T)d\lambda(T) + \int_{\mathcal{X}}f(x)\log(f(x))dx+\int_{\mathcal{S}}d\lambda(T)\int_T \tilde{g}(y)\log(\tilde{g}(y))d\mathcal{H}^1(y).
\end{align}

Let $\mathbb{T}_0:=\{t<-\beta(x):x+tV(x)\notin\mathcal{Y}\}$. By \eqref{P3.1.eq2} and \cref{P-lb-1}, for $\mathcal{L}^d$-a.e.\ $x\in\mathcal{T}_1^{*}$ and $\mathcal{L}^1$-a.e.\ $t\in\mathbb{R}\backslash\mathbb{T}_0$, as $j\rightarrow\infty$ we have 
\begin{align*}
   & \big|G_{\eps_j}(x,x+tV(x)) \vee (-M)-G(x,x+tV(x)) \vee (-M)\big|\psi_{\eps_j}(x,t)\mathbbm{1}_{\psi_{\eps_j}(x,t)\leq M'} \nonumber\\
    \leq\, & M'\big|G_{\eps_j}(x,x+tV(x))\vee(-M)-G(x,x+tV(x))\vee(-M)\big|\mathbbm{1}_{\max\{|t|,\|(x+tV(x))-a(x)\|\}\leq 2D}\rightarrow 0.
\end{align*}
Similarly, using \Cref{NewGepsconvergence}, we obtain that for any $x\in\mathcal{T}_1^*$ and $t\in \mathbb{T}_0$, as $j\rightarrow\infty$, 
\begin{equation*}
     \big|\bar{G}_{\eps_j}(x,x+tV(x)) \vee (-M)-G(x,x+tV(x)) \vee (-M)\big|\psi_{\eps_j}(x,t)\mathbbm{1}_{\psi_{\eps_j}(x,t)\leq M'} \rightarrow 0.
\end{equation*}
By \cref{P-lb-2} and \Cref{barGequation},
for any $x\in\mathcal{X}  \cap \mathcal{T}_1^{*}\cap\mathcal{A}_{\delta}$, if $t\in\mathbb{R}\backslash\mathbb{T}_0$, then
\begin{align*}
    &\big|G_{\eps_j}(x,x+tV(x)) \vee (-M)-G(x,x+tV(x)) \vee (-M)\big|\psi_{\eps_j}(x,t)\mathbbm{1}_{\psi_{\eps_j}(x,t)\leq M'}\nonumber\\
    \leq\, & C\max\{\log\big(e\delta^{-1}\big), M\}M'\mathbbm{1}_{|t|\leq 2D};
\end{align*}
if $t\in\mathbb{T}_0$, then 
\begin{align*}
    &\big|\bar{G}_{\eps_j}(x,x+tV(x)) \vee (-M)-G(x,x+tV(x)) \vee (-M)\big|\psi_{\eps_j}(x,t)\mathbbm{1}_{\psi_{\eps_j}(x,t)\leq M'}\nonumber\\
    \leq\, & C\max\{\log\big(e\delta^{-1}\big), M\}M'\mathbbm{1}_{|t|\leq 2D}.
\end{align*}
Hence by the dominated convergence theorem, as $j  \rightarrow \infty$, 
\begin{equation*}
    \int_{(\mathcal{X}\cap\mathcal{T}_1^{*}\cap\mathcal{A}_{\delta})\times(\mathbb{R}\backslash\mathbb{T}_0)}\big|G_{\eps_j}(x,x+tV(x))\vee (-M)-G(x,x+tV(x))\vee (-M)\big|\psi_{\eps_j}(x,t)\mathbbm{1}_{\psi_{\eps_j}(x,t)\leq M'}dxdt \rightarrow 0,
\end{equation*}
\begin{equation*}
    \int_{(\mathcal{X}\cap\mathcal{T}_1^{*}\cap\mathcal{A}_{\delta})\times\mathbb{T}_0}\big|\bar{G}_{\eps_j}(x,x+tV(x))\vee (-M)-G(x,x+tV(x))\vee (-M)\big|\psi_{\eps_j}(x,t)\mathbbm{1}_{\psi_{\eps_j}(x,t)\leq M'}dxdt \rightarrow 0.
\end{equation*}
Therefore, by \Cref{barGequation}, \eqref{L3.11.e1}, and \cref{L3.7}, we have
\begin{align}\label{Eqq3.2}
   & \limsup_{j\rightarrow \infty}\int_{(\mathcal{X}\cap\mathcal{T}_1^{*}\cap\mathcal{A}_{\delta})\times\mathbb{R}}\big(\bar{G}_{\eps_j}(x,x+tV(x))\vee (-M)+M\big)\psi_{\eps_j}(x,t)\mathbbm{1}_{\psi_{\eps_j}(x,t)\leq M'}dxdt\nonumber\\
   \leq\,& \limsup_{j\rightarrow \infty}\int_{(\mathcal{X}\cap\mathcal{T}_1^{*}\cap\mathcal{A}_{\delta})\times\mathbb{R}}\big(G(x,x+tV(x))\vee (-M)+M\big) \psi_{\eps_j}(x,t)\mathbbm{1}_{\psi_{\eps_j}(x,t)\leq M'}dxdt \nonumber\\
   \leq\,& \limsup_{j\rightarrow \infty}\int_{(\mathcal{X}\cap\mathcal{T}_1^{*}\cap\mathcal{A}_{\delta})\times\mathbb{R}}\big(G(x,x+tV(x))\vee (-M)+M\big)\psi_{\eps_j}(x,t)dxdt \nonumber\\
   =\,& \limsup_{j\rightarrow \infty}\int_{(\mathcal{X}\cap\mathcal{T}_1^{*}\cap\mathcal{A}_{\delta})\times\mathbb{R}} \big(G(x,x+tV(x))\vee (-M)\big)\psi_{\eps_j}(x,t)dxdt+M\mu(\mathcal{X}\cap\mathcal{T}_1^{*}\cap\mathcal{A}_{\delta}) \nonumber\\
   \leq\,& \int_{(\mathcal{X}\cap\mathcal{T}_1^{*}\cap\mathcal{A}_{\delta})\times\mathbb{R}^d}G(x,y)\vee (-M)d\gamma_0'(x,y)+M\mu(\mathcal{X}\cap\mathcal{T}_1^{*}\cap\mathcal{A}_{\delta})\nonumber\\
   =\,&  \int_{\mathcal{S}} d\lambda(T) \int_{T\times T} \mathbbm{1}_{\mathcal{X}\cap\mathcal{T}_1^{*}\cap\mathcal{A}_{\delta}}(x) \big(G(x,y)\vee (-M)\big) d\kappa'_T(x,y)+M\mu(\mathcal{X}\cap\mathcal{T}_1^{*}\cap\mathcal{A}_{\delta}).
\end{align}

For any $x\in\mathcal{T}_1^*$ and $y\in L(x)$, we denote
\begin{equation*}
    G_{+}(x,y):=\max\{G(x,y),0\},  \qquad  G_{-}(x,y):=\max\{-G(x,y),0\}.
\end{equation*}
By \cref{Lemma2.7} and \cref{P-lb-2}, we have
\begin{align}\label{Eq3.4n}
  &  \int_{\mathcal{S}} d\lambda(T) \int_{T\times T} \mathbbm{1}_{\mathcal{X}\cap\mathcal{T}_1^{*}}(x) G_{+}(x,y)d\kappa'_T(x,y)\nonumber\\
  \leq\,&  \int_{\mathcal{S}} d\lambda(T) \int_{T\times T} \mathbbm{1}_{\mathcal{X}\cap\mathcal{T}_1^{*}}(x)\Big(-\frac{d-1}{2}\log(\alpha(x))+C\Big) d\kappa'_T(x,y)\nonumber\\
  =\,& \int_{\mathcal{S}} d\lambda(T) \int_{T}\mathbbm{1}_{\mathcal{X}\cap\mathcal{T}_1^{*}}(x)\Big(-\frac{d-1}{2}\log(\alpha(x))+C\Big) \tilde{f}(x)d\mathcal{H}^1(x)\nonumber\\
  =\,& \int_{\mathcal{X}\cap\mathcal{T}_1^{*}}\Big(-\frac{d-1}{2}\log(\alpha(x))+C\Big)f(x)dx\nonumber\\
  \leq\,&-\frac{d-1}{2}\int_{\mathcal{X}\cap\mathcal{T}_1^{*}}f(x)\log(\alpha(x))dx+C<\infty.
\end{align}

By \eqref{Eqq3.1} and \eqref{Eqq3.2}, for $M'$ sufficiently large (depending on $M,\delta$), we have 
\begin{align*}
    & \liminf_{j\rightarrow\infty}\bigg\{\eps_j^{-1}\int_{\mathcal{X}\times\mathcal{Y}} \phi_{\eps_j}(x,y)\big(E(x,y)+\eps_j\log(\phi_{\eps_j}(x,y))\big)dxdy+\frac{d-1}{2}\log{\eps_j}\bigg\}\nonumber\\
    \geq\,& -\int_{\mathcal{S}} d\lambda(T) \int_{T\times T} \mathbbm{1}_{\mathcal{X}\cap\mathcal{T}_1^{*}\cap\mathcal{A}_{\delta}}(x) G(x,y)\vee (-M) d\kappa'_T(x,y)  \nonumber\\
    & +\bigg(1-\frac{C(\log\big(e\delta^{-1}\big)+M)}{\log(M')}\bigg)\liminf_{j\rightarrow\infty}\bigg\{\int_{(\mathcal{X}\cap\mathcal{T}_1^{*})\times\mathbb{R}}\psi_{\eps_j}(x,t) \log(\psi_{\eps_j}(x,t))dxdt\bigg\}\nonumber\\
    &-\sup_{j\in \mathbb{N}^*}\bigg\{\int_{((\mathcal{X}\cap\mathcal{T}_1^{*})\backslash\mathcal{A}_{\delta})\times\mathbb{R}}\bar{G}_{\eps_j}(x,x+tV(x))\psi_{\eps_j}(x,t)dxdt\bigg\}-\frac{C(\log(e\delta^{-1})+M)}{\log(M')}.
\end{align*}
Taking $M'\rightarrow\infty$, we get 
\begin{align}\label{Eqq3.3}
    & \liminf_{j\rightarrow\infty}\bigg\{\eps_j^{-1}\int_{\mathcal{X}\times\mathcal{Y}} \phi_{\eps_j}(x,y)\big(E(x,y)+\eps_j\log(\phi_{\eps_j}(x,y))\big)dxdy+\frac{d-1}{2}\log{\eps_j}\bigg\}\nonumber\\
    \geq\,& -\int_{\mathcal{S}} d\lambda(T) \int_{T\times T} \mathbbm{1}_{\mathcal{X}\cap\mathcal{T}_1^{*}\cap\mathcal{A}_{\delta}}(x) G(x,y)\vee (-M) d\kappa'_T(x,y)  \nonumber\\   &+\liminf_{j\rightarrow\infty}\bigg\{\int_{(\mathcal{X}\cap\mathcal{T}_1^{*})\times\mathbb{R}}\psi_{\eps_j}(x,t) \log(\psi_{\eps_j}(x,t))dxdt\bigg\}\nonumber\\
    &-\sup_{j\in \mathbb{N}^*}\bigg\{\int_{((\mathcal{X}\cap\mathcal{T}_1^{*})\backslash\mathcal{A}_{\delta})\times\mathbb{R}}\bar{G}_{\eps_j}(x,x+tV(x))\psi_{\eps_j}(x,t)dxdt\bigg\}.
\end{align}
Note that
for any $x\in\mathcal{X}\cap\mathcal{T}_1^{*}\cap\mathcal{A}_{\delta}$ and $y\in L(x)$, $(-G(x,y))\wedge M\geq -C\log\big(e\delta^{-1}\big)$, and $(-G(x,y))\wedge M$ is monotonically non-decreasing in $M$. Taking $M\rightarrow\infty$ in \eqref{Eqq3.3} and using the monotone convergence theorem, we obtain 
\begin{align}\label{Eqq3.4}
     & \liminf_{j\rightarrow\infty}\bigg\{\eps_j^{-1}\int_{\mathcal{X}\times\mathcal{Y}} \phi_{\eps_j}(x,y)\big(E(x,y)+\eps_j\log(\phi_{\eps_j}(x,y))\big)dxdy+\frac{d-1}{2}\log{\eps_j}\bigg\}\nonumber\\
    \geq\,& -\int_{\mathcal{S}} d\lambda(T) \int_{T\times T} \mathbbm{1}_{\mathcal{X}\cap\mathcal{T}_1^{*}\cap\mathcal{A}_{\delta}}(x) G(x,y) d\kappa_T'(x,y)\nonumber\\
    &+\liminf_{j\rightarrow\infty}\bigg\{\int_{(\mathcal{X}\cap\mathcal{T}_1^{*})\times\mathbb{R}}\psi_{\eps_j}(x,t) \log(\psi_{\eps_j}(x,t))dxdt\bigg\}\nonumber\\
    &-\sup_{j\in\mathbb{N}^*}\bigg\{\int_{((\mathcal{X}\cap\mathcal{T}_1^{*})\backslash\mathcal{A}_{\delta})\times\mathbb{R}}\bar{G}_{\eps_j}(x,x+tV(x))\psi_{\eps_j}(x,t)dxdt\bigg\}  \nonumber\\
    \geq\,& -\int_{\mathcal{S}} d\lambda(T) \int_{T\times T} \mathbbm{1}_{\mathcal{X}\cap\mathcal{T}_1^{*}}(x) G_{+}(x,y)d\kappa_T'(x,y)+\int_{\mathcal{S}} d\lambda(T) \int_{T\times T} \mathbbm{1}_{\mathcal{X}\cap\mathcal{T}_1^{*}\cap\mathcal{A}_{\delta}}(x) G_{-}(x,y)d\kappa_T'(x,y)\nonumber\\
    &+\liminf_{j\rightarrow\infty}\bigg\{\int_{(\mathcal{X}\cap\mathcal{T}_1^{*})\times\mathbb{R}}\psi_{\eps_j}(x,t) \log(\psi_{\eps_j}(x,t))dxdt\bigg\}\nonumber\\
    &-\sup_{j\in\mathbb{N}^*}\bigg\{\int_{((\mathcal{X}\cap\mathcal{T}_1^{*})\backslash\mathcal{A}_{\delta})\times\mathbb{R}}\bar{G}_{\eps_j}(x,x+tV(x))\psi_{\eps_j}(x,t)dxdt\bigg\}.
\end{align}
Note that by \eqref{Eq3.4n}, the first term on the right-hand side of \eqref{Eqq3.4} is not $-\infty$. By the monotone convergence theorem and noting \eqref{P3.1.eqnew8}, we have
\begin{align}\label{P3.1.eq8n}
    &\lim_{\delta\rightarrow 0^{+}}\int_{\mathcal{S}} d\lambda(T) \int_{T\times T} \mathbbm{1}_{\mathcal{X}\cap\mathcal{T}_1^{*}\cap\mathcal{A}_{\delta}}(x) G_{-}(x,y)d\kappa_T'(x,y)\nonumber\\
    =\,& \lim_{\delta\rightarrow 0^{+}}\int_{\mathcal{S}} d\lambda(T)\int_{T\times T}  \mathbbm{1}_{\mathcal{X}\cap\mathcal{T}_1^{*}\cap\tilde{\mathcal{A}}_{\delta}\cap\mathscr{G}}(x)G_{-}(x,y)d\kappa_T'(x,y)\nonumber\\
    =\,& \int_{\mathcal{S}} d\lambda(T)\int_{T\times T}  \mathbbm{1}_{\mathcal{X}\cap\mathcal{T}_1^{*}\cap\mathscr{G}}(x)G_{-}(x,y)d\kappa_T'(x,y)=\int_{\mathcal{S}} d\lambda(T)\int_{T\times T}  \mathbbm{1}_{\mathcal{X}\cap\mathcal{T}_1^{*}}(x)G_{-}(x,y)d\kappa_T'(x,y),
\end{align}
where the last equality uses the facts that $\mathcal{L}^d(\mathcal{T}_1^*\backslash\mathscr{G})=0$ and $\gamma_0'\in\Pi(\mu,\nu)$. By \eqref{Eqq3.4}, \eqref{P3.1.eq8n}, \eqref{P3.1.eq7}, and \eqref{P3.1.eq8}, taking $\delta\rightarrow 0^{+}$, we obtain that
\begin{align}\label{Eq3.1new}
     & \liminf_{j\rightarrow\infty}\bigg\{\eps_j^{-1}\int_{\mathcal{X}\times\mathcal{Y}} \phi_{\eps_j}(x,y)\big(E(x,y)+\eps_j\log(\phi_{\eps_j}(x,y))\big)dxdy+\frac{d-1}{2}\log{\eps_j}\bigg\}\nonumber\\
    \geq\,& -\int_{\mathcal{S}} d\lambda(T) \int_{T\times T}G(x,y)d\kappa'_T(x,y)+\int_{\mathcal{S}}H(\kappa'_T|\tilde{\mu}_T\otimes\tilde{\nu}_T)d\lambda(T) \nonumber\\
    & +\int_{\mathcal{X}}f(x)\log(f(x))dx+\int_{\mathcal{S}}d\lambda(T)\int_T \tilde{g}(y)\log(\tilde{g}(y))d\mathcal{H}^1(y) .
\end{align}
By \eqref{Eq2.20} and \eqref{Eq3.1new}, we have
\begin{align*}
   & \liminf_{j\rightarrow\infty}\bigg\{\frac{\mathcal{C}_{\eps_j}(\gamma'_{\eps_j})-\OT(\mu,\nu)}{\eps_j}+\frac{d-1}{2}\log{\eps_j}\bigg\}\nonumber\\
   \geq\,& -\int_{\mathcal{S}} d\lambda(T) \int_{T\times T}G(x,y)d\kappa'_T(x,y)+\int_{\mathcal{S}}H(\kappa'_T|\tilde{\mu}_T\otimes\tilde{\nu}_T)d\lambda(T)\nonumber\\
    & -\int_{\mathcal{Y}}g(y)\log(g(y))dy+\int_{\mathcal{S}}d\lambda(T)\int_T \tilde{g}(y)\log(\tilde{g}(y))d\mathcal{H}^1(y)\nonumber\\
    =\,& \int_{\mathcal{S}} \bigg(-\frac{d-1}{2}\int_{T\times T}\log(2\pi\|x-y\|) d\kappa'_T(x,y)+ H(\kappa'_T|\tilde{\mu}_T \otimes \tilde{\nu}_T)\bigg) d\lambda(T)\nonumber\\
    & +\frac{1}{2}\int_{\mathcal{S}}d\lambda(T)\int_{T\times T}\log\det\big(\mathbf{I}_d+\|x-y\|F(y)\big)d\kappa_T'(x,y)\nonumber\\
    & -\int_{\mathcal{Y}}g(y)\log(g(y))dy+\int_{\mathcal{S}}d\lambda(T)\int_T \tilde{g}(y)\log(\tilde{g}(y))d\mathcal{H}^1(y)\nonumber\\
     =\,&\int_{\mathcal{S}} \bigg(-\frac{d-1}{2}\int_{T\times T}\log(2\pi\|x-y\|) d\kappa'_T(x,y)+ H(\kappa'_T|\tilde{\mu}_T\otimes\tilde{\nu}_T)\bigg) d\lambda(T)\nonumber\\
    &-\frac{1}{2}\int_{\mathcal{X}}f(x)\log(f(x))dx-\frac{1}{2}\int_{\mathcal{Y}}g(y)\log(g(y))dy\nonumber\\
    &+\frac{1}{2}\int_{\mathcal{S}}d\lambda(T)\int_T\tilde{f}(x)\log(\tilde{f}(x))d\mathcal{H}^1(x)+\frac{1}{2}\int_{\mathcal{S}}d\lambda(T)\int_T\tilde{g}(y)\log(\tilde{g}(y))d\mathcal{H}^1(y)\nonumber\\
    =\,& \mathtt{S}(\mu,\nu;\gamma_0'),
\end{align*}
where the first and second equalities use \cref{P-lb-2}(b) and \cref{Lem3.11}, respectively. 
\end{proof}

\section{Proof of Theorem \ref{th:upperBound}}\label{Sec4}

In this section, we present the proof of the variational upper bound. We begin by establishing several properties of the constrained EOT problem~\eqref{opt.solver} in Section~\ref{Sect.4.1}. Then, recalling the domain decomposition $\{\mathcal{X}_k\}_{k\in [K]}$, $\{\mathcal{Y}_{k'}\}_{k'\in [K']}$ introduced in Section~\ref{sectdom}, we construct in Section~\ref{Sect.4.2} the measures $\tilde{\gamma}_{\eps;k,k'}$ whose sum will yield the coupling $\tilde{\gamma}_{\eps}\in\Pi(\mu,\nu)$ announced in~\eqref{ubb.estimates} of the proof overview. Section~\ref{Sec4.0} reduces the proof of Theorem~\ref{th:upperBound} to a local estimate for each $k\in[K],k'\in[K']$, formulated as Proposition~\ref{P4.2}. This local estimate is then deduced in Section~\ref{Sect.4.3} from two key ingredients: a total mass estimate (\Cref{P4.3}) and an entropy estimate (\Cref{P4.4}). Those estimates, in turn, are proved in \cref{Sect.4.4} and~\ref{Sect.4.5}, respectively.

\subsection{Properties of the constrained EOT problem~\eqref{opt.solver}}\label{Sect.4.1}

In this subsection, we establish existence, uniqueness, and regularity properties of the constrained EOT problem~\eqref{opt.solver}. We begin by defining a set of ``good'' rays on which the induced marginals admit a monotone coupling.

\begin{definition}[The set $\tilde{\mathcal{S}}^*$]\label{defstilde}
From the proof of \cite[Theorem 6.2]{MR2011032}, for $\lambda$-a.e.\ $T\in\mathcal{S}$, there exists $\zeta_T\in\Pi(\tilde{\mu}_T,\tilde{\nu}_T)$ such that 
\begin{equation}\label{cond_Zeta}
    \zeta_T\big((\mathbb{R}^d\times\mathbb{R}^d)\backslash \{(x,y)\in T\times T: u(x)\geq u(y)\}\big)=0.
\end{equation}
As $\lambda(\mathcal{S}\backslash\mathcal{S}^*)=0$ (recall Definition \ref{factors} and Lemma \ref{L3nnn}), there exists a Borel set $\tilde{\mathcal{S}}^*\subseteq\mathcal{S}^*$, such that $\lambda(\mathcal{S}\backslash\tilde{\mathcal{S}}^*)=0$, and for every $T\in \tilde{\mathcal{S}}^*$, there exists $\zeta_T\in\Pi(\tilde{\mu}_T,\tilde{\nu}_T)$ such that \eqref{cond_Zeta} holds.
\end{definition}

Recall the definitions of $\{\mathcal{X}_k\}_{k\in [K]},\{\mathcal{Y}_{k'}\}_{k'\in [K']},a_T,b_T,V_T$ from \eqref{eq:mathtthDefn} and Definition \ref{TransportRays}. The following theorem establishes the existence, uniqueness, and regularity of the solution $\kappa_T$ to the constrained EOT problem~\eqref{opt.solver}.

\begin{theorem}\label{Lem2}
For every $T\in\tilde{\mathcal{S}}^*$, the following holds:
\begin{itemize}
     \item[(a)] There exists $\kappa_T^{(0)}\in\Pi(\tilde{\mu}_T,\tilde{\nu}_T)$, such that $\kappa_T^{(0)}\big((\mathbb{R}^d\times\mathbb{R}^d)\backslash\{(x,y)\in T\times T: u(x)\geq u(y)\}\big)=0$ and $H\big(\kappa_T^{(0)}|\tilde{\mu}_T\otimes \tilde{\nu}_T\big)<\infty$.  
     \item[(b)] The minimization problem \eqref{opt.solver} admits a unique solution $\kappa_T$, and \begin{equation*}
    -\frac{d-1}{2}\int_{T\times T}\log(\|x-y\|)d\kappa_T(x,y)+H(\kappa_T|\tilde{\mu}_T\otimes\tilde{\nu}_T)<\infty.
\end{equation*} 
    \item[(c)] $\kappa_T$ is absolutely continuous with respect to $\mathcal{H}^1|_T\otimes\mathcal{H}^1|_T$. Moreover, $\kappa_T$ admits a density $h_T(\cdot,\cdot)$ with respect to $\mathcal{H}^1|_T\otimes\mathcal{H}^1|_T$ on $T\times T$, such that $h_T(x,y)\in [0,\infty)$ for every $(x,y)\in T\times T$, $h_T(\cdot,\cdot)$ is continuous on $\mathrm{int}(T)\times\mathrm{int}(T)$, and $h_T(x,y)=0$ for any $(x,y)\in T\times T$ such that $x=a_T$ or $y=b_T$ or $u(x)\leq u(y)$. (Note that $h_T(\cdot,\cdot)$ is uniquely determined by these criteria.) 
\item[(d)] There exists a nonnegative function $\tilde{h}_T(\cdot,\cdot)$ on $T\times T$, such that $h_T(x,y)=\tilde{f}(x)\tilde{g}(y)\tilde{h}_T(x, y)$ for all $x,y\in T$, and for any $k\in[K],k'\in[K']$ and $x,x',y,y'\in T$ such that (recall \eqref{dfrakdef}) $\max\{\dist(x,\mathcal{X}_k), \dist(x',\mathcal{X}_k), \dist(y,\mathcal{Y}_{k'}),\dist(y',\mathcal{Y}_{k'})\}\leq d_0$, we have  
\begin{equation}\label{inhtild}
 \tilde{h}_T(x,y)e^{-C(\|x-x'\|+\|y-y'\|)} \leq \tilde{h}_T(x',y')\leq \tilde{h}_T(x,y) e^{C(\|x-x'\|+\|y-y'\|)}. 
\end{equation}
\end{itemize}
\end{theorem}

\begin{proof}

Throughout the proof, we fix an arbitrary $T\in\tilde{\mathcal{S}}^*$. Note that $a_T\in\mathcal{X}$ and $b_T\in\mathcal{Y}$. By \eqref{distcondi}, for every $k\in[K],k'\in[K']$, we have $\dist(\mathcal{X}_k,\mathcal{Y}_{k'})\geq 4d_0$. Moreover, $\{\mathcal{X}_k\cap T\}_{k\in[K]}$ and $\{\mathcal{Y}_{k'}\cap T\}_{k'\in [K']}$ are intervals (including the degenerate case of a single point or the empty set). Hence there exist points $0=t_0<t_1<t_2<\cdots<t_{2\EuScript{N}}=\|b_T-a_T\|$ (where $\EuScript{N}\in\mathbb{N}^{*}$) that define a partition of $T$ into $2\EuScript{N}$ disjoint intervals $\{T_{\ell}\}_{\ell\in [2\EuScript{N}]}$, where $T_{\ell}:=[a_T-t_{\ell-1}V_T,a_T-t_{\ell}V_T)$ for every $\ell\in [2\EuScript{N}-1]$ and $T_{2\EuScript{N}}:=[a_T-t_{2\EuScript{N}-1}V_T,a_T-t_{2\EuScript{N}}V_T]=[a_T-t_{2\EuScript{N}-1}V_T,b_T]$, such that the following properties hold:
\begin{itemize}
    \item[(i)] For every $k\in [K]$, $\mathcal{X}_k\cap T$ is contained in $T_{2\ell-1}$ for some $\ell\in [\EuScript{N}]$. Moreover, $\dist(a_T-t_{2\ell-1} V_T,\mathcal{X}_k\cap T)\geq 2d_0$, and if $\ell\neq 1$, then $\dist(a_T-t_{2\ell-2}V_T,\mathcal{X}_k\cap T)\geq 2d_0$.
    \item[(ii)] For every $k'\in[K']$, $\mathcal{Y}_{k'}\cap T$ is contained in $T_{2\ell}$ for some $\ell\in [\EuScript{N}]$. Moreover, $\dist(a_T-t_{2\ell-1} V_T,\mathcal{Y}_{k'}\cap T)\geq 2d_0$, and if $\ell\neq \EuScript{N}$, then $\dist(a_T-t_{2\ell}V_T,\mathcal{Y}_{k'}\cap T)\geq 2d_0$.
    \item[(iii)] For every $\ell\in [\EuScript{N}]$, $T_{2\ell-1}\cap\big(\bigcup_{k'=1}^{K'}\mathcal{Y}_{k'}\big)=\emptyset$ and $T_{2\ell}\cap\big(\bigcup_{k=1}^{K}\mathcal{X}_{k}\big)=\emptyset$.
\end{itemize}
Note that property (iii) above implies that 
\begin{equation}\label{dich}
    \text{Either } T_{\ell}\cap\mathcal{X}=\emptyset \text{ or } T_{\ell}\cap\mathcal{Y}=\emptyset, \quad \text{ for every }\ell\in [ 2\EuScript{N}].
\end{equation} 

Let $\zeta_T\in\Pi(\tilde{\mu}_T,\tilde{\nu}_T)$ be such that \eqref{cond_Zeta} holds. For any $x,y\in T$, we define
\begin{equation}\label{defh0}
    h^{(0)}_T(x,y):=\sum_{\substack{1\leq \ell<\ell' \leq 2\EuScript{N}:\\ \tilde{\mu}_T(T_{\ell})>0,\tilde{\nu}_T(T_{\ell'})>0}} \frac{\zeta_T(T_{\ell}\times T_{\ell'})}{\tilde{\mu}_T(T_{\ell})\tilde{\nu}_T(T_{\ell'})}\cdot \mathbbm{1}_{T_{\ell}}(x)\tilde{f}(x)\mathbbm{1}_{T_{\ell'}}(y)\tilde{g}(y).
\end{equation}
Let $\kappa_T^{(0)}$ be the Borel measure on $\mathbb{R}^d\times\mathbb{R}^d$ with density $h^{(0)}_T(\cdot,\cdot)$ on $T\times T$ (and $0$ outside) with respect to $\mathcal{H}^1|_T\otimes \mathcal{H}^1|_T$. By \eqref{dich}, $\zeta_T(T_{\ell}\times T_{\ell})=0$ for every $\ell\in [ 2\EuScript{N}]$. Moreover, for any $\ell,\ell'\in [ 2\EuScript{N}]$ such that $\ell'< \ell$, we have $u(x)<u(y)$ for every $x\in T_{\ell},y\in T_{\ell'}$, and so $\zeta_T(T_{\ell}\times T_{\ell'})=0$ by \eqref{cond_Zeta}. Therefore, for any Borel set $A\subseteq\mathbb{R}^d$, we have (recall \eqref{ress})
\begin{align}\label{mar}
    &\kappa_T^{(0)}(A\times\mathbb{R}^d)=\sum_{\substack{1\leq \ell<\ell' \leq  2\EuScript{N}:\\ \tilde{\mu}_T(T_{\ell})>0,\tilde{\nu}_T(T_{\ell'})>0}} \frac{\zeta_T(T_{\ell}\times T_{\ell'})}{\tilde{\mu}_T(T_{\ell})\tilde{\nu}_T(T_{\ell'})}\cdot \tilde{\mu}_T(T_{\ell}\cap A)\tilde{\nu}_T(T_{\ell'})\nonumber\\
    =& \sum_{\ell\in[ 2\EuScript{N}]:\tilde{\mu}_T(T_{\ell})>0}\frac{\tilde{\mu}_T(T_{\ell}\cap A)}{\tilde{\mu}_T(T_{\ell})}\cdot
    \sum_{\ell< \ell'\leq 2\EuScript{N}}\zeta_T(T_{\ell}\times T_{\ell'})\nonumber\\
    =&\sum_{\ell\in[ 2\EuScript{N}]:\tilde{\mu}_T(T_{\ell})>0}\frac{\tilde{\mu}_T(T_{\ell}\cap A)}{\tilde{\mu}_T(T_{\ell})}\cdot
    \sum_{\ell'\in[ 2\EuScript{N}]}\zeta_T(T_{\ell}\times T_{\ell'})\nonumber\\
    =&\sum_{\ell\in[ 2\EuScript{N}]:\tilde{\mu}_T(T_{\ell})>0}\frac{\tilde{\mu}_T(T_{\ell}\cap A)}{\tilde{\mu}_T(T_{\ell})}\cdot \tilde{\mu}_T(T_{\ell})=\sum_{\ell\in[ 2\EuScript{N}]}\tilde{\mu}_T(T_{\ell}\cap A)=\tilde{\mu}_T(A).
\end{align}
Similarly, we can deduce that $\kappa_T^{(0)}(\mathbb{R}^d\times A)=\tilde{\nu}_T(A)$ for any Borel set $A\subseteq\mathbb{R}^d$. Hence $\kappa_T^{(0)}\in\Pi(\tilde{\mu}_T,\tilde{\nu}_T)$. By \eqref{defh0}, $\kappa_T^{(0)}\big((\mathbb{R}^d\times\mathbb{R}^d)\backslash\{(x,y)\in T\times T: u(x)\geq u(y)\}\big)=0$. Moreover,
\begin{align}\label{ent}
    H\big(\kappa_T^{(0)}|\tilde{\mu}_T\otimes \tilde{\nu}_T\big)=&\sum_{\substack{1\leq \ell<\ell' \leq 2\EuScript{N}:\\ \tilde{\mu}_T(T_{\ell})>0,\tilde{\nu}_T(T_{\ell'})>0}} \frac{\zeta_T(T_{\ell}\times T_{\ell'})}{\tilde{\mu}_T(T_{\ell})\tilde{\nu}_T(T_{\ell'})} \log\bigg(\frac{\zeta_T(T_{\ell}\times T_{\ell'})}{\tilde{\mu}_T(T_{\ell})\tilde{\nu}_T(T_{\ell'})}\bigg)\nonumber\\
    &\hspace{1in}\times\int_{T\times T}\mathbbm{1}_{T_{\ell}}(x)\tilde{f}(x)\mathbbm{1}_{T_{\ell'}}(y)\tilde{g}(y)d\mathcal{H}^1(x)d\mathcal{H}^1(y)\nonumber\\
    =&\sum_{\substack{1\leq \ell<\ell' \leq 2\EuScript{N}:\\ \tilde{\mu}_T(T_{\ell})>0,\tilde{\nu}_T(T_{\ell'})>0}} \zeta_T(T_{\ell}\times T_{\ell'})\log\bigg(\frac{\zeta_T(T_{\ell}\times T_{\ell'})}{\tilde{\mu}_T(T_{\ell})\tilde{\nu}_T(T_{\ell'})}\bigg)<\infty,
\end{align}
where we note that $\zeta_T(T_{\ell}\times T_{\ell'})\leq 1$ for every $\ell,\ell'\in[ 2\EuScript{N}]$ in the last inequality. This establishes part (a) of the theorem.

Now we proceed to the proofs of parts (b)--(d). By Definition \ref{De3.8},
\begin{equation}\label{supports}
    \tilde{\mu}_T\big(\mathbb{R}^d\backslash (T\cap\mathcal{X})\big)=\tilde{\nu}_T\big(\mathbb{R}^d\backslash (T\cap\mathcal{Y})\big)=0.
\end{equation}
For any $(x,y)\in T\times T$, we define
\begin{equation}\label{def_cT}
    c_T(x,y):=\begin{cases}
        -\frac{d-1}{2}\log(\|x-y\|), & \text{ if }u(x)> u(y),\\
        +\infty, & \text{ if }u(x)\leq u(y).
    \end{cases}
\end{equation}
Note that by \eqref{supports} and \eqref{domainconst}, 
\begin{align*}
    \int_{T\times T}&e^{-c_T(x,y)}d\tilde{\mu}_T\otimes\tilde{\nu}_T(x,y)=\int_{\{(x,y)\in (T\cap\mathcal{X})\times (T\cap\mathcal{Y}): u(x)>u(y)\}}e^{\frac{d-1}{2}\log(\|x-y\|)}d\tilde{\mu}_T\otimes\tilde{\nu}_T(x,y)\nonumber\\
&\begin{cases}
    \leq (2D)^{(d-1)\slash 2}\int_{\{(x,y)\in (T\cap\mathcal{X})\times (T\cap\mathcal{Y}): u(x)>u(y)\}}d\tilde{\mu}_T\otimes\tilde{\nu}_T(x,y)\leq (2D)^{(d-1)\slash 2}\leq C,\\
    \geq  d_0^{(d-1)\slash 2}\tilde{\mu}_T\otimes\tilde{\nu}_T\big(\{(x,y)\in (T\cap\mathcal{X})\times (T\cap\mathcal{Y}): u(x)>u(y)\}\big).
\end{cases}
\end{align*}
We have shown previously that $\kappa_T^{(0)}\in\Pi(\tilde{\mu}_T,\tilde{\nu}_T)$; in particular, $\kappa_T^{(0)}(T\times T)=1$ and by \eqref{defh0}, there exist $\ell,\ell'\in[2\EuScript{N}]$ such that $\ell<\ell'$ and $\tilde{\mu}_T(T_{\ell})>0,\tilde{\nu}_T(T_{\ell'})>0$, hence
\begin{equation*}
    \tilde{\mu}_T\otimes\tilde{\nu}_T\big(\{(x,y)\in (T\cap\mathcal{X})\times (T\cap\mathcal{Y}): u(x)>u(y)\}\big)\geq \tilde{\mu}_T\otimes\tilde{\nu}_T(T_{\ell}\times T_{\ell'})=\tilde{\mu}_T(T_{\ell})\tilde{\nu}_T(T_{\ell'})>0.
\end{equation*}
Therefore, $\int_{T\times T}e^{-c_T(x,y)}d\tilde{\mu}_T\otimes\tilde{\nu}_T(x,y)\in (0,\infty)$. Define $\EuScript{M}_T$ to be the Borel probability measure on $\mathbb{R}^d\times\mathbb{R}^d$ with density $\frac{e^{-c_T(\cdot,\cdot)}}{\int_{T\times T}e^{-c_T(x,y)}d\tilde{\mu}_T\otimes\tilde{\nu}_T(x,y)}$ on $T\times T$ (and $0$ outside) with respect to $\tilde{\mu}_T\otimes\tilde{\nu}_T$. 

Now we consider any Borel set $A\subseteq\mathbb{R}^d\times\mathbb{R}^d$. By \eqref{supports} and \eqref{domainconst},
\begin{align}\label{res_1}
   & \int_{T\times T}\mathbbm{1}_A(x,y) e^{-c_T(x,y)}d\tilde{\mu}_T\otimes\tilde{\nu}_T(x,y)\nonumber\\
   =& \int_{A\cap\{(x,y)\in (T\cap\mathcal{X})\times (T\cap\mathcal{Y}): u(x)>u(y)\}}e^{\frac{d-1}{2}\log(\|x-y\|)}d\tilde{\mu}_T\otimes\tilde{\nu}_T(x,y)\nonumber\\
   \geq& d_0^{(d-1)\slash 2}\tilde{\mu}_T\otimes\tilde{\nu}_T\big(A \cap \{(x,y)\in (T\cap\mathcal{X})\times (T\cap\mathcal{Y}): u(x)>u(y)\}\big)\nonumber\\
   =& d_0^{(d-1)\slash 2}\tilde{\mu}_T\otimes\tilde{\nu}_T\big(A \cap \{(x,y)\in T \times T: u(x)>u(y)\}\big)\nonumber\\
   =& d_0^{(d-1)\slash 2} \sum_{1\leq \ell<\ell'\leq 2\EuScript{N}}\tilde{\mu}_T\otimes \tilde{\nu}_T(A\cap (T_{\ell}\times T_{\ell'})),
\end{align}
where we use \eqref{dich} in the last line. If $\EuScript{M}_T(A)=0$, then $\int_{T\times T}\mathbbm{1}_A(x,y) e^{-c_T(x,y)}d\tilde{\mu}_T\otimes\tilde{\nu}_T(x,y)=0$, and by \eqref{res_1}, for every $1\leq \ell<\ell'\leq 2\EuScript{N}$, we have $\tilde{\mu}_T\otimes\tilde{\nu}_T(A\cap (T_{\ell}\times T_{\ell'}))=0$. Hence noting \eqref{defh0}, we have
\begin{equation*}
    \kappa_T^{(0)}(A)=\sum_{\substack{1\leq \ell<\ell' \leq 2\EuScript{N}:\\ \tilde{\mu}_T(T_{\ell})>0,\tilde{\nu}_T(T_{\ell'})>0}} \frac{\zeta_T(T_{\ell}\times T_{\ell'})}{\tilde{\mu}_T(T_{\ell})\tilde{\nu}_T(T_{\ell'})}\cdot \tilde{\mu}_T\otimes\tilde{\nu}_T(A\cap (T_{\ell}\times T_{\ell'}))=0. 
\end{equation*}
Therefore, $\kappa_T^{(0)}$ is absolutely continuous with respect to $\EuScript{M}_T$. 

Now note that for any $\chi_T\in\Pi(\tilde{\mu}_T,\tilde{\nu}_T)$, if $\chi_T\big((\mathbb{R}^d\times\mathbb{R}^d)\backslash\{(x,y)\in T\times T: u(x)\geq u(y)\}\big)>0$, then $\chi_T$ is not absolutely continuous with respect to $\EuScript{M}_T$ (as $\EuScript{M}_T\big((\mathbb{R}^d\times\mathbb{R}^d)\backslash\{(x,y)\in T\times T: u(x)\geq u(y)\}\big)=0$), hence $H(\chi_T|\EuScript{M}_T)=\infty$; if $\chi_T\big((\mathbb{R}^d\times\mathbb{R}^d)\backslash\{(x,y)\in T\times T: u(x)\geq u(y)\}\big)=0$, then 
\begin{align}\label{relH}
     H(\chi_T|\EuScript{M}_T)=&\log\bigg(\int_{T\times T}e^{-c_T(x,y)}d\tilde{\mu}_T\otimes\tilde{\nu}_T(x,y)\bigg)+\int_{T\times T}c_T(x,y)d\chi_T(x,y)+H(\chi_T|\tilde{\mu}_T\otimes\tilde{\nu}_T)\nonumber\\
    =&\log\bigg(\int_{T\times T}e^{-c_T(x,y)}d\tilde{\mu}_T\otimes\tilde{\nu}_T(x,y)\bigg)\nonumber\\
    &-\frac{d-1}{2}\int_{T\times T}\log(\|x-y\|)d\chi_T(x,y)+H(\chi_T|\tilde{\mu}_T\otimes\tilde{\nu}_T).
\end{align}
Hence the minimization problem \eqref{opt.solver} is equivalent to the minimization problem
\begin{equation}\label{new.solver}
    \min_{\chi_T\in\Pi(\tilde{\mu}_T,\tilde{\nu}_T)} H(\chi_T|\EuScript{M}_T). 
\end{equation}
Note that by \eqref{domainconst},
\begin{align*}
    \int_{T\times T}c_T(x,y)d\kappa_T^{(0)}(x,y)=&-\int_{\{(x,y)\in (T\cap\mathcal{X})\times (T\cap\mathcal{Y}): u(x)>u(y)\}}\frac{d-1}{2}\log(\|x-y\|)d\kappa_T^{(0)}(x,y)\nonumber\\
    \leq& -\frac{d-1}{2}\log d_0<\infty.
\end{align*}
Hence by part (a) and \eqref{relH}, we have $H\big(\kappa^{(0)}_T|\EuScript{M}_T\big)<\infty$. By \cite[Theorem 2.1]{Nutz.20}, the minimization problem \eqref{new.solver} admits a unique solution $\kappa_T$ and $H(\kappa_T|\EuScript{M}_T)<\infty$. Hence the minimization problem \eqref{opt.solver} also admits a unique solution (which is $\kappa_T$), and 
\begin{equation*}
    -\frac{d-1}{2}\int_{T\times T}\log(\|x-y\|)d\kappa_T(x,y)+H(\kappa_T|\tilde{\mu}_T\otimes\tilde{\nu}_T)<\infty.
\end{equation*}
This establishes part (b). In particular, $H(\kappa_T|\tilde{\mu}_T\otimes\tilde{\nu}_T)<\infty$ and $\kappa_T$ is absolutely continuous with respect to $\tilde{\mu}_T\otimes\tilde{\nu}_T$. As $\tilde{\mu}_T\otimes\tilde{\nu}_T$ is absolutely continuous with respect to $\mathcal{H}^1|_T\otimes\mathcal{H}^1|_T$ (see Definition~\ref{De3.8}), we conclude that $\kappa_T$ is absolutely continuous with respect to $\mathcal{H}^1|_T\otimes\mathcal{H}^1|_T$. 

We define $\EuScript{L}_T$ to be the set of $(\ell,\ell')$ such that $1\leq \ell<\ell'\leq  2\EuScript{N}$ and $\kappa_T(T_{\ell}\times T_{\ell'})>0$. Note that for any $(\ell,\ell')\in \EuScript{L}_T$, $\ell$ is odd and $\ell'$ is even (otherwise $\tilde{\mu}_T(T_{\ell})=0$ or $\tilde{\nu}_T(T_{\ell'})=0$, which implies $\kappa_T(T_{\ell}\times T_{\ell'})=0$). For any $(x,y)\in T\times T$, we define
\begin{equation}
    \tilde{c}_T(x,y):=\begin{cases}
        -\frac{d-1}{2}\log(\|x-y\|), & \text{ if }(x,y)\in\bigcup_{(\ell,\ell')\in\EuScript{L}_T}(T_{\ell}\times T_{\ell'}),\\
        +\infty, & \text{ otherwise}.
    \end{cases}
\end{equation}
Note that $\tilde{c}_T(x,y)\geq c_T(x,y)$ for every $(x,y)\in T\times T$. Hence for any $\chi_T\in\Pi(\tilde{\mu}_T,\tilde{\nu}_T)$ such that $\chi_T\big((\mathbb{R}^d\times\mathbb{R}^d)\backslash\{(x,y)\in T\times T: u(x)\geq u(y)\}\big)=0$, we have  
\begin{align*}
   & \int_{T\times T}\tilde{c}_T(x,y)d\chi_T(x,y)+H(\chi_T|\tilde{\mu}_T\otimes\tilde{\nu}_T)\geq \int_{T\times T}c_T(x,y)d\chi_T(x,y)+H(\chi_T|\tilde{\mu}_T\otimes\tilde{\nu}_T)\nonumber\\
   =& -\frac{d-1}{2} \int_{T\times T}\log(\|x-y\|)d\chi_T(x,y)+H(\chi_T|\tilde{\mu}_T\otimes\tilde{\nu}_T)\nonumber\\
   \geq& -\frac{d-1}{2} \int_{T\times T}\log(\|x-y\|)d\kappa_T(x,y)+H(\kappa_T|\tilde{\mu}_T\otimes\tilde{\nu}_T)\nonumber\\
   =& -\frac{d-1}{2}\sum_{1\leq \ell<\ell'\leq  2\EuScript{N}}\int_{T_{\ell}\times T_{\ell'}}\log(\|x-y\|)d\kappa_T(x,y) +H(\kappa_T|\tilde{\mu}_T\otimes\tilde{\nu}_T)\nonumber\\ 
   =& -\frac{d-1}{2}\sum_{(\ell,\ell')\in\EuScript{L}_T}\int_{T_{\ell}\times T_{\ell'}}\log(\|x-y\|)d\kappa_T(x,y) +H(\kappa_T|\tilde{\mu}_T\otimes\tilde{\nu}_T)\nonumber\\ 
   =& \int_{T\times T}\tilde{c}_T(x,y) d\kappa_T(x,y)+H(\kappa_T|\tilde{\mu}_T\otimes\tilde{\nu}_T),
\end{align*}
where we use the fact that $\kappa_T$ is the solution of \eqref{opt.solver} in the third line and the definition of $\EuScript{L}_T$ in the fifth line. Hence $\kappa_T$ is the solution to the following minimization problem:
\begin{equation}\label{opt.solver1}
    \min_{\substack{\chi_T\in\Pi(\tilde{\mu}_T,\tilde{\nu}_T):\\ \chi_T((\mathbb{R}^d\times\mathbb{R}^d)\backslash\{(x,y)\in T\times T: u(x)\geq u(y)\})=0}}\bigg\{ \int_{T\times T} \tilde{c}_T(x,y) d\chi_T(x,y) + H(\chi_T|\tilde{\mu}_T\otimes \tilde{\nu}_T)\bigg\}.
\end{equation}
Note that by \eqref{supports} and \eqref{domainconst}, 
\begin{align*}
    \int_{T\times T}&e^{-\tilde{c}_T(x,y)}d\tilde{\mu}_T\otimes\tilde{\nu}_T(x,y)=\sum_{(\ell,\ell')\in\EuScript{L}_T}\int_{(T_{\ell}\cap\mathcal{X})\times (T_{\ell'}\cap\mathcal{Y})}e^{\frac{d-1}{2}\log(\|x-y\|)}d\tilde{\mu}_T\otimes\tilde{\nu}_T(x,y) \nonumber\\
&\begin{cases}
    \leq (2D)^{(d-1)\slash 2}\sum_{(\ell,\ell')\in\EuScript{L}_T} \tilde{\mu}_T\otimes\tilde{\nu}_T\big((T_{\ell}\cap\mathcal{X})\times (T_{\ell'}\cap\mathcal{Y})\big)\leq (2D)^{(d-1)\slash 2}\leq C,\\
    \geq  d_0^{(d-1)\slash 2}\sum_{(\ell,\ell')\in\EuScript{L}_T} \tilde{\mu}_T\otimes\tilde{\nu}_T\big((T_{\ell}\cap\mathcal{X})\times (T_{\ell'}\cap\mathcal{Y})\big)=d_0^{(d-1)\slash 2}\sum_{(\ell,\ell')\in\EuScript{L}_T} \tilde{\mu}_T(T_{\ell})\tilde{\nu}_T(T_{\ell'}).
\end{cases}
\end{align*}
By \eqref{dich}, we have $1=\kappa_T(T\times T)=\sum_{(\ell,\ell')\in\EuScript{L}_T}\kappa_T(T_{\ell}\times T_{\ell'})$. Hence $\EuScript{L}_T\neq\emptyset$. For any $(\ell,\ell')\in\EuScript{L}_T$, as $\kappa_T\in\Pi(\tilde{\mu}_T,\tilde{\nu}_T)$, we have
\begin{equation}\label{pso}
    \min\{\tilde{\mu}_T(T_{\ell}),\tilde{\nu}_T(T_{\ell'})\}\geq \kappa_T(T_{\ell}\times T_{\ell'})>0.
\end{equation}
Therefore, we have $\int_{T\times T}e^{-\tilde{c}_T(x,y)}d\tilde{\mu}_T\otimes\tilde{\nu}_T(x,y)\in (0,\infty)$. Now define $\tilde{\EuScript{M}}_T$ to be the Borel probability measure on $\mathbb{R}^d\times\mathbb{R}^d$ with density $\frac{e^{-\tilde{c}_T(\cdot,\cdot)}}{\int_{T\times T}e^{-\tilde{c}_T(x,y)}d\tilde{\mu}_T\otimes\tilde{\nu}_T(x,y)}$ on $T\times T$ (and $0$ outside) with respect to $\tilde{\mu}_T\otimes\tilde{\nu}_T$. It can be checked that \eqref{opt.solver1} is equivalent to the following minimization problem: 
\begin{equation}\label{opt.solver1_1}
    \min_{\chi_T\in\Pi(\tilde{\mu}_T,\tilde{\nu}_T)}H(\chi_T|\tilde{\EuScript{M}}_T). 
\end{equation}
Hence $\kappa_T$ is the solution to \eqref{opt.solver1_1}.

For any $x,y\in T$, we define (note \eqref{pso})
\begin{equation}
    \EuScript{H}_T(x,y):=\sum_{(\ell,\ell')\in\EuScript{L}_T} \frac{\kappa_T(T_{\ell}\times T_{\ell'})}{\tilde{\mu}_T(T_{\ell})\tilde{\nu}_T(T_{\ell'})}\cdot \mathbbm{1}_{T_{\ell}}(x)\tilde{f}(x)\mathbbm{1}_{T_{\ell'}}(y)\tilde{g}(y).
\end{equation}
Let $\tilde{\kappa}_T$ be the Borel measure on $\mathbb{R}^d\times\mathbb{R}^d$ with density $\EuScript{H}_T(\cdot,\cdot)$ on $T\times T$ (and $0$ outside) with respect to $\mathcal{H}^1|_T\otimes\mathcal{H}^1|_T$. Arguing similarly as in \eqref{mar} and \eqref{ent}, we get $\tilde{\kappa}_T\in\Pi(\tilde{\mu}_T,\tilde{\nu}_T)$ and $H(\tilde{\kappa}_T|\tilde{\mu}_T\otimes\tilde{\nu}_T)<\infty$, from which we can deduce that $H(\tilde{\kappa}_T|\tilde{\EuScript{M}}_T)<\infty$. 

Now consider any Borel set $A\subseteq\mathbb{R}^d\times\mathbb{R}^d$. By \eqref{supports} and \eqref{domainconst},
\begin{align}\label{Eqq1.1}
    &\int_{T\times T}\mathbbm{1}_A(x,y) e^{-\tilde{c}_T(x,y)}d\tilde{\mu}_T\otimes\tilde{\nu}_T(x,y)\nonumber\\
    =& \sum_{(\ell,\ell')\in\EuScript{L}_T}\int_{(T_{\ell}\cap\mathcal{X})\times (T_{\ell'}\cap\mathcal{Y})}\mathbbm{1}_A(x,y) e^{\frac{d-1}{2}\log(\|x-y\|)}d\tilde{\mu}_T\otimes\tilde{\nu}_T(x,y)\nonumber\\
    &\begin{cases}
        \leq (2D)^{(d-1)  \slash 2}\sum_{(\ell,\ell')\in  \EuScript{L}_T}\tilde{\mu}_T\otimes\tilde{\nu}_T(A\cap(T_{\ell}\times T_{\ell'})),\\
        \geq d_0^{(d-1)\slash 2} \sum_{(\ell,\ell')\in  \EuScript{L}_T}\tilde{\mu}_T\otimes\tilde{\nu}_T(A\cap(T_{\ell}\times T_{\ell'})). 
    \end{cases}
\end{align}
Moreover, 
\begin{equation}\label{Eqq1.1.1}
    \tilde{\kappa}_T(A)=\sum_{(\ell,\ell')\in\EuScript{L}_T} \frac{\kappa_T(T_{\ell}\times T_{\ell'})}{\tilde{\mu}_T(T_{\ell})\tilde{\nu}_T(T_{\ell'})} \cdot   \tilde{\mu}_T\otimes\tilde{\nu}_T(A\cap (T_{\ell}\times T_{\ell'}))
\end{equation}
If $\tilde{\EuScript{M}}_T(A)=0$, then $\int_{T\times T}\mathbbm{1}_A(x,y) e^{-\tilde{c}_T(x,y)}d\tilde{\mu}_T\otimes\tilde{\nu}_T(x,y)=0$, and by \eqref{Eqq1.1}, for every $(\ell,\ell')\in\EuScript{L}_T$, $\tilde{\mu}_T\otimes \tilde{\nu}_T(A\cap (T_{\ell}\times T_{\ell'}))=0$. Hence by \eqref{Eqq1.1.1}, $\tilde{\kappa}_T(A)=0$. Now if $\tilde{\kappa}_T(A)=0$, as $\frac{\kappa_T(T_{\ell}\times T_{\ell'})}{\tilde{\mu}_T(T_{\ell})\tilde{\nu}_T(T_{\ell'})}>0$ for every $(\ell,\ell')\in\EuScript{L}_T$, we have $\tilde{\mu}_T\otimes \tilde{\nu}_T(A\cap (T_{\ell}\times T_{\ell'}))=0$ for every $(\ell,\ell')\in\EuScript{L}_T$. Hence by \eqref{Eqq1.1}, $\tilde{\EuScript{M}}_T(A)=0$. Therefore, the measures $\tilde{\kappa}_T$ and $\tilde{\EuScript{M}}_T$ are equivalent. 

As $\kappa_T$ is the solution to \eqref{opt.solver1_1}, by \cite[Equation (3.4) and Theorem 3.43]{FG} (note that $\tilde{\kappa}_T\in\Pi(\tilde{\mu}_T,\tilde{\nu}_T)$, $H(\tilde{\kappa}_T|\tilde{\EuScript{M}}_T)<\infty$, $\tilde{\kappa}_T$ and $\tilde{\EuScript{M}}_T$ are equivalent, and $\tilde{\EuScript{M}}_T$ is absolutely continuous with respect to $\tilde{\mu}_T\otimes\tilde{\nu}_T$), the measures $\kappa_T$ and $\tilde{\EuScript{M}}_T$ are equivalent, and there exist Borel measurable functions $\varphi_T,\psi_T:T\rightarrow [0,\infty)$, such that $\frac{d\kappa_T}{d\tilde{\EuScript{M}}_T}(x,y)=\varphi_T(x)\psi_T(y)$ for $\kappa_T$-a.e.\ (hence $\tilde{\EuScript{M}}_T$-a.e.) $(x,y)\in T\times T$. Therefore, for $\tilde{\mu}_T\otimes\tilde{\nu}_T$-a.e.\ $(x,y)\in\bigcup_{(\ell,\ell')\in\EuScript{L}_T}(T_{\ell}\times T_{\ell'})$, we have
\begin{equation*}
     \frac{d\kappa_T}{d\tilde{\mu}_T\otimes\tilde{\nu}_T}(x,y)=e^{\frac{d-1}{2}\log(\|x-y\|)} \bar{\varphi}_T(x)  \bar{\psi}_T(y)=\|x-y\|^{(d-1)\slash 2}\bar{\varphi}_T(x)  \bar{\psi}_T(y),
\end{equation*}
where $\bar{\varphi}_T,\bar{\psi}_T:T\rightarrow [0,\infty)$ are Borel measurable functions. As $\kappa_T$ is concentrated on $\bigcup_{(\ell,\ell')\in\EuScript{L}_T}(T_{\ell}\times T_{\ell'})$, for $\tilde{\mu}_T\otimes\tilde{\nu}_T$-a.e.\ $(x,y)\in T\times T$, we have 
\begin{equation}\label{factorize}
    \frac{d\kappa_T}{d\tilde{\mu}_T\otimes\tilde{\nu}_T}(x,y)=\mathbbm{1}_{\bigcup_{(\ell,\ell')\in\EuScript{L}_T}(T_{\ell}\times T_{\ell'})}(x,y)\|x-y\|^{(d-1)\slash 2}\bar{\varphi}_T(x)  \bar{\psi}_T(y).
\end{equation}
As $\kappa_T\in\Pi(\tilde{\mu}_T,\tilde{\nu}_T)$, for every $\ell\in [\EuScript{N}]$, we have
\begin{equation}\label{Eqqq1}
    \int_{\bigcup_{\ell'\in [\EuScript{N}]:(2\ell-1,2\ell')\in\EuScript{L}_T}T_{2\ell'}}\|x-y\|^{(d-1)\slash 2}\bar{\varphi}_T(x)\bar{\psi}_T(y)\tilde{g}(y)d\mathcal{H}^1(y)=1 ,\quad\text{ for }\tilde{\mu}_T\text{-a.e.\ }x\in T_{2\ell-1}, 
\end{equation}
\begin{equation}\label{Eqqq1.1}
    \int_{\bigcup_{\ell'\in [\EuScript{N}]:(2\ell'-1,2\ell)\in\EuScript{L}_T}T_{2\ell'-1}}\|x-y\|^{(d-1)\slash 2}\bar{\varphi}_T(x)\bar{\psi}_T(y)\tilde{f}(x)d\mathcal{H}^1(x)=1 ,\quad\text{ for }\tilde{\nu}_T\text{-a.e.\ }y\in T_{2\ell}.
\end{equation}

Recall that $T_{\ell}=[a_T-t_{\ell-1}V_T,a_T-t_{\ell}V_T)$ for every $\ell\in [2\EuScript{N}-1]$ and $T_{2\EuScript{N}}=[a_T-t_{2\EuScript{N}-1}V_T,b_T]$. Also recall the properties satisfied by $\{T_{\ell}\}_{\ell\in [2\EuScript{N}]}$. Note that by \eqref{supports} and \eqref{domainconst}, for any $\ell\in [\EuScript{N}]$, $\int_{\bigcup_{\ell'\in [\EuScript{N}]:(2\ell-1,2\ell')\in \EuScript{L}_T}T_{2\ell'}}\bar{\psi}_T(y)d\tilde{\nu}_T(y)>0$ if and only if
\begin{align*}
&\int_{\bigcup_{\ell'\in [\EuScript{N}]:(2\ell-1,2\ell')\in \EuScript{L}_T}T_{2\ell'}}\|x-y\|^{(d-1)\slash 2}\bar{\psi}_T(y)d\tilde{\nu}_T(y)\nonumber\\
=&\int_{\big(\bigcup_{\ell'\in [\EuScript{N}]:(2\ell-1,2\ell')\in \EuScript{L}_T}T_{2\ell'}\big)\cap\mathcal{Y}}\|x-y\|^{(d-1)\slash 2}\bar{\psi}_T(y)d\tilde{\nu}_T(y)>0
\end{align*}
for all $x\in [a_T-t_{2\ell-2}V_T,a_T-t_{2\ell-1}V_T]$,
and $\int_{\bigcup_{\ell'\in [\EuScript{N}]:(2\ell'-1,2\ell)\in \EuScript{L}_T}T_{2\ell'-1}}\bar{\varphi}_T(x)d\tilde{\mu}_T(x)>0$ if and only if
\begin{align*}
  &\int_{\bigcup_{\ell'\in [\EuScript{N}]:(2\ell'-1,2\ell)\in \EuScript{L}_T}T_{2\ell'-1}}\|x-y\|^{(d-1)\slash 2}\bar{\varphi}_T(x)d\tilde{\mu}_T(x)\nonumber\\
  =&\int_{\big(\bigcup_{\ell'\in [\EuScript{N}]:(2\ell'-1,2\ell)\in \EuScript{L}_T}T_{2\ell'-1}\big)\cap\mathcal{X}}\|x-y\|^{(d-1)\slash 2}\bar{\varphi}_T(x)d\tilde{\mu}_T(x)>0
\end{align*}
for all $y\in[a_T-t_{2\ell-1}V_T, a_T-t_{2\ell}V_T]$.

For any $\ell\in [\EuScript{N}]$ and $x\in [a_T-t_{2\ell-2}V_T,a_T-t_{2\ell-1}V_T]$, we define 
\begin{align*}
    \tilde{\varphi}_T(x):=\begin{cases}
        \frac{1}{\int_{\bigcup_{\ell'\in [\EuScript{N}]:(2\ell-1,2\ell')\in \EuScript{L}_T}T_{2\ell'}}\|x-y\|^{(d-1)\slash 2}\bar{\psi}_T(y)d\tilde{\nu}_T(y)}, & \text{ if denominator $>0$}, \\
        0, & \text{ otherwise}.
    \end{cases}
\end{align*}
Moreover, we define $\tilde{\varphi}_T(b_T):=0$. For any $\ell\in [\EuScript{N}]$ and $x\in (a_T-t_{2\ell-1}V_T, a_T-t_{2\ell}V_T)$, we define
\begin{equation*}
    \tilde{\varphi}_T(x):=\frac{\|x-(a_T-t_{2\ell}V_T)\|\tilde{\varphi}_T(a_T-t_{2\ell-1}V_T)+\|x-(a_T-t_{2\ell-1}V_T)\|\tilde{\varphi}_T(a_T-t_{2\ell}V_T)}{t_{2\ell}-t_{2\ell-1}}.
\end{equation*}

For any $\ell\in [\EuScript{N}]$ and $y\in [a_T-t_{2\ell-1}V_T, a_T-t_{2\ell}V_T]$, we define
\begin{equation*}
    \tilde{\psi}_T(y):=\begin{cases}
        \frac{1}{\int_{\bigcup_{\ell'\in [\EuScript{N}]:(2\ell'-1,2\ell)\in \EuScript{L}_T}T_{2\ell'-1}}\|x-y\|^{(d-1)\slash 2}\bar{\varphi}_T(x)d\tilde{\mu}_T(x)}, & \text{ if denominator $>0$}, \\
        0, & \text{ otherwise}.
    \end{cases}
\end{equation*}
Moreover, we define $\tilde{\psi}_T(a_T):=0$. For any $\ell\in [\EuScript{N}]$ and $y\in (a_T-t_{2\ell-2}V_T, a_T-t_{2\ell-1}V_T)$, we define
\begin{equation*}
    \tilde{\psi}_T(y):=\frac{\|y-(a_T-t_{2\ell-1}V_T)\|\tilde{\psi}_T(a_T-t_{2\ell-2}V_T)+\|y-(a_T-t_{2\ell-2}V_T)\|\tilde{\psi}_T(a_T-t_{2\ell-1}V_T)}{t_{2\ell-1}-t_{2\ell-2}}.
\end{equation*}

For any $\ell\in [\EuScript{N}]$, if $\int_{\bigcup_{\ell'\in [\EuScript{N}]:(2\ell-1,2\ell')\in \EuScript{L}_T}T_{2\ell'}}\bar{\psi}_T(y)d\tilde{\nu}_T(y)=0$, then we have
\begin{equation*}
    \int_{\bigcup_{\ell'\in [\EuScript{N}]:(2\ell-1,2\ell')\in\EuScript{L}_T}T_{2\ell'}}\|x-y\|^{(d-1)\slash 2}\bar{\psi}_T(y)\tilde{g}(y)d\mathcal{H}^1(y)=0, \quad\text{for every }x\in T_{2\ell-1}.
\end{equation*}
Hence by \eqref{Eqqq1}, we get $0=1$ for $\tilde{\mu}_T$-a.e.\ $x\in T_{2\ell-1}$, which implies $\tilde{\mu}_T(T_{2\ell-1})=0$. Similarly, if $\int_{\bigcup_{\ell'\in [\EuScript{N}]:(2\ell'-1,2\ell)\in \EuScript{L}_T}T_{2\ell'-1}}\bar{\varphi}_T(x)d\tilde{\mu}_T(x)=0$, then $\tilde{\nu}_T(T_{2\ell})=0$. Moreover, for every $\ell\in[\EuScript{N}]$, as $T_{2\ell-1}\cap\mathcal{Y}=T_{2\ell}\cap\mathcal{X}=\emptyset$, by \eqref{supports}, we have $\tilde{\mu}_T(T_{2\ell})=\tilde{\nu}_T(T_{2\ell-1})=0$. Hence by \eqref{Eqqq1}--\eqref{Eqqq1.1} and the definitions of $\tilde{\varphi}_T,\tilde{\psi}_T$, we have $\tilde{\varphi}_T(x)=\bar{\varphi}_T(x)$ for $\tilde{\mu}_T$-a.e.\ $x\in T$ and $\tilde{\psi}_T(y)=\bar{\psi}_T(y)$ for $\tilde{\nu}_T$-a.e.\ $y\in T$. Therefore, by \eqref{factorize}, for $\tilde{\mu}_T\otimes \tilde{\nu}_T$-a.e.\ $(x,y)\in T\times T$, we have 
\begin{equation}\label{RKd}
     \frac{d\kappa_T}{d\tilde{\mu}_T\otimes\tilde{\nu}_T}(x,y)=\mathbbm{1}_{\bigcup_{(\ell,\ell')\in\EuScript{L}_T}(T_{\ell}\times T_{\ell'})}(x,y)\|x-y\|^{(d-1)\slash 2}\tilde{\varphi}_T(x)\tilde{\psi}_T(y).
\end{equation}
Now for any $(x,y)\in T\times T$, we define 
\begin{equation}\label{eqtildeh}
    \tilde{h}_T(x,y):=\mathbbm{1}_{\bigcup_{(\ell,\ell')\in\EuScript{L}_T}(T_{\ell}\times T_{\ell'})}(x,y)\|x-y\|^{(d-1)\slash 2}\tilde{\varphi}_T(x)\tilde{\psi}_T(y), 
\end{equation}
\begin{equation}
    h_T(x,y):=\tilde{f}(x)\tilde{g}(y)\tilde{h}_T(x,y)=\mathbbm{1}_{\bigcup_{(\ell,\ell')\in\EuScript{L}_T}(T_{\ell}\times T_{\ell'})}(x,y)\|x-y\|^{(d-1)\slash 2}\tilde{\varphi}_T(x)\tilde{\psi}_T(y)\tilde{f}(x)\tilde{g}(y).
\end{equation}
Note that $h_T(x,y)\in [0,\infty)$ for any $(x,y)\in T\times T$. By \eqref{RKd}, $h_T(\cdot,\cdot)$ is a density of $\kappa_T$ with respect to $\mathcal{H}^1|_T\otimes \mathcal{H}^1|_T$ on $T\times T$. For any $x,y\in T$ such that $u(x)\leq u(y)$, by the definition of $\EuScript{L}_T$, we have $(x,y)\notin \bigcup_{(\ell,\ell')\in\EuScript{L}_T}(T_{\ell}\times T_{\ell'})$ and hence
$h_T(x,y)=\tilde{h}_T(x,y)=0$.

Below, we consider any $\ell\in [\EuScript{N}]$. For any $x,x'\in [a_T-t_{2\ell-2}V_T,a_T-t_{2\ell-1}V_T]$ and $y\in\big(\bigcup_{\ell'\in [\EuScript{N}]:(2\ell-1,2\ell')\in\EuScript{L}_T}T_{2\ell'}\big)\cap\mathcal{Y}$, we have $\min\{\|x-y\|,\|x'-y\|\}\geq 2d_0$, hence
\begin{equation*}
    \max\Big\{\frac{\|x'-y\|}{\|x-y\|},\frac{\|x-y\|}{\|x'-y\|}\Big\}\leq 1+\max\Big\{\frac{\|x'-x\|}{\|x-y\|},\frac{\|x'-x\|}{\|x'-y\|}\Big\}\leq 1+\frac{\|x-x'\|}{2d_0}\leq e^{\|x-x'\|\slash (2d_0)}.
\end{equation*}
Consequently, using the definition of $\tilde{\varphi}_T(\cdot)$, we obtain that
\begin{equation}\label{eqvarphi}
    e^{-C\|x-x'\|}\tilde{\varphi}_T(x)\leq \tilde{\varphi}_T(x')\leq e^{C\|x-x'\|}\tilde{\varphi}_T(x), \quad\text{for any }x,x'\in [a_T-t_{2\ell-2}V_T,a_T-t_{2\ell-1}V_T].
\end{equation}
Similarly, we can deduce that  
\begin{equation}\label{eqpsi}
    e^{-C\|y-y'\|}\tilde{\psi}_T(y)\leq \tilde{\psi}_T(y')\leq e^{C\|y-y'\|}\tilde{\psi}_T(y), \quad\text{for any }y,y'\in [a_T-t_{2\ell-1}V_T, a_T-t_{2\ell}V_T].
\end{equation}
The above two displays, combined with the definitions of $\tilde{\varphi}_T(\cdot)$ and $\tilde{\psi}_T(\cdot)$, imply that $\tilde{\varphi}_T(\cdot)$ and $\tilde{\psi}_T(\cdot)$ are continuous on $T$.

Note that by Definitions \ref{factors} and \ref{De3.8}, $\tilde{f}(a_T)=\tilde{g}(b_T)=0$, hence $h_T(x,y)=0$ for any $(x,y)\in T\times T$ such that $x=a_T$ or $y=b_T$. Below we show that $h_T(\cdot,\cdot)$ is continuous on $\mathrm{int}(T)\times\mathrm{int}(T)$. By Definition \ref{factors}, $\mathfrak{F}(\cdot)$ is continuous on $\mathrm{int}(T)$, hence $\tilde{f}(\cdot)$ and $\tilde{g}(\cdot)$ are continuous on $\mathrm{int}(T)$. Now for any $(x,y)\in \mathrm{int}(T)\times \mathrm{int}(T)$, if $x\in\mathrm{int}(T_{\ell})$ and $y\in\mathrm{int}(T_{\ell'})$ for some $\ell,\ell'\in [2\EuScript{N}]$, then $\mathbbm{1}_{\bigcup_{(\ell,\ell')\in\EuScript{L}_T}(T_{\ell}\times T_{\ell'})}(\cdot,\cdot)$ is continuous at $(x,y)$, hence $h_T(\cdot,\cdot)$ is continuous at $(x,y)$; otherwise, for any $(x',y')\in \mathrm{int}(T)\times \mathrm{int}(T)$ that are sufficiently close to $(x,y)$, by the properties satisfied by $\{T_{\ell}\}_{\ell\in[2\EuScript{N}]}$, either $x'\notin \mathcal{X}$ or $y'\notin\mathcal{Y}$ (which implies $\tilde{f}(x')=0$ or $\tilde{g}(y')=0$), hence $h_T(x',y')=0$, and $h_T(\cdot,\cdot)$ is continuous at $(x,y)$. Therefore, $h_T(\cdot,\cdot)$ is continuous on $\mathrm{int}(T)\times\mathrm{int}(T)$. This completes the proof of part (c).  

Proceeding to the proof of part (d), we consider any $k\in[K],k'\in[K']$ and $x,x',y,y'\in T$ such that $\max\{\dist(x,\mathcal{X}_k), \dist(x',\mathcal{X}_k), \dist(y,\mathcal{Y}_{k'}),\dist(y',\mathcal{Y}_{k'})\}\leq d_0$. Note that by the properties satisfied by $\{T_{\ell}\}_{\ell\in[2\EuScript{N}]}$, there exist $\ell_1,\ell_2\in[\EuScript{N}]$ such that $x,x'\in T_{2\ell_1-1}(\subseteq [a_T-t_{2\ell_1-2}V_T,a_T-t_{2\ell_1-1}V_T])$ and $y,y'\in T_{2\ell_2}(\subseteq [a_T-t_{2\ell_2-1}V_T,a_T-t_{2\ell_2}V_T])$. Moreover, we have $\min\{\|x-y\|, \|x'-y'\|\}\geq d_0$, hence 
\begin{align}\label{bddnew}
    \max\bigg\{\frac{\|x'-y'\|}{\|x-y\|},\frac{\|x-y\|}{\|x'-y'\|}\bigg\}&\leq 1+(\|x'-x\|+\|y'-y\|)\max\bigg\{\frac{1}{\|x-y\|}, \frac{1}{\|x'-y'\|}\bigg\}\nonumber\\
    &\leq 1+\frac{\|x'-x\|+\|y'-y\|}{d_0}\leq e^{(\|x'-x\|+\|y'-y\|)\slash d_0}. 
\end{align}
If $(2\ell_1-1,2\ell_2)\notin\EuScript{L}_T$, by \eqref{eqtildeh}, we have $\tilde{h}_T(x',y')=\tilde{h}_T(x,y)=0$. If $(2\ell_1-1,2\ell_2)\in\EuScript{L}_T$, by \eqref{eqtildeh}, we have $\tilde{h}_T(x,y)=\|x-y\|^{(d-1)\slash 2}\tilde{\varphi}_T(x)\tilde{\psi}_T(y)$ and $\tilde{h}_T(x',y')=\|x'-y'\|^{(d-1)\slash 2}\tilde{\varphi}_T(x')\tilde{\psi}_T(y')$. Hence using \eqref{eqvarphi}--\eqref{bddnew}, we obtain \eqref{inhtild}. This completes the proof of part (d).
\end{proof}

With $h_T$ as in \Cref{Lem2}(c), we define $h:\mathbb{R}^d\times\mathbb{R}^d\rightarrow [0,\infty)$ as follows. 

\begin{definition}[The function $h$]\label{DefhH}
For any $T\in\tilde{\mathcal{S}}^*$ and $x,y\in\mathrm{int}(T)$, we define $h(x,y):=h_T(x,y)$. For any $(x,y)\in(\mathbb{R}^d\times\mathbb{R}^d)\big\backslash\big(\bigcup_{T\in\tilde{\mathcal{S}}^*}(\mathrm{int}(T)\times \mathrm{int}(T))\big)$, we define $h(x,y):=0$. 
\end{definition}

The following lemma establishes the Borel measurability of $h$, which, roughly speaking, follows by combining the measurability of the mappings from $(x,y)$ to the ray $T$ containing $(x,y)$, from~$T$ to the marginals $(\tilde\mu_T,\tilde\nu_T)$, and from those marginals to the solution $h_T$ of the associated EOT problem.

\begin{lemma}\label{Lem4.2n} 
The function $h: \mathbb{R}^d\times\mathbb{R}^d\rightarrow [0,\infty)$ is Borel measurable. 
\end{lemma}

\begin{proof}

For any Polish space $\mathfrak{X}$, we denote by $\mathcal{P}(\mathfrak{X})$ the space of Borel probability measures on $\mathfrak{X}$ equipped with the weak convergence topology. Note that $\mathcal{P}(\mathfrak{X})$ is a Polish space (see, e.g., \cite[Theorem 17.23]{MR1321597}). We also recall the definitions of $a_T,b_T,V_T$ for $T\in\mathcal{S}$ from Definition \ref{TransportRays}. For any $T\in\tilde{\mathcal{S}}^*$, let $\bar{\mu}_T,\bar{\nu}_T\in\mathcal{P}([0,1])$ be such that for any Borel set $A\subseteq [0,1]$,
\begin{equation*}
    \bar{\mu}_T(A)=\tilde{\mu}_T\bigg(\bigg\{x\in T:\frac{\|x-b_T\|}{\|a_T-b_T\|}\in A\bigg\}\bigg), \qquad \bar{\nu}_T(A)=\tilde{\nu}_T\bigg(\bigg\{x\in T:\frac{\|x-b_T\|}{\|a_T-b_T\|}\in A\bigg\}\bigg).
\end{equation*}
Note that by Definition \ref{De3.8},
\begin{equation*}
    \bar{\mu}_T(A)=\int_{\{ta_T+(1-t)b_T:t\in A\}}\tilde{f}(x)d\mathcal{H}^1(x)=\|b_T-a_T\|\int_A \tilde{f}(ta_T+(1-t)b_T)dt,
\end{equation*}
\begin{equation*}
    \bar{\nu}_T(A)=\int_{\{ta_T+(1-t)b_T:t\in A\}}\tilde{g}(x)d\mathcal{H}^1(x)=\|b_T-a_T\|\int_A\tilde{g}(ta_T+(1-t)b_T)dt.
\end{equation*}
Using the above two displays, the Borel measurability of $\tilde{f},\tilde{g}$, and Fubini's theorem, we conclude that
for any Borel set $A\subseteq [0,1]$, the mappings $T\mapsto \bar{\mu}_T(A)$ and $T\mapsto \bar{\nu}_T(A)$ (where $T\in\tilde{\mathcal{S}}^*$) are Borel measurable. Hence by \cite[Theorem 17.24]{MR1321597},
\begin{equation}\label{measuremunu}
    \text{the mappings }T\mapsto \bar{\mu}_T  \text{ and } T\mapsto\bar{\nu}_T \text{ from }\tilde{\mathcal{S}}^*\text{ to }\mathcal{P}([0,1])\text{ are Borel measurable}.
\end{equation}

For any $T\in\tilde{\mathcal{S}}^*$, with $\kappa_T,h_T(\cdot,\cdot)$ defined as in \Cref{Lem2}, we define $\rho_T\in\mathcal{P}([0,1]^2)$ such that for any Borel set $A\subseteq [0,1]^2$,
\begin{equation*}
  \rho_T(A)=\kappa_T\bigg(\bigg\{(x,y)\in T\times T:\bigg(\frac{\|x-b_T\|}{\|a_T-b_T\|},\frac{\|y-b_T\|}{\|a_T-b_T\|}\bigg)\in A\bigg\}\bigg),
\end{equation*}
and define the function $\bar{h}_T:[0,1]^2\rightarrow [0,\infty)$ such that for every $(s,t)\in [0,1]^2$,
\begin{equation*}
    \bar{h}_T(s,t)=\|b_T-a_T\|^2 h_T\big(s a_T+(1-s)b_T, ta_T+(1-t)b_T\big).
\end{equation*}
Note that for any $(x,y)\in\mathrm{int}(T)\times\mathrm{int}(T)$,
\begin{equation}\label{formulashT}
    h_T(x,y)=\frac{1}{\|b_T-a_T\|^2}\bar{h}_T\bigg(\frac{\|x-b_T\|}{\|a_T-b_T\|},\frac{\|y-b_T\|}{\|a_T-b_T\|}\bigg).
\end{equation}
As $h_T(\cdot,\cdot)$ is the density of $\kappa_T$ with respect to $\mathcal{H}^1|_T\otimes\mathcal{H}^1|_T$ on $T\times T$,
\begin{equation}\label{htt}
    \bar{h}_T(\cdot,\cdot) \text{ is the density of }\rho_T  \text{ with respect to }\mathcal{L}^2.
\end{equation}
By \Cref{Lem2}(b), $\rho_T$ is the unique solution to the following minimization problem:
\begin{equation}\label{E4.29}
    \inf_{\substack{\chi\in \Pi(\bar{\mu}_T,\bar{\nu}_T):\\ \chi(\{(s,t)\in [0,1]^2:s< t\})=0}}\bigg\{-\frac{d-1}{2}\int_{[0,1]^2}\log(|s-t|)d\chi(s,t)+H(\chi|\bar{\mu}_T\otimes\bar{\nu}_T)\bigg\},
\end{equation}
with
\begin{equation}\label{finiteness}
    -\frac{d-1}{2}\int_{[0,1]^2}\log(|s-t|)d\rho_T(s,t)+H(\rho_T|\bar{\mu}_T\otimes\bar{\nu}_T)<\infty.
\end{equation}
By \Cref{Lem2}(c), we have
\begin{equation}\label{hconti}
    \bar{h}_T(s,t)=0\text{ if }s=1\text{ or }t=0\text{ or }s\leq t, \text{ and }\bar{h}_T(\cdot,\cdot)\text{ is continuous on }(0,1)^2.
\end{equation}

\paragraph{Step 1.} In this step, we show that the mapping $T\mapsto \rho_T$ from $\tilde{\mathcal{S}}^*$ to $\mathcal{P}([0,1]^2)$ is Borel measurable.

As $\tilde{\mathcal{S}}^*$ is a Borel subset of the standard Borel space $\mathcal{S}$ (recall Lemma \ref{Le2.2} and Definition \ref{defstilde}), by \cite[Corollary 13.4]{MR1321597}, $\tilde{\mathcal{S}}^*$ is a standard Borel space. For any $T\in\tilde{\mathcal{S}}^*$ and $\chi\in\mathcal{P}([0,1]^2)$, we define
\begin{equation}\label{OTC}
    \EuScript{O}(T,\chi):=-\frac{d-1}{2}\int_{[0,1]^2}\log(|s-t|)d\chi(s,t)+H(\chi|\bar{\mu}_T\otimes\bar{\nu}_T)\in (-\infty,\infty].
\end{equation}
We also define 
\begin{equation}\label{GDEf}
\mathcal{G}:=\big\{(T,\chi)\in \tilde{\mathcal{S}}^*\times \mathcal{P}([0,1]^2):\chi\in \Pi(\bar{\mu}_T,\bar{\nu}_T),\chi(\{(s,t)\in [0,1]^2:s< t\})=0\big\}.
\end{equation}  

Consider any fixed $T\in\tilde{\mathcal{S}}^*$. Note that the mapping $\chi\mapsto H(\chi|\bar{\mu}_T\otimes\bar{\nu}_T)$ (where $\chi\in\mathcal{P}([0,1]^2)$) is lower semicontinuous (see, e.g., \cite[Lemma 1.3]{Nutz.20}). Moreover, for any $M\geq 0$, the function $(s,t)\mapsto \min\big\{-\frac{d-1}{2}\log(|s-t|),M\big\}$ is bounded and continuous on $[0,1]^2$. Hence for any sequence $(\chi_n)_{n\geq 1}$ in $\mathcal{P}([0,1]^2)$ such that $\chi_n$ converges weakly to $\chi\in \mathcal{P}([0,1]^2)$ as $n\rightarrow\infty$, we have 
\begin{align*}
    &\liminf_{n\rightarrow\infty}\bigg\{-\frac{d-1}{2}\int_{[0,1]^2}\log(|s-t|)d\chi_n(s,t)\bigg\}\nonumber\\
    \geq\,& \liminf_{n\rightarrow\infty}\bigg\{\int_{[0,1]^2}\min\Big\{-\frac{d-1}{2}\log(|s-t|),M \Big\}d\chi_n(s,t)\bigg\}\nonumber\\
    =\,&\int_{[0,1]^2}\min\bigg\{-\frac{d-1}{2}\log(|s-t|),M\bigg\}d\chi(s,t).
\end{align*}
Taking $M\rightarrow\infty$ in the above display and using the monotone convergence theorem, we obtain that 
\begin{equation*}
    \liminf_{n\rightarrow\infty}\bigg\{-\frac{d-1}{2}\int_{[0,1]^2}\log(|s-t|)d\chi_n(s,t)\bigg\}\geq -\frac{d-1}{2}  
\int_{[0,1]^2}\log(|s-t|)d\chi(s,t).
\end{equation*}
Hence the mapping $\chi\mapsto -\frac{d-1}{2}\int_{[0,1]^2}\log(|s-t|)d\chi(s,t)$ is lower semicontinuous. Thus, recalling \eqref{OTC}, we conclude that for any fixed $T\in\tilde{\mathcal{S}}^*$,
\begin{equation}\label{Eq4.34}
    \text{the mapping }\chi\mapsto \EuScript{O}(T,\chi)\text{ is lower semicontinuous}.
\end{equation}

Using \eqref{measuremunu} and the continuity of the mapping $(\chi_1,\chi_2)\mapsto \chi_1 \otimes \chi_2$ (where $\chi_1,\chi_2\in\mathcal{P}([0,1])$; see, e.g., \cite[Theorem 2.8(ii)]{MR1700749}), we obtain the Borel measurability of the mapping $T\mapsto \bar{\mu}_T\otimes\bar{\nu}_T$ from $\tilde{\mathcal{S}}^*$ to $\mathcal{P}([0,1]^2)$. Consequently, the mapping $(T,\chi)\mapsto (\chi,\bar{\mu}_T\otimes\bar{\nu}_T)$ from $\tilde{\mathcal{S}}^*\times \mathcal{P}([0,1]^2)$ to $\mathcal{P}([0,1]^2)\times \mathcal{P}([0,1]^2)$ is Borel measurable. Using this and the joint lower semicontinuity of the mapping $(\chi,\chi')\mapsto H(\chi|\chi')$ (where $\chi,\chi'\in\mathcal{P}([0,1]^2)$; see, e.g., \cite[Lemma 1.3]{Nutz.20}), we obtain the Borel measurability of the mapping $(T,\chi)\mapsto H(\chi|\bar{\mu}_T\otimes\bar{\nu}_T)$. Moreover, using the lower semicontinuity of the mapping $\chi\mapsto -\frac{d-1}{2}\int_{[0,1]^2}\log(|s-t|)d\chi(s,t)$, we obtain the Borel measurability of the mapping $(T,\chi)\mapsto -\frac{d-1}{2}\int_{[0,1]^2}\log(|s-t|)d\chi(s,t)$. Therefore, 
\begin{equation}\label{Eq4.35}
    \text{the mapping } (T,\chi)\mapsto \EuScript{O}(T,\chi) \text{ is Borel measurable}.
\end{equation}

Now we show that $\mathcal{G}$ as defined in \eqref{GDEf} is a Borel subset of $\tilde{\mathcal{S}}^*\times \mathcal{P}([0,1]^2)$. Note that
\begin{align*}
    \mathcal{G}=\,&\big\{(T,\chi)\in \tilde{\mathcal{S}}^*\times \mathcal{P}([0,1]^2):\chi\in \Pi(\bar{\mu}_T,\bar{\nu}_T)\big\}\nonumber\\
    &\cap \big(\tilde{\mathcal{S}}^*\times \{\chi\in\mathcal{P}([0,1]^2):\chi(\{(s,t)\in [0,1]^2:s<t\})=0\}\big).
\end{align*}
By \cite[Theorem 17.24]{MR1321597}, $\{\chi\in\mathcal{P}([0,1]^2):\chi(\{(s,t)\in [0,1]^2:s<t\})=0\}$ is a Borel subset of $\mathcal{P}([0,1]^2)$. Hence $\tilde{\mathcal{S}}^*\times \{\chi\in\mathcal{P}([0,1]^2):\chi(\{(s,t)\in [0,1]^2:s<t\})=0\}$ is a Borel subset of $\tilde{\mathcal{S}}^*\times \mathcal{P}([0,1]^2)$. For any $\chi\in\mathcal{P}([0,1]^2)$, we define $\mathtt{P}_1(\chi)$ and $\mathtt{P}_2(\chi)$ to be the two marginals of $\chi$; note that the mappings $\mathtt{P}_1,\mathtt{P}_2:\mathcal{P}([0,1]^2)\rightarrow \mathcal{P}([0,1])$ are continuous. Hence by \eqref{measuremunu}, the mapping $\mathfrak{E}:\tilde{\mathcal{S}}^*\times\mathcal{P}([0,1]^2)\rightarrow \mathcal{P}([0,1])^4$ with $\mathfrak{E}(T,\chi)=(\bar{\mu}_T,\bar{\nu}_T,\mathtt{P}_1(\chi),\mathtt{P}_2(\chi))$ for every $(T,\chi)\in\tilde{\mathcal{S}}^*\times\mathcal{P}([0,1]^2)$ is Borel measurable. As $\Delta_0:=\{(\chi_1,\chi_2,\chi_1,\chi_2):\chi_1,\chi_2\in\mathcal{P}([0,1])\}$ is a closed subset of $\mathcal{P}([0,1])^4$, we conclude that $\big\{(T,\chi)\in \tilde{\mathcal{S}}^*\times \mathcal{P}([0,1]^2):\chi\in \Pi(\bar{\mu}_T,\bar{\nu}_T)\big\}=\mathfrak{E}^{-1}(\Delta_0)$ is a Borel subset of $\tilde{\mathcal{S}}^*\times \mathcal{P}([0,1]^2)$. Therefore, $\mathcal{G}$ is a Borel subset of $\tilde{\mathcal{S}}^*\times \mathcal{P}([0,1]^2)$. 

Define the projection map $\mathtt{proj}_{\tilde{\mathcal{S}}^*}:\tilde{\mathcal{S}}^*\times\mathcal{P}([0,1]^2) \rightarrow \tilde{\mathcal{S}}^*$ with $\mathtt{proj}_{\tilde{\mathcal{S}}^*}(T,\chi)=T$ for every $(T,\chi)\in \tilde{\mathcal{S}}^*\times\mathcal{P}([0,1]^2)$. For any $T\in\tilde{\mathcal{S}}^*$, we define
\begin{equation*}
    \EuScript{F}(T):=\inf_{\substack{\chi\in \Pi(\bar{\mu}_T,\bar{\nu}_T):\\ \chi(\{(s,t)\in [0,1]^2:s< t\})=0}}\{\EuScript{O}(T,\chi)\}.
\end{equation*}
Note that by \eqref{finiteness}, $\EuScript{F}(T)\in (-\infty,\infty)$. Below we consider any $a_0\in\mathbb{R}$. As the minimization problem \eqref{E4.29} admits the solution $\rho_T$ that satisfies \eqref{finiteness}, for any $T\in\tilde{\mathcal{S}}^*$, $\EuScript{F}(T)\leq a_0$ if and only if there exists $\chi\in \Pi(\bar{\mu}_T,\bar{\nu}_T)$ such that $\chi(\{(s,t)\in [0,1]^2:s< t\})=0$ and $\EuScript{O}(T,\chi)\leq a_0$. Hence with $\mathcal{G}$ as in \eqref{GDEf}, we have  
\begin{equation}\label{Eq4.36}
    \{T\in\tilde{\mathcal{S}}^*:\EuScript{F}(T)\leq a_0\}=\mathrm{proj}_{\tilde{\mathcal{S}}^*}(\{(T,\chi)\in\mathcal{G}:\EuScript{O}(T,\chi)\leq a_0\}). 
\end{equation}
Note that by \eqref{Eq4.35} and the Borel measurability of $\mathcal{G}$,
\begin{equation}\label{Eq4.37}
    \{(T,\chi)\in\mathcal{G}:\EuScript{O}(T,\chi)\leq a_0\}\text{ is a Borel subset of }\tilde{\mathcal{S}}^*\times \mathcal{P}([0,1]^2).
\end{equation}
Let $A_0:=\{(T,\chi)\in\mathcal{G}:\EuScript{O}(T,\chi)\leq a_0\}$. For any $T\in\tilde{\mathcal{S}}^*$, we have  
\begin{align}\label{chide}
    \{\chi\in\mathcal{P}([0,1]^2):(T,\chi)\in A_0\}=\,&\big\{\chi\in\mathcal{P}([0,1]^2):\chi(\{(s,t)\in [0,1]^2:s<t\})=0\big\}\nonumber\\
    &\cap\{\chi\in\Pi(\bar{\mu}_T,\bar{\nu}_T):\EuScript{O}(T,\chi)\leq a_0\}.
\end{align}
By \eqref{Eq4.34}, the set $\{\chi\in\mathcal{P}([0,1]^2):\EuScript{O}(T,\chi)\leq a_0\}$ is closed. As $\Pi(\bar{\mu}_T,\bar{\nu}_T)$ is compact (hence closed), $\{\chi\in\Pi(\bar{\mu}_T,\bar{\nu}_T):\EuScript{O}(T,\chi)\leq a_0\}=\Pi(\bar{\mu}_T,\bar{\nu}_T)\cap \{\chi\in\mathcal{P}([0,1]^2):\EuScript{O}(T,\chi)\leq a_0\}$ is closed. Moreover, as $\{(s,t)\in [0,1]^2:s<t\}$ is an open subset of $[0,1]^2$, for any $(\chi_n)_{n\geq 1},\chi$ in $\mathcal{P}([0,1]^2)$ such that $\chi_n(\{(s,t)\in [0,1]^2:s<t\})=0$ for all $n\in\mathbb{N}^*$ and $\chi_n$ converges weakly to $\chi$ as $n\rightarrow\infty$, we have 
\begin{equation*}
    \chi(\{(s,t)\in [0,1]^2:s<t\})\leq\liminf\limits_{n\rightarrow\infty}\chi_n(\{(s,t)\in [0,1]^2:s<t\})=0,
\end{equation*}
hence $\chi(\{(s,t)\in [0,1]^2:s<t\})=0$. Consequently, $\big\{\chi\in\mathcal{P}([0,1]^2):\chi(\{(s,t)\in [0,1]^2:s<t\})=0\big\}$ is closed. Therefore, by \eqref{chide}, $\{\chi\in\mathcal{P}([0,1]^2):(T,\chi)\in A_0\}$ is closed. As $\{\chi\in\mathcal{P}([0,1]^2):(T,\chi)\in A_0\}$ is a subset of the compact set $\Pi(\bar{\mu}_T,\bar{\nu}_T)$, we conclude that 
\begin{equation}\label{Eq4.38}
    \text{for any }T\in\tilde{\mathcal{S}}^*, \text{ the set } \{\chi\in\mathcal{P}([0,1]^2):(T,\chi)\in A_0\} \text{ is compact}.
\end{equation}
As $\tilde{\mathcal{S}}^*$ is a standard Borel space and $\mathcal{P}([0,1]^2)$ is a Polish space, by the Arsenin--Kunugui theorem (\cite[Theorem 18.18]{MR1321597}), noting \eqref{Eq4.36}, \eqref{Eq4.37}, and \eqref{Eq4.38}, we conclude that $\{T\in\tilde{\mathcal{S}}^*:\EuScript{F}(T)\leq a_0\}$ is a Borel subset of $\tilde{\mathcal{S}}^*$. As $a_0\in\mathbb{R}$ is arbitrary, the mapping $\EuScript{F}:\tilde{\mathcal{S}}^*\rightarrow \mathbb{R}$ is Borel measurable.     

As $\rho_T$ (where $T\in\tilde{\mathcal{S}}^*$) is the unique solution to the minimization problem \eqref{E4.29} and satisfies \eqref{finiteness}, we have  
\begin{equation*}
    \{(T,\chi)\in\tilde{\mathcal{S}}^*\times\mathcal{P}([0,1]^2):\chi=\rho_T\}=\{(T,\chi)\in\tilde{\mathcal{S}}^*\times\mathcal{P}([0,1]^2):\EuScript{O}(T,\chi)=\EuScript{F}(T)\}\cap\mathcal{G},
\end{equation*}
which is a Borel subset of $\tilde{\mathcal{S}}^*\times\mathcal{P}([0,1]^2)$ by \eqref{Eq4.35} and the Borel measurability of $\EuScript{F}$ and $\mathcal{G}$. Hence by \cite[Theorem 14.12]{MR1321597}, 
\begin{equation}\label{borelrhot}
    \text{the mapping } T \mapsto \rho_T \text{ from } \tilde{\mathcal{S}}^* \text{ to } \mathcal{P}([0,1]^2) \text{ is Borel measurable}.
\end{equation} 

\paragraph{Step 2.} In this step, we complete the proof of the Borel measurability of $h$.

By the continuity of $\bar{h}_T(\cdot,\cdot)$ on $(0,1)^2$ (see \eqref{hconti}) and \eqref{htt}, for any $T\in\tilde{\mathcal{S}}^*$ and $(s,t)\in (0,1)^2$, 
\begin{align}\label{Eqnnn1}
    \bar{h}_T(s,t)=\,&\lim_{n\rightarrow\infty}\frac{\int_{\{(s',t')\in(0,1)^2:|s'-s|\leq n^{-1},|t'-t|\leq n^{-1}\}}\bar{h}_T(s',t')ds'dt'}{4n^{-2}}\nonumber\\
=\,&\limsup_{n\rightarrow\infty} \bigg\{\frac{n^2}{4}\cdot \rho_T(\{(s',t')\in(0,1)^2:|s'-s|\leq n^{-1},|t'-t|\leq n^{-1}\})\bigg\}  .
\end{align}
By \eqref{borelrhot} and \cite[Theorem 17.24]{MR1321597}, for any $n\in\mathbb{N}^*$ and $(s,t)\in (0,1)^2$, the mapping $T\mapsto \rho_T(\{(s',t')\in(0,1)^2:|s'-s|\leq n^{-1},|t'-t|\leq n^{-1}\})$ (where $T\in\tilde{\mathcal{S}}^*$) is Borel measurable. Hence by \eqref{Eqnnn1}, for any $(s,t)\in (0,1)^2$, the mapping $T\mapsto \bar{h}_T(s,t)$ (where $T\in\tilde{\mathcal{S}}^*$) is Borel measurable. Combined with the continuity of $\bar{h}_T(\cdot,\cdot)$ on $(0,1)^2$, this implies that the mapping $\mathtt{H}:\tilde{\mathcal{S}}^*\times (0,1)^2\rightarrow [0,\infty)$ with $\mathtt{H}(T,s,t)=\bar{h}_T(s,t)$ for every $(T,s,t)\in \tilde{\mathcal{S}}^*\times(0,1)^2$ is a Carath\'eodory function and hence is Borel measurable (see \cite[Theorem 4.51]{MR2378491}).  

Let $\EuScript{Z}:=\bigcup_{T\in\tilde{\mathcal{S}}^*}(\mathrm{int}(T)\times \mathrm{int}(T))$. Below we show that $\EuScript{Z}$ is a Borel set. Define the mapping $\EuScript{W}:\mathcal{T}_1^*\times\mathcal{T}_1^*\rightarrow \mathcal{S}\times\mathcal{S}$ such that $\EuScript{W}(x,y)=(T(x),T(y))$ for every $(x,y)\in\mathcal{T}_1^*\times\mathcal{T}_1^*$. By Lemma \ref{L2.0nn}, $\EuScript{W}$ is Borel measurable. Moreover, $\Delta_{\mathcal{S}}:=\{(T,T')\in\mathcal{S}\times\mathcal{S}: T=T'\}$ is a Borel subset of $\mathcal{S}\times\mathcal{S}$. Now note that 
\begin{equation*}
\EuScript{Z}=\big\{(x,y)\in\mathcal{T}_1^*\times \mathcal{T}_1^*:T(x)\in \tilde{\mathcal{S}}^*, T(x)=T(y)\big\}=\EuScript{W}^{-1}\big((\tilde{\mathcal{S}}^*\times\mathcal{S})\cap\Delta_{\mathcal{S}}\big).
\end{equation*}
As $\tilde{\mathcal{S}}^*$ is a Borel subset of $\mathcal{S}$, $(\tilde{\mathcal{S}}^*\times\mathcal{S})\cap\Delta_{\mathcal{S}}$ is a Borel subset of $\mathcal{S}\times\mathcal{S}$, and so $\EuScript{Z}$ is a Borel set.  

Now note that by \eqref{formulashT}, for any $(x,y)\in\EuScript{Z}$, we have 
\begin{align*}
    h(x,y)=h_{T(x)}(x,y)=\,&\frac{1}{\|b(x)-a(x)\|^2}\bar{h}_{T(x)}\bigg(\frac{\|x-b(x)\|}{\|a(x)-b(x)\|},\frac{\|y-b(x)\|}{\|a(x)-b(x)\|}\bigg)\nonumber\\
    =\,& \frac{1}{\|b(x)-a(x)\|^2}\mathtt{H}\bigg(T(x),\frac{\|x-b(x)\|}{\|a(x)-b(x)\|},\frac{\|y-b(x)\|}{\|a(x)-b(x)\|}\bigg),
\end{align*}
where $a(\cdot),b(\cdot)$ are as in Definition \ref{Def2.6}. For any $(x,y)\in (\mathbb{R}^d\times\mathbb{R}^d)\backslash\EuScript{Z}$, we have $h(x,y)=0$. Recall that $a(\cdot),b(\cdot),T(\cdot)$ are Borel measurable (see Definition \ref{Def2.6} and Lemma \ref{L2.0nn}). Consequently, using the Borel measurability of $\EuScript{Z}$ and $\mathtt{H}$, we conclude that $h$ is Borel measurable.
\end{proof}

\subsection{Construction of the measure $\tilde{\gamma}_{\eps;k,k'}$}\label{Sect.4.2}

In this subsection, we construct Borel measures $\tilde{\gamma}_{\eps;k,k'}$ on $\mathbb{R}^d\times\mathbb{R}^d$ for each $k\in[K],k'\in[K']$. These measures will later be used (see \eqref{def_gamma_eps}) to construct the coupling $\tilde{\gamma}_{\eps}\in\Pi(\mu,\nu)$ announced in~\eqref{ubb.estimates} of the proof overview in Section~\ref{se:upper}.  

Section~\ref{sec:preliminary_definition} introduces preliminary notation and definitions. In Sections~\ref{Sect.4.2.1} and~\ref{Sect.4.2.2}, we introduce two functions $p(\cdot,\cdot)$ and $r(\cdot,\cdot)$ on $\mathbb{R}^d\times\mathbb{R}^d$, corresponding to Steps~1--2 and Steps 4--5 of the proof overview, respectively. Section~\ref{Sect.4.2.3} then defines $\tilde{\gamma}_{\eps;k,k'}$ in terms of $p(\cdot,\cdot)$ and $r(\cdot,\cdot)$.

\subsubsection{Preliminary notation and definitions}\label{sec:preliminary_definition}

Recall the domain decomposition scheme as described in Section \ref{sectdom}. For any $k\in[K],k'\in[K']$ and $x\in\mathbb{R}^d$, with $h$ as in \Cref{DefhH}, we define
\begin{equation}\label{fkkdef}
    f_{k,k'}(x):=\mathbbm{1}_{\mathcal{X}_k^{\circ}\cap\mathfrak{T}^*\cap\tilde{\mathcal{T}}_{1;k,k'}^*}(x)
      \mathfrak{F}(x)^{-1}  
    \int_{\mathrm{int}(T(x))\cap \mathcal{Y}_{k'}^{\circ}}h(x,z)d\mathcal{H}^1(z),
\end{equation} 
\begin{equation}\label{gkkdef}
    g_{k,k'}(x):=\mathbbm{1}_{\mathcal{Y}_{k'}^{\circ}\cap\mathfrak{T}^*\cap\tilde{\mathcal{T}}_{1;k,k'}^*}(x)
     \mathfrak{F}(x)^{-1}  
    \int_{\mathrm{int}(T(x))\cap \mathcal{X}_{k}^{\circ}}h(z,x)d\mathcal{H}^1(z),
\end{equation}
where $\tilde{\mathcal{T}}_{1;k,k'}^*$, $\mathfrak{T}^*$, and $\mathfrak{F}(\cdot)$ are as in \eqref{deftildeTkk} and Definition~\ref{factors}, and we recall from Remark~\ref{positivityF} that $\mathfrak{F}(x)>0$ for every $x\in\mathfrak{T}^*$. Define $\mu_{k,k'},\nu_{k,k'}$ to be the Borel measures on $\mathbb{R}^d$ with densities $f_{k,k'},g_{k,k'}$ (respectively) with respect to $\mathcal{L}^d$, and
\begin{equation}\label{def:OT}
    \OT_{k,k'}(\mu,\nu):=\int_{\mathbb{R}^d}u(x)d\mu_{k,k'}(x)-\int_{\mathbb{R}^d}u(x)d\nu_{k,k'}(x).
\end{equation}

Note that by \Cref{Lem2} and Definition~\ref{DefhH}, for any $T\in\tilde{\mathcal{S}}^*$ and any Borel set $A\subseteq T$, $\int_{A}f_{k,k'}(x)\mathfrak{F}(x)d\mathcal{H}^1(x)\leq \int_{A\times T}h(x,z)d\mathcal{H}^1(x)d\mathcal{H}^1(z)=\int_{A}f(x)\mathfrak{F}(x)d\mathcal{H}^1(x)$. Hence for any $T\in\tilde{\mathcal{S}}^*$, we have $f_{k,k'}(x)\mathfrak{F}(x)\leq f(x)\mathfrak{F}(x)$ for $\mathcal{H}^1$-a.e.\ $x\in T$. Similarly, $g_{k,k'}(x)\mathfrak{F}(x)\leq g(x)\mathfrak{F}(x)$ for $\mathcal{H}^1$-a.e.\ $x\in T$. Moreover, for any $T\in\mathcal{S}\backslash\tilde{\mathcal{S}}^*$ and $x\in \mathrm{int}(T)$, $f_{k,k'}(x)=g_{k,k'}(x)=0$. Hence for any $T\in\mathcal{S}$,
\begin{equation}\label{fkkgkkb}
    f_{k,k'}(x)\mathfrak{F}(x)\leq f(x)\mathfrak{F}(x), \quad g_{k,k'}(x)\mathfrak{F}(x)\leq g(x)\mathfrak{F}(x), \quad\text{ for }\mathcal{H}^1\text{-a.e.\ }x\in T.
\end{equation}

The following lemma provides a decomposition of the marginals $\mu,\nu$ and of the optimal cost $\OT(\mu,\nu)$ in terms of $\mu_{k,k'}$, $\nu_{k,k'}$, and $\OT_{k,k'}(\mu,\nu)$.

\begin{lemma}\label{Ln4.1}
We have $\mu=\sum_{k=1}^K\sum_{k'=1}^{K'}\mu_{k,k'}$ and $\nu=\sum_{k=1}^K\sum_{k'=1}^{K'}\nu_{k,k'}$. Moreover,
\begin{equation*}
    \OT(\mu,\nu)=\sum_{k=1}^K\sum_{k'=1}^{K'}\OT_{k,k'}(\mu,\nu).
\end{equation*} 
\end{lemma}
\begin{proof}

Note that by Lemma \ref{Lem3.18nn},
\begin{align}
0=\mathcal{L}^d\Big(\Big(\mathcal{X}\backslash\bigcup_{k=1}^K\mathcal{X}_k^{\circ}\Big)\cap\mathfrak{T}^*\Big)=\int_{\mathcal{S}}d\lambda(T)\int_{T}\mathbbm{1}_{\big(\mathcal{X}\backslash\bigcup_{k=1}^K\mathcal{X}_k^{\circ}\big)\cap\mathfrak{T}^*}(x)\mathfrak{F}(x)d\mathcal{H}^1(x).
\end{align}
Hence for $\lambda$-a.e.\ $T\in\mathcal{S}$ and $\mathcal{H}^1$-a.e.\ $x\in T$, we have $\mathbbm{1}_{\big(\mathcal{X}\backslash\bigcup_{k=1}^K\mathcal{X}_k^{\circ}\big)\cap\mathfrak{T}^*}(x) \mathfrak{F}(x)=0$. Moreover, by Remark~\ref{positivityF}, if $x\in \big(\mathcal{X}\backslash\bigcup_{k=1}^K\mathcal{X}_k^{\circ}\big)\cap\mathfrak{T}^*$, then $\mathfrak{F}(x)>0$. Hence
\begin{equation}\label{generx}
    x\notin \Big(\mathcal{X}\backslash\bigcup_{k=1}^K\mathcal{X}_k^{\circ}\Big)\cap\mathfrak{T}^*,\quad\text{ for }\lambda\text{-a.e.\ }T\in\mathcal{S},\text{ }\mathcal{H}^1\text{-a.e.\ }x\in T.
\end{equation}
Similarly,
\begin{equation}\label{genery}
    x\notin \Big(\mathcal{Y}\backslash\bigcup_{k'=1}^{K'}\mathcal{Y}_{k'}^{\circ}\Big)\cap\mathfrak{T}^*,\quad\text{ for }\lambda\text{-a.e.\ }T\in\mathcal{S},\text{ }\mathcal{H}^1\text{-a.e.\ }x\in T.
\end{equation}

For any Borel set $A\subseteq \mathbb{R}^d$, by \eqref{fkkdef} and Lemma \ref{Lem3.18nn}, we have 
\begin{align}
&\sum_{k=1}^K\sum_{k'=1}^{K'}\mu_{k,k'}(A)=\sum_{k=1}^K\sum_{k'=1}^{K'}\int_A f_{k,k'}(x)dx=\sum_{k=1}^K\sum_{k'=1}^{K'}\int_{A\cap\mathfrak{T}^*} f_{k,k'}(x)dx\nonumber\\
=\, & \sum_{k=1}^K\sum_{k'=1}^{K'}\int_{\mathcal{S}}d\lambda(T)\int_T \mathbbm{1}_A(x) f_{k,k'}(x)\mathfrak{F}(x)d\mathcal{H}^1(x)\nonumber\\
=\, &\sum_{k=1}^K\sum_{k'=1}^{K'}\int_{\mathcal{S}}d\lambda(T)\int_{(T\cap\mathcal{X}_k^{\circ})\times (T\cap\mathcal{Y}_{k'}^{\circ})}\mathbbm{1}_{A\cap\mathfrak{T}^*\cap \tilde{\mathcal{T}}_{1;k,k'}^*}(x)h(x,y)d\mathcal{H}^1(x) d\mathcal{H}^1(y)\nonumber\\
=\, &\sum_{k=1}^K\sum_{k'=1}^{K'}\int_{\tilde{\mathcal{S}}^*}d\lambda(T)\int_{(T\cap\mathcal{X}_k^{\circ})\times (T\cap\mathcal{Y}_{k'}^{\circ})}\mathbbm{1}_{A\cap\mathfrak{T}^*}(x)h(x,y)d\mathcal{H}^1(x) d\mathcal{H}^1(y)\nonumber\\
=\, & \int_{\tilde{\mathcal{S}}^*}d\lambda(T)\int_{(T\cap\mathcal{X})\times (T\cap\mathcal{Y})}\mathbbm{1}_{A\cap\mathfrak{T}^*}(x)h_T(x,y)d\mathcal{H}^1(x) d\mathcal{H}^1(y)\nonumber\\
=\, & \int_{\tilde{\mathcal{S}}^*}d\lambda(T)\int_{T\cap\mathcal{X}}\mathbbm{1}_{A\cap\mathfrak{T}^*}(x)\tilde{f}(x)d\mathcal{H}^1(x)=\int_{\tilde{\mathcal{S}}^*}d\lambda(T)\int_{T\cap\mathcal{X}}\mathbbm{1}_{A\cap\mathfrak{T}^*}(x)f(x)\mathfrak{F}(x)d\mathcal{H}^1(x)\nonumber\\
=\, & \int_{\mathcal{S}}d\lambda(T)\int_{T}\mathbbm{1}_{A\cap\mathfrak{T}^*}(x)f(x)\mathfrak{F}(x)d\mathcal{H}^1(x)=\int_{A\cap\mathfrak{T}^*}f(x)dx=\int_{A}f(x)dx,
\end{align}
where in the fourth line we note that $\lambda(\mathcal{S}\backslash\tilde{\mathcal{S}}^*)=0$ (by \Cref{defstilde}) and $\mathrm{int}(T)\subseteq \tilde{\mathcal{T}}_{1;k,k'}^*$ for any $T\in\mathcal{S}$ such that $T\cap\mathcal{X}_k^{\circ},T\cap\mathcal{Y}_{k'}^{\circ}\neq\emptyset$. Moreover, we use  \eqref{generx}--\eqref{genery} and \Cref{DefhH} in the fifth line, $\kappa_T\in\Pi(\tilde{\mu}_T,\tilde{\nu}_T)$ in the sixth line, and Lemmas~\ref{Lem3.18nn} and~\ref{Lem3.15} in the last two equalities. Hence $\sum_{k=1}^K\sum_{k'=1}^{K'}\mu_{k,k'}=\mu$. Similarly, we can deduce that $\sum_{k=1}^K\sum_{k'=1}^{K'}\nu_{k,k'}=\nu$. By \eqref{def:OT} and \eqref{OTrela},
\begin{align*}
    \sum_{k=1}^K\sum_{k'=1}^{K'}\OT_{k,k'}(\mu,\nu)=\,&\sum_{k=1}^K\sum_{k'=1}^{K'}\int_{\mathbb{R}^d}u(x)d\mu_{k,k'}(x)-\sum_{k=1}^K\sum_{k'=1}^{K'}\int_{\mathbb{R}^d}u(x)d\nu_{k,k'}(x)\nonumber\\
    =\,&\int_{\mathbb{R}^d}u(x)f(x) dx-\int_{\mathbb{R}^d}u(x)g(x) dx=\OT(\mu,\nu). \qedhere
\end{align*}
\end{proof}

Throughout the rest of this subsection, we fix $k\in[K],k'\in[K']$. We also fix $M_1\geq 10$ and $\eps,\delta_1,\delta_2,\delta_3\in (0,1\slash 100)$ such that $\delta_2\sqrt{\eps}\leq \mathtt{h}\slash 10$, where $\mathtt{h}$ is as defined in~\eqref{eq:mathtthDefn} (note that $\mathtt{h}$ is of constant order).

\subsubsection{Definition of $p(\cdot,\cdot)$}\label{Sect.4.2.1}

For any $x,y\in\mathbb{R}^d$, we define
\begin{eqnarray}\label{defz}
    \mathtt{z}(x,y):=x+\langle y-x, V(x)\rangle V(x),\qquad
   \mathtt{w}(x,y):=y-x-\langle y-x, V(x)\rangle V(x),
\end{eqnarray}
\begin{eqnarray}\label{defzp}
    \mathtt{z}'(x,y):=y+\langle x-y, V(y)\rangle V(y),\qquad \mathtt{w}'(x,y):=x-y-\langle x-y,V(y)\rangle V(y).    
\end{eqnarray}
Note that for any $x,y\in\mathbb{R}^d$, $\mathtt{z}(x,y)+\mathtt{w}(x,y)=y$ and $\mathtt{z}'(x,y)+\mathtt{w}'(x,y)=x$. Moreover, if $x\in\mathcal{T}_1^{*}$, then $\mathtt{w}(x,y)\in O(V(x))$; if $y\in\mathcal{T}_1^{*}$, then $\mathtt{w}'(x,y)\in O(V(y))$. For any $x,y\in\mathbb{R}^d$, with $\tilde{\mathcal{T}}_{1;k,k'}^*,\mathfrak{T}^*\subseteq\mathcal{T}_1^{*}$ and $\mathfrak{F}(\cdot)$ as in \eqref{deftildeTkk} and Definition \ref{factors}, we define (recall from Remark \ref{positivityF} that $\mathfrak{F}(x)>0$ for any $x\in\mathfrak{T}^*$)
\begin{align}\label{p1}
   p_1(x,y):=\,& \mathbbm{1}_{x\in\mathcal{X}_k^{\circ}\cap\mathfrak{T}^*\cap\tilde{\mathcal{T}}_{1;k,k'}^{*}}\mathbbm{1}_{\min\{\alpha(x),\beta(x)\}\geq\delta_1}\mathbbm{1}_{\mathtt{z}(x,y)\in\mathrm{int}(T(x))\cap\mathcal{Y}_{k'}^{\circ}}\mathbbm{1}_{\langle x-\mathtt{z}(x,y),V(x)\rangle>0}\nonumber\\
   &\times\mathbbm{1}_{\min\{\alpha(\mathtt{z}(x,y)),\beta(\mathtt{z}(x,y))\}\geq\delta_1}
   \mathfrak{F}(x)^{-1}h(x,\mathtt{z}(x,y))\nonumber\\
    & 
\times(2\pi\|x-\mathtt{z}(x,y)\|\eps)^{-(d-1)\slash 2}\sqrt{\det\big(\mathbf{I}_d+\|x-\mathtt{z}(x,y)\|F(\mathtt{z}(x,y))\big)}\nonumber\\
&\times e^{-\frac{1}{2\eps}\mathtt{w}(x,y)^{\top}\big(\frac{1}{\|x-\mathtt{z}(x,y)\|}\mathbf{I}_d+F(\mathtt{z}(x,y))\big)\mathtt{w}(x,y)}\cdot  \mathbbm{1}_{\|\mathtt{w}(x,y)\|\leq M_1\sqrt{\eps}},
\end{align}
\begin{align}\label{p2}
  p_2(x,y):=\,& \mathbbm{1}_{y\in\mathcal{Y}_{k'}^{\circ}\cap\mathfrak{T}^*\cap\tilde{\mathcal{T}}_{1;k,k'}^{*}}\mathbbm{1}_{\min\{\alpha(y),\beta(y)\}\geq\delta_1}\mathbbm{1}_{\mathtt{z}'(x,y)\in \mathrm{int}(T(y))\cap\mathcal{X}_k^{\circ}}\mathbbm{1}_{\langle y-\mathtt{z}'(x,y),V(y)\rangle<0}\nonumber\\
  &\times\mathbbm{1}_{\min\{\alpha(\mathtt{z}'(x,y)),\beta(\mathtt{z}'(x,y))\}\geq \delta_1}\mathfrak{F}(y)^{-1}h(\mathtt{z}'(x,y),y)\nonumber\\
  &\times(2\pi\|y-\mathtt{z}'(x,y)\|\eps)^{-(d-1)\slash 2}\sqrt{\det\big(\mathbf{I}_d-\|y-\mathtt{z}'(x,y)\|F(\mathtt{z}'(x,y))\big)}\nonumber\\
    &\times e^{-\frac{1}{2\eps}\mathtt{w}'(x,y)^{\top}\big(\frac{1}{\|y-\mathtt{z}'(x,y)\|}\mathbf{I}_d-F(\mathtt{z}'(x,y))\big)\mathtt{w}'(x,y)}.
\end{align}
Note that by \Cref{L2.7,Defn4.1.1n}, $p_1(x,y),p_2(x,y)\geq 0$ for every $(x,y)\in\mathbb{R}^d\times\mathbb{R}^d$.  

Recall the definition of $\mathfrak{q}_{k,k'}(\cdot)$ from Definition~\ref{defkk}. For any $x\in\tilde{\mathcal{T}}^{*}_{1;k,k'}$, we define 
\begin{equation}\label{def_Gkk}
    \mathfrak{G}_{k,k'}(x):=\det\big(\mathbf{I}_d+\langle x-\mathfrak{q}_{k,k'}(x), V(x)\rangle F(\mathfrak{q}_{k,k'}(x))\big).
\end{equation}
As $\mathfrak{q}_{k,k'}:\tilde{\mathcal{T}}^{*}_{1;k,k'}\rightarrow H_{k,k'}$ is Borel measurable,  the mapping $\mathfrak{G}_{k,k'}:\tilde{\mathcal{T}}^{*}_{1;k,k'}\rightarrow\mathbb{R}$ is also Borel measurable. Note that by \Cref{L2.7,Defn4.1.1n},
\begin{equation}\label{Gkk}
    \mathfrak{G}_{k,k'}(x)>0,\quad\text{ for all }x\in\mathfrak{T}\cap \tilde{\mathcal{T}}^{*}_{1;k,k'}. 
\end{equation}

For any $x,y\in\mathcal{T}_1^{*}$, if $V(x)+V(y)\neq 0$, we define
\begin{equation}\label{x'd}
  \mathtt{x}'(x,y) := y-\frac{\langle y-x,V(x)+V(y)\rangle}{\langle V(y), V(x)+V(y)\rangle}V(y),
\end{equation}
\begin{equation}\label{y'd}
  \mathtt{y}'(x,y) := x+\frac{\langle y-x,V(x)+V(y)\rangle}{\langle V(x), V(x)+V(y)\rangle}V(x);
\end{equation}
if $V(x)+V(y)=0$, we define
\begin{equation}\label{alt_def}
   \mathtt{x}'(x,y):= y,\qquad 
 \mathtt{y}'(x,y):= x.
\end{equation}
Note that for any $x,y\in\mathcal{T}_1^*$ such that $V(x)+V(y)\neq 0$, $\mathtt{x}'(x,y)\in\mathrm{int}(T(y))$, $\mathtt{y}'(x,y)\in\mathrm{int}(T(x))$, we have $V(\mathtt{x}'(x,y))=V(y)$ and $V(\mathtt{y}'(x,y))=V(x)$, hence
{\small
\begin{align}\label{zwzw}
    \mathtt{z}(\mathtt{x}'(x,y),\mathtt{y}'(x,y))&=y-\langle y-x,V(y)\rangle V(y)+\frac{\langle y-x, V(x)+V(y)\rangle \langle V(x), V(y)\rangle}{\langle V(x), V(x)+V(y)\rangle} V(y),\nonumber\\
    \mathtt{w}(\mathtt{x}'(x,y),\mathtt{y}'(x,y))&=x-y-\langle x-y,V(y)\rangle V(y)+\frac{\langle y-x, V(x)+V(y)\rangle}{\langle V(x), V(x)+V(y)\rangle} \big(V(x)-\langle V(x),V(y)\rangle V(y)\big),\nonumber\\ \mathtt{z}'(\mathtt{x}'(x,y),\mathtt{y}'(x,y))&=x-\langle x-y, V(x)\rangle V(x)+\frac{\langle x-y, V(x)+V(y)\rangle \langle V(x), V(y)\rangle}{\langle V(y), V(x)+V(y)\rangle} V(x),\nonumber\\ \mathtt{w}'(\mathtt{x}'(x,y),\mathtt{y}'(x,y))&=y-x-\langle y-x,V(x)\rangle V(x)+\frac{\langle x-y, V(x)+V(y)\rangle}{\langle V(y), V(x)+V(y)\rangle} \big(V(y)-\langle V(x),V(y)\rangle V(x)\big).
\end{align}
}

For any $(x,y)\in(\mathfrak{T}\cap\tilde{\mathcal{T}}_{1;k,k'}^{*})\times (\mathfrak{T}\cap\tilde{\mathcal{T}}_{1;k,k'}^{*})$, by \eqref{Gkk}, we have $\mathfrak{G}_{k,k'}(x),\mathfrak{G}_{k,k'}(y)>0$; moreover, if we further assume that $\mathtt{x}'(x,y)\in\mathrm{int}(T(y))$ and $\mathtt{y}'(x,y)\in\mathrm{int}(T(x))$, then $\mathtt{x}'(x,y),\mathtt{y}'(x,y)\in\mathfrak{T}\cap\tilde{\mathcal{T}}_{1;k,k'}^{*}$ (recall \eqref{deftildeTkk} and \cref{Defn4.1.1n}), which via \eqref{Gkk} implies that $\mathfrak{G}_{k,k'}(\mathtt{x}'(x,y)),\mathfrak{G}_{k,k'}(\mathtt{y}'(x,y))>0$. Now for any $(x,y)\in(\mathcal{X}_k^{\circ}\cap\mathfrak{T}\cap\tilde{\mathcal{T}}_{1;k,k'}^{*})\times (\mathcal{Y}_{k'}^{\circ}\cap\mathfrak{T}\cap\tilde{\mathcal{T}}_{1;k,k'}^{*})$, we define
\begin{align}\label{p3}
p_3(x,y):=\,& \mathbbm{1}_{\mathtt{x}'(x,y)\in\mathrm{int}(T(y))\cap\mathcal{X}_k^{\circ}}\mathbbm{1}_{\mathtt{y}'(x,y)\in\mathrm{int}(T(x))\cap\mathcal{Y}_{k'}^{\circ}}\mathbbm{1}_{V(x)+V(y)\neq 0}\nonumber\\
&\times p_2(\mathtt{x}'(x,y),\mathtt{y}'(x,y))\frac{\mathfrak{G}_{k,k'}(\mathtt{x}'(x,y))\mathfrak{G}_{k,k'}(\mathtt{y}'(x,y))}{\mathfrak{G}_{k,k'}(x)\mathfrak{G}_{k,k'}(y)},
\end{align}
\begin{align}\label{p4}
p_4(x,y):=\,&\mathbbm{1}_{\mathtt{x}'(x,y)\in\mathrm{int}(T(y))\cap\mathcal{X}_k^{\circ}}\mathbbm{1}_{\mathtt{y}'(x,y)\in\mathrm{int}(T(x))\cap\mathcal{Y}_{k'}^{\circ}}\mathbbm{1}_{V(x)+V(y)\neq 0}\nonumber\\
&\times  p_1(\mathtt{x}'(x,y),\mathtt{y}'(x,y))\frac{\mathfrak{G}_{k,k'}(\mathtt{x}'(x,y))\mathfrak{G}_{k,k'}(\mathtt{y}'(x,y))}{\mathfrak{G}_{k,k'}(x)\mathfrak{G}_{k,k'}(y)}.
\end{align}
For any $(x,y)\in(\mathbb{R}^d\times\mathbb{R}^d)\big\backslash\big((\mathcal{X}_k^{\circ}\cap\mathfrak{T}\cap\tilde{\mathcal{T}}_{1;k,k'}^{*})\times (\mathcal{Y}_{k'}^{\circ}\cap\mathfrak{T}\cap\tilde{\mathcal{T}}_{1;k,k'}^{*})\big)$, we define $p_3(x,y)=p_4(x,y):=0$. Note that $p_3(x,y),p_4(x,y)\geq 0$ for every $(x,y)\in\mathbb{R}^d\times\mathbb{R}^d$.    

Finally, for any $x,y\in\mathbb{R}^d$, we define
\begin{equation}\label{def:p}
    p(x,y):=(p_1\wedge p_2\wedge p_3\wedge p_4)(x,y)\geq 0.
\end{equation}
Note that by \eqref{p1} and \eqref{p2}, we have
\begin{equation}\label{domains}
    \{(x,y)\in\mathbb{R}^d\times \mathbb{R}^d:p(x,y)>0\}\subseteq (\mathcal{X}_k^{\circ}\cap\mathfrak{T}^*)\times(\mathcal{Y}_{k'}^{\circ}\cap\mathfrak{T}^*). 
\end{equation}

The following lemma shows that $p(x,y)$ and $p(\mathtt{x}'(x,y),\mathtt{y}'(x,y))$ coincide up to multiplicative factors. 

\begin{lemma}\label{L4.6}
For any $x,y\in\tilde{\mathcal{T}}_{1;k,k'}^{*}$ such that $\mathtt{x}'(x,y)\in \mathrm{int}(T(y))$ and $\mathtt{y}'(x,y)\in \mathrm{int}(T(x))$ (which implies $\mathtt{x}'(x,y),\mathtt{y}'(x,y)\in\tilde{\mathcal{T}}_{1;k,k'}^{*}$), we have
\begin{equation}\label{conl4.6}
   p(x,y)\mathfrak{G}_{k,k'}(x)\mathfrak{G}_{k,k'}(y) = p(\mathtt{x}'(x,y),\mathtt{y}'(x,y))\mathfrak{G}_{k,k'}(\mathtt{x}'(x,y))\mathfrak{G}_{k,k'}(\mathtt{y}'(x,y)).
\end{equation}
\end{lemma}
\begin{proof}
If $(\mathtt{x}'(x,y),\mathtt{y}'(x,y))\notin (\mathcal{X}_k^{\circ}\cap\mathfrak{T}\cap\tilde{\mathcal{T}}_{1;k,k'}^{*})\times (\mathcal{Y}_{k'}^{\circ}\cap\mathfrak{T}\cap\tilde{\mathcal{T}}_{1;k,k'}^{*})$, by \eqref{p1} and \eqref{p2}, we have $(p_1\wedge p_2)(\mathtt{x}'(x,y),\mathtt{y}'(x,y))=0$. Hence by \eqref{p3}--\eqref{def:p}, we have $p(x,y)=p(\mathtt{x}'(x,y),\mathtt{y}'(x,y))=0$, which establishes \eqref{conl4.6}.

Below we assume that $(\mathtt{x}'(x,y),\mathtt{y}'(x,y))\in (\mathcal{X}_k^{\circ}\cap\mathfrak{T}\cap\tilde{\mathcal{T}}_{1;k,k'}^{*})\times (\mathcal{Y}_{k'}^{\circ}\cap\mathfrak{T}\cap\tilde{\mathcal{T}}_{1;k,k'}^{*})$. As $\mathtt{x}'(x,y)\in \mathrm{int}(T(y))$ and $\mathtt{y}'(x,y)\in \mathrm{int}(T(x))$, we have $V(\mathtt{x}'(x,y))=V(y)$ and $V(\mathtt{y}'(x,y))=V(x)$. If $V(x)+V(y)=0$, then  $V(\mathtt{x}'(x,y))+V(\mathtt{y}'(x,y))=0$, hence by \eqref{p3} and \eqref{p4}, $(p_3\wedge p_4)(x,y)=(p_3\wedge p_4)(\mathtt{x}'(x,y),\mathtt{y}'(x,y))=0$, and so $p(x,y)=p(\mathtt{x}'(x,y),\mathtt{y}'(x,y))=0$. This yields \eqref{conl4.6}. Below we assume that $V(x)+V(y)\neq 0$. This implies $V(\mathtt{x}'(x,y))+V(\mathtt{y}'(x,y))\neq 0$, hence by \eqref{x'd} and \eqref{y'd}, 
\begin{equation}\label{neq4.1}
    \mathtt{x}'(\mathtt{x}'(x,y),\mathtt{y}'(x,y))=\mathtt{y}'(x,y)-\frac{\langle \mathtt{y}'(x,y)-\mathtt{x}'(x,y), V(x)+V(y)\rangle}{\langle V(x), V(x)+V(y)\rangle}V(x)=x,
\end{equation}
\begin{equation}\label{neq4.2}
    \mathtt{y}'(\mathtt{x}'(x,y),\mathtt{y}'(x,y))=\mathtt{x}'(x,y)+\frac{\langle \mathtt{y}'(x,y)-\mathtt{x}'(x,y), V(x)+V(y)\rangle}{\langle V(y), V(x)+V(y)\rangle}V(y)=y. 
\end{equation}

We first consider the case where $(x,y)\notin (\mathcal{X}_k^{\circ}\cap\mathfrak{T}\cap\tilde{\mathcal{T}}_{1;k,k'}^{*})\times (\mathcal{Y}_{k'}^{\circ}\cap\mathfrak{T}\cap\tilde{\mathcal{T}}_{1;k,k'}^{*})$. By \eqref{p1} and \eqref{p2}, $(p_1\wedge p_2)(x,y)=0$, hence $p(x,y)=0$. By \eqref{p3}--\eqref{p4} and \eqref{neq4.1}--\eqref{neq4.2}, we have $(p_3\wedge p_4)(\mathtt{x}'(x,y),\mathtt{y}'(x,y))=0$, hence $p(\mathtt{x}'(x,y),\mathtt{y}'(x,y))=0$. Thus \eqref{conl4.6} holds.  

Now we consider the case where $(x,y)\in (\mathcal{X}_k^{\circ}\cap\mathfrak{T}\cap\tilde{\mathcal{T}}_{1;k,k'}^{*})\times (\mathcal{Y}_{k'}^{\circ}\cap\mathfrak{T}\cap\tilde{\mathcal{T}}_{1;k,k'}^{*})$. By \eqref{p3} and \eqref{p4}, we have
\begin{equation}\label{EEq4.1}
     (p_3\wedge p_4)(x,y)=(p_1\wedge p_2)(\mathtt{x}'(x,y),\mathtt{y}'(x,y))\frac{\mathfrak{G}_{k,k'}(\mathtt{x}'(x,y))\mathfrak{G}_{k,k'}(\mathtt{y}'(x,y))}{\mathfrak{G}_{k,k'}(x)\mathfrak{G}_{k,k'}(y)}.
\end{equation}
Moreover, by \eqref{p3}--\eqref{p4} and \eqref{neq4.1}--\eqref{neq4.2}, we have
\begin{equation*}
    (p_3\wedge p_4)(\mathtt{x}'(x,y),\mathtt{y}'(x,y))=(p_1\wedge p_2)(x,y)\frac{\mathfrak{G}_{k,k'}(x)\mathfrak{G}_{k,k'}(y)}{\mathfrak{G}_{k,k'}(\mathtt{x}'(x,y))\mathfrak{G}_{k,k'}(\mathtt{y}'(x,y))},
\end{equation*}
and so
\begin{equation}\label{EEq4.2}
    (p_1\wedge p_2)(x,y)=(p_3\wedge p_4)(\mathtt{x}'(x,y),\mathtt{y}'(x,y))\frac{\mathfrak{G}_{k,k'}(\mathtt{x}'(x,y))\mathfrak{G}_{k,k'}(\mathtt{y}'(x,y))}{\mathfrak{G}_{k,k'}(x)\mathfrak{G}_{k,k'}(y)}.
\end{equation}
Combining \eqref{EEq4.1} and \eqref{EEq4.2}, we obtain that
\begin{equation*}
    p(x,y)=p(\mathtt{x}'(x,y),\mathtt{y}'(x,y))\frac{\mathfrak{G}_{k,k'}(\mathtt{x}'(x,y))\mathfrak{G}_{k,k'}(\mathtt{y}'(x,y))}{\mathfrak{G}_{k,k'}(x)\mathfrak{G}_{k,k'}(y)},
\end{equation*}
which yields \eqref{conl4.6}.
\end{proof}

\subsubsection{Definition of $r(\cdot,\cdot)$}\label{Sect.4.2.2}

Recall the definitions of $f_{k,k'}(\cdot)$ and $g_{k,k'}(\cdot)$ from \eqref{fkkdef} and \eqref{gkkdef}. For any $x\in\mathbb{R}^d$, we define
\begin{equation}\label{Eq.4.1new}
    \bar{f}_1(x):=(1-\delta_3)\int_{\mathbb{R}^d}p(x,z)dz,  \quad \bar{g}_1(x):=(1-\delta_3)\int_{\mathbb{R}^d}p(z,x)dz, 
\end{equation}
\begin{equation}\label{Eq.4.2new}
    \bar{f}_2(x):=f_{k,k'}(x)-\bar{f}_1(x), \quad \bar{g}_2(x):=g_{k,k'}(x)-\bar{g}_1(x).
\end{equation}

The following lemma provides a lower bound on $\bar{f}_2$ and $\bar{g}_2$ in terms of $f_{k,k'}$ and $g_{k,k'}$, respectively.    

\begin{lemma}\label{Lem4.7}
For every $x\in\mathbb{R}^d$, we have $\bar{f}_1(x)\leq (1-\delta_3) f_{k,k'}(x)$ and $\bar{g}_1(x)\leq (1-\delta_3) g_{k,k'}(x)$. Note that this implies $\bar{f}_2(x)\geq \delta_3 f_{k,k'}(x)$ and $\bar{g}_2(x)\geq \delta_3 g_{k,k'}(x)$ for every $x\in\mathbb{R}^d$. 
\end{lemma}
\begin{proof}
By \eqref{p1}, for any $x\in\mathcal{X}_k^{\circ}\cap\mathfrak{T}^*\cap\tilde{\mathcal{T}}_{1;k,k'}^{*}$, we have 
\begin{align*}
  \bar{f}_1(x)=\,&(1-\delta_3)\int_{\mathbb{R}^d}p(x,y)dy\leq (1-\delta_3)\int_{\mathbb{R}^d}p_1(x,y)dy\nonumber\\
 \leq\,& (1-\delta_3)\int_{\{z\in \mathrm{int}(T(x))\cap\mathcal{Y}_{k'}^{\circ}:\langle x-z,V(x)\rangle   >0\}}(2\pi\|x-z\|\eps)^{-(d-1)\slash 2}\mathfrak{F}(x)^{-1}h(x,z) d\mathcal{H}^1(z)\nonumber\\
 &\hspace{0.6in}\int_{O(V(x))}\sqrt{\det\big(\mathbf{I}_d+\|x-z\|F(z)\big)}e^{-\frac{1}{2\eps}w^{\top}\big(\frac{1}{\|x-z\|}\mathbf{I}_d+F(z)\big)w} d\mathcal{H}^{d-1}(w)\nonumber\\
 \leq\,& (1-\delta_3)\int_{\{z\in\mathrm{int}(T(x))\cap\mathcal{Y}_{k'}^{\circ}:\langle x-z,V(x)\rangle   >0\}}\mathfrak{F}(x)^{-1}h(x,z) d\mathcal{H}^1(z)\nonumber\\
 \leq\,& (1-\delta_3)\mathfrak{F}(x)^{-1}\int_{\mathrm{int}(T(x))\cap\mathcal{Y}_{k'}^{\circ}}h(x,z)d\mathcal{H}^1(z)
 = (1-\delta_3) f_{k,k'}(x),
\end{align*}
where we use \Cref{L2.7} (note the definition of $\mathfrak{T}$ in \Cref{Defn4.1.1n}) in the third inequality. For any $x\in\mathbb{R}^d\backslash(\mathcal{X}_k^{\circ}\cap\mathfrak{T}^*\cap\tilde{\mathcal{T}}_{1;k,k'}^{*})$, by \eqref{p1}, we have $p_1(x,y)=0$ for every $y\in\mathbb{R}^d$, hence
\begin{equation*}
    \bar{f}_1(x)\leq\int_{\mathbb{R}^d}p_1(x,y)dy=0\leq (1-\delta_3) f_{k,k'}(x).
\end{equation*}
Similarly, using \eqref{p2} we can deduce that $\bar{g}_1(x)\leq (1-\delta_3) g_{k,k'}(x)$ for every $x\in\mathbb{R}^d$.
\end{proof}

For any $T,T'\in\mathcal{S}$ such that $V_T+V_{T'}\neq 0$ and any $x\in T$, we define $\EuScript{R}(x;T,T')$ to be the unique point on $L_{T'}$ (recall Definition \ref{TransportRays}) such that
\begin{equation}\label{defvarphi}
    \langle \EuScript{R}(x;T,T')-x, V_T+V_{T'}\rangle =0.
\end{equation}
For any $T\in \mathcal{S}$ and $x\in \mathrm{int}(T)\cap\mathcal{X}_k^{\circ},y\in\mathrm{int}(T)\cap\mathcal{Y}_{k'}^{\circ}$, with $\mathfrak{F}(\cdot)$ as in \Cref{factors}, we define
\begin{equation}\label{defxiq}
    \xi_T(x,y):=\mathfrak{F}(x)\int_{\{T' \in\mathcal{S}: V_{T}+V_{T'}\neq 0,\EuScript{R}(y;T,T')\in\mathrm{int}(T')\cap\mathcal{Y}_{k'}^{\circ}\}} p(x,\EuScript{R}(y;T,T'))\mathfrak{F}(\EuScript{R}(y;T,T'))d\lambda(T').
\end{equation}

Now we define $\bar{\gamma}_{\eps, 1}$ and $\bar{\gamma}_{\eps, 2}$ to be (positive) Borel measures on $\mathbb{R}^d\times\mathbb{R}^d$ such that for any Borel set $A\subseteq \mathbb{R}^d\times\mathbb{R}^d$, 
\begin{equation}\label{def.gamma.1}
    \bar{\gamma}_{\eps, 1}(A)=(1-\delta_3)\int_{\mathcal{S}}d\lambda(T)\int_{(\mathrm{int}(T)\cap\mathcal{X}_k^{\circ})\times (\mathrm{int}(T)\cap\mathcal{Y}_{k'}^{\circ})}\mathbbm{1}_A(x,y)\xi_T(x,y)d\mathcal{H}^1(x)d\mathcal{H}^1(y),
\end{equation}
\begin{align}\label{def.gamma.2}
 \bar{\gamma}_{\eps, 2}(A)=\,&\int_{\mathcal{S}}\mathbbm{1}_{\int_{\mathrm{int}(T)\cap \mathcal{X}_k^{\circ}}\bar{f}_2(z)\mathfrak{F}(z)d\mathcal{H}^1(z)>0}\nonumber\\
 &\hspace{0.2in}\times\frac{\int_{(\mathrm{int}(T)\cap\mathcal{X}_k^{\circ})\times (\mathrm{int}(T)\cap\mathcal{Y}_{k'}^{\circ})}\mathbbm{1}_A(x,y)\bar{f}_2(x)\bar{g}_2(y)\mathfrak{F}(x)\mathfrak{F}(y) d\mathcal{H}^1(x)d\mathcal{H}^1(y)}{\int_{\mathrm{int}(T)\cap \mathcal{X}_k^{\circ}}\bar{f}_2(z)\mathfrak{F}(z)d\mathcal{H}^1(z)}
 d\lambda(T).
\end{align}

The next lemma shows that $\bar{\gamma}_{\eps,1}$ and $\bar{\gamma}_{\eps,2}$ only transport mass in the direction where $u$ decreases. 

\begin{lemma}\label{Lem4.6}
We have 
\begin{equation*}
    \bar{\gamma}_{\eps,1}(\{(x,y)\in\mathbb{R}^d\times\mathbb{R}^d:u(x)\leq u(y)\})=0, \qquad \bar{\gamma}_{\eps,2}(\{(x,y)\in\mathbb{R}^d\times\mathbb{R}^d:u(x)\leq u(y)\})=0.
\end{equation*}
\end{lemma}
\begin{proof}

Suppose that there exists $T\in\mathcal{S}$, such that there exist $x\in T\cap \mathcal{X}_k$ and $y\in T\cap \mathcal{Y}_{k'}$ satisfying $u(x)\leq u(y)$. As $x,y\in T$ and $\dist(\mathcal{X}_k, \mathcal{Y}_{k'})\geq 2d_0$, we have $u(y)-u(x)=\|y-x\|\geq 2d_0$. As $u\in \text{Lip}_1(\mathbb{R}^d)$, for any $x'\in \mathcal{X}_k$ and $y'\in \mathcal{Y}_{k'}$, we have (recall \eqref{eq:mathtthDefn})
\begin{align*}
    u(y')-u(x')\geq & u(y)-u(x)-|u(y')-u(y)|-|u(x')-u(x)|\geq u(y)-u(x)-\|y'-y\|-\|x'-x\|\nonumber\\
    \geq& 2d_0-2\sqrt{d}\mathtt{h}-2\sqrt{d}\mathtt{h}\geq 2d_0-4\sqrt{d}\cdot\frac{d_0}{10\sqrt{d}} \geq d_0,
\end{align*}
where we use the fact that $\mathtt{h}<d_0\slash (10\sqrt{d})$. Hence by \Cref{Lem2}(c) and Definition \ref{DefhH}, for any $x'\in\mathcal{X}_k$ and $y'\in\mathcal{Y}_{k'}$, we have $h(x',y')=0$. Therefore, by \eqref{Eq.4.2new} and \eqref{fkkdef}, for any $x'\in\mathcal{X}_k$, we have $\bar{f}_2(x')\leq f_{k,k'}(x')=0$. By \eqref{def.gamma.2}, this implies $\bar{\gamma}_{\eps,2}(\mathbb{R}^d\times\mathbb{R}^d)=0$. Moreover, by \eqref{p1}, for any $x,y\in\mathbb{R}^d$, we have $p(x,y)\leq p_1(x,y)=0$ (in order for $p_1(x,y)>0$, we must have $x\in\mathcal{X}_k^{\circ}$ and $\mathtt{z}(x,y)\in \mathcal{Y}_{k'}^{\circ}$, but then $h(x,\mathtt{z}(x,y))=0$ yielding $p_1(x,y)=0$). Hence by \eqref{defxiq}, $\xi_T(x,y)=0$ for any $T\in\mathcal{S}$ and $x\in \mathrm{int}(T)\cap\mathcal{X}_k^{\circ},y\in\mathrm{int}(T)\cap\mathcal{Y}_{k'}^{\circ}$. By \eqref{def.gamma.1}, this implies $\bar{\gamma}_{\eps,1}(\mathbb{R}^d\times\mathbb{R}^d)=0$.

Below we assume that for any $T\in\mathcal{S}$ and any $x\in T\cap\mathcal{X}_k,y\in T\cap\mathcal{Y}_{k'}$, we have $u(x)>u(y)$. Let $A:=\{(x,y)\in\mathbb{R}^d\times\mathbb{R}^d:u(x)\leq u(y)\}$. Note that for any $T\in\mathcal{S}$,
\begin{align*}
   & \int_{(\mathrm{int}(T)\cap\mathcal{X}_k^{\circ})\times (\mathrm{int}(T)\cap\mathcal{Y}_{k'}^{\circ})}\mathbbm{1}_A(x,y)\xi_T(x,y)d\mathcal{H}^1(x)d\mathcal{H}^1(y)\nonumber\\
   =\,&\int_{(\mathrm{int}(T)\cap\mathcal{X}_k^{\circ})\times (\mathrm{int}(T)\cap\mathcal{Y}_{k'}^{\circ})}\mathbbm{1}_{u(x)\leq u(y)}\xi_T(x,y)d\mathcal{H}^1(x)d\mathcal{H}^1(y)=0,
\end{align*}
\begin{align*}
   & \int_{(\mathrm{int}(T)\cap\mathcal{X}_k^{\circ})\times (\mathrm{int}(T)\cap\mathcal{Y}_{k'}^{\circ})}\mathbbm{1}_A(x,y)\bar{f}_2(x)\bar{g}_2(y)\mathfrak{F}(x)\mathfrak{F}(y) d\mathcal{H}^1(x)d\mathcal{H}^1(y)\nonumber\\
   =\,& \int_{(\mathrm{int}(T)\cap\mathcal{X}_k^{\circ})\times (\mathrm{int}(T)\cap\mathcal{Y}_{k'}^{\circ})}\mathbbm{1}_{u(x)\leq u(y)}\bar{f}_2(x)\bar{g}_2(y)\mathfrak{F}(x)\mathfrak{F}(y) d\mathcal{H}^1(x)d\mathcal{H}^1(y)=0.
\end{align*}
Hence by \eqref{def.gamma.1} and \eqref{def.gamma.2}, $\bar{\gamma}_{\eps,1}(A)=\bar{\gamma}_{\eps,2}(A)=0$.
\end{proof}

The following lemma shows that the marginal densities of $\bar{\gamma}_{\eps,1}$ are given by $\bar{f}_1$ and $\bar{g}_1$.

\begin{lemma}\label{Lem4.11}
For any Borel set $A\subseteq\mathbb{R}^d$, we have
\begin{eqnarray*}
    \bar{\gamma}_{\eps,1}(A\times \mathbb{R}^d)=\int_A \bar{f}_1(x)dx,\qquad  \bar{\gamma}_{\eps,1}(\mathbb{R}^d\times A)=\int_A \bar{g}_1(y)dy. 
\end{eqnarray*}
\end{lemma}
\begin{proof}
For any $T\in\mathcal{S}$ and any Borel set $A\subseteq\mathbb{R}^d$, we have 
\begin{align}\label{neq4.3}
   & \int_{(\mathrm{int}(T)\cap\mathcal{X}_k^{\circ})\times (\mathrm{int}(T)\cap\mathcal{Y}_{k'}^{\circ})}\mathbbm{1}_A(x)\xi_T(x,y)d\mathcal{H}^1(x)d\mathcal{H}^1(y)\nonumber\\
   =\,&  \int_{\{T'\in\mathcal{S}: V_{T}+V_{T'}\neq 0\}} d\lambda(T')\int_{(\mathrm{int}(T)\cap\mathcal{X}_k^{\circ})\times (\mathrm{int}(T)\cap\mathcal{Y}_{k'}^{\circ})}\mathbbm{1}_A(x) \mathbbm{1}_{\EuScript{R}(y;T,T')\in\mathrm{int}(T')\cap\mathcal{Y}_{k'}^{\circ}}p(x,\EuScript{R}(y;T,T'))
   \nonumber\\
   &\hspace{3.1in}\times\mathfrak{F}(x)\mathfrak{F}(\EuScript{R}(y;T,T'))d\mathcal{H}^1(x)d\mathcal{H}^1(y)\nonumber\\
   =\,&\int_{\{T'\in\mathcal{S}: V_{T}+V_{T'}\neq 0\}} d\lambda(T')\int_{(\mathrm{int}(T)\cap\mathcal{X}_k^{\circ})\times (\mathrm{int}(T')\cap\mathcal{Y}_{k'}^{\circ})}\mathbbm{1}_A(x) p(x,y')\mathfrak{F}(x)\mathfrak{F}(y')d\mathcal{H}^1(x)d\mathcal{H}^1(y')\nonumber\\
   =\,& \int_{\mathrm{int}(T)\cap\mathcal{X}_k^{\circ}} \mathbbm{1}_A(x)\mathfrak{F}(x)  d\mathcal{H}^1(x)
   \int_{\{T'\in\mathcal{S}: V_{T}+V_{T'}\neq 0\}} d\lambda(T')\int_{\mathrm{int}(T')\cap\mathcal{Y}_{k'}^{\circ}} p(x,y')\mathfrak{F}(y')d\mathcal{H}^1(y')\nonumber\\
   =\,&  \int_{\mathrm{int}(T)\cap\mathcal{X}_k^{\circ}} \mathbbm{1}_A(x)\mathfrak{F}(x)  d\mathcal{H}^1(x)
   \int_{\mathcal{S}} d\lambda(T')\int_{\mathrm{int}(T')\cap\mathcal{Y}_{k'}^{\circ}} p(x,y')\mathfrak{F}(y')d\mathcal{H}^1(y')\nonumber\\
   =\,& \int_{\mathrm{int}(T)\cap\mathcal{X}_k^{\circ}} \mathbbm{1}_A(x)\mathfrak{F}(x)\bigg(\int_{\mathcal{Y}_{k'}^{\circ}\cap\mathfrak{T}^*}p(x,y')dy'\bigg)d\mathcal{H}^1(x)\nonumber\\
   =\,&\int_{\mathrm{int}(T)\cap\mathcal{X}_k^{\circ}} \mathbbm{1}_A(x)\mathfrak{F}(x)\bigg(\int_{\mathbb{R}^d}p(x,y)dy\bigg)d\mathcal{H}^1(x),
\end{align}
where the second equality uses the fact that for any $T,T'\in \mathcal{S}$ such that $V_T+V_{T'}\neq 0$ and any $x\in\mathrm{int}(T)\cap\mathcal{X}_{k}^{\circ},y'\in \mathrm{int}(T')\cap\mathcal{Y}_{k'}^{\circ}$, if $\mathtt{y}'(x,y')\notin \mathrm{int}(T)\cap\mathcal{Y}_{k'}^{\circ}$, then by \eqref{p3}, $p(x,y')\leq p_3(x,y')=0$. Moreover, the fourth equality uses the fact that for any $T,T'\in \mathcal{S}$ and $x\in\mathrm{int}(T)\cap\mathcal{X}_k^{\circ},y'\in \mathrm{int}(T')\cap\mathcal{Y}_{k'}^{\circ}$, if $V_T+V_{T'}=0$, then $p(x,y')=0$ (by \eqref{p3}); the last two equalities use Lemma \ref{Lem3.18nn} and \eqref{domains}. Hence by Lemma \ref{Lem3.18nn} and \eqref{domains}, for any Borel set $A\subseteq\mathbb{R}^d$,
\begin{align*}
 & \bar{\gamma}_{\eps,1}(A\times\mathbb{R}^d)=(1-\delta_3)\int_{\mathcal{S}}d\lambda(T)
    \int_{(\mathrm{int}(T)\cap\mathcal{X}_k^{\circ})\times (\mathrm{int}(T)\cap\mathcal{Y}_{k'}^{\circ})}\mathbbm{1}_A(x)\xi_T(x,y)d\mathcal{H}^1(x)d\mathcal{H}^1(y)\nonumber\\
    =\,&
    (1-\delta_3)\int_{\mathcal{S}}d\lambda(T)\int_{\mathrm{int}(T)\cap\mathcal{X}_k^{\circ}} \mathbbm{1}_A(x)\mathfrak{F}(x)\bigg(\int_{\mathbb{R}^d}p(x,y)dy\bigg)d\mathcal{H}^1(x)\nonumber\\
    =\,& (1-\delta_3)\int_{\mathcal{X}_k^{\circ}\cap\mathfrak{T}^*}\mathbbm{1}_A(x)\bigg(\int_{\mathbb{R}^d}p(x,y)dy\bigg)dx  
    = (1-\delta_3)\int_{A\times\mathbb{R}^d}p(x,y)dxdy=\int_A \bar{f}_1(x)dx.
\end{align*}

Consider any $T,T'\in\mathcal{S}$ such that $V_T+V_{T'}\neq 0$ and any $x\in \mathrm{int}(T)\cap\mathcal{X}_k^{\circ},y\in\mathrm{int}(T)\cap\mathcal{Y}_{k'}^{\circ}$. Note that the existence of $x,y$ implies that $\mathrm{int}(T)\cap\mathcal{X}_k^{\circ},\mathrm{int}(T)\cap\mathcal{Y}_{k'}^{\circ}\neq\emptyset$, hence $\mathrm{int}(T)\subseteq \tilde{\mathcal{T}}_{1;k,k'}^{*}$ (recall \eqref{deftildeTkk}) and $x,y\in  \tilde{\mathcal{T}}_{1;k,k'}^{*}$. First consider the case where $\EuScript{R}(x;T,T')\in\mathrm{int}(T')\cap\mathcal{X}_k^{\circ}, \EuScript{R}(y;T,T')\in\mathrm{int}(T')\cap\mathcal{Y}_{k'}^{\circ}$. Note that this implies $\mathrm{int}(T')\cap\mathcal{X}_k^{\circ}, \mathrm{int}(T')\cap\mathcal{Y}_{k'}^{\circ}\neq\emptyset$, hence $\mathrm{int}(T')\subseteq \tilde{\mathcal{T}}^{*}_{1;k,k'}$ and $\EuScript{R}(x;T,T'), \EuScript{R}(y;T,T')  \in  \tilde{\mathcal{T}}^{*}_{1;k,k'}$. Hence by Lemma \ref{L4.6}, we have  
\begin{equation}\label{L4.6nne1}
    p(x,\EuScript{R}(y;T,T'))\mathfrak{G}_{k,k'}(x)\mathfrak{G}_{k,k'}(\EuScript{R}(y;T,T')) = p(\EuScript{R}(x;T,T'),y)\mathfrak{G}_{k,k'}(\EuScript{R}(x;T,T'))\mathfrak{G}_{k,k'}(y).
\end{equation}
Now if $T,T'\in\mathcal{S}^*$, then by Definition \ref{factors} and \eqref{def_Gkk}, we have (note that $\mathfrak{q}_{k,k'}(x)=\mathfrak{q}_{k,k'}(y)$ and $\mathfrak{q}_{k,k'}(\EuScript{R}(x;T,T'))=\mathfrak{q}_{k,k'}(\EuScript{R}(y;T,T'))$)
\begin{align}\label{L4.6nne2}
    \mathfrak{F}(x)&=\frac{\mathfrak{G}_{k,k'}(x)}{\int_{T}f(z)\det(\mathbf{I}_d+\langle z-\mathfrak{q}_{k,k'}(x),V_T\rangle F(\mathfrak{q}_{k,k'}(x)))d\mathcal{H}^1(z)},\nonumber\\
     \mathfrak{F}(y)&=\frac{\mathfrak{G}_{k,k'}(y)}{\int_{T}f(z)\det(\mathbf{I}_d+\langle z-\mathfrak{q}_{k,k'}(x),V_T\rangle F(\mathfrak{q}_{k,k'}(x)))d\mathcal{H}^1(z)},\nonumber\\
     \mathfrak{F}(\EuScript{R}(x;T,T'))&=\frac{\mathfrak{G}_{k,k'}(\EuScript{R}(x;T,T'))}{\int_{T'}f(z)\det(\mathbf{I}_d+\langle z-\mathfrak{q}_{k,k'}(\EuScript{R}(x;T,T')),V_{T'}\rangle F(\mathfrak{q}_{k,k'}(\EuScript{R}(x;T,T'))))d\mathcal{H}^1(z)},\nonumber\\
      \mathfrak{F}(\EuScript{R}(y;T,T'))&=\frac{\mathfrak{G}_{k,k'}(\EuScript{R}(y;T,T'))}{\int_{T'}f(z)\det(\mathbf{I}_d+\langle z-\mathfrak{q}_{k,k'}(\EuScript{R}(x;T,T')),V_{T'}\rangle F(\mathfrak{q}_{k,k'}(\EuScript{R}(x;T,T'))))d\mathcal{H}^1(z)}.
\end{align}
Combining \eqref{L4.6nne1} and \eqref{L4.6nne2}, we get
\begin{equation}\label{nneq4.1}
    p(x,\EuScript{R}(y;T,T'))\mathfrak{F}(x)\mathfrak{F}(\EuScript{R}(y;T,T'))= p(\EuScript{R}(x;T,T'),y)\mathfrak{F}(\EuScript{R}(x;T,T'))\mathfrak{F}(y). 
\end{equation}
If at least one of $T$ and $T'$ is not in $\mathcal{S}^*$, then by Remark \ref{positivityF}, either $\mathfrak{F}(x)=\mathfrak{F}(y)=0$ or $\mathfrak{F}(\EuScript{R}(x;T,T'))=\mathfrak{F}(\EuScript{R}(y;T,T'))=0$, hence \eqref{nneq4.1} also holds. If $\EuScript{R}(x;T,T')\in\mathrm{int}(T')\cap\mathcal{X}_{k}^{\circ}$ and $\EuScript{R}(y;T,T')\notin\mathrm{int}(T')\cap\mathcal{Y}_{k'}^{\circ}$, then by \eqref{p3}, we have
\begin{equation}
    p(\EuScript{R}(x;T,T'),y)\leq p_3(\EuScript{R}(x;T,T'),y) =0.
\end{equation}
If $\EuScript{R}(x;T,T')\notin\mathrm{int}(T')\cap\mathcal{X}_{k}^{\circ}$ and $\EuScript{R}(y;T,T')\in\mathrm{int}(T')\cap\mathcal{Y}_{k'}^{\circ}$, then by \eqref{p3}, we have
\begin{equation}\label{nneq4.2}
 p(x,\EuScript{R}(y;T,T'))\leq p_3(x,\EuScript{R}(y;T,T'))  =0.
\end{equation}
By \eqref{nneq4.1}--\eqref{nneq4.2}, noting \eqref{defxiq}, for any $T\in\mathcal{S}$ and $x\in \mathrm{int}(T)\cap\mathcal{X}_k^{\circ},y\in\mathrm{int}(T)\cap\mathcal{Y}_{k'}^{\circ}$, we get 
\begin{equation*}
 \xi_T(x,y) =  \mathfrak{F}(y)\int_{\{T'\in \mathcal{S}: V_{T}+V_{T'}\neq 0,\EuScript{R}(x;T,T')\in\mathrm{int}(T')\cap\mathcal{X}_{k}^{\circ}\}}p(\EuScript{R}(x;T,T'),y)\mathfrak{F}(\EuScript{R}(x;T,T'))d\lambda(T').
\end{equation*}
Hence for any $T\in\mathcal{S}$ and any Borel set $A\subseteq\mathbb{R}^d$, following a similar argument as in \eqref{neq4.3}, we get
\begin{align}\label{meq4.4}
& \int_{(\mathrm{int}(T)\cap\mathcal{X}_k^{\circ})\times (\mathrm{int}(T)\cap\mathcal{Y}_{k'}^{\circ})}\mathbbm{1}_A(y)\xi_T(x,y)d\mathcal{H}^1(x)d\mathcal{H}^1(y)\nonumber\\
=\,& \int_{\{T'\in\mathcal{S}: V_{T}+V_{T'}\neq 0\}} d\lambda(T')\int_{(\mathrm{int}(T)\cap\mathcal{X}_k^{\circ})\times (\mathrm{int}(T)\cap\mathcal{Y}_{k'}^{\circ})}\mathbbm{1}_A(y) \mathbbm{1}_{\EuScript{R}(x;T,T')\in\mathrm{int}(T')\cap\mathcal{X}_{k}^{\circ}}p(\EuScript{R}(x;T,T'),y)
   \nonumber\\
   &\hspace{3.1in}\times\mathfrak{F}(y)\mathfrak{F}(\EuScript{R}(x;T,T'))d\mathcal{H}^1(x)d\mathcal{H}^1(y)\nonumber\\
   =\,& \int_{\{T'\in\mathcal{S}: V_{T}+V_{T'}\neq 0\}} d\lambda(T')\int_{(\mathrm{int}(T')\cap\mathcal{X}_k^{\circ})\times (\mathrm{int}(T)\cap\mathcal{Y}_{k'}^{\circ})}\mathbbm{1}_A(y) p(x',y) \mathfrak{F}(y)\mathfrak{F}(x')d\mathcal{H}^1(x')d\mathcal{H}^1(y)\nonumber\\
   =\,& \int_{\mathrm{int}(T)\cap\mathcal{Y}_{k'}^{\circ}} \mathbbm{1}_A(y)\mathfrak{F}(y)  d\mathcal{H}^1(y)
   \int_{\{T'\in\mathcal{S}: V_{T}+V_{T'}\neq 0\}} d\lambda(T')\int_{\mathrm{int}(T')\cap\mathcal{X}_{k}^{\circ}} p(x',y)\mathfrak{F}(x')d\mathcal{H}^1(x')\nonumber\\
   =\,&  \int_{\mathrm{int}(T)\cap\mathcal{Y}_{k'}^{\circ}} \mathbbm{1}_A(y)\mathfrak{F}(y)  d\mathcal{H}^1(y)
   \int_{\mathcal{S}} d\lambda(T')\int_{\mathrm{int}(T')\cap\mathcal{X}_{k}^{\circ}} p(x',y)\mathfrak{F}(x')d\mathcal{H}^1(x'),\nonumber\\
   =\,& \int_{\mathrm{int}(T)\cap\mathcal{Y}_{k'}^{\circ}} \mathbbm{1}_A(y)\mathfrak{F}(y)\bigg(\int_{\mathcal{X}_{k}^{\circ}\cap\mathfrak{T}^*}p(x',y)dx'\bigg)d\mathcal{H}^1(y)\nonumber\\
   =\,&\int_{\mathrm{int}(T)\cap\mathcal{Y}_{k'}^{\circ}} \mathbbm{1}_A(y)\mathfrak{F}(y)\bigg(\int_{\mathbb{R}^d}p(x,y)dx\bigg)d\mathcal{H}^1(y).
\end{align}
Consequently, by Lemma \ref{Lem3.18nn} and \eqref{domains}, for any Borel set $A\subseteq\mathbb{R}^d$, we have
\begin{align*}
  & \bar{\gamma}_{\eps,1}(\mathbb{R}^d\times A)=(1-\delta_3)\int_{\mathcal{S}}d\lambda(T)
    \int_{(\mathrm{int}(T)\cap\mathcal{X}_k^{\circ})\times (\mathrm{int}(T)\cap\mathcal{Y}_{k'}^{\circ})}\mathbbm{1}_A(y)\xi_T(x,y)d\mathcal{H}^1(x)d\mathcal{H}^1(y)\nonumber\\
    =\,&
    (1-\delta_3)\int_{\mathcal{S}}d\lambda(T)\int_{\mathrm{int}(T)\cap\mathcal{Y}_{k'}^{\circ}} \mathbbm{1}_A(y)\mathfrak{F}(y)\bigg(\int_{\mathbb{R}^d}p(x,y)dx\bigg)d\mathcal{H}^1(y)\nonumber\\
    =\,& (1-\delta_3)\int_{\mathcal{Y}_{k'}^{\circ}\cap\mathfrak{T}^*}\mathbbm{1}_A(y)\bigg(\int_{\mathbb{R}^d}p(x,y)dx\bigg)dy
    =(1-\delta_3)\int_{\mathbb{R}^d\times A}p(x,y)dxdy=\int_A \bar{g}_1(x)dx.\qedhere
\end{align*}
\end{proof} 

The next lemma shows that the pairs of densities $(\bar{f}_1,\bar{g}_1)$ and $(\bar{f}_2,\bar{g}_2)$, weighted by $\mathfrak{F}$, are balanced along each transport ray. 

\begin{lemma}\label{Lem4.7nn}
For any $T\in\mathcal{S}$, we have
\begin{equation*}
    \int_{T}\bar{f}_1(x)\mathfrak{F}(x)d\mathcal{H}^1(x) = \int_{T}\bar{g}_1(x)\mathfrak{F}(x)d\mathcal{H}^1(x),\quad \int_{T}\bar{f}_2(x)\mathfrak{F}(x)d\mathcal{H}^1(x) = \int_{T}\bar{g}_2(x)\mathfrak{F}(x)d\mathcal{H}^1(x).
\end{equation*}
\end{lemma}
\begin{proof}

Recall Definition \ref{factors}. For any $T\in\mathcal{S}\backslash\mathcal{S}^*$ and $x\in T$, we have $x\notin \mathfrak{T}^*$, hence by Remark \ref{positivityF}, $\mathfrak{F}(x)=0$. Therefore, for any $T\in\mathcal{S}\backslash\mathcal{S}^*$, $\int_{T}\bar{f}_1(x)\mathfrak{F}(x)d\mathcal{H}^1(x) = \int_{T}\bar{g}_1(x)\mathfrak{F}(x)d\mathcal{H}^1(x)=0$ and $\int_{T}\bar{f}_2(x)\mathfrak{F}(x)d\mathcal{H}^1(x) = \int_{T}\bar{g}_2(x)\mathfrak{F}(x)d\mathcal{H}^1(x)=0$. 

Below, we consider any $T\in\mathcal{S}^*$. Note that $\mathrm{int}(T)\subseteq \mathfrak{T}^*$. By \eqref{fkkdef} and \eqref{gkkdef}, we have  
\begin{align*}
    \int_{T}f_{k,k'}(x)\mathfrak{F}(x)d\mathcal{H}^1(x)&=\int_{\mathrm{int}(T)} \mathbbm{1}_{\mathcal{X}_k^{\circ}\cap\tilde{\mathcal{T}}_{1;k,k'}^*}(x)
     d\mathcal{H}^1(x)\int_{\mathrm{int}(T)\cap \mathcal{Y}_{k'}^{\circ}}h(x,z)d\mathcal{H}^1(z)\nonumber\\
     &=\int_{(\mathrm{int}(T)\cap\mathcal{X}_k^{\circ})\times (\mathrm{int}(T)\cap \mathcal{Y}_{k'}^{\circ})} \mathbbm{1}_{\tilde{\mathcal{T}}_{1;k,k'}^*}(x)h(x,z)
     d\mathcal{H}^1(x)d\mathcal{H}^1(z),
\end{align*}
\begin{align*}
    \int_{T}g_{k,k'}(x)\mathfrak{F}(x)d\mathcal{H}^1(x)&=\int_{\mathrm{int}(T)}\mathbbm{1}_{\mathcal{Y}_{k'}^{\circ}\cap\tilde{\mathcal{T}}_{1;k,k'}^*}(x)d\mathcal{H}^1(x)\int_{\mathrm{int}(T)\cap \mathcal{X}_{k}^{\circ}}h(z,x)d\mathcal{H}^1(z)\nonumber\\
    &=\int_{(\mathrm{int}(T)\cap\mathcal{X}_k^{\circ})\times (\mathrm{int}(T)\cap \mathcal{Y}_{k'}^{\circ})} \mathbbm{1}_{\tilde{\mathcal{T}}_{1;k,k'}^*}(z)h(x,z)
     d\mathcal{H}^1(x)d\mathcal{H}^1(z).
\end{align*}
As $\tilde{\mathcal{T}}^{*}_{1;k,k'}$ is a transport set (recall \eqref{deftildeTkk}), for any $x,z\in\mathrm{int}(T)$, we have $\mathbbm{1}_{\tilde{\mathcal{T}}^{*}_{1;k,k'}}(x)=\mathbbm{1}_{\tilde{\mathcal{T}}^{*}_{1;k,k'}}(z)$. Hence 
\begin{equation}\label{Lem4.7n1}
    \int_{T}f_{k,k'}(x)\mathfrak{F}(x)d\mathcal{H}^1(x)=\int_{T}g_{k,k'}(x)\mathfrak{F}(x)d\mathcal{H}^1(x).
\end{equation}
By \eqref{neq4.3} and \eqref{meq4.4} (taking $A=\mathbb{R}^d$), and noting \eqref{Eq.4.1new}, we have
\begin{align}\label{Lem4.7n2}
    \int_{T}\bar{f}_1(x)\mathfrak{F}(x)d\mathcal{H}^1(x)=\,&(1-\delta_3)\int_{\mathrm{int}(T)}\mathfrak{F}(x)\bigg(\int_{\mathbb{R}^d}p(x,y)dy\bigg)d\mathcal{H}^1(x)\nonumber\\
    =\,&(1-\delta_3)\int_{\mathrm{int}(T)\cap\mathcal{X}_k^{\circ}}\mathfrak{F}(x)\bigg(\int_{\mathbb{R}^d}p(x,y)dy\bigg)d\mathcal{H}^1(x)\nonumber\\
    =\,&(1-\delta_3)\int_{(\mathrm{int}(T)\cap\mathcal{X}_k^{\circ})\times (\mathrm{int}(T)\cap\mathcal{Y}_{k'}^{\circ})}\xi_T(x,y)d\mathcal{H}^1(x)d\mathcal{H}^1(y)\nonumber\\
    =\,& (1-\delta_3)\int_{\mathrm{int}(T)\cap\mathcal{Y}_{k'}^{\circ}} \mathfrak{F}(y)\bigg(\int_{\mathbb{R}^d}p(x,y)dx\bigg)d\mathcal{H}^1(y)\nonumber\\
    =\,&(1-\delta_3)\int_{\mathrm{int}(T)} \mathfrak{F}(y)\bigg(\int_{\mathbb{R}^d}p(x,y)dx\bigg)d\mathcal{H}^1(y)=\int_T \bar{g}_1(y)\mathfrak{F}(y)d\mathcal{H}^1(y),
\end{align}
where we use \eqref{domains} in the second and fifth equalities. The conclusion of the lemma now follows from \eqref{Eq.4.2new}, \eqref{Lem4.7n1}, and \eqref{Lem4.7n2}.
\end{proof}

The following lemma shows that the marginal densities of $\bar{\gamma}_{\eps,2}$ are given by $\bar{f}_2$ and $\bar{g}_2$.

\begin{lemma}\label{Lem4.9}
For any Borel set $A\subseteq \mathbb{R}^d$, we have
\begin{equation*}
   \bar{\gamma}_{\eps,2}(A\times\mathbb{R}^d)=\int_A \bar{f}_2(x) dx,\qquad \bar{\gamma}_{\eps,2}(\mathbb{R}^d\times A)=\int_A\bar{g}_2(x)dx.
\end{equation*}
\end{lemma}
\begin{proof}

By Lemma \ref{Lem4.7}, for any $x\in\mathbb{R}^d$, $\max\{\bar{f}_1(x),\bar{f}_2(x)\}\leq f_{k,k'}(x)$ and $\max\{\bar{g}_1(x),\bar{g}_2(x)\}\leq g_{k,k'}(x)$. Hence by \eqref{fkkdef} and \eqref{gkkdef},
\begin{equation}\label{subinclud}
    \big\{x\in\mathbb{R}^d:\bar{f}_2(x)>0\big\}\subseteq \mathcal{X}_k^{\circ}, \qquad \big\{x\in\mathbb{R}^d:\bar{g}_2(x)>0\big\}\subseteq \mathcal{Y}_{k'}^{\circ}. 
\end{equation}
By Lemmas \ref{Lem3.15} and \ref{Ln4.1},
\begin{equation}\label{subinclud1}
    \mu_{k,k'}(\mathbb{R}^d\backslash\mathfrak{T}^*)\leq \mu(\mathbb{R}^d\backslash\mathfrak{T}^*)=0, \qquad \nu_{k,k'}(\mathbb{R}^d\backslash\mathfrak{T}^*)\leq \nu(\mathbb{R}^d\backslash\mathfrak{T}^*)=0.
\end{equation}

By \eqref{def.gamma.2}, for any Borel set $A\subseteq \mathbb{R}^d$, we have
\begin{align*}
     \bar{\gamma}_{\eps,2}(\mathbb{R}^d\times A)
    =\,& \int_{\mathcal{S}}\mathbbm{1}_{\int_{\mathrm{int}(T)\cap \mathcal{X}_k^{\circ}}\bar{f}_2(z)\mathfrak{F}(z)d\mathcal{H}^1(z)>0}\nonumber\\
    & \hspace{0.2in}\times\frac{\int_{(\mathrm{int}(T)\cap\mathcal{X}_k^{\circ})\times (\mathrm{int}(T)\cap\mathcal{Y}_{k'}^{\circ})}\mathbbm{1}_A(y)\bar{f}_2(x)\bar{g}_2(y)\mathfrak{F}(x)\mathfrak{F}(y) d\mathcal{H}^1(x)d\mathcal{H}^1(y)}{\int_{\mathrm{int}(T)\cap \mathcal{X}_k^{\circ}}\bar{f}_2(z)\mathfrak{F}(z)d\mathcal{H}^1(z)}
 d\lambda(T)\nonumber\\
 =\,& \int_{\mathcal{S}}\mathbbm{1}_{\int_{\mathrm{int}(T)\cap \mathcal{X}_k^{\circ}}\bar{f}_2(z)\mathfrak{F}(z)d\mathcal{H}^1(z)>0}d\lambda(T)\int_{\mathrm{int}(T)\cap\mathcal{Y}_{k'}^{\circ}}\mathbbm{1}_A(y)\bar{g}_2(y)\mathfrak{F}(y)d\mathcal{H}^1(y)\nonumber\\
 =\,& \int_{\mathcal{S}}d\lambda(T)\int_{\mathrm{int}(T)\cap\mathcal{Y}_{k'}^{\circ}}\mathbbm{1}_A(y)\bar{g}_2(y)\mathfrak{F}(y)d\mathcal{H}^1(y)\nonumber\\
 =\,&\int_{\mathcal{S}}d\lambda(T)\int_{T}\mathbbm{1}_A(y)\bar{g}_2(y)\mathfrak{F}(y)d\mathcal{H}^1(y)
 = \int_{A\cap\mathfrak{T}^*}\bar{g}_2(y)dy = \int_A \bar{g}_2(y)dy,
\end{align*}
where in the third equality we note that by \eqref{subinclud} and Lemma \ref{Lem4.7nn}, for any $T\in\mathcal{S}$, 
\begin{equation}\label{subinclud1.1}
    \int_{\mathrm{int}(T)\cap \mathcal{X}_k^{\circ}}\bar{f}_2(z)\mathfrak{F}(z)d\mathcal{H}^1(z)=\int_{\mathrm{int}(T)\cap\mathcal{Y}_{k'}^{\circ}}\bar{g}_2(z)\mathfrak{F}(z)d\mathcal{H}^1(z).
\end{equation}
Moreover, in the last two equalities we use Lemma \ref{Lem3.18nn} and \eqref{subinclud1}. Similarly, using \eqref{def.gamma.2} and \eqref{subinclud1.1}, we obtain that
\begin{align*}
     \bar{\gamma}_{\eps,2}(A\times\mathbb{R}^d)
    =\,& \int_{\mathcal{S}}\mathbbm{1}_{\int_{\mathrm{int}(T)\cap \mathcal{X}_k^{\circ}}\bar{f}_2(z)\mathfrak{F}(z)d\mathcal{H}^1(z)>0}\nonumber\\
    & \hspace{0.2in}\times\frac{\int_{(\mathrm{int}(T)\cap\mathcal{X}_k^{\circ})\times (\mathrm{int}(T)\cap\mathcal{Y}_{k'}^{\circ})}\mathbbm{1}_A(x)\bar{f}_2(x)\bar{g}_2(y)\mathfrak{F}(x)\mathfrak{F}(y) d\mathcal{H}^1(x)d\mathcal{H}^1(y)}{\int_{\mathrm{int}(T)\cap \mathcal{X}_k^{\circ}}\bar{f}_2(z)\mathfrak{F}(z)d\mathcal{H}^1(z)}
 d\lambda(T)\nonumber\\
 =\,& \int_{\mathcal{S}}\mathbbm{1}_{\int_{\mathrm{int}(T)\cap \mathcal{X}_k^{\circ}}\bar{f}_2(z)\mathfrak{F}(z)d\mathcal{H}^1(z)>0}\nonumber\\
    & \hspace{0.2in}\times\frac{\int_{(\mathrm{int}(T)\cap\mathcal{X}_k^{\circ})\times (\mathrm{int}(T)\cap\mathcal{Y}_{k'}^{\circ})}\mathbbm{1}_A(x)\bar{f}_2(x)\bar{g}_2(y)\mathfrak{F}(x)\mathfrak{F}(y) d\mathcal{H}^1(x)d\mathcal{H}^1(y)}{\int_{\mathrm{int}(T)\cap \mathcal{Y}_{k'}^{\circ}}\bar{g}_2(z)\mathfrak{F}(z)d\mathcal{H}^1(z)}
 d\lambda(T)\nonumber\\
 =\,& \int_{\mathcal{S}}\mathbbm{1}_{\int_{\mathrm{int}(T)\cap \mathcal{X}_k^{\circ}}\bar{f}_2(z)\mathfrak{F}(z)d\mathcal{H}^1(z)>0}d\lambda(T)\int_{\mathrm{int}(T)\cap\mathcal{X}_k^{\circ}}\mathbbm{1}_A(x)\bar{f}_2(x)\mathfrak{F}(x)d\mathcal{H}^1(x)\nonumber\\
 =\,& \int_{\mathcal{S}}d\lambda(T)\int_{\mathrm{int}(T)\cap\mathcal{X}_{k}^{\circ}}\mathbbm{1}_A(x)\bar{f}_2(x)\mathfrak{F}(x)d\mathcal{H}^1(x)\nonumber\\
 =\,&\int_{\mathcal{S}}d\lambda(T)\int_{T}\mathbbm{1}_A(x)\bar{f}_2(x)\mathfrak{F}(x)d\mathcal{H}^1(x)
 = \int_{A\cap\mathfrak{T}^*}\bar{f}_2(x)dx = \int_A \bar{f}_2(x)dx. \qedhere
\end{align*}
\end{proof}

Recall from \eqref{eq:mathtthDefn} that $\mathcal{X}_k^{\circ}=\prod_{\ell=1}^d(\mathtt{x}_{k,\ell}-\mathtt{h},\mathtt{x}_{k,\ell}+\mathtt{h})$ and $\mathcal{Y}_{k'}^{\circ}=\prod_{\ell=1}^d(\mathtt{y}_{k',\ell}-\mathtt{h},\mathtt{y}_{k',\ell}+\mathtt{h})$, where $\mathtt{h}$ is of constant order. We now partition $\mathcal{X}_k^{\circ}$ into the union of smaller, disjoint $d$-dimensional cubes $\{L_i\}_{i=1}^I$ (each cube is formed by the product of $d$ open or half-open intervals), requiring that the number of such cubes is of order $(\delta_2 \sqrt{\eps})^{-d}$ and the side lengths of each cube are constrained to be within $[\delta_2 \sqrt{\eps}\slash 2, \delta_2\sqrt{\eps}]$. Thus $\mathcal{X}_k^{\circ}=\bigcup_{i=1}^{I} L_i$. We do the same for $\mathcal{Y}_{k'}^{\circ}$, so that $\mathcal{Y}_{k'}^{\circ}= \bigcup_{j=1}^{J} R_j$, where $\{R_j\}_{j=1}^J$ are the corresponding (disjoint) cubes.

For every $i\in [I]$ and $j\in [J]$, we denote by $\mathtt{c}_i=(\mathtt{c}_{i,1},\cdots,\mathtt{c}_{i,d})$ the center of $L_i$ and $\mathtt{c}_j'=(\mathtt{c}'_{j,1},\cdots,\mathtt{c}'_{j,d})$ the center of $R_j$. Let
\begin{equation}\label{defU}
    U_{i,j}:=\frac{\mathtt{c}_i-\mathtt{c}_j'}{\|\mathtt{c}_i-\mathtt{c}_j'\|}.
\end{equation}
Note that for any $x\in L_i$ and $y\in R_j$,
\begin{equation*}
    |\langle U_{i,j}, x-\mathtt{c}_i\rangle|\leq \|U_{i,j}\|\|x-\mathtt{c}_i\|=\|x-\mathtt{c}_i\|\leq \frac{1}{2}\delta_2\sqrt{d\eps},
\end{equation*}
\begin{equation*}
    |\langle U_{i,j}, y-\mathtt{c}'_j\rangle|\leq \|U_{i,j}\|\|y-\mathtt{c}'_j\|=\|y-\mathtt{c}'_j\|\leq \frac{1}{2}\delta_2\sqrt{d\eps}.
\end{equation*}
Let $N_0:=\lceil 2\sqrt{d}\eps^{-1\slash 2}\rceil$. For every $\ell\in [N_0]$, we define 
\begin{equation}\label{Lijk}
    L_{i,j;\ell}:=\Big\{x\in L_i: \langle U_{i,j},x-\mathtt{c}_i\rangle \in \big(-\delta_2\sqrt{d\eps}+(\ell-1)\delta_2 \eps, -\delta_2\sqrt{d\eps}+\ell\delta_2 \eps\big]\Big\}, 
\end{equation}
\begin{equation}\label{Rijk}
    R_{i,j;\ell}:=\Big\{y\in R_j: \langle U_{i,j},y-\mathtt{c}'_j\rangle \in \big(-\delta_2\sqrt{d\eps}+(\ell-1)\delta_2 \eps, -\delta_2\sqrt{d\eps}+\ell\delta_2 \eps\big]\Big\}. 
\end{equation}
Note that $L_i=\bigcup_{\ell=1}^{N_0}L_{i,j;\ell}$ and $R_j=\bigcup_{\ell=1}^{N_0}R_{i,j;\ell}$. Moreover, 
\begin{equation}\label{relation}
    c(\delta_2\sqrt{\eps})^{-d}\leq I,J\leq C(\delta_2\sqrt{\eps})^{-d}, \qquad c\eps^{-1\slash 2}\leq N_0\leq C\eps^{-1\slash2}.
\end{equation}

For any $j\in [J]$ and $x\in\mathbb{R}^d$, we define
\begin{equation}\label{deffx}
    f_j(x):=\mathbbm{1}_{\mathcal{X}_k^{\circ}\cap\mathfrak{T}^{*}}(x)
    \bar{f}_2(x) \mathbbm{1}_{\int_{\mathrm{int}(T(x))\cap \mathcal{Y}_{k'}^{\circ}}\bar{g}_2(z)\mathfrak{F}(z)d\mathcal{H}^1(z)>0}\cdot\frac{\int_{\mathrm{int}(T(x))\cap\mathcal{Y}_{k'}^{\circ}\cap R_j}\bar{g}_2(z)\mathfrak{F}(z)d\mathcal{H}^1(z)}{\int_{\mathrm{int}(T(x))\cap \mathcal{Y}_{k'}^{\circ}}\bar{g}_2(z)\mathfrak{F}(z)d\mathcal{H}^1(z)}.
\end{equation}
For any $i\in [I]$ and $x\in\mathbb{R}^d$, we define
\begin{equation}
    g_i(x):=\mathbbm{1}_{\mathcal{Y}_{k'}^{\circ}\cap\mathfrak{T}^{*}}(x)
    \bar{g}_2(x) \mathbbm{1}_{\int_{\mathrm{int}(T(x))\cap \mathcal{X}_k^{\circ}}\bar{f}_2(z)\mathfrak{F}(z)d\mathcal{H}^1(z)>0}\cdot\frac{\int_{\mathrm{int}(T(x))\cap\mathcal{X}_{k}^{\circ} \cap  L_i}\bar{f}_2(z)\mathfrak{F}(z)d\mathcal{H}^1(z)}{\int_{\mathrm{int}(T(x))\cap \mathcal{X}_k^{\circ}}\bar{f}_2(z)\mathfrak{F}(z)d\mathcal{H}^1(z)}.
\end{equation}
Then, $f_j$ and $g_i$ are convenient versions of the densities of $\bar{\gamma}_{\eps,2}(\cdot\times R_j)$ and $\bar{\gamma}_{\eps,2}(L_i\times \cdot)$, respectively.

\begin{lemma}\label{Lem4.9n}
For any $i\in [I],j\in [J]$, $f_j$ is a density of $\bar{\gamma}_{\eps,2}(\cdot\times R_j)$ with respect to $\mathcal{L}^d$, and $g_i$ is a density of $\bar{\gamma}_{\eps,2}(L_i\times \cdot)$ with respect to $\mathcal{L}^d$. 
\end{lemma}
\begin{proof}

Fix any Borel set $A\subseteq \mathbb{R}^d$. By \eqref{def.gamma.2}, \eqref{subinclud1.1}, and Lemma \ref{Lem3.18nn}, for any $j\in [J]$, we have 
\begin{align*}
   & \bar{\gamma}_{\eps,2}(A\times R_j)= \int_{\mathcal{S}}\mathbbm{1}_{\int_{\mathrm{int}(T)\cap \mathcal{X}_k^{\circ}}\bar{f}_2(z)\mathfrak{F}(z)d\mathcal{H}^1(z)>0}\nonumber\\
&\hspace{1.2in}\times\frac{\int_{\mathrm{int}(T)\cap\mathcal{X}_k^{\circ}}\mathbbm{1}_A(x)\bar{f}_2(x)\mathfrak{F}(x)d\mathcal{H}^1(x)
 \int_{\mathrm{int}(T)\cap\mathcal{Y}_{k'}^{\circ}\cap R_j}\bar{g}_2(y)\mathfrak{F}(y)d\mathcal{H}^1(y)}{\int_{\mathrm{int}(T)\cap \mathcal{X}_k^{\circ}}\bar{f}_2(z)\mathfrak{F}(z)d\mathcal{H}^1(z)}
 d\lambda(T)\nonumber\\
 =\,&\int_{\mathcal{S}}d\lambda(T)\int_{\mathrm{int}(T)}\mathbbm{1}_{A\cap\mathcal{X}_k^{\circ}}(x)\bar{f}_2(x)\mathfrak{F}(x)\mathbbm{1}_{\int_{\mathrm{int}(T(x))\cap \mathcal{X}_k^{\circ}}\bar{f}_2(z)\mathfrak{F}(z)d\mathcal{H}^1(z)>0}\nonumber\\
 &\hspace{0.95in}\times\frac{
 \int_{\mathrm{int}(T(x))\cap\mathcal{Y}_{k'}^{\circ}\cap R_j}\bar{g}_2(z)\mathfrak{F}(z)d\mathcal{H}^1(z)}{\int_{\mathrm{int}(T(x))\cap \mathcal{X}_k^{\circ}}\bar{f}_2(z)\mathfrak{F}(z)d\mathcal{H}^1(z)}d\mathcal{H}^1(x)\nonumber\\
 =\,&\int_{\mathcal{S}}d\lambda(T)\int_{\mathrm{int}(T)}\mathbbm{1}_{A\cap\mathcal{X}_k^{\circ}}(x)\bar{f}_2(x)\mathfrak{F}(x)\mathbbm{1}_{\int_{\mathrm{int}(T(x))\cap \mathcal{Y}_{k'}^{\circ}}\bar{g}_2(z)\mathfrak{F}(z)d\mathcal{H}^1(z)>0}\nonumber\\
 &\hspace{0.95in}\times\frac{
 \int_{\mathrm{int}(T(x))\cap\mathcal{Y}_{k'}^{\circ}\cap R_j}\bar{g}_2(z)\mathfrak{F}(z)d\mathcal{H}^1(z)}{\int_{\mathrm{int}(T(x))  \cap \mathcal{Y}_{k'}^{\circ}}\bar{g}_2(z)\mathfrak{F}(z)d\mathcal{H}^1(z)}d\mathcal{H}^1(x)\nonumber\\
 =\,& \int_{\mathfrak{T}^*}\mathbbm{1}_{A\cap\mathcal{X}_k^{\circ}}(x)\bar{f}_2(x)\mathbbm{1}_{\int_{\mathrm{int}(T(x))\cap \mathcal{Y}_{k'}^{\circ}}\bar{g}_2(z)\mathfrak{F}(z)d\mathcal{H}^1(z)>0}\cdot \frac{
 \int_{\mathrm{int}(T(x))\cap\mathcal{Y}_{k'}^{\circ}\cap R_j}\bar{g}_2(z)\mathfrak{F}(z)d\mathcal{H}^1(z)}{\int_{\mathrm{int}(T(x))  \cap \mathcal{Y}_{k'}^{\circ}}\bar{g}_2(z)\mathfrak{F}(z)d\mathcal{H}^1(z)}dx\nonumber\\
 =\,& \int_A f_j(x)dx.
\end{align*}
Similarly, for any $i\in[I]$, we have 
\begin{align*}
   &\bar{\gamma}_{\eps,2}(L_i\times A) = \int_{\mathcal{S}}\mathbbm{1}_{\int_{\mathrm{int}(T)\cap \mathcal{X}_k^{\circ}}\bar{f}_2(z)\mathfrak{F}(z)d\mathcal{H}^1(z)>0}\nonumber\\
&\hspace{1.2in}\times\frac{\int_{\mathrm{int}(T)\cap\mathcal{X}_k^{\circ}\cap L_i}\bar{f}_2(x)\mathfrak{F}(x)d\mathcal{H}^1(x)
 \int_{\mathrm{int}(T)\cap\mathcal{Y}_{k'}^{\circ}}\mathbbm{1}_A(y)\bar{g}_2(y)\mathfrak{F}(y)d\mathcal{H}^1(y)}{\int_{\mathrm{int}(T)\cap \mathcal{X}_k^{\circ}}\bar{f}_2(z)\mathfrak{F}(z)d\mathcal{H}^1(z)}
 d\lambda(T)\nonumber\\
 =\,& \int_{\mathcal{S}}d\lambda(T)\int_{\mathrm{int}(T)}\mathbbm{1}_{A\cap\mathcal{Y}_{k'}^{\circ}}(y)\bar{g}_2(y)\mathfrak{F}(y)\mathbbm{1}_{\int_{\mathrm{int}(T(y))\cap \mathcal{X}_k^{\circ}}\bar{f}_2(z)\mathfrak{F}(z)d\mathcal{H}^1(z)>0}\nonumber\\
 &\hspace{0.95in}\times\frac{
 \int_{\mathrm{int}(T(y))\cap\mathcal{X}_k^{\circ}\cap L_i}\bar{f}_2(z)\mathfrak{F}(z)d\mathcal{H}^1(z)}{\int_{\mathrm{int}(T(y))  \cap \mathcal{X}_k^{\circ}}\bar{f}_2(z)\mathfrak{F}(z)d\mathcal{H}^1(z)}d\mathcal{H}^1(y)\nonumber\\
 =\,&\int_{\mathfrak{T}^*}\mathbbm{1}_{A\cap\mathcal{Y}_{k'}^{\circ}}(y)\bar{g}_2(y)\mathbbm{1}_{\int_{\mathrm{int}(T(y))\cap \mathcal{X}_k^{\circ}}\bar{f}_2(z)\mathfrak{F}(z)d\mathcal{H}^1(z)>0}\cdot \frac{
 \int_{\mathrm{int}(T(y))\cap\mathcal{X}_k^{\circ}\cap L_i}\bar{f}_2(z)\mathfrak{F}(z)d\mathcal{H}^1(z)}{\int_{\mathrm{int}(T(y))  \cap \mathcal{X}_k^{\circ}}\bar{f}_2(z)\mathfrak{F}(z)d\mathcal{H}^1(z)}dy\nonumber\\
 =\,&\int_A g_i(y)dy. \qedhere
\end{align*}
\end{proof}

For any $x,y\in\mathbb{R}^d$, we define
\begin{equation}\label{def:r}
    r(x,y):=\sum_{i=1}^I\sum_{j=1}^J
   \sum_{\substack{\ell,\ell'\in [N_0]:\\\bar{\gamma}_{\eps,2}(L_{i,j;\ell}\times R_{i,j;\ell'})>0}}
    \frac{\bar{\gamma}_{\eps,2}(L_{i,j;\ell}\times R_{i,j;\ell'})}{\bar{\gamma}_{\eps,2}(L_{i,j;\ell}\times R_j) \bar{\gamma}_{\eps,2}(L_i\times R_{i,j;\ell'})}\cdot f_{j}(x)g_{i}(y)\mathbbm{1}_{L_{i,j;\ell}}(x)\mathbbm{1}_{R_{i,j;\ell'}}(y). 
\end{equation}

\subsubsection{Definition of $\tilde{\gamma}_{\eps;k,k'}$}\label{Sect.4.2.3}

With the functions $p(\cdot,\cdot)$ and $r(\cdot,\cdot)$ in hand (see \eqref{def:p} and \eqref{def:r}), we are ready to define $\tilde{\gamma}_{\eps;k,k'}$.

\begin{definition}\label{def4.4}
For any $x,y\in\mathbb{R}^d$, we define
\begin{equation*}
    \tilde{\phi}_{\eps;k,k'}(x,y):=(1-\delta_3)p(x,y)+r(x,y).
\end{equation*}
We define $\tilde{\gamma}_{\eps;k,k'}$ to be the Borel measure on $\mathbb{R}^d\times\mathbb{R}^d$ with density $\tilde{\phi}_{\eps;k,k'}(\cdot,\cdot)$ with respect to $\mathcal{L}^d\otimes\mathcal{L}^d$.
\end{definition}

The following lemma shows that $\tilde{\gamma}_{\eps;k,k'}$ has marginals $\mu_{k,k'}$ and $\nu_{k,k'}$, and is concentrated on the set $\mathcal{X}_k^{\circ}\times\mathcal{Y}_{k'}^{\circ}$.   

\begin{lemma}\label{lemmacoup}
We have $\tilde{\gamma}_{\eps;k,k'}\in\Pi(\mu_{k,k'},\nu_{k,k'})$. Moreover, for any $(x,y)\in (\mathbb{R}^d\times \mathbb{R}^d)\backslash (\mathcal{X}_k^{\circ}\times\mathcal{Y}_{k'}^{\circ})$, we have   $\tilde{\phi}_{\eps;k,k'}(x,y)=0$. 
\end{lemma}
\begin{proof}

By Lemma \ref{Lem4.9} and \eqref{subinclud}, we have 
\begin{equation}\label{Lem4.9ne1}
     \bar{\gamma}_{\eps,2}((\mathbb{R}^d\backslash\mathcal{X}_k^{\circ})\times\mathbb{R}^d)=\int_{\mathbb{R}^d\backslash\mathcal{X}_k^{\circ}}\bar{f}_2(x)dx=0,\qquad \bar{\gamma}_{\eps,2}(\mathbb{R}^d\times(\mathbb{R}^d\backslash\mathcal{Y}_{k'}^{\circ}))=\int_{\mathbb{R}^d\backslash\mathcal{Y}_{k'}^{\circ}}\bar{g}_2(x)dx=0.
\end{equation}
By Lemma \ref{Lem4.9n}, for any Borel set $A\subseteq \mathbb{R}^d$, we have 
\begin{align*}
   & \int_{A\times\mathbb{R}^d} r(x,y)dxdy 
    = \sum_{i=1}^I\sum_{j=1}^J
    \sum_{\substack{\ell,\ell'\in [N_0]:\\\bar{\gamma}_{\eps,2}(L_{i,j;\ell}\times R_{i,j;\ell'})>0}}
    \frac{\bar{\gamma}_{\eps,2}(L_{i,j;\ell}\times R_{i,j;\ell'})}{\bar{\gamma}_{\eps,2}(L_{i,j;\ell}\times R_j) \bar{\gamma}_{\eps,2}(L_i\times R_{i,j;\ell'})}\nonumber\\
    &\hspace{3in}\times \int_{A\times\mathbb{R}^d}f_j(x)g_i(y)\mathbbm{1}_{L_{i,j;\ell}}(x)\mathbbm{1}_{R_{i,j;\ell'}}(y)dxdy\nonumber\\
    =\,& \sum_{i=1}^I\sum_{j=1}^J
    \sum_{\substack{\ell,\ell'\in [N_0]:\\\bar{\gamma}_{\eps,2}(L_{i,j;\ell}\times R_{i,j;\ell'})>0}}
    \frac{\bar{\gamma}_{\eps,2}(L_{i,j;\ell}\times R_{i,j;\ell'})\bar{\gamma}_{\eps,2}((A\cap L_{i,j;\ell})\times R_j)}{\bar{\gamma}_{\eps,2}(L_{i,j;\ell}\times R_j)}\nonumber\\
    =\,& \sum_{i=1}^I\sum_{j=1}^J
    \sum_{\substack{\ell\in [N_0]:\\\bar{\gamma}_{\eps,2}(L_{i,j;\ell}\times R_j)>0}}\frac{\sum_{\ell'=1}^{N_0}\bar{\gamma}_{\eps,2}(L_{i,j;\ell}\times R_{i,j;\ell'})\bar{\gamma}_{\eps,2}((A\cap L_{i,j;\ell})\times R_j)}{\bar{\gamma}_{\eps,2}(L_{i,j;\ell}\times R_j)}\nonumber\\
    =\,& \sum_{i=1}^I\sum_{j=1}^J
    \sum_{\substack{\ell\in[N_0]:\\\bar{\gamma}_{\eps,2}(L_{i,j;\ell}\times R_j)>0}}\bar{\gamma}_{\eps,2}((A\cap L_{i,j;\ell})\times R_j)=\sum_{i=1}^I\sum_{j=1}^J
\sum_{\ell=1}^{N_0}\bar{\gamma}_{\eps,2}((A\cap L_{i,j;\ell})\times R_j)\nonumber\\
=\,&\sum_{i=1}^I\sum_{j=1}^J
    \bar{\gamma}_{\eps,2}((A\cap L_i)\times R_j)=\bar{\gamma}_{\eps,2}((A\cap\mathcal{X}_k^{\circ})\times\mathcal{Y}_{k'}^{\circ})=\bar{\gamma}_{\eps,2}(A\times\mathbb{R}^d),
\end{align*}
where we use \eqref{Lem4.9ne1} in the last equality. Similarly, we can deduce that $\int_{\mathbb{R}^d\times A} r(x,y)dxdy=\bar{\gamma}_{\eps,2}(\mathbb{R}^d\times A)$. Hence by Lemma \ref{Lem4.9}, the marginal densities of $r(\cdot,\cdot)$ are given by $\bar{f}_2$ and $\bar{g}_2$. Therefore, on noting \eqref{Eq.4.1new}, \eqref{Eq.4.2new}, and Lemma \ref{Lem4.11}, we conclude that $\tilde{\gamma}_{\eps;k,k'}\in\Pi(\mu_{k,k'},\nu_{k,k'})$. By \eqref{domains} and \eqref{def:r}, we have $\tilde{\phi}_{\eps;k,k'}(x,y)=0$ for any $(x,y)\in (\mathbb{R}^d\times \mathbb{R}^d)\backslash (\mathcal{X}_k^{\circ}\times\mathcal{Y}_{k'}^{\circ})$. 
\end{proof}

\subsection{Reduction to a local estimate}\label{Sec4.0}

In this subsection, we reduce Theorem~\ref{th:upperBound} to a local estimate, stated as Proposition~\ref{P4.2} below, for each $k\in[K], k'\in[K']$. We begin with the following definition.

\begin{definition}
For each $k\in[K],k'\in[K']$, we define
\begin{align*}
    & \mathtt{S}_{k,k'}(\mu,\nu;\gamma_0)\nonumber\\
    :=\,& \int_{\mathcal{X}_k^{\circ}\cap\mathfrak{T}^*\cap \tilde{\mathcal{T}}_{1;k,k'}^{*}}dx\int_{\mathrm{int}(T(x))\cap\mathcal{Y}_{k'}^{\circ}}\mathfrak{F}(x)^{-1}h(x,z)\bigg(\log\big(\mathfrak{F}(x)^{-1}h(x,z)\big)-\frac{d-1}{2}\log(2\pi\|x-z\|)\nonumber\\
&\hspace{2.9in}+\frac{1}{2}\log\det\big(\mathbf{I}_d+\|x-z\|F(z)\big)\bigg)d\mathcal{H}^1(z). 
\end{align*}
\end{definition}

With $\tilde{\gamma}_{\eps;k,k'}$ and $\tilde{\phi}_{\eps;k,k'}$ as in Definition~\ref{def4.4}, we can now state the local estimate.

\begin{proposition}[Local cost estimate]\label{P4.2}
For any $k\in[K],k'\in[K']$, $\mathtt{S}_{k,k'}(\mu,\nu;\gamma_0)$ is well-defined and finite. Moreover, for any $k\in[K],k'\in[K']$ and $\delta_2\in (0,1\slash 100)$, we have
\begin{multline*}
   \limsup_{\delta_3\rightarrow 0^{+}} \limsup_{\delta_1\rightarrow 0^{+}}\limsup_{M_1\rightarrow\infty}\limsup_{\eps\rightarrow 0^+}\bigg\{\eps^{-1}\int_{\mathcal{X}_k^{\circ}\times\mathcal{Y}_{k'}^{\circ}}\|x-y\|d\tilde{\gamma}_{\eps;k,k'}(x,y)-\eps^{-1}\OT_{k,k'}(\mu,\nu) \\
   +\int_{\mathcal{X}_k^{\circ}\times\mathcal{Y}_{k'}^{\circ}}\tilde{\phi}_{\eps;k,k'}(x,y)\log\big(\eps^{(d-1)\slash 2}\tilde{\phi}_{\eps;k,k'}(x,y)\big)dxdy\bigg\}\leq \mathtt{S}_{k,k'}(\mu,\nu;\gamma_0).
\end{multline*}
\end{proposition}

The proof of Proposition~\ref{P4.2} will be given in Sections~\ref{Sect.4.3}--\ref{Sect.4.5}. In the remainder of this subsection, we show how Theorem~\ref{th:upperBound} follows from Proposition~\ref{P4.2}. We begin with the following integrability lemma.

\begin{lemma}\label{L4.11}
We have 
\begin{equation*}
    \int_{\mathcal{S}}d\lambda(T)\int_{T\times T}\bigg(\frac{d-1}{2}\big|\log(2\pi\|x-z\|)\big|h(x,z)+\big|h(x,z)\log(h(x,z))\big|\bigg)d\mathcal{H}^1(x)d\mathcal{H}^1(z)\leq C.
\end{equation*}
\end{lemma}
\begin{proof}
By Lemma \ref{Lem3.11},
\begin{equation*}
    \int_{\mathcal{S}}d\lambda(T)\int_{T\times T}h(x,z)\big|\log(\mathfrak{F}(x))\big|d\mathcal{H}^1(x)d\mathcal{H}^1(z)\leq \int_{\mathcal{S}}d\lambda(T)\int_T\tilde{f}(x)\big|\log(\mathfrak{F}(x))\big|d\mathcal{H}^1(x)\leq C.
\end{equation*}
By \eqref{generx} and \eqref{genery}, noting that $\lambda(\mathcal{S}\backslash\tilde{\mathcal{S}}^*)=0$ (recall Definition \ref{defstilde}) and $\mathrm{int}(T)\subseteq \tilde{\mathcal{T}}_{1;k,k'}^*$ for any $T\in\mathcal{S}$ such that $\mathrm{int}(T)\cap\mathcal{X}_k^{\circ}\neq\emptyset,\mathrm{int}(T)\cap\mathcal{Y}_{k'}^{\circ}\neq\emptyset$, we obtain that
\begin{align*}
&\int_{\mathcal{S}}d\lambda(T)\int_{T\times T}h(x,z)\bigg(\frac{d-1}{2}\big|\log(2\pi\|x-z\|)\big|+\big|\log\big(\mathfrak{F}(x)^{-1} h(x,z)\big)\big|\bigg)d\mathcal{H}^1(x)d\mathcal{H}^1(z)\nonumber\\
=\,& \int_{\tilde{\mathcal{S}}^*}d\lambda(T)\int_{T\times T}h(x,z)\bigg(\frac{d-1}{2}\big|\log(2\pi\|x-z\|)\big|+\big|\log\big(\mathfrak{F}(x)^{-1} h(x,z)\big)\big|\bigg)d\mathcal{H}^1(x)d\mathcal{H}^1(z)\nonumber\\
\leq\,& \sum_{k=1}^K\sum_{k'=1}^{K'} \int_{\tilde{\mathcal{S}}^*}\mathbbm{1}_{\mathrm{int}(T)\subseteq\tilde{\mathcal{T}}_{1;k,k'}^*}d\lambda(T)\nonumber\\
&\hspace{0.1in}\int_{(\mathrm{int}(T)\cap\mathcal{X}_k^{\circ})\times (\mathrm{int}(T)\cap\mathcal{Y}_{k'}^{\circ})}h(x,z)\bigg(\frac{d-1}{2}\big|\log(2\pi\|x-z\|)\big|+\big|\log\big(\mathfrak{F}(x)^{-1} h(x,z)\big)\big|\bigg)d\mathcal{H}^1(x)d\mathcal{H}^1(z)\nonumber\\
\leq\,& \sum_{k=1}^K\sum_{k'=1}^{K'}\int_{\mathcal{X}_k^{\circ}\cap\mathfrak{T}^*\cap \tilde{\mathcal{T}}_{1;k,k'}^{*}}dx\nonumber\\
&\hspace{0.5in}\int_{\mathrm{int}(T(x))\cap\mathcal{Y}_{k'}^{\circ}}\mathfrak{F}(x)^{-1}h(x,z)\bigg(\frac{d-1}{2}\big|\log(2\pi\|x-z\|)\big|+\big|\log\big(\mathfrak{F}(x)^{-1} h(x,z)\big)\big|\bigg)d\mathcal{H}^1(z)\nonumber\\
\leq\,& C,
\end{align*}
where we use Lemma \ref{Lem3.18nn} and the estimate~\eqref{eq3L4.15} (derived independently in Section~\ref{Sect.4.3} below) in the second and third inequalities, respectively. The conclusion of the lemma follows from the above two displays. 
\end{proof}

The next lemma expresses $\mathtt{S}(\mu,\nu;\gamma_0)$ in terms of $\mathtt{S}_{k,k'}(\mu,\nu;\gamma_0)$.

\begin{lemma}\label{L4.12}
With $\mathtt{S}(\mu,\nu;\gamma_0)$ as in Definition~\ref{de:Sfunctional}, we have
\begin{equation*}
\mathtt{S}(\mu,\nu;\gamma_0)=\sum_{k=1}^K\sum_{k'=1}^{K'}\mathtt{S}_{k,k'}(\mu,\nu;\gamma_0)-\int_{\mathcal{X}}f(x)\log(f(x))dx-\int_{\mathcal{Y}}g(y)\log(g(y))dy.
\end{equation*}
\end{lemma}
\begin{proof}

By Lemma \ref{Lem3.18nn}, for any $k\in[K],k'\in [K']$, we have
\begin{align*}
    & \mathtt{S}_{k,k'}(\mu,\nu;\gamma_0)\nonumber\\
    =\,&\int_{\mathcal{S}}d\lambda(T)\int_{(\mathrm{int}(T)\cap\mathcal{X}_k^{\circ})\times (\mathrm{int}(T) \cap\mathcal{Y}_{k'}^{\circ})}\mathbbm{1}_{\tilde{\mathcal{T}}_{1;k,k'}^*}(x) h(x,z)\nonumber\\
   &\hspace{0.08in}\times\bigg(\log\big(\mathfrak{F}(x)^{-1}h(x,z)\big)-\frac{d-1}{2}\log(2\pi\|x-z\|)+\frac{1}{2}\log\det\big(\mathbf{I}_d+\|x-z\|F(z)\big)\bigg)d\mathcal{H}^1(x)d\mathcal{H}^1(z)\nonumber\\
   =\,&\int_{\tilde{\mathcal{S}}^*}d\lambda(T)\int_{(\mathrm{int}(T)\cap\mathcal{X}_k^{\circ})\times (\mathrm{int}(T) \cap\mathcal{Y}_{k'}^{\circ})}h(x,z)\nonumber\\
   &\hspace{0.08in}\times\bigg(\log\big(\mathfrak{F}(x)^{-1}h(x,z)\big)-\frac{d-1}{2}\log(2\pi\|x-z\|)+\frac{1}{2}\log\det\big(\mathbf{I}_d+\|x-z\|F(z)\big)\bigg)d\mathcal{H}^1(x)d\mathcal{H}^1(z),
\end{align*}
where in the second equality, we use the facts that $\lambda(\mathcal{S}\backslash\tilde{\mathcal{S}}^*)=0$ (recall Definition \ref{defstilde}) and that $\mathrm{int}(T)\subseteq \tilde{\mathcal{T}}_{1;k,k'}^*$ for any $T\in\mathcal{S}$ such that $\mathrm{int}(T)\cap\mathcal{X}_k^{\circ}\neq\emptyset,\mathrm{int}(T)\cap\mathcal{Y}_{k'}^{\circ}\neq\emptyset$. Hence by \eqref{generx}--\eqref{genery} and the facts that $\lambda(\mathcal{S}\backslash\tilde{\mathcal{S}}^*)=0$ and $\tilde{\mathcal{S}}^*\subseteq\mathcal{S}^*\subseteq\mathcal{S}$, we have   
\begin{align*}
    &\sum_{k=1}^K\sum_{k'=1}^{K'}\mathtt{S}_{k,k'}(\mu,\nu;\gamma_0)=\int_{\tilde{\mathcal{S}}^*}d\lambda(T)\int_{T\times T}h(x,z)\bigg(\log\big(\mathfrak{F}(x)^{-1}h(x,z)\big)-\frac{d-1}{2}\log(2\pi\|x-z\|)\nonumber\\
    & \hspace{3.22in}+\frac{1}{2}\log\det\big(\mathbf{I}_d+\|x-z\|F(z)\big)\bigg)d\mathcal{H}^1(x)d\mathcal{H}^1(z)\nonumber\\
    =\,& \int_{\mathcal{S}}d\lambda(T)\int_{T\times T}\bigg(-\frac{d-1}{2}\log(2\pi\|x-z\|)h(x,z)+h(x,z)\log(h(x,z))\bigg)d\mathcal{H}^1(x)d\mathcal{H}^1(z)\nonumber\\
    &-\int_{\mathcal{S}^*}d\lambda(T)\int_{T}\tilde{f}(x)\log(\mathfrak{F}(x))d\mathcal{H}^1(x)\nonumber\\
    &+\frac{1}{2}\int_{\mathcal{S}}d\lambda(T)\int_{T\times T}\log\det\big(\mathbf{I}_d+\|x-z\|F(z)\big)h(x,z)d\mathcal{H}^1(x)d\mathcal{H}^1(z)\nonumber\\
    =\,& \int_{\mathcal{S}}d\lambda(T)\int_{T\times T}\bigg(-\frac{d-1}{2}\log(2\pi\|x-y\|)h(x,y)+h(x,y)\log(h(x,y))\bigg)d\mathcal{H}^1(x)d\mathcal{H}^1(y)\nonumber\\
    &+\frac{1}{2}\int_{\mathcal{X}}f(x)\log(f(x))dx+\frac{1}{2}\int_{\mathcal{Y}}g(y)\log(g(y))dy\nonumber\\
    &-\frac{1}{2}\int_{\mathcal{S}}d\lambda(T)\int_T\tilde{f}(x)\log(\tilde{f}(x))d\mathcal{H}^1(x)-\frac{1}{2}\int_{\mathcal{S}}d\lambda(T)\int_T\tilde{g}(y)\log(\tilde{g}(y))d\mathcal{H}^1(y)\nonumber\\
    =\,& \mathtt{S}(\mu,\nu;\gamma_0)+\int_{\mathcal{X}}f(x)\log(f(x))dx+\int_{\mathcal{Y}}g(y)\log(g(y))dy,
\end{align*}
where the second to last equality uses \Cref{Lem3.11}.
\end{proof}

We now give the proof of Theorem~\ref{th:upperBound}, assuming Proposition~\ref{P4.2}.

\begin{proof}[Proof of Theorem~\ref{th:upperBound}]

We define
\begin{equation}\label{def_gamma_eps}
    \tilde{\gamma}_{\eps}:=\sum_{k=1}^K\sum_{k'=1}^{K'}\tilde{\gamma}_{\eps;k,k'},
\end{equation}
where $\tilde{\gamma}_{\eps;k,k'}$ is as in Definition~\ref{def4.4}. By Lemmas~\ref{Ln4.1} and~\ref{lemmacoup}, $\tilde{\gamma}_{\eps}\in\Pi(\mu,\nu)$. Hence 
\begin{align*}
    \mathcal{C}_{\eps}(\gamma_{\eps})\leq\,& \mathcal{C}_{\eps}(\tilde{\gamma}_{\eps})=\int_{\mathbb{R}^d\times \mathbb{R}^d} \|x-y\| d\tilde{\gamma}_{\eps}(x,y)+\eps H(\tilde{\gamma}_{\eps}|\mu\otimes \nu)\nonumber\\
   =\,& \sum_{k=1}^K\sum_{k'=1}^{K'}\bigg(\int_{\mathcal{X}_k^{\circ}\times\mathcal{Y}_{k'}^{\circ}}\|x-y\|d\tilde{\gamma}_{\eps;k,k'}(x,y)+\eps\int_{\mathcal{X}_k^{\circ}\times\mathcal{Y}_{k'}^{\circ}}\tilde{\phi}_{\eps;k,k'}(x,y)\log\big(\tilde{\phi}_{\eps;k,k'}(x,y)\big)dxdy\bigg)\nonumber\\
   &-\eps\int_{\mathcal{X}}f(x)\log(f(x))dx-\eps\int_{\mathcal{Y}}g(y)\log(g(y))dy,
\end{align*}
where we use Lemma~\ref{lemmacoup} in the second line. Moreover,
\begin{align*}
& \sum_{k=1}^K\sum_{k'=1}^{K'}\int_{\mathcal{X}_k^{\circ}\times\mathcal{Y}_{k'}^{\circ}}\tilde{\phi}_{\eps;k,k'}(x,y)\log\big(\eps^{(d-1)\slash 2}\big)dxdy=\frac{d-1}{2}\log\eps\sum_{k=1}^K\sum_{k'=1}^{K'}\tilde{\gamma}_{\eps;k,k'}(\mathbb{R}^d\times\mathbb{R}^d)\nonumber\\
=\,& \frac{d-1}{2}\log\eps\sum_{k=1}^K\sum_{k'=1}^{K'}\mu_{k,k'}(\mathbb{R}^d)=\frac{d-1}{2}\log\eps \cdot\mu(\mathbb{R}^d)=\frac{d-1}{2}\log\eps.
\end{align*}
Combining the above two displays and noting Lemma~\ref{Ln4.1}, we get
\begin{align*}
    & \frac{\mathcal{C}_{\eps}(\gamma_{\eps})-\OT(\mu,\nu)}{\eps}+\frac{d-1}{2}\log\eps\nonumber\\
    \leq\,& \sum_{k=1}^K\sum_{k'=1}^{K'}\bigg(\eps^{-1}\int_{\mathcal{X}_k^{\circ}\times\mathcal{Y}_{k'}^{\circ}}\|x-y\|d\tilde{\gamma}_{\eps;k,k'}(x,y)-\eps^{-1}\OT_{k,k'}(\mu,\nu)\nonumber\\
    &\hspace{0.6in}+\int_{\mathcal{X}_k^{\circ}\times\mathcal{Y}_{k'}^{\circ}}\tilde{\phi}_{\eps;k,k'}(x,y)\log\big(\eps^{(d-1)\slash 2}\tilde{\phi}_{\eps;k,k'}(x,y)\big)dxdy\bigg)\nonumber\\
   &-\int_{\mathcal{X}}f(x)\log(f(x))dx-\int_{\mathcal{Y}}g(y)\log(g(y))dy.
\end{align*}
Therefore, by Proposition~\ref{P4.2} and Lemma~\ref{L4.12}, we have  
\begin{align*}
  & \limsup_{\eps\rightarrow 0^+}\bigg\{\frac{\mathcal{C}_{\eps}(\gamma_{\eps})-\OT(\mu,\nu)}{\eps}+\frac{d-1}{2}\log{\eps}\bigg\}\nonumber\\
  \leq\,& \sum_{k=1}^K\sum_{k'=1}^{K'}\mathtt{S}_{k,k'}(\mu,\nu;\gamma_0)-\int_{\mathcal{X}}f(x)\log(f(x))dx-\int_{\mathcal{Y}}g(y)\log(g(y))dy= \mathtt{S}(\mu,\nu;\gamma_0).\qedhere
\end{align*}
\end{proof}  

\subsection{Mass and entropy estimates, and proof of the local estimate}\label{Sect.4.3}

Throughout the rest of this section, we fix $k\in[K],k'\in[K']$. Recall the definitions of $p(\cdot,\cdot)$ and $r(\cdot,\cdot)$ from \eqref{def:p} and \eqref{def:r}, respectively. We now state two key propositions that play a central role in the proof of the desired local estimate, \cref{P4.2}. The first one shows that the mass of the sub-coupling with density $p(\cdot,\cdot)$ is close to $\mu_{k,k'}(\mathbb{R}^d)$ (cf.\ Step~3 in Section~\ref{se:upper}).

\begin{proposition}[Mass estimate]\label{P4.3}
For any fixed $\delta_2,\delta_3\in (0,1\slash 100)$, we have 
\begin{equation*}
    \liminf_{\delta_1\rightarrow 0^{+}}\liminf_{M_1\rightarrow\infty}\liminf_{\eps\rightarrow 0^{+}}\int_{\mathbb{R}^d\times\mathbb{R}^d}p(x,y)dxdy=\mu_{k,k'}(\mathbb{R}^d).
\end{equation*}
\end{proposition}

The second proposition provides tight control on the entropy of $r(\cdot,\cdot)$ (cf.\ Step~5 in Section~\ref{se:upper}).

\begin{proposition}[Entropy estimate]\label{P4.4}
For any fixed $\delta_2\in (0,1\slash 100)$, we have 
\begin{equation*}
   \limsup_{\delta_3\rightarrow 0^{+}} \limsup_{\delta_1\rightarrow 0^{+}}\limsup_{M_1\rightarrow\infty}\limsup_{\eps\rightarrow 0^{+}}\int_{\mathbb{R}^d\times\mathbb{R}^d}r(x,y)\max\big\{\log\big(\eps^{(d-1)\slash 2}r(x,y)\big),0\big\}dxdy \leq 0.
\end{equation*}
\end{proposition}

The proofs of Propositions~\ref{P4.3} and~\ref{P4.4} will be presented in Sections~\ref{Sect.4.4} and~\ref{Sect.4.5}, respectively.

In the remainder of this subsection, we detail the proof of Proposition~\ref{P4.2}, assuming Propositions~\ref{P4.3} and~\ref{P4.4}. We begin with the following three lemmas. The first lemma provides an upper bound on the transport cost of $\tilde{\gamma}_{\eps;k,k'}$ (recall Definition~\ref{def4.4}).  

\begin{lemma}\label{lemma4.14}
For any fixed $\delta_2\in (0,1\slash 100)$, we have 
\begin{align*}
   & \limsup_{\delta_3\rightarrow 0^+}\limsup_{\delta_1\rightarrow 0^+}\limsup_{M_1\rightarrow\infty}\limsup_{\eps\rightarrow 0^+}\eps^{-1}\bigg(\int_{\mathcal{X}_{k}^{\circ}\times\mathcal{Y}_{k'}^{\circ}}\|x-y\|d\tilde{\gamma}_{\eps;k,k'}(x,y)-\OT_{k,k'}(\mu,\nu)\bigg)\nonumber\\
   \leq\,&\frac{d-1}{2}\int_{\mathcal{X}_k^{\circ}\cap\mathfrak{T}^*\cap\tilde{\mathcal{T}}_{1;k,k'}^*} dx \int_{\mathrm{int}(T(x))\cap\mathcal{Y}_{k'}^{\circ}}\mathfrak{F}(x)^{-1}h(x,z) d\mathcal{H}^1(z).
\end{align*}
\end{lemma}

\begin{proof}
By \eqref{p3}, for any $(x,y)\in \mathbb{R}^d\times\mathbb{R}^d$, if $V(x)+V(y)= 0$, then $p(x,y)=0$. Hence by \eqref{domains}, 
\begin{align}\label{lem.neq.e12}
    &\int_{\mathcal{X}_{k}^{\circ}\times\mathcal{Y}_{k'}^{\circ}}\|x-y\|d\tilde{\gamma}_{\eps;k,k'}(x,y)=\int_{\mathcal{X}_{k}^{\circ}\times\mathcal{Y}_{k'}^{\circ}}\|x-y\|\tilde{\phi}_{\eps;k,k'}(x,y)dxdy\nonumber\\
    =\,& (1-\delta_3)\int_{(\mathcal{X}_k^{\circ}\cap\mathfrak{T}^*)\times(\mathcal{Y}_{k'}^{\circ}\cap\mathfrak{T}^*)}\mathbbm{1}_{V(x)+V(y)\neq 0}\|x-y\|p(x,y)dxdy+\int_{\mathcal{X}_{k}^{\circ}\times\mathcal{Y}_{k'}^{\circ}}\|x-y\|r(x,y)dxdy\nonumber\\
    \leq\,& (1-\delta_3)\int_{(\mathcal{X}_k^{\circ}\cap\mathfrak{T}^*)\times(\mathcal{Y}_{k'}^{\circ}\cap\mathfrak{T}^*)}\mathbbm{1}_{V(x)+V(y)\neq 0}(\|x-y\|-\|x-\mathtt{y}'(x,y)\|)p(x,y)dxdy\nonumber\\
    &+(1-\delta_3)\int_{(\mathcal{X}_k^{\circ}\cap\mathfrak{T}^*)\times(\mathcal{Y}_{k'}^{\circ}\cap\mathfrak{T}^*)}\mathbbm{1}_{V(x)+V(y)\neq 0}\|x-\mathtt{y}'(x,y)\|p(x,y)dxdy\nonumber\\
    &+\int_{\mathbb{R}^d\times\mathbb{R}^d}\|x-y\|r(x,y)dxdy,
\end{align}
where $\mathtt{y}'(\cdot,\cdot)$ is as in \eqref{y'd}.

\paragraph{Step 1: Bounding $\int_{(\mathcal{X}_k^{\circ}\cap\mathfrak{T}^*)\times(\mathcal{Y}_{k'}^{\circ}\cap\mathfrak{T}^*)}\mathbbm{1}_{V(x)+V(y)\neq 0}(\|x-y\|-\|x-\mathtt{y}'(x,y)\|)p(x,y)dxdy$.} We have
\begin{align}\label{Lem.4.14.e10}
    & \int_{(\mathcal{X}_k^{\circ}\cap\mathfrak{T}^*)\times(\mathcal{Y}_{k'}^{\circ}\cap\mathfrak{T}^*)}\mathbbm{1}_{V(x)+V(y)\neq 0}(\|x-y\|-\|x-\mathtt{y}'(x,y)\|)p(x,y)dxdy\nonumber\\
    =\,& \int_{\mathcal{X}_k^{\circ}\cap\mathfrak{T}^{*} \cap \tilde{\mathcal{T}}_{1;k,k'}^{*}} \mathbbm{1}_{\min\{\alpha(x),\beta(x)\}\geq\delta_1} dx\int_{\mathrm{int}(T(x))\cap\mathcal{Y}_{k'}^{\circ}}\mathbbm{1}_{\min\{\alpha(z),\beta(z)\}\geq\delta_1} \mathbbm{1}_{\langle x-z, V(x)\rangle>0} d\mathcal{H}^1(z)\nonumber\\    &\hspace{0.2in}\int_{O(V(x))}\mathbbm{1}_{z+w\in\mathcal{Y}_{k'}^{\circ}\cap\mathfrak{T}^*\cap \tilde{\mathcal{T}}_{1;k,k'}^{*}}\mathbbm{1}_{V(x)+V(z+w)\neq 0}\mathbbm{1}_{\min\{\alpha(z+w),\beta(z+w)\}\geq\delta_1} \mathbbm{1}_{\|w\|\leq M_1\sqrt{\eps}}\nonumber\\
    &\hspace{0.7in}\times(\|x-(z+w)\|-\|x-\mathtt{y}'(x,z+w)\|)p(x,z+w)d\mathcal{H}^{d-1}(w),
\end{align}
where we note \eqref{p1}--\eqref{p2} and \eqref{def:p}. 

Below, we fix any $x\in\mathcal{X}_k^{\circ}\cap\mathfrak{T}^{*} \cap \tilde{\mathcal{T}}_{1;k,k'}^*$ and $z\in\mathrm{int}(T(x))\cap\mathcal{Y}_{k'}^{\circ}$ such that $\min\{\alpha(x),\beta(x)\}\geq \delta_1$, $\min\{\alpha(z),\beta(z)\}\geq \delta_1$, and $\langle x-z, V(x)\rangle>0$. We also fix any $\delta\in (0,1)$, and assume that $\eps>0$ is sufficiently small (depending on $M_1,\delta_1,\delta$ and $d,\mu,\nu$). Let
\begin{align*}
    \mathscr{W}_0(x,z;\eps):=\,&\Big\{w\in O(V(x)):\|w\|\leq M_1 \sqrt{\eps}, z+w\in \mathcal{Y}_{k'}^{\circ}\cap\mathfrak{T}^*\cap \tilde{\mathcal{T}}_{1;k,k'}^{*}, 
 V(x)+V(z+w) \neq 0,\nonumber\\
    & \hspace{1.1in}  \min\{\alpha(z+w), 
 \beta(z+w)\}\geq\delta_1\Big\},  
\end{align*}
\begin{equation*}
    \mathscr{W}(x,z;\delta,\eps):=\Big\{w\in\mathscr{W}_0(x,z;\eps):w\in \mathcal{W}'(z;M_1\sqrt{\eps},\delta)\Big\}, 
\end{equation*}
where $\mathcal{W}'(z;M_1\sqrt{\eps},\delta)$ is as in Definition \ref{Defd}. Note that $\mathscr{W}_0(x,z;\eps)\backslash\mathscr{W}(x,z;\delta,\eps)\subseteq\mathcal{W}(z;M_1\sqrt{\eps},\delta)$. As $x\in\mathfrak{T}$ and $\mathfrak{T}$ is a transport set (recall \cref{Defn4.1.1n}), we have $z\in\mathfrak{T}$. Hence by \Cref{L2.4}, 
\begin{equation}\label{Lem4.14.e10}
    \lim_{\eps\rightarrow 0^+}\frac{\mathcal{H}^{d-1}\big(\mathcal{W}(z;M_1\sqrt{\eps},\delta)\big)}{\eps^{(d-1)\slash 2}}=0 \,\Rightarrow\, \lim_{\eps\rightarrow 0^+}\frac{\mathcal{H}^{d-1}\big(\mathscr{W}_0(x,z;\eps)\backslash\mathscr{W}(x,z;\delta,\eps)\big)}{\eps^{(d-1)\slash 2}}=0.
\end{equation}
For any $w\in\mathscr{W}_0(x,z;\eps)$, noting that $z-x=-\|z-x\|V(x)$ (as $\langle x-z, V(x)\rangle>0$), we have  
\begin{eqnarray}\label{Lem4.14.e1}
     \mathtt{y}'(x,z+w)&=& x+\frac{\langle z-x+w,V(x)+V(z+w)\rangle}{\langle V(x), V(x)+V(z+w)\rangle}V(x)\nonumber\\
    &=&x-\|z-x\|V(x)+\frac{\langle w,V(x)+V(z+w)\rangle}{\langle V(x), V(x)+V(z+w)\rangle}V(x).
\end{eqnarray}
Moreover,
\begin{equation}\label{Lem4.14.e2}
    \|x-(z+w)\|=\sqrt{\|x-z\|^2+\|w\|^2}\leq \|x-z\|+\frac{\|w\|^2}{2\|x-z\|}. 
\end{equation}
As $\min\{\alpha(z),\beta(z)\}\geq \delta_1$ and $\min\{\alpha(z+w),\beta(z+w)\}\geq\delta_1$, by Lemma \ref{L2.0},
\begin{equation*}
    \|V(z+w)-V(x)\|=\|V(z+w)-V(z)\|\leq 4\|w\|\slash \delta_1 \leq 4M_1\delta_1^{-1}\sqrt{\eps}.
\end{equation*}
Hence noting that $\eps$ is sufficiently small, we have \begin{equation}\label{Lem.q.4.1}
    \langle V(x), V(z+w)\rangle= 1-\frac{\|V(x)-V(z+w)\|^2}{2} \begin{cases}
        \leq 1, \\
        \geq 1-8M_1^2\delta_1^{-2}\eps\geq \frac{1}{2},
    \end{cases}
\end{equation}
\begin{align*}
    |\langle w,V(x)+V(z+w)\rangle|=\,&|\langle w, V(z+w)-V(x)\rangle|\leq \|w\|\|V(z+w)-V(x)\|\nonumber\\
    \leq\,& M_1\sqrt{\eps}\cdot 4M_1\delta_1^{-1}\sqrt{\eps}=4M_1^2\delta_1^{-1}\eps,
\end{align*}
where we note that $\langle w, V(x)\rangle=0$ in the second display. Therefore, for any $w\in\mathscr{W}_0(x,z;\eps)$, we have  
\begin{equation}\label{Lem4.14.e3}
    \bigg|\frac{\langle w,V(x)+V(z+w)\rangle}{\langle V(x), V(x)+V(z+w)\rangle}\bigg|\leq |\langle w, V(x)+V(z+w)\rangle|\leq 4M_1^2\delta_1^{-1}\eps,
\end{equation}
which combined with \eqref{Lem4.14.e1} and the fact that $\|z-x\|\geq 2d_0>4M_1^2\delta_1^{-1}\eps$ (since $x\in\mathcal{X}_k^{\circ}$, $z\in\mathcal{Y}_{k'}^{\circ}$, and $\eps$ is sufficiently small) yields
\begin{equation}\label{Lem4.14.e5}
    \|x-\mathtt{y}'(x,z+w)\|=\|z-x\|-\frac{\langle w,V(x)+V(z+w)\rangle}{\langle V(x), V(x)+V(z+w)\rangle}.
\end{equation}
Moreover, by \Cref{L2.4} (as $x\in\mathfrak{T}$ and $\min\{\alpha(z),\beta(z)\}\geq\delta_1$), we have $\|F(z)\|_2\leq 4\delta_1^{-1}$, hence
\begin{equation}\label{Lem4.14.e6}
     |w^{\top}F(z)w|\leq \|F(z)\|_2\|w\|^2\leq 4M_1^2\delta_1^{-1}\eps.
\end{equation}
Now note that for any $w\in \mathscr{W}(x,z;\delta,\eps)$, we have $\|V(z+w)-V(z)-F(z)w\|\leq\delta M_1\sqrt{\eps}$, hence  
\begin{align}\label{Lem.q.4.2}
   &|\langle w,V(x)+V(z+w)\rangle-w^{\top}F(z)w|=|\langle w,V(z+w)-V(z)-F(z)w\rangle|\nonumber\\
   \leq\,& \|w\|\|V(z+w)-V(z)-F(z)w\|\leq M_1\sqrt{\eps}\cdot \delta M_1\sqrt{\eps} =\delta M_1^2 \eps.
\end{align}
By \eqref{Lem.q.4.1} and \eqref{Lem.q.4.2}, for any $w\in\mathscr{W}(x,z;\delta,\eps)$, we have 
\begin{align}\label{Lem4.14.e4}
    &\bigg|\frac{\langle w,V(x)+V(z+w)\rangle}{\langle V(x), V(x)+V(z+w)\rangle}-\frac{1}{2}w^{\top}F(z)w\bigg|\nonumber\\
    \leq\,& \bigg|\frac{\langle w,V(x)+V(z+w)\rangle}{\langle V(x), V(x)+V(z+w)\rangle}-\frac{1}{2}\langle w,V(x)+V(z+w)\rangle\bigg|+\frac{1}{2}\big|\langle w,V(x)+V(z+w)\rangle-w^{\top}F(z)w\big|\nonumber\\
    \leq\,& \frac{\|w\|\|V(x)+V(z+w)\||\langle V(x), V(x)+V(z+w)\rangle-2|}{2\langle V(x), V(x)+V(z+w)\rangle} +\frac{1}{2}\delta M_1^2\eps\nonumber\\
    \leq\,& M_1\sqrt{\eps}\cdot 8M_1^2\delta_1^{-2}\eps+\frac{1}{2}\delta M_1^2\eps = 8M_1^3\delta_1^{-2}\eps^{3\slash 2}+\frac{1}{2}\delta M_1^2\eps \leq  \delta M_1^2\eps,
\end{align}
where we note that $\eps$ is sufficiently small in the last inequality. By \eqref{Lem4.14.e2},  \eqref{Lem4.14.e3}--\eqref{Lem4.14.e6}, and \eqref{Lem4.14.e4}, 
\begin{align}\label{Eq.4.14.e1}
&\int_{O(V(x))}\mathbbm{1}_{z+w\in\mathcal{Y}_{k'}^{\circ}\cap\mathfrak{T}^*\cap \tilde{\mathcal{T}}_{1;k,k'}^{*}}\mathbbm{1}_{V(x)+V(z+w)\neq 0}\mathbbm{1}_{\min\{\alpha(z+w),\beta(z+w)\}\geq\delta_1} \mathbbm{1}_{\|w\|\leq M_1\sqrt{\eps}}\nonumber\\
    &\hspace{0.52in}\times(\|x-(z+w)\|-\|x-\mathtt{y}'(x,z+w)\|)p(x,z+w)d\mathcal{H}^{d-1}(w)\nonumber\\
    =\,& \int_{\mathscr{W}_0(x,z;\eps)}(\|x-(z+w)\|-\|x-\mathtt{y}'(x,z+w)\|)p(x,z+w)d\mathcal{H}^{d-1}(w)\nonumber\\
    =\,& \int_{\mathscr{W}(x,z;\delta,\eps)}(\|x-(z+w)\|-\|x-\mathtt{y}'(x,z+w)\|)p(x,z+w)d\mathcal{H}^{d-1}(w)\nonumber\\
    & +\int_{\mathscr{W}_0(x,z;\eps)\backslash\mathscr{W}(x,z;\delta,\eps)}(\|x-(z+w)\|-\|x-\mathtt{y}'(x,z+w)\|)p(x,z+w)d\mathcal{H}^{d-1}(w)\nonumber\\
    \leq\,& \int_{\mathscr{W}(x,z;\delta,\eps)}\bigg(\frac{\|w\|^2}{2\|x-z\|}+\frac{1}{2}w^{\top}F(z)w+\delta M_1^2\eps\bigg)p(x,z+w)d\mathcal{H}^{d-1}(w)\nonumber\\
    & + \int_{\mathscr{W}_0(x,z;\eps)\backslash\mathscr{W}(x,z;\delta,\eps)}\bigg(\frac{\|w\|^2}{2\|x-z\|}+\frac{1}{2}w^{\top}F(z)w+6M_1^2\delta_1^{-1}\eps\bigg)p(x,z+w)d\mathcal{H}^{d-1}(w)\nonumber\\
    \leq\,& \eps\int_{\eps^{-1\slash 2}\mathscr{W}(x,z;\delta,\eps)}\bigg(\frac{\|w\|^2}{2\|x-z\|}+\frac{1}{2}w^{\top}F(z)w+\delta M_1^2\bigg)\eps^{(d-1)\slash 2}p(x,z+\sqrt{\eps}w)d\mathcal{H}^{d-1}(w)\nonumber\\
    & + \eps\int_{\eps^{-1\slash 2}(\mathscr{W}_0(x,z;\eps)\backslash\mathscr{W}(x,z;\delta,\eps))}\bigg(\frac{\|w\|^2}{2\|x-z\|}+\frac{1}{2}w^{\top}F(z)w+6M_1^2\delta_1^{-1}\bigg)\eps^{(d-1)\slash 2}p(x,z+\sqrt{\eps}w)d\mathcal{H}^{d-1}(w)\nonumber\\
    =\,& \eps\int_{\eps^{-1\slash 2}\mathscr{W}_0(x,z;\eps)}\bigg(\frac{\|w\|^2}{2\|x-z\|}+\frac{1}{2}w^{\top}F(z)w\bigg)\eps^{(d-1)\slash 2}p(x,z+\sqrt{\eps}w)d\mathcal{H}^{d-1}(w)\nonumber\\
    &+\eps\delta M_1^2\int_{\eps^{-1\slash 2}\mathscr{W}(x,z;\delta,\eps)}\eps^{(d-1)\slash 2}p(x,z+\sqrt{\eps}w)d\mathcal{H}^{d-1}(w)\nonumber\\
    &+6\eps M_1^2\delta_1^{-1}\int_{\eps^{-1\slash 2}(\mathscr{W}_0(x,z;\eps)\backslash\mathscr{W}(x,z;\delta,\eps))}\eps^{(d-1)\slash 2}p(x,z+\sqrt{\eps}w)d\mathcal{H}^{d-1}(w).
\end{align}
By Lemma \ref{L2.7} (as $x\in\mathfrak{T}$), $\mathbf{I}_d+\|x-z\| F(z)$ is positive definite, hence for any $w\in O(V(x))$, 
\begin{equation}\label{Eq.4.14.e2}
    \frac{\|w\|^2}{2\|x-z\|}+\frac{1}{2}w^{\top}F(z)w=\frac{1}{2\|x-z\|}\cdot w^{\top}\big(\mathbf{I}_d+\|x-z\|F(z))w\geq 0.
\end{equation}  
Consequently, using \eqref{p1}, we obtain that
\begin{align}\label{Eq.4.14.e7}
        &\int_{\eps^{-1\slash 2}\mathscr{W}_0(x,z;\eps)}\bigg(\frac{\|w\|^2}{2\|x-z\|}+\frac{1}{2}w^{\top}F(z)w\bigg)\eps^{(d-1)\slash 2}p(x,z+\sqrt{\eps}w)d\mathcal{H}^{d-1}(w)\nonumber\\
    \leq\,& \int_{O(V(x))}\bigg(\frac{\|w\|^2}{2\|x-z\|}+\frac{1}{2}w^{\top}F(z)w\bigg)\eps^{(d-1)\slash 2}p_1(x,z+\sqrt{\eps}w)d\mathcal{H}^{d-1}(w)\nonumber\\
    \leq\,& (2\pi\|x-z\|)^{-(d-1) \slash  2}\mathfrak{F}(x)^{-1}h(x,z) \sqrt{\det\big(\mathbf{I}_d+\|x-z\|F(z)\big)}\nonumber\\
    &\times \int_{O(V(x))}   \bigg(\frac{\|w\|^2}{2\|x-z\|}+\frac{1}{2}w^{\top}F(z)w\bigg) e^{-\frac{1}{2}w^{\top}\big(\frac{1}{\|x-z\|}\mathbf{I}_d+F(z)\big)w}d\mathcal{H}^{d-1}(w)\nonumber\\
    =\,& \frac{1}{2}(2\pi\|x-z\|)^{-(d-1) \slash  2}\mathfrak{F}(x)^{-1}h(x,z) \sqrt{\det\big(\mathbf{I}_d+\|x-z\|F(z)\big)}\nonumber\\
    &\times \int_{O(V(x))}   \bigg(w^{\top}\Big(\frac{1}{\|x-z\|}\mathbf{I}_d+F(z)\Big)w\bigg) e^{-\frac{1}{2}w^{\top}\big(\frac{1}{\|x-z\|}\mathbf{I}_d+F(z)\big)w}d\mathcal{H}^{d-1}(w)\nonumber\\
    =\,& \frac{d-1}{2}\mathfrak{F}(x)^{-1}h(x,z).
\end{align}
By \eqref{p1} and \eqref{Eq.4.14.e2}, for any $w\in O(V(x))$, we have
\begin{align*}
   &\eps^{(d-1)\slash 2}p(x,z+\sqrt{\eps}w)\leq \eps^{(d-1)\slash 2}p_1(x,z+\sqrt{\eps}w)\nonumber\\
   \leq\,& (2\pi\|x-z\|)^{-(d-1) \slash  2}\mathfrak{F}(x)^{-1}h(x,z) \sqrt{\det\big(\mathbf{I}_d+\|x-z\|F(z)\big)}\leq C\delta_1^{-d\slash 2}\mathfrak{F}(x)^{-1}h(x,z),
\end{align*} 
where in the last inequality we note that $\|x-z\|\geq 2d_0$ and that by Hadamard's inequality  (see, e.g., \cite[Lemma 2.5]{ipsen2008perturbation}) and \Cref{L2.4} (note that $\min\{\alpha(z), \beta(z)\}\geq\delta_1$),
\begin{equation*}
    \big|\det\big(\mathbf{I}_d+\|x-z\|F(z)\big)\big|\leq \big\|\mathbf{I}_d+\|x-z\|F(z)\big\|_2 ^d \leq \big(1+\|x-z\|\|F(z)\|_2\big)^d
    \leq \big(1+2D\cdot 4\delta_1^{-1}\big)^d\leq C\delta_1^{-d}.
\end{equation*}
Hence we have 
\begin{align}\label{Eq.4.14.e8}
   &\int_{\eps^{-1\slash 2}\mathscr{W}(x,z;\delta,\eps)}\eps^{(d-1)\slash 2}p(x,z+\sqrt{\eps}w)d\mathcal{H}^{d-1}(w) \nonumber\\
   \leq\,&  C\delta_1^{-d\slash 2}\mathfrak{F}(x)^{-1}h(x,z)   \mathcal{H}^{d-1}\big(\big\{w\in O(V(x)):\|w\|\leq M_1\big\}\big)\leq CM_1^{d-1} \delta_1^{-d\slash 2} \mathfrak{F}(x)^{-1}h(x,z),
\end{align}
\begin{align}\label{Eq.4.14.e9}
  &\int_{\eps^{-1\slash 2}(\mathscr{W}_0(x,z;\eps)\backslash\mathscr{W}(x,z;\delta,\eps))}\eps^{(d-1)\slash 2}p(x,z+\sqrt{\eps}w)d\mathcal{H}^{d-1}(w)\nonumber\\
  \leq\,& C\delta_1^{-d\slash 2}\mathfrak{F}(x)^{-1}h(x,z) \cdot  \eps^{-(d-1)\slash 2}\mathcal{H}^{d-1}\big(\mathscr{W}_0(x,z;\eps)\backslash\mathscr{W}(x,z;\delta,\eps)\big),
\end{align}
where the second display uses
\begin{equation*}
   \mathcal{H}^{d-1}\big(\eps^{-1\slash 2}(\mathscr{W}_0(x,z;\eps)\backslash\mathscr{W}(x,z;\delta,\eps))\big)=\eps^{-(d-1)\slash 2}\mathcal{H}^{d-1}\big(\mathscr{W}_0(x,z;\eps)\backslash\mathscr{W}(x,z;\delta,\eps)\big).
\end{equation*}
Combining \eqref{Eq.4.14.e1} and \eqref{Eq.4.14.e7}--\eqref{Eq.4.14.e9}, we get
\begin{align}
    & \eps^{-1}  \int_{O(V(x))}\mathbbm{1}_{z+w\in\mathcal{Y}_{k'}^{\circ}\cap\mathfrak{T}^*\cap \tilde{\mathcal{T}}_{1;k,k'}^{*}}\mathbbm{1}_{V(x)+V(z+w)\neq 0}\mathbbm{1}_{\min\{\alpha(z+w),\beta(z+w)\}\geq\delta_1} \mathbbm{1}_{\|w\|\leq M_1\sqrt{\eps}}\nonumber\\
    &\hspace{0.6in}\times(\|x-(z+w)\|-\|x-\mathtt{y}'(x,z+w)\|)p(x,z+w)d\mathcal{H}^{d-1}(w)\nonumber\\
    \leq\,  &  \bigg(\frac{d-1}{2}+ C \delta M_1^{d+1}\delta_1^{-d\slash 2}+CM_1^2\delta_1^{-(d\slash 2+1)}\cdot \eps^{-(d-1)\slash 2}\mathcal{H}^{d-1}\big(\mathscr{W}_0(x,z;\eps)\backslash\mathscr{W}(x,z;\delta,\eps)\big)\bigg)\label{Lem.4.14.e8}\mathfrak{F}(x)^{-1}h(x,z)\\
    \leq\,& \bigg(\frac{d-1}{2}+ C M_1^{d+1}\delta_1^{-(d\slash 2+1)}\bigg) \mathfrak{F}(x)^{-1}h(x,z).\label{Lem.4.14.e9}
\end{align}

Note that by \eqref{Lem.4.14.e8}, for any $x\in\mathcal{X}_k^{\circ}\cap\mathfrak{T}^{*} \cap \tilde{\mathcal{T}}_{1;k,k'}^*$ and $z\in\mathrm{int}(T(x))\cap\mathcal{Y}_{k'}^{\circ}$ such that $\min\{\alpha(x),\beta(x)\}\geq \delta_1$, $\min\{\alpha(z),\beta(z)\}\geq \delta_1$, and $\langle x-z, V(x)\rangle>0$ and any $\delta\in (0,1)$, we have 
\begin{align}\label{Lem.4.14.e11}
    & \limsup_{\eps\rightarrow 0^+}\bigg\{\eps^{-1}\int_{O(V(x))}\mathbbm{1}_{z+w\in\mathcal{Y}_{k'}^{\circ}\cap\mathfrak{T}^*\cap \tilde{\mathcal{T}}_{1;k,k'}^{*}}\mathbbm{1}_{V(x)+V(z+w)\neq 0}\mathbbm{1}_{\min\{\alpha(z+w),\beta(z+w)\}\geq\delta_1} \mathbbm{1}_{\|w\|\leq M_1\sqrt{\eps}}\nonumber\\
    &\hspace{1in}\times(\|x-(z+w)\|-\|x-\mathtt{y}'(x,z+w)\|)p(x,z+w)d\mathcal{H}^{d-1}(w)\bigg\}\nonumber\\
    \leq\,& \mathfrak{F}(x)^{-1}h(x,z)\bigg(\frac{d-1}{2}+ C \delta M_1^{d+1}\delta_1^{-d\slash 2}+CM_1^2\delta_1^{-(d\slash 2+1)}\cdot \limsup_{\eps\rightarrow 0^+}\frac{\mathcal{H}^{d-1}\big(\mathscr{W}_0(x,z;\eps)\backslash\mathscr{W}(x,z;\delta,\eps)\big)}{\eps^{(d-1)\slash 2}}\bigg)\nonumber\\
    =\,&  \mathfrak{F}(x)^{-1}h(x,z)\bigg(\frac{d-1}{2}+ C \delta M_1^{d+1}\delta_1^{-d\slash 2}\bigg),
\end{align}
where we use \eqref{Lem4.14.e10} in the last line. Moreover, note that 
\begin{align}\label{Lem4.14.e12}
    &\int_{\mathcal{X}_k^{\circ}\cap\mathfrak{T}^{*} \cap \tilde{\mathcal{T}}_{1;k,k'}^{*}} \mathbbm{1}_{\min\{\alpha(x),\beta(x)\}\geq\delta_1} dx\nonumber\\
    &\hspace{0.2in}\int_{\mathrm{int}(T(x))\cap\mathcal{Y}_{k'}^{\circ}}\mathbbm{1}_{\min\{\alpha(z),\beta(z)\}\geq\delta_1}\mathbbm{1}_{\langle x-z, V(x)\rangle>0}\bigg(\frac{d-1}{2}+ C M_1^{d+1}\delta_1^{-(d \slash 2+1)}\bigg)\mathfrak{F}(x)^{-1}h(x,z) d\mathcal{H}^1(z)\nonumber\\
     \leq\,  &  \frac{d-1}{2}+ C  M_1^{d+1}\delta_1^{-(d\slash 2+1)}<\infty.
\end{align}
By \eqref{Lem.4.14.e10} and the reverse Fatou lemma (noting \eqref{Lem.4.14.e9} and \eqref{Lem4.14.e12}), we have 
\begin{align*}
    & \limsup_{\eps\rightarrow 0^{+}}\bigg\{\eps^{-1}\int_{(\mathcal{X}_k^{\circ}\cap\mathfrak{T}^*)\times(\mathcal{Y}_{k'}^{\circ}\cap\mathfrak{T}^*)}\mathbbm{1}_{V(x)+V(y)\neq 0}(\|x-y\|-\|x-\mathtt{y}'(x,y)\|)p(x,y)dxdy\bigg\} \nonumber\\
    =\,&  \limsup_{\eps\rightarrow 0^{+}}\bigg\{\int_{\mathcal{X}_k^{\circ}\cap\mathfrak{T}^{*} \cap \tilde{\mathcal{T}}_{1;k,k'}^{*}} \mathbbm{1}_{\min\{\alpha(x),\beta(x)\}\geq\delta_1} dx\int_{\mathrm{int}(T(x))\cap\mathcal{Y}_{k'}^{\circ}}\mathbbm{1}_{\min\{\alpha(z),\beta(z)\}\geq\delta_1} \mathbbm{1}_{\langle x-z, V(x)\rangle>0} d\mathcal{H}^1(z)\nonumber\\    &\hspace{0.8in}\eps^{-1}\int_{O(V(x))}\mathbbm{1}_{z+w\in\mathcal{Y}_{k'}^{\circ}\cap\mathfrak{T}^*\cap \tilde{\mathcal{T}}_{1;k,k'}^{*}}\mathbbm{1}_{V(x)+V(z+w)\neq 0}\mathbbm{1}_{\min\{\alpha(z+w),\beta(z+w)\}\geq\delta_1} \mathbbm{1}_{\|w\|\leq M_1\sqrt{\eps}}\nonumber\\
    &\hspace{1in}\times(\|x-(z+w)\|-\|x-\mathtt{y}'(x,z+w)\|)p(x,z+w)d\mathcal{H}^{d-1}(w)\bigg\}\nonumber\\
    \leq\,& \int_{\mathcal{X}_k^{\circ}\cap\mathfrak{T}^{*} \cap \tilde{\mathcal{T}}_{1;k,k'}^{*}} \mathbbm{1}_{\min\{\alpha(x),\beta(x)\}\geq\delta_1} dx\int_{\mathrm{int}(T(x))\cap\mathcal{Y}_{k'}^{\circ}}\mathbbm{1}_{\min\{\alpha(z),\beta(z)\}\geq\delta_1} \mathbbm{1}_{\langle x-z, V(x)\rangle>0} d\mathcal{H}^1(z)\nonumber\\
    &\hspace{0.2in}\limsup_{\eps\rightarrow 0^+}\bigg\{\eps^{-1}\int_{O(V(x))}\mathbbm{1}_{z+w\in\mathcal{Y}_{k'}^{\circ}\cap\mathfrak{T}^*\cap \tilde{\mathcal{T}}_{1;k,k'}^{*}}\mathbbm{1}_{V(x)+V(z+w)\neq 0}\mathbbm{1}_{\min\{\alpha(z+w),\beta(z+w)\}\geq\delta_1} \mathbbm{1}_{\|w\|\leq M_1\sqrt{\eps}}\nonumber\\
    &\hspace{1.11in}\times(\|x-(z+w)\|-\|x-\mathtt{y}'(x,z+w)\|)p(x,z+w)d\mathcal{H}^{d-1}(w)\bigg\} \nonumber\\
    \leq\,& \bigg(\frac{d-1}{2}+ C \delta M_1^{d+1}\delta_1^{-d   
     \slash 2}\bigg)\int_{\mathcal{X}_k^{\circ}\cap\mathfrak{T}^{*} \cap \tilde{\mathcal{T}}_{1;k,k'}^{*}} \mathbbm{1}_{\min\{\alpha(x),\beta(x)\}\geq\delta_1} dx\nonumber\\
    &\hspace{1.85in}\int_{\mathrm{int}(T(x))\cap\mathcal{Y}_{k'}^{\circ}}\mathbbm{1}_{\min\{\alpha(z),\beta(z)\}\geq\delta_1}\mathbbm{1}_{\langle x-z, V(x)\rangle>0} \mathfrak{F}(x)^{-1}h(x,z) d\mathcal{H}^1(z)\nonumber\\
    \leq\,&\bigg(\frac{d-1}{2}+ C \delta M_1^{d+1}\delta_1^{-d\slash 2}\bigg)\int_{\mathcal{X}_k^{\circ}\cap\mathfrak{T}^{*} \cap \tilde{\mathcal{T}}_{1;k,k'}^{*}}dx\int_{\mathrm{int}(T(x))\cap\mathcal{Y}_{k'}^{\circ}}\mathfrak{F}(x)^{-1}h(x,z) d\mathcal{H}^1(z),
\end{align*}
where we use \eqref{Lem.4.14.e11} in the second to last inequality. Taking $\delta\rightarrow 0^+$, we conclude that
\begin{align}\label{lem.neq.e9}
    & \limsup_{\eps\rightarrow 0^{+}}\bigg\{\eps^{-1}\int_{(\mathcal{X}_k^{\circ}\cap\mathfrak{T}^*)\times(\mathcal{Y}_{k'}^{\circ}\cap\mathfrak{T}^*)}\mathbbm{1}_{V(x)+V(y)\neq 0}(\|x-y\|-\|x-\mathtt{y}'(x,y)\|)p(x,y)dxdy\bigg\}\nonumber\\
    \leq \,&\frac{d-1}{2}\int_{\mathcal{X}_k^{\circ}\cap\mathfrak{T}^{*} \cap \tilde{\mathcal{T}}_{1;k,k'}^{*}}dx\int_{\mathrm{int}(T(x))\cap\mathcal{Y}_{k'}^{\circ}}\mathfrak{F}(x)^{-1}h(x,z) d\mathcal{H}^1(z). 
\end{align}

\paragraph{Step 2: Bounding $\int_{(\mathcal{X}_k^{\circ}\cap\mathfrak{T}^*)\times(\mathcal{Y}_{k'}^{\circ}\cap\mathfrak{T}^*)}\mathbbm{1}_{V(x)+V(y)\neq 0}\|x-\mathtt{y}'(x,y)\|p(x,y)dxdy$.} 

By Lemma \ref{Lem3.18nn}, 
\begin{align}\label{lem4.14.ee3}
    & \int_{(\mathcal{X}_k^{\circ}\cap\mathfrak{T}^*)\times(\mathcal{Y}_{k'}^{\circ}\cap\mathfrak{T}^*)}\mathbbm{1}_{V(x)+V(y)\neq 0}\|x-\mathtt{y}'(x,y)\|p(x,y)dxdy\nonumber\\
    =\,& \int_{\mathcal{S}}d\lambda(T)\int_{\mathrm{int}(T)\cap\mathcal{X}_k^{\circ}}\mathfrak{F}(x)d\mathcal{H}^1(x)\int_{\mathcal{Y}_{k'}^{\circ}\cap\mathfrak{T}^*}\mathbbm{1}_{V(x)+V(y)\neq 0}\|x-\mathtt{y}'(x,y)\|p(x,y)dy.
\end{align}
Moreover, by Lemma \ref{Lem3.18nn}, for any $T\in\mathcal{S}$ and $x\in\mathrm{int}(T)\cap\mathcal{X}_k^{\circ}$, we have
\begin{align}\label{lem4.14.ee4}
    &\int_{\mathcal{Y}_{k'}^{\circ}\cap\mathfrak{T}^*}\mathbbm{1}_{V(x)+V(y)\neq 0}\|x-\mathtt{y}'(x,y)\|p(x,y)dy\nonumber\\
    =\,& \int_{\mathcal{S}}d\lambda(T')\int_{\mathrm{int}(T')\cap\mathcal{Y}_{k'}^{\circ}}\mathbbm{1}_{V(x)+V(y)\neq 0}\|x-\mathtt{y}'(x,y)\|p(x,y)\mathfrak{F}(y)d\mathcal{H}^1(y)\nonumber\\
    =\,& \int_{\{T'\in\mathcal{S}:V_T+V_{T'}\neq 0\}}d\lambda(T')\int_{\mathrm{int}(T')\cap\mathcal{Y}_{k'}^{\circ}} \|x-\mathtt{y}'(x,y)\|p(x,y)\mathfrak{F}(y)d\mathcal{H}^1(y)\nonumber\\
    =\,& \int_{\{T'\in\mathcal{S}:V_T+V_{T'}\neq 0\}}d\lambda(T')\nonumber\\
    &\hspace{0.2in}\int_{\mathrm{int}(T)\cap\mathcal{Y}_{k'}^{\circ}}\mathbbm{1}_{\EuScript{R}(y';T,T')\in\mathrm{int}(T') \cap \mathcal{Y}_{k'}^{\circ}}\|x-y'\|p(x,\EuScript{R}(y';T,T'))\mathfrak{F}(\EuScript{R}(y';T,T'))d\mathcal{H}^1(y'),
\end{align}
where in the last line $\EuScript{R}(\cdot;T,T')$ is as defined in \eqref{defvarphi}, and we use the fact that for any $T,T'\in \mathcal{S}$ such that $V_T+V_{T'}\neq 0$ and any $x\in\mathrm{int}(T)\cap\mathcal{X}_k^{\circ},y\in \mathrm{int}(T')\cap\mathcal{Y}_{k'}^{\circ}$, if $\mathtt{y}'(x,y)\notin \mathrm{int}(T)\cap\mathcal{Y}_{k'}^{\circ}$, then by \eqref{p3}, $p(x,y)\leq p_3(x,y)=0$. Combining \eqref{lem4.14.ee3} and \eqref{lem4.14.ee4}, we get
\begin{align}\label{lem4.14.ee1}
  & \int_{(\mathcal{X}_k^{\circ}\cap\mathfrak{T}^*)\times(\mathcal{Y}_{k'}^{\circ}\cap\mathfrak{T}^*)}\mathbbm{1}_{V(x)+V(y)\neq 0}\|x-\mathtt{y}'(x,y)\|p(x,y)dxdy\nonumber\\
  \leq\, & \int_{\mathcal{S}}d\lambda(T)\int_{\mathrm{int}(T)\cap\mathcal{X}_k^{\circ}}\mathfrak{F}(x)d\mathcal{H}^1(x)  \int_{\{T'\in\mathcal{S}:V_T+V_{T'}\neq 0\}}d\lambda(T')\nonumber\\
    &\hspace{0.2in}\int_{\mathrm{int}(T)\cap\mathcal{Y}_{k'}^{\circ}}\mathbbm{1}_{\EuScript{R}(y';T,T')\in\mathrm{int}(T') \cap \mathcal{Y}_{k'}^{\circ}}\|x-y'\|p(x,\EuScript{R}(y';T,T'))\mathfrak{F}(\EuScript{R}(y';T,T'))d\mathcal{H}^1(y') \nonumber\\
    =\,& \int_{\mathcal{S}}d\lambda(T)\int_{(\mathrm{int}(T)\cap\mathcal{X}_k^{\circ})\times(\mathrm{int}(T)\cap\mathcal{Y}_{k'}^{\circ})}\|x-y'\|d\mathcal{H}^1(x)d\mathcal{H}^1(y')\nonumber\\
    &\hspace{0.5in}\mathfrak{F}(x)\int_{\{T'\in\mathcal{S}:V_T+V_{T'}\neq 0,
\EuScript{R}(y';T,T')\in\mathrm{int}(T')\cap\mathcal{Y}_{k'}^{\circ}\}}p(x,\EuScript{R}(y';T,T'))\mathfrak{F}(\EuScript{R}(y';T,T'))d\lambda(T')\nonumber\\
 = \, & \int_{\mathcal{S}}d\lambda(T)\int_{(\mathrm{int}(T)\cap\mathcal{X}_k^{\circ})\times(\mathrm{int}(T)\cap\mathcal{Y}_{k'}^{\circ})}\|x-y\|\xi_T(x,y)d\mathcal{H}^1(x)d\mathcal{H}^1(y),
\end{align}
where in the last line we use the definition of $\xi_T(\cdot,\cdot)$ as in \eqref{defxiq}. By \eqref{lem4.14.ee1} and the definition of $\bar{\gamma}_{\eps,1}$ (see \eqref{def.gamma.1}), we have
\begin{equation}\label{lem4.14.ee2}
    (1-\delta_3)\int_{(\mathcal{X}_k^{\circ}\cap\mathfrak{T}^*)\times(\mathcal{Y}_{k'}^{\circ}\cap\mathfrak{T}^*)}\mathbbm{1}_{V(x)+V(y)\neq 0}\|x-\mathtt{y}'(x,y)\|p(x,y)dxdy\leq\int_{\mathbb{R}^d\times\mathbb{R}^d}\|x-y\|d\bar{\gamma}_{\eps,1}(x,y). 
\end{equation}
By \eqref{def.gamma.1} and Lemma \ref{Lem4.6}, for $\bar{\gamma}_{\eps,1}$-a.e.\ $(x,y)\in\mathbb{R}^d\times\mathbb{R}^d$, $x$ and $y$ lie on the same transport ray and $\|x-y\|=u(x)-u(y)$. Hence by \eqref{lem4.14.ee2} and Lemma \ref{Lem4.11}, we have
\begin{align}\label{lem.neq.e11}
   & (1-\delta_3)\int_{(\mathcal{X}_k^{\circ}\cap\mathfrak{T}^*)\times(\mathcal{Y}_{k'}^{\circ}\cap\mathfrak{T}^*)}\mathbbm{1}_{V(x)+V(y)\neq 0}\|x-\mathtt{y}'(x,y)\|p(x,y)dxdy\nonumber\\
   \leq\, & \int_{\mathbb{R}^d\times\mathbb{R}^d}(u(x)-u(y))d\bar{\gamma}_{\eps,1}(x,y)=\int_{\mathbb{R}^d}u(x)\bar{f}_1(x)dx-\int_{\mathbb{R}^d}u(y)\bar{g}_1(y)dy.
\end{align}

\paragraph{Step 3: Bounding $\int_{\mathbb{R}^d\times\mathbb{R}^d}\|x-y\|r(x,y)dxdy$.} By \eqref{def.gamma.2} and Lemma \ref{Lem4.6}, for $\bar{\gamma}_{\eps,2}$-a.e.\ $(x,y)\in\mathbb{R}^d\times\mathbb{R}^d$, $x$ and $y$ lie on the same transport ray and $u(x)-u(y)=\|x-y\|$. Hence by Lemma \ref{Lem4.9}, 
\begin{equation}\label{lem.neq.e15}
    \int_{\mathbb{R}^d\times\mathbb{R}^d}\|x-y\|d\bar{\gamma}_{\eps,2}(x,y)=\int_{\mathbb{R}^d\times\mathbb{R}^d}(u(x)-u(y))d\bar{\gamma}_{\eps,2}(x,y)
    =\int_{\mathbb{R}^d}u(x)\bar{f}_2(x)dx-\int_{\mathbb{R}^d}u(y)\bar{g}_2(y)dy. 
\end{equation}

By the definition of $r(\cdot,\cdot)$ as in \eqref{def:r} and noting that $\bar{\gamma}_{\eps,2}((\mathbb{R}^d\times\mathbb{R}^d)\backslash (\mathcal{X}_k^{\circ}\times\mathcal{Y}_{k'}^{\circ}))=0$ (by \eqref{def.gamma.2}), we have 
\begin{align}\label{lem.neq.e10}
&\int_{\mathbb{R}^d\times\mathbb{R}^d}\|x-y\|r(x,y)dxdy-\int_{\mathbb{R}^d\times\mathbb{R}^d}\|x-y\|d\bar{\gamma}_{\eps,2}(x,y) \nonumber\\
=\,&\sum_{i=1}^I\sum_{j=1}^J\sum_{\substack{\ell,\ell'\in [N_0]:\\ 
 \bar{\gamma}_{\eps,2}(L_{i,j;\ell}\times R_{i,j;\ell'})>0}}\bigg(\int_{L_{i,j;\ell}\times R_{i,j;\ell'}}\|x-y\|r(x,y)dxdy-\int_{L_{i,j;\ell}\times R_{i,j;\ell'}}\|x-y\|d\bar{\gamma}_{\eps,2}(x,y)\bigg).
\end{align}
For any $i\in[I],j\in [J]$ and $\ell,\ell'\in[N_0]$ such that $\bar{\gamma}_{\eps,2}(L_{i,j;\ell}\times R_{i,j;\ell'})>0$, let $\Upsilon_{i,j;\ell,\ell'}$ be the Borel measure on $L_{i,j;\ell}\times R_{i,j;\ell'}\times L_{i,j;\ell}\times R_{i,j;\ell'}$ such that for any Borel sets $A,A'\subseteq L_{i,j;\ell}\times R_{i,j;\ell'}$,
\begin{equation*}
     \Upsilon_{i,j;\ell,\ell'}(A\times A')=\frac{\bar{\gamma}_{\eps,2}(A')\int_{A}f_{j}(x)g_{i}(y)dxdy}{\bar{\gamma}_{\eps,2}(L_{i,j;\ell}\times R_j) \bar{\gamma}_{\eps,2}(L_i\times R_{i,j;\ell'})}.
\end{equation*}
Note that by Lemma \ref{Lem4.9n}, for any Borel sets $A,A'\subseteq L_{i,j;\ell}\times R_{i,j;\ell'}$,
\begin{equation*}
    \Upsilon_{i,j;\ell,\ell'}((L_{i,j;\ell}\times R_{i,j;\ell'})\times A')=\bar{\gamma}_{\eps,2}(A'),\qquad \Upsilon_{i,j;\ell,\ell'}(A\times (L_{i,j;\ell}\times R_{i,j;\ell'}))=\int_{A} r(x,y) dxdy.
\end{equation*}
Hence
\begin{align}\label{lem.neq.e7}
    & \int_{L_{i,j;\ell}\times R_{i,j;\ell'}}\|x-y\|r(x,y)dxdy-\int_{L_{i,j;\ell}\times R_{i,j;\ell'}}\|x-y\|d\bar{\gamma}_{\eps,2}(x,y)\nonumber\\
    =\,& \int_{L_{i,j;\ell}\times R_{i,j;\ell'}\times L_{i,j;\ell}\times R_{i,j;\ell'}}(\|x-y\|-\|x'-y'\|)d\Upsilon_{i,j;\ell,\ell'}(x,y,x',y'). 
\end{align}
Now note that by \eqref{Lijk} and \eqref{Rijk}, for any $(x,y),(x',y')\in L_{i,j;\ell}\times R_{i,j;\ell'}$, we have
\begin{equation}\label{lem.neq.e1}
    |\langle U_{i,j}, x-x'\rangle|\leq \delta_2\eps, \qquad |\langle U_{i,j}, y-y'\rangle|\leq \delta_2\eps. 
\end{equation}
As $x,x'\in L_i$, $y,y'\in R_j$, $\dist(L_i,R_j)\geq 2d_0$, and $L_i,R_j\subseteq B_d(0,D+1)$ (by \eqref{distcondi}), we have
\begin{equation}\label{lem.neq.e2}
    \max\{\|x-x'\|,\|y-y'\|\}\leq \delta_2\sqrt{d\eps}, \quad d_0\leq\|x-y\|,\|x'-y'\|\leq 2(D+1).
\end{equation}
Moreover, by the definition of $U_{i,j}$ as in \eqref{defU}, we have
\begin{equation*}
    x-y-\langle x-y, U_{i,j}\rangle U_{i,j}=(x-\mathtt{c}_i)-(y-\mathtt{c}'_j)-\langle (x-\mathtt{c}_i)-(y-\mathtt{c}'_j), U_{i,j}\rangle U_{i,j},
\end{equation*}
hence
\begin{equation}\label{lem.neq.e3}
    \|x-y-\langle x-y, U_{i,j}\rangle U_{i,j}\|\leq 2(\|x-\mathtt{c}_i\|+\|y-\mathtt{c}'_j\|)\leq 4\delta_2\sqrt{d\eps}.
\end{equation}
Similarly, we have
\begin{equation}\label{lem.neq.e4}
    \|x'-y'-\langle x'-y', U_{i,j}\rangle U_{i,j}\|\leq 4\delta_2\sqrt{d\eps}.
\end{equation}
By \eqref{lem.neq.e1}--\eqref{lem.neq.e4}, we have 
\begin{align}\label{lem.neq.e5}
  & \big||\langle x-y, U_{i,j}\rangle |^2-|\langle x'-y', U_{i,j}\rangle|^2\big|\nonumber\\
  =\,& \big||\langle x-y, U_{i,j}\rangle|-|\langle x'-y', U_{i,j}\rangle|\big|(|\langle x-y, U_{i,j}\rangle|+|\langle x'-y', U_{i,j}\rangle|)\nonumber\\
  \leq\,& 2|\langle (x-y)-(x'-y'),U_{i,j}\rangle|(\|x-y\|+\|x'-y'\|)\nonumber\\
  \leq\,& 2(|\langle x-x',U_{i,j}\rangle|+|\langle y-y',U_{i,j}\rangle|)(\|x-y\|+\|x'-y'\|)\leq C\delta_2\eps,
\end{align}
\begin{align}\label{lem.neq.e6}
    & \big|\|x-y-\langle x-y, U_{i,j}\rangle U_{i,j}\|^2-\|x'-y'-\langle x'-y', U_{i,j}\rangle U_{i,j}\|^2\big| \nonumber\\
    =\,& \big|\|x-y-\langle x-y, U_{i,j}\rangle U_{i,j}\|-\|x'-y'-\langle x'-y', U_{i,j}\rangle U_{i,j}\|\big|\nonumber\\
    &\times (\|x-y-\langle x-y, U_{i,j}\rangle U_{i,j}\|+\|x'-y'-\langle x'-y', U_{i,j}\rangle U_{i,j}\|)\nonumber\\
    \leq\,& 8\delta_2\sqrt{d\eps}\|(x-y)-(x'-y')-\langle (x-y)-(x'-y'),U_{i,j}\rangle U_{i,j}\| \nonumber\\
    \leq\,& 16\delta_2\sqrt{d\eps}(\|x-x'\|+\|y-y'\|)\leq C\delta_2^2\eps.
\end{align}
As $x-y=\langle x-y, U_{i,j}\rangle U_{i,j}+\big(x-y-\langle x-y, U_{i,j}\rangle U_{i,j}\big)$ and $\big\langle \langle x-y, U_{i,j}\rangle U_{i,j},x-y-\langle x-y, U_{i,j}\rangle U_{i,j}\big\rangle =0$, 
\begin{align*}
    \|x-y\|=\,&\sqrt{\|\langle x-y, U_{i,j}\rangle U_{i,j}\|^2+\|x-y-\langle x-y, U_{i,j}\rangle U_{i,j}\|^2}\nonumber\\
    =\,& \sqrt{|\langle x-y, U_{i,j}\rangle|^2+\|x-y-\langle x-y, U_{i,j}\rangle U_{i,j}\|^2}.
\end{align*}
Similarly, 
\begin{equation*}
    \|x'-y'\|=\sqrt{|\langle x'-y', U_{i,j}\rangle|^2+\|x'-y'-\langle x'-y', U_{i,j}\rangle U_{i,j}\|^2}.
\end{equation*}
Hence
\begin{align}\label{lem.neq.e8}
   & \big|\|x-y\|-\|x'-y'\|\big|=\frac{\big|\|x-y\|^2-\|x'-y'\|^2\big|}{\|x-y\|+\|x'-y'\|}\nonumber\\
   \leq\,&\frac{\big|\|x-y-\langle x-y, U_{i,j}\rangle U_{i,j}\|^2-\|x'-y'-\langle x'-y', U_{i,j}\rangle U_{i,j}\|^2\big|}{2d_0}\nonumber\\
   &+\frac{\big||\langle x-y, U_{i,j}\rangle |^2-|\langle x'-y', U_{i,j}\rangle|^2\big|}{2d_0}\leq C\delta_2\eps,
\end{align}
where we use \eqref{lem.neq.e2} in the first inequality and \eqref{lem.neq.e5}--\eqref{lem.neq.e6} in the second inequality. Combining \eqref{lem.neq.e7} and \eqref{lem.neq.e8}, we get
\begin{align*}
& \bigg|\int_{L_{i,j;\ell}\times R_{i,j;\ell'}}\|x-y\|r(x,y)dxdy-\int_{L_{i,j;\ell}\times R_{i,j;\ell'}}\|x-y\|d\bar{\gamma}_{\eps,2}(x,y)\bigg|\nonumber\\
\leq\,&  C\delta_2\eps\Upsilon_{i,j;\ell,\ell'}(L_{i,j;\ell}\times R_{i,j;\ell'}\times L_{i,j;\ell}\times R_{i,j;\ell'})=C\delta_2\eps\int_{L_{i,j;\ell}\times R_{i,j;\ell'}}r(x,y)dxdy.
\end{align*}
Hence by \eqref{lem.neq.e10}, we have
\begin{align}\label{lem.neq.e14}
&  \bigg|\int_{\mathbb{R}^d\times\mathbb{R}^d}\|x-y\|r(x,y)dxdy-\int_{\mathbb{R}^d\times\mathbb{R}^d}\|x-y\|d\bar{\gamma}_{\eps,2}(x,y)\bigg|
\nonumber\\
\leq\,& C\delta_2\eps\sum_{i=1}^I\sum_{j=1}^J\sum_{\ell,\ell'\in[N_0]}\int_{L_{i,j;\ell}\times R_{i,j;\ell'}}r(x,y)dxdy\nonumber\\
\leq\,& C\delta_2\eps\int_{\mathbb{R}^d\times\mathbb{R}^d}r(x,y)dxdy=C\delta_2\eps\bigg(\mu_{k,k'}(\mathbb{R}^d)-(1-\delta_3)\int_{\mathbb{R}^d\times\mathbb{R}^d}p(x,y)dxdy\bigg),
\end{align}
where we note Lemma \ref{lemmacoup} in the last equality.

\paragraph{Step 4: Conclusion.} Combining \eqref{lem.neq.e12}, \eqref{lem.neq.e11}, \eqref{lem.neq.e15}, and \eqref{lem.neq.e14}, we obtain 
\begin{align*}
    &\int_{\mathcal{X}_{k}^{\circ}\times\mathcal{Y}_{k'}^{\circ}}\|x-y\|d\tilde{\gamma}_{\eps;k,k'}(x,y)\nonumber\\
    \leq\,& (1-\delta_3)\int_{(\mathcal{X}_k^{\circ}\cap\mathfrak{T}^*)\times(\mathcal{Y}_{k'}^{\circ}\cap\mathfrak{T}^*)}\mathbbm{1}_{V(x)+V(y)\neq 0}(\|x-y\|-\|x-\mathtt{y}'(x,y)\|)p(x,y)dxdy\nonumber\\
    &+\int_{\mathbb{R}^d}u(x)\bar{f}_1(x)dx-\int_{\mathbb{R}^d}u(y)\bar{g}_1(y)dy+\int_{\mathbb{R}^d}u(x)\bar{f}_2(x)dx-\int_{\mathbb{R}^d}u(y)\bar{g}_2(y)dy\nonumber\\
    &+C\delta_2\eps\bigg(\mu_{k,k'}(\mathbb{R}^d)-(1-\delta_3)\int_{\mathbb{R}^d\times\mathbb{R}^d}p(x,y)dxdy\bigg).
\end{align*}
By \eqref{Eq.4.2new}, recalling \eqref{def:OT}, we have
\begin{align*}
    & \int_{\mathbb{R}^d}u(x)\bar{f}_1(x)dx-\int_{\mathbb{R}^d}u(y)\bar{g}_1(y)dy+\int_{\mathbb{R}^d}u(x)\bar{f}_2(x)dx-\int_{\mathbb{R}^d}u(y)\bar{g}_2(y)dy\nonumber\\
    =\,& \int_{\mathbb{R}^d}u(x)f_{k,k'}(x)dx-\int_{\mathbb{R}^d}u(y)g_{k,k'}(y)dy=\OT_{k,k'}(\mu,\nu).
\end{align*}
Combining the above two displays, we get
\begin{align*}
&\eps^{-1}\bigg(\int_{\mathcal{X}_k^{\circ}\times\mathcal{Y}_{k'}^{\circ}}\|x-y\|d\tilde{\gamma}_{\eps;k,k'}(x,y)-\OT_{k,k'}(\mu,\nu)\bigg)\nonumber\\
\leq\,& (1-\delta_3)\eps^{-1}\int_{(\mathcal{X}_k^{\circ}\cap\mathfrak{T}^*)\times(\mathcal{Y}_{k'}^{\circ}\cap\mathfrak{T}^*)}\mathbbm{1}_{V(x)+V(y)\neq 0}(\|x-y\|-\|x-\mathtt{y}'(x,y)\|)p(x,y)dxdy\nonumber\\
&+C\delta_2\bigg(\mu_{k,k'}(\mathbb{R}^d)-(1-\delta_3)\int_{\mathbb{R}^d\times\mathbb{R}^d}p(x,y)dxdy\bigg).
\end{align*}
Hence by \eqref{lem.neq.e9} and Proposition \ref{P4.3}, we get (note that $\mu_{k,k'}(\mathbb{R}^d)\leq 1$ by Lemma \ref{Ln4.1}) 
\begin{align*}
 & \limsup_{\delta_1\rightarrow 0^+}\limsup_{M_1\rightarrow\infty}\limsup_{\eps\rightarrow 0^+}\eps^{-1}\bigg(\int_{\mathcal{X}_k^{\circ}\times\mathcal{Y}_{k'}^{\circ}}\|x-y\|d\tilde{\gamma}_{\eps;k,k'}(x,y)-\OT_{k,k'}(\mu,\nu)\bigg)  \nonumber\\
 \leq \,  &\frac{(1-\delta_3)(d-1)}{2} \int_{\mathcal{X}_k^{\circ}\cap\mathfrak{T}^{*} \cap \tilde{\mathcal{T}}_{1;k,k'}^{*}}dx\int_{\mathrm{int}(T(x))\cap\mathcal{Y}_{k'}^{\circ}}\mathfrak{F}(x)^{-1}h(x,z) d\mathcal{H}^1(z)+C\delta_2\delta_3.
\end{align*}
Taking $\delta_3\rightarrow 0^+$ in the above display, we conclude that 
\begin{align*}
   & \limsup_{\delta_3\rightarrow 0^+}\limsup_{\delta_1\rightarrow 0^+}\limsup_{M_1\rightarrow\infty}\limsup_{\eps\rightarrow 0^+}\eps^{-1}\bigg(\int_{\mathcal{X}_k^{\circ}\times\mathcal{Y}_{k'}^{\circ}}\|x-y\|d\tilde{\gamma}_{\eps;k,k'}(x,y)-\OT_{k,k'}(\mu,\nu)\bigg)\nonumber\\
   \leq\,& \frac{d-1}{2}\int_{\mathcal{X}_k^{\circ}\cap\mathfrak{T}^*\cap\tilde{\mathcal{T}}_{1;k,k'}^*} dx \int_{\mathrm{int}(T(x))\cap\mathcal{Y}_{k'}^{\circ}}\mathfrak{F}(x)^{-1}h(x,z) d\mathcal{H}^1(z). \qedhere
\end{align*}
\end{proof}

For any $x,y\in\mathbb{R}^d$, we define 
\begin{align}\label{p1til}
   \tilde{p}_1 (x,y)   :=\,& \mathbbm{1}_{x\in\mathcal{X}_k^{\circ}\cap \mathfrak{T}^*
\cap\tilde{\mathcal{T}}_{1;k,k'}^{*}}\mathbbm{1}_{\mathtt{z}(x,y)\in\mathrm{int}(T(x))\cap\mathcal{Y}_{k'}^{\circ}}\mathbbm{1}_{\langle x-\mathtt{z}(x,y),V(x)\rangle>0}\mathfrak{F}(x)^{-1}h(x,\mathtt{z}(x,y))\nonumber\\
   & \times (2\pi\|x-\mathtt{z}(x,y)\|\eps)^{-(d-1)\slash 2} \sqrt{\det\big(\mathbf{I}_d+\|x-\mathtt{z}(x,y)\|F(\mathtt{z}(x,y))\big)}\nonumber\\
   & \times e^{-\frac{1}{2\eps}\mathtt{w}(x,y)^{\top}\big(\frac{1}{\|x-\mathtt{z}(x,y)\|}\mathbf{I}_d+F(\mathtt{z}(x,y))\big)\mathtt{w}(x,y)}.
\end{align}
Note that $p_1(x,y)\leq\tilde{p}_1(x,y)$ (recall the definition of $p_1$ in \eqref{p1}). The next lemma computes several integrals involving $\tilde{p}_1$ for later reference. 

\begin{lemma}\label{Lemma4.15}
We have 
\begin{equation}\label{eq0L4.15}
    \int_{\mathbb{R}^d\times\mathbb{R}^d}p_1(x,y)dxdy \leq \int_{\mathbb{R}^d\times\mathbb{R}^d}\tilde{p}_1(x,y)dxdy=\mu_{k,k'}(\mathbb{R}^d),
\end{equation}
\begin{align}\label{eq2L4.15}
    &\int_{\mathbb{R}^d\times\mathbb{R}^d}\tilde{p}_1(x,y)\log\big(\eps^{(d-1)\slash 2}\tilde{p}_1(x,y)\big)dxdy\nonumber\\
    =\,&\int_{\mathcal{X}_k^{\circ}\cap\mathfrak{T}^*\cap \tilde{\mathcal{T}}_{1;k,k'}^{*}}dx\int_{\mathrm{int}(T(x))\cap\mathcal{Y}_{k'}^{\circ}}\mathfrak{F}(x)^{-1}h(x,z)\bigg(\log\big(\mathfrak{F}(x)^{-1}h(x,z)\big)-\frac{d-1}{2}\log(2\pi e \|x-z\|)\nonumber\\
& \hspace{2.95in}+\frac{1}{2}\log\det\big(\mathbf{I}_d+\|x-z\| F(z)\big)\bigg)d\mathcal{H}^1(z),
\end{align}
\begin{align}\label{eq3L4.15}
&\int_{\mathcal{X}_k^{\circ}\cap\mathfrak{T}^*\cap \tilde{\mathcal{T}}_{1;k,k'}^{*}}dx\int_{\mathrm{int}(T(x))\cap\mathcal{Y}_{k'}^{\circ}}\mathfrak{F}(x)^{-1}h(x,z)\bigg(\big|\log\big(\mathfrak{F}(x)^{-1}h(x,z)\big)\big|+\frac{d-1}{2}\big|\log(2\pi\|x-z\|)\big| \nonumber\\
&\hspace{2.9in}+\frac{1}{2}\big|\log\det\big(\mathbf{I}_d+\|x-z\| F(z)\big)\big|\bigg)d\mathcal{H}^1(z)\leq C.
\end{align}
Moreover, for any $M\geq 0$, we have
\begin{align}\label{eq1L4.15}
    & \int_{\mathbb{R}^d\times\mathbb{R}^d}\tilde{p}_1(x,y)\max\big\{\log\big(\eps^{(d-1)\slash 2}\tilde{p}_1(x,y)\big), -M\big\}dxdy\nonumber\\
    =\,& \int_{\mathcal{X}_k^{\circ}\cap\mathfrak{T}^*\cap \tilde{\mathcal{T}}_{1;k,k'}^{*}}dx\int_{\mathrm{int}(T(x))\cap\mathcal{Y}_{k'}^{\circ}}\mathfrak{F}(x)^{-1}h(x,z)(2\pi\|x-z\|)^{-(d-1)\slash 2}\sqrt{\det\big(\mathbf{I}_d+\|x-z\|F(z)\big)}d\mathcal{H}^1(z)\nonumber\\
&\hspace{0.1in}\int_{O(V(x))}\max\bigg\{\log\big(\mathfrak{F}(x)^{-1}h(x,z)\big)-\frac{d-1}{2}\log(2\pi\|x-z\|) +\frac{1}{2}\log\det\big(\mathbf{I}_d+\|x-z\| F(z)\big)\nonumber\\
&\hspace{1.05in}-\frac{1}{2}w^{\top}\Big(\frac{1}{\|x-z\|}\mathbf{I}_d+F(z)\Big)w,-M\bigg\} e^{-\frac{1}{2}w^{\top}\big(\frac{1}{\|x-z\|}\mathbf{I}_d+F(z)\big)w}d\mathcal{H}^{d-1}(w).
\end{align}
\end{lemma}
\begin{proof}

For any $x\in\mathcal{X}_k^{\circ}\cap\mathfrak{T}^*\cap \tilde{\mathcal{T}}_{1;k,k'}^{*}$ and $z\in \mathrm{int}(T(x))\cap\mathcal{Y}_{k'}^{\circ}$ such that $\langle x-z,V(x)\rangle>0$ (which is the case if $h(x,z)>0$), $\frac{1}{\|x-z\|}\mathbf{I}_d+F(z)$ is positive definite (by Lemma \ref{L2.7}), hence
\begin{equation}\label{Lem4.15.n1}
     (2\pi\|x-z\|)^{-(d-1)\slash 2}\sqrt{\det\big(\mathbf{I}_d+\|x-z\|F(z)\big)} \int_{O(V(x))} e^{-\frac{1}{2}w^{\top}\big(\frac{1}{\|x-z\|}\mathbf{I}_d+F(z)\big)w}d\mathcal{H}^{d-1}(w)=1,
\end{equation}
\begin{align}\label{Lem4.15.n2}
  & (2\pi\|x-z\|)^{-(d-1)\slash 2}\sqrt{\det\big(\mathbf{I}_d+\|x-z\|F(z)\big)}\nonumber\\
  &\hspace{0.1in} \times \int_{O(V(x))}\frac{1}{2}w^{\top}\Big(\frac{1}{\|x-z\|}\mathbf{I}_d+F(z)\Big)w \cdot  e^{-\frac{1}{2}w^{\top}\big(\frac{1}{\|x-z\|}\mathbf{I}_d+F(z)\big)w}d\mathcal{H}^{d-1}(w)=\frac{d-1}{2}.
\end{align}
By \eqref{fkkdef} and \eqref{Lem4.15.n1}, we have
\begin{align*}
&\int_{\mathbb{R}^d\times\mathbb{R}^d}p_1(x,y)dxdy\leq \int_{\mathbb{R}^d\times\mathbb{R}^d}\tilde{p}_1(x,y)dxdy\nonumber\\
=\,&\int_{\mathcal{X}_k^{\circ}\cap\mathfrak{T}^*\cap \tilde{\mathcal{T}}_{1;k,k'}^{*}}dx\int_{\mathrm{int}(T(x))\cap\mathcal{Y}_{k'}^{\circ}}\mathfrak{F}(x)^{-1}h(x,z) (2\pi\|x-z\|)^{-(d-1)\slash 2}\sqrt{\det\big(\mathbf{I}_d+\|x-z\|F(z)\big)} d\mathcal{H}^1(z)\nonumber\\
    &\hspace{2in} \int_{O(V(x))}e^{-\frac{1}{2}w^{\top}\big(\frac{1}{\|x-z\|}\mathbf{I}_d+F(z)\big)w}d\mathcal{H}^{d-1}(w)\nonumber\\
    =\,&\int_{\mathcal{X}_k^{\circ}\cap\mathfrak{T}^*\cap \tilde{\mathcal{T}}_{1;k,k'}^{*}}dx\int_{\mathrm{int}(T(x))\cap\mathcal{Y}_{k'}^{\circ}}\mathfrak{F}(x)^{-1}h(x,z)d\mathcal{H}^1(z)=\int_{\mathbb{R}^d}f_{k,k'}(x)dx=\mu_{k,k'}(\mathbb{R}^d),
\end{align*}
which establishes \eqref{eq0L4.15}. Moreover, 
\begin{align}\label{Lem.4.15.ee1}
&\int_{\mathbb{R}^d\times\mathbb{R}^d}\tilde{p}_1(x,y)\log\big(\eps^{(d-1)\slash 2}\tilde{p}_1(x,y)\big)dxdy\nonumber\\
=\,& \int_{\mathcal{X}_k^{\circ}\cap\mathfrak{T}^*\cap \tilde{\mathcal{T}}_{1;k,k'}^{*}}dx\int_{\mathrm{int}(T(x))\cap\mathcal{Y}_{k'}^{\circ}} \mathfrak{F}(x)^{-1}h(x,z)(2\pi\|x-z\|)^{-(d-1)\slash 2}\sqrt{\det\big(\mathbf{I}_d+\|x-z\|F(z)\big)} d\mathcal{H}^1(z)\nonumber\\
&\hspace{0.55in}\int_{O(V(x))}\bigg(\log\big(\mathfrak{F}(x)^{-1}h(x,z)\big)-\frac{d-1}{2}\log(2\pi\|x-z\|) +\frac{1}{2}\log\det\big(\mathbf{I}_d+\|x-z\|F(z)\big)\nonumber\\
&\hspace{1.25in}-\frac{1}{2}w^{\top}\Big(\frac{1}{\|x-z\|}\mathbf{I}_d+F(z)\Big)w\bigg) e^{-\frac{1}{2}w^{\top}\big(\frac{1}{\|x-z\|}\mathbf{I}_d+F(z)\big)w}d\mathcal{H}^{d-1}(w).
\end{align}
Similarly, for any $M\geq 0$, \eqref{eq1L4.15} holds. By \eqref{Lem4.15.n1}--\eqref{Lem.4.15.ee1}, we obtain \eqref{eq2L4.15}.

Now note that by Lemma \ref{Lem3.18nn}, we have 
\begin{align}\label{Eqnn4.0}
    & \int_{\mathcal{X}_k^{\circ}\cap\mathfrak{T}^*\cap \tilde{\mathcal{T}}_{1;k,k'}^{*}}dx\int_{\mathrm{int}(T(x))\cap\mathcal{Y}_{k'}^{\circ}}\mathfrak{F}(x)^{-1}h(x,z)\bigg(\big|\log\big(\mathfrak{F}(x)^{-1}h(x,z)\big)\big|+\frac{d-1}{2}\big|\log(2\pi\|x-z\|)\big| \nonumber\\
&\hspace{2.9in}+\frac{1}{2}\big|\log\det\big(\mathbf{I}_d+\|x-z\| F(z)\big)\big|\bigg)d\mathcal{H}^1(z)\nonumber\\
=\,& \int_{\mathcal{S}}d\lambda(T)\int_{(\mathrm{int}(T)\cap\mathcal{X}_k^{\circ})\times(\mathrm{int}(T)\cap\mathcal{Y}_{k'}^{\circ})}\mathbbm{1}_{x\in\tilde{\mathcal{T}}_{1;k,k'}^*} h(x,z)\bigg(\big|\log\big(\mathfrak{F}(x)^{-1}h(x,z)\big)\big|+\frac{d-1}{2}\big|\log(2\pi\|x-z\|)\big| \nonumber\\
&\hspace{3.15in}+\frac{1}{2}\big|\log\det\big(\mathbf{I}_d+\|x-z\| F(z)\big)\big|\bigg)d\mathcal{H}^1(x)d\mathcal{H}^1(z)\nonumber\\
\leq\,& \int_{\tilde{\mathcal{S}}^*}
d\lambda(T)\int_{(\mathrm{int}(T)\cap\mathcal{X}_k^{\circ})\times(\mathrm{int}(T)\cap\mathcal{Y}_{k'}^{\circ})}
h(x,z)\bigg(\big|\log\big(\mathfrak{F}(x)^{-1}h(x,z)\big)\big|+\frac{d-1}{2}\big|\log(2\pi\|x-z\|)\big| \nonumber\\
&\hspace{2.7in}+\frac{1}{2}\big|\log\det\big(\mathbf{I}_d+\|x-z\| F(z)\big)\big|\bigg)d\mathcal{H}^1(x)d\mathcal{H}^1(z),
\end{align}
where we use the fact that $\lambda(\mathcal{S}\backslash\tilde{\mathcal{S}}^*)=0$ (recall Definition \ref{defstilde}) in the last inequality.

Fix any $T\in\tilde{\mathcal{S}}^*$. If $h(x,z)=0$ for $\mathcal{H}^1\otimes\mathcal{H}^1$-a.e.\ $(x,z)\in (\mathrm{int}(T)\cap\mathcal{X}_k^{\circ})\times(\mathrm{int}(T)\cap\mathcal{Y}_{k'}^{\circ})$, we have 
\begin{equation}\label{Eqnn4.1}
  \int_{(\mathrm{int}(T)\cap\mathcal{X}_k^{\circ})\times(\mathrm{int}(T)\cap\mathcal{Y}_{k'}^{\circ})}
h(x,z)\big|\log\big(\mathfrak{F}(x)^{-1}h(x,z)\big)\big|d\mathcal{H}^1(x)d\mathcal{H}^1(z)=0.
\end{equation}
Otherwise, we have
\begin{equation}\label{bddinte}
    \int_{(\mathrm{int}(T)\cap\mathcal{X}_k^{\circ})\times(\mathrm{int}(T)\cap\mathcal{Y}_{k'}^{\circ})}h(x,z) d\mathcal{H}^1(x)d\mathcal{H}^1(z)\in (0,1].
\end{equation}
This implies the existence of some $(x_0,z_0)\in(\mathrm{int}(T)\cap\mathcal{X}_k^{\circ})\times(\mathrm{int}(T)\cap\mathcal{Y}_{k'}^{\circ})$ with $h_T(x_0,z_0)>0$ (recall Definition \ref{DefhH}). By \Cref{Lem2}(d) and \eqref{bddsalphabeta}, for any $(x,z)\in(\mathrm{int}(T)\cap\mathcal{X}_k^{\circ})\times(\mathrm{int}(T)\cap\mathcal{Y}_{k'}^{\circ})$, we have $c\tilde{h}_T(x_0,z_0)\leq \tilde{h}_T(x,z)\leq C\tilde{h}_T(x_0,z_0)$, hence
\begin{align*}
   & \int_{\mathrm{int}(T)\cap\mathcal{X}_k^{\circ}}\tilde{f}(x')d\mathcal{H}^1(x')\int_{\mathrm{int}(T)\cap\mathcal{Y}_{k'}^{\circ}}\tilde{g}(z')d\mathcal{H}^1(z')\nonumber\\
   \geq\,&\frac{c \int_{(\mathrm{int}(T)\cap\mathcal{X}_k^{\circ})\times(\mathrm{int}(T)\cap\mathcal{Y}_{k'}^{\circ})}\tilde{f}(x')\tilde{g}(z')\tilde{h}_T(x',z') d\mathcal{H}^1(x')d\mathcal{H}^1(z')}{\tilde{h}_T(x_0,z_0)}\nonumber\\
   =\,&\frac{c \int_{(\mathrm{int}(T)\cap\mathcal{X}_k^{\circ})\times(\mathrm{int}(T)\cap\mathcal{Y}_{k'}^{\circ})}h(x',z') d\mathcal{H}^1(x')d\mathcal{H}^1(z')}{\tilde{h}_T(x_0,z_0)}>0\nonumber\\
   &\Rightarrow\, \int_{\mathrm{int}(T)\cap\mathcal{X}_k^{\circ}}\tilde{f}(x')d\mathcal{H}^1(x'),\int_{\mathrm{int}(T)\cap\mathcal{Y}_{k'}^{\circ}}\tilde{g}(z')d\mathcal{H}^1(z')>0,
\end{align*}
where we note that $h_T(x_0,z_0)>0$ implies $\tilde{h}_T(x_0,z_0)>0$. Moreover, noting that $x\in\mathfrak{T}^*$ and $\mathfrak{F}(x)>0$ (by Remark \ref{positivityF}), we have 
\begin{align}\label{bddnew1}
& \frac{\mathfrak{F}(x)^{-1}h(x,z)}{\int_{(\mathrm{int}(T)\cap\mathcal{X}_k^{\circ})\times(\mathrm{int}(T)\cap\mathcal{Y}_{k'}^{\circ})}h(x',z') d\mathcal{H}^1(x')d\mathcal{H}^1(z')}\nonumber\\
=\,&\frac{\mathfrak{F}(x)^{-1}\tilde{f}(x)\tilde{g}(z)\tilde{h}_T(x,z)}{\int_{(\mathrm{int}(T)\cap\mathcal{X}_k^{\circ})\times(\mathrm{int}(T)\cap\mathcal{Y}_{k'}^{\circ})}\tilde{f}(x')\tilde{g}(z')\tilde{h}_T(x',z')d\mathcal{H}^1(x')d\mathcal{H}^1(z')}\nonumber\\
& \begin{cases}
   \geq  \frac{cf(x)\tilde{g}(z)}{\int_{\mathrm{int}(T)\cap\mathcal{X}_k^{\circ}}\tilde{f}(x')d\mathcal{H}^1(x')\int_{\mathrm{int}(T)\cap\mathcal{Y}_{k'}^{\circ}}\tilde{g}(z')d\mathcal{H}^1(z')},\\
   \leq \frac{Cf(x)\tilde{g}(z)}{\int_{\mathrm{int}(T)\cap\mathcal{X}_k^{\circ}}\tilde{f}(x')d\mathcal{H}^1(x')\int_{\mathrm{int}(T)\cap\mathcal{Y}_{k'}^{\circ}}\tilde{g}(z')d\mathcal{H}^1(z')}.
\end{cases}
\end{align}
Hence, with $\Phi(\cdot)$ as in \eqref{Phi_def}, we have 
\begin{align}\label{Eqnn4.2}
  &  \int_{(\mathrm{int}(T)\cap\mathcal{X}_k^{\circ})\times(\mathrm{int}(T)\cap\mathcal{Y}_{k'}^{\circ})}
h(x,z)\big|\log\big(\mathfrak{F}(x)^{-1}h(x,z)\big)\big|d\mathcal{H}^1(x)d\mathcal{H}^1(z)\nonumber\\
\leq\,& \int_{(\mathrm{int}(T)\cap\mathcal{X}_k^{\circ})\times(\mathrm{int}(T)\cap\mathcal{Y}_{k'}^{\circ})}
h(x,z)\bigg|\log\bigg(\frac{\mathfrak{F}(x)^{-1}h(x,z)}{\int_{(\mathrm{int}(T)\cap\mathcal{X}_k^{\circ})\times(\mathrm{int}(T)\cap\mathcal{Y}_{k'}^{\circ})}h(x',z') d\mathcal{H}^1(x')d\mathcal{H}^1(z')}\bigg)\bigg|d\mathcal{H}^1(x)d\mathcal{H}^1(z)\nonumber\\
&+\bigg|\Phi\bigg(\int_{(\mathrm{int}(T)\cap\mathcal{X}_k^{\circ})\times(\mathrm{int}(T)\cap\mathcal{Y}_{k'}^{\circ})}
h(x,z)d\mathcal{H}^1(x)d\mathcal{H}^1(z)\bigg)\bigg|\nonumber\\
\leq\,& \int_{(\mathrm{int}(T)\cap\mathcal{X}_k^{\circ})\times(\mathrm{int}(T)\cap\mathcal{Y}_{k'}^{\circ})}h(x,z)\bigg|\log\bigg(\frac{f(x)\tilde{g}(z)}{\int_{\mathrm{int}(T)\cap\mathcal{X}_k^{\circ}}\tilde{f}(x')d\mathcal{H}^1(x')\int_{\mathrm{int}(T)\cap\mathcal{Y}_{k'}^{\circ}}\tilde{g}(z')d\mathcal{H}^1(z')}\bigg)\bigg|d\mathcal{H}^1(x)d\mathcal{H}^1(z)\nonumber\\
&+C \int_{(\mathrm{int}(T)\cap\mathcal{X}_k^{\circ})\times(\mathrm{int}(T)\cap\mathcal{Y}_{k'}^{\circ})}h(x,z)d\mathcal{H}^1(x)d\mathcal{H}^1(z)+C\nonumber\\
\leq\,& \int_{(\mathrm{int}(T)\cap\mathcal{X}_k^{\circ})\times(\mathrm{int}(T)\cap\mathcal{Y}_{k'}^{\circ})}h(x,z)|\log(f(x))|d\mathcal{H}^1(x)d\mathcal{H}^1(z)\nonumber\\
&+\int_{(\mathrm{int}(T)\cap\mathcal{X}_k^{\circ})\times(\mathrm{int}(T)\cap\mathcal{Y}_{k'}^{\circ})}h(x,z)|\log(\tilde{g}(z))|d\mathcal{H}^1(x)d\mathcal{H}^1(z)\nonumber\\
&+\bigg|\log\bigg(\int_{\mathrm{int}(T)\cap\mathcal{X}_k^{\circ}}\tilde{f}(x')d\mathcal{H}^1(x')\bigg)\bigg|\int_{(\mathrm{int}(T)\cap\mathcal{X}_k^{\circ})\times(\mathrm{int}(T)\cap\mathcal{Y}_{k'}^{\circ})}h(x,z)d\mathcal{H}^1(x)d\mathcal{H}^1(z)\nonumber\\
&+\bigg|\log\bigg(\int_{\mathrm{int}(T)\cap\mathcal{Y}_{k'}^{\circ}}\tilde{g}(z')d\mathcal{H}^1(z')\bigg)\bigg|\int_{(\mathrm{int}(T)\cap\mathcal{X}_k^{\circ})\times(\mathrm{int}(T)\cap\mathcal{Y}_{k'}^{\circ})}h(x,z)d\mathcal{H}^1(x)d\mathcal{H}^1(z)+C \nonumber\\
\leq\,& \int_{\mathrm{int}(T)\cap\mathcal{X}_k^{\circ}}\tilde{f}(x)|\log(f(x))|d\mathcal{H}^1(x)+\int_{\mathrm{int}(T)\cap\mathcal{Y}_{k'}^{\circ}}\tilde{g}(z)|\log(\tilde{g}(z))|d\mathcal{H}^1(z)+\bigg|\Phi\bigg(\int_{\mathrm{int}(T)\cap\mathcal{X}_k^{\circ}}\tilde{f}(x')d\mathcal{H}^1(x')\bigg)\bigg|\nonumber\\
&+\bigg|\Phi\bigg(\int_{\mathrm{int}(T)\cap\mathcal{Y}_{k'}^{\circ}}\tilde{g}(z')d\mathcal{H}^1(z')\bigg)\bigg|+C\nonumber\\
\leq\,& \int_{\mathrm{int}(T)\cap\mathcal{X}_k^{\circ}}\tilde{f}(x)|\log(f(x))|d\mathcal{H}^1(x)+\int_{\mathrm{int}(T)\cap\mathcal{Y}_{k'}^{\circ}}\tilde{g}(z)|\log(\tilde{g}(z))|d\mathcal{H}^1(z)+C,
\end{align}
where we use \eqref{bddinte} and \eqref{bddnew1} in the second and third inequalities, note the fact that $h_T(\cdot,\cdot)$ is the density of $\kappa_T\in\Pi(\tilde{\mu}_T,\tilde{\nu}_T)$ on $T\times T$ in the fourth inequality, and use the facts that $\int_{\mathrm{int}(T)\cap\mathcal{X}_k^{\circ}}\tilde{f}(x')d\mathcal{H}^1(x')\leq 1$ and $\int_{\mathrm{int}(T)\cap\mathcal{Y}_{k'}^{\circ}}\tilde{g}(z')d\mathcal{H}^1(z')\leq 1$ in the last inequality. By \eqref{Eqnn4.1} and \eqref{Eqnn4.2}, we have 
\begin{align}\label{Eqnn4.3}
    & \int_{\tilde{\mathcal{S}}^*}d\lambda(T)\int_{(\mathrm{int}(T)\cap\mathcal{X}_k^{\circ})\times(\mathrm{int}(T)\cap\mathcal{Y}_{k'}^{\circ})}
h(x,z)\big|\log\big(\mathfrak{F}(x)^{-1}h(x,z)\big)\big|d\mathcal{H}^1(x)d\mathcal{H}^1(z)\nonumber\\
\leq\, &\int_{\mathcal{S}}d\lambda(T)\int_{\mathrm{int}(T)\cap\mathcal{X}_k^{\circ}}\tilde{f}(x)|\log(f(x))|d\mathcal{H}^1(x)+\int_{\mathcal{S}}d\lambda(T)\int_{\mathrm{int}(T)\cap\mathcal{Y}_{k'}^{\circ}}\tilde{g}(z)|\log(\tilde{g}(z))|d\mathcal{H}^1(z)+C\nonumber\\
\leq\,& \int_{\mathfrak{T}^*}f(x)|\log(f(x))|dx+C\leq C,
\end{align}
where we use Lemmas \ref{Lem3.18nn} and \ref{L3.18} in the second inequality and \eqref{bddLip} in the third inequality. 

By \eqref{distcondi}, for any $x\in\mathcal{X}_k$ and $z\in\mathcal{Y}_{k'}$, we have $2d_0\leq\|x-z\|\leq 2(D+1)$, hence
\begin{align}\label{Eqnn4.4}
& \int_{\tilde{\mathcal{S}}^*}
d\lambda(T)\int_{(\mathrm{int}(T)\cap\mathcal{X}_k^{\circ})\times(\mathrm{int}(T)\cap\mathcal{Y}_{k'}^{\circ})}
h(x,z)\big|\log(2\pi\|x-z\|)\big|d\mathcal{H}^1(x)d\mathcal{H}^1(z)\nonumber\\
\leq\, & C\int_{\tilde{\mathcal{S}}^*}
d\lambda(T)\int_{(\mathrm{int}(T)\cap\mathcal{X}_k^{\circ})\times(\mathrm{int}(T)\cap\mathcal{Y}_{k'}^{\circ})}h(x,z)d\mathcal{H}^1(x)d\mathcal{H}^1(z)\leq C\lambda(\tilde{\mathcal{S}}^*)\leq C.
\end{align}
By Lemma \ref{Lem3.11}, we have
\begin{align}\label{Eqnn4.5}
& \int_{\tilde{\mathcal{S}}^*}
d\lambda(T)\int_{(\mathrm{int}(T)\cap\mathcal{X}_k^{\circ})\times(\mathrm{int}(T)\cap\mathcal{Y}_{k'}^{\circ})}
h(x,z)\big|\log\det\big(\mathbf{I}_d+\|x-z\| F(z)\big)\big|d\mathcal{H}^1(x)d\mathcal{H}^1(z)\nonumber\\
\leq\, & \int_{\mathcal{S}}d\lambda(T)\int_{T\times T}h(x,z)\big|\log\det\big(\mathbf{I}_d+\|x-z\| F(z)\big)\big|d\mathcal{H}^1(x)d\mathcal{H}^1(z) \leq C.
\end{align}
By \eqref{Eqnn4.0} and \eqref{Eqnn4.3}--\eqref{Eqnn4.5}, we obtain \eqref{eq3L4.15}.
\end{proof}

The last lemma provides an upper bound on the entropy of $\tilde{\gamma}_{\eps;k,k'}$. 

\begin{lemma}\label{lemma4.16}
For any fixed $\delta_2\in (0,1\slash 100)$, we have 
\begin{align*}
      &\limsup_{\delta_3\rightarrow 0^+} \limsup_{\delta_1\rightarrow 0^+}\limsup_{M_1\rightarrow\infty}\limsup_{\eps\rightarrow 0^{+}}\int_{\mathbb{R}^d\times\mathbb{R}^d}\tilde{\phi}_{\eps;k,k'}(x,y)\log\big(\eps^{(d-1)\slash 2}\tilde{\phi}_{\eps;k,k'}(x,y)\big)dxdy \nonumber\\
    \leq\,& \int_{\mathcal{X}_k^{\circ}\cap\mathfrak{T}^*\cap \tilde{\mathcal{T}}_{1;k,k'}^{*}}dx\int_{\mathrm{int}(T(x))\cap\mathcal{Y}_{k'}^{\circ}}\mathfrak{F}(x)^{-1}h(x,z)\bigg(\log\big(\mathfrak{F}(x)^{-1}h(x,z)\big)-\frac{d-1}{2}\log(2\pi\|x-z\|)\nonumber\\
&\hspace{2.9in}+\frac{1}{2}\log\det\big(\mathbf{I}_d+\|x-z\|F(z)\big)-\frac{d-1}{2}\bigg)d\mathcal{H}^1(z).
\end{align*}
\end{lemma}
\begin{proof}

We fix any $M\geq 1$ and $ \delta\in (0,1)$. With $\tilde{p}_1(\cdot,\cdot)$ as in \eqref{p1til}, we define 
\begin{align}
   & \mathscr{D}_{\delta}:=\big\{(x,y)\in\mathbb{R}^d\times\mathbb{R}^d: r(x,y)\leq \delta(1-\delta_3)\tilde{p}_1(x,y)\big\}, \nonumber\\
   &  \mathscr{D}'_{\delta}:=\big\{(x,y)\in\mathbb{R}^d\times\mathbb{R}^d: r(x,y)> \delta(1-\delta_3)\tilde{p}_1(x,y)\big\}.
\end{align}

For any $(x,y)\in \mathscr{D}'_{\delta}$, we have
\begin{equation*}
    \tilde{\phi}_{\eps;k,k'}(x,y)=(1-\delta_3)p(x,y)+r(x,y)\leq (1-\delta_3)\tilde{p}_1(x,y)+r(x,y)\leq (1+\delta^{-1})r(x,y),
\end{equation*}
hence
\begin{align*}
    & \tilde{\phi}_{\eps;k,k'}(x,y)\log\big(\eps^{(d-1)\slash 2}\tilde{\phi}_{\eps;k,k'}(x,y)\big) \leq \tilde{\phi}_{\eps;k,k'}(x,y)\max\big\{\log\big(\eps^{(d-1)\slash 2}\tilde{\phi}_{\eps;k,k'}(x,y)\big),0\big\}\nonumber\\
    \leq\,& (1+\delta^{-1})r(x,y) \max\big\{ \log\big((1+\delta^{-1})\eps^{(d-1)\slash 2}r(x,y)\big) , 0\big\}\nonumber\\
    \leq\,& (1+\delta^{-1})\log(1+\delta^{-1})r(x,y)+(1+\delta^{-1})r(x,y)\max\big\{\log\big(\eps^{(d-1)\slash 2}r(x,y)\big),0\big\}.
 \end{align*}  
Consequently, we have 
\begin{align}\label{Lem4.16.1}
   &  \int_{\mathscr{D}_{\delta}'}\tilde{\phi}_{\eps;k,k'}(x,y)\log\big(\eps^{(d-1)\slash 2}\tilde{\phi}_{\eps;k,k'}(x,y)\big) dxdy\nonumber\\
    \leq\,& (1+\delta^{-1})\log(1+\delta^{-1})\int_{\mathscr{D}_{\delta}'}r(x,y)dxdy+(1+\delta^{-1})\int_{\mathscr{D}_{\delta}'}r(x,y)\max\big\{ \log\big(\eps^{(d-1)\slash 2} r(x,y)\big),0\big\}dxdy\nonumber\\
    \leq\,& (1+\delta^{-1})\log(1+\delta^{-1})\int_{\mathbb{R}^d\times\mathbb{R}^d}r(x,y)dxdy\nonumber\\
    &+(1+\delta^{-1})\int_{\mathbb{R}^d\times\mathbb{R}^d}r(x,y)\max\big\{\log\big(\eps^{(d-1)\slash 2}r(x,y)\big),0\big\}dxdy.
\end{align}
 
For any $(x,y)\in  \mathscr{D}_{\delta}$, we have  
\begin{equation*}
     \tilde{\phi}_{\eps;k,k'}(x,y)=(1-\delta_3)p(x,y)+r(x,y) \leq  (1+\delta)(1-\delta_3)\tilde{p}_1(x,y), 
\end{equation*}  
hence
\begin{align*}
     &\tilde{\phi}_{\eps;k,k'}(x,y)\log\big(\eps^{(d-1)\slash 2}\tilde{\phi}_{\eps;k,k'}(x,y)\big)\nonumber\\
     =\,& \tilde{\phi}_{\eps;k,k'}(x,y)\big(\log\big(\eps^{(d-1)\slash 2}\tilde{\phi}_{\eps;k,k'}(x,y)\big)\big)_{+}-\tilde{\phi}_{\eps;k,k'}(x,y)\big(\log\big(\eps^{(d-1)\slash 2}\tilde{\phi}_{\eps;k,k'}(x,y)\big)\big)_{-} \nonumber\\
     \leq\,& \tilde{\phi}_{\eps;k,k'}(x,y)\big(\log\big(\eps^{(d-1)\slash 2}\tilde{\phi}_{\eps;k,k'}(x,y)\big)\big)_{+}-\tilde{\phi}_{\eps;k,k'}(x,y)\min\big\{\big(\log\big(\eps^{(d-1)\slash 2}\tilde{\phi}_{\eps;k,k'}(x,y)\big)\big)_{-},M \big\}\nonumber\\
     \leq\,& (1+\delta)(1-\delta_3)\tilde{p}_1(x,y)\big(\log\big( (1+\delta)(1-\delta_3)\eps^{(d-1)\slash 2}\tilde{p}_1(x,y)\big)\big)_{+}\nonumber\\
     &-\tilde{\phi}_{\eps;k,k'}(x,y)\min\big\{\big(\log\big((1+\delta)(1-\delta_3)\eps^{(d-1)\slash 2}\tilde{p}_1(x,y)\big)\big)_{-}, M  \big\}\nonumber\\
     =\,& (1+\delta)(1-\delta_3)\tilde{p}_1(x,y)\max\big\{\log\big( (1+\delta)(1-\delta_3)\eps^{(d-1)\slash 2}\tilde{p}_1(x,y)\big),-M  \big\}\nonumber\\
     &+\big((1+\delta)(1-\delta_3)\tilde{p}_1(x,y)-\tilde{\phi}_{\eps;k,k'}(x,y)\big)\min\big\{\big(\log\big((1+\delta)(1-\delta_3)\eps^{(d-1)\slash 2}\tilde{p}_1(x,y)\big)\big)_{-}, M\big\}\nonumber\\
     \leq\,& (1+\delta)(1-\delta_3)\tilde{p}_1(x,y)\max\big\{\log\big((1+\delta)\eps^{(d-1)\slash 2}\tilde{p}_1(x,y)\big),-M  \big\}\nonumber\\
     &+ M \big((1+\delta)(1-\delta_3)\tilde{p}_1(x,y)-\tilde{\phi}_{\eps;k,k'}(x,y)\big)_{+}\nonumber\\
     \leq\,& (1+\delta)\log(1+\delta)(1-\delta_3)\tilde{p}_1(x,y)+(1+\delta)(1-\delta_3)\tilde{p}_1(x,y)\max\big\{\log\big(\eps^{(d-1)\slash 2}\tilde{p}_1(x,y)\big),-M  \big\}\nonumber\\
     &+ M (1-\delta_3) \big((1+\delta)\tilde{p}_1(x,y)-p(x,y)\big),
\end{align*}
where the last inequality uses the fact that $\tilde{\phi}_{\eps;k,k'}(x,y)\geq (1-\delta_3)p(x,y)$ and $p(x,y)\leq \tilde{p}_1(x,y)\leq (1+\delta)\tilde{p}_1(x,y)$ for any $x,y\in\mathbb{R}^d$. Therefore, we have 
\begin{align*}
& \int_{\mathscr{D}_{\delta}}\tilde{\phi}_{\eps;k,k'}(x,y)\log\big(\eps^{(d-1)\slash 2}\tilde{\phi}_{\eps;k,k'}(x,y)\big)dxdy\nonumber\\
\leq\,&  (1+\delta)\log(1+\delta)(1-\delta_3)\int_{\mathscr{D}_{\delta}}\tilde{p}_1(x,y)dxdy+M(1-\delta_3)\int_{\mathscr{D}_{\delta}}\big((1+\delta)\tilde{p}_1(x,y)-p(x,y)\big)dxdy\nonumber\\
&+(1+\delta)(1-\delta_3)\int_{\mathscr{D}_{\delta}}\tilde{p}_1(x,y)\max\big\{\log\big(\eps^{(d-1)\slash 2}\tilde{p}_1(x,y)\big),-M  \big\}dxdy\nonumber\\
\leq\,&  (1+\delta)\log(1+\delta)(1-\delta_3)+M(1-\delta_3)\int_{\mathbb{R}^d\times\mathbb{R}^d}\big((1+\delta)\tilde{p}_1(x,y)-p(x,y)\big)dxdy\nonumber\\
&+(1+\delta)(1-\delta_3)\int_{\mathscr{D}_{\delta}}\tilde{p}_1(x,y)\max\big\{\log\big(\eps^{(d-1)\slash 2}\tilde{p}_1(x,y)\big),-M  \big\}dxdy\nonumber\\
\leq\,& (1+\delta)\log(1+\delta)+M\bigg( (1+\delta)\mu_{k,k'}(\mathbb{R}^d)-\int_{\mathbb{R}^d\times\mathbb{R}^d}p(x,y)dxdy \bigg)\nonumber\\
&+(1+\delta)(1-\delta_3)\int_{\mathscr{D}_{\delta}}\tilde{p}_1(x,y)\max\big\{\log\big(\eps^{(d-1)\slash 2}\tilde{p}_1(x,y)\big),-M  \big\}dxdy,
\end{align*}
where we use \eqref{eq0L4.15} in the last two inequalities. Moreover, note that
\begin{align*}
    &(1-\delta_3)\int_{\mathscr{D}_{\delta}}\tilde{p}_1(x,y)\max\big\{\log\big(\eps^{(d-1)\slash 2}\tilde{p}_1(x,y)\big),-M  \big\}dxdy\nonumber\\
    =\,& (1-\delta_3)\int_{\mathbb{R}^d\times\mathbb{R}^d}\tilde{p}_1(x,y)\max\big\{\log\big(\eps^{(d-1)\slash 2}\tilde{p}_1(x,y)\big),-M  \big\}dxdy\nonumber\\
    &- (1-\delta_3)\int_{\mathscr{D}_{\delta}'}\tilde{p}_1(x,y)\max\big\{\log\big(\eps^{(d-1)\slash 2}\tilde{p}_1(x,y)\big),-M  \big\}dxdy\nonumber\\
    \leq\,&  (1-\delta_3)\int_{\mathbb{R}^d\times\mathbb{R}^d}\tilde{p}_1(x,y)\max\big\{\log\big(\eps^{(d-1)\slash 2}\tilde{p}_1(x,y)\big),-M  \big\}dxdy+M(1-\delta_3)\int_{\mathscr{D}_{\delta}'}\tilde{p}_1(x,y)dxdy\nonumber\\
    \leq\,& (1-\delta_3)\int_{\mathbb{R}^d\times\mathbb{R}^d}\tilde{p}_1(x,y)\max\big\{\log\big(\eps^{(d-1)\slash 2}\tilde{p}_1(x,y)\big),-M  \big\}dxdy+M\delta^{-1}\int_{\mathscr{D}_{\delta}'}r(x,y)dxdy,
\end{align*}
where we use the definition of $\mathscr{D}_{\delta}'$ in the last inequality. Combining the above two displays, we get
\begin{align}\label{Lem4.16.2}
    & \int_{\mathscr{D}_{\delta}}\tilde{\phi}_{\eps;k,k'}(x,y)\log\big(\eps^{(d-1)\slash 2}\tilde{\phi}_{\eps;k,k'}(x,y)\big)dxdy\nonumber\\
    \leq\,& (1+\delta)(1-\delta_3)\int_{\mathbb{R}^d\times\mathbb{R}^d}\tilde{p}_1(x,y)\max\big\{\log\big(\eps^{(d-1)\slash 2}\tilde{p}_1(x,y)\big),-M  \big\}dxdy\nonumber\\
    & +M(1+\delta^{-1})\int_{\mathbb{R}^d\times\mathbb{R}^d}r(x,y)dxdy+M\bigg((1+\delta)\mu_{k,k'}(\mathbb{R}^d)-\int_{\mathbb{R}^d\times\mathbb{R}^d}p(x,y)dxdy   \bigg)\nonumber\\
    & +(1+\delta)\log(1+\delta).
\end{align}

Note that by Lemma \ref{lemmacoup},
\begin{equation*}
    \int_{\mathbb{R}^d\times\mathbb{R}^d}r(x,y)dxdy=\mu_{k,k'}(\mathbb{R}^d)-(1-\delta_3)\int_{\mathbb{R}^d\times\mathbb{R}^d}p(x,y)dxdy.
\end{equation*}
Moreover, denoting by $\mathfrak{J}_M$ the right-hand side of \eqref{eq1L4.15} (note that $\mathfrak{J}_M$ is independent of $\eps$), by \eqref{eq1L4.15} we have
\begin{equation*}
    \int_{\mathbb{R}^d\times\mathbb{R}^d}\tilde{p}_1(x,y)\max\big\{\log\big(\eps^{(d-1)\slash 2}\tilde{p}_1(x,y)\big),-M  \big\}dxdy=\mathfrak{J}_M.
\end{equation*}
Hence by \eqref{Lem4.16.1} and \eqref{Lem4.16.2}, we have
\begin{align*}
& 
 \int_{\mathbb{R}^d\times\mathbb{R}^d}\tilde{\phi}_{\eps;k,k'}(x,y)\log\big(\eps^{(d-1)\slash 2}\tilde{\phi}_{\eps;k,k'}(x,y)\big)dxdy\nonumber\\
\leq\,& (1+\delta)(1-\delta_3)\mathfrak{J}_M +\big(M+\log(1+\delta^{-1})\big)(1+\delta^{-1})\bigg(\mu_{k,k'}(\mathbb{R}^d)-(1-\delta_3)\int_{\mathbb{R}^d\times\mathbb{R}^d}p(x,y)dxdy\bigg)\nonumber\\
&+(1+\delta^{-1})\int_{\mathbb{R}^d\times\mathbb{R}^d}r(x,y)\max\big\{\log\big(\eps^{(d-1)\slash 2}r(x,y)\big),0\big\}dxdy\nonumber\\
&+M\bigg((1+\delta)\mu_{k,k'}(\mathbb{R}^d)-\int_{\mathbb{R}^d\times\mathbb{R}^d}p(x,y)dxdy \bigg)+(1+\delta)\log(1+\delta) .
\end{align*}
Sequentially taking $\eps\rightarrow 0^+,M_1\rightarrow\infty,\delta_1\rightarrow 0^+,\delta_3\rightarrow 0^+$ in the above display and using Propositions \ref{P4.3} and \ref{P4.4}, we obtain that 
\begin{align*}
    &  \limsup_{\delta_3\rightarrow 0^+} \limsup_{\delta_1\rightarrow 0^+}\limsup_{M_1\rightarrow\infty}\limsup_{\eps\rightarrow 0^{+}}\int_{\mathbb{R}^d\times\mathbb{R}^d}\tilde{\phi}_{\eps;k,k'}(x,y)\log\big(\eps^{(d-1)\slash 2}\tilde{\phi}_{\eps;k,k'}(x,y)\big)dxdy \nonumber\\
   \leq\,& (1+\delta)\mathfrak{J}_M+M\delta+(1+\delta)\log(1+\delta).
\end{align*}
Taking $\delta\rightarrow 0^+$, we get
\begin{equation}\label{Eq4.102}
    \limsup_{\delta_3\rightarrow 0^+} \limsup_{\delta_1\rightarrow 0^+}\limsup_{M_1\rightarrow\infty}\limsup_{\eps\rightarrow 0^{+}}\int_{\mathbb{R}^d\times\mathbb{R}^d}\tilde{\phi}_{\eps;k,k'}(x,y)\log\big(\eps^{(d-1)\slash 2}\tilde{\phi}_{\eps;k,k'}(x,y)\big)dxdy 
   \leq \mathfrak{J}_M.
\end{equation}

For any $x\in\mathcal{X}_k^{\circ}\cap\mathfrak{T}^*\cap \tilde{\mathcal{T}}_{1;k,k'}^{*}$ and $z\in\mathrm{int}(T(x))\cap\mathcal{Y}_{k'}^{\circ}$ such that $\langle x-z, V(x)\rangle >0$ (note that if $\langle x-z, V(x)\rangle \leq 0$, then by \Cref{Lem2}(c), $h(x,z)=0$) and any $w\in O(V(x))$, as $\frac{1}{\|x-z\|}\mathbf{I}_d+F(z)$ is positive definite (by Lemma \ref{L2.7}), we have
\begin{align*}
    &\bigg|\max\bigg\{\log\big(\mathfrak{F}(x)^{-1}h(x,z)\big)-\frac{d-1}{2}\log(2\pi\|x-z\|) +\frac{1}{2}\log\det\big(\mathbf{I}_d+\|x-z\| F(z)\big)\nonumber\\
&\hspace{0.5in}-\frac{1}{2}w^{\top}\Big(\frac{1}{\|x-z\|}\mathbf{I}_d+F(z)\Big)w,-M\bigg\}\bigg|\nonumber\\
\leq\, & \bigg|\log\big(\mathfrak{F}(x)^{-1}h(x,z)\big)-\frac{d-1}{2}\log(2\pi\|x-z\|) +\frac{1}{2}\log\det\big(\mathbf{I}_d+\|x-z\| F(z)\big)\nonumber\\
&\hspace{0.1in}-\frac{1}{2}w^{\top}\Big(\frac{1}{\|x-z\|}\mathbf{I}_d+F(z)\Big)w\bigg|\nonumber\\
\leq\, & \big|\log\big(\mathfrak{F}(x)^{-1}h(x,z)\big)\big|+\frac{d-1}{2}\big|\log(2\pi\|x-z\|)\big|+\frac{1}{2}\big|\log\det\big(\mathbf{I}_d+\|x-z\|F(z)\big)\big|\nonumber\\
&+\frac{1}{2}w^{\top}\Big(\frac{1}{\|x-z\|}\mathbf{I}_d+F(z)\Big)w.
\end{align*}
Note that by \eqref{Lem4.15.n1}, \eqref{Lem4.15.n2}, and \eqref{eq3L4.15}, we have
\begin{align*}
    &\int_{\mathcal{X}_k^{\circ}\cap\mathfrak{T}^*\cap \tilde{\mathcal{T}}_{1;k,k'}^{*}}dx\int_{\mathrm{int}(T(x))\cap\mathcal{Y}_{k'}^{\circ}}\mathfrak{F}(x)^{-1}h(x,z)(2\pi\|x-z\|)^{-(d-1)\slash 2}\sqrt{\det\big(\mathbf{I}_d+\|x-z\|F(z)\big)}d\mathcal{H}^1(z)\nonumber\\
&\hspace{0.2in}\int_{O(V(x))}\bigg(\big|\log\big(\mathfrak{F}(x)^{-1}h(x,z)\big)\big|+\frac{d-1}{2}\big|\log(2\pi\|x-z\|)\big|+\frac{1}{2}\big|\log\det\big(\mathbf{I}_d+\|x-z\|F(z)\big)\big|\nonumber\\
&\hspace{0.9in}+\frac{1}{2}w^{\top}\Big(\frac{1}{\|x-z\|}\mathbf{I}_d+F(z)\Big)w\bigg)e^{-\frac{1}{2}w^{\top}\big(\frac{1}{\|x-z\|}\mathbf{I}_d+F(z)\big)w}d\mathcal{H}^{d-1}(w)\nonumber\\
\leq\,& \int_{\mathcal{X}_k^{\circ}\cap\mathfrak{T}^*\cap \tilde{\mathcal{T}}_{1;k,k'}^{*}}dx\int_{\mathrm{int}(T(x))\cap\mathcal{Y}_{k'}^{\circ}}
\mathfrak{F}(x)^{-1}h(x,z)
\bigg(\big|\log\big(\mathfrak{F}(x)^{-1}h(x,z)\big)\big|+\frac{d-1}{2}\big|\log(2\pi\|x-z\|)\big|\nonumber\\
& \hspace{2.71in}+\frac{1}{2}\big|\log\det\big(\mathbf{I}_d+\|x-z\|F(z)\big)\big|\bigg)d\mathcal{H}^1(z)+\frac{d-1}{2}  \leq C.
\end{align*}
By the above two displays and the dominated convergence theorem, we get
\begin{align*}
&\lim_{M\rightarrow\infty}\mathfrak{J}_M \nonumber\\
=\, &\int_{\mathcal{X}_k^{\circ}\cap\mathfrak{T}^*\cap \tilde{\mathcal{T}}_{1;k,k'}^{*}}dx\int_{\mathrm{int}(T(x))\cap\mathcal{Y}_{k'}^{\circ}}\mathfrak{F}(x)^{-1}h(x,z)(2\pi\|x-z\|)^{-(d-1)\slash 2}\sqrt{\det\big(\mathbf{I}_d+\|x-z\|F(z)\big)}d\mathcal{H}^1(z)\nonumber\\
&\hspace{0.2in}\int_{O(V(x))}\bigg(\log\big(\mathfrak{F}(x)^{-1}h(x,z)\big)-\frac{d-1}{2}\log(2\pi\|x-z\|) +\frac{1}{2}\log\det\big(\mathbf{I}_d+\|x-z\|F(z)\big)\nonumber\\
&\hspace{0.82in}-\frac{1}{2}w^{\top}\Big(\frac{1}{\|x-z\|}\mathbf{I}_d+F(z)\Big)w\bigg) e^{-\frac{1}{2}w^{\top}\big(\frac{1}{\|x-z\|}\mathbf{I}_d+F(z)\big)w}d\mathcal{H}^{d-1}(w)\nonumber\\
=\, & \int_{\mathcal{X}_k^{\circ}\cap\mathfrak{T}^*\cap \tilde{\mathcal{T}}_{1;k,k'}^{*}}dx\int_{\mathrm{int}(T(x))\cap\mathcal{Y}_{k'}^{\circ}}\mathfrak{F}(x)^{-1}h(x,z)\bigg(\log\big(\mathfrak{F}(x)^{-1}h(x,z)\big)-\frac{d-1}{2}\log(2\pi\|x-z\|)\nonumber\\
& \hspace{2.9in}+\frac{1}{2}\log\det\big(\mathbf{I}_d+\|x-z\|F(z)\big)-\frac{d-1}{2} \bigg)d\mathcal{H}^1(z),
\end{align*}
where the last equality uses \eqref{Lem4.15.n1} and \eqref{Lem4.15.n2}. Therefore, taking $M\rightarrow \infty$ in \eqref{Eq4.102}, we conclude that
\begin{align*}
    &\limsup_{\delta_3\rightarrow 0^+} \limsup_{\delta_1\rightarrow 0^+}\limsup_{M_1\rightarrow\infty}\limsup_{\eps\rightarrow 0^{+}}\int_{\mathbb{R}^d\times\mathbb{R}^d}\tilde{\phi}_{\eps;k,k'}(x,y)\log\big(\eps^{(d-1)\slash 2}\tilde{\phi}_{\eps;k,k'}(x,y)\big)dxdy \nonumber\\
    \leq\,& \int_{\mathcal{X}_k^{\circ}\cap\mathfrak{T}^*\cap \tilde{\mathcal{T}}_{1;k,k'}^{*}}dx\int_{\mathrm{int}(T(x))\cap\mathcal{Y}_{k'}^{\circ}}\mathfrak{F}(x)^{-1}h(x,z)\bigg(\log\big(\mathfrak{F}(x)^{-1}h(x,z)\big)-\frac{d-1}{2}\log(2\pi\|x-z\|)\nonumber\\
& \hspace{2.9in}+\frac{1}{2}\log\det\big(\mathbf{I}_d+\|x-z\|F(z)\big)-\frac{d-1}{2}\bigg)d\mathcal{H}^1(z).\qedhere
\end{align*}
\end{proof}

We can now complete the proof of Proposition \ref{P4.2}. 

\begin{proof}[Proof of Proposition \ref{P4.2}]

By \eqref{eq3L4.15}, $\mathtt{S}_{k,k'}(\mu,\nu;\gamma_0)$ is well-defined and finite. By \Cref{lemma4.14,lemma4.16}, and noting \Cref{lemmacoup}, for any $\delta_2\in (0,1\slash 100)$, we have
\begin{align*}
   & \limsup_{\delta_3\rightarrow 0^{+}} \limsup_{\delta_1\rightarrow 0^{+}}\limsup_{M_1\rightarrow\infty}\limsup_{\eps\rightarrow 0^+}\bigg\{\eps^{-1}\int_{\mathcal{X}_k^{\circ}\times\mathcal{Y}_{k'}^{\circ}}\|x-y\|d\tilde{\gamma}_{\eps;k,k'}(x,y)-\eps^{-1}\OT_{k,k'}(\mu,\nu)\nonumber\\
   & \hspace{2.05in}+\int_{\mathcal{X}_k^{\circ}\times\mathcal{Y}_{k'}^{\circ}}\tilde{\phi}_{\eps;k,k'}(x,y)\log\big(\eps^{(d-1)\slash 2}\tilde{\phi}_{\eps;k,k'}(x,y)\big)dxdy\bigg\}\nonumber\\
   \leq\,& \limsup_{\delta_3\rightarrow 0^{+}} \limsup_{\delta_1\rightarrow 0^{+}}\limsup_{M_1\rightarrow\infty}\limsup_{\eps\rightarrow 0^+}\bigg\{\eps^{-1}\int_{\mathcal{X}_k^{\circ}\times\mathcal{Y}_{k'}^{\circ}}\|x-y\|d\tilde{\gamma}_{\eps;k,k'}(x,y)-\eps^{-1}\OT_{k,k'}(\mu,\nu)\bigg\}\nonumber\\
   &+\limsup_{\delta_3\rightarrow 0^{+}} \limsup_{\delta_1\rightarrow 0^{+}}\limsup_{M_1\rightarrow\infty}\limsup_{\eps\rightarrow 0^+}\bigg\{\int_{\mathcal{X}_k^{\circ}\times\mathcal{Y}_{k'}^{\circ}}\tilde{\phi}_{\eps;k,k'}(x,y)\log\big(\eps^{(d-1)\slash 2}\tilde{\phi}_{\eps;k,k'}(x,y)\big)dxdy\bigg\}\nonumber\\
   \leq\,&\mathtt{S}_{k,k'}(\mu,\nu;\gamma_0),
\end{align*}
as desired.
\end{proof}

\subsection{Proof of the mass estimate (Proposition \ref{P4.3})}\label{Sect.4.4}

In this subsection, we give the proof of Proposition \ref{P4.3}. Recall that we have fixed $k\in[K],k'\in[K']$.

\subsubsection{Preparatory results}

We begin by introducing several functions and sets. 

\begin{definition}[The functions $\tilde{\mathscr{F}},\tilde{\mathscr{G}},\mathscr{H}$]\label{deffgh}
For any $q\in \tilde{H}_{k,k'}$ (recall \eqref{deftildehat}) and $t\in\mathbb{R}$, we define
\begin{equation*}
    \tilde{\mathscr{F}}(q,t):=\begin{cases}
        \tilde{f}(q+tV(q)), & \mbox{ if } t\in (-\beta(q),\alpha(q)),\\
        0, & \mbox{ if }t\notin (-\beta(q),\alpha(q)),
    \end{cases}
\end{equation*}
\begin{equation*}
    \tilde{\mathscr{G}}(q,t):=\begin{cases}
        \tilde{g}(q+tV(q)), & \mbox{ if } t\in (-\beta(q),\alpha(q)),\\
        0, & \mbox{ if }t\notin (-\beta(q),\alpha(q)).
    \end{cases}
\end{equation*}
For any $q\in \tilde{H}_{k,k'}$ and $s,t\in\mathbb{R}$, we define
\begin{equation*}
    \mathscr{H}(q,s,t):=\begin{cases}
        h(q+sV(q), q+tV(q)),&\mbox{ if }(s,t)\in (-\beta(q),\alpha(q))^2,\\
        0,&\mbox{ if }(s,t)\notin (-\beta(q),\alpha(q))^2.
    \end{cases}
\end{equation*}
For any $q\in H_{k,k'}\backslash\tilde{H}_{k,k'}$ and $s,t\in\mathbb{R}$, we define $\tilde{\mathscr{F}}(q,t)=\tilde{\mathscr{G}}(q,t)=\mathscr{H}(q,s,t):=0$. Note that by Definition \ref{De3.8} and Lemma \ref{Lem4.2n}, the functions $\tilde{\mathscr{F}},\tilde{\mathscr{G}},\mathscr{H}$ are Borel measurable.   
\end{definition}

\begin{definition}[The sets $\Gamma_0,\Gamma$ and the functions $\mathscr{I}_0,\mathscr{I}_1$]\label{def.4.6}
Recall \cref{Defn4.1.1n}. We define $\Gamma_0$ to be the set of $(q,s,t)\in \mathfrak{H}_{k,k'}\times\mathbb{R}\times\mathbb{R}$ such that $s,t \in(-\beta(q),\alpha(q))$, $s-t\geq 2d_0$, $q+sV(q)\in\mathcal{X}_k^{\circ}$, and $q+tV(q)\in\mathcal{Y}_{k'}^{\circ}$. For $q\in \mathfrak{H}_{k,k'}$, we define   
\begin{equation*}
    \mathscr{I}_0(q):=\int_{T(q)}f(x)\det\big(\mathbf{I}_d+\langle x-q, V(q)\rangle F(q)\big)d\mathcal{H}^1(x). 
\end{equation*}
For $q\in H_{k,k'}\backslash \mathfrak{H}_{k,k'}$, we define $\mathscr{I}_0(q):=0$. For any $(q,s,t)\in\Gamma_0$, we define
\begin{equation*}
    \mathscr{I}_1(q,s,t):=\mathbbm{1}_{q\in\mathfrak{T}^*} (2\pi(s-t))^{-(d-1)\slash 2}\sqrt{\det\big(\mathbf{I}_d+(s-t)F(q+tV(q))\big)} \mathfrak{F}(q+sV(q))^{-1}\mathscr{H}(q,s,t),
\end{equation*}
where we note that $q\in\mathfrak{T}^*$ implies $q+sV(q)\in \mathfrak{T}^*$ and $\mathfrak{F}(q+sV(q))>0$ (see Remark \ref{positivityF}). For any $(q,s,t)\in (H_{k,k'}\times\mathbb{R}\times\mathbb{R})\backslash\Gamma_0$, we define $\mathscr{I}_1(q,s,t):=0$. We define $\Gamma$ to be the set of $(q,s,t)\in\Gamma_0$ such that the functions $\mathscr{I}_0(\cdot)$, $\mathscr{I}_1(\cdot,s,t)$, $F(\cdot)$, $F(\cdot+sV(\cdot))$, $F(\cdot+tV(\cdot))$, $\mathbbm{1}_{\mathfrak{T}^*\cap\mathfrak{H}_{k,k'}}(\cdot)$, and $\mathbbm{1}_{T^{-1}(\tilde{\mathcal{S}}^*)}(\cdot)$ on $H_{k,k'}$ are approximately continuous at $q$. 
\end{definition}

The following lemma shows that $\Gamma$ has full $\mathcal{H}^{d-1}\otimes\mathcal{L}^1\otimes\mathcal{L}^1$-measure in $\Gamma_0$.

\begin{lemma}\label{Lem4.15}
We have $\mathcal{H}^{d-1}\otimes\mathcal{L}^1\otimes\mathcal{L}^1(\Gamma_0\backslash\Gamma)=0$.
\end{lemma}
\begin{proof}

We have
\begin{equation*}
    \mathcal{H}^{d-1}\otimes\mathcal{L}^1\otimes\mathcal{L}^1(\Gamma_0\backslash\Gamma)=\int_{\mathbb{R}\times\mathbb{R}}dsdt\int_{H_{k,k'}}\mathbbm{1}_{\Gamma_0\backslash \Gamma}(q,s,t)d\mathcal{H}^{d-1}(q). 
\end{equation*}
For any fixed $(s,t)\in \mathbb{R}^2$, as the functions $\mathscr{I}_0(\cdot)$, $\mathscr{I}_1(\cdot,s,t)$, $F(\cdot)$, $F(\cdot+sV(\cdot))$, $F(\cdot+tV(\cdot))$, $\mathbbm{1}_{\mathfrak{T}^*\cap\mathfrak{H}_{k,k'}}(\cdot)$, and $\mathbbm{1}_{T^{-1}(\tilde{\mathcal{S}}^*)}(\cdot)$ (on $H_{k,k'}$) are Borel measurable, by \cite[Theorem 1.37]{MR3409135}, these functions are approximately continuous at $\mathcal{H}^{d-1}$-a.e.\ $q\in H_{k,k'}$. Hence for any fixed $(s,t)\in \mathbb{R}^2$, 
\begin{equation*}
    \mathcal{H}^{d-1}(\{q\in H_{k,k'}: (q,s,t)\in\Gamma_0\backslash \Gamma\})=0.
\end{equation*}
Combining the above two displays yields the desired conclusion. 
\end{proof}

We next isolate the region where $f(x)$ or $g(\mathtt{z}(x,y))$ (recall \eqref{defz}) is small.

\begin{definition}[The set $\mathfrak{R}_{\delta}$]\label{def_Rdelta}
For any $\delta>0$, we define $\mathfrak{R}_{\delta}$ to be the set of $(x , y)\in\mathbb{R}^d\times\mathbb{R}^d$ such that $f(x)\leq \delta$ or $g(\mathtt{z}(x,y))\leq \delta$. 
\end{definition}

The following lemma shows that the $p_1$-mass of this region is small.

\begin{lemma}\label{Lem.4.10}
For any $\delta>0$, we have $\int_{\mathfrak{R}_{\delta}}p_1(x,y)dxdy\leq C\delta$.
\end{lemma}

\begin{proof}
By \eqref{p1} and \Cref{Lem3.18nn}, we have
\begin{align*}
&\int_{\mathfrak{R}_{\delta}}p_1(x,y)dxdy\nonumber\\ 
\leq\,& \int_{\mathfrak{T}^*}\mathbbm{1}_{\min\{\alpha(x),\beta(x)\}\geq\delta_1} \mathfrak{F}(x)^{-1}dx\nonumber\\
&\hspace{0.1in}\int_{\mathrm{int}(T(x))}\mathbbm{1}_{\langle x-z,V(x)\rangle>0} \mathbbm{1}_{\min\{\alpha(z),\beta(z)\}\geq\delta_1}(\mathbbm{1}_{f(x)\leq\delta}+\mathbbm{1}_{g(z)\leq\delta})h(x,z)d\mathcal{H}^1(z)\nonumber\\
&\hspace{0.3in}(2\pi\|x-z\|\eps)^{-(d-1)\slash 2}\sqrt{\det\big(\mathbf{I}_d+\|x-z\|F(z)\big)}\int_{O(V(x))}e^{-\frac{1}{2\eps}w^{\top}\big(\frac{1}{\|x-z\|}\mathbf{I}_d+F(z)\big)w}d\mathcal{H}^{d-1}(w)\nonumber\\
=\,& \int_{\mathfrak{T}^*}\mathbbm{1}_{\min\{\alpha(x),\beta(x)\}\geq\delta_1} \mathfrak{F}(x)^{-1}dx\nonumber\\
&\hspace{0.2in}\int_{\mathrm{int}(T(x))}\mathbbm{1}_{\langle x-z,V(x)\rangle>0} \mathbbm{1}_{\min\{\alpha(z),\beta(z)\}\geq\delta_1}(\mathbbm{1}_{f(x)\leq\delta}+\mathbbm{1}_{g(z)\leq\delta})h(x,z)d\mathcal{H}^1(z)\nonumber\\
\leq\,& \int_{\mathcal{S}} d\lambda(T)\int_{T\times T}h(x,z)(\mathbbm{1}_{f(x)\leq\delta}+\mathbbm{1}_{g(z)\leq\delta})d\mathcal{H}^1(x)d\mathcal{H}^1(z)\nonumber\\
\leq\,&\int_{\mathcal{X}}f(x)\mathbbm{1}_{f(x)\leq\delta}dx+\int_{\mathcal{Y}}g(z)\mathbbm{1}_{g(z)\leq\delta}dz\leq \delta(\mathcal{L}^d(\mathcal{X})+\mathcal{L}^d(\mathcal{Y}))\leq  C\delta, 
\end{align*}
where the last line uses \Cref{L3nnn}.
\end{proof}

The next lemma bounds the contribution to the total mass of $p_1(x,z+w)$ coming from the region where $x$ and $z$ lie on the same transport ray and $w\in O(V(x))$ satisfies $\|w\|\geq M\sqrt{\eps}$.

\begin{lemma}\label{Lem.4.11}
For any $M>0$, define
\begin{align}
    \mathfrak{B}_M:=\,&\int_{\mathcal{X}_k^{\circ}\cap\mathfrak{T}^*\cap\tilde{\mathcal{T}}_{1;k,k'}^{*}}dx\int_{\mathrm{int}(T(x))}\mathfrak{F}(x)^{-1}h(x,z)\mathbbm{1}_{\langle x-z,V(x)\rangle>0}(2\pi\|x-z\|)^{-(d-1)\slash 2}d\mathcal{H}^1(z)\nonumber\\
    &\quad\int_{O(V(x))}\sqrt{\det\big(\mathbf{I}_d+\|x-z\|F(z)\big)}e^{-\frac{1}{2}w^{\top}\big(\frac{1}{\|x-z\|}\mathbf{I}_d+F(z)\big)w}\mathbbm{1}_{\|w\|\geq M}d\mathcal{H}^{d-1}(w).
\end{align}
Then for any $M>0$, we have 
\begin{equation*}
    \int_{\mathcal{X}_k^{\circ}\cap\mathfrak{T}^*\cap\tilde{\mathcal{T}}_{1;k,k'}^{*}}dx\int_{T(x)}d\mathcal{H}^1(z)\int_{O(V(x))}p_1(x,z+w)\mathbbm{1}_{\|w\|\geq M\sqrt{\eps}} d\mathcal{H}^{d-1}(w) \leq\mathfrak{B}_M.
\end{equation*}
Moreover, $\lim_{M\rightarrow \infty}\mathfrak{B}_M=0$. 
\end{lemma}

\begin{proof}

By \eqref{p1}, for any $M>0$,
\begin{align*}
    &
\int_{\mathcal{X}_k^{\circ}\cap\mathfrak{T}^*\cap\tilde{\mathcal{T}}_{1;k,k'}^{*}}dx\int_{T(x)}d\mathcal{H}^1(z)\int_{O(V(x))}p_1(x,z+w)\mathbbm{1}_{\|w\|\geq M\sqrt{\eps}}d\mathcal{H}^{d-1}(w)\nonumber\\
    =\,& \int_{\mathcal{X}_k^{\circ}\cap\mathfrak{T}^*\cap\tilde{\mathcal{T}}_{1;k,k'}^{*}}dx\int_{T(x)}d\mathcal{H}^1(z)\int_{O(V(x))}\eps^{(d-1)\slash 2}p_1(x,z+\sqrt{\eps}w)\mathbbm{1}_{\|w\|\geq M} d\mathcal{H}^{d-1}(w)\nonumber\\
    \leq\,& \int_{\mathcal{X}_k^{\circ}\cap\mathfrak{T}^*\cap\tilde{\mathcal{T}}_{1;k,k'}^{*}}dx\int_{\mathrm{int}(T(x))}\mathfrak{F}(x)^{-1}h(x,z)\mathbbm{1}_{\langle x-z,V(x)\rangle>0}(2\pi\|x-z\|)^{-(d-1)\slash 2}d\mathcal{H}^1(z)\nonumber\\
    &\quad\int_{O(V(x))}\sqrt{\det\big(\mathbf{I}_d+\|x-z\|F(z)\big)}e^{-\frac{1}{2}w^{\top}\big(\frac{1}{\|x-z\|}\mathbf{I}_d+F(z)\big)w}\mathbbm{1}_{\|w\|\geq M}d\mathcal{H}^{d-1}(w) = \mathfrak{B}_M.
\end{align*}
Note that
\begin{align*}
&\int_{\mathcal{X}_k^{\circ}\cap\mathfrak{T}^*\cap\tilde{\mathcal{T}}_{1;k,k'}^{*}}dx\int_{\mathrm{int}(T(x))}\mathfrak{F}(x)^{-1}h(x,z)\mathbbm{1}_{\langle x-z,V(x)\rangle>0}(2\pi\|x-z\|)^{-(d-1)\slash 2}d\mathcal{H}^1(z)\nonumber\\
    &\quad\int_{O(V(x))}\sqrt{\det\big(\mathbf{I}_d+\|x-z\|F(z)\big)}e^{-\frac{1}{2}w^{\top}\big(\frac{1}{\|x-z\|}\mathbf{I}_d+F(z)\big)w}d\mathcal{H}^{d-1}(w)\nonumber\\
    \leq \,& \int_{\mathcal{X}_k^{\circ}\cap\mathfrak{T}^*\cap\tilde{\mathcal{T}}_{1;k,k'}^{*}}dx\int_{T(x)}\mathfrak{F}(x)^{-1}h(x,z)d\mathcal{H}^1(z)\leq \int_{\mathcal{X}_k^{\circ}\cap\mathfrak{T}^*\cap\tilde{\mathcal{T}}_{1;k,k'}^{*}}f(x)dx
    \leq  1<\infty.
\end{align*}
As $\lim_{M\rightarrow\infty}\mathbbm{1}_{\|w\|\geq M}=0$ for any fixed $w\in O(V(x))$ (where $x\in \mathcal{X}_k^{\circ}\cap\mathfrak{T}^*\cap\tilde{\mathcal{T}}_{1;k,k'}^*$), by the above display and the dominated convergence theorem, we have $\lim_{M\rightarrow \infty}\mathfrak{B}_M = 0$.
\end{proof}

The next lemma bounds the total mass of $p_1(x,y)$ over the set of $(x,y)$ for which both $x$ and $\mathtt{z}(x,y)$ stay a positive distance away from the ray ends.  

\begin{lemma}\label{Lem4.18}
The value of
\begin{equation*}
     \int_{\mathbb{R}^d\times\mathbb{R}^d} p_1(x,y) \mathbbm{1}_{\min\{\alpha(x),\beta(x)\}\geq 2\delta_1}\mathbbm{1}_{\min\{\alpha(\mathtt{z}(x,y)),\beta(\mathtt{z}(x,y))\}\geq 2\delta_1}dxdy
\end{equation*}
is independent of $\eps$, and
\begin{equation*}
    \lim_{\delta_1\rightarrow 0^+}\lim_{M_1\rightarrow\infty}\int_{\mathbb{R}^d\times\mathbb{R}^d} p_1(x,y) \mathbbm{1}_{\min\{\alpha(x),\beta(x)\}\geq 2\delta_1} \mathbbm{1}_{\min\{\alpha(\mathtt{z}(x,y)),\beta(\mathtt{z}(x,y))\}\geq 2\delta_1}dxdy=\mu_{k,k'}(\mathbb{R}^d).
\end{equation*}
\end{lemma}
\begin{proof}
We have 
\begin{align*}
    &\int_{\mathbb{R}^d\times\mathbb{R}^d} p_1(x,y) \mathbbm{1}_{\min\{\alpha(x),\beta(x)\}\geq 2\delta_1}\mathbbm{1}_{\min\{\alpha(\mathtt{z}(x,y)),\beta(\mathtt{z}(x,y))\}\geq 2\delta_1}dxdy\nonumber\\
    =\,& \int_{\mathcal{X}_k^{\circ}\cap\mathfrak{T}^*\cap\tilde{\mathcal{T}}_{1;k,k'}^{*}}\mathbbm{1}_{\min\{\alpha(x),\beta(x)\}\geq 2\delta_1}dx\int_{\mathrm{int}(T(x))}\mathbbm{1}_{\min\{\alpha(z),\beta(z)\}\geq 2\delta_1}d\mathcal{H}^1(z)\int_{O(V(x))}p_1(x,z+w)d\mathcal{H}^{d-1}(w)\nonumber\\
    =\,& \int_{\mathcal{X}_k^{\circ}\cap\mathfrak{T}^*\cap\tilde{\mathcal{T}}_{1;k,k'}^{*}}\mathbbm{1}_{\min\{\alpha(x),\beta(x)\}\geq 2\delta_1}dx\int_{\mathrm{int}(T(x))\cap\mathcal{Y}_{k'}^{\circ}}\mathbbm{1}_{\min\{\alpha(z),\beta(z)\}\geq 2\delta_1}\mathbbm{1}_{\langle x-z,V(x)\rangle>0} \mathfrak{F}(x)^{-1}h(x,z)
    d\mathcal{H}^1(z)\nonumber\\
    &\hspace{0.1in} 
    (2\pi\|x-z\|\eps)^{-(d-1)\slash 2}\sqrt{\det\big(\mathbf{I}_d+\|x-z\|F(z)\big)}\int_{O(V(x))}   e^{-\frac{1}{2\eps}w^{\top}\big(\frac{1}{\|x-z\|}\mathbf{I}_d+F(z)\big)w}\mathbbm{1}_{\|w\|\leq M_1\sqrt{\eps}} d\mathcal{H}^{d-1}(w)\nonumber\\
    =\,& \int_{\mathcal{X}_k^{\circ}\cap\mathfrak{T}^*\cap\tilde{\mathcal{T}}_{1;k,k'}^{*}}\mathbbm{1}_{\min\{\alpha(x),\beta(x)\}\geq 2\delta_1}dx\int_{\mathrm{int}(T(x))\cap\mathcal{Y}_{k'}^{\circ}}\mathbbm{1}_{\min\{\alpha(z),\beta(z)\}\geq 2\delta_1}\mathbbm{1}_{\langle x-z,V(x)\rangle>0}\mathfrak{F}(x)^{-1}h(x,z)d\mathcal{H}^1(z)  \nonumber\\
    &\hspace{0.1in} 
    (2\pi\|x-z\|)^{-(d-1)\slash 2}\sqrt{\det\big(\mathbf{I}_d+\|x-z\|F(z)\big)} \int_{O(V(x))}e^{-\frac{1}{2} w^{\top}\big(\frac{1}{\|x-z\|}\mathbf{I}_d+F(z)\big)w}\mathbbm{1}_{\|w\|\leq M_1} d\mathcal{H}^{d-1}(w),
\end{align*}
which is independent of $\eps$. By the monotone convergence theorem, we have 
\begin{align}\label{L4.18.eq1}
    & \lim_{M_1\rightarrow\infty}\int_{\mathbb{R}^d\times\mathbb{R}^d} p_1(x,y) \mathbbm{1}_{\min\{\alpha(x),\beta(x)\}\geq 2\delta_1}\mathbbm{1}_{\min\{\alpha(\mathtt{z}(x,y)),\beta(\mathtt{z}(x,y))\}\geq 2\delta_1}dxdy\nonumber\\
  =\,  & \int_{\mathcal{X}_k^{\circ}\cap\mathfrak{T}^*\cap\tilde{\mathcal{T}}_{1;k,k'}^{*}}\mathbbm{1}_{\min\{\alpha(x),\beta(x)\}\geq 2\delta_1}dx\int_{\mathrm{int}(T(x))\cap\mathcal{Y}_{k'}^{\circ}}\mathbbm{1}_{\min\{\alpha(z),\beta(z)\}\geq 2\delta_1}\mathbbm{1}_{\langle x-z,V(x)\rangle>0}\mathfrak{F}(x)^{-1}h(x,z) d\mathcal{H}^1(z)  \nonumber\\
    &\hspace{0.12in}(2\pi\|x-z\|)^{-(d-1)\slash 2}\sqrt{\det\big(\mathbf{I}_d+\|x-z\|F(z)\big)}\int_{O(V(x))}e^{-\frac{1}{2}w^{\top}\big(\frac{1}{\|x-z\|}\mathbf{I}_d+F(z)\big)w} d\mathcal{H}^{d-1}(w) \nonumber\\
 =\,  & \int_{\mathcal{X}_k^{\circ}\cap\mathfrak{T}^*\cap\tilde{\mathcal{T}}_{1;k,k'}^{*}}dx\int_{\mathrm{int}(T(x))\cap\mathcal{Y}_{k'}^{\circ}}\mathbbm{1}_{\min\{\alpha(x),\beta(x)\}\geq 2\delta_1}\mathbbm{1}_{\min\{\alpha(z),\beta(z)\}\geq 2\delta_1}\mathfrak{F}(x)^{-1}h(x,z)d\mathcal{H}^1(z). 
\end{align}
Note that by \eqref{fkkdef}, we have  
\begin{equation*}
    \int_{\mathcal{X}_{k}^{\circ}\cap\mathfrak{T}^*\cap\tilde{\mathcal{T}}_{1;k,k'}^*}dx\int_{\mathrm{int}(T(x))\cap\mathcal{Y}_{k'}^{\circ}}\mathfrak{F}(x)^{-1}h(x,z) d\mathcal{H}^1(z)=\int_{\mathbb{R}^d}f_{k,k'}(x)dx=\mu_{k,k'}(\mathbb{R}^d).
\end{equation*}
Moreover, for any $x\in \mathcal{X}_{k}^{\circ}\cap\mathfrak{T}^*\cap\tilde{\mathcal{T}}_{1;k,k'}^*$ and $z\in\mathrm{int}(T(x))\cap\mathcal{Y}_{k'}^{\circ}$, we have
\begin{equation*}
     \lim\limits_{\delta_1\rightarrow 0^+}\mathbbm{1}_{\min\{\alpha(x),\beta(x)\}\geq 2\delta_1}\mathbbm{1}_{\min\{\alpha(z),\beta(z)\}\geq 2\delta_1}=1.
\end{equation*}
Hence by \eqref{L4.18.eq1} and the dominated convergence theorem, we have 
\begin{align*}
   & \lim_{\delta_1\rightarrow 0^+}\lim_{M_1\rightarrow\infty}\int_{\mathbb{R}^d\times\mathbb{R}^d} p_1(x,y) \mathbbm{1}_{\min\{\alpha(x),\beta(x)\}\geq 2\delta_1} \mathbbm{1}_{\min\{\alpha(\mathtt{z}(x,y)),\beta(\mathtt{z}(x,y))\}\geq 2\delta_1}dxdy\nonumber\\
   =\,& \lim_{\delta_1\rightarrow 0^+}\int_{\mathcal{X}_{k}^{\circ}\cap\mathfrak{T}^*\cap\tilde{\mathcal{T}}_{1;k,k'}^*}dx\int_{\mathrm{int}(T(x))\cap\mathcal{Y}_{k'}^{\circ}}\mathbbm{1}_{\min\{\alpha(x),\beta(x)\}\geq 2\delta_1}\mathbbm{1}_{\min\{\alpha(z),\beta(z)\}\geq 2\delta_1}\mathfrak{F}(x)^{-1}h(x,z) d\mathcal{H}^1(z)\nonumber\\
   =\,&\mu_{k,k'}(\mathbb{R}^d). \qedhere
\end{align*}
\end{proof}

The next two definitions introduce sets encoding the continuity properties of several functions.

\begin{definition}[The set $\mathcal{Q}(q,s,t;\delta,\eta)$]
For any $(q,s,t)\in\Gamma$ and $\delta,\eta>0$, we define $\mathcal{Q}(q,s,t;\delta,\eta)$ to be the set of $q'\in H_{k,k'}$ such that $\|q'-q\|\leq \eta$ and the following conditions hold:
\begin{itemize}
    \item[(a)] $\big|\mathscr{I}_0(q')-\mathscr{I}_0(q)\big|\leq\delta$, $\big|\mathscr{I}_1(q',s,t)-\mathscr{I}_1(q,s,t)\big|\leq\delta$;
    \item[(b)] $\|F(q')-F(q)\|_2\leq\delta$, $\|F(q'+sV(q'))-F(q+sV(q))\|_2\leq\delta$, $\|F(q'+tV(q'))-F(q+tV(q))\|_2\leq\delta$;
    \item[(c)] $\mathbbm{1}_{\mathfrak{T}^*\cap\mathfrak{H}_{k,k'}}(q')=\mathbbm{1}_{\mathfrak{T}^*\cap\mathfrak{H}_{k,k'}}(q)$, $\mathbbm{1}_{T^{-1}(\tilde{\mathcal{S}}^*)}(q')=\mathbbm{1}_{T^{-1}(\tilde{\mathcal{S}}^*)}(q)$.
\end{itemize}
\end{definition}

\begin{definition}[The set $\mathscr{B}(q,s,t;\delta,\eta)$]\label{def4.9}
For any $(q,s,t)\in\Gamma$ and $\delta,\eta>0$, we define $\mathscr{B}(q,s,t;\delta,\eta)$ to be the set of $w\in O(V(q))$ such that $\|w\|\leq\eta$, $q+tV(q)+w\in\tilde{\mathcal{T}}_{1;k,k'}^*$, and the following conditions hold with $q'(w):=\mathfrak{q}_{k,k'}(q+tV(q)+w)$ (recall Definition \ref{defkk}):
\begin{itemize}
    \item[(a)] $\max\big\{|\alpha(q+tV(q)+w)-\alpha(q+tV(q))|,|\beta(q+tV(q)+w)-\beta(q+tV(q))|\big\}\leq \delta$;
    \item[(b)] $\|V(q+tV(q)+w)-V(q+tV(q))-F(q+tV(q))w\|\leq \delta\eta$;
     \item[(c)] $\big|\mathscr{I}_0(q'(w))-\mathscr{I}_0(q)\big|\leq\delta$, $\big|\mathscr{I}_1(q'(w),s,t)-\mathscr{I}_1(q,s,t)\big|\leq\delta$;
    \item[(d)] $\|F(q'(w))-F(q)\|_2\leq\delta$, $\|F(q'(w)+sV(q'(w)))-F(q+sV(q))\|_2\leq\delta$,\\ $\|F(q'(w)+tV(q'(w)))-F(q+tV(q))\|_2\leq\delta$;
    \item[(e)] $\mathbbm{1}_{\mathfrak{T}^*\cap\mathfrak{H}_{k,k'}}(q'(w))=\mathbbm{1}_{\mathfrak{T}^*\cap\mathfrak{H}_{k,k'}}(q)$, $\mathbbm{1}_{T^{-1}(\tilde{\mathcal{S}}^*)}(q'(w))=\mathbbm{1}_{T^{-1}(\tilde{\mathcal{S}}^*)}(q)$.
\end{itemize}
\end{definition}

The following lemma shows that ``most'' points of $H_{k,k'}$ in a small ball around $q$ belong to $\mathcal{Q}(q,s,t;\delta,\eta)$. 

\begin{lemma}\label{Lem4.22}
For any $(q,s,t)\in\Gamma$ and $\delta>0$, we have  
\begin{equation*}
    \lim_{\eta\rightarrow 0^{+}}\frac{\mathcal{H}^{d-1}(\{q' \in H_{k,k'}:\|q'-q\|\leq\eta\}\backslash\mathcal{Q}(q,s,t;\delta,\eta))}{\eta^{d-1}}=0.
\end{equation*}
\end{lemma}
\begin{proof}
This follows from the definition of $\Gamma$ (see Definition \ref{def.4.6}).
\end{proof}

Using \Cref{Lem4.22}, we derive the following result, which shows that ``most'' $w\in O(V(q))$ in a small ball belong to $\mathscr{B}(q,s,t;\delta,\eta)$.

\begin{lemma}\label{Lemma4.23}
For any $(q,s,t)\in\Gamma$ and  $\delta>0$, we have
\begin{equation*}
    \lim_{\eta\rightarrow 0^{+}}\frac{\mathcal{H}^{d-1}\big(\big\{w\in O(V(q)):\|w\|\leq\eta\big\}\backslash\mathscr{B}(q,s,t;\delta,\eta)\big)}{\eta^{d-1}} = 0.
\end{equation*}
\end{lemma}

\begin{proof} 

We denote $x:=q+sV(q)$ and $z:=q+tV(q)$. As $(q,s,t)\in\Gamma$, we have $q\in\mathfrak{H}_{k,k'}$, $s>t$, $x\in\mathcal{T}_1^*\cap\mathcal{X}_k^{\circ}$, and $z\in\mathcal{T}_1^*\cap\mathcal{Y}_{k'}^{\circ}$. Note that this implies $\langle V(q),\mathtt{a}_{k,k'} \rangle\neq 0$ and $\min\{\alpha(q),\beta(q)\}\geq 2d_0$ (recall the setup at the beginning of Section \ref{sectdom}). Without loss of generality, we assume that $\langle V(q),\mathtt{a}_{k,k'} \rangle>0$. As $\mathscr{B}(q,s,t;\delta,\eta)\subseteq \mathscr{B}(q,s,t;\delta',\eta)$ for any $(q,s,t)\in\Gamma$ and  $\delta,\delta',\eta>0$ with $\delta<\delta'$, we assume, without loss of generality, that $\delta\in (0,\min\{\alpha(z),\beta(z),d_0\}\slash 10]$. We also assume that $\eta\in \big(0,\frac{\langle V(q), \mathtt{a}_{k,k'}\rangle\delta^2d_0}{1000D}\big)$.

Below, we consider any $w\in O(V(q))$ such that
\begin{equation}\label{assumpcl}
  \|w\|\leq\eta, \qquad  \max\{|\alpha(z+w)-\alpha(z)|,|\beta(z+w)-\beta(z)|\}\leq \delta.
\end{equation}
As $\min\{\alpha(z+w),\beta(z+w)\}\geq \min\{\alpha(z),\beta(z)\}-\delta\geq\delta>0$, we have $z+w\in \mathcal{T}_1^{*}$. Hence by Lemma~\ref{L2.0}, 
\begin{equation}\label{Lem4.23.eq2}
    \|V(z+w)-V(z)\|
   \leq \frac{4\|w\|}{\min\{\alpha(z+w),\beta(z+w),\alpha(z),\beta(z)\}}\leq \frac{4\eta}{\delta}.
\end{equation}
Note that this implies $|\langle V(z+w),\mathtt{a}_{k,k'}\rangle-\langle V(z),\mathtt{a}_{k,k'}\rangle|\leq \|V(z+w)-V(z)\|\leq 4\eta\slash \delta$. Consequently,
\begin{equation}\label{Lem4.23.eq1}
   \langle V(z+w),\mathtt{a}_{k,k'}\rangle\geq \langle V(z),\mathtt{a}_{k,k'}\rangle-\frac{4\eta}{\delta}=\langle V(q),\mathtt{a}_{k,k'}\rangle-\frac{4\eta}{\delta}\geq  \frac{1}{2}\langle V(q),\mathtt{a}_{k,k'}\rangle>0.
\end{equation}
Note that $L(z+w)$ intersects $H_{k,k'}$ at (recall \eqref{Hkkdef}) 
\begin{align}\label{E4.new.1}
    \tilde{q}(w):=\,&z+w-\frac{\langle z+w,\mathtt{a}_{k,k'}\rangle -\mathtt{b}_{k,k'}}{\langle V(z+w), \mathtt{a}_{k,k'}\rangle }  V(z+w)\nonumber\\
    =\,& z+w- \frac{\langle q+tV(q)+w,\mathtt{a}_{k,k'}\rangle -\mathtt{b}_{k,k'}}{\langle V(z+w), \mathtt{a}_{k,k'}\rangle }  V(z+w)\nonumber\\
    =\,& z+w-\frac{\langle tV(z)+w,\mathtt{a}_{k,k'}\rangle}{\langle V(z+w), \mathtt{a}_{k,k'}\rangle }V(z+w),
\end{align}
where we note that $\langle q, \mathtt{a}_{k,k'}\rangle =\mathtt{b}_{k,k'}$ (as $q\in H_{k,k'}$) in the third equality. By \eqref{Lem4.23.eq2}--\eqref{Lem4.23.eq1},
\begin{align}\label{nbdd}
  & \big|\langle z+w-\tilde{q}(w), V(z+w)\rangle-t\big|=
     \bigg|\frac{\langle tV(z)+w,\mathtt{a}_{k,k'}\rangle}{\langle V(z+w), \mathtt{a}_{k,k'}\rangle}-t\bigg|=\frac{|\langle t(V(z)-V(z+w))+w,\mathtt{a}_{k,k'}\rangle|}{|\langle V(z+w), \mathtt{a}_{k,k'}\rangle|}\nonumber\\
     \leq\,& \frac{2}{\langle V(q),\mathtt{a}_{k,k'}\rangle} \cdot(2D\|V(z+w)-V(z)\|+\|w\|)\leq \frac{20D\eta}{\delta\langle V(q),\mathtt{a}_{k,k'}\rangle}\leq \frac{\delta}{5}.
\end{align}
Hence by \eqref{assumpcl} and the fact that $\min\{\alpha(q),\beta(q)\}\geq 2d_0\geq 20\delta$, we get $\tilde{q}(w)\in\mathrm{int}(T(z+w))$ and
\begin{align*}
    &\max\{|\alpha(\tilde{q}(w))-\alpha(q)|,|\beta(\tilde{q}(w))-\beta(q)|\}\nonumber\\
    \leq\,& \big|\langle z+w-\tilde{q}(w), V(z+w)\rangle-t\big|+ \max\{|\alpha(z+w)-\alpha(z)|,|\beta(z+w)-\beta(z)|\}\leq \frac{6\delta}{5}<d_0.
\end{align*}
As $\min\{\alpha(q),\beta(q)\}\geq 2d_0$, we have $\min\{\alpha(\tilde{q}(w)),\beta(\tilde{q}(w))\}\geq d_0$; by \eqref{Lem4.23.eq1},  $\langle V(\tilde{q}(w)), \mathtt{a}_{k,k'}\rangle=\langle V(z+w), \mathtt{a}_{k,k'}\rangle>0$. Therefore, noting \eqref{deftildehat} and \eqref{deftildeTkk}, we have (recall Definition \ref{defkk})
\begin{equation}\label{zpw}
    \tilde{q}(w)\in \tilde{H}_{k,k'},\qquad z+w\in\tilde{\mathcal{T}}^*_{1;k,k'}, \qquad \mathfrak{q}_{k,k'}(z+w)=\tilde{q}(w).
\end{equation}
Further, by \eqref{E4.new.1}, \eqref{Lem4.23.eq2}, \eqref{nbdd}, and the fact that $|\langle V(q),\mathtt{a}_{k,k'}\rangle|\leq \|V(q)\|\|\mathtt{a}_{k,k'}\|= 1$,
\begin{align}\label{Lem4.23.eq4}
\|\tilde{q}(w)-q\|\leq\,&\|tV(z)- tV(z+w)\|+\bigg\|tV(z+w)-\frac{\langle tV(z)+w,\mathtt{a}_{k,k'}\rangle }{\langle V(z+w),\mathtt{a}_{k,k'}\rangle} V(z+w)\bigg\|+\|w\|\nonumber\\
\leq\,& |t|\|V(z+w)-V(z)\|+\bigg|\frac{\langle tV(z)+w,\mathtt{a}_{k,k'}\rangle }{\langle V(z+w),\mathtt{a}_{k,k'}\rangle}-t\bigg|+\|w\|\nonumber\\   
    \leq\,&2D\cdot\frac{4\eta}{\delta}+\frac{20D\eta}{\delta\langle V(q),\mathtt{a}_{k,k'}\rangle}+\eta\leq \frac{30D\eta}{\delta\langle V(q),\mathtt{a}_{k,k'}\rangle}.
\end{align}
Below we denote $\eta':=\frac{30D\eta}{\delta\langle V(q),\mathtt{a}_{k,k'}\rangle}$. By \eqref{zpw} and \eqref{Lem4.23.eq4}, 
\begin{align}\label{E.n.q}
 &\,\big\{w\in  O(V(q)):\|w\|\leq\eta,\max\{|\alpha(z+w)-\alpha(z)|,|\beta(z+w)-\beta(z)|\}\leq \delta\big\}\nonumber\\
 \subseteq\,& \big\{w\in O(V(q)):\|w\|\leq \eta,\max\{|\alpha(z+w)-\alpha(z)|,|\beta(z+w)-\beta(z)|\}\leq \delta,\nonumber\\
  &\hspace{1.1in} z+w\in \tilde{\mathcal{T}}_{1;k,k'}^*, \mathfrak{q}_{k,k'}(z+w)=\tilde{q}(w)\in \{q'\in \tilde{H}_{k,k'}:\|q'-q\|\leq \eta'\}\big\}.
\end{align}

Let 
\begin{align*}
    \mathscr{D}:=\,&\big\{w\in O(V(q)):\|w\|\leq \eta,\max\{|\alpha(z+w)-\alpha(z)|,|\beta(z+w)-\beta(z)|\}\leq \delta,\nonumber\\
  &\hspace{1.1in}\tilde{q}(w)\in \{q'\in \tilde{H}_{k,k'}:\|q'-q\|\leq \eta'\}\backslash\mathcal{Q}(q,s,t;\delta,\eta')\big\}.
\end{align*}
By \eqref{E.n.q}, recalling Definition \ref{Defd}, we have
\begin{equation}\label{Lem4.23.eq7}
  \big\{w\in O(V(q)):\|w\|\leq\eta\big\}\backslash \mathscr{B}(q,s,t;\delta,\eta)\subseteq\mathcal{D}(z;\eta,\delta)\cup \mathcal{W}(z;\eta,\delta)\cup\mathscr{D}.
\end{equation}
As $q\in \mathfrak{H}_{k,k'}$, by \Cref{L2.2,L2.4} and noting \cref{Defn4.1.1n}, we have 
\begin{equation}\label{Lem4.23.eq8}
    \lim_{\eta\rightarrow 0^+}\frac{\mathcal{H}^{d-1}(\mathcal{D}(z;\eta,\delta))}{\eta^{d-1}}=0, \qquad \lim_{\eta\rightarrow 0^+}\frac{\mathcal{H}^{d-1}(\mathcal{W}(z;\eta,\delta))}{\eta^{d-1}}=0.
\end{equation}

For any $q',q''\in \tilde{H}_{k,k'}$ such that $\max\{\|q'-q\|,\|q''-q\|\}\leq \eta'$, we have
\begin{equation*}
    \min\{\alpha(q'),\beta(q'),\alpha(q''),\beta(q'')\}\geq d_0,
\end{equation*}
hence by Lemma \ref{L2.0}, 
\begin{equation}\label{Lem4.23eq5}
    \|V(q')-V(q'')\|\leq \frac{4\|q'-q''\|}{\min\{\alpha(q'),\beta(q'),\alpha(q''),\beta(q'')\}}\leq\frac{4\|q'-q''\|}{d_0}\leq\frac{8\eta'}{d_0}=\frac{240D\eta}{\delta d_0\langle V(q),\mathtt{a}_{k,k'}\rangle}.
\end{equation}
In particular, for any $q'\in \tilde{H}_{k,k'}$ such that $\|q'-q\|\leq\eta'$, noting that $\eta<\frac{\langle V(q), \mathtt{a}_{k,k'}\rangle\delta^2d_0}{1000D}$, we have
\begin{equation}\label{Lem4.23eq6}
    \langle V(q'), V(q)\rangle=1+\langle V(q')-V(q),V(q)\rangle\geq 1-\|V(q')-V(q)\|\geq 1-\frac{240D\eta}{\delta d_0\langle V(q),\mathtt{a}_{k,k'}\rangle}\geq\frac{1}{2}.
\end{equation}
Now for any $q'\in\tilde{H}_{k,k'}$ such that $\|q'-q\|\leq\eta'$, we define
\begin{equation*}
    \varphi(q'):=q'-z-\frac{\langle q'-z,V(q)\rangle}{\langle V(q'),V(q)\rangle}V(q').
\end{equation*}
Note that $\varphi(q')\in O(V(q))$. By \eqref{Lem4.23eq5} and \eqref{Lem4.23eq6}, for any $q',q''\in \tilde{H}_{k,k'}$ such that $\max\{\|q'-q\|,\|q''-q\|\}\leq\eta'$, we have 
\begin{align*}
   &\bigg|\frac{\langle q'-z,V(q)\rangle}{\langle V(q'),V(q)\rangle}-\frac{\langle q''-z,V(q)\rangle}{\langle V(q''),V(q)\rangle}\bigg|\nonumber\\
   \leq\,&\frac{|\langle q'-z,V(q)\rangle-\langle q''-z,V(q)\rangle||\langle V(q''),V(q)\rangle|+|\langle V(q'')-V(q'),V(q)\rangle||\langle q''-z,V(q)\rangle|}{\langle V(q'),V(q)\rangle 
 \langle V(q''),V(q)\rangle}\nonumber\\
 \leq\,& 4\big(\|q'-q''\|+\|q''-z\|\|V(q'')-V(q')\|\big)
 \leq (4+32D d_0^{-1})\|q'-q''\|,
\end{align*}
where we note that $\|q''-z\|\leq 2D$ (note that $q'',z\in\mathcal{T}_1^*$) in the last inequality. Hence
\begin{align*}
    &\|\varphi(q')-\varphi(q'')\|\nonumber\\
    \leq\,& \|q'-q''\|+\bigg|\frac{\langle q'-z,V(q)\rangle}{\langle V(q'),V(q)\rangle}-\frac{\langle q''-z,V(q)\rangle}{\langle V(q''),V(q)\rangle}\bigg|\|V(q')\|+\|V(q')-V(q'')\|  \bigg|\frac{\langle q''-z,V(q)\rangle}{\langle V(q''),V(q)\rangle}\bigg|\nonumber\\
    \leq\,& (5+32Dd_0^{-1})\|q'-q''\|+2\|q''-z\|\|V(q')-V(q'')\|\leq (5+48Dd_0^{-1})\|q'-q''\|. 
\end{align*}
By Kirszbraun's theorem, there exists a Lipschitz mapping $\tilde{\varphi}:H_{k,k'}\rightarrow  O(V(q))$ such that $\|\tilde{\varphi}\|_{\mathrm{Lip}}\leq 5+48Dd_0^{-1}$ and $\tilde{\varphi}(q')=\varphi(q')$ for every $q'\in \tilde{H}_{k,k'}$ such that $\|q'-q\|\leq \eta'$. 

Note that for any $w\in O(V(q))$ such that $\|w\|\leq \eta$ and $\max\{|\alpha(z+w)-\alpha(z)|,|\beta(z+w)-\beta(z)|\}\leq \delta$, by \eqref{E4.new.1}, we have $\varphi(\tilde{q}(w))=w$. Hence
\begin{align*}
    \mathscr{D} \subseteq\,& \varphi\big(\{q'\in \tilde{H}_{k,k'}:\|q'-q\|\leq \eta'\}\backslash\mathcal{Q}(q,s,t;\delta,\eta')\big)\nonumber\\
    =\,& \tilde{\varphi}\big(\{q'\in \tilde{H}_{k,k'}:\|q'-q\|\leq \eta'\}\backslash\mathcal{Q}(q,s,t;\delta,\eta')\big).
\end{align*}
By \cite[Theorem 2.8]{MR3409135},
\begin{align*}
   &\frac{\mathcal{H}^{d-1}(\mathscr{D})}{\eta^{d-1}}\leq \frac{(5+48Dd_0^{-1})^{d-1}\mathcal{H}^{d-1}(\{q'\in \tilde{H}_{k,k'}:\|q'-q\|\leq \eta'\}\backslash\mathcal{Q}(q,s,t;\delta,\eta'))}{\eta^{d-1}}\nonumber\\
    =\,& (5+48Dd_0^{-1})^{d-1}\bigg(\frac{30D }{\delta\langle V(q),\mathtt{a}_{k,k'}\rangle}\bigg)^{d-1}\cdot \frac{\mathcal{H}^{d-1}(\{q'\in \tilde{H}_{k,k'}:\|q'-q\|\leq \eta'\}\backslash\mathcal{Q}(q,s,t;\delta,\eta'))}{(\eta')^{d-1}}. 
\end{align*}
Hence by Lemma \ref{Lem4.22}, 
\begin{equation}\label{Lem4.23.eq9}
    \lim_{\eta\rightarrow 0^+}\frac{ \mathcal{H}^{d-1}(\mathscr{D})}{\eta^{d-1}}=0. 
\end{equation}
Combining \eqref{Lem4.23.eq7}, \eqref{Lem4.23.eq8}, and \eqref{Lem4.23.eq9}, we conclude that  
\begin{equation*}
    \lim_{\eta\rightarrow 0^+}\frac{\mathcal{H}^{d-1}\big(\big\{w\in O(V(q)):\|w\|\leq\eta\big\}\backslash \mathscr{B}(q,s,t;\delta,\eta)\big)}{\eta^{d-1}} = 0. \qedhere
\end{equation*}    
\end{proof}

\subsubsection{Proof of Proposition \ref{P4.3}}

We fix any $\delta_0 \in (0,1)$, and let $\mathfrak{R}_{\delta_0}$ be defined as in Definition \ref{def_Rdelta}. We also define
\begin{align}\label{Gdeltadef}
    \mathtt{G}_{\delta_1}:=\,&\big\{(x,y)\in\mathbb{R}^d\times \mathbb{R}^d:x\in\mathcal{X}_k^{\circ}\cap\mathfrak{T}^*\cap\tilde{\mathcal{T}}_{1;k,k'}^{*},  \mathtt{z}(x,y)\in\mathrm{int}(T(x))\cap\mathcal{Y}_{k'}^{\circ},\nonumber\\
&\hspace{1.27in}\min\{\alpha(x),\beta(x)\}\geq 2\delta_1,\min\{\alpha(\mathtt{z}(x,y)),\beta(\mathtt{z}(x,y))\}\geq 2\delta_1 \big\}.
\end{align} 
By Lemma \ref{Lem.4.10}, we have
\begin{align}\label{Gdelta_1}
 \int_{\mathtt{G}_{\delta_1}}(p_1(x,y)-p(x,y))dxdy
   \leq\,&\int_{\mathfrak{R}_{\delta_0}}p_1(x,y)dxdy+\int_{\mathtt{G}_{\delta_1}\backslash \mathfrak{R}_{\delta_0}}(p_1(x,y)-p(x,y))dxdy\nonumber\\
   \leq\,& C\delta_0+\sum_{i=2}^4\int_{\mathtt{G}_{\delta_1}\backslash\mathfrak{R}_{\delta_0}}(p_1(x,y)-p_i(x,y))_{+}dxdy,
\end{align}
where we note $p_1(x,y)-p(x,y)\leq \sum_{i=2}^4 (p_1(x,y)-p_i(x,y))_{+}$ for every $x,y\in\mathbb{R}^d$ in the last inequality.

Consider any $i\in \{2,3,4\}$. Define
\begin{equation*}
    \Pi_i:=\int_{\mathtt{G}_{\delta_1}\backslash\mathfrak{R}_{\delta_0}}(p_1(x,y)-p_i(x,y))_{+}\mathbbm{1}_{\|\mathtt{w}(x,y)\|\leq M_1\sqrt{\eps}\slash 2}dxdy.
\end{equation*}
With $\mathfrak{B}_{M_1\slash 2}$ defined as in Lemma \ref{Lem.4.11}, we have 
\begin{align}\label{G_result}
&\int_{\mathtt{G}_{\delta_1}\backslash\mathfrak{R}_{\delta_0}}(p_1(x,y)-p_i(x,y))_{+}dxdy\nonumber\\
\leq\,& \int_{\mathtt{G}_{\delta_1}\backslash\mathfrak{R}_{\delta_0}}p_1(x,y) \mathbbm{1}_{\|\mathtt{w}(x,y)\|>M_1\sqrt{\eps}\slash 2}dxdy+\int_{\mathtt{G}_{\delta_1}\backslash\mathfrak{R}_{\delta_0}}(p_1(x,y)-p_i(x,y))_{+}\mathbbm{1}_{\|\mathtt{w}(x,y)\|\leq M_1\sqrt{\eps}\slash 2}dxdy\nonumber\\
 \leq\,&
\int_{\mathcal{X}_k^{\circ}\cap\mathfrak{T}^*\cap\tilde{\mathcal{T}}_{1;k,k'}^{*}}dx\int_{T(x)}d\mathcal{H}^1(z)\int_{O(V(x))} p_1(x,z+w)\mathbbm{1}_{\|w\| > M_1\sqrt{\eps}\slash 2}d\mathcal{H}^{d-1}(w)+\Pi_i\nonumber\\
\leq\,&\mathfrak{B}_{M_1\slash 2}+\Pi_i,
\end{align}
where we use Lemma \ref{Lem.4.11} in the last inequality. 

Let $\mathfrak{S}_0$ be the set of $T\in\tilde{\mathcal{S}}^*$ such that $T$ intersects $H_{k,k'}$ at some $q\in \tilde{H}_{k,k'}$. Note that $T^{-1}(\mathfrak{S}_0)\subseteq\tilde{\mathcal{T}}_{1;k,k'}^*$ and $\lambda(\mathcal{S}\backslash\tilde{\mathcal{S}}^*)=0$ (recall Definition \ref{defstilde}). Moreover, for any $T\in\tilde{\mathcal{S}}^*\backslash \mathfrak{S}_0$, we have $\mathrm{int}(T)\cap\tilde{\mathcal{T}}_{1;k,k'}^*=\emptyset$. Hence by Lemmas \ref{Lem3.18nn} and \ref{L3.14n}, noting Definition \ref{factors}, for any $i\in\{ 2, 3, 4 \}$, 
\begin{align*}
\Pi_i=\,&\int_{\mathfrak{T}^*}\mathbbm{1}_{x\in \mathcal{X}_k^{\circ}\cap\mathfrak{T}^*\cap\tilde{\mathcal{T}}_{1;k,k'}^{*}}dx\nonumber\\
&\hspace{0.35in}\int_{\mathbb{R}^d}\mathbbm{1}_{(x,y)\in\mathtt{G}_{\delta_1}\backslash\mathfrak{R}_{\delta_0}}(p_1(x,y)-p_i(x,y))_{+}\mathbbm{1}_{\|\mathtt{w}(x,y)\|\leq M_1\sqrt{\eps}\slash 2}dy\nonumber\\
=\,& \int_{\mathfrak{S}_0}d\lambda(T)\int_T\bigg(\int_{\mathbb{R}^d}\mathbbm{1}_{(x,y)\in\mathtt{G}_{\delta_1}\backslash\mathfrak{R}_{\delta_0}}(p_1(x,y)-p_i(x,y))_{+}\mathbbm{1}_{\|\mathtt{w}(x,y)\|\leq M_1\sqrt{\eps}\slash 2}dy\bigg)\nonumber\\
&\hspace{0.9in}\times\mathbbm{1}_{x\in \mathcal{X}_k^{\circ}\cap\mathfrak{T}^*\cap\tilde{\mathcal{T}}_{1;k,k'}^{*}}\mathfrak{F}(x)d\mathcal{H}^1(x)\nonumber\\
\leq\, &\int_{\mathfrak{H}_{k,k'}\cap T^{-1}(\tilde{\mathcal{S}}^*)}|\langle V(q), \mathtt{a}_{k,k'}\rangle|d\mathcal{H}^{d-1}(q)\nonumber\\
&\hspace{0.14in}\int_{T(q)}\bigg(\int_{\mathbb{R}^d}\mathbbm{1}_{(x,y)\in\mathtt{G}_{\delta_1}\backslash\mathfrak{R}_{\delta_0}}(p_1(x,y)-p_i(x,y))_{+}\mathbbm{1}_{\|\mathtt{w}(x,y)\|\leq M_1\sqrt{\eps}\slash 2}dy\bigg)\nonumber\\
&\hspace{0.7in}\times\mathbbm{1}_{x\in \mathcal{X}_k^{\circ}\cap\mathfrak{T}^*\cap\tilde{\mathcal{T}}_{1;k,k'}^{*}}\det\big(\mathbf{I}_d+\langle x-q,V(q)\rangle F(q)\big) d\mathcal{H}^1(x).
\end{align*}
Recalling Definition \ref{def_Rdelta} and \eqref{Gdeltadef}, we deduce
\begin{align*}
\Pi_i\leq\,& \eps^{(d-1)\slash 2}\int_{\mathfrak{H}_{k,k'}\cap T^{-1}(\tilde{\mathcal{S}}^*)}|\langle V(q), \mathtt{a}_{k,k'}\rangle|d\mathcal{H}^{d-1}(q)\nonumber\\
&\hspace{0.6in}\int_{T(q)\times T(q)}\mathbbm{1}_{x\in\mathcal{X}_k^{\circ}\cap\mathfrak{T}^*\cap\tilde{\mathcal{T}}_{1;k,k'}^{*}}\mathbbm{1}_{z\in \mathcal{Y}_{k'}^{\circ}}\mathbbm{1}_{\min\{\alpha(x),\beta(x)\}\geq2\delta_1,\min\{\alpha(z),\beta(z)\}\geq2\delta_1,  \min\{f(x),g(z)\}\geq\delta_0}\nonumber\\
&\hspace{1.2in}\times\det\big(\mathbf{I}_d+\langle x-q,V(q)\rangle F(q)\big)d\mathcal{H}^1(x)d\mathcal{H}^1(z)\nonumber\\
&\hspace{0.7in} \int_{O(V(q))}(p_1(x,z+\sqrt{\eps}w)-p_i(x,z+\sqrt{\eps}w))_{+}\mathbbm{1}_{\|w\|\leq M_1\slash 2}d\mathcal{H}^{d-1}(w).
\end{align*}
Note that by \eqref{p1} and Definition \ref{def.4.6}, for any $(q,s,t)\in (\mathfrak{H}_{k,k'}\times\mathbb{R}\times\mathbb{R})\backslash \Gamma_0$ such that $s,t\in (-\beta(q),\alpha(q))$, we have $p_1(q+sV(q),q+tV(q)+\sqrt{\eps}w)=0$ for every $w\in O(V(q))$ (using $\dist(\mathcal{X}_k,\mathcal{Y}_{k'})\geq 2d_0$). Hence by Lemma \ref{Lem4.15}, we have
\begin{align}\label{defPi1}
    \Pi_i \leq\,& \eps^{(d-1)\slash 2}\int_{\mathfrak{H}_{k,k'}\cap T^{-1}(\tilde{\mathcal{S}}^*)}|\langle V(q), \mathtt{a}_{k,k'}\rangle|d\mathcal{H}^{d-1}(q)\nonumber\\
&\hspace{0.42in}\int_{\mathbb{R}\times \mathbb{R}}\mathbbm{1}_{(q,s,t)\in\Gamma_0}\mathbbm{1}_{q+sV(q)\in\mathcal{X}_k^{\circ}\cap\mathfrak{T}^*\cap\tilde{\mathcal{T}}_{1;k,k'}^{*}}
\mathbbm{1}_{\substack{\min\{\alpha(q+sV(q)),\beta(q+sV(q))\}\geq2\delta_1,\\\min\{\alpha(q+tV(q)),\beta(q+tV(q))\}\geq2\delta_1,\\\min\{f(q+sV(q)),g(q+tV(q))\}\geq\delta_0}}\det(\mathbf{I}_d+s F(q))dsdt\nonumber\\
&\hspace{0.7in}\int_{O(V(q))}(p_1(q+sV(q),q+tV(q)+\sqrt{\eps}w)-p_i(q+sV(q),q+tV(q)+\sqrt{\eps}w))_{+}\nonumber\\
&\hspace{1.2in}\times\mathbbm{1}_{\|w\|\leq M_1\slash 2}d\mathcal{H}^{d-1}(w)\nonumber\\
=\,& \int_{\mathfrak{H}_{k,k'}\cap T^{-1}(\tilde{\mathcal{S}}^*)
}|\langle V(q), \mathtt{a}_{k,k'}\rangle|d\mathcal{H}^{d-1}(q)\nonumber\\
&\hspace{0.3in}\int_{\mathbb{R}\times \mathbb{R}}\mathbbm{1}_{(q,s,t)\in\Gamma}\mathbbm{1}_{q+sV(q)\in\mathcal{X}_k^{\circ}\cap\mathfrak{T}^*\cap\tilde{\mathcal{T}}_{1;k,k'}^{*}}
\mathbbm{1}_{\substack{\min\{\alpha(q+sV(q)),\beta(q+sV(q))\}\geq2\delta_1,\\\min\{\alpha(q+tV(q)),\beta(q+tV(q))\}\geq2\delta_1,\\\min\{f(q+sV(q)),g(q+tV(q))\}\geq\delta_0}}\det(\mathbf{I}_d+s F(q))dsdt\nonumber\\
&\hspace{0.5in} \eps^{(d-1)\slash 2}\int_{O(V(q))}(p_1(q+sV(q),q+tV(q)+\sqrt{\eps}w)-p_i(q+sV(q),q+tV(q)+\sqrt{\eps}w))_{+}\nonumber\\
&\hspace{1.45in}\times\mathbbm{1}_{\|w\|\leq M_1\slash 2} d\mathcal{H}^{d-1}(w).
\end{align}     

Note that for any $(q,s,t)\in\Gamma$, $q+sV(q)\in \mathfrak{T}^*$ implies $q\in\mathfrak{T}^*$. Below, we fix any $(q,s,t)\in\Gamma$ such that $q\in \mathfrak{H}_{k,k'}\cap T^{-1}(\tilde{\mathcal{S}}^*)$, $q\in\mathfrak{T}^*$, and
\begin{equation}\label{alpha_beta_delta1}
    \min\{\alpha(q+sV(q)),\beta(q+sV(q))\}\geq2\delta_1, \quad \min\{\alpha(q+tV(q)),\beta(q+tV(q))\}\geq2\delta_1,
\end{equation}
\begin{equation}\label{fgdeltab}
    \min\{f(q+sV(q)),g(q+tV(q))\}\geq\delta_0   .
\end{equation}
We also fix any $\delta\in (0,\min\{\delta_0,\delta_1\}\slash 10)$, and set $x:=q+sV(q),z:=q+tV(q)$. As $(q,s,t)\in\Gamma$ and $q\in\mathfrak{H}_{k,k'}\cap\mathfrak{T}^*$, we have $x=q+sV(q)\in\mathcal{X}_k^{\circ}\cap\mathfrak{T}^*\cap\tilde{\mathcal{T}}_{1;k,k'}^{*}$ and $z=q+tV(q)\in\mathrm{int}(T(x))\cap\mathcal{Y}_{k'}^{\circ}$. Hence by \eqref{p1}, for any $w\in O(V(q))$, we have (recall Definition \ref{def.4.6})
\begin{equation}\label{p1exp}
   \eps^{(d-1)\slash 2}p_1(x,z+\sqrt{\eps}w)=\mathscr{I}_1(q,s,t)e^{-\frac{1}{2}w^{\top}\big(\frac{1}{s-t} \mathbf{I}_d+F(q+tV(q))\big)w}\mathbbm{1}_{\|w\|\leq M_1}.
\end{equation}

By \eqref{p1exp}, recalling Definition \ref{def.4.6}, we have for any $i\in\{2,3,4\}$, 
\begin{align}\label{conc1}
    &\eps^{(d-1)\slash 2}\int_{O(V(q))}(p_1(q+sV(q),q+tV(q)+\sqrt{\eps}w)-p_i(q+sV(q),q+tV(q)+\sqrt{\eps}w))_{+}\nonumber\\
&\hspace{1in}\times\mathbbm{1}_{\|w\|\leq M_1\slash 2} d\mathcal{H}^{d-1}(w)\nonumber\\
\leq\,& \eps^{(d-1)\slash 2}\int_{\{w\in O(V(q)):\|w\|\leq M_1\slash 2\}}p_1(q+sV(q),q+tV(q)+\sqrt{\eps}w)d\mathcal{H}^{d-1}(w)\nonumber\\
\leq\,& \int_{O(V(q))}\mathbbm{1}_{\|w\|\leq M_1\slash 2}\mathscr{I}_1(q,s,t)e^{-\frac{1}{2}w^{\top}\big(\frac{1}{s-t}\mathbf{I}_d+F(q+tV(q))\big)w}d\mathcal{H}^{d-1}(w). 
\end{align}
As $q+sV(q)\in\mathfrak{T}^*$, by Remark \ref{positivityF}, $\mathfrak{F}(q+sV(q))>0$. Consequently, 
\begin{equation}\label{conc5} 
    \int_{O(V(q))}\mathscr{I}_1(q,s,t)e^{-\frac{1}{2}w^{\top}\big(\frac{1}{s-t}\mathbf{I}_d+F(q+tV(q))\big)w}d\mathcal{H}^{d-1}(w) \leq \mathfrak{F}(q+sV(q))^{-1}h(q+sV(q),q+tV(q))<\infty.
\end{equation}

For each $i\in\{2,3,4\}$, the following lemma provides a lower bound on $p_i(x,z+w)$ for $w\in\mathscr{B}(q,s,t;\delta,M_1\sqrt{\eps}\slash 2)$ (recall Definition \ref{def4.9}). 

\begin{lemma}\label{Lem4.25}
Suppose that $\delta$ is sufficiently small (depending on $d,\mu,\nu$) and $\eps$ is sufficiently small (depending on $d,\mu,\nu$ and $x,z,\delta,\delta_0,\delta_1,M_1$). Then, for any $w\in\mathscr{B}(q,s,t;\delta,M_1\sqrt{\eps}\slash 2)$, we have
\begin{itemize}
    \item[(a)] $\eps^{(d-1)\slash 2}p_2(x,z+w)\geq e^{-C\delta M_1^2\delta_1^{-2}}(\mathscr{I}_1(q,s,t)-\delta)_{+}e^{-\frac{1}{2\eps}w^{\top}\big(\frac{1}{s-t} \mathbf{I}_d+F(z)\big)w} $;
    \item[(b)] $\eps^{(d-1)\slash 2} p_3(x,z+w) 
 \geq  e^{-C\delta M_1^2\delta_1^{-2}}  (1+C\delta\delta_1^{-d})^{-2}\mathscr{I}_1(q,s,t) e^{-\frac{1}{2\eps}w^{\top}\big(\frac{1}{s-t}\mathbf{I}_d+F(z)\big)w}$;
  \item[(c)] $ \eps^{(d-1)\slash 2} p_4(x,z+w) 
 \geq e^{-C\delta M_1^2\delta_1^{-2}}  (1+C\delta\delta_1^{-d})^{-2}(\mathscr{I}_1(q,s,t)-\delta)_{+} e^{-\frac{1}{2\eps}w^{\top}\big(\frac{1}{s-t}\mathbf{I}_d+F(z)\big)w}$.
\end{itemize}
\end{lemma} 

The proof of the lemma is deferred to \cref{se:proof-of-Lem4.25}.

Continuing with the proof of Proposition~\ref{P4.3}, we assume that $\delta$ is sufficiently small (depending on $d,\mu,\nu$) and $\eps$ is sufficiently small (depending on $d,\mu,\nu$ and $x,z,\delta,\delta_0,\delta_1,M_1$), so that the conclusions of Lemma \ref{Lem4.25} hold. Consider any $i\in\{2,3,4\}$. By Lemma \ref{Lem4.25} and \eqref{p1exp}, we have
\begin{align}\label{conclu2}
    &\eps^{(d-1)\slash 2}\int_{O(V(q))}(p_1(q+sV(q),q+tV(q)+\sqrt{\eps}w)-p_i(q+sV(q),q+tV(q)+\sqrt{\eps}w))_{+}\nonumber\\
&\hspace{1.06in}\times\mathbbm{1}_{\|w\|\leq M_1\slash 2}d\mathcal{H}^{d-1}(w)\nonumber\\
\leq\,& \eps^{(d-1)\slash 2}\int_{\substack{\{w\in O(V(q)):\|w\|\leq M_1\slash 2,\\\quad\sqrt{\eps}w\notin \mathscr{B}(q,s,t;\delta,M_1\sqrt{\eps}\slash 2)\}}}p_1(q+sV(q),q+tV(q)+\sqrt{\eps}w)d\mathcal{H}^{d-1}(w)\nonumber\\
&+\int_{\substack{\{w\in O(V(q)):\|w\|\leq M_1\slash 2,\\\quad\sqrt{\eps}w\in \mathscr{B}(q,s,t;\delta,M_1\sqrt{\eps}\slash 2)\}}}
\Big(\mathscr{I}_1(q,s,t)-e^{-C\delta M_1^2\delta_1^{-2}}  (1+C\delta\delta_1^{-d})^{-2}(\mathscr{I}_1(q,s,t)-\delta)_{+}\Big)\nonumber\\
&\hspace{1.8in}\times e^{-\frac{1}{2}w^{\top}\big(\frac{1}{s-t} \mathbf{I}_d+F(q+tV(q))\big)w}d\mathcal{H}^{d-1}(w)\nonumber\\
\leq\,& \int_{\substack{\{w\in O(V(q)):\|w\|\leq M_1\slash 2,\\\quad\sqrt{\eps}w\notin \mathscr{B}(q,s,t;\delta,M_1\sqrt{\eps}\slash 2)\}}}\mathscr{I}_1(q,s,t)e^{-\frac{1}{2}w^{\top}\big(\frac{1}{s-t} \mathbf{I}_d+F(q+tV(q))\big)w}d\mathcal{H}^{d-1}(w)\nonumber\\
&+\int_{\{w\in O(V(q)):\|w\|\leq M_1\slash 2\}}\Big(\mathscr{I}_1(q,s,t)-e^{-C\delta M_1^2\delta_1^{-2}}  (1+C\delta\delta_1^{-d})^{-2}(\mathscr{I}_1(q,s,t)-\delta)_{+}\Big)
\nonumber\\
&\hspace{1.6in}\times e^{-\frac{1}{2}w^{\top}\big(\frac{1}{s-t} \mathbf{I}_d+F(q+tV(q))\big)w}d\mathcal{H}^{d-1}(w).
\end{align}
By \eqref{conc5}, the function $w\mapsto \mathscr{I}_1(q,s,t)e^{-\frac{1}{2}w^{\top}\big(\frac{1}{s-t}\mathbf{I}_d+F(q+tV(q))\big)w}$ is $\mathcal{H}^{d-1}$-integrable on $O(V(q))$. Moreover, by Lemma \ref{Lemma4.23} (note that $(q,s,t)\in \Gamma$), we have
\begin{align*}
   & \lim_{\eps\rightarrow 0^+}\mathcal{H}^{d-1}\big(\big\{w\in O(V(q)):\|w\|\leq M_1\slash 2,\sqrt{\eps}w\notin \mathscr{B}(q,s,t;\delta,M_1\sqrt{\eps}\slash 2)\big\}\big)\nonumber\\
   =\,& \lim_{\eps\rightarrow 0^+}\frac{\mathcal{H}^{d-1}\big(\big\{w\in O(V(q)):\|w\|\leq M_1\sqrt{\eps}\slash 2,w\notin \mathscr{B}(q,s,t;\delta,M_1\sqrt{\eps}\slash 2)\big\}\big)}{\eps^{(d-1)\slash 2}} = 0.
\end{align*}
Hence
\begin{equation}\label{conclu3}
   \lim_{\eps\rightarrow 0^{+}}\int_{\substack{\{w\in O(V(q)):\|w\|\leq M_1\slash 2,\\\quad\sqrt{\eps}w\notin \mathscr{B}(q,s,t;\delta,M_1\sqrt{\eps}\slash 2)\}}}  \mathscr{I}_1(q,s,t)e^{-\frac{1}{2}w^{\top}\big(\frac{1}{s-t}\mathbf{I}_d+F(q+tV(q))\big)w}d\mathcal{H}^{d-1}(w)=0.
\end{equation}
By \eqref{conclu2} and \eqref{conclu3}, we have
\begin{align*}
    &\limsup_{\eps\rightarrow 0^{+}}\bigg\{\eps^{(d-1)\slash 2}\int_{O(V(q))}(p_1(q+sV(q),q+tV(q)+\sqrt{\eps}w)-p_i(q+sV(q),q+tV(q)+\sqrt{\eps}w))_{+}\nonumber\\
&\hspace{1.6in}\times\mathbbm{1}_{\|w\|\leq M_1\slash 2}d\mathcal{H}^{d-1}(w)\bigg\}\nonumber\\
\leq\,&\int_{\{w\in O(V(q)):\|w\|\leq M_1\slash 2\}}\Big(\mathscr{I}_1(q,s,t)-e^{-C\delta M_1^2\delta_1^{-2}}  (1+C\delta\delta_1^{-d})^{-2}(\mathscr{I}_1(q,s,t)-\delta)_{+}\Big)\nonumber\\
&\hspace{1.5in}\times e^{-\frac{1}{2}w^{\top}\big(\frac{1}{s-t} \mathbf{I}_d+F(q+tV(q))\big)w}d\mathcal{H}^{d-1}(w).
\end{align*}
Taking $\delta\rightarrow 0^{+}$ in the above display, using the dominated convergence theorem (note \eqref{conc5} and the fact that $0\leq \mathscr{I}_1(q,s,t)-e^{-C\delta M_1^2\delta_1^{-2}}  (1+C\delta\delta_1^{-d})^{-2}(\mathscr{I}_1(q,s,t)-\delta)_{+}\leq \mathscr{I}_1(q,s,t)$), we get 
\begin{align}\label{conc6}
    &\lim_{\eps\rightarrow 0^{+}}\bigg\{\eps^{(d-1)\slash 2}\int_{O(V(q))}(p_1(q+sV(q),q+tV(q)+\sqrt{\eps}w)-p_i(q+sV(q),q+tV(q)+\sqrt{\eps}w))_{+}\nonumber\\
&\hspace{1.5in}\times\mathbbm{1}_{\|w\|\leq M_1\slash 2}d\mathcal{H}^{d-1}(w)\bigg\}=0.
\end{align}

Note that for any $q\in\mathfrak{H}_{k,k'}$, with $\mathscr{I}_0(q)$ as in Definition \ref{def.4.6}, by Hadamard's inequality (see, e.g., \cite[Lemma 2.5]{ipsen2008perturbation}) and \eqref{bddsalphabeta} we have 
\begin{align*}
    |\mathscr{I}_0(q)|\leq\,& C\sup_{x\in T(q)}\big\{\big|\det\big(\mathbf{I}_d+\langle x-q, V(q)\rangle F(q)\big)\big|\big\}\leq C\sup_{x\in T(q)}\big\{\|\mathbf{I}_d+\langle x-q, V(q)\rangle F(q)\|^d\big\}\nonumber\\
    \leq\,& C\bigg(1+2D\cdot \frac{4}{\min\{\alpha(q),\beta(q)\}}\bigg)^d\leq C\bigg(1+2D\cdot \frac{4}{2d_0}\bigg)^d\leq C,
\end{align*}
where we use \Cref{L2.4} (see \eqref{L2.4.enew3}) in the second line. Hence by Definition \ref{factors}, we have
\begin{align}\label{conc7}
   & \int_{\mathfrak{H}_{k,k'}\cap T^{-1}(\tilde{\mathcal{S}}^*)}|\langle V(q), \mathtt{a}_{k,k'}\rangle|d\mathcal{H}^{d-1}(q)\nonumber\\
&\hspace{0.2in}\int_{\mathbb{R}\times \mathbb{R}}\mathbbm{1}_{(q,s,t)\in\Gamma}\mathbbm{1}_{q+sV(q)\in\mathcal{X}_k^{\circ}\cap\mathfrak{T}^*\cap\tilde{\mathcal{T}}_{1;k,k'}^{*}}
\mathbbm{1}_{\substack{\min\{\alpha(q+sV(q)),\beta(q+sV(q))\}\geq2\delta_1,\\\min\{\alpha(q+tV(q)),\beta(q+tV(q))\}\geq2\delta_1,\\\min\{f(q+sV(q)),g(q+tV(q))\}\geq\delta_0}}\det(\mathbf{I}_d+s F(q))\nonumber\\
&\hspace{0.5in}\times\mathfrak{F}(q+sV(q))^{-1}h(q+sV(q),q+tV(q))dsdt\nonumber\\
\leq\,& \int_{\mathfrak{H}_{k,k'}\cap T^{-1}(\tilde{\mathcal{S}}^*)}|\langle V(q), \mathtt{a}_{k,k'}\rangle|d\mathcal{H}^{d-1}(q)\int_{\mathrm{int}(T(q))\times\mathrm{int}(T(q))}\mathscr{I}_0(q)h(x,z)d\mathcal{H}^1(x)d\mathcal{H}^1(z)\nonumber\\
\leq\,& \int_{\mathfrak{H}_{k,k'}\cap T^{-1}(\tilde{\mathcal{S}}^*)} \mathscr{I}_0(q)d\mathcal{H}^{d-1}(q)\leq C\mathcal{H}^{d-1}(\mathfrak{H}_{k,k'})\leq C\mathcal{H}^{d-1}(\tilde{H}_{k,k'})\leq C, 
\end{align}
where we use \eqref{BddHtildekk} in the last inequality. By \eqref{defPi1}, \eqref{conc1}, \eqref{conc5}, \eqref{conc7}, and \eqref{conc6}, using the dominated convergence theorem, we get for any $i\in\{2,3,4\}$,
\begin{align*}
   & \limsup_{\eps  \rightarrow 0^+}\Pi_i\nonumber\\
   \leq\,& \int_{\mathfrak{H}_{k,k'}\cap T^{-1}(\tilde{\mathcal{S}}^*)}|\langle V(q), \mathtt{a}_{k,k'}\rangle|d\mathcal{H}^{d-1}(q)\nonumber\\
  & \hspace{0.1in} \int_{\mathbb{R}\times \mathbb{R}}\mathbbm{1}_{(q,s,t)\in\Gamma}\mathbbm{1}_{q+sV(q)\in\mathcal{X}_k^{\circ}\cap\mathfrak{T}^*\cap\tilde{\mathcal{T}}_{1;k,k'}^{*}}
\mathbbm{1}_{\substack{\min\{\alpha(q+sV(q)),\beta(q+sV(q))\}\geq2\delta_1,\\\min\{\alpha(q+tV(q)),\beta(q+tV(q))\}\geq2\delta_1,\\\min\{f(q+sV(q)),g(q+tV(q))\}\geq\delta_0}}\det(\mathbf{I}_d+s F(q))dsdt\nonumber\\
&\hspace{0.2in} \lim_{\eps\rightarrow 0^+}\bigg\{\eps^{(d-1)\slash 2}\int_{O(V(q))}(p_1(q+sV(q),q+tV(q)+\sqrt{\eps}w)-p_i(q+sV(q),q+tV(q)+\sqrt{\eps}w))_{+}\nonumber\\
&\hspace{1.6in}\times\mathbbm{1}_{\|w\|\leq M_1\slash 2} d\mathcal{H}^{d-1}(w)\bigg\}=0.
\end{align*}
Hence
\begin{equation}\label{P2limits}
    \lim_{\eps\rightarrow 0^+}\Pi_i=0,\quad\text{ for every }i\in\{2,3,4\}.
\end{equation}

Sequentially taking $\eps\rightarrow 0^+,  M_1\rightarrow\infty$ in \eqref{G_result}, and noting Lemma \ref{Lem.4.11} (with $M=M_1\slash 2$) and \eqref{P2limits}, we get for every $i\in \{2,3,4\}$,
\begin{equation*}
    \limsup_{M_1\rightarrow \infty}\limsup_{\eps\rightarrow  0^{+}}\int_{\mathtt{G}_{\delta_1}\backslash\mathfrak{R}_{\delta_0}}(p_1(x,y)-p_i(x,y))_{+}dxdy=0.
\end{equation*}
Combining this with \eqref{Gdelta_1}, we get
\begin{equation*}
    \limsup_{M_1\rightarrow\infty}\limsup_{\eps\rightarrow 0^+}\int_{\mathtt{G}_{\delta_1}}(p_1(x,y)-p(x,y))dxdy\leq C\delta_0.
\end{equation*}
Taking $\delta_0\rightarrow 0^+$, we obtain that 
\begin{equation*}
   \limsup_{M_1\rightarrow\infty}\limsup_{\eps\rightarrow 0^+}\int_{\mathtt{G}_{\delta_1}}(p_1(x,y)-p(x,y))dxdy=0.
\end{equation*}
By Lemma \ref{Lem4.18} (recall the definition of $\mathtt{G}_{\delta_1}$ as in \eqref{Gdeltadef}, and note that by \eqref{p1}, $p_1(x,y)>0$ implies $x\in\mathcal{X}_k^{\circ}\cap\mathfrak{T}^*\cap\tilde{\mathcal{T}}_{1;k,k'}^{*},  \mathtt{z}(x,y)\in\mathrm{int}(T(x))\cap\mathcal{Y}_{k'}^{\circ}$), 
$\int_{\mathtt{G}_{\delta_1}}p_1(x,y)dxdy$ is independent of $\eps$, and
\begin{equation*}
    \lim\limits_{\delta_1\rightarrow 0^{+}}\lim\limits_{M_1\rightarrow\infty}\int_{\mathtt{G}_{\delta_1}}p_1(x,y)dxdy =\mu_{k,k'}(\mathbb{R}^d).
\end{equation*}
Hence
\begin{align*}
   & \liminf_{M_1\rightarrow\infty}\liminf_{\eps\rightarrow 0^+}\int_{\mathbb{R}^d\times\mathbb{R}^d}   p(x,y)dxdy \geq \liminf_{M_1\rightarrow\infty}\liminf_{\eps\rightarrow 0^+}\int_{\mathtt{G}_{\delta_1}} p(x,y)dxdy \nonumber\\
   \geq\,& \liminf_{M_1\rightarrow\infty}\liminf_{\eps\rightarrow 0^+}\int_{\mathtt{G}_{\delta_1}}p_1(x,y)dxdy-\limsup_{M_1\rightarrow\infty}\limsup_{\eps\rightarrow 0^+}\int_{\mathtt{G}_{\delta_1}}(p_1(x,y)-p(x,y))dxdy\nonumber\\
   =\,& \liminf_{M_1\rightarrow\infty} 
   \int_{\mathtt{G}_{\delta_1}}p_1(x,y)dxdy.
\end{align*}
Taking $\delta_1\rightarrow 0^+$ in the above display, we get
\begin{equation*}
    \liminf_{\delta_1\rightarrow 0^+}\liminf_{M_1\rightarrow\infty} \liminf_{\eps\rightarrow 0^{+}}\int_{\mathbb{R}^d\times\mathbb{R}^d}   p(x,y)dxdy\geq \mu_{k,k'}(\mathbb{R}^d). 
\end{equation*}
By \eqref{eq0L4.15}, we have
\begin{equation*}
     \liminf_{\delta_1\rightarrow 0^+}\liminf_{M_1\rightarrow\infty} \liminf_{\eps\rightarrow 0^{+}}\int_{\mathbb{R}^d\times\mathbb{R}^d}   p(x,y)dxdy\leq  \liminf_{\delta_1\rightarrow 0^+}\liminf_{M_1\rightarrow\infty} \liminf_{\eps\rightarrow 0^{+}}\int_{\mathbb{R}^d\times\mathbb{R}^d}   p_1(x,y)dxdy \leq\mu_{k,k'}(\mathbb{R}^d).
\end{equation*}
Therefore, we conclude that
\begin{equation*}
    \liminf_{\delta_1\rightarrow 0^+}\liminf_{M_1\rightarrow\infty} \liminf_{\eps\rightarrow 0^{+}}\int_{\mathbb{R}^d\times\mathbb{R}^d}   p(x,y)dxdy=\mu_{k,k'}(\mathbb{R}^d),
\end{equation*}
as desired. This completes the proof of Proposition \ref{P4.3}, modulo the proof of Lemma~\ref{Lem4.25} which is given next.

\subsubsection{Proof of Lemma \ref{Lem4.25}}\label{se:proof-of-Lem4.25}

Throughout the proof, we assume that $\delta$ is sufficiently small (depending on $d,\mu,\nu$) and $\eps$ is sufficiently small (depending on $d,\mu,\nu$ and $x,z,\delta,\delta_0,\delta_1,M_1$). We fix any $w\in \mathscr{B}(q,s,t;\delta,M_1\sqrt{\eps}\slash 2)$, and denote $y:=z+w$. Note that $y\in \tilde{\mathcal{T}}_{1;k,k'}^*$ (recall Definition \ref{def4.9}). We define $ q':=\mathfrak{q}_{k,k'}(y)$ (recall Definition \ref{defkk}). By Definition \ref{def4.9}, we have (recall that $q\in\mathfrak{T}^*\cap\mathfrak{H}_{k,k'}$ and $T(q)\in\tilde{\mathcal{S}}^*$)
\begin{equation}\label{alphabdd}
    \max\{|\alpha(y)-\alpha(z)|,|\beta(y)-\beta(z)|\}\leq \delta,\quad \|V(y)-V(z)-F(z)w\|\leq \frac{1}{2}\delta M_1\sqrt{\eps},
\end{equation}
\begin{equation}\label{Fbound}
    \max\{\|F(q')-F(q)\|_2,\|F(q'+sV(q'))-F(q+sV(q))\|_2,\|F(q'+tV(q'))-F(q+tV(q))\|_2\}\leq\delta,
\end{equation}
\begin{equation}\label{mathIbound}
  q'\in\mathfrak{T}^*\cap\mathfrak{H}_{k,k'},\quad T(q')\in\tilde{\mathcal{S}}^*,\quad   \big|\mathscr{I}_0(q')-\mathscr{I}_0(q)\big|\leq\delta,\quad \big|\mathscr{I}_1(q',s,t)-\mathscr{I}_1(q,s,t)\big|\leq\delta.
\end{equation}
By \eqref{alpha_beta_delta1} and \eqref{alphabdd}, recalling that $\delta\in (0,\delta_1\slash 10)$, we have
\begin{equation}\label{alphald}
    \alpha(y)\geq \alpha(z)-\delta\geq \delta_1, \quad \beta(y) \geq\beta(z)-\delta\geq \delta_1. 
\end{equation}
Note that by Lemma \ref{L2.0}, \eqref{alpha_beta_delta1}, and \eqref{alphald},  
\begin{equation}\label{Vdds}
    \|V(y)-V(z)\|\leq \frac{4\|y-z\|}{\min\{\alpha(z),\beta(z),\alpha(y),\beta(y)\}}\leq 4\|w\|\delta_1^{-1}\leq 2M_1\delta_1^{-1}\sqrt{\eps}.
\end{equation}

\medskip

\noindent\textbf{Proof of (a).} We define (recall \eqref{defzp}) 
\begin{equation}\label{zwdef0}
    z_0 :=\mathtt{z}'(x,y)=y+\langle x-y, V(y)\rangle V(y), \quad w_0 :=\mathtt{w}'(x,y)=x-y-\langle x-y,V(y)\rangle V(y), 
\end{equation}
\begin{equation}\label{stdef0}
   s_0:=\langle z_0-q', V(y)\rangle, \quad t_0:=\langle y-q', V(y)\rangle.
\end{equation}
Note that $V(q')=V(y)$, $z_0=q'+s_0V(q')$, $y=q'+t_0 V(q')$, and $\|y-z_0\|=|s_0-t_0|$. As $(q,s,t)\in\Gamma$ and $\min\{f(q+sV(q)),g(q+tV(q))\}\geq\delta_0$, we have
\begin{equation}\label{xzres}
    x\in\mathcal{X}_k^{\circ},\quad z\in\mathcal{Y}_{k'}^{\circ},\quad \min\{f(x),g(z)\}\geq\delta_0.
\end{equation}
As $\|y-z\|=\|w\|\leq M_1\sqrt{\eps}\slash 2$ is sufficiently small, we have $y\in \mathcal{Y}_{k'}^{\circ}$. Combining this with the facts that $q'\in\mathfrak{T}^*$ (see \eqref{mathIbound}) and $y\in\tilde{\mathcal{T}}_{1;k,k'}^*$ implies that $y\in\mathcal{Y}_{k'}^{\circ}\cap\mathfrak{T}^*\cap\tilde{\mathcal{T}}_{1;k,k'}^{*}$. Hence by \eqref{p2} and \eqref{alphald}, we have
\begin{align}\label{pbdds}
   \eps^{(d-1)\slash 2}p_2(x,y)=\,&\mathbbm{1}_{z_0\in\mathrm{int}(T(y))\cap\mathcal{X}_{k}^{\circ}}\mathbbm{1}_{s_0>t_0}\mathbbm{1}_{\min\{\alpha(z_0),\beta(z_0)\}\geq \delta_1}\mathfrak{F}(y)^{-1}h(z_0,y)\nonumber\\
   &\times (2\pi\|y-z_0\|)^{-(d-1)\slash 2} \sqrt{\det\big(\mathbf{I}_d-\|y-z_0\|F(z_0)\big)} e^{-\frac{1}{2\eps}w_0^{\top}\big(\frac{1}{\|y-z_0\|}\mathbf{I}_d-F(z_0)\big)w_0}.
\end{align}

Next, we establish three claims which together will imply~(a).

\begin{claim}\label{Claim.a1}
We have $z_0\in \mathrm{int}(T(y))\cap\mathcal{X}_{k}^{\circ}$, $s_0>t_0$, and $\min\{\alpha(z_0),\beta(z_0)\}\geq\delta_1$.
\end{claim}
\begin{proof}

By \eqref{zwdef0} and \eqref{stdef0}, we have
\begin{equation*}
    s_0-t_0=\langle z_0-y, V(y)\rangle =\langle x-y, V(y)\rangle =\langle x-z, V(y)\rangle - \langle w, V(y)\rangle. 
\end{equation*}
As $w\in O(V(q))$, by the above display, \eqref{bddsalphabeta}, and \eqref{Vdds}, we have (note that $\langle x-z,V(z)\rangle =s-t$)
\begin{align}\label{s0t0}
   & |(s_0-t_0)-(s-t)|=|\langle x-y, V(y)\rangle -(s-t)|\leq\, |\langle x-z, V(y)\rangle - \langle x-z, V(z)\rangle|+|\langle w, V(y)\rangle|\nonumber\\
    \leq\,&\|x-z\|\|V(y)-V(z)\|+\|w\|\leq 2D\cdot 2M_1\delta_1^{-1}\sqrt{\eps}+\frac{1}{2}M_1\sqrt{\eps}\leq CM_1\delta_1^{-1}\sqrt{\eps}.
\end{align}
As $s-t\geq 2d_0$ (recall Definition \ref{def.4.6}) and $\eps$ is sufficiently small, we have  
\begin{equation}\label{s0t0bd}
    s_0-t_0\geq d_0\,\Rightarrow\, s_0>t_0. 
\end{equation}

By \eqref{alphabdd} and \eqref{s0t0}, noting that $\delta\in (0,\delta_1\slash 10)$ and $\eps$ is sufficiently small, we have
\begin{align}\label{boundingalpha1}
   & |\alpha(y)-\langle z_0-y,V(y)\rangle-\alpha(x)|=|\alpha(y)-\langle x-y,V(y)\rangle-\alpha(x)|\nonumber\\
   \leq\,&  |\alpha(y)-\alpha(z)|+|\alpha(z)-\alpha(x)-(s-t)|+|\langle x-y,V(y)\rangle-(s-t)|\nonumber\\
   \leq\,& \delta+CM_1\delta_1^{-1}\sqrt{\eps}\leq \delta_1,
\end{align}
\begin{align}\label{boundingbeta1}
     & |\beta(y)+\langle z_0-y,V(y)\rangle-\beta(x)|=|\beta(y)+\langle x-y,V(y)\rangle-\beta(x)|\nonumber\\
   \leq\,&  |\beta(y)-\beta(z)|+|\beta(z)-\beta(x)+(s-t)|+|\langle x-y,V(y)\rangle-(s-t)|\nonumber\\
   \leq\,& \delta+CM_1\delta_1^{-1}\sqrt{\eps}\leq \delta_1.
\end{align}
By \eqref{alpha_beta_delta1} and the above two displays, we get
\begin{equation*}
    \min\{\alpha(y)-\langle z_0-y,V(y)\rangle, \beta(y)+\langle z_0-y,V(y)\rangle\}\geq \min\{\alpha(x),\beta(x)\}-\delta_1\geq \delta_1>0.
\end{equation*}
Hence $z_0\in\mathrm{int}(T(y))$, which implies $\alpha(z_0)=\alpha(y)-\langle z_0-y,V(y)\rangle$ and $\beta(z_0)=\beta(y)+\langle z_0-y,V(y)\rangle$. Hence $\min\{\alpha(z_0),\beta(z_0)\}\geq\delta_1$. Now note that by \eqref{zwdef0}, \eqref{bddsalphabeta}, \eqref{s0t0}, and \eqref{Vdds}, we have  
\begin{align}\label{bddz}
   \|w_0\|=\,& \|z_0-x\|=\|y-x+\langle x-y, V(y)\rangle V(y)\|\nonumber\\
   \leq\,&\|w\|+\|z-x+(s-t)V(y)\|+\|(\langle x-y,V(y)\rangle-(s-t))V(y)\|\nonumber\\
   \leq\,&\frac{1}{2}M_1\sqrt{\eps}+|s-t|\|V(y)-V(z)\|+|\langle x-y,V(y)\rangle-(s-t)|\nonumber\\
   \leq\,& \frac{1}{2}M_1\sqrt{\eps}+2D \cdot 2M_1\delta_1^{-1}\sqrt{\eps} + CM_1\delta_1^{-1}\sqrt{\eps} \leq CM_1\delta_1^{-1}\sqrt{\eps}. 
\end{align}
As $x\in\mathcal{X}_k^{\circ}$ (recall \eqref{xzres}) and $\eps$ is sufficiently small, we have $z_0\in\mathcal{X}_k^{\circ}$. 
\end{proof}

\begin{claim}\label{Claim.a2}
We have
\begin{equation*}
    \mathfrak{F}(y)^{-1}h(z_0,y)(2\pi\|y-z_0\|)^{-(d-1)\slash 2} \sqrt{\det\big(\mathbf{I}_d-\|y-z_0\|F(z_0)\big)}\geq e^{-CM_1\delta_1^{-2}\delta_0^{-1}\sqrt{\eps}}(\mathscr{I}_1(q,s,t)-\delta)_{+}.
\end{equation*}
\end{claim}

\begin{proof}

Note that by Claim \ref{Claim.a1}, $z_0,y\in\mathrm{int}(T(q'))$. By Definitions \ref{factors} and \ref{def.4.6}, we have 
\begin{align*}
   & \mathfrak{F}(z_0)\sqrt{\det\big(\mathbf{I}_d-\|y-z_0\|F(z_0)\big)}\nonumber\\
   =\,& \frac{\det(\mathbf{I}_d+\langle z_0-q',V(q')\rangle F(q'))}{\int_{T(q')} f(\theta)\det(\mathbf{I}_d+\langle \theta-q',V(q')\rangle F(q'))d\mathcal{H}^1(\theta)}\cdot  \sqrt{\det\big(\mathbf{I}_d +\langle y-z_0, V(q')\rangle F(z_0) \big)}\nonumber\\
   =\,& \mathscr{I}_0(q')^{-1} \sqrt{\det\big(\mathbf{I}_d+\langle z_0-q',V(q')\rangle F(q')\big)}\nonumber\\
   &\times  \sqrt{\det\big(\mathbf{I}_d+\langle z_0-q',V(q')\rangle F(q')\big) \det\big(\mathbf{I}_d +\langle y-z_0, V(q')\rangle F(z_0) \big)}\nonumber\\
   =\,& \mathscr{I}_0(q')^{-1}\sqrt{\det\big(\mathbf{I}_d+\langle z_0-q',V(q')\rangle F(q')\big)}\cdot \sqrt{\det\big(\mathbf{I}_d+\langle y-q',V(q')\rangle F(q')\big)} 
\end{align*}
where we use the fact that $s_0>t_0$ (see Claim \ref{Claim.a1}) in the first equality and use \Cref{L2.4} (see \eqref{Eq2.6.4}; note that $q'\in\mathfrak{T}^*$ by \eqref{mathIbound}) in the third equality. Hence by Definitions \ref{De3.8} and \ref{DefhH},
\begin{align}\label{E4.183}
   & \mathfrak{F}(y)^{-1}h(z_0,y)(2\pi\|y-z_0\|)^{-(d-1)\slash 2} \sqrt{\det\big(\mathbf{I}_d-\|y-z_0\|F(z_0)\big)}\nonumber\\
   =\,& \mathfrak{F}(y)^{-1}\tilde{f}(z_0)\tilde{g}(y)\tilde{h}_{T(q')}(z_0,y)(2\pi\|y-z_0\|)^{-(d-1)\slash 2} \sqrt{\det\big(\mathbf{I}_d-\|y-z_0\|F(z_0)\big)}\nonumber\\
   =\,& f(z_0)g(y)\tilde{h}_{T(q')}(z_0,y)(2\pi\|y-z_0\|)^{-(d-1)\slash 2}\mathfrak{F}(z_0) \sqrt{\det\big(\mathbf{I}_d-\|y-z_0\|F(z_0)\big)}\nonumber\\
   =\,&  \mathscr{I}_0(q')^{-1}f(z_0)g(y)\tilde{h}_{T(q')}(z_0,y)(2\pi\|y-z_0\|)^{-(d-1)\slash 2}\nonumber\\
   &\times \sqrt{\det\big(\mathbf{I}_d+\langle z_0-q',V(q')\rangle F(q')\big)}\cdot \sqrt{\det\big(\mathbf{I}_d+\langle y-q',V(q')\rangle F(q')\big)}.
\end{align}

Note that by Claim \ref{Claim.a1}, $\|y-z_0\|=|s_0-t_0|=s_0-t_0$. Hence by \eqref{s0t0} and the fact that $s-t\geq 2d_0$,
\begin{equation*}
    \|y-z_0\|\leq s-t+ CM_1\delta_1^{-1}\sqrt{\eps}\leq (1+CM_1\delta_1^{-1}\sqrt{\eps})(s-t),
\end{equation*}
which implies 
\begin{equation}\label{E4.n.1}
    (2\pi\|y-z_0\|)^{-(d-1)\slash 2}\geq (1+CM_1\delta_1^{-1}\sqrt{\eps})^{-(d-1)\slash 2}(2\pi(s-t))^{-(d-1)\slash 2}. 
\end{equation}

As $V(z)=\frac{x-z}{\|x-z\|}$, by \eqref{xzres}, \eqref{bddakk}, and \eqref{bddsalphabeta}, we have $\langle V(z),\mathtt{a}_{k,k'}\rangle = \frac{\langle x-z,\mathtt{a}_{k,k'}\rangle}{\|x-z\|}\geq \frac{4d_0}{2D} \geq c$. As $\eps$ is sufficiently small, by \eqref{Vdds}, $\langle V(y),\mathtt{a}_{k,k'}\rangle\geq c$. Consequently, noting that $q\in H_{k,k'}$, we have
\begin{align*}
   & \bigg|t-\frac{\langle y, \mathtt{a}_{k,k'}\rangle-\mathtt{b}_{k,k'}}{\langle V(y), \mathtt{a}_{k,k'}\rangle}\bigg|=\bigg|t-\frac{\langle q+tV(q)+w, \mathtt{a}_{k,k'}\rangle-\mathtt{b}_{k,k'}}{\langle V(y), \mathtt{a}_{k,k'}\rangle}\bigg|=\bigg|t-\frac{\langle tV(q)+w, \mathtt{a}_{k,k'}\rangle}{\langle V(y), \mathtt{a}_{k,k'}\rangle}\bigg|\nonumber\\
   =\,& \frac{|\langle t(V(y)-V(z))-w, \mathtt{a}_{k,k'}\rangle|}{\langle V(y),\mathtt{a}_{k,k'}\rangle}\leq   C(|t|\|V(y)-V(z)\|+\|w\|)\nonumber\\
   \leq\,& C\bigg(2D\cdot 2M_1\delta_1^{-1}\sqrt{\eps}+\frac{1}{2}M_1\sqrt{\eps}\bigg)\leq CM_1\delta_1^{-1} \sqrt{\eps}, 
\end{align*}
where we use \eqref{Vdds} in the last line. By Definition \ref{defkk} and \eqref{Vdds}, we thus have
\begin{align}\label{qbd}
   & \|q'-q\| = \|\mathfrak{q}_{k,k'}(y)-q\|=\bigg\|y-\frac{\langle y, \mathtt{a}_{k,k'}\rangle-\mathtt{b}_{k,k'}}{\langle V(y), \mathtt{a}_{k,k'}\rangle}V(y)-q\bigg\|\nonumber\\
   \leq\,& \|w\|+\|tV(q)-tV(y)\|+\bigg\|tV(y)-\frac{\langle y, \mathtt{a}_{k,k'}\rangle-\mathtt{b}_{k,k'}}{\langle V(y), \mathtt{a}_{k,k'}\rangle}V(y)\bigg\|\nonumber\\
   \leq\,&\frac{1}{2}M_1\sqrt{\eps}+|t|\|V(z)-V(y)\|+\bigg|t-\frac{\langle y, \mathtt{a}_{k,k'}\rangle-\mathtt{b}_{k,k'}}{\langle V(y), \mathtt{a}_{k,k'}\rangle}\bigg|\leq CM_1\delta_1^{-1}\sqrt{\eps}.
\end{align}
By \eqref{stdef0}, \eqref{Vdds}, and \eqref{qbd}, we have  
\begin{align}\label{tbdd}
    |t_0-t|=\,&|\langle y-q',V(y)\rangle -\langle z-q,V(z)\rangle |\nonumber\\
    \leq\,&|\langle z+w-q', V(y)\rangle-\langle z-q,V(y)\rangle|+|\langle z-q, V(y)-V(z)\rangle|\nonumber\\
    \leq\,& \|w\|+\|q'-q\|+\|z-q\|\|V(y)-V(z)\|\nonumber\\
    \leq\,& \frac{1}{2}M_1\sqrt{\eps}+CM_1\delta_1^{-1}\sqrt{\eps}+2D\cdot 2M_1\delta_1^{-1}\sqrt{\eps}\leq  CM_1\delta_1^{-1}\sqrt{\eps}.
\end{align}
By \eqref{s0t0} and \eqref{tbdd}, we get
\begin{equation}\label{sbdd}
    |s_0-s|\leq CM_1\delta_1^{-1}\sqrt{\eps}.
\end{equation}
By \eqref{alphald} and Claim \ref{Claim.a1}, noting that $\eps$ is sufficiently small, we have 
\begin{equation}\label{intTq}
    q'+s V(q'), q'+t  V(q')\in\mathrm{int}(T(q')),
\end{equation}
\begin{equation}\label{ldalpha}
    \min\{\alpha(q'+s V(q')),\beta(q'+s V(q'))\}  \geq \frac{\delta_1}{2},\quad \min\{\alpha(q'+t V(q')),\beta(q'+t V(q'))\}    \geq \frac{\delta_1}{2}. 
\end{equation}
By \eqref{bddz}, \eqref{tbdd}, and \eqref{sbdd}, 
\begin{equation}\label{bddn1}
    \|q'+sV(q')-x\|\leq \|q'+sV(q')-z_0\|+\|z_0-x\|=|s-s_0|+\|z_0-x\|\leq CM_1\delta_1^{-1}\sqrt{\eps},
\end{equation}
\begin{equation}\label{bddn0}
    \|q'+tV(q')-z\|\leq \|q'+tV(q')-y\|+\|w\|\leq |t-t_0|+\frac{1}{2}M_1\sqrt{\eps}\leq CM_1\delta_1^{-1}\sqrt{\eps}.
\end{equation}
By \eqref{xzres}, noting that $\eps$ is sufficiently small, we get
\begin{equation}\label{xkykk}
    q'+sV(q')\in\mathcal{X}_k^{\circ},\quad q'+tV(q')\in\mathcal{Y}_{k'}^{\circ}.
\end{equation}
As $T(q')\in\tilde{\mathcal{S}}^*$ (see \eqref{mathIbound}),
by \Cref{Lem2} and \eqref{tbdd}--\eqref{sbdd},  
\begin{align}\label{E4.n.2}
    \tilde{h}_{T(q')}(z_0,y)=\,& \tilde{h}_{T(q')}(q'+s_0 V(q'),q'+t_0 V(q'))\geq\,  e^{-C(|s_0-s|+|t_0-t|)}\tilde{h}_{T(q')}(q'+s V(q'),q'+t V(q'))\nonumber\\
    \geq\,& e^{-CM_1\delta_1^{-1}\sqrt{\eps}}\tilde{h}_{T(q')}(q'+sV(q'),q'+tV(q')).
\end{align}

By \eqref{bddLip} and \eqref{bddz},
\begin{equation*}
    |f(z_0)-f(x)|\leq M_0\|z_0-x\|\leq CM_1\delta_1^{-1}\sqrt{\eps}, \quad |g(y)-g(z)|\leq M_0\|w\|\leq CM_1\sqrt{\eps}.
\end{equation*}
As $\eps$ is sufficiently small, noting \eqref{xzres}, we get $\min\{f(z_0),g(y)\}\geq \delta_0\slash 2$. Hence by \eqref{bddz} and \eqref{bddn1}--\eqref{bddn0}, 
\begin{equation*}
    |f(q'+sV(q'))-f(z_0)|\leq M_0\|q'+sV(q')-z_0\|\leq CM_1\delta_1^{-1}\sqrt{\eps}\leq  CM_1\delta_1^{-1}\delta_0^{-1}\sqrt{\eps}f(z_0),
\end{equation*}
\begin{equation*}
    |g(q'+tV(q'))-g(y)|\leq M_0\|q'+tV(q')-y\|\leq CM_1\delta_1^{-1}\sqrt{\eps}\leq CM_1\delta_1^{-1}\delta_0^{-1}\sqrt{\eps}g(y),
\end{equation*}
which yields
\begin{equation}\label{E4.n.3}
    f(z_0)\geq (1+CM_1\delta_1^{-1}\delta_0^{-1}\sqrt{\eps})^{-1}f(q'+sV(q')), \quad g(y)\geq (1+CM_1\delta_1^{-1}\delta_0^{-1}\sqrt{\eps})^{-1} g(q'+tV(q')).
\end{equation}

For any $\theta,\theta'\in\mathrm{int}(T(q'))$, by \Cref{L2.4} (see \eqref{Eq2.6.4}; with $x,x',x_0$ replaced by $\theta,q',\theta'$), 
\begin{equation}
    \det\big(\mathbf{I}_d+\langle \theta'-q', V(q')\rangle F(q')\big)=  \det\big(\mathbf{I}_d+\langle \theta-q', V(q')\rangle F(q')\big)\det\big(\mathbf{I}_d+\langle \theta'-\theta, V(q')\rangle F(\theta)\big).
\end{equation}
By Lemma \ref{L2.7}, we have $ \det\big(\mathbf{I}_d+\langle \theta'-q', V(q')\rangle F(q')\big)>0$, $\det\big(\mathbf{I}_d+\langle \theta-q', V(q')\rangle F(q')\big)>0$, and $\det\big(\mathbf{I}_d+\langle \theta'-\theta, V(q')\rangle F(\theta)\big)>0$. Now by Hadamard's inequality (see, e.g., \cite[Lemma 2.5]{ipsen2008perturbation}), \eqref{bddsalphabeta}, and \Cref{L2.4} (see \eqref{L2.4.enew3}), 
\begin{align*}
    \big|\det\big(\mathbf{I}_d+\langle \theta'-\theta, V(q')\rangle F(\theta)\big)\big|
  \leq\,& \|\mathbf{I}_d+\langle \theta'-\theta, V(q')\rangle F(\theta)\|_2^d  \nonumber\\
  \leq\,&(1+\|\theta'-\theta\| \|F(\theta)\|_2)^d\leq\bigg(1+\frac{C\|\theta-\theta'\|}{\min\{\alpha(\theta),\beta(\theta)\}}\bigg)^d
\end{align*}
Hence
\begin{equation}\label{E4.195}
    \det\big(\mathbf{I}_d+\langle \theta-q', V(q')\rangle F(q')\big)\geq \bigg(1+\frac{C\|\theta-\theta'\|}{\min\{\alpha(\theta),\beta(\theta)\}}\bigg)^{-d} \det\big(\mathbf{I}_d+\langle \theta'-q', V(q')\rangle F(q')\big).
\end{equation}
Taking $\theta=z_0,\theta'=q'+sV(q')$ in \eqref{E4.195} and noting \eqref{intTq}, Claim \ref{Claim.a1} and \eqref{sbdd}, we obtain
\begin{align}\label{E4.n.4}
    \det\big(\mathbf{I}_d+\langle z_0-q', V(q')\rangle F(q')\big)\geq\,& \bigg(1+\frac{C|s_0-s|}{\delta_1}\bigg)^{-d}\det(\mathbf{I}_d+sF(q'))\nonumber\\
    \geq\,& (1+CM_1\delta_1^{-2}\sqrt{\eps})^{-d}\det(\mathbf{I}_d+sF(q')).
\end{align}
Taking $\theta=y,\theta'=q'+tV(q')$ in \eqref{E4.195}, and noting \eqref{intTq}, \eqref{alphald}, and \eqref{tbdd}, we obtain
\begin{align}\label{E4.n.5}
    \det\big(\mathbf{I}_d+\langle y-q', V(q')\rangle F(q')\big)\geq\,& \bigg(1+\frac{C|t_0-t|}{\delta_1}\bigg)^{-d}\det(\mathbf{I}_d+tF(q'))\nonumber\\
    \geq\,& (1+CM_1\delta_1^{-2}\sqrt{\eps})^{-d}\det(\mathbf{I}_d+tF(q')).
\end{align}

By \eqref{E4.183}, \eqref{E4.n.1}, \eqref{E4.n.2}, \eqref{E4.n.3}, \eqref{E4.n.4}, and \eqref{E4.n.5}, we have 
\begin{align*}
    & \mathfrak{F}(y)^{-1}h(z_0,y)(2\pi\|y-z_0\|)^{-(d-1)\slash 2} \sqrt{\det\big(\mathbf{I}_d-\|y-z_0\|F(z_0)\big)}\nonumber\\
    \geq\,&(1+CM_1\delta_1^{-1}\delta_0^{-1}\sqrt{\eps})^{-2}e^{-CM_1\delta_1^{-1}\sqrt{\eps}} (1+CM_1\delta_1^{-1}\sqrt{\eps})^{-(d-1)\slash 2}(1+CM_1\delta_1^{-2}\sqrt{\eps})^{-d}\nonumber\\
    &\times\mathscr{I}_0(q')^{-1}f(q'+sV(q'))g(q'+tV(q'))\tilde{h}_{T(q')}(q'+sV(q'),q'+tV(q'))\nonumber\\
    &\times(2\pi(s-t))^{-(d-1)\slash 2}\sqrt{\det(\mathbf{I}_d+sF(q'))}\cdot \sqrt{\det(\mathbf{I}_d+t F(q'))}\nonumber\\
    \geq\,& e^{-CM_1\delta_1^{-2}\delta_0^{-1}\sqrt{\eps}}\mathfrak{F}(q'+sV(q'))^{-1} h(q'+sV(q'),q'+tV(q'))  \nonumber\\
    & \times\mathscr{I}_0(q')^{-1}\mathfrak{F}(q'+tV(q'))^{-1}(2\pi(s-t))^{-(d-1)\slash 2}\sqrt{\det(\mathbf{I}_d+sF(q'))}\cdot \sqrt{\det(\mathbf{I}_d+t F(q'))},
\end{align*}
where we use Definitions \ref{De3.8} and \ref{DefhH} in the second inequality. By \eqref{intTq}, Definitions \ref{factors} and \ref{def.4.6}, and \Cref{L2.4} (see \eqref{Eq2.6.4}, with $x,x',x_0$ replaced by $q'+tV(q'), q', q'+sV(q')$), we have 
\begin{align}\label{rel1}
&\mathscr{I}_0(q')^{-1}\mathfrak{F}(q'+tV(q'))^{-1}\sqrt{\det(\mathbf{I}_d+sF(q'))}\cdot \sqrt{\det(\mathbf{I}_d+t F(q'))}\nonumber\\
=\, & \sqrt{\frac{\det(\mathbf{I}_d+sF(q'))}{\det(\mathbf{I}_d+t F(q'))}}=\sqrt{\det\big(\mathbf{I}_d+(s-t)F(q'+tV(q'))\big)}.
\end{align}
Combining the above two displays and using Definition \ref{def.4.6} (note that by \eqref{mathIbound}, \eqref{intTq}, and \eqref{xkykk}, $(q',s,t)\in\Gamma_0$), we have
\begin{align*}
    & \mathfrak{F}(y)^{-1}h(z_0,y)(2\pi\|y-z_0\|)^{-(d-1)\slash 2} \sqrt{\det\big(\mathbf{I}_d-\|y-z_0\|F(z_0)\big)}\nonumber\\
    \geq\, & e^{-CM_1\delta_1^{-2}\delta_0^{-1}\sqrt{\eps}}(2\pi(s-t))^{-(d-1)\slash 2}\sqrt{\det\big(\mathbf{I}_d+(s-t)F(q'+tV(q'))\big)}\nonumber\\
    &\times\mathfrak{F}(q'+sV(q'))^{-1} h(q'+sV(q'),q'+tV(q'))
    =  e^{-CM_1\delta_1^{-2}\delta_0^{-1}\sqrt{\eps}}  \mathscr{I}_1(q',s,t). 
\end{align*}
As $w\in \mathscr{B}(q,s,t;\delta,M_1\sqrt{\eps}\slash 2)$ (recall Definition \ref{def4.9}), we have
\begin{equation}\label{rel1.1}
    \big|\mathscr{I}_1(q',s,t)-\mathscr{I}_1(q,s,t)\big|\leq\delta.
\end{equation}
Hence
\begin{equation*}
    \mathfrak{F}(y)^{-1}h(z_0,y)(2\pi\|y-z_0\|)^{-(d-1)\slash 2} \sqrt{\det\big(\mathbf{I}_d-\|y-z_0\|F(z_0)\big)}\geq e^{-CM_1\delta_1^{-2}\delta_0^{-1}\sqrt{\eps}}(\mathscr{I}_1(q,s,t)-\delta)_{+}.
\end{equation*}
\end{proof}

\begin{claim}\label{Claim.a3}
We have 
\begin{equation}\label{Claim.a3.eq1}
    \bigg| w_0^{\top}\bigg(\frac{1}{\|y-z_0\|}\mathbf{I}_d-F(z_0)\bigg)w_0-w^{\top}\bigg(\frac{1}{s-t} \mathbf{I}_d+F(z)\bigg)w  \bigg|\leq C\delta M_1^2\delta_1^{-2}\eps.
\end{equation}
\end{claim}

\begin{proof}

As $w\in \mathscr{B}(q,s,t;\delta,M_1\sqrt{\eps}\slash 2)$ (recall Definition \ref{def4.9}), we have
\begin{equation}\label{BdVy}
    \|V(y)-V(z)-F(z)w\|\leq \frac{1}{2}\delta M_1\sqrt{\eps}.
\end{equation}
As $w\in O(V(q))$ and $F(z)^{\top}V(z)=0$ (see \eqref{Eq2.6.5} in \Cref{L2.4}), we have
\begin{equation*}
    \langle x-y, V(z)+F(z)w\rangle=\langle (s-t)V(z)-w, V(z)+F(z)w\rangle=s-t-w^{\top}F(z)w.
\end{equation*}
By \Cref{L2.4} (see \eqref{L2.4.enew3}) and \eqref{alpha_beta_delta1},
\begin{equation}\label{Fbdds}
  \|F(x)\|_2\leq \frac{4}{\min\{\alpha(x),\beta(x)\}}\leq 2\delta_1^{-1},\quad  \|F(z)\|_2\leq \frac{4}{\min\{\alpha(z),\beta(z)\}}\leq 2\delta_1^{-1}.
\end{equation}
By the above three displays, noting that $\eps$ is sufficiently small, we get 
\begin{align}\label{xyV}
   & |\langle x-y, V(y)\rangle -(s-t)|\leq |\langle x-y, V(y)-V(z)-F(z) w \rangle|+|\langle x-y, V(z)+F(z)w\rangle -(s-t)|\nonumber\\
   \leq\,& (\|x-z\|+\|w\|)\|V(y)-V(z)-F(z) w\|+|w^{\top} F(z) w| \nonumber\\
   \leq\, & \Big(2D+\frac{1}{2}M_1\sqrt{\eps} \Big)\cdot \frac{1}{2}\delta M_1\sqrt{\eps}+\|F(z)\|_2\|w\|^2 \leq C \delta M_1\sqrt{\eps}+\delta_1^{-1} M_1^2 \eps \leq C\delta M_1\sqrt{\eps}.
\end{align}
Consequently, by \eqref{zwdef0} and \eqref{BdVy}, we have
\begin{align}\label{E4.199}
   & \|w_0+(\mathbf{I}_d+(s-t)F(z))w\|=\|x-y-\langle x-y,V(y)\rangle V(y)+(\mathbf{I}_d+(s-t)F(z))w\|\nonumber\\
   =\,& \|x-z-\langle x-y,V(y)\rangle V(y)+(s-t)F(z)w\|\nonumber\\
   \leq\,& \|x-z-(s-t)V(y)+(s-t)F(z)w\|+|\langle x-y, V(y)\rangle -(s-t)|\nonumber\\
   \leq\,& \|x-z-(s-t)(V(z)+F(z)w)+(s-t)F(z)w\|+|s-t|\|V(y)-V(z)-F(z)w\|+C\delta M_1\sqrt{\eps}\nonumber\\
   \leq\,& 2D\cdot \frac{1}{2}\delta M_1\sqrt{\eps}+C\delta M_1\sqrt{\eps}\leq C\delta M_1\sqrt{\eps}.
\end{align}
By \eqref{E4.199} and \eqref{Fbdds},
\begin{align}
   & \max\{\|w_0\|,\|(\mathbf{I}_d+(s-t)F(z))w\|\}\leq \|(\mathbf{I}_d+(s-t)F(z))w\|+C\delta M_1\sqrt{\eps}\nonumber\\
   \leq\,& (1+|s-t|\|F(z)\|_2)\|w\|+C\delta M_1\sqrt{\eps}\leq (1+2D\cdot 2\delta_1^{-1})\cdot\frac{1}{2}M_1\sqrt{\eps}+C\delta M_1\sqrt{\eps} \leq C\delta_1^{-1}M_1\sqrt{\eps}.
\end{align}

As $s-t\geq 2d_0$ and $s_0-t_0\geq d_0$ (recall \eqref{s0t0bd}), by \eqref{tbdd} and \eqref{sbdd}, we have  
\begin{equation}
    \bigg\|\frac{1}{s_0-t_0}\mathbf{I}_d-\frac{1}{s-t}\mathbf{I}_d\bigg\|_2=\bigg|\frac{1}{s_0-t_0}-\frac{1}{s-t}\bigg|\leq\frac{|s_0-s|+|t_0-t|}{(s_0-t_0)(s-t)}\leq CM_1\delta_1^{-1}\sqrt{\eps},
\end{equation}
\begin{equation}\label{stid}
    \bigg\|\frac{1}{s-t}\mathbf{I}_d-F(x)\bigg\|_2\leq \frac{1}{s-t}+\|F(x)\|_2\leq C\delta_1^{-1},
\end{equation}
where we use \eqref{Fbdds} in the second display. By \Cref{L2.4} (see \eqref{L2.4.enew4}), noting Claim \ref{Claim.a1} and \eqref{sbdd}--\eqref{ldalpha}, we get 
\begin{align*}
\|F(z_0)-F(q'+sV(q'))\|_2\leq\,& \frac{128\|z_0-(q'+sV(q'))\|}{\min\{\alpha(z_0),\beta(z_0)\}\min\{\alpha(q'+sV(q')),\beta(q'+sV(q'))\}}\nonumber\\
\leq\,&  C\delta_1^{-2}|s_0-s|\leq CM_1\delta_1^{-3}\sqrt{\eps}.
\end{align*}
As $w\in \mathscr{B}(q,s,t;\delta,M_1\sqrt{\eps}\slash 2)$ (recall Definition \ref{def4.9}), we have
\begin{equation*}
    \|F(q'+sV(q')) -F(x)\|_2=\|F(q'+sV(q')) -F(q+sV(q))\|_2\leq \delta.
\end{equation*}
Combining the above two displays, noting that $\eps$ is sufficiently small, we get
\begin{equation}\label{E4.202}
    \|F(z_0)-F(x)\|_2 \leq CM_1\delta_1^{-3}\sqrt{\eps}+\delta\leq 2\delta.
\end{equation}

By Remark \ref{Rmk.L2.4} (with $x,x'$ replaced by $z=q+tV(q), x=q+sV(q)$), we have
\begin{equation*}
    (\mathbf{I}_d-(s-t)F(x))(\mathbf{I}_d+(s-t)F(z))=\mathbf{I}_d,
\end{equation*}
which implies (note that $F(z)^{\top}=F(z)$ by \Cref{L2.4}(c))
\begin{align}\label{E4.206}
    & \big(-(\mathbf{I}_d+(s-t)F(z))w\big)^{\top}\bigg(\frac{1}{s-t}\mathbf{I}_d-F(x)\bigg)\big(-(\mathbf{I}_d+(s-t)F(z))w\big)\nonumber\\
    =\,& w^{\top}(\mathbf{I}_d+(s-t)F(z))^{\top}\bigg(\frac{1}{s-t}\mathbf{I}_d-F(x)\bigg)(\mathbf{I}_d+(s-t)F(z))w\nonumber\\
    =\,&  w^{\top}\bigg(\frac{1}{s-t} \mathbf{I}_d+F(z)\bigg)w.
\end{align}
Hence by \eqref{E4.199}--\eqref{E4.202}, noting that $\|y-z_0\|=s_0-t_0$ and $\eps$ is sufficiently small, we have 
\begin{align*}
    &\bigg| w_0^{\top}\bigg(\frac{1}{\|y-z_0\|}\mathbf{I}_d-F(z_0)\bigg)w_0-w^{\top}\bigg(\frac{1}{s-t} \mathbf{I}_d+F(z)\bigg)w\bigg|\nonumber\\
    =  \,& \bigg|w_0^{\top}\bigg(\frac{1}{s_0-t_0}\mathbf{I}_d-F(z_0)\bigg)w_0\nonumber\\
    &\hspace{0.1in}-\big(-(\mathbf{I}_d+(s-t)F(z))w\big)^{\top}\bigg(\frac{1}{s-t}\mathbf{I}_d-F(x)\bigg)\big(-(\mathbf{I}_d+(s-t)F(z))w\big)\bigg|\nonumber\\
    \leq\,& \bigg|w_0^{\top}\bigg(\frac{1}{s_0-t_0}\mathbf{I}_d-F(z_0)\bigg)w_0-w_0^{\top}\bigg(\frac{1}{s-t}\mathbf{I}_d-F(x)\bigg)w_0\bigg|\nonumber\\
    &+\bigg|\big(w_0+(\mathbf{I}_d+(s-t)F(z))w\big)^{\top}\bigg(\frac{1}{s-t}\mathbf{I}_d-F(x)\bigg)w_0\bigg|\nonumber\\
    &+\bigg|\big(-(\mathbf{I}_d+(s-t)F(z))w\big)^{\top}\bigg(\frac{1}{s-t}\mathbf{I}_d-F(x)\bigg)\big(w_0+(\mathbf{I}_d+(s-t)F(z))w\big)\bigg|\nonumber\\
    \leq\,& \bigg(\bigg\|\frac{1}{s_0-t_0}\mathbf{I}_d-\frac{1}{s-t}\mathbf{I}_d\bigg\|_2+\|F(z_0)-F(x)\|_2\bigg)\|w_0\|^2\nonumber\\
    &+2\max\{\|w_0\|,\|(\mathbf{I}_d+(s-t)F(z))w\|\}\bigg\|\frac{1}{s-t}\mathbf{I}_d-F(x)\bigg\|_2\|w_0+(\mathbf{I}_d+(s-t)F(z))w\|\nonumber\\
    \leq\,& C\delta\cdot (\delta_1^{-1}M_1\sqrt{\eps})^2 +C\delta_1^{-1}M_1\sqrt{\eps}\cdot \delta_1^{-1}\cdot \delta M_1 \sqrt{\eps}\leq C\delta M_1^2\delta_1^{-2}\eps. \qedhere
\end{align*}
\end{proof}

Combining \eqref{pbdds} and Claims \ref{Claim.a1}--\ref{Claim.a3}, noting that $\eps$ is sufficiently small, we conclude that
\begin{align*}
    \eps^{(d-1)\slash 2}p_2(x,y)\geq\,&e^{-CM_1\delta_1^{-2}\delta_0^{-1}\sqrt{\eps}}(\mathscr{I}_1(q,s,t)-\delta)_{+}\cdot e^{-C\delta M_1^2\delta_1^{-2}}e^{-\frac{1}{2\eps}w^{\top}\big(\frac{1}{s-t}\mathbf{I}_d+F(z)\big)w}\nonumber\\
    \geq\,& e^{-C\delta M_1^2\delta_1^{-2}}(\mathscr{I}_1(q,s,t)-\delta)_{+}e^{-\frac{1}{2\eps}w^{\top}\big(\frac{1}{s-t}\mathbf{I}_d+F(z)\big)w}.
\end{align*}
This completes the proof of (a).

\medskip

\noindent\textbf{Proof of (b) and (c).} By \eqref{Vdds}, noting that $\eps$ is sufficiently small, we have 
\begin{equation*}
    \|V(y)+V(z)\|^2=4-\|V(y)-V(z)\|^2\geq 4-4M_1^2\delta_1^{-2}\eps>0,
\end{equation*}
hence $V(x)+V(y)=V(z)+V(y)\neq 0$. We define (recall \eqref{x'd} and \eqref{y'd})
\begin{equation}\label{xpp1}
    x':=\mathtt{x}'(x,y)=y-\frac{\langle y-x,V(z)+V(y)\rangle}{\langle V(y),V(z)+V(y)\rangle} V(y),
\end{equation}
\begin{equation}\label{ypp1}
    y':=\mathtt{y}'(x,y)= x+\frac{\langle y-x,V(z)+V(y)\rangle}{\langle V(z), V(z)+V(y)\rangle}V(z).
\end{equation}
By \eqref{Vdds}, we have
\begin{align*}
   & |\langle V(y), V(z)+V(y)\rangle -2|=|\langle V(y), V(z)-V(y)\rangle|\leq\|V(y)-V(z)\|\leq 2M_1\delta_1^{-1}\sqrt{\eps},\nonumber\\
   & |\langle V(z), V(z)+V(y)\rangle -2|=|\langle V(z), V(y)-V(z)\rangle|\leq\|V(y)-V(z)\|\leq 2M_1\delta_1^{-1}\sqrt{\eps}.
\end{align*}
Hence noting that $\eps$ is sufficiently small, we get
\begin{equation}\label{E4.205}
     \min\{\langle V(y), V(z)+V(y)\rangle, \langle V(z), V(z)+V(y)\rangle\}\geq 1.
\end{equation}
By \eqref{Vdds} and \eqref{E4.205}, we get
\begin{align}\label{EEEn1}
  &  \bigg|\frac{\langle y-x,V(z)+V(y)\rangle}{\langle V(y),V(z)+V(y)\rangle}+(s-t)\bigg|\nonumber\\
  =\,& \bigg|\frac{\langle -(s-t)V(z), V(z)+V(y)\rangle +(s-t)\langle V(y),V(z)+V(y)\rangle+\langle w, V(z)+V(y)\rangle}{\langle V(y),V(z)+V(y)\rangle}\bigg|\nonumber\\
  \leq\,& |s-t|\|V(z)+V(y)\|\|V(y)-V(z)\|+\|w\|\|V(z)+V(y)\|\nonumber\\
  \leq\,& 2D\cdot 2\cdot 2M_1\delta_1^{-1}\sqrt{\eps}+\frac{1}{2}M_1\sqrt{\eps}\cdot 2 \leq CM_1\delta_1^{-1}\sqrt{\eps}.
\end{align}
By \eqref{alphabdd} and the above display, noting that $\delta\in (0,\delta_1\slash 10)$ and $\eps$ is sufficiently small, we get
\begin{align*}
   & | \alpha(y) - \langle x'-y,V(y)\rangle -\alpha(x)|= \bigg|\alpha(y)-\alpha(x)+\frac{\langle y-x,V(z)+V(y)\rangle}{\langle V(y),V(z)+V(y)\rangle}\bigg| \nonumber\\
   \leq\,& |\alpha(y)-\alpha(z)|+|\alpha(z)-\alpha(x)-(s-t)|+\bigg|\frac{\langle y-x,V(z)+V(y)\rangle}{\langle V(y),V(z)+V(y)\rangle}+(s-t)\bigg|\nonumber\\
   \leq\,& \delta + CM_1\delta_1^{-1}\sqrt{\eps}
   \leq \delta_1,
\end{align*}
\begin{align*}
& | \beta(y) + \langle x'-y,V(y)\rangle -\beta(x)|=\bigg|\beta(y)-\beta(x)-\frac{\langle y-x,V(z)+V(y)\rangle}{\langle V(y),V(z)+V(y)\rangle}\bigg|\nonumber\\
\leq\,& |\beta(y)-\beta(z)|+|\beta(z)-\beta(x)+s-t|+\bigg|\frac{\langle y-x,V(z)+V(y)\rangle}{\langle V(y),V(z)+V(y)\rangle}+(s-t)\bigg|\nonumber\\
   \leq\,& \delta + CM_1\delta_1^{-1}\sqrt{\eps}
   \leq \delta_1.
\end{align*}
By \eqref{alpha_beta_delta1} and the above two displays, we get
\begin{equation*}
    \min\{\alpha(y) - \langle x'-y,V(y)\rangle,  \beta(y) + \langle x'-y,V(y)\rangle\}\geq \min\{\alpha(x),\beta(x)\}-\delta_1\geq\delta_1.
\end{equation*}
Hence $x'\in\mathrm{int}(T(y))$, which implies $\alpha(x')=\alpha(y) - \langle x'-y,V(y)\rangle$ and $\beta(x')=\beta(y) + \langle x'-y,V(y)\rangle$. Consequently, we have
\begin{equation}\label{xp1}
    x'\in\mathrm{int}(T(y)),\quad \min\{\alpha(x'),\beta(x')\}\geq\delta_1.
\end{equation}
Now by \eqref{E4.205}, we get
\begin{align}\label{E4.208}
    \bigg|\frac{\langle y-x,V(z)+V(y)\rangle}{\langle V(z),V(z)+V(y)\rangle}+(s-t)\bigg|=\,&\bigg|\frac{\langle w, V(z)+V(y)\rangle}{\langle V(z),V(z)+V(y)\rangle}\bigg|\leq\|w\|\|V(z)+V(y)\|\nonumber\\
     \leq\,& \frac{1}{2}M_1\sqrt{\eps}\cdot 2 \leq M_1\sqrt{\eps}.
\end{align}
Hence noting that $\eps$ is sufficiently small, we get
\begin{align*}
   & | \alpha(x) - \langle y'-x,V(z)\rangle -\alpha(z)|= \bigg|\alpha(x)-\alpha(z)-\frac{\langle y-x,V(z)+V(y)\rangle}{\langle V(z),V(z)+V(y)\rangle}\bigg| \nonumber\\
  =\,& \bigg|\frac{\langle y-x,V(z)+V(y)\rangle}{\langle V(z),V(z)+V(y)\rangle}+(s-t)\bigg|
   \leq M_1\sqrt{\eps}
   \leq \delta_1,
\end{align*}
\begin{align*}
& | \beta(x) + \langle y'-x,V(z)\rangle -\beta(z)|=\bigg|\beta(x)-\beta(z)+\frac{\langle y-x,V(z)+V(y)\rangle}{\langle V(z),V(z)+V(y)\rangle}\bigg|\nonumber\\
=\,& \bigg|\frac{\langle y-x,V(z)+V(y)\rangle}{\langle V(z),V(z)+V(y)\rangle}+(s-t)\bigg|\leq M_1\sqrt{\eps}
   \leq \delta_1.
\end{align*}
By \eqref{alpha_beta_delta1} and the above two displays, we get
\begin{equation*}
    \min\{\alpha(x) - \langle y'-x,V(z)\rangle,  \beta(x) + \langle y'-x,V(z)\rangle\}  \geq \min\{\alpha(z),\beta(z)\}-\delta_1\geq\delta_1.
\end{equation*}
Hence $y'\in\mathrm{int}(T(x))$, which implies $\alpha(y')=\alpha(x) - \langle y'-x,V(z)\rangle$ and $\beta(y')=\beta(x) + \langle y'-x,V(z)\rangle$. Consequently, we have
\begin{equation}\label{yp1}
    y'\in\mathrm{int}(T(x)),\quad \min\{\alpha(y'),\beta(y')\}\geq\delta_1.
\end{equation}

Noting \eqref{xp1}--\eqref{yp1} and \eqref{zwzw}, we define
\begin{align}\label{z1w1z2w2}
    z_1&:=\mathtt{z}(x',y') = y-\langle y-x,V(y)\rangle V(y)+\frac{\langle y-x, V(z)+V(y)\rangle \langle V(z), V(y)\rangle}{\langle V(z), V(z)+V(y)\rangle} V(y) ,\nonumber\\
    w_1&:=\mathtt{w}(x',y')=x-y-\langle x-y,V(y)\rangle V(y)+\frac{\langle y-x, V(z)+V(y)\rangle}{\langle V(z), V(z)+V(y)\rangle} \big(V(z)-\langle V(z),V(y)\rangle V(y)\big) , \nonumber\\
    z_2&:=\mathtt{z}'(x',y')=x-\langle x-y, V(z)\rangle V(z)+\frac{\langle x-y, V(z)+V(y)\rangle \langle V(z), V(y)\rangle}{\langle V(y), V(z)+V(y)\rangle} V(z),\nonumber\\
    w_2&:=\mathtt{w}'(x',y')=y-x-\langle y-x,V(z)\rangle V(z)+\frac{\langle x-y, V(z)+V(y)\rangle}{\langle V(y), V(z)+V(y)\rangle} \big(V(y)-\langle V(z),V(y)\rangle V(z)\big).
\end{align}
We also define (recall that $q'=\mathfrak{q}_{k,k'}(y)$ and $V(q')=V(y)$)
\begin{equation}\label{s1s2}
   s_1:=\langle x'-q', V(y)\rangle, \quad     
t_1:=\langle z_1-q', V(y)\rangle, \quad
s_2:=\langle z_2-q, V(z)\rangle, \quad t_2:=\langle y'-q, V(z)\rangle.
\end{equation}
Note that $x'=q'+s_1V(q')$, $z_1=q'+t_1 V(q')$, $z_2=q+s_2 V(q)$, and $y'=q +t_2 V(q)$. 

By \eqref{Vdds}, \eqref{xpp1}--\eqref{EEEn1}, and \eqref{E4.208}, noting that $\eps$ is sufficiently small, we have
\begin{align}\label{x'bdd}
   & \|x'-x\|=\bigg\|y-x-\frac{\langle y-x,V(z)+V(y)\rangle}{\langle V(y),V(z)+V(y)\rangle} V(y)\bigg\|\nonumber\\
   \leq\,& \|w\|+\|z-x+(s-t)V(z)\|+\bigg\|\bigg(s-t+\frac{\langle y-x,V(z)+V(y)\rangle}{\langle V(y),V(z)+V(y)\rangle}\bigg) V(z)\bigg\|\nonumber\\
   &+\frac{|\langle y-x,V(z)+V(y)\rangle|}{|\langle V(y),V(z)+V(y)\rangle|}\|V(y)-V(z)\|\nonumber\\
   \leq\,& \frac{1}{2}M_1\sqrt{\eps}+\bigg|s-t+\frac{\langle y-x,V(z)+V(y)\rangle}{\langle V(y),V(z)+V(y)\rangle}\bigg|+(\|z-x\|+\|w\|)\|V(z)+V(y)\|\|V(y)-V(z)\|\nonumber\\
   \leq\,& \frac{1}{2}M_1\sqrt{\eps}+CM_1\delta_1^{-1}\sqrt{\eps}+\Big(2D+\frac{1}{2}M_1\sqrt{\eps}\Big)\cdot 2 \cdot 2M_1\delta_1^{-1}\sqrt{\eps}\leq CM_1\delta_1^{-1}\sqrt{\eps},
\end{align}
\begin{equation}\label{y'bdd}
    \|y'-z\|=\bigg\|x-z+\frac{\langle y-x,V(z)+V(y)\rangle}{\langle V(z), V(z)+V(y)\rangle}V(z)\bigg\|=\bigg|s-t+\frac{\langle y-x,V(z)+V(y)\rangle}{\langle V(z),V(z)+V(y)\rangle}\bigg|\leq M_1\sqrt{\eps}.
\end{equation}
Similarly, using \eqref{Vdds} and \eqref{z1w1z2w2} we can deduce that
\begin{equation}\label{z1z2ub}
    \|z_1-z\|\leq CM_1\delta_1^{-1}\sqrt{\eps}, \qquad \|z_2-x\|\leq CM_1\delta_1^{-1}\sqrt{\eps}.
\end{equation}
By the above displays, noting that $\eps$ is sufficiently small and $x\in\mathcal{X}_k^{\circ}, z\in\mathcal{Y}_{k'}^{\circ}$ (see \eqref{xzres}), we get
\begin{equation}\label{x'y'}
    x',z_2\in\mathcal{X}_k^{\circ}, \qquad y', z_1\in\mathcal{Y}_{k'}^{\circ}.
\end{equation}
Moreover, by \eqref{z1z2ub}, $\|z_1-y\|\leq \|z_1-z\|+\|w\|\leq CM_1\delta_1^{-1}\sqrt{\eps}\leq \delta_1\slash 10$ and $\|z_2-x\|<\delta_1\slash 10$ (as $\eps$ is sufficiently small); by \eqref{alpha_beta_delta1} and \eqref{alphabdd}, we have $\min\{\alpha(x),\beta(x),\alpha(y),\beta(y)\}\geq 1.9\delta_1$ (where we use $\delta\in (0,\delta_1\slash 10)$). Consequently,
\begin{equation}\label{z1z2w1w2res}
    z_1\in\mathrm{int}(T(y)), \quad z_2\in\mathrm{int}(T(x)), \quad \min\{\alpha(z_1),\beta(z_1)\}\geq \delta_1, \quad \min\{\alpha(z_2),\beta(z_2)\}\geq \delta_1.
\end{equation}
By \eqref{s1s2}, \eqref{Vdds}, \eqref{qbd}, and \eqref{x'bdd}--\eqref{z1z2ub} we have
\begin{align}\label{s1t1}
    |s_1-s|=\,& |\langle (x'-q')-(x-q), V(y)\rangle +\langle x-q, V(y)-V(z)\rangle|\nonumber\\
    \leq\,& \|x'-x\|+\|q'-q\|+2D\|V(y)-V(z)\|\leq CM_1\delta_1^{-1}\sqrt{\eps},
\end{align}
\begin{align}\label{s1t11}
    |t_1-t|=\,&|\langle (z_1-q')-(z-q),V(y)\rangle+\langle z-q,V(y)-V(z)\rangle|\nonumber\\
    \leq\,& \|z_1-z\|+\|q'-q\|+2D\|V(y)-V(z)\|\leq  CM_1\delta_1^{-1}\sqrt{\eps},
\end{align}
\begin{equation}\label{s2t2}
    |s_2-s|=|\langle (z_2-q)-(x-q),V(z)\rangle|\leq \|z_2-x\| \leq CM_1\delta_1^{-1}\sqrt{\eps},
\end{equation}
\begin{equation}\label{s2t21}
    |t_2-t|=|\langle (y'-q)-(z-q), V(z)\rangle|\leq \|y'-z\|\leq M_1\sqrt{\eps}.
\end{equation}
As $s-t\geq 2d_0$ (recall Definition \ref{def.4.6}) and $\eps$ is sufficiently small, we have  
\begin{equation}\label{s1t1s2t2b}
    s_1-t_1\geq d_0, \qquad s_2-t_2\geq d_0.
\end{equation}

By \eqref{p1}--\eqref{p2} and \eqref{p3}--\eqref{p4}, we have (recall that $y=q'+t_0 V(q')$)
\begin{align}\label{E4.226}
   \eps^{(d-1)\slash 2} p_3(x,y) =\,&\eps^{(d-1)\slash 2} p_2(x',y')\cdot \frac{\det(\mathbf{I}_d+s_1 F(q'))\det(\mathbf{I}_d+t_2F(q))}{\det(\mathbf{I}_d+t_0  F(q'))\det(\mathbf{I}_d+sF(q))}\nonumber\\
   =\,& \mathfrak{F}(y')^{-1}h(z_2,y')(2\pi(s_2-t_2))^{-(d-1)\slash 2}\sqrt{\det\big(\mathbf{I}_d-(s_2-t_2)F(z_2)\big)}\nonumber\\
   &\times e^{-\frac{1}{2\eps} w_2^{\top} \big(\frac{1}{s_2-t_2}\mathbf{I}_d-F(z_2)\big) w_2}\cdot \frac{\det(\mathbf{I}_d+s_1 F(q'))\det(\mathbf{I}_d+t_2F(q))}{\det(\mathbf{I}_d+t_0  F(q'))\det(\mathbf{I}_d+sF(q))},
\end{align}
\begin{align}\label{E4.227}
    \eps^{(d-1)\slash 2} p_4(x,y)=\,&\eps^{(d-1)\slash 2} p_1(x',y')\cdot \frac{\det(\mathbf{I}_d+s_1 F(q'))\det(\mathbf{I}_d+t_2F(q))}{\det(\mathbf{I}_d+t_0  F(q'))\det(\mathbf{I}_d+sF(q))}\nonumber\\
    =\,& \mathfrak{F}(x')^{-1}h(x',z_1)(2\pi(s_1-t_1))^{-(d-1)\slash 2}\sqrt{\det\big(\mathbf{I}_d+(s_1-t_1)F(z_1)\big)} \cdot \mathbbm{1}_{\|w_1\|\leq M_1\sqrt{\eps}}\nonumber\\
   &\times e^{-\frac{1}{2\eps} w_1^{\top} \big(\frac{1}{s_1-t_1}\mathbf{I}_d+F(z_1)\big) w_1}\cdot \frac{\det(\mathbf{I}_d+s_1 F(q'))\det(\mathbf{I}_d+t_2F(q))}{\det(\mathbf{I}_d+t_0  F(q'))\det(\mathbf{I}_d+sF(q))},
\end{align}
where we use \eqref{xp1}, \eqref{yp1}, \eqref{x'y'}, \eqref{z1z2w1w2res}, and \eqref{s1t1s2t2b}.

Next, we prove three further claims which together will imply (b) and (c). 

\begin{claim}\label{Claim.a4}
We have
\begin{equation*}
    \mathfrak{F}(x')^{-1}h(x',z_1)(2\pi(s_1-t_1))^{-(d-1)\slash 2}\sqrt{\det\big(\mathbf{I}_d+(s_1-t_1)F(z_1)\big)} \geq e^{-CM_1\delta_1^{-2}\delta_0^{-1}\sqrt{\eps}}(\mathscr{I}_1(q,s,t)-\delta)_{+},
\end{equation*}
\begin{equation*}
    \mathfrak{F}(y')^{-1}h(z_2,y')(2\pi(s_2-t_2))^{-(d-1)\slash 2}\sqrt{\det\big(\mathbf{I}_d-(s_2-t_2)F(z_2)\big)}\geq e^{-CM_1\delta_1^{-2}\delta_0^{-1}\sqrt{\eps}}\mathscr{I}_1(q,s,t).
\end{equation*}
\end{claim}

\begin{proof}

By Definitions \ref{factors} and \ref{def.4.6}, noting \eqref{yp1} and \eqref{z1z2w1w2res}, we have
\begin{align*}
&\mathfrak{F}(z_2)\sqrt{\det\big(\mathbf{I}_d-(s_2-t_2)F(z_2)\big)}\nonumber\\
=\,& \mathscr{I}_0(q)^{-1}\det\big(\mathbf{I}_d+\langle z_2-q,V(q)\rangle F(q)\big)\sqrt{\det\big(\mathbf{I}_d+\langle y'-z_2,V(q)\rangle 
F(z_2)\big)}\nonumber\\
=\,& \mathscr{I}_0(q)^{-1}\sqrt{\det\big(\mathbf{I}_d+\langle z_2-q,V(q)\rangle F(q)\big)}\nonumber\\
&\times\sqrt{\det\big(\mathbf{I}_d+\langle z_2-q,V(q)\rangle F(q)\big)\det\big(\mathbf{I}_d+\langle y'-z_2,V(q)\rangle 
F(z_2)\big)}\nonumber\\
=\,& \mathscr{I}_0(q)^{-1}\sqrt{\det\big(\mathbf{I}_d+\langle z_2-q,V(q)\rangle F(q)\big)}\cdot\sqrt{\det\big(\mathbf{I}_d+\langle y'-q,V(q)\rangle F(q)\big)},
\end{align*}
where we use \eqref{Eq2.6.4} in \Cref{L2.4} (with $x,x',x_0$ replaced by $z_2,q,y'$) in the last equality. Similarly,
\begin{align*}
   & \mathfrak{F}(z_1)  \sqrt{\det\big(\mathbf{I}_d+(s_1-t_1)F(z_1)\big)}\nonumber\\
   =\,& \mathscr{I}_0(q')^{-1}\det\big(\mathbf{I}_d+\langle z_1-q',V(q')\rangle F(q')\big)\sqrt{\det\big(\mathbf{I}_d+\langle x'-z_1,V(q')\rangle F(z_1)\big)}\nonumber\\
   =\,& \mathscr{I}_0(q')^{-1}\sqrt{\det\big(\mathbf{I}_d+\langle z_1-q',V(q')\rangle F(q')\big)}\cdot\sqrt{\det\big(\mathbf{I}_d+\langle x'-q',V(q')\rangle F(q')\big)}.
\end{align*}
Hence by Definitions \ref{De3.8} and \ref{DefhH},
\begin{align}\label{series1}
  &  \mathfrak{F}(y')^{-1}h(z_2,y')(2\pi(s_2-t_2))^{-(d-1)\slash 2}\sqrt{\det\big(\mathbf{I}_d-(s_2-t_2)F(z_2)\big)}\nonumber\\
  =\,& \mathfrak{F}(y')^{-1}\tilde{f}(z_2)\tilde{g}(y')\tilde{h}_{T(q)}(z_2,y')(2\pi(s_2-t_2))^{-(d-1)\slash 2}\sqrt{\det\big(\mathbf{I}_d-(s_2-t_2)F(z_2)\big)}\nonumber\\
  =\,& f(z_2)g(y')\tilde{h}_{T(q)}(z_2,y')(2\pi(s_2-t_2))^{-(d-1)\slash 2}\mathfrak{F}(z_2)\sqrt{\det\big(\mathbf{I}_d-(s_2-t_2)F(z_2)\big)}\nonumber\\
  =\,& \mathscr{I}_0(q)^{-1}f(z_2)g(y')\tilde{h}_{T(q)}(z_2,y')(2\pi(s_2-t_2))^{-(d-1)\slash 2}\nonumber\\
  &\times\sqrt{\det\big(\mathbf{I}_d+\langle z_2-q,V(q)\rangle F(q)\big)}\cdot\sqrt{\det\big(\mathbf{I}_d+\langle y'-q,V(q)\rangle F(q)\big)},
\end{align}
\begin{align}
&\mathfrak{F}(x')^{-1}h(x',z_1)(2\pi(s_1-t_1))^{-(d-1)\slash 2}\sqrt{\det\big(\mathbf{I}_d+(s_1-t_1)F(z_1)\big)}\nonumber\\
=\,&  \mathfrak{F}(x')^{-1}\tilde{f}(x')\tilde{g}(z_1)\tilde{h}_{T(q')}(x',z_1)(2\pi(s_1-t_1))^{-(d-1)\slash 2}\sqrt{\det\big(\mathbf{I}_d+(s_1-t_1)F(z_1)\big)}\nonumber\\
=\,& f(x')g(z_1)\tilde{h}_{T(q')}(x',z_1)(2\pi(s_1-t_1))^{-(d-1)\slash 2}\mathfrak{F}(z_1)\sqrt{\det\big(\mathbf{I}_d+(s_1-t_1)F(z_1)\big)}\nonumber\\
=\,& \mathscr{I}_0(q')^{-1}f(x')g(z_1)\tilde{h}_{T(q')}(x',z_1)(2\pi(s_1-t_1))^{-(d-1)\slash 2}\nonumber\\
&\times \sqrt{\det\big(\mathbf{I}_d+\langle z_1-q',V(q')\rangle F(q')\big)}\cdot\sqrt{\det\big(\mathbf{I}_d+\langle x'-q',V(q')\rangle F(q')\big)}.
\end{align}

As $s-t\geq 2d_0$, by \eqref{s1t1}--\eqref{s2t21}, we have
\begin{equation*}
    \max\{s_1-t_1,s_2-t_2\}\leq s-t+CM_1\delta_1^{-1}\sqrt{\eps}\leq (1+CM_1\delta_1^{-1}\sqrt{\eps})(s-t).
\end{equation*}
Hence
\begin{equation}
    \min\{(2\pi(s_1-t_1))^{-(d-1)\slash 2}, (2\pi(s_2-t_2))^{-(d-1)\slash 2}\}\geq (1+CM_1\delta_1^{-1}\sqrt{\eps})^{-(d-1)\slash 2}(2\pi(s-t))^{-(d-1)\slash 2}.
\end{equation}

By \eqref{bddLip} and \eqref{x'bdd}, we have $|f(x')-f(x)|\leq M_0\|x'-x\|\leq CM_1\delta_1^{-1}\sqrt{\eps}$. As $\eps$ is sufficiently small, noting \eqref{xzres}, we get $f(x')\geq \delta_0\slash 2$. Similarly, we can deduce that $f(z_2),g(y'),g(z_1)\geq \delta_0\slash 2$. Now by \eqref{bddn1}--\eqref{bddn0} and \eqref{x'bdd}--\eqref{z1z2ub}, we have  
\begin{align*}
    |f(x')-f(q'+sV(q'))|&\leq M_0(\|x'-x\|+\|q'+sV(q')-x\|)\leq CM_1\delta_1^{-1}\sqrt{\eps}\leq CM_1\delta_1^{-1}\delta_0^{-1}\sqrt{\eps}f(x'),\nonumber\\
    |g(z_1)-g(q'+tV(q'))|&\leq M_0(\|z_1-z\|+\|q'+tV(q')-z\|)\leq CM_1\delta_1^{-1}\sqrt{\eps}\leq CM_1\delta_1^{-1}\delta_0^{-1}\sqrt{\eps}g(z_1),\nonumber\\
    |f(z_2)-f(x)|&\leq M_0\|z_2-x\|\leq CM_1\delta_1^{-1}\sqrt{\eps}\leq CM_1\delta_1^{-1}\delta_0^{-1}\sqrt{\eps}f(z_2),\nonumber\\
    |g(y')-g(z)|&\leq M_0\|y'-z\|\leq CM_1\delta_1^{-1}\sqrt{\eps}\leq CM_1\delta_1^{-1}\delta_0^{-1}\sqrt{\eps}g(y').
\end{align*}
Hence
\begin{align}
    & f(x')\geq (1+CM_1\delta_1^{-1}\delta_0^{-1}\sqrt{\eps})^{-1}f(q'+sV(q')), \quad g(z_1)\geq (1+CM_1\delta_1^{-1}\delta_0^{-1}\sqrt{\eps})^{-1}g(q'+tV(q')),\nonumber\\
    & f(z_2)\geq (1+CM_1\delta_1^{-1}\delta_0^{-1}\sqrt{\eps})^{-1}f(x), \quad g(y')\geq (1+CM_1\delta_1^{-1}\delta_0^{-1}\sqrt{\eps})^{-1}g(z).
\end{align}

As $T(q')\in\tilde{\mathcal{S}}^*$ (see \eqref{mathIbound}),
by \Cref{Lem2} and noting \eqref{xzres}, \eqref{xkykk}, and \eqref{x'y'}, we get
\begin{align}
    \tilde{h}_{T(q')}(x',z_1)=\,& \tilde{h}_{T(q')}(q'+s_1V(q'), q'+t_1 V(q'))\geq e^{-C(|s_1-s|+|t_1-t|)}\tilde{h}_{T(q')}(q'+sV(q'),q'+tV(q'))\nonumber\\
    \geq\,&  e^{-CM_1\delta_1^{-1}\sqrt{\eps}}\tilde{h}_{T(q')}(q'+sV(q'),q'+tV(q')),
\end{align}
\begin{align}
    \tilde{h}_{T(q)}(z_2,y')=\,&\tilde{h}_{T(q)}(q+s_2V(q),q+t_2 V(q))\geq e^{-C(|s_2-s|+|t_2-t|)}\tilde{h}_{T(q)}(q+sV(q),q+tV(q))\nonumber\\
    \geq\,& e^{-CM_1\delta_1^{-1}\sqrt{\eps}}\tilde{h}_{T(q)}(q+sV(q),q+tV(q))=e^{-CM_1\delta_1^{-1}\sqrt{\eps}}\tilde{h}_{T(q)}(x,z),
\end{align}
where we use \eqref{s1t1}--\eqref{s2t21}.

By \eqref{E4.195} (taking $(\theta,\theta')=(x',q'+sV(q'))$ or $(\theta,\theta')=(z_1,q'+t V(q'))$), noting \eqref{xp1} and \eqref{z1z2w1w2res}--\eqref{s1t11}, we get
\begin{align}
    \det\big(\mathbf{I}_d+\langle x'-q',V(q')\rangle F(q')\big)\geq\,& \bigg(1+\frac{C|s_1-s|}{\delta_1}\bigg)^{-d}\det(\mathbf{I}_d+s F(q'))\nonumber\\
    \geq\,& (1+CM_1\delta_1^{-2}\sqrt{\eps})^{-d}\det(\mathbf{I}_d+s F(q')),
\end{align}
\begin{align}
   \det\big(\mathbf{I}_d+\langle z_1-q',V(q')\rangle F(q')\big)\geq\,& \bigg(1+\frac{C|t_1-t|}{\delta_1}\bigg)^{-d}\det(\mathbf{I}_d+t F(q'))\nonumber\\
   \geq\,&(1+CM_1\delta_1^{-2}\sqrt{\eps})^{-d}\det(\mathbf{I}_d+ t F(q')).
\end{align}
Similarly, we can deduce that
\begin{equation}
    \det\big(\mathbf{I}_d+\langle z_2-q, V(q)\rangle F(q)\big)\geq (1+CM_1\delta_1^{-2}\sqrt{\eps})^{-d}\det(\mathbf{I}_d+ s F(q)),
\end{equation}
\begin{equation}\label{series1.1}
    \det\big(\mathbf{I}_d+\langle y'-q, V(q)\rangle F(q)\big)\geq (1+CM_1\delta_1^{-2}\sqrt{\eps})^{-d}\det(\mathbf{I}_d+ t F(q)).
\end{equation}

Combining \eqref{series1}--\eqref{series1.1}, we obtain that
\begin{align}
    &\mathfrak{F}(x')^{-1}h(x',z_1)(2\pi(s_1-t_1))^{-(d-1)\slash 2}\sqrt{\det\big(\mathbf{I}_d+(s_1-t_1)F(z_1)\big)}\nonumber\\
    \geq\,&  (1+CM_1\delta_1^{-1}\delta_0^{-1}\sqrt{\eps})^{-2}e^{-CM_1\delta_1^{-1}\sqrt{\eps}}(1+CM_1\delta_1^{-1}\sqrt{\eps})^{-(d-1)\slash 2}(1+CM_1\delta_1^{-2}\sqrt{\eps})^{-d}\nonumber\\
    &\times\mathscr{I}_0(q')^{-1}f(q'+sV(q'))g(q'+tV(q'))\tilde{h}_{T(q')}(q'+sV(q'),q'+tV(q'))\nonumber\\
    &\times(2\pi(s-t))^{-(d-1)\slash 2}\sqrt{\det(\mathbf{I}_d+s F(q'))}\cdot\sqrt{\det(\mathbf{I}_d+ t F(q'))}\nonumber\\
    \geq\,& e^{-CM_1\delta_1^{-2}\delta_0^{-1}\sqrt{\eps}}\mathfrak{F}(q'+sV(q'))^{-1} h(q'+sV(q'),q'+tV(q'))  \nonumber\\
    & \times\mathscr{I}_0(q')^{-1}\mathfrak{F}(q'+tV(q'))^{-1}(2\pi(s-t))^{-(d-1)\slash 2}\sqrt{\det(\mathbf{I}_d+sF(q'))}\cdot \sqrt{\det(\mathbf{I}_d+t F(q'))}\nonumber\\
    \geq\,& e^{-CM_1\delta_1^{-2}\delta_0^{-1}\sqrt{\eps}}(\mathscr{I}_1(q,s,t)-\delta)_{+},
\end{align}
where we use Definitions \ref{De3.8} and \ref{DefhH} in the second inequality and use \eqref{rel1}--\eqref{rel1.1} in the third inequality. Similarly, we have
\begin{align}
    &\mathfrak{F}(y')^{-1}h(z_2,y')(2\pi(s_2-t_2))^{-(d-1)\slash 2}\sqrt{\det\big(\mathbf{I}_d-(s_2-t_2)F(z_2)\big)}\nonumber\\
  \geq\,&  (1+CM_1\delta_1^{-1}\delta_0^{-1}\sqrt{\eps})^{-2}e^{-CM_1\delta_1^{-1}\sqrt{\eps}}(1+CM_1\delta_1^{-1}\sqrt{\eps})^{-(d-1)\slash 2}(1+CM_1\delta_1^{-2}\sqrt{\eps})^{-d}\nonumber\\
    &\times\mathscr{I}_0(q)^{-1}f(x)g(z)\tilde{h}_{T(q)}(x,z)(2\pi(s-t))^{-(d-1)\slash 2}\sqrt{\det(\mathbf{I}_d+s F(q))}\cdot\sqrt{\det(\mathbf{I}_d+ t F(q))}\nonumber\\
    \geq\,& e^{-CM_1\delta_1^{-2}\delta_0^{-1}\sqrt{\eps}}\mathscr{I}_1(q,s,t),
\end{align}
where in the last inequality, we note that by Definitions \ref{factors}, \ref{De3.8}, \ref{DefhH}, and \ref{def.4.6} as well as \Cref{L2.4} (see \eqref{Eq2.6.4}, with $x,x',x_0$ replaced by $q+tV(q), q, q+sV(q)$), 
\begin{align*}
 &\mathscr{I}_0(q)^{-1}f(x)g(z)\tilde{h}_{T(q)}(x,z)(2\pi(s-t))^{-(d-1)\slash 2}\sqrt{\det(\mathbf{I}_d+s F(q))}\cdot\sqrt{\det(\mathbf{I}_d+ t F(q))}\nonumber\\
 =\,& \mathscr{I}_0(q)^{-1}\mathfrak{F}(x)^{-1}\mathfrak{F}(z)^{-1}h(x,z)(2\pi(s-t))^{-(d-1)\slash 2}\sqrt{\det(\mathbf{I}_d+s F(q))}\cdot\sqrt{\det(\mathbf{I}_d+ t F(q))}\nonumber\\
 =\,& \mathfrak{F}(x)^{-1}h(x,z)(2\pi(s-t))^{-(d-1)\slash 2}\sqrt{\frac{\det(\mathbf{I}_d+sF(q))}{\det(\mathbf{I}_d+tF(q))}}\nonumber\\
 =\,&\mathfrak{F}(x)^{-1}h(x,z)(2\pi(s-t))^{-(d-1)\slash 2}\sqrt{\det\big(\mathbf{I}_d+(s-t)F(z)\big)}=\mathscr{I}_1(q,s,t). \qedhere
\end{align*}
\end{proof}

\begin{claim}\label{Claim.a5}
We have
\begin{equation*}
  \frac{\det(\mathbf{I}_d+s_1 F(q'))\det(\mathbf{I}_d+t_2F(q))}{\det(\mathbf{I}_d+t_0  F(q'))\det(\mathbf{I}_d+sF(q))}\geq e^{-CM_1\delta_1^{-2}\sqrt{\eps}} (1+C\delta\delta_1^{-d})^{-2}.
\end{equation*}
\end{claim}

\begin{proof}
By \eqref{E4.195} (taking $\theta=x'=q'+s_1V(q')$ and $\theta'=q'+sV(q')$), noting \eqref{xp1} and \eqref{s1t1}, we have 
\begin{equation*}
    \det(\mathbf{I}_d+s_1 F(q'))\geq (1+C\delta_1^{-1}|s_1-s|)^{-d}\det(\mathbf{I}_d+s F(q')) \geq (1+CM_1\delta_1^{-2}\sqrt{\eps})^{-d}\det(\mathbf{I}_d+s F(q')).
\end{equation*}
Similarly, by \eqref{yp1}, \eqref{s2t21}, \eqref{ldalpha}, and \eqref{tbdd}, 
\begin{equation*}
    \det(\mathbf{I}_d+t_2 F(q))\geq (1+C\delta_1^{-1}|t_2-t|)^{-d}\det(\mathbf{I}_d+t F(q)) \geq (1+CM_1\delta_1^{-2}\sqrt{\eps})^{-d}\det(\mathbf{I}_d+t F(q)),
\end{equation*}
\begin{equation*}
    \det(\mathbf{I}_d+t F(q'))\geq (1+C\delta_1^{-1}|t_0-t|)^{-d}\det(\mathbf{I}_d+ t_0 F(q')) \geq (1+CM_1\delta_1^{-2}\sqrt{\eps})^{-d}\det(\mathbf{I}_d+t_0 F(q')).
\end{equation*}
Hence noting Lemma \ref{L2.7}, we have
\begin{equation}\label{E4.238}
     \frac{\det(\mathbf{I}_d+s_1 F(q'))\det(\mathbf{I}_d+t_2F(q))}{\det(\mathbf{I}_d+t_0  F(q'))\det(\mathbf{I}_d+sF(q))}\geq e^{-CM_1\delta_1^{-2}\sqrt{\eps}}\cdot\frac{\det(\mathbf{I}_d+s F(q'))\det(\mathbf{I}_d+tF(q))}{\det(\mathbf{I}_d+t  F(q'))\det(\mathbf{I}_d+sF(q))}.
\end{equation}

As $w\in \mathscr{B}(q,s,t;\delta,M_1\sqrt{\eps}\slash 2)$ (recall Definition \ref{def4.9}), we have $\|F(q')-F(q)\|_2\leq \delta$. By \eqref{L2.4.enew3}, noting that $q,q'\in\mathfrak{H}_{k,k'}\subseteq\tilde{H}_{k,k'}$, we have $\max\{\|F(q)\|_2,\|F(q')\|_2\}\leq C$. Hence by \cite[Theorem 2.12]{ipsen2008perturbation},
\begin{align}\label{E4.239}
 &\big|\det(\mathbf{I}_d+sF(q'))-\det(\mathbf{I}_d+sF(q))\big|\nonumber\\
 \leq\,&d\|s(F(q')-F(q))\|_2\max\{\|\mathbf{I}_d+sF(q')\|_2,\|\mathbf{I}_d+sF(q)\|_2\}^{d-1} 
    \leq C\|F(q')-F(q)\|_2\leq C\delta.
\end{align}
By Remark \ref{Rmk.L2.4} (with $x,x'$ replaced by $q,q+sV(q)$),
\begin{equation*}
    \det\big(\mathbf{I}_d-sF(q+sV(q))\big)\det(\mathbf{I}_d+sF(q))=1.
\end{equation*}
By Lemma \ref{L2.7}, $\det\big(\mathbf{I}_d-sF(q+sV(q))\big),\det(\mathbf{I}_d+sF(q))>0$. By \eqref{L2.4.enew3} in \Cref{L2.4} and \eqref{alpha_beta_delta1}, $\|F(q+sV(q))\|_2\leq C   \delta_1^{-1}$, hence by Hadamard's inequality, we have
\begin{equation*}
    \big|\det\big(\mathbf{I}_d-sF(q+sV(q))\big)\big|\leq (1+2D\|F(q+sV(q))\|_2)^d\leq C\delta_1^{-d}.
\end{equation*}
Combining the above two displays, we get $\det(\mathbf{I}_d+sF(q))\geq c\delta_1^d$. Similarly, we can deduce that $\det(\mathbf{I}_d+sF(q'))\geq c\delta_1^d$. Noting \eqref{E4.239}, we get
\begin{equation*}
    \big|\det(\mathbf{I}_d+sF(q'))-\det(\mathbf{I}_d+sF(q))\big|\leq C\delta\delta_1^{-d}\min\big\{\det(\mathbf{I}_d+sF(q)),\det(\mathbf{I}_d+sF(q'))\big\}.
\end{equation*}
Similarly, we can deduce that 
\begin{equation*}
    \big|\det(\mathbf{I}_d+tF(q'))-\det(\mathbf{I}_d+tF(q))\big|\leq C\delta\delta_1^{-d}\min\big\{\det(\mathbf{I}_d+tF(q)),\det(\mathbf{I}_d+tF(q'))\big\}.
\end{equation*}
The above two displays imply that
\begin{align}\label{E4.240}
   & \det(\mathbf{I}_d+sF(q')) \geq (1+C\delta\delta_1^{-d})^{-1}\det(\mathbf{I}_d+sF(q)),\nonumber\\
   & \det(\mathbf{I}_d+tF(q)) \geq (1+C\delta\delta_1^{-d})^{-1}\det(\mathbf{I}_d+tF(q')).
\end{align}

Combining \eqref{E4.238} and \eqref{E4.240}, we conclude that
\begin{equation*}
    \frac{\det(\mathbf{I}_d+s_1 F(q'))\det(\mathbf{I}_d+t_2F(q))}{\det(\mathbf{I}_d+t_0  F(q'))\det(\mathbf{I}_d+sF(q))} \geq e^{-CM_1\delta_1^{-2}\sqrt{\eps}} (1+C\delta\delta_1^{-d})^{-2}. \qedhere
\end{equation*}
\end{proof}

\begin{claim}\label{Claim.a6}
We have $\|w_1\|\leq M_1\sqrt{\eps}$ and
\begin{equation*}
    \bigg|w_1^{\top} \bigg(\frac{1}{s_1-t_1}\mathbf{I}_d+F(z_1)\bigg) w_1-w^{\top}\bigg(\frac{1}{s-t} \mathbf{I}_d+F(z)\bigg)w\bigg|\leq C\delta M_1^2\delta_1^{-1}\eps,
\end{equation*}
\begin{equation*}
    \bigg|w_2^{\top} \bigg(\frac{1}{s_2-t_2}\mathbf{I}_d-F(z_2)\bigg) w_2-w^{\top}\bigg(\frac{1}{s-t} \mathbf{I}_d+F(z)\bigg)w\bigg|\leq C\delta M_1^2\delta_1^{-2}\eps.
\end{equation*}
\end{claim}

\begin{proof}
By \eqref{Vdds}, as $\eps$ is sufficiently small, we have  
\begin{equation}\label{E4.242}
    \langle V(z),V(y)\rangle=\frac{2-\|V(z)-V(y)\|^2}{2}\,\begin{cases}
       \, \leq 1, \\
       \, \geq 1-CM_1^2\delta_1^{-2}\eps \geq \frac{1}{2}.
    \end{cases}
\end{equation}
Hence
\begin{equation}\label{E4.243}
  |\langle V(z),V(z)+V(y)\rangle-2|= |\langle V(y),V(z)+V(y)\rangle-2|=|\langle V(z),V(y)\rangle-1|\leq CM_1^2\delta_1^{-2}\eps.
\end{equation}
Using \eqref{Vdds} and \eqref{E4.243}, noting that $\langle w, V(z)\rangle=0$, we get 
\begin{align}
 &\big|\langle x-y, V(y)\rangle-(s-t)\big|=\big|\langle x-y, V(y)\rangle -\langle x-z, V(z) \rangle \big|\nonumber\\
 \leq\,& \big|\langle w, V(y)\rangle|+\big|\langle x-z, V(y)-V(z)\rangle\big|=\big|\langle w, V(y)-V(z)\rangle|+\big|\langle x-z, V(y)-V(z)\rangle\big|\nonumber\\
 \leq\, & \|w\|\|V(y)-V(z)\|+\|x-z\||\langle V(z),V(y)-V(z)\rangle|\nonumber\\
 =\,&  \|w\|\|V(y)-V(z)\|+\|x-z\||\langle V(z),V(y)\rangle-1|\leq CM_1^2\delta_1^{-2}\eps,
\end{align}
\begin{equation}\label{xyv}
    \langle x-y, V(z)\rangle =\langle x-z,V(z)\rangle -\langle w, V(z)\rangle =s-t.
\end{equation}
By \eqref{BdVy} and \eqref{E4.243}, as $\eps$ is sufficiently small, we get
\begin{align}
   & \|V(z)-\langle V(z), V(y)\rangle V(y)+F(z)w\|\leq\|V(z)+F(z)w-V(y)\|+\|(1-\langle V(z),V(y)\rangle)V(y)\|\nonumber\\
    \leq\,& \frac{1}{2}\delta M_1\sqrt{\eps}+|1-\langle V(z),V(y)\rangle|\leq \frac{1}{2}\delta M_1\sqrt{\eps}+CM_1^2\delta_1^{-2}\eps\leq \delta M_1\sqrt{\eps},
\end{align}
\begin{align}\label{E4.247}
   & \|V(y)-\langle V(z),V(y)\rangle V(z) -F(z)w \|\leq\|V(y)-V(z)-F(z)w\|+\|(1-\langle V(z),V(y)\rangle) V(z)\|\nonumber\\
    \leq \,& \frac{1}{2}\delta M_1\sqrt{\eps}+|1-\langle V(z),V(y)\rangle|\leq \frac{1}{2}\delta M_1\sqrt{\eps}+CM_1^2\delta_1^{-2}\eps\leq \delta M_1\sqrt{\eps}.
\end{align}
By \eqref{E4.242}--\eqref{E4.247}, noting that $\eps$ is sufficiently small, we have
\begin{equation*}
    \bigg\|\frac{\langle y-x, V(z)+V(y)\rangle}{\langle V(z), V(z)+V(y)\rangle} \big(V(z)-\langle V(z),V(y)\rangle V(y)\big)-(s-t)F(z)w\bigg\|\leq C\delta M_1\sqrt{\eps},
\end{equation*}
\begin{equation*}
 \bigg\|\frac{\langle x-y, V(z)+V(y)\rangle}{\langle V(y), V(z)+V(y)\rangle} \big(V(y)-\langle V(z),V(y)\rangle V(z)\big)-(s-t)F(z)w\bigg\|\leq C\delta M_1\sqrt{\eps},
\end{equation*}
By \eqref{BdVy} and \eqref{xyV},
\begin{align*}
    & \big\|x-z-\langle x-y, V(y)\rangle V(y)+(s-t)  F(z)w\big\|\nonumber\\
   =\,& \big\|(s-t)V(z)-\langle x-y,V(y)\rangle V(y)+(s-t) F(z)w \big\|\nonumber\\
   \leq\,& (s-t)\|V(z)-V(y)+F(z)w\|+\big\|((s-t)-\langle x-y,V(y)\rangle)V(y)\big\|\nonumber\\
   \leq\,& 2D\cdot \frac{1}{2} \delta M_1 \sqrt{\eps}+\big|(s-t)-\langle x-y, V(y)\rangle\big|\leq C\delta M_1\sqrt{\eps} . 
\end{align*}
Using \eqref{z1w1z2w2} and the above three displays, we get
\begin{align}\label{w1b}
    &\|w_1+w\|=\,\bigg\|x-z-\langle x-y,V(y)\rangle V(y)+\frac{\langle y-x, V(z)+V(y)\rangle}{\langle V(z), V(z)+V(y)\rangle} \big(V(z)-\langle V(z),V(y)\rangle V(y)\big)\bigg\|   \nonumber\\
    \leq\, & \big\|x-z-\langle x-y, V(y)\rangle V(y)+(s-t)  F(z)w\big\|\nonumber\\
    &+\bigg\|\frac{\langle y-x, V(z)+V(y)\rangle}{\langle V(z), V(z)+V(y)\rangle} \big(V(z)-\langle V(z),V(y)\rangle V(y)\big)-(s-t)F(z)w\bigg\|\leq C\delta M_1\sqrt{\eps},
\end{align}
\begin{align}\label{w2b}
 &\|w_2-(\mathbf{I}_d+(s-t)F(z))w\|\nonumber\\
 =\,&\bigg\|z-x+(s-t)V(z)+\frac{\langle x-y, V(z)+V(y)\rangle}{\langle V(y), V(z)+V(y)\rangle} \big(V(y)-\langle V(z),V(y)\rangle V(z)\big)-(s-t)F(z)w\bigg\|\nonumber\\
 =\,& \bigg\|\frac{\langle x-y, V(z)+V(y)\rangle}{\langle V(y), V(z)+V(y)\rangle} \big(V(y)-\langle V(z),V(y)\rangle V(z)\big)-(s-t)F(z)w\bigg\|\leq C\delta M_1\sqrt{\eps},
\end{align}
where in the second display we note \eqref{xyv}. By \eqref{w1b}, noting that $\delta$ is sufficiently small, we get 
\begin{equation}\label{w1ubb}
    \|w_1\|\leq \|w\|+\|w_1+w\| \leq \frac{1}{2}M_1\sqrt{\eps}+C\delta M_1\sqrt{\eps}\leq M_1\sqrt{\eps}.
\end{equation}
By \eqref{w2b} and \eqref{Fbdds},
\begin{align}\label{w2ubb}
    & \max\{\|w_2\|,\|(\mathbf{I}_d+(s-t)F(z))w\|\}\leq \|(\mathbf{I}_d+(s-t)F(z))w\|+C\delta M_1\sqrt{\eps}\nonumber\\
    \leq\,& (1+|s-t|\|F(z)\|_2)\|w\|+C\delta M_1\sqrt{\eps}\leq (1+2D\cdot 2\delta_1^{-1})\cdot\frac{1}{2}M_1\sqrt{\eps}+C\delta M_1\sqrt{\eps} \leq C\delta_1^{-1}M_1\sqrt{\eps}.
\end{align}

As $s-t\geq 2d_0$ and $\min\{s_1-t_1,s_2-t_2\}\geq d_0$ (recall \eqref{s1t1s2t2b}), by \eqref{s1t1}--\eqref{s2t21}, we have
\begin{equation}\label{E4.251}
    \bigg\|\frac{1}{s_1-t_1}\mathbf{I}_d-\frac{1}{s-t}\mathbf{I}_d\bigg\|_2=\bigg|\frac{1}{s_1-t_1}-\frac{1}{s-t}\bigg|\leq\frac{|s_1-s|+|t_1-t|}{(s_1-t_1)(s-t)}\leq CM_1\delta_1^{-1}\sqrt{\eps},
\end{equation}
\begin{equation}\label{E4.252}
    \bigg\|\frac{1}{s_2-t_2}\mathbf{I}_d-\frac{1}{s-t}\mathbf{I}_d\bigg\|_2=\bigg|\frac{1}{s_2-t_2}-\frac{1}{s-t}\bigg|\leq\frac{|s_2-s|+|t_2-t|}{(s_2-t_2)(s-t)}\leq CM_1\delta_1^{-1}\sqrt{\eps}.
\end{equation}
Moreover, by \eqref{Fbdds} we have
\begin{equation}\label{E4.253}
    \bigg\|\frac{1}{s-t} \mathbf{I}_d+F(z)\bigg\|_2\leq \frac{1}{s-t}+\|F(z)\|_2\leq C\delta_1^{-1}.
\end{equation}
By Theorem~\ref{L2.4} (see \eqref{L2.4.enew4}), noting \eqref{ldalpha}, \eqref{z1z2w1w2res}, and \eqref{s1t11}, we get 
\begin{align*}
&\|F(z_1)-F(q'+tV(q'))\|_2=\|F(q'+t_1V(q'))-F(q'+tV(q'))\|_2\nonumber\\
\leq\,&\frac{C|t_1-t|}{\min\{\alpha(z_1),\beta(z_1)\}\min\{\alpha(q'+tV(q')),\beta(q'+tV(q'))\}}
\leq C\delta_1^{-2}|t_1-t|\leq CM_1\delta_1^{-3}\sqrt{\eps}.
\end{align*}
As $w\in \mathscr{B}(q,s,t;\delta,M_1\sqrt{\eps}\slash 2)$ (recall Definition \ref{def4.9}), we have
\begin{equation*}
    \|F(q'+tV(q')) -F(z)\|_2=\|F(q'+tV(q')) -F(q+tV(q))\|_2\leq \delta.
\end{equation*}
Combining the above two displays, noting that $\eps$ is sufficiently small, we get
\begin{equation}\label{E4.254}
    \|F(z_1)-F(z)\|_2 \leq CM_1\delta_1^{-3}\sqrt{\eps}+\delta\leq 2\delta.
\end{equation}
By Theorem~\ref{L2.4} (see \eqref{L2.4.enew4}), noting \eqref{alpha_beta_delta1}, \eqref{z1z2w1w2res}, and \eqref{s2t2}, we get 
\begin{align}\label{E4.255}
&\|F(z_2)-F(x)\|_2=\|F(q+s_2V(q))-F(q+sV(q))\|_2\nonumber\\
\leq\,&\frac{C|s_2-s|}{\min\{\alpha(z_2),\beta(z_2)\}\min\{\alpha(x),\beta(x)\}}
\leq C\delta_1^{-2}|s_2-s|\leq CM_1\delta_1^{-3}\sqrt{\eps}.
\end{align}

By \eqref{w1b}, \eqref{w1ubb}, \eqref{E4.251}, \eqref{E4.253}, and \eqref{E4.254}, noting that $\eps$ is sufficiently small, we have 
\begin{align*}
    &\bigg|w_1^{\top} \bigg(\frac{1}{s_1-t_1}\mathbf{I}_d+F(z_1)\bigg) w_1-w^{\top}\bigg(\frac{1}{s-t} \mathbf{I}_d+F(z)\bigg)w\bigg|\nonumber\\
    \leq\,& \bigg|w_1^{\top}\bigg(\frac{1}{s_1-t_1}\mathbf{I}_d+F(z_1)\bigg)w_1-w_1^{\top}\bigg(\frac{1}{s-t}\mathbf{I}_d+F(z)\bigg)w_1\bigg|\nonumber\\
    &+\bigg|(w_1+w)^{\top}\bigg(\frac{1}{s-t}\mathbf{I}_d+F(z)\bigg)w_1\bigg|+\bigg|(-w)^{\top}\bigg(\frac{1}{s-t}\mathbf{I}_d+F(z)\bigg)(w_1+w)\bigg|\nonumber\\
    \leq\,& \bigg(\bigg\|\frac{1}{s_1-t_1}\mathbf{I}_d-\frac{1}{s-t}\mathbf{I}_d\bigg\|_2+\|F(z_1)-F(z)\|_2\bigg)\|w_1\|^2\nonumber\\
    &+2\max\{\|w_1\|,\|w\|\}\bigg\|\frac{1}{s-t}\mathbf{I}_d+F(z)\bigg\|_2\|w_1+w\|\nonumber\\
    \leq\,& C\delta\cdot  M_1^2\eps +CM_1\sqrt{\eps}\cdot \delta_1^{-1}\cdot \delta M_1 \sqrt{\eps}\leq C\delta M_1^2\delta_1^{-1}\eps.
\end{align*}
By \eqref{stid}, \eqref{E4.206}, \eqref{w2b}, \eqref{w2ubb}, \eqref{E4.252}, and \eqref{E4.255}, noting that $\eps$ is sufficiently small, we have
\begin{align*}
    &\bigg| w_2^{\top}\bigg(\frac{1}{s_2-t_2}\mathbf{I}_d-F(z_2)\bigg)w_2-w^{\top}\bigg(\frac{1}{s-t} \mathbf{I}_d+F(z)\bigg)w\bigg|\nonumber\\
    =  \,& \bigg|w_2^{\top}\bigg(\frac{1}{s_2-t_2}\mathbf{I}_d-F(z_2)\bigg)w_2\nonumber\\
    &\hspace{0.1in}-\big((\mathbf{I}_d+(s-t)F(z))w\big)^{\top}\bigg(\frac{1}{s-t}\mathbf{I}_d-F(x)\bigg)(\mathbf{I}_d+(s-t)F(z))w\bigg|\nonumber\\
    \leq\,& \bigg|w_2^{\top}\bigg(\frac{1}{s_2-t_2}\mathbf{I}_d-F(z_2)\bigg)w_2-w_2^{\top}\bigg(\frac{1}{s-t}\mathbf{I}_d-F(x)\bigg)w_2\bigg|\nonumber\\
    &+\bigg|\big(w_2-(\mathbf{I}_d+(s-t)F(z))w\big)^{\top}\bigg(\frac{1}{s-t}\mathbf{I}_d-F(x)\bigg)w_2\bigg|\nonumber\\
    &+\bigg|\big((\mathbf{I}_d+(s-t)F(z))w\big)^{\top}\bigg(\frac{1}{s-t}\mathbf{I}_d-F(x)\bigg)\big(w_2-(\mathbf{I}_d+(s-t)F(z))w\big)\bigg|\nonumber\\
    \leq\,& \bigg(\bigg\|\frac{1}{s_2-t_2}\mathbf{I}_d-\frac{1}{s-t}\mathbf{I}_d\bigg\|_2+\|F(z_2)-F(x)\|_2\bigg)\|w_2\|^2\nonumber\\
    &+2\max\{\|w_2\|,\|(\mathbf{I}_d+(s-t)F(z))w\|\}\bigg\|\frac{1}{s-t}\mathbf{I}_d-F(x)\bigg\|_2\|w_2-(\mathbf{I}_d+(s-t)F(z))w\|\nonumber\\
    \leq\,& C\delta\cdot \delta_1^{-2}M_1^2\eps +C\delta_1^{-1}M_1\sqrt{\eps}\cdot \delta_1^{-1}\cdot \delta M_1 \sqrt{\eps}\leq C\delta M_1^2\delta_1^{-2}\eps.\qedhere
\end{align*}
\end{proof}

Combining \eqref{E4.226}--\eqref{E4.227} and Claims~\ref{Claim.a4}--\ref{Claim.a6}, noting that $\eps$ is sufficiently small, we conclude that
\begin{align*}
    \eps^{(d-1)\slash 2} p_3(x,y)\geq\,& e^{-CM_1\delta_1^{-2}\delta_0^{-1}\sqrt{\eps}}\mathscr{I}_1(q,s,t)\cdot e^{-CM_1\delta_1^{-2}\sqrt{\eps}} (1+C\delta\delta_1^{-d})^{-2}\nonumber\\
    & \times e^{-C\delta M_1^2\delta_1^{-2}}e^{-\frac{1}{2\eps}w^{\top}\big(\frac{1}{s-t}\mathbf{I}_d+F(z)\big)w}\nonumber\\
    \geq\,& e^{-C\delta M_1^2\delta_1^{-2}}  (1+C\delta\delta_1^{-d})^{-2}\mathscr{I}_1(q,s,t) e^{-\frac{1}{2\eps}w^{\top}\big(\frac{1}{s-t}\mathbf{I}_d+F(z)\big)w},
\end{align*}
\begin{align*}
    \eps^{(d-1)\slash 2}p_4(x,y)\geq\,& e^{-CM_1\delta_1^{-2}\delta_0^{-1}\sqrt{\eps}}(\mathscr{I}_1(q,s,t)-\delta)_{+}\cdot e^{-CM_1\delta_1^{-2}\sqrt{\eps}} (1+C\delta\delta_1^{-d})^{-2}\nonumber\\
    &\times e^{-C\delta M_1^2\delta_1^{-2}}e^{-\frac{1}{2\eps}w^{\top}\big(\frac{1}{s-t}\mathbf{I}_d+F(z)\big)w}\nonumber\\
    \geq\,& e^{-C\delta M_1^2\delta_1^{-2}}  (1+C\delta\delta_1^{-d})^{-2}(\mathscr{I}_1(q,s,t)-\delta)_{+} e^{-\frac{1}{2\eps}w^{\top}\big(\frac{1}{s-t}\mathbf{I}_d+F(z)\big)w}. \qedhere
\end{align*}
This completes the proof of (b) and (c), and hence of \cref{Lem4.25}.

\subsection{Proof of the entropy estimate (Proposition \ref{P4.4})}\label{Sect.4.5}

In this subsection, we give the proof of Proposition \ref{P4.4}. Recall that we have fixed $k\in[K],k'\in [K']$. For any $i\in[I], j\in[J]$ and $\ell,\ell'\in [N_{0}]$, we define (recall \eqref{def.gamma.2} and note that $\lambda(\mathcal{S}\backslash\tilde{\mathcal{S}}^*)=0$ from Definition \ref{defstilde})
\begin{align}\label{def:M}
    \mathcal{M}_{i,j;\ell,\ell'}:=& \bar{\gamma}_{\eps,2}(L_{i,j;\ell}\times R_{i,j;\ell'})\nonumber\\
 =& \int_{\tilde{\mathcal{S}}^*}\mathbbm{1}_{\int_{\mathrm{int}(T)\cap \mathcal{X}_k^{\circ}}\bar{f}_2(z)\mathfrak{F}(z)d\mathcal{H}^1(z)>0}\nonumber\\
 &\hspace{0.2in}\times\frac{\int_{\mathrm{int}(T)\cap L_{i,j;\ell}}\bar{f}_2(z)\mathfrak{F}(z)d\mathcal{H}^1(z)\int_{\mathrm{int}(T)\cap R_{i,j;\ell'}}\bar{g}_2(z)\mathfrak{F}(z) d\mathcal{H}^1(z)}{\int_{\mathrm{int}(T)\cap \mathcal{X}_k^{\circ}}\bar{f}_2(z)\mathfrak{F}(z)d\mathcal{H}^1(z)}
 d\lambda(T).
\end{align}

We begin with an upper bound on $\int_{\mathbb{R}^d\times\mathbb{R}^d}r(x,y)\max\big\{\log\big(\eps^{(d-1)\slash 2}r(x,y)\big),0\big\}dxdy$. By the definition of $r(\cdot,\cdot)$ in~\eqref{def:r}, we have  
\begin{align*}
    &\int_{\mathbb{R}^d\times\mathbb{R}^d}r(x,y)\max\big\{\log\big(\eps^{(d-1)\slash 2}r(x,y)\big),0\big\}dxdy\nonumber\\
    =& \sum_{i=1}^{I}\sum_{j=1}^{J}\sum_{\substack{\ell,\ell'\in[N_0]:\\\mathcal{M}_{i,j;\ell,\ell'}>0}}\frac{\mathcal{M}_{i,j;\ell,\ell'}}{\bar{\gamma}_{\eps,2}(L_{i,j;\ell}\times R_j) \bar{\gamma}_{\eps,2}(L_i\times R_{i,j;\ell'})}\nonumber\\
    &\hspace{0.9in} \times\int_{L_{i,j;\ell}\times R_{i,j;\ell'}}f_j(x)g_i(y)\bigg(\log\bigg(\frac{\eps^{(d-1)\slash 2}\mathcal{M}_{i,j;\ell,\ell'}f_j(x)g_i(y)}{\bar{\gamma}_{\eps,2}(L_{i,j;\ell}\times R_j) \bar{\gamma}_{\eps,2}(L_i\times R_{i,j;\ell'})}\bigg)\bigg)_{+} dxdy\nonumber\\
    \leq& \sum_{i=1}^{I}\sum_{j=1}^{J}\sum_{\ell,\ell'\in[N_0]} \Big(\mathcal{M}_{i,j;\ell,\ell'}\big(\log\big(\eps^{-(d+3)\slash 2}\mathcal{M}_{i,j;\ell,\ell'}\big)\big)_{+}+\mathcal{M}_{i,j;\ell,\ell'}\big(\log\big(\eps^{-d\slash 2-1}\bar{\gamma}_{\eps,2}(L_{i,j;\ell}\times R_j)\big)\big)_{-}\nonumber\\
    &\hspace{1.2in}+\mathcal{M}_{i,j;\ell,\ell'}\big(\log\big(\eps^{-d\slash 2-1}\bar{\gamma}_{\eps,2}(L_i\times R_{i,j;\ell'})\big)\big)_{-}\Big)\nonumber\\
    &+\sum_{i=1}^{I}\sum_{j=1}^{J}\sum_{\substack{\ell,\ell'\in[N_0]:\\\mathcal{M}_{i,j;\ell,\ell'}>0}}\frac{\mathcal{M}_{i,j;\ell,\ell'}}{\bar{\gamma}_{\eps,2}(L_{i,j;\ell}\times R_j)}\int_{L_{i,j;\ell}}f_j(x)\big(\log\big(\eps^{-1\slash 2}f_j(x)\big)\big)_{+}dx\nonumber\\
    &+\sum_{i=1}^{I}\sum_{j=1}^{J}\sum_{\substack{\ell,\ell'\in[N_0]:\\\mathcal{M}_{i,j;\ell,\ell'}>0}}\frac{\mathcal{M}_{i,j;\ell,\ell'}}{\bar{\gamma}_{\eps,2}(L_{i}\times R_{i,j;\ell'})}\int_{R_{i,j;\ell'}}g_i(y)\big(\log\big(\eps^{-1\slash 2}g_i(y)\big)\big)_{+}dy,
\end{align*}
where for the inequality, we use that (by Lemma \ref{Lem4.9n})
\begin{equation*}
    \int_{L_{i,j;\ell}\times R_{i,j;\ell'}}f_j(x)g_i(y)dxdy=\bar{\gamma}_{\eps,2}(L_{i,j;\ell}\times R_j)\bar{\gamma}_{\eps,2}(L_{i}\times R_{i,j;\ell'})
\end{equation*}
and
\begin{align*}
   &\log\bigg(\frac{\eps^{(d-1)\slash 2}\mathcal{M}_{i,j;\ell,\ell'}f_j(x)g_i(y)}{\bar{\gamma}_{\eps,2}(L_{i,j;\ell}\times R_j) \bar{\gamma}_{\eps,2}(L_i\times R_{i,j;\ell'})}\bigg)\nonumber\\
   =\,& \log\big(\eps^{-(d+3)\slash 2}\mathcal{M}_{i,j;\ell,\ell'}\big)+\log\big(\eps^{-1\slash 2}f_j(x)\big)+\log\big(\eps^{-1\slash 2}g_i(y)\big)\nonumber\\
   & -\log\big(\eps^{-d\slash 2-1}\bar{\gamma}_{\eps,2}(L_{i,j;\ell}\times R_j)\big)-\log\big(\eps^{-d\slash 2-1}\bar{\gamma}_{\eps,2}(L_i\times R_{i,j;\ell'})\big).
\end{align*}
Rewriting the right-hand side of the initial display using~\eqref{def:M}, we arrive at
\begin{align}\label{Pro4.4.Eq1}
    &\int_{\mathbb{R}^d\times\mathbb{R}^d}r(x,y)\max\big\{\log\big(\eps^{(d-1)\slash 2}r(x,y)\big),0\big\}dxdy\nonumber\\
    \leq& \sum_{i=1}^{I}\sum_{j=1}^{J}\sum_{\ell=1}^{N_0}\sum_{\ell'=1}^{N_0}\mathcal{M}_{i,j;\ell,\ell'}\big(\log\big(\eps^{-(d+3)\slash 2}\mathcal{M}_{i,j;\ell,\ell'}\big)\big)_{+}\nonumber\\
    &+\sum_{i=1}^{I}\sum_{j=1}^{J}\sum_{\ell=1}^{N_0}\bar{\gamma}_{\eps,2}(L_{i,j;\ell}\times R_{j})\big(\log\big(\eps^{-d\slash 2-1}\bar{\gamma}_{\eps,2}(L_{i,j;\ell}\times R_{j})\big)\big)_{-}\nonumber\\
    &+\sum_{i=1}^{I}\sum_{j=1}^{J}\sum_{\ell'=1}^{N_0}\bar{\gamma}_{\eps,2}(L_i\times R_{i,j;\ell'})\big(\log\big(\eps^{-d\slash 2-1}\bar{\gamma}_{\eps,2}(L_i\times R_{i,j;\ell'})\big)\big)_{-}\nonumber\\
    &+\sum_{i=1}^I\sum_{j=1}^J \int_{L_i}f_j(x)\big(\log\big(\eps^{-1\slash 2}f_j(x)\big)\big)_{+}dx+\sum_{i=1}^I\sum_{j=1}^J \int_{R_j}g_i(y)\big(\log\big(\eps^{-1\slash 2}g_i(y)\big)\big)_{+}dy.
\end{align}

Next, we state and prove three lemmas bounding each term appearing on the right-hand side of \eqref{Pro4.4.Eq1}. 

\begin{lemma}\label{Lem4.24}
For any $K_0\in\mathbb{N}^*$ such that $2^{K_0}\geq 4d_0^{-1}$, we have 
\begin{align}\label{Lem4.28.1}
    & \sum_{i=1}^{I}\sum_{j=1}^{J}\sum_{\ell=1}^{N_0}\bar{\gamma}_{\eps,2}(L_{i,j;\ell}\times R_{j})\big(\log\big(\eps^{-d\slash 2-1}\bar{\gamma}_{\eps,2}(L_{i,j;\ell}\times R_{j})\big)\big)_{-}\nonumber\\
    \leq& C\big(\log(\delta_2^{-1})+K_0\big)\bar{\gamma}_{\eps,2}(\mathbb{R}^d\times\mathbb{R}^d)-\bar{\gamma}_{\eps,2}(\mathbb{R}^d\times\mathbb{R}^d)\log\big(\bar{\gamma}_{\eps,2}(\mathbb{R}^d\times\mathbb{R}^d)\big)+CK_02^{-K_0},
\end{align}
\begin{align}\label{Lem4.28.2}
    &\sum_{i=1}^{I}\sum_{j=1}^{J}\sum_{\ell=1}^{N_0}\bar{\gamma}_{\eps,2}(L_{i}\times R_{i,j;\ell})\big(\log\big(\eps^{-d\slash 2-1}\bar{\gamma}_{\eps,2}(L_{i}\times R_{i,j;\ell})\big)\big)_{-}\nonumber\\
    \leq& C\big(\log(\delta_2^{-1})+K_0\big)\bar{\gamma}_{\eps,2}(\mathbb{R}^d\times\mathbb{R}^d)-\bar{\gamma}_{\eps,2}(\mathbb{R}^d\times\mathbb{R}^d)\log\big(\bar{\gamma}_{\eps,2}(\mathbb{R}^d\times\mathbb{R}^d)\big)+CK_02^{-K_0}.
\end{align}
\end{lemma}

\begin{proof}

For any $i\in[I]$ and $j\in[J]$, we define
\begin{equation*}
    q_{i,j}^{(L)}:=\sum_{\ell=1}^{N_0}\min\big\{\bar{\gamma}_{\eps,2}(L_{i,j;\ell}\times R_{j}),\eps^{d\slash 2+1}\big\}, \qquad q_{i,j}^{(R)}:=\sum_{\ell=1}^{N_0}\min\big\{\bar{\gamma}_{\eps,2}(L_{i}\times R_{i,j; \ell}),\eps^{d\slash 2+1}\big\},
\end{equation*}
and set $q_{i,j}:=\bar{\gamma}_{\eps,2}(L_i\times R_j)$. Note that
\begin{equation}\label{E4.278}
  q_{i,j}^{(L)}\leq \sum_{\ell=1}^{N_0}\bar{\gamma}_{\eps,2}(L_{i,j;\ell}\times R_{j})=q_{i,j}, \qquad q_{i,j}^{(R)}\leq \sum_{\ell=1}^{N_0}\bar{\gamma}_{\eps,2}(L_{i}\times R_{i,j;\ell})=q_{i,j}.
\end{equation}
By Jensen's inequality and the convexity of $\Phi(\cdot)$ on $[0,\infty)$ (recall \eqref{Phi_def}), for any $i\in [I]$ and $j\in [J]$, 
\begin{align}\label{E4.279}
&-\sum_{\ell=1}^{N_0}\bar{\gamma}_{\eps,2}(L_{i,j;\ell}\times R_{j})\big(\log\big(\eps^{-d\slash 2-1}\bar{\gamma}_{\eps,2}(L_{i,j;\ell}\times R_{j})\big)\big)_{-}\nonumber\\
=& \eps^{d\slash 2+1}\sum_{\ell=1}^{N_0}\min\big\{\eps^{-d\slash 2-1}\bar{\gamma}_{\eps,2}(L_{i,j;\ell}\times R_{j}),1\big\}\log\big(\min\big\{\eps^{-d\slash 2-1}\bar{\gamma}_{\eps,2}(L_{i,j;\ell}\times R_{j}),1\big\}\big)\nonumber\\
=& \eps^{d\slash 2+1}N_0\cdot\frac{1}{N_0}\sum_{\ell=1}^{N_0}\Phi\Big(\min\big\{\eps^{-d\slash 2-1}\bar{\gamma}_{\eps,2}(L_{i,j;\ell}\times R_{j}),1\big\}\Big)\nonumber\\
\geq& \eps^{d\slash 2+1}N_0\Phi\bigg(\frac{1}{N_0}\sum_{\ell=1}^{N_0}\min\big\{\eps^{-d\slash 2-1}\bar{\gamma}_{\eps,2}(L_{i,j;\ell}\times R_{j}),1\big\}\bigg) \nonumber\\
=& q_{i,j}^{(L)}\log\bigg(\frac{\eps^{-d\slash 2-1}}{N_0}q_{i,j}^{(L)}\bigg)\geq q_{i,j}^{(L)}\log\Big(c\eps^{-(d+1)\slash 2}q_{i,j}^{(L)}\Big),
\end{align}
where we use \eqref{relation} in the last inequality. Similarly, we have
\begin{equation}
    -\sum_{\ell=1}^{N_0}\bar{\gamma}_{\eps,2}(L_{i}\times R_{i,j;\ell})\big(\log\big(\eps^{-d\slash 2-1}\bar{\gamma}_{\eps,2}(L_{i}\times R_{i,j;\ell})\big)\big)_{-}\geq q_{i,j}^{(R)}\log\Big(c\eps^{-(d+1)\slash 2} q_{i,j}^{(R)}\Big).
\end{equation}

Let $\mathcal{I}_0$ be the set of $i\in [I]$ such that $L_i\cap\mathcal{T}_1^{*}\neq \emptyset$. For any $i\in\mathcal{I}_0$, we take $x_i\in L_i\cap\mathcal{T}_1^{*}$ such that $\min\{\alpha(x_i),\beta(x_i)\}\geq\frac{1}{2}\sup\limits_{x\in L_i\cap\mathcal{T}_1^{*}}\min\{\alpha(x),\beta(x)\}$. For any $i\in\mathcal{I}_0$ and $t> 0$, we define 
\begin{equation*}
    \mathscr{D}_{i,t}:=\big\{x\in L_i\cap\mathcal{T}_1^{*}:\|V(x)-V(x_i)\|> t\sqrt{d \eps}\big\}, \quad \mathscr{F}_{i,t}:=\big\{x\in L_i\cap\mathcal{T}_1^{*}:\min\{\alpha(x),\beta(x)\}< t^{-1}\big\}.
\end{equation*}
If $x\in L_i\cap\mathcal{T}_1^{*}$ satisfies $\min\{\alpha(x),\beta(x)\}\geq t^{-1}$, then $\min\{\alpha(x_i),\beta(x_i)\}\geq \frac{1}{2}\min\{\alpha(x),\beta(x)\}\geq \frac{t^{-1}}{2}$, hence by Lemma \ref{L2.0} (note that $x,x_i\in L_i\cap\mathcal{T}_1^*$ and $\delta_2\in (0,1\slash 100)$), 
\begin{equation*}
    \|V(x)-V(x_i)\|\leq \frac{4\|x-x_i\|}{\min\{\alpha(x),\beta(x),\alpha(x_i),\beta(x_i)\}} \leq  \frac{4\delta_2\sqrt{d\eps}}{t^{-1}\slash 2}= 8t \delta_2\sqrt{d\eps}\leq t\sqrt{d\eps}.
\end{equation*}
Hence $\mathscr{D}_{i,t}\subseteq \mathscr{F}_{i,t}$. Thus, by Lemma \ref{Lem4.4} and the fact that $\{L_i\}_{i=1}^I$ are disjoint, for any $t>0$,
\begin{eqnarray}\label{neqq4.6}
   \sum_{i\in\mathcal{I}_0} \mathcal{L}^d(\mathscr{D}_{i,t})&\leq& \sum_{i\in\mathcal{I}_0}\mathcal{L}^d(\mathscr{F}_{i,t})\leq \mathcal{L}^d(\{x\in\mathcal{T}_1^{*}:\min\{\alpha(x),\beta(x)\}<t^{-1}\})\nonumber\\
    &\leq&\mathcal{L}^d(\{x\in\mathcal{T}_1^{*}:\alpha(x)<t^{-1}\})+\mathcal{L}^d(\{x\in\mathcal{T}_1^{*}:\beta(x)<t^{-1}\})\leq Ct^{-1}.
\end{eqnarray}

Recall the definition of $U_{i,j}$ for $i\in [I],j\in[J]$ from \eqref{defU}. For every $i\in \mathcal{I}_0$, we define
\begin{equation*}
    \mathcal{E}_{i,0}:=\big\{j\in [J]:\|U_{i,j}-V(x_i)\|\leq \sqrt{d\eps}\big\},
\end{equation*}
\begin{equation*}
    \mathcal{E}_{i,\ell}:=\big\{j\in[J]: 2^{\ell-1}\sqrt{d\eps}<\|U_{i,j}-V(x_i)\|\leq 2^{\ell}\sqrt{d\eps}\big\},\quad\text{ for all }\ell\in\mathbb{N}^{*}.
\end{equation*}

Fix any $i\in \mathcal{I}_0$, $\ell\in\mathbb{N}$, and $j\in\mathcal{E}_{i,\ell}$. For any $x\in L_i$ and $y\in R_j$, we have $\|x-y\|\geq d_0$ (as $\dist(L_i,R_j)\geq d_0$), $\|x-\mathtt{c}_i\|\leq \delta_2\sqrt{d\eps}$, and $\|y-\mathtt{c}'_j\|\leq \delta_2\sqrt{d\eps}$. Hence (note that $\delta_2\in (0,1\slash 100)$)
\begin{eqnarray}\label{bddUij}
   && \bigg\|\frac{x-y}{\|x-y\|}-U_{i,j}\bigg\|\leq \bigg\|\frac{x-y}{\|x-y\|}-\frac{\mathtt{c}_i-\mathtt{c}'_j}{\|x-y\|}\bigg\|+\bigg\|\frac{\mathtt{c}_i-\mathtt{c}'_j}{\|x-y\|}-\frac{\mathtt{c}_i-\mathtt{c}'_j}{\|\mathtt{c}_i-\mathtt{c}'_j\|}\bigg\|\nonumber\\
   &\leq& \frac{\|(x-y)-(\mathtt{c}_i-\mathtt{c}'_j)\|}{\|x-y\|}+\frac{|\|x-y\|-\|\mathtt{c}_i-\mathtt{c}'_j\||}{\|x-y\|}\nonumber\\
   &\leq& \frac{2(\|x-\mathtt{c}_i\|+\|y-\mathtt{c}'_j\|)}{\|x-y\|}
   \leq 4\delta_2d_0^{-1}\sqrt{d\eps}\leq d_0^{-1}\sqrt{d\eps}.
\end{eqnarray}
Therefore, using the fact that $\|U_{i,j}-V(x_i)\|\leq 2^{\ell}\sqrt{d\eps}$ (as $j\in\mathcal{E}_{i,\ell}$), we obtain
\begin{equation*}
    \bigg\|\frac{x-y}{\|x-y\|}-V(x_i)\bigg\|\leq \bigg\|\frac{x-y}{\|x-y\|}-U_{i,j}\bigg\|+\|U_{i,j}-V(x_i)\|\leq  (d_0^{-1}+2^{\ell})\sqrt{d\eps}.
\end{equation*}
Taking $x=x_i$ in the above display, we obtain that for any $y\in \bigcup_{j\in\mathcal{E}_{i,\ell}} R_j$,
\begin{equation*}
    \bigg\|\frac{x_i-y}{\|x_i-y\|}-V(x_i)\bigg\|\leq (d_0^{-1}+2^{\ell})\sqrt{d\eps}\leq C\cdot 2^{\ell}\sqrt{\eps}.
\end{equation*}
As $L_i,\bigcup_{j\in\mathcal{E}_{i,\ell}} R_j\subseteq B_d(0,D+1)$ (by \eqref{distcondi}), for any $i\in\mathcal{I}_0$ and $\ell\in\mathbb{N}$, using the above display, we get  $\mathcal{L}^d\big(\bigcup_{j\in\mathcal{E}_{i,\ell}} R_j\big)\leq C\cdot 2^{(d-1)\ell}\eps^{(d-1)\slash 2}$. As $\{R_j\}_{j=1}^J$ are disjoint and $\mathcal{L}^d(R_j)\geq (\delta_2\sqrt{\eps}\slash 2)^{d}$ for every $j\in [J]$, we have
\begin{equation}\label{bddEil}
    |\mathcal{E}_{i,\ell}|\leq \frac{\mathcal{L}^d\big(\bigcup_{j\in\mathcal{E}_{i,\ell}} R_j\big)}{(\delta_2\sqrt{\eps}\slash   2  )^{d}}\leq C\delta_2^{-d}2^{(d-1)\ell}\eps^{-1\slash 2},\quad \text{for every } i\in\mathcal{I}_0,\ell\in\mathbb{N}.  
\end{equation}

Now fix any $i\in\mathcal{I}_0$ and $\ell\in\mathbb{N}^{*}$. For any $x\in L_i\cap\mathcal{T}_1^{*}$, if there exists $j\in\mathcal{E}_{i,\ell}$ and $y\in R_j$ such that $y\in T(x)$ and $u(x)\geq u(y)$, then $\frac{x-y}{\|x-y\|}=V(x)$, and   
\begin{equation*}
    \|V(x)-V(x_i)\|=\bigg\|\frac{x-y}{\|x-y\|}-V(x_i)\bigg\|\geq \|U_{i,j}-V(x_i)\|-\bigg\|\frac{x-y}{\|x-y\|}-U_{i,j}\bigg\|> (2^{\ell-1}-d_0^{-1})\sqrt{d\eps},
\end{equation*}
where we use the definition of $\mathcal{E}_{i,\ell}$ and \eqref{bddUij} in the last inequality. Hence if $2^{\ell-1}>d_0^{-1}$, then  
\begin{equation}\label{neqq4.5}
    \bigg\{x\in L_i\cap\mathcal{T}_1^{*}:  \Big\{y\in \Big(\bigcup_{j\in\mathcal{E}_{i,\ell}} R_j\Big)\cap T(x):u(x)\geq u(y)\Big\}\neq\emptyset\bigg\}\subseteq \mathscr{D}_{i,2^{\ell-1}-d_0^{-1}}.
\end{equation}
Note that by \eqref{Eq.4.2new} and Lemma \ref{Ln4.1}, for $\mathcal{L}^d$-a.e.\ $x\in\mathbb{R}^d$, $\bar{f}_2(x)\leq f_{k,k'}(x)\leq f(x)$. Hence if $2^{\ell-1}>d_0^{-1}$, then by \eqref{def.gamma.2} and \Cref{Lem4.6}, we have 
\begin{align*}\label{neqq4.12}
   & \sum_{j\in\mathcal{E}_{i,\ell}}q_{i,j} = \sum_{j\in\mathcal{E}_{i,\ell}}\bar{\gamma}_{\eps,2}(L_i\times R_j)=\bar{\gamma}_{\eps,2}\Big(L_i\times \Big(\bigcup_{j\in\mathcal{E}_{i,\ell}}R_j\Big)\Big)\nonumber\\
   =& \bar{\gamma}_{\eps,2}\Big(\Big\{(x,y)\in\mathcal{T}_1^*\times\mathcal{T}_1^{*}:u(x)\geq u(y), x\in L_i,y\in \Big(\bigcup_{j\in\mathcal{E}_{i,\ell}} R_j\Big)\cap T(x)\Big\}\Big)\nonumber\\
   \leq& \bar{\gamma}_{\eps,2}\bigg(\bigg\{x\in L_i\cap\mathcal{T}_1^{*}:\Big\{y\in \Big(\bigcup_{j\in\mathcal{E}_{i,\ell}} R_j\Big)\cap T(x):u(x)\geq u(y)\Big\}\neq\emptyset\bigg\}\times \mathbb{R}^d\bigg)\nonumber\\
   =&\int_{\big\{x\in L_i\cap\mathcal{T}_1^{*}:\big\{y\in\big(\bigcup_{j\in\mathcal{E}_{i,\ell}}  R_j\big)\cap T(x):u(x)\geq u(y)\big\}\neq\emptyset\big\}}\bar{f}_2(x)dx\nonumber\\
   \leq&\int_{\big\{x\in L_i\cap\mathcal{T}_1^{*}:\big\{y\in\big(\bigcup_{j\in\mathcal{E}_{i,\ell}}  R_j\big)\cap T(x):u(x)\geq u(y)\big\}\neq\emptyset\big\}}f(x)dx\nonumber\\
   \leq& C\mathcal{L}^d\bigg(\bigg\{x\in L_i\cap\mathcal{T}_1^{*}:  \Big\{y\in \Big(\bigcup_{j\in\mathcal{E}_{i,\ell}} R_j\Big)\cap T(x):u(x)\geq u(y)\Big\}\neq\emptyset\bigg\}\bigg)\leq C \mathcal{L}^d\big(\mathscr{D}_{i,2^{\ell-1}-d_0^{-1}}\big), 
\end{align*}
where we use Lemma \ref{Lem4.9} in the fourth line and \eqref{neqq4.5} in the last inequality. By the above display and \eqref{neqq4.6}, if $2^{\ell-1}>d_0^{-1}$, then 
\begin{equation}\label{E4.285}
    \sum_{i\in\mathcal{I}_0}\sum_{j\in\mathcal{E}_{i,\ell}}q_{i,j}\leq C\sum_{i\in\mathcal{I}_0}\mathcal{L}^d\big(\mathscr{D}_{i,2^{\ell-1}-d_0^{-1}}\big)\leq \frac{C}{2^{\ell-1}-d_0^{-1}}.
\end{equation}

Now note that for any $i\in [I]\backslash\mathcal{I}_0$, $L_i\cap\mathcal{T}_1^{*}=\emptyset$. Hence by Lemma \ref{Lem4.9}, 
\begin{equation*}
    \bar{\gamma}_{\eps,2}(L_i\times\mathbb{R}^d)=
    \int_{L_i}\bar{f}_2(x)dx=\int_{L_i\cap\mathcal{X}}\bar{f}_2(x)dx=\int_{L_i\cap\mathcal{X}\cap\mathcal{T}_1^*}\bar{f}_2(x)dx=0,
\end{equation*}
where the second equality uses $\bar{f}_2(x)\leq f(x)$ for $\mathcal{L}^d$-a.e.\ $x\in\mathbb{R}^d$, and the third equality uses \eqref{new1.1}. Consequently, for any $i\in [I]\backslash\mathcal{I}_0$ and $j\in [J]$, $q_{i,j}^{(L)}\leq q_{i,j}\leq \bar{\gamma}_{\eps,2}(L_i\times\mathbb{R}^d)=0$ (recall \eqref{E4.278}). Hence
\begin{align}\label{neqq.4.21}
    \sum_{i=1}^I\sum_{j=1}^J q_{i,j}^{(L)}\log\Big(c\eps^{-(d+1)\slash 2}q_{i,j}^{(L)}\Big)=&\sum_{i\in \mathcal{I}_0} \sum_{j=1}^J q_{i,j}^{(L)}\log\Big(c\eps^{-(d+1)\slash 2}q_{i,j}^{(L)}\Big)\nonumber\\
    =& \sum_{\ell=0}^{\infty}\sum_{i\in \mathcal{I}_0} \sum_{j\in\mathcal{E}_{i,\ell}}
    q_{i,j}^{(L)}\log\Big(c\eps^{-(d+1)\slash 2}q_{i,j}^{(L)}\Big).
\end{align}

In the following, we fix any $K_0\in\mathbb{N}^{*}$ such that $2^{K_0}\geq 4d_0^{-1}$. For any $\ell\in\mathbb{N}^{*}$ such that $\ell\geq K_0+1$ and $\sum_{i\in \mathcal{I}_0}|\mathcal{E}_{i,\ell}|>0$, by Jensen's inequality, we have 
\begin{eqnarray*}
\sum_{i\in\mathcal{I}_0}\sum_{j\in\mathcal{E}_{i,\ell}}q_{i,j}^{(L)}\log\Big(c\eps^{-(d+1)\slash 2}q_{i,j}^{(L)}\Big)
&=&\bigg(c^{-1}\eps^{(d+1)\slash 2}\sum_{i\in \mathcal{I}_0}|\mathcal{E}_{i,\ell}|\bigg)\cdot \frac{\sum_{i\in \mathcal{I}_0}\sum_{j\in\mathcal{E}_{i,\ell}}\Phi\Big(c\eps^{-(d+1)\slash 2}q_{i,j}^{(L)}\Big)}{\sum_{i\in \mathcal{I}_0}|\mathcal{E}_{i,\ell}|}\nonumber\\
&\geq& \bigg(c^{-1}\eps^{(d+1)\slash 2}\sum_{i\in \mathcal{I}_0}|\mathcal{E}_{i,\ell}|\bigg)\cdot\Phi\Bigg(\frac{c\eps^{-(d+1)\slash 2}\sum_{i\in \mathcal{I}_0}\sum_{j\in\mathcal{E}_{i,\ell}}q_{i,j}^{(L)}}{\sum_{i\in \mathcal{I}_0}|\mathcal{E}_{i,\ell}|}\Bigg)\nonumber\\
&=& \bigg(\sum_{i\in \mathcal{I}_0}\sum_{j\in\mathcal{E}_{i,\ell}}q_{i,j}^{(L)}\bigg)\log\Bigg(\frac{c\eps^{-(d+1)\slash 2}\sum_{i\in \mathcal{I}_0}\sum_{j\in\mathcal{E}_{i,\ell}}q_{i,j}^{(L)}}{\sum_{i\in \mathcal{I}_0}|\mathcal{E}_{i,\ell}|}\Bigg).
\end{eqnarray*}
By \eqref{bddEil} and \eqref{relation},
\begin{equation*}
    \sum_{i\in \mathcal{I}_0}|\mathcal{E}_{i,\ell}|\leq C\delta_2^{-d}2^{(d-1)\ell}\eps^{-1\slash 2}I\leq C\delta_2^{-2d}2^{(d-1)\ell}\eps^{-(d+1)\slash 2}.
\end{equation*}
By the above two displays, we have
\begin{equation}\label{neqq4.15}
 \sum_{i\in\mathcal{I}_0}\sum_{j\in\mathcal{E}_{i,\ell}}q_{i,j}^{(L)}\log\Big(c\eps^{-(d+1)\slash 2}q_{i,j}^{(L)}\Big)\geq \bigg(\sum_{i\in \mathcal{I}_0}\sum_{j\in\mathcal{E}_{i,\ell}}q_{i,j}^{(L)}\bigg)\log\Bigg(\frac{\sum_{i\in \mathcal{I}_0}\sum_{j\in\mathcal{E}_{i,\ell}}q_{i,j}^{(L)}}{C\delta_2^{-2d}2^{(d-1)\ell}}\Bigg).
\end{equation}
By \eqref{E4.278} and \eqref{E4.285}, noting that $2^{\ell-1}\geq 2^{K_0}\geq 4d_0^{-1}$, we have
\begin{equation}\label{neqq4.16}
     \sum_{i\in \mathcal{I}_0}\sum_{j\in\mathcal{E}_{i,\ell}}q_{i,j}^{(L)} \leq \sum_{i\in \mathcal{I}_0}\sum_{j\in\mathcal{E}_{i,\ell}}q_{i,j} \leq \frac{C}{2^{\ell-1}-d_0^{-1}}\leq C_1\cdot 2^{-\ell}, 
\end{equation}
where $C_1$ is a positive constant that depends only on $d,\mu,\nu$. Note that for any $\alpha>0$, the function $\mathfrak{g}_{\alpha}(t):=t\log\big(\frac{t}{\alpha}\big)$ is monotonically decreasing on $[0,e^{-1}\alpha]$ and monotonically increasing on $[e^{-1}\alpha,\infty)$. Hence by \eqref{neqq4.15} and \eqref{neqq4.16}, 
\begin{align}\label{neqq4.22}
&\sum_{i\in\mathcal{I}_0}\sum_{j\in\mathcal{E}_{i,\ell}}q_{i,j}^{(L)}\log\Big(c\eps^{-(d+1)\slash 2}q_{i,j}^{(L)}\Big)\nonumber\\
\geq& -\log(e^{-1}CC_1^{-1}\delta_2^{-2d})\bigg(\sum_{i\in \mathcal{I}_0}\sum_{j\in\mathcal{E}_{i,\ell}}q_{i,j}^{(L)}\bigg)+\bigg(\sum_{i\in \mathcal{I}_0}\sum_{j\in\mathcal{E}_{i,\ell}}q_{i,j}^{(L)}\bigg)\log\Bigg(\frac{\sum_{i\in \mathcal{I}_0}\sum_{j\in\mathcal{E}_{i,\ell}}q_{i,j}^{(L)}}{eC_1\cdot 2^{(d-1)\ell}}\Bigg)\nonumber\\
\geq& (2d\log\delta_2-C)\bigg(\sum_{i\in \mathcal{I}_0}\sum_{j\in\mathcal{E}_{i,\ell}}q_{i,j}^{(L)}\bigg)+C_1\cdot 2^{-\ell}\log\big(e^{-1}2^{-d\ell}\big)\nonumber\\
\geq& (2d\log\delta_2-C)\bigg(\sum_{i\in \mathcal{I}_0}\sum_{j\in\mathcal{E}_{i,\ell}}q_{i,j}^{(L)}\bigg)-C\cdot\ell 2^{-\ell}.
\end{align}
Note that \eqref{neqq4.22} also holds if $\sum_{i\in \mathcal{I}_0}|\mathcal{E}_{i,\ell}|=0$.

By Jensen's inequality, if $\sum_{\ell=0}^{K_0}\sum_{i\in \mathcal{I}_0}|\mathcal{E}_{i,\ell}|>0$, then
\begin{align}\label{neqq4.10}
&\sum_{\ell=0}^{K_0}\sum_{i\in\mathcal{I}_0}\sum_{j\in\mathcal{E}_{i,\ell}}q_{i,j}^{(L)}\log\Big(c \eps^{-(d+1)\slash 2}q_{i,j}^{(L)}\Big)\nonumber\\
=&\bigg(c^{-1}\eps^{(d+1)\slash 2}\sum_{\ell=0}^{K_0}\sum_{i\in \mathcal{I}_0}|\mathcal{E}_{i,\ell}|\bigg)\cdot \frac{\sum_{\ell=0}^{K_0}\sum_{i\in \mathcal{I}_0}\sum_{j\in\mathcal{E}_{i,\ell}}\Phi\Big(c\eps^{-(d+1)\slash 2}q_{i,j}^{(L)}\Big)}{\sum_{\ell=0}^{K_0}\sum_{i\in \mathcal{I}_0}|\mathcal{E}_{i,\ell}|}\nonumber\\
\geq& \bigg(c^{-1}\eps^{(d+1)\slash 2}\sum_{\ell=0}^{K_0}\sum_{i\in \mathcal{I}_0}|\mathcal{E}_{i,\ell}|\bigg)\cdot\Phi\Bigg(\frac{c\eps^{-(d+1)\slash 2}\sum_{\ell=0}^{K_0}\sum_{i\in \mathcal{I}_0}\sum_{j\in\mathcal{E}_{i,\ell}}q_{i,j}^{(L)}}{\sum_{\ell=0}^{K_0}\sum_{i\in \mathcal{I}_0}|\mathcal{E}_{i,\ell}|}\Bigg)\nonumber\\
=& \bigg(\sum_{\ell=0}^{K_0}\sum_{i\in \mathcal{I}_0}\sum_{j\in\mathcal{E}_{i,\ell}}q_{i,j}^{(L)}\bigg)\log\Bigg(\frac{c\eps^{-(d+1)\slash 2}\sum_{\ell=0}^{K_0}\sum_{i\in \mathcal{I}_0}\sum_{j\in\mathcal{E}_{i,\ell}}q_{i,j}^{(L)}}{\sum_{\ell=0}^{K_0}\sum_{i\in \mathcal{I}_0}|\mathcal{E}_{i,\ell}|}\Bigg)\nonumber\\
\geq& \bigg(\sum_{\ell=0}^{K_0}\sum_{i\in \mathcal{I}_0}\sum_{j\in\mathcal{E}_{i,\ell}}q_{i,j}^{(L)}\bigg)\log\Bigg(\frac{\sum_{\ell=0}^{K_0}\sum_{i\in \mathcal{I}_0}\sum_{j\in\mathcal{E}_{i,\ell}}q_{i,j}^{(L)}}{C\delta_2^{-2d}2^{(d-1)K_0}}\Bigg),
\end{align}
where for the last inequality, we note that by \eqref{bddEil} and \eqref{relation},
\begin{equation*}
    \sum_{\ell=0}^{K_0}\sum_{i\in \mathcal{I}_0}|\mathcal{E}_{i,\ell}|\leq C\sum_{\ell=0}^{K_0}\sum_{i\in\mathcal{I}_0}\delta_2^{-d}2^{(d-1)\ell}\eps^{-1\slash 2}\leq C\delta_2^{-d}2^{(d-1)K_0}\eps^{-1\slash 2}I\leq C\delta_2^{-2d}2^{(d-1)K_0}\eps^{-(d+1)\slash 2}. 
\end{equation*}
Moreover, by \eqref{E4.278},
\begin{equation}\label{neqq4.9}
    \sum_{\ell=0}^{K_0}\sum_{i\in \mathcal{I}_0}\sum_{j\in\mathcal{E}_{i,\ell}}q_{i,j}^{(L)}\leq \sum_{i=1}^I\sum_{j=1}^J q_{i,j}^{(L)}\leq \sum_{i=1}^I\sum_{j=1}^J q_{i,j}\leq \bar{\gamma}_{\eps,2}(\mathbb{R}^d\times\mathbb{R}^d)\leq 1.
\end{equation}
By \eqref{neqq4.10} and \eqref{neqq4.9}, we have  
\begin{eqnarray}\label{neqq.4.23}
&&\sum_{\ell=0}^{K_0}\sum_{i\in\mathcal{I}_0}\sum_{j\in\mathcal{E}_{i,\ell}}q_{i,j}^{(L)}\log\Big(c\eps^{-(d+1)\slash 2}q_{i,j}^{(L)}\Big)\nonumber\\
&\geq& (2d\log\delta_2-C)
\bigg(\sum_{\ell=0}^{K_0}\sum_{i\in \mathcal{I}_0}\sum_{j\in\mathcal{E}_{i,\ell}}q_{i,j}^{(L)}\bigg)+\bigg(\sum_{\ell=0}^{K_0}\sum_{i\in \mathcal{I}_0}\sum_{j\in\mathcal{E}_{i,\ell}}q_{i,j}^{(L)}\bigg)\log\Bigg(\frac{\sum_{\ell=0}^{K_0}\sum_{i\in \mathcal{I}_0}\sum_{j\in\mathcal{E}_{i,\ell}}q_{i,j}^{(L)}}{e\cdot 2^{(d-1)K_0}}\Bigg)\nonumber\\
&\geq& (2d\log\delta_2-C)
\bar{\gamma}_{\eps,2}(\mathbb{R}^d\times\mathbb{R}^d)+\bar{\gamma}_{\eps,2}(\mathbb{R}^d\times\mathbb{R}^d)\log\Bigg(\frac{\bar{\gamma}_{\eps,2}(\mathbb{R}^d\times\mathbb{R}^d)}{e\cdot 2^{(d-1)K_0}}\Bigg)\nonumber\\
&\geq& (2d\log\delta_2-CK_0)\bar{\gamma}_{\eps,2}(\mathbb{R}^d\times\mathbb{R}^d)+\bar{\gamma}_{\eps,2}(\mathbb{R}^d\times\mathbb{R}^d)\log\big(\bar{\gamma}_{\eps,2}(\mathbb{R}^d\times\mathbb{R}^d)\big).
\end{eqnarray}
Note that \eqref{neqq.4.23} also holds if $\sum_{\ell=0}^{K_0}\sum_{i\in \mathcal{I}_0}|\mathcal{E}_{i,\ell}|=0$. 

Combining \eqref{E4.279}, \eqref{neqq.4.21}, \eqref{neqq4.22}, and \eqref{neqq.4.23}, we obtain that
\begin{align*}
    &-\sum_{i=1}^{I}\sum_{j=1}^{J}\sum_{\ell=1}^{N_0}\bar{\gamma}_{\eps,2}(L_{i,j;\ell}\times R_{j})\big(\log\big(\eps^{-d\slash 2-1}\bar{\gamma}_{\eps,2}(L_{i,j;\ell}\times R_{j})\big)\big)_{-}
    \nonumber\\
    \geq&\sum_{i=1}^I\sum_{j=1}^Jq_{i,j}^{(L)}\log\Big(c\eps^{-(d+1)\slash 2}q_{i,j}^{(L)}\Big)=\sum_{\ell=0}^{\infty}\sum_{i\in \mathcal{I}_0} \sum_{j\in\mathcal{E}_{i,\ell}}
    q_{i,j}^{(L)}\log\Big(c\eps^{-(d+1)\slash 2}q_{i,j}^{(L)}\Big)\nonumber\\
    \geq&(2d\log\delta_2-CK_0)\bar{\gamma}_{\eps,2}(\mathbb{R}^d\times\mathbb{R}^d)+\bar{\gamma}_{\eps,2}(\mathbb{R}^d\times\mathbb{R}^d)\log\big(\bar{\gamma}_{\eps,2}(\mathbb{R}^d\times\mathbb{R}^d)\big)\nonumber\\
    &+(2d\log\delta_2-C)\bigg(\sum_{\ell=K_0+1}^{\infty}\sum_{i\in \mathcal{I}_0}\sum_{j\in\mathcal{E}_{i,\ell}}q_{i,j}^{(L)}\bigg)-C\sum_{\ell=K_0+1}^{\infty}  \ell 2^{-\ell}\nonumber\\
    \geq& (4d\log\delta_2-CK_0)\bar{\gamma}_{\eps,2}(\mathbb{R}^d\times\mathbb{R}^d)+\bar{\gamma}_{\eps,2}(\mathbb{R}^d\times\mathbb{R}^d)\log\big(\bar{\gamma}_{\eps,2}(\mathbb{R}^d\times\mathbb{R}^d)\big)-CK_02^{-K_0},
\end{align*}
where we use $\sum_{\ell=K_0+1}^{\infty}\sum_{i\in \mathcal{I}_0}\sum_{j\in\mathcal{E}_{i,\ell}}q_{i,j}^{(L)}\leq \sum_{i=1}^I\sum_{j=1}^Jq_{i,j}^{(L)}\leq \sum_{i=1}^I\sum_{j=1}^Jq_{i,j}\leq \bar{\gamma}_{\eps,2}(\mathbb{R}^d\times\mathbb{R}^d)$ in the last inequality. This establishes \eqref{Lem4.28.1}. The bound \eqref{Lem4.28.2} can be similarly deduced.
\end{proof}

\begin{lemma}\label{Lem4.26}
For any $\delta\in(0,e^{-1}]$, we have
\begin{align*}
   & \sum_{i=1}^I\sum_{j=1}^J
    \sum_{\ell=1}^{N_0}\sum_{\ell'=1}^{N_0}\mathcal{M}_{i,j;\ell,\ell'}\big(\log\big(\eps^{-(d+3)\slash 2}\mathcal{M}_{i,j;\ell,\ell'}\big)\big)_{+}\nonumber\\
    \leq& \big(C\log(\delta^{-1})+\log(\delta_3^{-1})\big)\bar{\gamma}_{\eps,2}(\mathbb{R}^d\times\mathbb{R}^d)+C\delta\log(\delta^{-1}).
\end{align*}
\end{lemma}
\begin{proof}

For any $i\in[I],j\in[J]$, we define $\mathfrak{S}_{i,j}:=\{T\in\tilde{\mathcal{S}}^*: T\cap L_i\neq\emptyset,T\cap R_j\neq\emptyset\}$. For any $i\in[I],j\in[J]$, $\ell,\ell'\in[N_0]$, and $T\in\tilde{\mathcal{S}}^*$, we define
\begin{equation}\label{def:mijl}
    \mathfrak{m}_{i,j;\ell,\ell'}(T):=\mathbbm{1}_{\int_{\mathrm{int}(T)\cap \mathcal{X}_k^{\circ}}\bar{f}_2(z)\mathfrak{F}(z)d\mathcal{H}^1(z)>0}\cdot\frac{\int_{\mathrm{int}(T)\cap L_{i,j;\ell}}\bar{f}_2(z)\mathfrak{F}(z)d\mathcal{H}^1(z)\int_{\mathrm{int}(T)\cap R_{i,j;\ell'}}\bar{g}_2(z)\mathfrak{F}(z) d\mathcal{H}^1(z)}{\int_{\mathrm{int}(T)\cap \mathcal{X}_k^{\circ}}\bar{f}_2(z)\mathfrak{F}(z)d\mathcal{H}^1(z)}.
\end{equation}
Note that by \eqref{def:M}, $\mathcal{M}_{i,j;\ell,\ell'}=\int_{\mathfrak{S}_{i,j}}\mathfrak{m}_{i,j;\ell,\ell'}(T)d\lambda(T)$. Hence by Jensen's inequality, with $\Phi(\cdot)$ as in \eqref{Phi_def}, we have
\begin{align}\label{Eq4.294}
 & \sum_{i=1}^I\sum_{j=1}^J
    \sum_{\ell=1}^{N_0}\sum_{\ell'=1}^{N_0}\mathcal{M}_{i,j;\ell,\ell'}\big(\log\big(\eps^{-(d+3)\slash 2}\mathcal{M}_{i,j;\ell,\ell'}\big)\big)_{+}\nonumber\\
    =& \sum_{i=1}^I\sum_{j=1}^J
    \sum_{\ell=1}^{N_0}\sum_{\ell'=1}^{N_0} \bigg(\int_{\mathfrak{S}_{i,j}}\mathfrak{m}_{i,j;\ell,\ell'}(T)d\lambda(T)\bigg)\bigg(\log\bigg(\eps^{-(d+3)\slash 2}\int_{\mathfrak{S}_{i,j}}\mathfrak{m}_{i,j;\ell,\ell'}(T)d\lambda(T)\bigg)\bigg)_{+}   \nonumber\\
    =& \sum_{i=1}^I\sum_{j=1}^J
    \sum_{\substack{\ell,\ell'\in[N_0]:\\\int_{\mathfrak{S}_{i,j}}\mathfrak{m}_{i,j;\ell,\ell'}(T)d\lambda(T)\geq \eps^{(d+3)\slash 2}}}\bigg(\int_{\mathfrak{S}_{i,j}}\mathfrak{m}_{i,j;\ell,\ell'}(T)d\lambda(T)\bigg)\log\bigg(\eps^{-(d+3)\slash 2}\int_{\mathfrak{S}_{i,j}}\mathfrak{m}_{i,j;\ell,\ell'}(T)d\lambda(T)\bigg)\nonumber\\
    =& \eps^{(d+3)\slash 2}\sum_{i=1}^I\sum_{j=1}^J
    \sum_{\substack{\ell,\ell'\in[N_0]:\\\int_{\mathfrak{S}_{i,j}}\mathfrak{m}_{i,j;\ell,\ell'}(T)d\lambda(T)\geq \eps^{(d+3)\slash 2}}}\Phi\bigg(\frac{1}{\lambda(\mathfrak{S}_{i,j})}\int_{\mathfrak{S}_{i,j}}\eps^{-(d+3)\slash 2}\lambda(\mathfrak{S}_{i,j})\mathfrak{m}_{i,j;\ell,\ell'}(T)d\lambda(T)\bigg)\nonumber\\
    \leq& \eps^{(d+3)\slash 2}\sum_{i=1}^I\sum_{j=1}^J
    \sum_{\substack{\ell,\ell'\in[N_0]:\\\int_{\mathfrak{S}_{i,j}}\mathfrak{m}_{i,j;\ell,\ell'}(T)d\lambda(T)\geq \eps^{(d+3)\slash 2}}}\frac{1}{\lambda(\mathfrak{S}_{i,j})}\int_{\mathfrak{S}_{i,j}}\Phi\Big(\eps^{-(d+3)\slash 2}\lambda(\mathfrak{S}_{i,j})\mathfrak{m}_{i,j;\ell,\ell'}(T)\Big)d\lambda(T)\nonumber\\
    \leq& \sum_{i=1}^I\sum_{j=1}^J
    \sum_{\ell,\ell'\in[N_0]}\int_{\mathfrak{S}_{i,j}}\mathfrak{m}_{i,j;\ell,\ell'}(T)\Big(\log\Big(\eps^{-(d+3)\slash 2}\lambda(\mathfrak{S}_{i,j})\mathfrak{m}_{i,j;\ell,\ell'}(T)\Big)\Big)_{+}d\lambda(T)\nonumber\\
    =& \sum_{i=1}^I\sum_{j=1}^J\sum_{\ell,\ell'\in [N_0]}\int_{\mathfrak{H}_{k,k'}\cap T^{-1}(\mathfrak{S}_{i,j})}|\langle V(q),\mathtt{a}_{k,k'}\rangle|\mathfrak{m}_{i,j;\ell,\ell'}(T(q))\Big(\log\Big(\eps^{-(d+3)\slash 2}\lambda(\mathfrak{S}_{i,j})\mathfrak{m}_{i,j;\ell,\ell'}(T(q))\Big)\Big)_{+}\nonumber\\
&\hspace{1.64in}\times\bigg(\int_{T(q)}f(x)\det\big(\mathbf{I}_d+\langle x-q, V(q)\rangle F(q)\big)d\mathcal{H}^1(x)\bigg)d\mathcal{H}^{d-1}(q),
\end{align}
where the last equality uses Lemma \ref{L3.14n} (note that by \eqref{deftildeTkk}, $T^{-1}(\mathfrak{S}_{i,j})\subseteq \tilde{\mathcal{T}}_{1;k,k'}^*$). 

Below, we consider any $i\in [I],j\in[J]$, $\ell,\ell'\in[N_0]$, and $T\in\mathfrak{S}_{i,j}$. Note that $\mathrm{int}(T)\subseteq \mathfrak{T}^*\cap\tilde{\mathcal{T}}_{1;k,k'}^*$. Let $q\in \tilde{H}_{k,k'}$ be the intersection point of $T$ and $H_{k,k'}$. By \eqref{fkkdef} and \eqref{gkkdef}, we have
\begin{align*}
   & \int_{\mathrm{int}(T)\cap \mathcal{X}_k^{\circ}}f_{k,k'}(z)\mathfrak{F}(z)d\mathcal{H}^1(z)=\int_{(\mathrm{int}(T)\cap\mathcal{X}_k^{\circ})\times(\mathrm{int}(T)\cap\mathcal{Y}_{k'}^{\circ})}h(x,y)d\mathcal{H}^1(x)d\mathcal{H}^1(y),\nonumber\\
   & \int_{\mathrm{int}(T)\cap L_{i,j;\ell}}f_{k,k'}(z)\mathfrak{F}(z)d\mathcal{H}^1(z)=\int_{(\mathrm{int}(T)\cap L_{i,j;\ell})\times(\mathrm{int}(T)\cap\mathcal{Y}_{k'}^{\circ})} h(x,y)d\mathcal{H}^1(x)d\mathcal{H}^1(y),\nonumber\\
   & \int_{\mathrm{int}(T)\cap R_{i,j;\ell'}}g_{k,k'}(z)\mathfrak{F}(z)d\mathcal{H}^1(z)=\int_{(\mathrm{int}(T)\cap\mathcal{X}_k^{\circ})\times(\mathrm{int}(T)\cap R_{i,j;\ell'})}h(x,y)d\mathcal{H}^1(x)d\mathcal{H}^1(y).
\end{align*}
Consequently, by Lemma \ref{Lem4.7} and Definition \ref{DefhH},
\begin{align}\label{Lem4.29.E1}
   & \mathfrak{m}_{i,j;\ell,\ell'}(T)\nonumber\\
  \leq & \delta_3^{-1}\cdot\mathbbm{1}_{\int_{\mathrm{int}(T)\cap \mathcal{X}_k^{\circ}}f_{k,k'}(z)\mathfrak{F}(z)d\mathcal{H}^1(z)>0}\cdot\frac{\int_{\mathrm{int}(T)\cap L_{i,j;\ell}}f_{k,k'}(z)\mathfrak{F}(z)d\mathcal{H}^1(z)\int_{\mathrm{int}(T)\cap R_{i,j;\ell'}}g_{k,k'}(z)\mathfrak{F}(z)d\mathcal{H}^1(z)}{\int_{\mathrm{int}(T)\cap \mathcal{X}_k^{\circ}}f_{k,k'}(z)\mathfrak{F}(z)d\mathcal{H}^1(z)}\nonumber\\
   =&\delta_3^{-1}\cdot\mathbbm{1}_{\int_{(\mathrm{int}(T)\cap\mathcal{X}_k^{\circ})\times(\mathrm{int}(T)\cap\mathcal{Y}_{k'}^{\circ})}h_T(x,y)d\mathcal{H}^1(x)d\mathcal{H}^1(y)>0}\nonumber\\
   &\times \frac{\int_{(\mathrm{int}(T)\cap L_{i,j;\ell})\times(\mathrm{int}(T)\cap\mathcal{Y}_{k'}^{\circ})} h_T(x,y)d\mathcal{H}^1(x)d\mathcal{H}^1(y)\int_{(\mathrm{int}(T)\cap\mathcal{X}_k^{\circ})\times(\mathrm{int}(T)\cap R_{i,j;\ell'})}h_T(x,y)d\mathcal{H}^1(x)d\mathcal{H}^1(y)}{\int_{(\mathrm{int}(T)\cap\mathcal{X}_k^{\circ})\times(\mathrm{int}(T)\cap\mathcal{Y}_{k'}^{\circ})}h_T(x,y)d\mathcal{H}^1(x)d\mathcal{H}^1(y)}\nonumber\\
   \leq&\delta_3^{-1}\cdot \mathbbm{1}_{\int_{(\mathrm{int}(T)\cap\mathcal{X}_k^{\circ})\times(\mathrm{int}(T)\cap\mathcal{Y}_{k'}^{\circ})}h_T(x,y)d\mathcal{H}^1(x)d\mathcal{H}^1(y)>0}\nonumber\\
   &\times \frac{\int_{(\mathrm{int}(T)\cap L_{i,j;\ell})\times(\mathrm{int}(T)\cap\mathcal{Y}_{k'}^{\circ})} h_T(x,y)d\mathcal{H}^1(x)d\mathcal{H}^1(y)\int_{(\mathrm{int}(T)\cap\mathcal{X}_k^{\circ})\times(\mathrm{int}(T)\cap R_{i,j;\ell'})}h_T(x,y)d\mathcal{H}^1(x)d\mathcal{H}^1(y)}{\big(\int_{(\mathrm{int}(T)\cap\mathcal{X}_k^{\circ})\times(\mathrm{int}(T)\cap\mathcal{Y}_{k'}^{\circ})}h_T(x,y)d\mathcal{H}^1(x)d\mathcal{H}^1(y)\big)^2},
\end{align}
where we use $\int_{(\mathrm{int}(T)\cap\mathcal{X}_k^{\circ})\times(\mathrm{int}(T)\cap\mathcal{Y}_{k'}^{\circ})}h_T(x,y)d\mathcal{H}^1(x)d\mathcal{H}^1(y)\leq \int_{T\times T}h_T(x,y)d\mathcal{H}^1(x)d\mathcal{H}^1(y)\leq 1$ in the last inequality. Note that by \Cref{Lem2}(d) and \eqref{bddsalphabeta}, for any $x,x_0\in \mathrm{int}(T)\cap \mathcal{X}_k^{\circ}$ and $y,y_0\in \mathrm{int}(T)\cap\mathcal{Y}_{k'}^{\circ}$, we have $c\tilde{h}_T(x_0,y_0)\leq \tilde{h}_T(x,y)\leq C\tilde{h}_T(x_0,y_0)$. Hence by \eqref{Lem4.29.E1}, 
\begin{align}\label{Lem4.29.NE1}
     & \mathfrak{m}_{i,j;\ell,\ell'}(T)\nonumber\\
     \leq& C\delta_3^{-1}\cdot \mathbbm{1}_{\int_{\mathrm{int}(T)\cap\mathcal{X}_k^{\circ}}\tilde{f}(x)d\mathcal{H}^1(x)>0}\mathbbm{1}_{\int_{\mathrm{int}(T)\cap\mathcal{Y}_{k'}^{\circ}}\tilde{g}(y)d\mathcal{H}^1(y)>0}\nonumber\\
     &\times \frac{\int_{(\mathrm{int}(T)\cap L_{i,j;\ell})\times(\mathrm{int}(T)\cap\mathcal{Y}_{k'}^{\circ})} \tilde{f}(x)\tilde{g}(y)d\mathcal{H}^1(x)d\mathcal{H}^1(y)\int_{(\mathrm{int}(T)\cap\mathcal{X}_k^{\circ})\times(\mathrm{int}(T)\cap R_{i,j;\ell'})}\tilde{f}(x)\tilde{g}(y)d\mathcal{H}^1(x)d\mathcal{H}^1(y)}{\big(\int_{(\mathrm{int}(T)\cap\mathcal{X}_k^{\circ})\times(\mathrm{int}(T)\cap\mathcal{Y}_{k'}^{\circ})}\tilde{f}(x)\tilde{g}(y)d\mathcal{H}^1(x)d\mathcal{H}^1(y)\big)^2}\nonumber\\
     =& C\delta_3^{-1}\cdot \mathbbm{1}_{\int_{\mathrm{int}(T)\cap\mathcal{X}_k^{\circ}}\tilde{f}(x)d\mathcal{H}^1(x)>0}\mathbbm{1}_{\int_{\mathrm{int}(T)\cap\mathcal{Y}_{k'}^{\circ}}\tilde{g}(x)d\mathcal{H}^1(x)>0}\nonumber\\
     &\times\frac{\int_{\mathrm{int}(T)\cap L_{i,j;\ell}}\tilde{f}(x)d\mathcal{H}^1(x)\int_{\mathrm{int}(T)\cap R_{i,j;\ell'}}\tilde{g}(x)d\mathcal{H}^1(x)}{\int_{\mathrm{int}(T)\cap\mathcal{X}_k^{\circ}}\tilde{f}(x)d\mathcal{H}^1(x)\int_{\mathrm{int}(T)\cap\mathcal{Y}_{k'}^{\circ}}\tilde{g}(x)d\mathcal{H}^1(x)}\nonumber\\
     =& C\delta_3^{-1} \cdot  \mathbbm{1}_{\int_{T(q)\cap\mathcal{X}_k^{\circ}}f(x)\det(\mathbf{I}_d+\langle x-q, V(q)\rangle F(q))d\mathcal{H}^1(x)>0}\mathbbm{1}_{\int_{T(q)\cap\mathcal{Y}_{k'}^{\circ}}g(x)\det(\mathbf{I}_d+\langle x-q, V(q)\rangle F(q))d\mathcal{H}^1(x)>0}\nonumber\\
    &\times \frac{\int_{T(q)\cap L_{i,j;\ell}}f(x)\det(\mathbf{I}_d+\langle x-q, V(q)\rangle F(q))d\mathcal{H}^1(x)}{\int_{T(q)\cap\mathcal{X}_k^{\circ}}f(x)\det(\mathbf{I}_d+\langle x-q, V(q)\rangle F(q))d\mathcal{H}^1(x)}\nonumber\\
    &\times  \frac{\int_{T(q)\cap R_{i,j;\ell'}}g(x)\det(\mathbf{I}_d+\langle x-q, V(q)\rangle F(q))d\mathcal{H}^1(x)}{\int_{T(q)\cap\mathcal{Y}_{k'}^{\circ}}g(x)\det(\mathbf{I}_d+\langle x-q, V(q)\rangle F(q))d\mathcal{H}^1(x)},
\end{align}
where the last equality uses Definition \ref{factors} and \eqref{defft}. If $\mathrm{int}(T)\cap   L_{i,j;\ell}=\emptyset$ or $\mathrm{int}(T)\cap R_{i,j;\ell'}=\emptyset$, then by \eqref{Lem4.29.NE1} (note that $T=T(q)$), $\mathfrak{m}_{i,j;\ell,\ell'}(T)=0$. Below we assume that $\mathrm{int}(T)\cap   L_{i,j;\ell}\neq\emptyset$ and $\mathrm{int}(T)\cap R_{i,j;\ell'}\neq\emptyset$. Suppose that $x_0\in \mathrm{int}(T)\cap L_{i,j;\ell}$ and $y_0\in \mathrm{int}(T)\cap R_{i,j;\ell'}$. If $u(x_0)\leq u(y_0)$, then $u(y_0)-u(x_0)=\|x_0-y_0\|\geq d_0$ (as $\dist(L_i,R_j)\geq\dist(\mathcal{X}_k,\mathcal{Y}_{k'})\geq d_0$), and so for any $x\in\mathrm{int}(T)\cap\mathcal{X}_k^{\circ}$ and $y\in\mathrm{int}(T)\cap\mathcal{Y}_{k'}^{\circ}$,
\begin{equation*}
    u(y)-u(x)\geq u(y_0)-u(x_0)-|u(x)-u(x_0)|-|u(y)-u(y_0)|\geq d_0-\|x-x_0\|-\|y-y_0\|>0,
\end{equation*}
where the last inequality uses $\max\{\|x-x_0\|,\|y-y_0\|\}\leq 2\sqrt{d}\mathtt{h}\leq d_0\slash 5$ (recall \eqref{eq:mathtthDefn}); consequently, by \Cref{Lem2}(c), $h_T(x,y)=0$ for all $x\in\mathrm{int}(T)\cap\mathcal{X}_k^{\circ},y\in\mathrm{int}(T)\cap\mathcal{Y}_{k'}^{\circ}$, and so $\mathfrak{m}_{i,j;\ell,\ell'}(T)=0$ (by \eqref{Lem4.29.E1}). Below we further assume that $u(x_0)> u(y_0)$, which implies $V(q)=\frac{x_0-y_0}{\|x_0-y_0\|}$. Thus, with $U_{i,j}$ as in \eqref{defU}, we have
\begin{align}\label{eq.neq.4.1}
    \|V(q)-U_{i,j}\|=&\bigg\|\frac{x_0-y_0}{\|x_0-y_0\|}-\frac{\mathtt{c}_i-\mathtt{c}_j'}{\|\mathtt{c}_i-\mathtt{c}_j'\|}\bigg\|\leq \frac{\|(x_0-y_0)-(\mathtt{c}_i-\mathtt{c}'_j)\|}{\|x_0-y_0\|}+\frac{\big|\|x_0-y_0\|-\|\mathtt{c}_i-\mathtt{c}'_j\|\big|}{\|x_0-y_0\|}\nonumber\\
    \leq& \frac{2(\|x_0-\mathtt{c}_i\|+\|y_0-\mathtt{c}'_j\|)}{\|x_0-y_0\|}\leq \frac{4\delta_2\sqrt{d\eps}}{d_0}\leq C\delta_2\sqrt{\eps},
\end{align}
where the second to last inequality uses $\max\{\|x_0-\mathtt{c}_i\|,\|y_0-\mathtt{c}'_j\|\}\leq\delta_2\sqrt{d\eps}$ (as $x_0,\mathtt{c}_i\in L_i$ and $y_0,\mathtt{c}_j'\in R_j$) and $\min\{\|x_0-y_0\|,\|\mathtt{c}_i-\mathtt{c}'_j\|\}\geq d_0$. By the definitions of $L_{i,j;\ell}$ and $R_{i,j;\ell'}$ as in \eqref{Lijk} and \eqref{Rijk}, for any $x\in T(q)\cap L_{i,j;\ell}$ and $y\in T(q)\cap R_{i,j;\ell'}$, we have $|\langle U_{i,j}, x-x_0\rangle|\leq \delta_2\eps$ and $|\langle U_{i,j}, y-y_0\rangle|\leq \delta_2\eps$, hence by \eqref{eq.neq.4.1},
\begin{eqnarray}\label{eq.neq.4.2}
     |\langle V(q), x-x_0\rangle| &\leq&  |\langle U_{i,j}, x-x_0\rangle|+|\langle V(q)-U_{i,j}, x-x_0\rangle|\leq \delta_2\eps+\|V(q)-U_{i,j}\|\|x-x_0\| \nonumber\\
     &\leq& \delta_2\eps+C\delta_2\sqrt{\eps}\cdot \delta_2\sqrt{d\eps}\leq C\delta_2\eps,
\end{eqnarray}
\begin{eqnarray}\label{eq.new.4.3}
    |\langle V(q), y-y_0\rangle| &\leq&  |\langle U_{i,j}, y-y_0\rangle|+|\langle V(q)-U_{i,j}, y-y_0\rangle|\leq \delta_2\eps+\|V(q)-U_{i,j}\|\|y-y_0\| \nonumber\\
     &\leq& \delta_2\eps+C\delta_2\sqrt{\eps}\cdot \delta_2\sqrt{d\eps}\leq C\delta_2\eps,
\end{eqnarray}
where we use the fact that $\max\{\|x-x_0\|,\|y-y_0\|\}\leq \delta_2\sqrt{d\eps}$ (as $x,x_0\in L_i$ and $y,y_0\in R_j$). By \eqref{eq.neq.4.2} and \eqref{eq.new.4.3},
\begin{equation}\label{eq.new.4.0}
    \mathcal{H}^1(T(q)\cap L_{i,j;\ell})\leq C\delta_2\eps, \quad \mathcal{H}^1(T(q)\cap R_{i,j;\ell'})\leq C\delta_2\eps.
\end{equation}
For any $x\in T(q)$, as $\min\{\alpha(q),\beta(q)\}\geq d_0$ (recall that $q\in \tilde{H}_{k,k'}$ and \eqref{deftildehat}), by Hadamard's inequality (see, e.g., \cite[Lemma 2.5]{ipsen2008perturbation}), \eqref{bddsalphabeta}, and \Cref{L2.4} (see \eqref{L2.4.enew3}), we have
\begin{equation}\label{eq.new.4.1n}
    \big|\det\big(\mathbf{I}_d+\langle x-q, V(q)\rangle F(q)\big)\big|\leq \big\|\mathbf{I}_d+\langle x-q, V(q)\rangle F(q)\big\|_2^d\leq (1+2D\|F(q)\|_2)^d\leq C.
\end{equation}
By \eqref{Lem4.29.NE1} and \eqref{eq.new.4.0}--\eqref{eq.new.4.1n}, we have
\begin{align}\label{Eq4.302}
    \mathfrak{m}_{i,j;\ell,\ell'}(T)\leq& C\delta_2^2\delta_3^{-1}\eps^2\cdot \frac{\mathbbm{1}_{\int_{T(q)\cap\mathcal{X}_k^{\circ}}f(x)\det(\mathbf{I}_d+\langle x-q, V(q)\rangle F(q))d\mathcal{H}^1(x)>0}}{\int_{T(q)\cap\mathcal{X}_k^{\circ}}f(x)\det(\mathbf{I}_d+\langle x-q, V(q)\rangle F(q))d\mathcal{H}^1(x)}\nonumber\\
    &\times\frac{\mathbbm{1}_{\int_{T(q)\cap\mathcal{Y}_{k'}^{\circ}}g(x)\det(\mathbf{I}_d+\langle x-q, V(q)\rangle F(q))d\mathcal{H}^1(x)>0}}{\int_{T(q)\cap\mathcal{Y}_{k'}^{\circ}}g(x)\det(\mathbf{I}_d+\langle x-q, V(q)\rangle F(q))d\mathcal{H}^1(x)}.
\end{align}

Next we bound $\lambda(\mathfrak{S}_{i,j})$. By Lemma \ref{L3.14n} (note that $T^{-1}(\mathfrak{S}_{i,j})\subseteq \tilde{\mathcal{T}}_{1;k,k'}^*$), we have 
\begin{equation*}
    \lambda(\mathfrak{S}_{i,j})= \int_{\mathfrak{H}_{k,k'}\cap T^{-1}(\mathfrak{S}_{i,j})}|\langle V(q), \mathtt{a}_{k,k'}\rangle| d\mathcal{H}^{d-1}(q)\int_{T(q)}f(x)\det\big(\mathbf{I}_d+\langle x-q, V(q)\rangle F(q)\big)d\mathcal{H}^1(x).
\end{equation*}
Consequently, noting \eqref{eq.new.4.1n}, we get
\begin{equation}\label{lambdaSij}
    \lambda(\mathfrak{S}_{i,j})\leq C\mathcal{H}^{d-1}(\mathfrak{H}_{k,k'}\cap T^{-1}(\mathfrak{S}_{i,j})).
\end{equation}
Let $q_{i,j}\in H_{k,k'}$ be the intersection point of $[\mathtt{c}_i,\mathtt{c}_j']$ and $H_{k,k'}$. For any $q\in \mathfrak{H}_{k,k'}\cap T^{-1}(\mathfrak{S}_{i,j})$, there exist $x\in T(q)\cap L_i$ and $y\in T(q)\cap R_j$. Note that $\|x-\mathtt{c}_i\|\leq \delta_2\sqrt{d\eps}$ and $\|y-\mathtt{c}'_j\|\leq \delta_2\sqrt{d\eps}$. Moreover, by elementary computation, we have 
\begin{equation*}
    q_{i,j}=\frac{\langle \mathtt{c}_i,\mathtt{a}_{k,k'}\rangle\mathtt{c}'_j-\langle \mathtt{c}'_j,\mathtt{a}_{k,k'}\rangle\mathtt{c}_i+\mathtt{b}_{k,k'}(\mathtt{c}_i-\mathtt{c}_j')}{\langle \mathtt{c}_i-\mathtt{c}'_j, \mathtt{a}_{k,k'}\rangle}, \quad q=\frac{\langle x, \mathtt{a}_{k,k'}\rangle y-\langle y,\mathtt{a}_{k,k'}\rangle x+\mathtt{b}_{k,k'}(x-y)}{\langle x-y, \mathtt{a}_{k,k'}\rangle}.
\end{equation*}
Hence on noting $\langle \mathtt{c}_i-\mathtt{c}'_j, \mathtt{a}_{k,k'}\rangle,\langle x-y, \mathtt{a}_{k,k'}\rangle\geq c$ (by \eqref{bddakk}), we get $\|q-q_{i,j}\|\leq C\delta_2\sqrt{d\eps} \leq C\sqrt{\eps}$. Consequently, $\mathcal{H}^{d-1}(\mathfrak{H}_{k,k'}\cap T^{-1}(\mathfrak{S}_{i,j}))\leq C\eps^{(d-1)\slash 2}$, which combined with \eqref{lambdaSij} yields
\begin{equation}\label{eq.neq.4.8}
    \lambda(\mathfrak{S}_{i,j})\leq C\eps^{(d-1)\slash 2}. 
\end{equation}

For any $q\in\mathfrak{H}_{k,k'}\cap  T^{-1}(\tilde{\mathcal{S}}^*)$, we define
\begin{align}\label{Eqnn1}
    \EuScript{J}(q):=&\bigg(\int_{T(q)\cap\mathcal{X}_k^{\circ}}f(x)\det\big(\mathbf{I}_d+\langle x-q, V(q)\rangle F(q)\big)d\mathcal{H}^1(x)\bigg)\nonumber\\
    &\times \bigg(\int_{T(q)\cap\mathcal{Y}_{k'}^{\circ}}g(x)\det\big(\mathbf{I}_d+\langle x-q, V(q)\rangle F(q)\big)d\mathcal{H}^1(x)\bigg)\geq 0. 
\end{align}
By \eqref{Eq4.294}, \eqref{Eq4.302}, and \eqref{eq.neq.4.8}, noting that $\mathfrak{S}_{i,j}\subseteq\tilde{\mathcal{S}}^*$, we get
\begin{align}
    & \sum_{i=1}^I\sum_{j=1}^J
    \sum_{\ell=1}^{N_0}\sum_{\ell'=1}^{N_0}\mathcal{M}_{i,j;\ell,\ell'}\big(\log\big(\eps^{-(d+3)\slash 2}\mathcal{M}_{i,j;\ell,\ell'}\big)\big)_{+}\nonumber\\
    \leq& \sum_{i=1}^I\sum_{j=1}^J
    \sum_{\ell,\ell'\in[N_0]}\int_{\mathfrak{H}_{k,k'}\cap T^{-1}(\tilde{\mathcal{S}}^*)}\mathfrak{m}_{i,j;\ell,\ell'}(T(q))\bigg(\int_{T(q)}f(x)\det\big(\mathbf{I}_d+\langle x-q, V(q)\rangle F(q)\big)d\mathcal{H}^1(x)\bigg)\nonumber\\
    &\hspace{1.9in}\times|\langle V(q),\mathtt{a}_{k,k'}\rangle|\big((\log(C\delta_3^{-1}))_{+}+(\log\EuScript{J}(q))_{-}\big)d\mathcal{H}^{d-1}(q).
\end{align}
Note that by \eqref{def:mijl} and \eqref{subinclud1.1}, for any $q\in \mathfrak{H}_{k,k'}\cap T^{-1}(\tilde{\mathcal{S}}^*)$, we have
\begin{align*}
    \sum_{i=1}^I\sum_{j=1}^J
    \sum_{\ell,\ell'\in[N_0]}\mathfrak{m}_{i,j;\ell,\ell'}(T(q))=&\mathbbm{1}_{\int_{T(q)\cap \mathcal{X}_k^{\circ}}\bar{f}_2(z)\mathfrak{F}(z)d\mathcal{H}^1(z)>0}\cdot\int_{T(q)\cap\mathcal{Y}_{k'}^{\circ}}\bar{g}_2(z)\mathfrak{F}(z) d\mathcal{H}^1(z)\nonumber\\
    =&  \int_{T(q)\cap \mathcal{X}_k^{\circ}}\bar{f}_2(z)\mathfrak{F}(z)d\mathcal{H}^1(z),
\end{align*}
hence by Definition \ref{factors} and \eqref{subinclud1.1},
\begin{align}
 &\bigg(\sum_{i=1}^I\sum_{j=1}^J
    \sum_{\ell,\ell'\in[N_0]}\mathfrak{m}_{i,j;\ell,\ell'}(T(q))\bigg)\bigg(\int_{T(q)}f(x)\det\big(\mathbf{I}_d+\langle x-q, V(q)\rangle F(q)\big)d\mathcal{H}^1(x)\bigg)\nonumber\\
    =& \bigg(\int_{T(q)\cap \mathcal{X}_k^{\circ}}\bar{f}_2(x)\mathfrak{F}(x)d\mathcal{H}^1(x)\bigg)\bigg(\int_{T(q)}f(x)\det\big(\mathbf{I}_d+\langle x-q, V(q)\rangle F(q)\big)d\mathcal{H}^1(x)\bigg)\nonumber\\
    =& \int_{T(q)\cap\mathcal{X}_k^{\circ}}\bar{f}_2(x)\det\big(\mathbf{I}_d+\langle x-q, V(q)\rangle F(q)\big)d\mathcal{H}^1(x)\\
    =& \int_{T(q)\cap\mathcal{Y}_{k'}^{\circ}}\bar{g}_2(x)\det\big(\mathbf{I}_d+\langle x-q, V(q)\rangle F(q)\big)d\mathcal{H}^1(x).
\end{align}
By \eqref{fkkgkkb} and the fact that $\bar{f}_2(x)\leq f_{k,k'}(x)$ and $\bar{g}_2(x)\leq g_{k,k'}(x)$ for every $x\in\mathbb{R}^d$ (see \eqref{Eq.4.2new}),
\begin{align}\label{Eqnn2}
&\int_{T(q)\cap\mathcal{X}_k^{\circ}}\bar{f}_2(x)\det\big(\mathbf{I}_d+\langle x-q, V(q)\rangle F(q)\big)d\mathcal{H}^1(x)\leq \int_{T(q)\cap\mathcal{X}_k^{\circ}}f(x)\det\big(\mathbf{I}_d+\langle x-q, V(q)\rangle F(q)\big)d\mathcal{H}^1(x),\nonumber\\
& \int_{T(q)\cap\mathcal{Y}_{k'}^{\circ}}\bar{g}_2(x)\det\big(\mathbf{I}_d+\langle x-q, V(q)\rangle F(q)\big)d\mathcal{H}^1(x)\leq \int_{T(q)\cap\mathcal{Y}_{k'}^{\circ}}g(x)\det\big(\mathbf{I}_d+\langle x-q, V(q)\rangle F(q)\big)d\mathcal{H}^1(x).
\end{align}
Combining \eqref{Eqnn1}--\eqref{Eqnn2}, we get
\begin{align*}
    & \sum_{i=1}^I\sum_{j=1}^J
    \sum_{\ell=1}^{N_0}\sum_{\ell'=1}^{N_0}\mathcal{M}_{i,j;\ell,\ell'}\big(\log\big(\eps^{-(d+3)\slash 2}\mathcal{M}_{i,j;\ell,\ell'}\big)\big)_{+}\nonumber\\
    \leq& \int_{\mathfrak{H}_{k,k'}\cap T^{-1}(\tilde{\mathcal{S}}^*)}|\langle V(q),\mathtt{a}_{k,k'}\rangle|\min\bigg\{\EuScript{J}(q)^{1\slash 2},   \int_{T(q)\cap\mathcal{X}_k^{\circ}}\bar{f}_2(x)\det\big(\mathbf{I}_d+\langle x-q, V(q)\rangle F(q)\big)d\mathcal{H}^1(x)\bigg\}\nonumber\\
    &\hspace{0.96in}\times\big((\log(C\delta_3^{-1}))_{+}+(\log\EuScript{J}(q))_{-}\big)d\mathcal{H}^{d-1}(q).
\end{align*}
For any $\delta\in (0,e^{-1}]$ and $q\in\mathfrak{H}_{k,k'}\cap T^{-1}(\tilde{\mathcal{S}}^*)$ such that $\EuScript{J}(q)\leq \delta^2$, as $x\log x$ is monotonically decreasing on $[0,e^{-1}]$, we have $\EuScript{J}(q)^{1\slash 2}(\log \EuScript{J}(q))_{-}=-2\EuScript{J}(q)^{1\slash 2}\log(\EuScript{J}(q)^{1\slash 2})\leq 2\delta\log(\delta^{-1})$. Hence for any $\delta\in (0,e^{-1}]$ and $q\in\mathfrak{H}_{k,k'}\cap T^{-1}(\tilde{\mathcal{S}}^*)$,
\begin{align*}
   & \min\bigg\{\EuScript{J}(q)^{1\slash 2},   \int_{T(q)\cap\mathcal{X}_k^{\circ}}\bar{f}_2(x)\det\big(\mathbf{I}_d+\langle x-q, V(q)\rangle F(q)\big)d\mathcal{H}^1(x)\bigg\}\big((\log(C\delta_3^{-1}))_{+}+(\log \EuScript{J}(q))_{-}\big)\nonumber\\
   \leq& (\log(C\delta_3^{-1}))_{+}\int_{T(q)\cap\mathcal{X}_k^{\circ}}\bar{f}_2(x)\det\big(\mathbf{I}_d+\langle x-q, V(q)\rangle F(q)\big)d\mathcal{H}^1(x)+2\delta\log(\delta^{-1})\nonumber\\
   &+2\log(\delta^{-1})\int_{T(q)\cap\mathcal{X}_k^{\circ}}\bar{f}_2(x)\det\big(\mathbf{I}_d+\langle x-q, V(q)\rangle F(q)\big)d\mathcal{H}^1(x)\nonumber\\
   \leq& \big(C\log(\delta^{-1})+\log(\delta_3^{-1})\big)\int_{T(q)\cap\mathcal{X}_k^{\circ}}\bar{f}_2(x)\det\big(\mathbf{I}_d+\langle x-q, V(q)\rangle F(q)\big)d\mathcal{H}^1(x)+2\delta\log(\delta^{-1}).
\end{align*}
By the above two displays, we conclude that for any $\delta\in (0,e^{-1}]$,
\begin{align*}
    &\sum_{i=1}^I\sum_{j=1}^J
    \sum_{\ell=1}^{N_0}\sum_{\ell'=1}^{N_0}\mathcal{M}_{i,j;\ell,\ell'}\big(\log\big(\eps^{-(d+3)\slash 2}\mathcal{M}_{i,j;\ell,\ell'}\big)\big)_{+}\nonumber\\
    \leq& \big(C\log(\delta^{-1})+\log(\delta_3^{-1})\big)\nonumber\\
    &\times\int_{\mathfrak{H}_{k,k'}\cap T^{-1}(\tilde{\mathcal{S}}^*)}|\langle V(q),\mathtt{a}_{k,k'}\rangle|d\mathcal{H}^{d-1}(q)\int_{T(q)\cap\mathcal{X}_k^{\circ}}\bar{f}_2(x)\det\big(\mathbf{I}_d+\langle x-q, V(q)\rangle F(q)\big)d\mathcal{H}^1(x)\nonumber\\
    &+2\delta\log(\delta^{-1})\mathcal{H}^{d-1}\big(\mathfrak{H}_{k,k'}\cap T^{-1}(\tilde{\mathcal{S}}^*)\big)\nonumber\\
    \leq& \big(C\log(\delta^{-1})+\log(\delta_3^{-1})\big)\int_{\mathcal{S}}d\lambda(T)\int_{T\cap\mathcal{X}_k^{\circ}}\bar{f}_2(x)\mathfrak{F}(x)d\mathcal{H}^1(x)+2\delta\log(\delta^{-1})\mathcal{H}^{d-1}(\tilde{H}_{k,k'})\nonumber\\
    \leq& \big(C\log(\delta^{-1})+\log(\delta_3^{-1})\big)\int_{\mathbb{R}^d}\bar{f}_2(x)dx+C\delta\log(\delta^{-1})\nonumber\\
    =& \big(C\log(\delta^{-1})+\log(\delta_3^{-1})\big)\bar{\gamma}_{\eps,2}(\mathbb{R}^d\times\mathbb{R}^d)+C\delta\log(\delta^{-1}),
\end{align*}
where we use Lemma \ref{L3.14n} in the second inequality, Lemma \ref{Lem3.18nn} and \eqref{BddHtildekk} in the third inequality, and Lemma \ref{Lem4.9} in the last line. 
\end{proof}

\begin{lemma}\label{Lem4.27}
For any $\delta\in (0,e^{-1}]$, we have 
\begin{equation}\label{Lem4.30.e1}
    \sum_{i=1}^I \sum_{j=1}^J \int_{L_i} f_j(x)\big(\log\big(\eps^{-1\slash 2} f_j(x)\big)\big)_{+}dx\leq C\log(\delta^{-1})\bar{\gamma}_{\eps,2}(\mathbb{R}^d\times\mathbb{R}^d)+C\delta\log(\delta^{-1}),
\end{equation}
\begin{equation}\label{Lem4.30.e2}
    \sum_{i=1}^I\sum_{j=1}^J\int_{R_j}g_i(y)\big(\log\big(\eps^{-1\slash 2}g_i(y)\big)\big)_{+}dy\leq C\log(\delta^{-1})\bar{\gamma}_{\eps,2}(\mathbb{R}^d\times\mathbb{R}^d)+C\delta\log(\delta^{-1}).
\end{equation}
\end{lemma}
\begin{proof}

By \eqref{deffx}, recalling the definition of $\tilde{\Phi}(\cdot)$ from \eqref{Phi_def}, we have
\begin{align}\label{Eq4.312}
    & \sum_{i=1}^I \sum_{j=1}^J \int_{L_i} f_j(x)\big(\log\big(\eps^{-1\slash 2} f_j(x)\big)\big)_{+}dx \nonumber\\
    =& \sum_{i=1}^I\sum_{j=1}^J \int_{L_i\cap\mathfrak{T}^*}\mathbbm{1}_{\int_{\mathrm{int}(T(x))\cap \mathcal{Y}_{k'}^{\circ}}\bar{g}_2(z)\mathfrak{F}(z)d\mathcal{H}^1(z)>0} \bar{f}_2(x)\cdot\frac{\int_{\mathrm{int}(T(x))\cap R_j}\bar{g}_2(z)\mathfrak{F}(z)d\mathcal{H}^1(z)}{\int_{\mathrm{int}(T(x))\cap \mathcal{Y}_{k'}^{\circ}}\bar{g}_2(z)\mathfrak{F}(z)d\mathcal{H}^1(z)}\nonumber\\
    &\hspace{1in}\times\Bigg(\log\Bigg(\bar{f}_2(x)\cdot\frac{\int_{\mathrm{int}(T(x))\cap R_j}\bar{g}_2(z)\mathfrak{F}(z)d\mathcal{H}^1(z)}{\sqrt{\eps}\int_{\mathrm{int}(T(x))\cap \mathcal{Y}_{k'}^{\circ}}\bar{g}_2(z)\mathfrak{F}(z)d\mathcal{H}^1(z)}\Bigg)\Bigg)_{+}dx\nonumber\\
    \leq& \sum_{i=1}^I\sum_{j=1}^J \int_{L_i\cap\mathfrak{T}^*}\mathbbm{1}_{\int_{\mathrm{int}(T(x))\cap \mathcal{Y}_{k'}^{\circ}}\bar{g}_2(z)\mathfrak{F}(z)d\mathcal{H}^1(z)>0} \bar{f}_2(x)(\log\bar{f}_2(x))_{+}\cdot\frac{\int_{\mathrm{int}(T(x))\cap R_j}\bar{g}_2(z)\mathfrak{F}(z)d\mathcal{H}^1(z)}{\int_{\mathrm{int}(T(x))\cap \mathcal{Y}_{k'}^{\circ}}\bar{g}_2(z)\mathfrak{F}(z)d\mathcal{H}^1(z)}dx\nonumber\\
    &+\sum_{i=1}^I\sum_{j=1}^J \int_{L_i\cap\mathfrak{T}^*}\mathbbm{1}_{\int_{\mathrm{int}(T(x))\cap \mathcal{Y}_{k'}^{\circ}}\bar{g}_2(z)\mathfrak{F}(z)d\mathcal{H}^1(z)>0} \sqrt{\eps}\bar{f}_2(x)\cdot\tilde{\Phi}\Bigg(\frac{\int_{\mathrm{int}(T(x))\cap R_j}\bar{g}_2(z)\mathfrak{F}(z)d\mathcal{H}^1(z)}{\sqrt{\eps}\int_{\mathrm{int}(T(x))\cap \mathcal{Y}_{k'}^{\circ}}\bar{g}_2(z)\mathfrak{F}(z)d\mathcal{H}^1(z)}\Bigg)dx\nonumber\\
    =& \int_{\mathcal{X}_k^{\circ}\cap\mathfrak{T}^*}\mathbbm{1}_{\int_{\mathrm{int}(T(x))\cap \mathcal{Y}_{k'}^{\circ}}\bar{g}_2(z)\mathfrak{F}(z)d\mathcal{H}^1(z)>0} \bar{f}_2(x)(\log\bar{f}_2(x))_{+}dx\nonumber\\
    &+\sum_{j=1}^J\int_{\mathcal{X}_k^{\circ}\cap\mathfrak{T}^*}\mathbbm{1}_{\int_{\mathrm{int}(T(x))\cap \mathcal{Y}_{k'}^{\circ}}\bar{g}_2(z)\mathfrak{F}(z)d\mathcal{H}^1(z)>0} \sqrt{\eps}\bar{f}_2(x)\cdot\tilde{\Phi}\Bigg(\frac{\int_{\mathrm{int}(T(x))\cap R_j}\bar{g}_2(z)\mathfrak{F}(z)d\mathcal{H}^1(z)}{\sqrt{\eps}\int_{\mathrm{int}(T(x))\cap \mathcal{Y}_{k'}^{\circ}}\bar{g}_2(z)\mathfrak{F}(z)d\mathcal{H}^1(z)}\Bigg)dx.
\end{align}
By \eqref{Eq.4.2new}, Lemma \ref{Ln4.1}, and \eqref{bddLip}, $\bar{f}_2(x)\leq f_{k,k'}(x)\leq f(x)\leq M_0$ for $\mathcal{L}^d$-a.e.\ $x\in\mathbb{R}^d$. Hence by Lemma \ref{Lem4.9},
\begin{align}\label{Eq4.313}
&\int_{\mathcal{X}_k^{\circ}\cap\mathfrak{T}^*}\mathbbm{1}_{\int_{\mathrm{int}(T(x))\cap \mathcal{Y}_{k'}^{\circ}}\bar{g}_2(z)\mathfrak{F}(z)d\mathcal{H}^1(z)>0} \bar{f}_2(x)(\log\bar{f}_2(x))_{+}dx\nonumber\\
\leq& \int_{\mathfrak{T}^*}\bar{f}_2(x)(\log f(x))_{+}dx\leq \log(\max\{M_0,1\})\int_{\mathfrak{T}^*}\bar{f}_2(x)dx\leq C\bar{\gamma}_{\eps,2}(\mathbb{R}^d\times\mathbb{R}^d).
\end{align}
By Lemma \ref{Lem3.18nn} and \eqref{subinclud1.1}, noting that $\lambda(\mathcal{S}\backslash\tilde{\mathcal{S}}^*)=0$ (see Definition \ref{defstilde}), for any $j\in[J]$, we get 
\begin{align}\label{E4.314}
    & \int_{\mathcal{X}_k^{\circ}\cap\mathfrak{T}^*}\mathbbm{1}_{\int_{\mathrm{int}(T(x))\cap \mathcal{Y}_{k'}^{\circ}}\bar{g}_2(z)\mathfrak{F}(z)d\mathcal{H}^1(z)>0} \sqrt{\eps}\bar{f}_2(x)\cdot\tilde{\Phi}\Bigg(\frac{\int_{\mathrm{int}(T(x))\cap R_j}\bar{g}_2(z)\mathfrak{F}(z)d\mathcal{H}^1(z)}{\sqrt{\eps}\int_{\mathrm{int}(T(x))\cap \mathcal{Y}_{k'}^{\circ}}\bar{g}_2(z)\mathfrak{F}(z)d\mathcal{H}^1(z)}\Bigg)dx\nonumber\\
    =&\sqrt{\eps} \int_{\mathcal{S}}\mathbbm{1}_{\int_{\mathrm{int}(T)\cap \mathcal{Y}_{k'}^{\circ}}\bar{g}_2(z)\mathfrak{F}(z)d\mathcal{H}^1(z)>0}\bigg(\int_{\mathrm{int}(T)\cap\mathcal{X}_k^{\circ}}\bar{f}_2(x)\mathfrak{F}(x)d\mathcal{H}^1(x)\bigg)\nonumber\\
    &\hspace{0.32in}\times \tilde{\Phi}\Bigg(\frac{\int_{\mathrm{int}(T)\cap R_j}\bar{g}_2(z)\mathfrak{F}(z)d\mathcal{H}^1(z)}{\sqrt{\eps}\int_{\mathrm{int}(T)\cap \mathcal{Y}_{k'}^{\circ}}\bar{g}_2(z)\mathfrak{F}(z)d\mathcal{H}^1(z)}\Bigg) d\lambda(T)\nonumber\\
    =&\sqrt{\eps} \int_{\tilde{\mathcal{S}}^*}\mathbbm{1}_{\int_{\mathrm{int}(T)\cap \mathcal{Y}_{k'}^{\circ}}\bar{g}_2(z)\mathfrak{F}(z)d\mathcal{H}^1(z)>0}\bigg(\int_{\mathrm{int}(T)\cap\mathcal{Y}_{k'}^{\circ}}\bar{g}_2(z)\mathfrak{F}(z)d\mathcal{H}^1(z)\bigg)\nonumber\\
    &\hspace{0.32in}\times \tilde{\Phi}\Bigg(\frac{\int_{\mathrm{int}(T)\cap R_j}\bar{g}_2(z)\mathfrak{F}(z)d\mathcal{H}^1(z)}{\sqrt{\eps}\int_{\mathrm{int}(T)\cap \mathcal{Y}_{k'}^{\circ}}\bar{g}_2(z)\mathfrak{F}(z)d\mathcal{H}^1(z)}\Bigg) d\lambda(T)\nonumber\\
    =& \int_{\tilde{\mathcal{S}}^*} \mathbbm{1}_{\int_{\mathrm{int}(T)\cap \mathcal{Y}_{k'}^{\circ}}\bar{g}_2(z)\mathfrak{F}(z)d\mathcal{H}^1(z)>0}\bigg(\int_{\mathrm{int}(T)\cap\mathcal{Y}_{k'}^{\circ}\cap R_j}\bar{g}_2(z)\mathfrak{F}(z)d\mathcal{H}^1(z)\bigg)\nonumber\\
    &\hspace{0.25in}\times\Bigg(\log\Bigg(\frac{\int_{\mathrm{int}(T)\cap R_j}\bar{g}_2(z)\mathfrak{F}(z)d\mathcal{H}^1(z)}{\sqrt{\eps}\int_{\mathrm{int}(T)\cap \mathcal{Y}_{k'}^{\circ}}\bar{g}_2(z)\mathfrak{F}(z)d\mathcal{H}^1(z)}\Bigg)\Bigg)_{+}d\lambda(T).
\end{align}
We denote by $\mathfrak{S}$ the set of $T\in\tilde{\mathcal{S}}^*$ such that $T\cap\mathcal{X}_k^{\circ}\neq\emptyset,T\cap\mathcal{Y}_{k'}^{\circ}\neq\emptyset$. Note that $T^{-1}(\mathfrak{S})\subseteq\tilde{\mathcal{T}}_{1;k,k'}^*$, and for any $T\in\tilde{\mathcal{S}}^*\backslash\mathfrak{S}$, by \eqref{gkkdef},
\begin{align*}
    \int_{\mathrm{int}(T)\cap \mathcal{Y}_{k'}^{\circ}}\bar{g}_2(z)\mathfrak{F}(z)d\mathcal{H}^1(z)\leq& \int_{\mathrm{int}(T)\cap \mathcal{Y}_{k'}^{\circ}} g_{k,k'}(z)\mathfrak{F}(z)d\mathcal{H}^1(z)\nonumber\\
    \leq&\int_{(\mathrm{int}(T)\cap \mathcal{X}_{k}^{\circ})\times (\mathrm{int}(T)\cap\mathcal{Y}_{k'}^{\circ})}h(x,y)d\mathcal{H}^1(x)d\mathcal{H}^1(y)=0.
\end{align*}
Hence by \eqref{E4.314} and Lemma \ref{L3.14n}, and noting Definition \ref{factors}, we get
\begin{align}\label{E4.315}
    & \int_{\mathcal{X}_k^{\circ}\cap\mathfrak{T}^*}\mathbbm{1}_{\int_{\mathrm{int}(T(x))\cap \mathcal{Y}_{k'}^{\circ}}\bar{g}_2(z)\mathfrak{F}(z)d\mathcal{H}^1(z)>0} \sqrt{\eps}\bar{f}_2(x)\cdot\tilde{\Phi}\Bigg(\frac{\int_{\mathrm{int}(T(x))\cap R_j}\bar{g}_2(z)\mathfrak{F}(z)d\mathcal{H}^1(z)}{\sqrt{\eps}\int_{\mathrm{int}(T(x))\cap \mathcal{Y}_{k'}^{\circ}}\bar{g}_2(z)\mathfrak{F}(z)d\mathcal{H}^1(z)}\Bigg)dx\nonumber\\
    =& \int_{\mathfrak{S}} \mathbbm{1}_{\int_{\mathrm{int}(T)\cap \mathcal{Y}_{k'}^{\circ}}\bar{g}_2(z)\mathfrak{F}(z)d\mathcal{H}^1(z)>0}\bigg(\int_{\mathrm{int}(T)\cap\mathcal{Y}_{k'}^{\circ}\cap R_j}\bar{g}_2(z)\mathfrak{F}(z)d\mathcal{H}^1(z)\bigg)\nonumber\\
    &\hspace{0.25in}\times\Bigg(\log\Bigg(\frac{\int_{\mathrm{int}(T)\cap R_j}\bar{g}_2(z)\mathfrak{F}(z)d\mathcal{H}^1(z)}{\sqrt{\eps}\int_{\mathrm{int}(T)\cap \mathcal{Y}_{k'}^{\circ}}\bar{g}_2(z)\mathfrak{F}(z)d\mathcal{H}^1(z)}\Bigg)\Bigg)_{+}d\lambda(T)\nonumber\\
    =& \int_{\mathfrak{H}_{k,k'}\cap T^{-1}(\mathfrak{S})}|\langle V(q), \mathtt{a}_{k,k'}\rangle|\mathbbm{1}_{\int_{T(q)\cap \mathcal{Y}_{k'}^{\circ}}\bar{g}_2(x)\det(\mathbf{I}_d+\langle x-q,V(q)\rangle F(q))d\mathcal{H}^1(x)>0}\nonumber\\
    &\hspace{1in}\times
    \bigg(\int_{T(q)\cap R_j}\bar{g}_2(x)\det\big(\mathbf{I}_d+\langle x-q,V(q)\rangle F(q)\big)d\mathcal{H}^1(x)\bigg)\nonumber\\
    &\hspace{1in}\times \Bigg(\log\Bigg(\frac{\int_{T(q)\cap R_j}\bar{g}_2(x)\det(\mathbf{I}_d+\langle x-q,V(q)\rangle F(q))d\mathcal{H}^1(x)}{\sqrt{\eps}\int_{T(q)\cap \mathcal{Y}_{k'}^{\circ}}\bar{g}_2(x)\det(\mathbf{I}_d+\langle x-q,V(q)\rangle F(q))d\mathcal{H}^1(x)}\Bigg)\Bigg)_{+}d\mathcal{H}^{d-1}(q).
\end{align}
For any $q\in\mathfrak{H}_{k,k'}\cap T^{-1}(\mathfrak{S})$ and $x\in T(q)$, as $\min\{\alpha(q),\beta(q)\}\geq d_0$ (recall that $\mathfrak{H}_{k,k'}\subseteq\tilde{H}_{k,k'}$ and \eqref{deftildehat}), by Hadamard's inequality (see, e.g., \cite[Lemma 2.5]{ipsen2008perturbation}), \eqref{bddsalphabeta}, and \Cref{L2.4} (see \eqref{L2.4.enew3}), we have
\begin{equation*}
    \big|\det\big(\mathbf{I}_d+\langle x-q, V(q)\rangle F(q)\big)\big|\leq \big\|\mathbf{I}_d+\langle x-q, V(q)\rangle F(q)\big\|_2^d\leq (1+2D\|F(q)\|_2)^d\leq C.
\end{equation*}
Hence by \eqref{Eq.4.2new} and \eqref{fkkgkkb}, for any $q\in\mathfrak{H}_{k,k'}\cap T^{-1}(\mathfrak{S})$,
\begin{equation}\label{E4.316}
\int_{T(q)\cap R_j}\bar{g}_2(x)\det\big(\mathbf{I}_d+\langle x-q,V(q)\rangle F(q)\big)d\mathcal{H}^1(x)\leq C\mathcal{H}^1(T(q)\cap R_j)\leq C\delta_2\sqrt{\eps}\leq C\sqrt{\eps}.
\end{equation}
For any $q\in\mathfrak{H}_{k,k'}\cap T^{-1}(\mathfrak{S})$, define $\EuScript{C}(q):=\int_{T(q)\cap \mathcal{Y}_{k'}^{\circ}}\bar{g}_2(x)\det\big(\mathbf{I}_d+\langle x-q,V(q)\rangle F(q)\big)d\mathcal{H}^1(x)$. By \eqref{E4.315} and \eqref{E4.316}, we have
\begin{align*}
    & \sum_{j=1}^J \int_{\mathcal{X}_k^{\circ}\cap\mathfrak{T}^*}\mathbbm{1}_{\int_{\mathrm{int}(T(x))\cap \mathcal{Y}_{k'}^{\circ}}\bar{g}_2(z)\mathfrak{F}(z)d\mathcal{H}^1(z)>0} \sqrt{\eps}\bar{f}_2(x)\cdot\tilde{\Phi}\Bigg(\frac{\int_{\mathrm{int}(T(x))\cap R_j}\bar{g}_2(z)\mathfrak{F}(z)d\mathcal{H}^1(z)}{\sqrt{\eps}\int_{\mathrm{int}(T(x))\cap \mathcal{Y}_{k'}^{\circ}}\bar{g}_2(z)\mathfrak{F}(z)d\mathcal{H}^1(z)}\Bigg)dx\nonumber\\
    \leq& \int_{\mathfrak{H}_{k,k'}\cap T^{-1}(\mathfrak{S})}|\langle V(q), \mathtt{a}_{k,k'}\rangle|\bigg(\int_{T(q)\cap\mathcal{Y}_{k'}^{\circ}}\bar{g}_2(x)\det\big(\mathbf{I}_d+\langle x-q,V(q)\rangle F(q)\big)d\mathcal{H}^1(x)\bigg)\nonumber\\
    &\hspace{0.9in}\times \Bigg(\log\Bigg(\frac{C}{\int_{T(q)\cap \mathcal{Y}_{k'}^{\circ}}\bar{g}_2(x)\det(\mathbf{I}_d+\langle x-q,V(q)\rangle F(q))d\mathcal{H}^1(x)}\Bigg)\Bigg)_{+}d\mathcal{H}^{d-1}(q)\nonumber\\
    \leq & \int_{\mathfrak{H}_{k,k'}\cap T^{-1}(\mathfrak{S})}|\langle V(q), \mathtt{a}_{k,k'}\rangle|\EuScript{C}(q)\big(C+(\log{\EuScript{C}(q)})_{-}\big)d\mathcal{H}^{d-1}(q).
\end{align*}
For any $\delta\in (0,e^{-1}]$ and $q\in\mathfrak{H}_{k,k'}\cap T^{-1}(\mathfrak{S})$ such that $\EuScript{C}(q)\leq \delta$, as $x\log x$ is monotonically decreasing on $[0,e^{-1}]$, we have $\EuScript{C}(q)(\log{\EuScript{C}(q)})_{-}=-\EuScript{C}(q)\log(\EuScript{C}(q))\leq \delta\log(\delta^{-1})$. Hence for any $\delta\in (0,e^{-1}]$ and $q\in\mathfrak{H}_{k,k'}\cap T^{-1}(\mathfrak{S})$,  
\begin{equation*}
    \EuScript{C}(q)\big(C+(\log{\EuScript{C}(q)})_{-}\big)\leq (C+\log(\delta^{-1}))\EuScript{C}(q)+\delta\log(\delta^{-1})\leq C\log(\delta^{-1})\EuScript{C}(q) + \delta\log(\delta^{-1}). 
\end{equation*}
By the above two displays, and noting that $\mathfrak{H}_{k,k'}\subseteq\tilde{H}_{k,k'}$, we obtain that for any $\delta\in(0,e^{-1}]$,
\begin{align}\label{Eq4.317}
    &\sum_{j=1}^J \int_{\mathcal{X}_k^{\circ}\cap\mathfrak{T}^*}\mathbbm{1}_{\int_{\mathrm{int}(T(x))\cap \mathcal{Y}_{k'}^{\circ}}\bar{g}_2(z)\mathfrak{F}(z)d\mathcal{H}^1(z)>0} \sqrt{\eps}\bar{f}_2(x)\cdot\tilde{\Phi}\Bigg(\frac{\int_{\mathrm{int}(T(x))\cap R_j}\bar{g}_2(z)\mathfrak{F}(z)d\mathcal{H}^1(z)}{\sqrt{\eps}\int_{\mathrm{int}(T(x))\cap \mathcal{Y}_{k'}^{\circ}}\bar{g}_2(z)\mathfrak{F}(z)d\mathcal{H}^1(z)}\Bigg)dx\nonumber\\
    \leq& C\log(\delta^{-1})\int_{\mathfrak{H}_{k,k'}\cap T^{-1}(\mathfrak{S})}|\langle V(q), \mathtt{a}_{k,k'}\rangle|\EuScript{C}(q)d\mathcal{H}^{d-1}(q)+\delta\log(\delta^{-1})\mathcal{H}^{d-1}(\tilde{H}_{k,k'})\nonumber\\
    \leq& C\log(\delta^{-1})\int_{\mathcal{S}}d\lambda(T)\int_T \bar{g}_2(x)\mathfrak{F}(x)d\mathcal{H}^1(x)+C\delta\log(\delta^{-1})\nonumber\\
    \leq& C\log(\delta^{-1})\int_{\mathbb{R}^d}\bar{g}_2(x)dx+C\delta\log(\delta^{-1})=  C\log(\delta^{-1})\bar{\gamma}_{\eps,2}(\mathbb{R}^d\times\mathbb{R}^d)+C\delta\log(\delta^{-1}),
\end{align}
where we use \Cref{L3.14n} and \eqref{BddHtildekk} in the second inequality, and \Cref{Lem3.18nn,Lem4.9} in the last line. Combining \eqref{Eq4.312}, \eqref{Eq4.313}, and \eqref{Eq4.317}, we obtain \eqref{Lem4.30.e1}. The bound \eqref{Lem4.30.e2} can be similarly deduced.
\end{proof}

We now complete the proof of Proposition \ref{P4.4}.

\begin{proof}[Proof of Proposition~\ref{P4.4}]
We fix any $K_0\in\mathbb{N}^{*}$ such that $2^{K_0}\geq 4d_0^{-1}$ and any $\delta\in (0,e^{-1}]$. By \eqref{Pro4.4.Eq1} and Lemmas \ref{Lem4.24}--\ref{Lem4.27}, we have 
\begin{eqnarray}\label{Pro4.4.Eq7}
    &&\int_{\mathbb{R}^d\times\mathbb{R}^d}r(x,y)\max\big\{\log\big(\eps^{(d-1)\slash 2}r(x,y)\big),0\big\}  dxdy\nonumber\\
    &\leq& C\big(\log(\delta^{-1})+\log(\delta_3^{-1})+\log(\delta_2^{-1})+K_0\big)\bar{\gamma}_{\eps,2}(\mathbb{R}^d\times\mathbb{R}^d)+C\delta\log(\delta^{-1})\nonumber\\
    &&-2\bar{\gamma}_{\eps,2}(\mathbb{R}^d\times\mathbb{R}^d)\log\big(\bar{\gamma}_{\eps,2}(\mathbb{R}^d\times\mathbb{R}^d)\big)+CK_02^{-K_0}.
\end{eqnarray}
Note that by Lemma \ref{Lem4.9} and \eqref{Eq.4.1new}--\eqref{Eq.4.2new},
\begin{equation*}
    \bar{\gamma}_{\eps,2}(\mathbb{R}^d\times\mathbb{R}^d)=\int_{\mathbb{R}^d}\bar{f}_2(x)dx=\int_{\mathbb{R}^d}f_{k,k'}(x)dx-\int_{\mathbb{R}^d}\bar{f}_1(x)dx=\mu_{k,k'}(\mathbb{R}^d)-(1-\delta_3)\int_{\mathbb{R}^d\times\mathbb{R}^d}p(x,y)dxdy.
\end{equation*}
Hence by Proposition \ref{P4.3},
\begin{equation}\label{Pro4.4.Eq6}
    \limsup_{\delta_1\rightarrow 0^+}\limsup_{M_1\rightarrow\infty}\limsup_{\eps\rightarrow 0^+}\bar{\gamma}_{\eps,2}(\mathbb{R}^d\times\mathbb{R}^d)=\delta_3 \mu_{k,k'}(\mathbb{R}^d).
\end{equation}
Sending $\eps\rightarrow 0^+$, $M_1\rightarrow\infty$, and $\delta_1\rightarrow 0^+$ sequentially in \eqref{Pro4.4.Eq7}, and using \eqref{Pro4.4.Eq6}, we obtain that 
\begin{align*}
&\limsup_{\delta_1\rightarrow 0^+}\limsup_{M_1\rightarrow\infty}\limsup_{\eps\rightarrow 0^+}\int_{\mathbb{R}^d\times\mathbb{R}^d}r(x,y)\max\big\{\log\big(\eps^{(d-1)\slash 2}r(x,y)\big),0\big\} dxdy\nonumber\\
\leq \,& C\big(\log(\delta^{-1})+\log(\delta_3^{-1})+\log(\delta_2^{-1})+K_0\big)\delta_3\mu_{k,k'}(\mathbb{R}^d)
+C\delta\log(\delta^{-1})\nonumber\\
&-2\delta_3\mu_{k,k'}(\mathbb{R}^d)\log\big(\delta_3\mu_{k,k'}(\mathbb{R}^d)\big)+CK_02^{-K_0}. 
\end{align*}
Sending $\delta_3\rightarrow 0^+$, $\delta\rightarrow 0^+$, and $K_0\rightarrow\infty$ sequentially in the above display, we conclude that 
\begin{equation*}
   \limsup_{\delta_3\rightarrow 0^{+}} \limsup_{\delta_1\rightarrow 0^{+}}\limsup_{M_1\rightarrow\infty}\limsup_{\eps\rightarrow 0^{+}}\int_{\mathbb{R}^d\times\mathbb{R}^d}r(x,y)\max\big\{\log\big(\eps^{(d-1)\slash 2}r(x,y)\big),0\big\}dxdy \leq 0. \qedhere
\end{equation*}
\end{proof}

\section{Counterexamples}\label{sec:counterex}

\subsection{Sharpness of the error tolerance in \cref{Thm1.1}}\label{sec:error_tolerance}

In this section, we demonstrate the sharpness of the error tolerance $o(\eps)$ in \cref{Thm1.1}. Assuming that the marginal supports $\mathcal{X}$ and $\mathcal{Y}$ are strictly separated by a hyperplane, we show that there exist $O(\eps)$-optimal couplings $\gamma_{\eps}'\in \Pi(\mu,\nu)$ which converge to an optimal transport plan $\gamma_{0}'$ different from the entropic Monge plan~$\gamma_0$.

We say that $\mathcal{X}$ and $\mathcal{Y}$ are strictly separated by a hyperplane if there exist $\mathtt{a}_0\in\mathbb{R}^d$ and $\mathtt{b}_0\in \mathbb{R}$ with $\|\mathtt{a}_0\|=1$ such that $\langle x, \mathtt{a}_0\rangle > \mathtt{b}_0$ for all $x\in\mathcal{X}$ and $\langle y, \mathtt{a}_0\rangle < \mathtt{b}_0$ for all $y\in\mathcal{Y}$. 

\begin{theorem}\label{thmcounter}
Suppose that $\mathcal{X}$ and $\mathcal{Y}$ are strictly separated by a hyperplane. Then there exist $\gamma_{\eps}'\in \Pi(\mu,\nu)$ such that as $\eps\rightarrow 0^+$, we have $\mathcal{C}_{\eps}(\gamma_{\eps}')\leq \EOT_\eps(\mu,\nu)+O(\eps)$ and $\gamma_{\eps}'$ converges weakly to $\gamma_0':=\lambda\otimes \kappa_T'$, where $\kappa_T':=\tilde{\mu}_T\otimes \tilde{\nu}_T$ for each transport ray $T$. In particular, $\gamma_{\eps}'$ does not converge to the entropic Monge plan~$\gamma_0$.
\end{theorem}

\begin{proof} 
Since $\mathcal{X}$ and $\mathcal{Y}$ are strictly separated by a hyperplane, we can take $K=K'=1$ in the domain decomposition in \Cref{sectdom}. Let $\hat{\gamma}_{\eps}$ be the Borel measure on $\mathbb{R}^d\times\mathbb{R}^d$ with density
\begin{align}\label{checks}
   \hat{p}_1 (x,y)   :=\,& \mathbbm{1}_{x\in\mathcal{X}_k^{\circ}\cap \mathfrak{T}^*
\cap\tilde{\mathcal{T}}_{1;k,k'}^{*}}\mathbbm{1}_{\mathtt{z}(x,y)\in\mathrm{int}(T(x))\cap\mathcal{Y}_{k'}^{\circ}}\mathbbm{1}_{\langle x-\mathtt{z}(x,y),V(x)\rangle>0}\mathfrak{F}(x)^{-1}\tilde{f}(x)\tilde{g}(\mathtt{z}(x,y))\nonumber\\
   & \times (2\pi\|x-\mathtt{z}(x,y)\|\eps)^{-(d-1)\slash 2} \sqrt{\det\big(\mathbf{I}_d+\|x-\mathtt{z}(x,y)\|F(\mathtt{z}(x,y))\big)}\nonumber\\
   & \times e^{-\frac{1}{2\eps}\mathtt{w}(x,y)^{\top}\big(\frac{1}{\|x-\mathtt{z}(x,y)\|}\mathbf{I}_d+F(\mathtt{z}(x,y))\big)\mathtt{w}(x,y)}
\end{align}
with respect to $\mathcal{L}^d\otimes\mathcal{L}^d$, where $k=k'=1$. Note that $\hat{p}_1 (x,y)$ is obtained from $\tilde{p}_1(x,y)$ in \eqref{p1til} by replacing $h(x,\mathtt{z}(x,y))$ with $\tilde{f}(x)\tilde{g}(\mathtt{z}(x,y))$. Thus, replacing $h_T(x,y)$ by $\tilde{f}(x)\tilde{g}(y)$ for $T\in\tilde{\mathcal{S}}^*$ and $(x,y)\in T\times T$ (cf.\ \cref{defstilde,Lem2}) in the arguments of \cref{Sec4} (cf.\ \Cref{th:upperBound}, \Cref{P4.3}, and \Cref{Lemma4.15}), we obtain couplings $\gamma_{\eps}'\in \Pi(\mu,\nu)$ for every $\eps>0$ such that
\begin{equation}\label{ubbdd}
     \limsup_{\eps\rightarrow 0^+}\bigg\{\frac{\mathcal{C}_{\eps}(\gamma_{\eps}')-\OT(\mu,\nu)}{\eps}+\frac{d-1}{2}\log{\eps}\bigg\}\leq \mathtt{S}(\mu,\nu;\gamma_0'),
\end{equation}
\begin{equation}\label{tvconver}
    d_{\mathrm{TV}}(\gamma_{\eps}',\hat{\gamma}_{\eps})\to 0 \quad\text{ as }\eps\to 0^{+},
\end{equation}
where $d_{\mathrm{TV}}(\cdot,\cdot)$ is the total variation distance. Using \cref{de:Sfunctional}, the assumption of the theorem, and \cref{L3.18}, we obtain that $\mathtt{S}(\mu,\nu;\gamma_0')<\infty$. Hence by \eqref{ubbdd} and the fact that $\mathcal{C}_{\eps}(\gamma_{\eps}')\geq \int_{\mathbb{R}^d\times\mathbb{R}^d}\|x-y\|d\gamma_{\eps}'(x,y)\geq \OT(\mu,\nu)$, we have
\begin{equation}\label{optimal _costs}
    \lim_{\eps\to 0^{+}}\mathcal{C}_{\eps}(\gamma_{\eps}')=\OT(\mu,\nu).
\end{equation}

Consider any sequence $\eps_j\to 0^+$. As $\Pi(\mu,\nu)$ is weakly compact, a subsequence of $(\gamma'_{\eps_j})$ (still denoted $(\gamma'_{\eps_j})$) converges weakly to some limit $\gamma_0''\in\Pi(\mu,\nu)$. By \eqref{optimal _costs}, $\gamma_0''$ is an optimal transport plan. Below we show that $\gamma_0''=\gamma_0'$. For any $(x,y)\in\mathbb{R}^d\times\mathbb{R}^d$, we define $\mathfrak{r}(x,y):=(x,\langle y-x, V(x)\rangle)$. Let $\mathfrak{P}',\mathfrak{P}'',\hat{\mathfrak{P}}_{\eps}$ be the pushforward measures of $\gamma_0', \gamma_0'',\hat{\gamma}_{\eps}$ (respectively) by $\mathfrak{r}$. By \eqref{tvconver}, $\hat{\gamma}_{\eps_j}$ converges weakly to $\gamma_0''$. For any continuous bounded function $\varphi$ on $\mathbb{R}^{d+1}$ and $\eps>0$, we have
\begin{align*}
    \int_{\mathbb{R}^{d+1}}\varphi d\hat{\mathfrak{B}}_{\eps}=\int_{\mathbb{R}^d\times\mathbb{R}^d}\varphi(x,\langle y-x,V(x)\rangle)d\hat{\gamma}_{\eps}(x,y).
\end{align*}
Define $\Psi(x,y):=\varphi(x,\langle y-x,V(x)\rangle)$ for any $x,y\in\mathbb{R}^d$. Note that $\Psi(\cdot,\cdot)$ is bounded and Borel measurable, and for each $x\in\mathbb{R}^d$, $\Psi(x,\cdot)$ is continuous. Hence by the above display and \Cref{L3.24}, and noting that the first marginals of $\hat{\gamma}_{\eps}$ and $\gamma_0''$ are all given by $\mu$, we get
\begin{equation*}
    \lim_{j\rightarrow\infty}\int_{\mathbb{R}^{d+1}}\varphi d\hat{\mathfrak{B}}_{\eps_j}=\lim_{j\rightarrow\infty}\int_{\mathbb{R}^d\times\mathbb{R}^d}\Psi(x,y)d\hat{\gamma}_{\eps_j}(x,y)=\int_{\mathbb{R}^d\times\mathbb{R}^d}\Psi(x,y)d\gamma_0''(x,y)=\int_{\mathbb{R}^{d+1}}\varphi d\mathfrak{B}''.
\end{equation*}
Consequently, $\mathfrak{P}_{\eps_j}$ converges weakly to $\mathfrak{P}''$. From \eqref{checks},  $\mathfrak{P}_{\eps_j}=\mathfrak{P}'$ for all $j$, and so $\mathfrak{P}''=\mathfrak{P}'$. This combined with the fact that $\gamma_0',\gamma_0''$ are both optimal transport plans yields $\gamma_0''=\gamma_0'$.   

Therefore, $\gamma_{\eps}'$ converges weakly to $\gamma_0'$. By \eqref{ubbdd} and \cref{Thm1.2}, $\mathcal{C}_{\eps}(\gamma_{\eps}')\leq \EOT_\eps(\mu,\nu)+O(\eps)$. From \cref{ruleskappa} below, we have $\gamma_0'\neq \gamma_0$. 
\end{proof}

\begin{lemma}\label{ruleskappa}
For $\lambda$-a.e.\ $T\in\mathcal{S}$, the solution $\kappa_T$ of the constrained EOT problem \eqref{opt.solver} cannot be $\tilde{\mu}_T\otimes \tilde{\nu}_T$.
\end{lemma}
\begin{proof}
Consider any $T\in\tilde{\mathcal{S}}^*$ (note from \Cref{defstilde} that $\lambda(\mathcal{S}\backslash\tilde{\mathcal{S}}^*)=0$). If $\tilde{\mu}_T\otimes \tilde{\nu}_T$ is not monotone on $T$, then it does not satisfy the constraint in \eqref{opt.solver} and hence cannot be a solution.

Below we assume that $\tilde{\mu}_T\otimes \tilde{\nu}_T$ is monotone on $T$. Note that this, together with \cref{assump1}, implies the existence of closed intervals $A,B\subseteq T$ such that $A\cap B=\emptyset$, $\tilde{\mu}_T$ is concentrated on $A$, $\tilde{\nu}_T$ is concentrated on $B$, and $u(x)>u(y)$ for all $x\in A, y\in B$. In particular, any $\chi_T\in \Pi(\tilde{\mu}_T,\tilde{\nu}_T)$ is monotone on $T$. Hence \eqref{opt.solver} is equivalent to
\begin{equation*}
\inf_{\chi_T\in\Pi(\tilde{\mu}_T,\tilde{\nu}_T)}\bigg\{-\frac{d-1}{2} \int_{T\times T} \log(\|x-y\|)d\chi_T(x,y) + H(\chi_T|\tilde{\mu}_T\otimes \tilde{\nu}_T)\bigg\}.
\end{equation*}
It follows that $\kappa_T$ satisfies (see, e.g., \cite[Theorem 4.2]{Nutz.20}) 
\begin{equation*}
    \frac{d\kappa_T(x,y)}{d(\tilde{\mu}_T\otimes\tilde{\nu}_T)}=e^{\varphi(x)+\psi(y)+\frac{d-1}{2}\log(\|x-y\|)}\quad\text{for } \tilde{\mu}_T\otimes \tilde{\nu}_T\text{-a.e.\ } (x,y)\in T\times T,
\end{equation*}
where $\varphi,\psi:T\rightarrow \mathbb{R}$ are Borel measurable. 

If $\kappa_T=\tilde{\mu}_T\otimes\tilde{\nu}_T$, then 
\begin{equation*}
    \varphi(x)+\psi(y)+\frac{d-1}{2}\log(\|x-y\|)=0\quad\text{for } \tilde{\mu}_T\otimes \tilde{\nu}_T\text{-a.e.\ } (x,y)\in T\times T.
\end{equation*}
Recall from \Cref{Sect.1.1} that $\tilde{\mu}_T$ and $\tilde{\nu}_T$ are absolutely continuous with respect to $\mathcal{H}^1|_T$. By Fubini's theorem, there exists a Borel set $A_0\subseteq A$ with $\tilde{\mu}_T(A_0)=1$ such that for every $x\in A_0$, the identity  $\varphi(x)+\psi(y)+\frac{d-1}{2}\log(\|x-y\|)=0$ holds for $\tilde{\nu}_T$-a.e.\ $y\in B$. Since $\tilde{\mu}_T$ is absolutely continuous with respect to $\mathcal{H}^1|_T$, $A_0$ must contain at least two distinct points; we pick $x_1\neq x_2$ in $A_0$. For each $i\in\{1,2\}$, let
\begin{equation*}
    B_i:=\bigg\{y\in B: \varphi(x_i)+\psi(y)+\frac{d-1}{2}\log(\|x_i-y\|)=0\bigg\}.
\end{equation*}
Then $\tilde{\nu}_T(B_i)=1$ for each $i\in\{1,2\}$, hence $\tilde{\nu}_T(B_1\cap B_2)=1$. Since $\tilde{\nu}_T$ is absolutely continuous with respect to $\mathcal{H}^1|_T$, it follows that $B_1\cap B_2$ contains at least two distinct points. We choose $y_1\neq y_2$ in $B_1\cap B_2$. Note that for all $i,j\in\{1,2\}$, we have
\begin{equation*}
    \varphi(x_i)+\psi(y_j)+\frac{d-1}{2}\log(\|x_i-y_j\|)=0.
\end{equation*}
Consequently, we have
\begin{align}\label{xyss}
  &  \log(\|x_1-y_1\|)+\log(\|x_2-y_2\|)-\log(\|x_1-y_2\|)-\log(\|x_2-y_1\|)=0, \nonumber\\
   & \Rightarrow \quad \|x_1-y_1\|\|x_2-y_2\|=\|x_1-y_2\|\|x_2-y_1\|. 
\end{align}
Suppose that $x_i=b_T+s_i V_T$ and $y_i=b_T+t_i V_T$ for every $i\in\{1,2\}$, 
where $b_T$ and $V_T$ are as in \Cref{TransportRays}. As $u(x)>u(y)$ for all $x\in A,y\in B$ and $x_1\neq x_2,y_1\neq y_2$, we have $s_i>t_j$ for all $i,j\in\{1,2\}$ and $s_1\neq s_2,t_1\neq  t_2$. Hence by \eqref{xyss}, we have $(s_1-t_1)(s_2-t_2)=(s_1-t_2)(s_2-t_1)$, and therefore $(s_1-s_2)(t_1-t_2)=0$. This contradicts the fact that  $s_1\neq s_2,t_1\neq  t_2$. Therefore, we conclude that $\kappa_T\neq \tilde{\mu}_T\otimes\tilde{\nu}_T$.   
\end{proof}

\subsection{Necessity of the directional constraint to define $\gamma_0$}\label{sec:counterex_directional_constraint}

We construct a counterexample to show that the directional constraint in~\eqref{opt.solver} cannot be omitted. That is, if a coupling $\gamma_0'$ is defined by the unconstrained analogue of~\eqref{opt.solver}, then $\gamma_0'$ is not the limit of the EOT optimizers $\gamma_\eps$. In fact, it is not even an optimal transport plan.

Let $d\geq 2$ and
\begin{equation*}
    A_1:=[0,2]\times [0,1]^{d-1},\quad A_2:=[5,6]\times [0,1]^{d-1},\quad B_1:=[3,4]\times [0,1]^{d-1},\quad B_2:=[7,9]\times [0,1]^{d-1}.
\end{equation*}
For every $x\in\mathbb{R}^d$, let
\begin{equation*}
  f(x):=\frac{\mathbbm{1}_{A_1}(x)}{3}+\frac{\mathbbm{1}_{A_2}(x)}{3},\quad g(x)=\frac{\mathbbm{1}_{B_1}(x)}{3}+\frac{\mathbbm{1}_{B_2}(x)}{3}.
\end{equation*}
Let $\mu$ and $\nu$ be Borel probability measures on $\mathbb{R}^d$ with densities $f$ and $g$, respectively, with respect to $\mathcal{L}^d$.

\begin{theorem}
Consider any Kantorovich potential $u\in\mathrm{Lip}_1(\mathbb{R}^d)$ for $(\mu,\nu)$ and its associated transport rays. For each transport ray $T$, define $\kappa_T'$ as the solution of the following problem:
\begin{equation}\label{optunconstrained}
\inf_{\chi_T\in\Pi(\tilde{\mu}_T,\tilde{\nu}_T)}\bigg\{-\frac{d-1}{2} \int_{T\times T} \log(\|x-y\|)d\chi_T(x,y) + H(\chi_T|\tilde{\mu}_T\otimes \tilde{\nu}_T)\bigg\},
\end{equation}
that is, the problem obtained from \eqref{opt.solver} by removing the directional constraint. Set $\gamma_0':=\lambda\otimes \kappa_T'$. Then $\gamma_0'$ is not an optimal transport plan. 
\end{theorem}
\begin{proof}
Note that $\mathcal{X}=\mathrm{supp}(\mu)=([0,2]\cup [5,6])\times [0,1]^{d-1}$ and $\mathcal{Y}=\mathrm{supp}(\nu)=([3,4]\cup[7,9])\times [0,1]^{d-1}$. We first show that for any Kantorovich potential $u$ and $(x_2,\cdots,x_d)\in [0,1]^{d-1}$, we have 
\begin{equation}\label{ray_identity}
    u(0,x_2,\cdots,x_d)-u(9,x_2,\cdots,x_d)=9.
\end{equation}
For any $x=(x_1,\cdots,x_d)\in\mathbb{R}^d$, define $u_0(x):=-x_1$ and  
\begin{equation*}
    s(x):=\begin{cases}
        (x_1+3,x_2,\cdots,x_d) & \text{ if }x_1\in [0,1]\text{ and }(x_2,\cdots,x_d)\in [0,1]^{d-1}\\
        (x_1+6,x_2,\cdots,x_d) & \text{ if }x_1\in (1,2]\text{ and }(x_2,\cdots,x_d)\in [0,1]^{d-1}\\
        (x_1+3,x_2,\cdots,x_d) & \text{ if }x_1\in [5,6]\text{ and }(x_2,\cdots,x_d)\in [0,1]^{d-1}\\
        0 & \text{ if }x_1\notin [0,2]\cup [5,6]\text{ or }(x_2,\cdots,x_d)\notin [0,1]^{d-1}.
    \end{cases}
\end{equation*}
It is easy to check that $s$ pushes $\mu$ forward to $\nu$ and $u_0(x)-u_0(s(x))=\|x-s(x)\|$ for $\mu$-a.e.\ $x\in \mathcal{X}$. Hence by \cite[Lemma 6]{MR1862796}, $s$ is an optimal transport map and $u_0$ is a Kantorovich potential. Now for any Kantorovich potential $u\in\mathrm{Lip}_1(\mathbb{R}^d)$, using \cite[Lemma 6]{MR1862796} again, we obtain that for any $(x_2,\cdots,x_d)\in [0,1]^{d-1}$,
\begin{align}
   & u(x_1+3,x_2,\cdots,x_d)-u(x_1,x_2,\cdots,x_d)=-3  \quad\text{for all }x_1\in [0,1],\label{display1}\\
&  u(x_1+6,x_2,\cdots,x_d)-u(x_1,x_2,\cdots,x_d)=-6 \quad\text{for all }x_1\in [1,2],\label{display2}\\
& u(x_1+3,x_2,\cdots,x_d)-u(x_1,x_2,\cdots,x_d)=-3\quad\text{for all }x_1\in [5,6]. \label{display3} 
\end{align}
Fix any $(x_2,\cdots,x_d)$ and let $\tilde{u}(x_1):=u(x_1,x_2,\cdots,x_d)$ for any $x_1\in [0,9]$. As $u\in\mathrm{Lip}_1(\mathbb{R}^d)$, $\tilde{u}$ is $1$-Lipschitz on $[0,9]$. From \eqref{display1}, for any $x_1\in [0,1]$, we have $\tilde{u}(x_1+3)-\tilde{u}(x_1)=-3$. Since $\tilde{u}$ is $1$-Lipschitz on $[0,9]$, it follows that $\tilde{u}$ is linear with slope $-1$ on $[x_1,x_1+3]$. As $x_1\in [0,1]$ is arbitrary, we conclude that $\tilde{u}$ is linear with slope $-1$ on $[0,4]$. Similarly, using \eqref{display2} and \eqref{display3}, we deduce that $\tilde{u}$ is linear with slope $-1$ on $[1,8]$ and $[5,9]$. Combining these facts, we conclude that $\tilde{u}$ is linear with slope $-1$ on $[0,9]$. Therefore, $u(0,x_2,\cdots,x_d)-u(9,x_2,\cdots,x_d)=\tilde{u}(0)-\tilde{u}(9)=9$, establishing \eqref{ray_identity}.     

By \eqref{ray_identity} and \Cref{TransportRays}, the family of transport rays is
\begin{equation*}
    \bigl\{\,T_{x_2,\ldots,x_d}\,\bigr\}_{(x_2,\ldots,x_d)\in[0,1]^{d-1}},
\qquad
T_{x_2,\ldots,x_d}:=\{(x_1,x_2,\ldots,x_d): x_1\in[0,9]\},
\end{equation*}
where $T_{x_2,\ldots,x_d}$ has upper end $(0,x_2,\ldots,x_d)$ and lower end $(9,x_2,\ldots,x_d)$. Now fix $(x_2,\cdots,x_d)\in [0,1]^{d-1}$ and set $T:=T_{x_2,\ldots,x_d}$. Disintegrating $\mu$ and $\nu$ along the transport rays, we obtain $\tilde{\mu}_T$ and $\tilde{\nu}_T$, which are the uniform measures on 
\[
\{(x_1,x_2,\ldots,x_d): x_1\in[0,2]\cup[5,6]\}
\quad\text{and}\quad
\{(x_1,x_2,\ldots,x_d): x_1\in[3,4]\cup[7,9]\},
\]
respectively. By \cite[Theorem 4.2]{Nutz.20}, the measures $\kappa_T'$ and $\tilde{\mu}_T\otimes \tilde{\nu}_T$ are equivalent, and so $\kappa_T'$ is not monotone on $T$. Therefore, $\gamma_0'=\lambda\otimes \kappa_T'$ is not an optimal transport plan.
\end{proof}

\appendix

\section{Proofs for Section~\ref{Sect.2.1}}\label{proof_Section2}

This section collects the proofs of the results in Section~\ref{Sect.2.1}.

\begin{proof}[Proof of \Cref{L2.0nn}]
Each transport ray $T\in\mathcal{S}$ is represented by its ends $(a_T,b_T)\in\mathbb{R}^d\times\mathbb{R}^d$, and then the topology on $\mathcal{S}$ (as introduced in Definition \ref{TransportRays}) is the same as the subspace topology on $\mathcal{S}$ induced from $\mathbb{R}^d\times\mathbb{R}^d$. As $a(\cdot),b(\cdot)$ are Borel measurable, $T:\mathcal{T}_1^{*}\rightarrow \mathcal{S}$ is also Borel measurable.
\end{proof}

\begin{proof}[Proof of \cref{Le2.2}]
As $\alpha(\cdot),\beta(\cdot)$ are Borel measurable and $u(\cdot)$ is continuous,
\begin{equation*}
    \EuScript{I}(\mathcal{S})= \big\{(x,y)\in\mathcal{X}\times\mathcal{Y}:x\neq y, u(x)-u(y)=\|x-y\|, \alpha(x)=0,\beta(x)>0,\alpha(y)>0,\beta(y)=0 \big\}
\end{equation*}
is a Borel subset of $\mathbb{R}^d\times\mathbb{R}^d$. As $\mathbb{R}^d\times\mathbb{R}^d$ is a Polish space, \cite[Corollary 13.4]{MR1321597} then yields that~$\mathcal{S}$ is a standard Borel space when equipped with its Borel $\sigma$-algebra.
\end{proof}

\begin{proof}[Proof of \Cref{L2.0}]

We denote $\sigma:=\min\{\alpha(x_1),\beta(x_1),\alpha(x_2),\beta(x_2)\}$. Let $T_1,T_2$ be the transport rays that contain $x_1,x_2$, respectively. We denote by $a_1,b_1$ the upper and lower ends of $T_1$, and by $a_2,b_2$ the upper and lower ends of $T_2$. 

First consider the case where $\|x_1-x_2\|\leq \sigma\slash 2$. As $|u(x_1)-u(x_2)|\leq \|x_1-x_2\|\leq\sigma\slash 2$, there exists $x_2'\in \mathrm{int}(T_2)$ such that $u(x_1)=u(x_2')$. Note that $\|x_2'-x_2\|=|u(x_1)-u(x_2)|\leq \|x_1-x_2\|$. Hence
\begin{align}\label{neqq2.1}
   &  \|x_2'-a_2\|\geq \|x_2-a_2\|-\|x_2-x_2'\|\geq \sigma-\|x_1-x_2\|\geq \sigma\slash 2,\nonumber\\
   &  \|x_2'-b_2\|\geq \|x_2-b_2\|-\|x_2-x_2'\|\geq \sigma-\|x_1-x_2\|\geq \sigma\slash 2.
\end{align}
By \cite[Lemma 16]{MR1862796} and \eqref{neqq2.1}, we have
\begin{align*}
    \|V(x_1)-V(x_2)\|&=\|V(x_1)-V(x_2')\|\leq \frac{\|x_1-x_2'\|}{\min\{\|x_1-a_1\|,\|x_1-b_1\|,\|x_2'-a_2\|,\|x_2'-b_2\|\}}\nonumber\\
    &\leq\frac{\|x_1-x_2'\|}{\sigma\slash 2}\leq \frac{\|x_1-x_2\|+\|x_2'-x_2\|}{\sigma\slash 2}\leq \frac{4\|x_1-x_2\|}{\sigma}.
\end{align*}

Now consider the case where $\|x_1-x_2\|>\sigma\slash 2$. We have
\begin{equation*}
    \|V(x_1)-V(x_2)\|\leq \|V(x_1)\|+\|V(x_2)\|=2\leq \frac{4\|x_1-x_2\|}{\sigma}, 
\end{equation*}
completing the proof.
\end{proof}

\begin{proof}[Proof of \Cref{L2.1n}]
As $\alpha(\cdot),\beta(\cdot)$ are upper semicontinuous, for any $n\in\mathbb{N}^{*}$, the sets $\{(q,t)\in H_{k,k'}\times\mathbb{R}:t\leq \alpha(q)-n^{-1}\}$ and $\{(q,t)\in H_{k,k'}\times\mathbb{R}:t\geq -\beta(q)+n^{-1}\}$ are closed in $H_{k,k'}\times\mathbb{R}$. Hence 
\begin{equation*}
    \{(q,t)\in H_{k,k'}\times\mathbb{R}:t< \alpha(q)\}=\bigcup_{n=1}^{\infty}\{(q,t)\in H_{k,k'}\times\mathbb{R}:t\leq \alpha(q)-n^{-1}\},
\end{equation*}
\begin{equation*}
    \{(q,t)\in H_{k,k'}\times\mathbb{R}:t> -\beta(q)\}=\bigcup_{n=1}^{\infty} \{(q,t)\in H_{k,k'}\times\mathbb{R}:t\geq -\beta(q)+n^{-1}\}
\end{equation*}
are Borel measurable. Consequently, noting that $\{q\in H_{k,k'}:\langle V(q),\mathtt{a}_{k,k'}\rangle\neq 0\}$ is a Borel subset of $H_{k,k'}$, we have  
\begin{align*}
    &\{(q,t)\in H_{k,k'}\times\mathbb{R}:q\in A\cap\mathcal{T}_1^{*},\langle V(q),\mathtt{a}_{k,k'}\rangle\neq 0,-\beta(q)<t<\alpha(q)\}\nonumber\\
    =\,& \big((\{q\in H_{k,k'}:\langle V(q),\mathtt{a}_{k,k'}\rangle\neq 0\}\cap A\cap\mathcal{T}_1^{*})\times\mathbb{R}\big)\nonumber\\
    &\cap\{(q,t)\in H_{k,k'}\times\mathbb{R}:t< \alpha(q)\}\cap\{(q,t)\in H_{k,k'}\times\mathbb{R}:t> -\beta(q)\}
\end{align*}
is Borel measurable. 

Now we define $\mathfrak{h}:\{(q,t)\in H_{k,k'}\times\mathbb{R}:q\in \mathcal{T}_1^{*},\langle V(q),\mathtt{a}_{k,k'}\rangle\neq 0,-\beta(q)<t<\alpha(q)\}\rightarrow \mathbb{R}^d$ by setting $\mathfrak{h}(q,t):=q+tV(q)\in\mathcal{T}_1^{*}$. By \cite[Lemma 10]{MR1862796}, $\mathfrak{h}$ is injective. For any $n\in\mathbb{N}^{*}$, we set 
\begin{equation*}
   \mathfrak{W}_{k,k';A;n}:=\big\{q+tV(q):q\in A\cap\mathcal{T}_1^*,\langle V(q),\mathtt{a}_{k,k'}\rangle\neq 0,-\beta(q)<t<\alpha(q),\min\{\alpha(q),\beta(q)\}\geq n^{-1}\big\},
\end{equation*}
\begin{align*}
   \tilde{\mathfrak{W}}_{k,k';A;n}
   :=\big\{(q,t)\in H_{k,k'}\times\mathbb{R}:\,& q\in A\cap\mathcal{T}_1^{*},\langle V(q),\mathtt{a}_{k,k'}\rangle\neq 0,\nonumber\\
   &-\beta(q)<t<\alpha(q),\min\{\alpha(q),\beta(q)\}\geq n^{-1}\big\}.
\end{align*}
From the previous paragraph, it follows that for any $n\in\mathbb{N}^{*}$, $\tilde{\mathfrak{W}}_{k,k';A;n}$ is a Borel subset of $H_{k,k'}\times\mathbb{R}$. Moreover, for any $(q,t),(q',t')\in \tilde{\mathfrak{W}}_{k,k';A;n}$, by Lemma~\ref{L2.0} and \eqref{bddsalphabeta}, we have $\|V(q)-V(q')\|\leq 4n\|q-q'\|$ and $\max\{|t|,|t'|\}\leq 2D$, hence
\begin{align*}
    \big\|(q+tV(q))-(q'+t'V(q'))\big\|&\leq \|q-q'\|+|t|\|V(q)-V(q')\|+|t-t'|\|V(q')\|\nonumber\\
    &\leq  (1+8nD)\|q-q'\|+|t-t'|.
\end{align*}
Hence $\mathfrak{h}|_{\tilde{\mathfrak{W}}_{k,k';A;n}}$ is Lipschitz continuous. By Kirszbraun's theorem, $\mathfrak{h}|_{\tilde{\mathfrak{W}}_{k,k';A;n}}$ can be extended to a Lipschitz continuous mapping $\tilde{\mathfrak{h}}: H_{k,k'}\times\mathbb{R}\rightarrow\mathbb{R}^d$. Note that $\tilde{\mathfrak{h}}|_{\tilde{\mathfrak{W}}_{k,k';A;n}}=\mathfrak{h}|_{\tilde{\mathfrak{W}}_{k,k';A;n}}$ is injective (as $\mathfrak{h}$ is injective). Consequently, using \cite[Section 2.2.10]{MR0257325} we conclude that $\mathfrak{W}_{k,k';A;n}=\mathfrak{h}(\tilde{\mathfrak{W}}_{k,k';A;n})=\tilde{\mathfrak{h}}(\tilde{\mathfrak{W}}_{k,k';A;n})$ is a Borel set. Therefore, $\mathfrak{W}_{k,k';A}=\bigcup_{n=1}^{\infty}\mathfrak{W}_{k,k';A;n}$ is a Borel set.
\end{proof}

\begin{proof}[Proof of \Cref{L2.4n}]
For any $a_T\in\mathcal{X}_k,b_T\in\mathcal{Y}_{k'}$, let $\psi(a_T,b_T)$ be the (unique) intersection point of $[a_T,b_T]$ and $H_{k,k'}$. As $\psi:\mathcal{X}_k\times\mathcal{Y}_{k'}\rightarrow H_{k,k'}$ is continuous, $\psi^{-1}(A)$ is a Borel subset of $\mathcal{X}_k\times \mathcal{Y}_{k'}$. Hence $\mathcal{S}_{k,k';A}=\{T\in \mathcal{S}: (a_T,b_T)\in \psi^{-1}(A)\}$ is a Borel subset of $\mathcal{S}$.
\end{proof}

\begin{proof}[Proof of \Cref{L2.3n}]
By the area formula (\cite[Theorem 3.9]{MR3409135}; note that the result holds for any nonnegative Lebesgue measurable function $g$, as shown in Case 1 of the proof of this result), for any nonnegative Lebesgue measurable function $\psi:\mathbb{R}^d\rightarrow [0,\infty)$, we have  
\begin{equation}\label{E2.3n.1.1}
    \int_{\mathbb{R}^d}\psi(x)J\EuScript{P}(x)dx=\int_{\mathbb{R}^d} \bigg(\sum_{x\in\EuScript{P}^{-1}(\{y\})}\psi(x)\bigg)dy. 
\end{equation}

Below, for any $x\in\mathbb{R}^d$, we take
\begin{equation}\label{defn.psi}
    \psi(x):=\varphi(\EuScript{P}(x))\mathbbm{1}_{\Omega_{k,k'}^{-1}(x_1,\cdots, x_{d-1})\in\tilde{H}_{k,k'}}\mathbbm{1}_{-\beta\big(\Omega_{k,k'}^{-1}(x_1,\cdots,x_{d-1})\big)<x_d<\alpha\big(\Omega_{k,k'}^{-1}(x_1,\cdots,x_{d-1})\big)}\geq 0.
\end{equation}
As $\EuScript{P}(\cdot)$ is Lipschitz continuous and $\tilde{H}_{k,k'}$ is a Borel subset of $H_{k,k'}$, $\psi(\cdot)$ is Lebesgue measurable. We have
\begin{align*}
 &\int_{\mathbb{R}^d}\psi(x)J\EuScript{P}(x)dx\nonumber\\
 =\,&    \int_{\mathbb{R}^d}\varphi(\EuScript{P}(x))\mathbbm{1}_{\Omega_{k,k'}^{-1}(x_1,\cdots,x_{d-1})\in\tilde{H}_{k,k'}}\mathbbm{1}_{-\beta\big(\Omega_{k,k'}^{-1}(x_1,\cdots,x_{d-1})\big)<x_d<\alpha\big(\Omega_{k,k'}^{-1}(x_1,\cdots,x_{d-1})\big)}J\EuScript{P}(x)dx\nonumber\\
  =\,& \int_{H_{k,k'}}d\mathcal{H}^{d-1}(q)\int_{\mathbb{R}}\varphi(\EuScript{P}(\Omega_{k,k'}(q),t))\mathbbm{1}_{q\in\tilde{H}_{k,k'}}\mathbbm{1}_{(-\beta(q),\alpha(q))}(t)J\EuScript{P}(\Omega_{k,k'}(q),t)dt\nonumber\\
  =\,& \int_{\tilde{H}_{k,k'}}d\mathcal{H}^{d-1}(q)\int_{\mathbb{R}}\mathbbm{1}_{(-\beta(q),\alpha(q))}(t)\varphi(\EuScript{P}(\Omega_{k,k'}(q),t))R_{k,k'}(q,t)dt,
\end{align*}
where we use \eqref{defRk} in the last line. Now note that for any $(q,t)\in\tilde{H}_{k,k'}\times\mathbb{R}$ such that $t\in (-\beta(q),\alpha(q))$, we have $\EuScript{P}(\Omega_{k,k'}(q),t)=\tilde{\Psi}_{k,k'}(q,t)=\Psi_{k,k'}(q,t)=q+t\tilde{V}_{k,k'}(q)=q+tV(q)$ (where we note that by \eqref{bddsalphabeta}, $|t|\leq 2D$). Hence
\begin{equation}\label{E2.3n.1.2}
    \int_{\mathbb{R}^d}\psi(x)J\EuScript{P}(x)dx=\int_{\tilde{H}_{k,k'}}d\mathcal{H}^{d-1}(q)\int_{\mathbb{R}}\mathbbm{1}_{(-\beta(q),\alpha(q))}(t)\varphi(q+tV(q))R_{k,k'}(q,t)dt.
\end{equation}

Now consider any $y\in\mathbb{R}^d$ for which there exists some $x=(x_1,\cdots,x_d)\in\EuScript{P}^{-1}(\{y\})$ such that $\psi(x)\neq 0$. Setting $q=\Omega_{k,k'}^{-1}(x_1,\cdots,x_{d-1})$ and $t=x_d$, by \eqref{defn.psi}, we have $q\in \tilde{H}_{k,k'}$ and $t\in (-\beta(q),\alpha(q))$ (note that by \eqref{bddsalphabeta}, $|t|\leq 2D$). Thus $y=\EuScript{P}(x)=\tilde{\Psi}_{k,k'}(q,t)=\Psi_{k,k'}(q,t)=q+t\tilde{V}_{k,k'}(q)=q+tV(q)$, which implies that $y\in \tilde{\mathcal{T}}^{*}_{1;k,k'}$. Moreover, $\psi(x)=\varphi(\EuScript{P}(x))=\varphi(y)$. 

From the previous paragraph, for any $y\in\mathbb{R}^d\backslash \tilde{\mathcal{T}}^{*}_{1;k,k'}$, we have $\sum_{x\in\EuScript{P}^{-1}(\{y\})}\psi(x)=0$. Now we consider any $y\in \tilde{\mathcal{T}}^{*}_{1;k,k'}$. Note that there exist $q\in \tilde{H}_{k,k'}$ and $t\in (-\beta(q),\alpha(q))$ such that $y=q+tV(q)$ (note that by \eqref{bddsalphabeta}, $|t|\leq 2D$). For any $x=(x_1,\cdots,x_d)\in\EuScript{P}^{-1}(\{y\})$ such that $\psi(x)\neq 0$, setting $q'=\Omega_{k,k'}^{-1}(x_1,\cdots,x_{d-1})$ and $t'= x_d$, we have (from the previous paragraph) $q'\in\tilde{H}_{k,k'}$, $t'\in (-\beta(q'),\alpha(q'))$, and $y=q'+t'V(q')$. By \cite[Lemma 10]{MR1862796}, this implies that $q'=q$ and $t'=t$, hence $x=(\Omega_{k,k'}(q),t)$. As $\EuScript{P}(\Omega_{k,k'}(q),t)=\tilde{\Psi}_{k,k'}(q,t)=\Psi_{k,k'}(q,t)=q+t\tilde{V}_{k,k'}(q)=q+tV(q)=y$, we have $\sum_{x\in\EuScript{P}^{-1}(\{y\})}\psi(x)=\psi(\Omega_{k,k'}(q),t)=\varphi(y)$. Therefore,
\begin{equation}\label{E2.3n.1.3}
     \int_{\mathbb{R}^d} \bigg(\sum_{x\in\EuScript{P}^{-1}(\{y\})}\psi(x)\bigg)dy=\int_{\tilde{\mathcal{T}}^{*}_{1;k,k'}}\varphi(y)dy.
\end{equation}
The desired conclusion \eqref{E2.3n.1} follows from \eqref{E2.3n.1.1}, \eqref{E2.3n.1.2}, and \eqref{E2.3n.1.3}.
\end{proof}

\begin{proof}[Proof of \Cref{L2.5n}]
Consider any Borel set $A\subseteq \tilde{H}_{k,k'}$. By Lemma \ref{L2.1n}, the set $\mathfrak{W}_{k,k';A}$ is Borel measurable. By the definition of $\tilde{H}_{k,k'}$ as in \eqref{deftildehat}, we have $\mathfrak{W}_{k,k';A}=\big\{q+tV(q):q\in A, -\beta(q)<t<\alpha(q)\big\}$. Note that $\mathfrak{W}_{k,k';A}$ is a transport set (recall Definition \ref{Def2.7}). Hence by \cite[Lemma 27]{MR1862796}, we have
\begin{equation*}
    \int_{\mathfrak{W}_{k,k';A}}f(x)dx=\int_{\mathfrak{W}_{k,k';A}}g(x)dx\leq 1.
\end{equation*}
By Lemma \ref{L2.3n}, we have
\begin{align*}
&\int_{\mathfrak{W}_{k,k';A}}f(x)dx=\int_{\tilde{\mathcal{T}}^{*}_{1;k,k'}}\mathbbm{1}_{\mathfrak{W}_{k,k';A}}(x)f(x)dx \nonumber\\
=\,& \int_{\tilde{H}_{k,k'}}d\mathcal{H}^{d-1}(q) \int_{\mathbb{R}} \mathbbm{1}_{(-\beta(q),\alpha(q))}(t) 
    \mathbbm{1}_{\mathfrak{W}_{k,k';A}}(q+tV(q))f(q+tV(q)) R_{k,k'}(q,t)dt\nonumber\\
    =\,&\int_A d\mathcal{H}^{d-1}(q)\int_{\mathbb{R}}\mathbbm{1}_{(-\beta(q),\alpha(q))}(t) 
      f(q+tV(q)) R_{k,k'}(q,t)dt.
\end{align*}
Similarly, 
\begin{equation*}
    \int_{\mathfrak{W}_{k,k';A}}g(x)dx=\int_A d\mathcal{H}^{d-1}(q)\int_{\mathbb{R}}\mathbbm{1}_{(-\beta(q),\alpha(q))}(t) 
      g(q+tV(q)) R_{k,k'}(q,t)dt.
\end{equation*}
Combining the above three displays, we obtain that for any Borel set $A\subseteq\tilde{H}_{k,k'}$,
\begin{align*}
   &  \int_A d\mathcal{H}^{d-1}(q)\int_{\mathbb{R}}\mathbbm{1}_{(-\beta(q),\alpha(q))}(t) 
      f(q+tV(q)) R_{k,k'}(q,t)dt\nonumber\\
      =\,&  \int_A d\mathcal{H}^{d-1}(q)\int_{\mathbb{R}}\mathbbm{1}_{(-\beta(q),\alpha(q))}(t) 
      g(q+tV(q)) R_{k,k'}(q,t)dt\leq 1.
\end{align*}
Hence for $\mathcal{H}^{d-1}$-a.e.\ $q\in \tilde{H}_{k,k'}$, \eqref{L2.5.e1} holds. By Lemma \ref{L2.3n}, this implies that for $\mathcal{L}^d$-a.e.\ $z_0\in \tilde{\mathcal{T}}_{1;k,k'}^*$, if we set $q$ to be the intersection point of $T(z_0)$ and $H_{k,k'}$, then \eqref{L2.5.e1} holds.   
\end{proof}

\begin{proof}[Proof of \Cref{Lemma2.7}]

We only present the proof of $\big|\int_{\mathcal{T}_1^{*}} f(x)\log(\alpha(x))dx\big|\leq C$. The proofs for the other results are similar. First note that from \eqref{bddsalphabeta} and the fact that $0\leq \int_{\mathcal{T}_1^{*}} f(x)dx \leq 1$, we have
\begin{equation}\label{L2.7n0}
    \int_{\mathcal{T}_1^{*}} f(x)\log(\alpha(x))dx\leq \log(2D)\int_{\mathcal{T}_1^{*}} f(x)dx \leq \log(2D)\leq C. 
\end{equation}
Moreover, 
\begin{equation}\label{L2.7n1}
    \int_{\mathcal{T}_1^{*}} f(x)\log(\alpha(x))dx=\log(2D)\int_{\mathcal{T}_1^{*}} f(x)dx+\int_{\mathcal{T}_1^{*}} f(x)\log\Big(\frac{\alpha(x)}{2D}\Big)dx\geq \int_{\mathcal{T}_1^{*}} f(x)\log\Big(\frac{\alpha(x)}{2D}\Big)dx.
\end{equation}
By \eqref{disj}, we have 
\begin{equation}\label{L2.7n2}
    \int_{\mathcal{T}_1^{*}} f(x)\log\Big(\frac{\alpha(x)}{2D}\Big)dx=\sum_{k=1}^K\sum_{k'=1}^{K'}\int_{\mathcal{T}_{1;k,k'}^{*}} f(x)\log\Big(\frac{\alpha(x)}{2D}\Big)dx. 
\end{equation}
Now for any $k\in[K], k'\in[K']$, taking 
\begin{equation*}
    \varphi(x)=\begin{cases}
    -f(x)\log\big(\frac{\alpha(x)}{2D}\big) & \text{ for every }x\in \mathcal{T}_{1;k,k'}^{*}, \\
    0 & \text{ for every }x \in \mathbb{R}^d\backslash \mathcal{T}_{1;k,k'}^{*}
    \end{cases}
\end{equation*}
in Lemma \ref{L2.3n} (note that $\varphi(x)\geq 0$ by \eqref{bddsalphabeta}), we get
\begin{align}\label{L2.7n3}
   &\int_{\mathcal{T}_{1;k,k'}^{*}} f(x)\log\Big(\frac{\alpha(x)}{2D}\Big)dx\nonumber\\
   =\,& \int_{\tilde{H}_{k,k'}}d\mathcal{H}^{d-1}(q)\int_{\mathbb{R}}\mathbbm{1}_{(-\beta(q),\alpha(q))}(t)f(q+tV(q))\log\bigg(\frac{\alpha(q+tV(q))}{2D}\bigg)\mathbbm{1}_{q+tV(q)\in\mathcal{T}_{1;k,k'}^{*}}R_{k,k'}(q,t)dt\nonumber\\
   =\,& \int_{\bar{H}_{k,k'}}d\mathcal{H}^{d-1}(q)\int_{\mathbb{R}}\mathbbm{1}_{(-\beta(q),\alpha(q))}(t)f(q+tV(q))\log\bigg(\frac{\alpha(q)-t}{2D}\bigg)R_{k,k'}(q,t)dt\nonumber\\
   \geq\,& -C\int_{\bar{H}_{k,k'}}d\mathcal{H}^{d-1}(q)\int_{\mathbb{R}}\mathbbm{1}_{(-\beta(q),\alpha(q))}(t)\bigg|\log\bigg(\frac{\alpha(q)-t}{2D}\bigg)\bigg|dt\geq -C, 
\end{align}
where we use \eqref{subs} and \eqref{deftildehat} in the first two equalities, \eqref{bddLip} and \eqref{Rkkb} in the first inequality, and \eqref{BddHtildekk} in the last inequality. By \eqref{L2.7n1}--\eqref{L2.7n3}, noting that $K,K'$ are of constant order, we get $\int_{\mathcal{T}_1^{*}} f(x)\log(\alpha(x))dx\geq -C$. Combining this with \eqref{L2.7n0} yields $\big|\int_{\mathcal{T}_1^{*}} f(x)\log(\alpha(x))dx\big|\leq C$.  
\end{proof}

\begin{proof}[Proof of \Cref{Lem4.4}]
By \eqref{disj}, for any $s\geq 0$, 
\begin{equation*}
    \mathcal{L}^d(\{x\in\mathcal{T}_1^{*}:\alpha(x)\leq s\})=\int_{\mathcal{T}_1^{*}}\mathbbm{1}_{\alpha(x)\leq s}dx=\sum_{k=1}^K\sum_{k'=1}^{K'} \int_{\mathcal{T}_{1;k,k'}^{*}}\mathbbm{1}_{\alpha(x)\leq s}dx.
\end{equation*}
By Lemma \ref{L2.3n}, for each $k\in[K],k'\in[K']$, we have
\begin{align*}
    \int_{\mathcal{T}_{1;k,k'}^{*}}\mathbbm{1}_{\alpha(x)\leq s}dx&=\int_{\bar{H}_{k,k'}}d\mathcal{H}^{d-1}(q)\int_{\mathbb{R}}\mathbbm{1}_{(-\beta(q),\alpha(q))}(t)\mathbbm{1}_{\alpha(q+tV(q))\leq s}R_{k,k'}(q,t)dt\nonumber\\
    &\leq Cs\mathcal{H}^{d-1}(\bar{H}_{k,k'})\leq Cs,
\end{align*}
where we use \eqref{Rkkb} and \eqref{BddHtildekk} in the last two inequalities. Combining the above two displays, noting that $K,K'$ are of constant order, we get $\mathcal{L}^d(\{x\in\mathcal{T}_1^{*}:\alpha(x)\leq s\})\leq Cs$. The second conclusion $\mathcal{L}^d(\{x\in\mathcal{T}_1^{*}:\beta(x)\leq s\})\leq Cs$ can be similarly deduced. 
\end{proof}

\bibliographystyle{acm}
\bibliography{OptimalTransport.bib}

\end{document}